\input amstex
\documentstyle{amsppt}
\magnification1200
\tolerance=10000
\overfullrule=0pt
\def\n#1{\Bbb #1}
\def\p{\Bbb C_{\infty}}
\def\fr{\hbox{fr}}

\def\im{\hbox{im }}
\def\invlim{\hbox{invlim}}
\def\tr{\hbox{tr }}

\def\Gal{\hbox{Gal }}

\def\Exp{\hbox{Exp}}

\def\Hom{\hbox{Hom}}

\def\End{\hbox{End}}
\def\Prin{\hbox{Prin}}
\def\Ker{\hbox{Ker }}
\def\Lie{\hbox{Lie}}

\def\Div{\hbox{Div}}
\def\Pic{\hbox{Pic}}

\def\ord{\hbox{ord}}

\def\Id{\hbox{Id}}
\def\Cl{\hbox{Cl}}
\def\Supp{\hbox{ Supp }}
\def\Spec{\hbox{ Spec }}

\def\diag{\hbox{ diag }}
\def\Diag{\hbox{ Diag }}

\def\e11{I_{11}}

\def\ga{\goth A}
\def\w{\hbox{\bf A}}
\def\x{\hbox{\bf K}}
\def\ve{\varepsilon}
\def\vf{\varphi}

\def\ve{\varepsilon}

\def\vf{\varphi}

\def\de{\delta}

\def\ga{\gamma}

\def\be{\beta}

\def\al{\alpha}

\def\om{\omega}
\def\g{\goth }

\topmatter
\title
Duality of Anderson $t$-motives
\endtitle
\author
A. Grishkov, D. Logachev\footnotemark \footnotetext{E-mails: shuragri{\@}gmail.com; logachev94{\@}gmail.com (corresponding author)\phantom{*******************}}
\endauthor
\thanks Thanks: The authors are grateful to FAPESP, S\~ao Paulo, Brazil for a financial support (process No. 2017/19777-6). The first author is grateful to SNPq, Brazil, to RFBR, Russia, grant 16-01-00577a (Secs. 1-4), and to Russian Science Foundation, project 16-11-10002 (Secs. 5-8) for a financial support. The second author is grateful to Gilles Lachaud for invitation to IML and to Laurent Lafforgue for invitation to
IHES where this paper was started, and to Vladimir
Drinfeld for invitation to the University of Chicago where
this paper was continued. Discussions with Greg Anderson, Vladimir
Drinfeld, Laurent Fargues, Alain Genestier, David Goss, Richard Pink, Yuichiro
Taguchi, Dinesh Thakur on
the subject of this paper were very important. Particularly, Alain
Genestier informed me about the paper of Taguchi where the notion of
the dual of a Drinfeld module is defined. Further, Richard Pink indicated me an important reference (see Section 6 for details); proof of the main theorem of the present paper grew from it. Finally, Vladimir
Drinfeld indicated me the proof of the Theorem 12.6 and Jorge Morales
gave me a reference on classification of quadratic forms over $\n
F_q[T]$ (Remark 7.8).
\endthanks
\NoRunningHeads
\address
First author: Departamento de Matem\'atica e estatistica
Universidade de S\~ao Paulo. Rua de Mat\~ao 1010, CEP 05508-090, S\~ao Paulo, Brasil, and Omsk State University n.a. F.M.Dostoevskii. Pr. Mira 55-A, Omsk 644077, Russia.
\medskip
Second author: Departamento de Matem\'atica, Universidade Federal do Amazonas, Manaus, Brasil
\endaddress

\abstract Let $M$ be a t-motive. We introduce the notion of duality for $M$. Main results of the paper (we consider uniformizable $M$ over $\n F_q[T]$ of rank $r$, dimension $n$, whose nilpotent operator $N$ is 0):

1. Algebraic duality implies analytic duality (Theorem 5). Explicitly, this means that the lattice of the dual of $M$ is the dual of the lattice of $M$, i.e. the transposed of a Siegel matrix of $M$ is a Siegel matrix of the dual of $M$.

2. Let $n=r-1$. There is a 1 -- 1 correspondence between pure t-motives (all they are uniformizable), and lattices of rank $r$ in $\p^{n}$ having dual (Corollary 8.4).

3. Pure t-motives have duals which are pure t-motives as well (Theorem 10.3).

4. Some explicit results are proved for $M$ having complete
multiplication. The CM-type of the dual of $M$ is the complement of
the CM-type of $M$. Moreover, for $M$ having multiplication by a
division algebra there exists a simple formula for the
CM-type of the dual of $M$ (Section 12).

5. We construct a class of non-pure t-motives (t-motives having the
completely non-pure row echelon form) for which duals
are explicitly calculated (Theorem 11.5). This is the first step
of the problem of description of all t-motives having duals.

6. If $M$ has good ordinary reduction then the kernels of reduction maps on groups of  torsion points for $M$ and its dual are complementary
with respect to a natural pairing (proof is given for a particular case, Conjecture 13.4.1).
\endabstract
\endtopmatter
\document
{\bf 0. Introduction.}
\nopagebreak
\medskip
t-motives are the function field analogs of abelian varieties (more exactly, of abelian varieties with multiplication by an imaginary quadratic field, see [L09]). Main references for t-motives are [A], [G]. Nevertheless, function field analogs of some basic results in the theory of abelian varieties are not known yet.

The present paper contains an analog of such result. Namely, we introduce the notion of duality for a t-motive $M$ (this is not the duality in
a Tannakian category!), and we prove some properties of this notion, see the abstract. Particularly, if $M$ is uniformizable and has dual then
the lattice of the dual of $M$ is the dual of the lattice of $M$ (Theorem 5)\footnotemark \footnotetext{Here this result is proved for $M$ having
the associated nilpotent operator $N$ (see (1.9.2))
equal to 0. The same result for $M$ having $N\ne0$ is proved in [GL18].}. An immediate corollary of the above theorem and the result of
Drinfeld on 1 -- 1 correspondence between Drinfeld modules and lattices in $\p$ (here $\p$ is the function field analog of $\n C$) is Corollary 8.4: there is a 1 -- 1 correspondence between pure t-motives of dimension $r-1$ and rank $r$, and lattices of rank $r$ in $\p^{r-1}$ having dual (not all such lattices have dual).

Let us give more details on the contents of the paper. For simplicity, most results are proved for t-motives over the ring $\hbox{\bf A}=\n F_q[T]$,
and we consider, with few exceptions, only the case $N=0$. The main definition of duality of t-motives
(definition 1.8 --- case $\hbox{\bf A}=\n F_q[T]$ and definition 1.13 --- general case) is given in
Section 1.\footnotemark \footnotetext{A version of the definition of duality is obtained independently in [Tae], 2.2.}
Lemma 1.10 gives the explicit matrix form of the definition of duality of t-motives. Since Taguchi in [T] gave a definition of dual to
a Drinfeld module, we prove in Proposition 1.12.3 that the definition of the present paper is equivalent to the original definition of Taguchi.
Section 1.14 contains a definition of duality for abelian $\tau$-sheaves ([BH], Definition 2.1), but we do not develop this subject.

Section 2 contains the definition of the dual lattice. Section 3 contains explicit formulas for the dual lattice.
Section 5 contains the statement and the proof of the main theorem 5 --- coincidence of algebraic and analytic duality for the case
$\hbox{\bf A}=\n F_q[T]$ (section 4 contains the statement of the corresponding conjecture for the case of general $\hbox{\bf A}$).
Section 6 contains the theorem 6 describing the lattice of the tensor product of two t-motives (case $N=0$; the proof for the general
case was obtained, but not published, by Anderson). Section 7 contains the notion of self-dual t-motives and polarization form on them.
Some examples are given. We discuss in Section 8 the problem of correspondence between uniformizable t-motives and lattices.
Section 9 gives the statement of the main result for the case $N\ne 0$ without proof and a reformulation of the theorem 5 in terms
of Hodge-Pink structures of constant weight.

Further on, we prove in Section 10 that pure t-motives have duals which are pure
t-motives as well, and some related results (a proof that the dual of an abelian
$\tau$-sheaf is also an abelian $\tau$-sheaf can be obtained using ideas
of Section 10). In Section 11 we consider t-motives having the
completely non-pure row echelon form, and we give an explicit formula
for their duals. In Section 12 we consider t-motives with complete multiplication, and we give for them a very simple version of the proof of the first part of the main theorem.
Section 13 contains some explicit formulas for t-motives of complete multiplication. In 13.1 we describe the dual lattice, in 13.2 we show that the results of Section 12 are compatible with (the first form of) the main theorem of complete multiplication. Section 13.3 contains an explicit proof of the main theorem for t-motives with complete multiplication by two types of simplest fields.
Section 13.4 gives us an application of the notion of duality to the
reduction of t-motives (subject in development, see [L]).
\medskip
{\bf Notations. }
\medskip
$q$ is a power of a prime $p$;
\medskip
{\it Case of $M$ over $\n F_q[T]$}:
\medskip
$\n Z_\infty:=\n F_q[\theta]$, $\n R_{\infty}:=\n F_q((1/\theta))$, $\p$ is the completion of its algebraic closure
($\n Z_\infty$, $\n R_{\infty}$, $\p$ are the function field analogs of $\n Z$, $\n R$, $\n C$ respectively);

$\w:=\n F_q[T]$, $\hbox{\bf K}:=\n F_q((1/T))$;

$\iota: \hbox{\bf A} \to \p$ ($\iota(T)=\theta$) is the standard map of generic characteristic (with one exception (1.16), we shall not consider the case of finite characteristic);
we extend $\iota$ to $\hbox{\bf K}$, and we have $\n Z_\infty=\iota(\w)\subset \p$, $\n R_\infty=\iota(\hbox{\bf K})\subset \p$.

$\goth C$ (resp. $\goth C_2$) is the Carlitz module over $\w=\n F_q[T]$ (resp. over $\n F_{q^2}[T]$).
\medskip
{\it Case of $M$ over an extension of $\n F_q[T]$}:
\medskip
$\n Q_\infty$ is a finite separable extension of $\n F_q(\theta)$;

$\infty$ is a fixed valuation of $\n Q_\infty$ over the infinity of $\n F_q(\theta)$;

${\n Z_\infty \subset \n Q_\infty}$ is the subring of elements which are regular outside $\infty$;

$\n R_\infty$ is the completion of $\n Q_\infty$ at infinity, and $\p$ --- the completion of its algebraic closure --- is the same as of the case of $M$ over $\n F_q[T]$.

$\w\supset \n F_q[T]$, $\hbox{\bf K}\supset \n F_q((1/T))$ are defined by the condition that $\iota: \w \to \n Z_\infty$, $\iota: \hbox{\bf K} \to \n R_\infty$ are isomorphisms.

$\w_C:=\w\underset{\n F_q}\to{\otimes} \p$ (i.e. $\w_C=\p[T]$ for the case of $M$ over $\n F_q[T]$).

$\goth C$ is a Drinfeld module of rank 1 over {\bf A}.
\medskip
If $P=\frac{\sum a_iT^i}{\sum b_iT^i}\in \p(T)$ then $P^{(k)}:=\frac{\sum a_i^{q^k}T^i}{\sum b_i^{q^k}T^i}$. For $x\in \w_C$, $x=a\otimes z$, $a\in \w$, $z\in \p$ we let $x^{(k)}:= a\otimes z^{q^k}$.

$M_r$ is the set of $r\times r$ matrices. If $C=\{c_{ij}\}$ is a matrix with entries $c_{ij}\in \p(T)$ then $C^{(k)}:=\{c_{ij}^{(k)}\}$, $C^t$ is the  transposed of $C$, $C^{(k)\ -1}=(C^{(k)})^{-1}$, $C^{t-1}=(C^t)^{-1}$.

If $M$ is an $\w_C$-module, we define $M^{(1)}$ as the tensor
product $M\otimes_{\w_C,*^{(1)}}\w_C$ with respect to the map $*^{(1)}: \w_C\to \w_C$ (this notation is concordant in the obvious sense with
the above notation $C^{(1)}$).

For a t-motive $M$ we denote by $E=E(M)$ the corresponding t-module (see [G], Theorem 5.4.11; Goss uses the inverse functor $E\mapsto M=M(E)$).

$\Lie(M)$ is $\Lie(E(M))$ ([G], 5.4).

$I_k$ is the unit matrix of size $k$.

Throughout the whole paper the word "canonical" will mean "canonical up to multiplication by elements of $\n F_q^*$".
\medskip
{\bf 1. Definitions. }

\nopagebreak
\medskip
If otherwise is not explicitly stated, throughout the whole paper we consider the case of t-motives $M$ over the ring $\hbox{\bf A}=\n F_q[T]$ such that $N=N(M)=0$. Exceptions: case of arbitrary $\w$ is treated in Sections 1.13, 1.14, 2, 4, 5.2. Case of arbitrary $N$ is treated in Sections 1, 10 and in statements of some results of Anderson in Sections 5, 6.
\medskip
In the present section we consider $M$ such that $N(M)$ is arbitary.
\medskip
Let $\p[T,\tau]$ be the Anderson ring, i.e. the ring of non-commutative polynomials satisfying the following relations (here
$a \in \p$):
$$Ta=aT, \ T\tau = \tau T, \ \tau a = a^q \tau \eqno{(1.1)}$$
We need also an extension of $\p[T,\tau]$ --- the ring $\p(T)[\tau]$
which is the ring of
non-commutative polynomials in $\tau$ over the field of rational
functions $\p(T)$ with
the same relations (1.1). For a left $\p[T,\tau]$-module $M$ we denote by $M_{\p[T]}$ the
same
$M$ treated as a
$\p[T]$-module with respect to the natural inclusion
$\p[T]\hookrightarrow\p[T,\tau]$.
Analogously, we define $M_{\p[\tau]}$; we shall use similar notations
also for the left
$\p(T)[\tau]$-modules.
\medskip
Obviously we have:
\medskip
{\bf (1.2)} For $C\in M_{r}(\p(T))$ operations $C^t$, $C^{-1}$
and $C^{(i)}$ commute.
\medskip
{\bf Definition 1.3.} ([G], 5.4.2, 5.4.12, 5.4.10). A t-motive $M$ is a left
$\p[T, \tau]$-module which is free and finitely generated as both $\p[T]$-,
$\p[\tau]$-module and such that
$$ \exists m=m(M) \ \hbox{such that}\ (T-\theta)^m M/\tau M=0\eqno{(1.3.1)}$$
\medskip
{\bf Remark.} The above object is called "abelian t-motive" (resp. "t-motive") in [G] (resp. [A]), while the name "t-motive" is used in [G] for a more general object ([G], Definition 5.4.2). Since we shall not use objects defined in [G], 5.4.2, I prefer to use a shorter name for the above $M$.
\medskip
t-motives are main objects of the present paper. If we affirm
that an object exists this means that it exists as a t-motive if otherwise
is not stated. We denote dimension of $M$ over $\p[\tau]$ (resp. $\p[T]$) by $n$ (resp.
$r$), these numbers are called dimension and rank of $M$. Morphisms of abelian
t-motives are morphisms of left $\p[T, \tau]$-modules.

To define a left $\p[T, \tau]$-module $M$ is the same as to define a left $\p[T]$-module
$M_{\p[T]}$
endowed by an action of $\tau$ satisfying $\tau(Pm)=P^{(1)}\tau(m)$, $P\in \p[T]$. In
this situation
we can also treate $\tau$ as a $\p[T]$-linear map $M^{(1)}\to M$. This interpretation is
necessary if
we consider the general case $\w\supset \n F_q[T]$.

We need two categories which are larger than the category of abelian
t-motives.
\medskip
{\bf Definition 1.4.} A pr\'e-t-motive is a left $\p[T,
\tau]$-module which is free
and finitely generated as $\p[T]$-module, and satisfies (1.3.1).
\medskip
{\bf Definition 1.5.} A rational pr\'e-t-motive is a left
$\p(T)[\tau]$-module which is
free and finitely generated as $\p(T)$-module.
\medskip
{\bf Remark 1.6.} An analog of (1.3.1) does not exist for them.
\medskip
There is an obvious functor from the category of t-motives to
the category of
pr\'e-t-motives which is fully faithful, and an obvious functor from the category of
pr\'e-t-motives to the
category of rational pr\'e-t-motives. We denote these functors by
$i_1$, $i_2$
respectively. It is easy to see (Remark 10.2.3) that if $M$ is a
pr\'e-t-motive then the
action of $\tau$ on $i_2(M)$ is invertible.

Let $M_1$, $M_2$ be rational pr\'e-t-motives such
that the action of
$\tau$ on $(M_1)_{\p(T)}$ is invertible.
\medskip
{\bf Definition 1.7. (1)} $\Hom(M_1,M_2)$ is a rational pr\'e-t-motive such that

$$\Hom(M_1,M_2)_{\p(T)}=\Hom_{\p(T)}((M_1)_{\p(T)}, (M_2)_{\p(T)})$$
and the action of $\tau$ is defined by the usual manner: for
$\varphi:M_1 \to M_2$, $m\in M_1$
$$(\tau\varphi)(m)=\tau(\varphi(\tau^{-1}(m)))$$

{\bf (2)} Let $M_1$, $M_2$ be t-motives. Their tensor product is defined by $M_1\otimes_{\p[T]}M_2$ where the action of $\tau$ is given by
$\tau(m_1\otimes m_2)=\tau(m_1)\otimes \tau(m_2)$. It is known (Anderson; see also [G]) that $M_1\otimes M_2$ is really a t-motive of rank
$r_1r_2$, of dimension $n_1r_2+n_2r_1$. $M_1\otimes M_2$ has $N\ne0$ even if $M_1$, $M_2$ have $N=0$.
\medskip
The Carlitz module $\goth C$ is the Anderson t-motive with $r=n=1$, it is unique over $\p$ (see, for example, [G], 3.3). The $\mu$-th tensor power
of $\goth C$ is denoted
by $\goth C^{\otimes \mu}$. Its rank $r$ is 1 and its dimension is $\mu$.
\medskip
{\bf Definition 1.8.} Let $M$ be a t-motive and
$\mu$ a positive
number. A t-motive
$M'={M'}^{\mu}$ is called the $\mu$-dual of $M$ (dual if $\mu=1$) if
$M'=\Hom(M,\goth
C^{\otimes \mu})$ as a rational pr\'e-t-motive, i.e. $$i_2\circ
i_1(M')=\Hom(i_2\circ i_1(M),\goth
C^{\otimes
\mu})\eqno{(1.8.1)}$$

{\bf Remark.} This definition generalizes the original one of Taguchi
([T], Section 5), see 1.12 below. A similar definition is in [F].
\medskip
{\bf 1.9.} We shall need the explicit matrix description of the above
objects. Let
$e_*=(e_1, ..., e_n)^t$ be the vector column of elements of a basis
of $M$ over $\p[\tau]$. There exists a matrix $\goth A\in M_n(\p[\tau])$ such that

$$T e_* = \goth A e_*, \ \ \goth A = \sum_{i=0}^l \goth A_i \tau^i \hbox{ where } \goth A_i
\in M_n(\p)\eqno{(1.9.1)}$$
Condition (1.3.1) is equivalent to the condition
$$\goth A_0=\theta I_n + N\eqno{(1.9.2)}$$
where $N$ is a
nilpotent matrix, and the
condition
$m(M)=1$ is equivalent to the condition $N=0$.
\medskip
Let $f_*=(f_1, ..., f_r)^t$ be the vector column of elements of a basis
of $M$ over
$\p[T]$. There exists a matrix $Q=Q(f_*)\in M_r(\p[T])$ such that
$$\tau f_* = Q f_*\eqno{(1.9.3)}$$

{\bf Lemma 1.10.} Let $M$ be as above. A t-motive
$M'$ is the $\mu$-dual of $M$ iff there exists a basis
$f'_*=(f'_1, ..., f'_r)^t$ of $M'$ over $\p[T]$ such that its matrix
$Q'=Q(f'_*)$ satisfies
$$Q'=(T-\theta)^{\mu}Q^{t-1} \eqno{(1.10.1)}$$

{\bf Proof.} The matrix $Q$ of $\goth C^{\otimes \mu}$ is $(T-\theta)^{\mu}$. This implies the formula. $\square$
\medskip
{\bf 1.10.2.} For further applications we shall need the following
lemma. The above
$f_*$, $f'_*$ are the dual bases (i.e. if we consider $f'_i$ as
elements of $\Hom(M,\goth
C)$ then $f'_i(f_j)=\delta^i_j\goth f$, where $\goth f$ is canonically defined by the condition that it generates $\goth C_{\p[T]}$ and satisfies $\tau \goth
f=(T-\theta)\goth f$). Let $\gamma$ be an endomorphism of $M$ and $D$ its
matrix in the basis
$f_*$ (i.e. $\gamma(f_*)=Df_*$). Let $\gamma'$ be the dual endomorphism.
\medskip
{\bf Lemma 1.10.3.} The matrix of $\gamma'$ in the basis $f'_*$ is $D^t$.
$\square$
\medskip
{\bf Remark 1.11.1.} For any $M$ having dual there exists a canonical homomorphism $\delta: \goth
C \to M\otimes M'$.
This is a well-known theorem of linear algebra. Really, in the above notations we have
$\goth f \mapsto \sum_i f_i \times f'_i$. It is obvious that $\delta$ is well-defined,
canonical and compatible
with the action of $\tau$.
\medskip
{\bf Remark 1.11.2.} The $\mu$-dual of $M$ --- if it exists --- is
unique, i.e. does not
depend on base change. This follows immediately from Definition 1.8,
but can be deduced
easily from 1.10.1. Really, let $g_*=(g_1, ..., g_r)^t$ be another basis
of $M$ over
$\p[T]$ and $C\in GL_r(\p[T])$ the matrix of base change (i.e. $g_*= C
f_*$). Then
$Q(g_*)=C^{(1)}QC^{-1}$. Let $g'_*=(g'_1, ..., g'_r)^t$ be a basis of
$M'$ over $\p[T]$
satisfying $g'_*= C^{t-1} f'_*$. Elementary calculation shows that
matrices $Q(g_*)$,
$Q(g'_*)$ satisfy (1.10.1).
\medskip
{\bf Remark 1.11.3.} The operation $M \mapsto {M'}^{\mu}$ is
obviously contravariant functorial. It is an exercise to the reader to give an exact definition of the corresponding category such that the functor of duality is defined on it, and is involutive (recall that not all t-motives have duals, and the dual of a map of t-motives is a priori a map of rational pr\'e-t-motives).
\medskip
{\bf 1.12.} The original definition of duality ([T], Definition 4.1; Theorem 5.1) from
the first
sight seems to be more
restrictive than the definition 1.8 of the present paper, but really
they are
equivalent. We recall some notations and definitions of [T] in a slightly less
general setting (rough statements; see [T] for the exact statements). Let $G$ be a finite affine group scheme over $\p$, i.e. $G=\Spec R$
where $R$ is a finite-dimensional $\p$-algebra. Let $\mu:R \to R\otimes R$ be the
comultiplication of $R$. Such group $G$ is called a finite $v$-module ([T], Definition
3.1) if there is a homomorphism $\psi: \w \to \End_{gr. \ sch.}(G)$ satisfying some natural conditions (for
example, an analog of 1.3.1). Further, let $\Cal E_G$ be a $\p$-subspace of $R$ defined as
follows: $$\Cal E_G=\{x\in R \ \vert \ \mu(x)=x\otimes 1+1\otimes x\}$$
The map $x\mapsto x^q$ is a $\p$-linear map $\fr: \Cal E_G^{(1)}\to \Cal E_G$. Further,
the map $\psi(T): G \to G$ can be defined on $\Cal E_G$. Let $v: \Cal E_G\to\Cal
E_G^{(1)}$ be a map satisfying $\fr \circ v= \psi(T)-\theta$.

We consider two finite $v$-modules $G$, $H$, the above objects $\fr$, $v$ etc. will carry
the respective subscript. Let * be the dual in the meaning of linear algebra.
\medskip
{\bf Definition 1.12.1} ([T], 4.1). Two finite $v$-modules $G$, $H$
are called dual if there exists an
isomorphism $\alpha: \Cal E^*_H\to\Cal E_G$ such that if we denote by
$\goth v: \Cal E_G\to \Cal E^{(1)}_G$ a map which enters in the commutative diagram
$$ \matrix \Cal E^*_H & \overset{\fr^*_H}\to{\longrightarrow} & \Cal E^{*(1)}_H \\ & & \\
\alpha \downarrow & & \alpha ^{(1)}\downarrow \\ & & \\ \Cal E_G & \overset{\goth
v}\to{\longrightarrow} & \Cal E^{(1)}_G \endmatrix $$
then we have:
$$\fr_G\circ \goth v= \psi_G(T)-\theta \eqno{(1.12.2)}$$
i.e. $\goth v=v_G$.
\medskip
Let $M$ be a t-motive having $m(M)=1$, $E=E(M)$ the corresponding t-module and
$a\in \hbox{\bf
A}$. We denote $E_a$ --- the set of $a$-torsion elements of $E$ --- by $M_a$. It is a
finite
$v$-module.
\medskip
{\bf Proposition 1.12.3.} Let $M$, $M'$ be t-motives which are
dual in the meaning
of Definition 1.8. Then $\forall a\in \hbox{\bf A}$, $a\ne 0$ we have: $M_a$,
$M'_a$ are dual in the meaning of 1.12.1 = [T], Definition 4.1.
\medskip
{\bf Proof.} Condition $a\in \n F_q[T]$ implies that multiplication by $\tau$ is
well-defined on $M/aM$.
\medskip
{\bf Lemma 1.12.3.1.} We have canonical isomorphisms $i: M/aM \to
\Cal E_{M_a}$,
$i^{(1)}: M/aM \to \Cal E^{(1)}_{M_a}$ such that the following
diagrams are commutative:
$$ \matrix M/aM & \overset{\tau}\to{\longrightarrow} & M/aM &&& M/aM &
\overset{T}\to{\longrightarrow} & M/aM
\\ & & &&&&&
\\ i^{(1)}\downarrow & & i\downarrow &&& i\downarrow & & i\downarrow \\ & & &&&&&
\\ \Cal E^{(1)}_{M_a} & \overset{fr}\to{\longrightarrow} & \Cal E_{M_a} &&&
\Cal E_{M_a} & \overset{\psi_T}\to{\longrightarrow} & \Cal E_{M_a}
\endmatrix $$
\medskip
{\bf Proof.} Let $R$ be a ring such that $\Spec R=M_a$. The pairing between $M$ and $E$
shows that there exists a map $M\to R$ which is obviously factorized via an inclusion
$M/aM\to R$. It is easy to see that the image of this inclusion is contained in $\Cal
E_{M_a}$, i.e. we get $i$. Since $\dim_{\p}(M/aM)=\deg a \cdot r(M)$ and
$\dim_{\p}(R)=q^{\deg a \cdot r(M)}$ we get from [T], Definition 1.3 that $i$ is an
isomorphism. Other statements of the lemma are obvious. $\square$
\medskip
This lemma means that we can rewrite Definition 1.12.1 for the case $G=M_a$, $H=N_a$ by
the following way:\footnotemark \footnotetext{Here and below a t-motive $N$ should not be confused with $N$ of 1.9.2.}
\medskip
{\bf 1.12.3.2.} Two finite $v$-modules $M_a$, $N_a$ are dual if there
exists an
isomorphism $\alpha: (N/aN)^* \to M/aM $ such that after identification
via $\alpha$ of
$\tau^*: (N/aN)^* \to (N/aN)^*$ with a map $\goth v: M/aM \to M/aM$ we
have on $M/aM$:
$$\tau\circ \goth v= t-\theta \eqno{(1.12.3.3)}$$
We need a
\medskip
{\bf Lemma 1.12.3.4.} For $i=1,2$ let $N_i$ be a free $\p[T]$-module of
dimension $r$ with a
base $f_{i*} =(f_{i1}, ... , f_{ir} )$, let $\varphi_i:N_i \to N_i$ be
$\p[T]$-linear
maps having matrices $\goth Q_i$ in $f_{i*} $ such that $\goth
Q_2=\goth Q_1^t$, and let
$a$ be as above. Let, further, $\varphi_{i,a}: N_i/aN_i \to N_i/aN_i $
be the natural
quotient of $\varphi_i$. Then there exist $\p$-bases $\tilde f_{i*} $
of $N_i/aN_i $ such
that the matrix of $\varphi_{1,a}$ in the base $\tilde f_{1*} $ is
transposed to the
matrix of $\varphi_{2,a}$ in the base $\tilde f_{2*} $.
\medskip
{\bf Proof.} We can identify elements of $N_2$ with $\p[T]$-linear
forms on $N_1$
(notation: for $x\in N_2$ the corresponding form is denoted by
$\chi_x$) such that
$\chi_{\varphi_2(x)}=\chi_x\circ\varphi_1$. Any $\p[T]$-linear form
$\chi$ on $N_i$
defines a $\p[T]/a\p[T]$-linear form on $N_i/aN_i$ which is denoted by
$\chi_a$. Let now
$x\in N_2/aN_2$, $\bar x$ its lift on $N_2$, then
$\chi_{x,a}=(\chi_{\bar x})_a$ is a
well-defined $\p[T]/a\p[T]$-linear form on $N_1/aN_1$. For $x\in
N_2/aN_2$ we have
$$\chi_{\vf_{2,a}(x),a}= \chi_{x,a}\circ \vf_{1,a}$$
Further, let $\lambda: \p[T] \to \p$ be a $\p$-linear map such that
\medskip
{\bf 1.12.3.5.} Its kernel does not contain any non-zero ideal of
$\p[T]/a\p[T]$.
\medskip
(such $\lambda$ obviously exist.) For $x\in N_2/aN_2$ we denote $\lambda \circ
\chi_{x,a}$ by $\psi_x$, it is a $\p$-linear form on $\p$-vector space
$N_1/aN_1$.
Obviously condition (1.12.3.5) implies that the map $x\mapsto \psi_x$
is an isomorphism
from $N_2/aN_2$ to the space of $\p$-linear forms on $\p$-vector space
$N_1/aN_1$, and we
have
$$\psi_{\vf_{2,a}(x)}=\psi_x\circ\vf_{1,a}$$
which is equivalent to the statement of the lemma. $\square$
\medskip
Finally, the proposition follows immediately from this lemma
multiplied by $T-\theta$,
formula 1.10.1 and 1.12.3.2. $\square$
\medskip
{\bf Remark.} Let $a=\sum_{i=0}^k g_iT^i$, $g_i\in \n F_q$, $g_k=1$.
Taguchi ([T], proof
of 5.1 (iv)) uses the following $\lambda$: $\lambda(T^j)=0$ for $j<k-1$,
$\lambda(T^{k-1})=1$. It is easy to check that for $x=(T^i+
T^{i-1}g_{k-1}+T^{i-2}g_{k-2}+ ...
+g_{k-i})f_{2j}$ for this $\lambda$ we have: $\psi_x(T^i f_{1j})=1$,
$\psi_x(T^{i'}
f_{1j'})=0$ for other $i'$, $j'$.
\medskip
{\bf 1.13.} We consider in Sections 1.13, 1.14 the case of arbitrary $\w\supset \n F_q[T]$.
\medskip
A t-motive over {\bf A} is defined for example in [BH], p.1. Let us reproduce this
definition for the case of characteristic 0. Let $J$ be an ideal of
$\w_C$ generated by the elements $a\otimes 1 - 1 \otimes \iota(a)$ for all $a\in \w$. The ring $\w_C[\tau]$ is defined by the formula $\tau\cdot (a\otimes z)=(a\otimes z^q)\cdot \tau$, $a\in \w$, $z\in \p$.
\medskip
{\bf Definition 1.13.1.} A t-motive $M$ over {\bf A} is a pair $(M, \tau)$
where $M$ is
a locally free $\hbox{\bf A}_C$-module and $\tau$ is an $\hbox{\bf A}_C$-linear map
$M^{(1)}\to M$ satisfying the following analog of 1.3.1, 1.9.2:
$$\exists m \hbox{ such that } J^m(M/\tau (M^{(1)}))=0\eqno{(1.13.2)}$$

{\bf Remark 1.13.3.} We can consider $M$ as an $\w_C[\tau]$-module using the following formula for the product $\tau \cdot m$:
$$\tau \cdot m = \tau(m\otimes 1)$$
where $m\in M$, $m\otimes 1\in M^{(1)}$.
\medskip
The rank of $M$ as a locally free $\hbox{\bf A}_C$-module is called the rank of the
corresponding t-motive $(M, \tau)$. If $\w=\n F_q[T]$ then $M^{(1)}$ is isomorphic to
$M$, we can
consider $M$ as a $\p[T,\tau]$-module, and it is possible to show that in this case
1.13.2 implies
that $M_{\p[\tau]}$ is a free $\p[\tau]$-module. In the general case, the dimension $n$
of $(M, \tau)$
is defined as $\dim_{\p}(M/\tau (M^{(1)}))$.

Let us fix $\goth C=(\goth C, \tau_{\goth C})$ --- a t-motive of rank 1 over
{\bf A}.
For a t-motive $M=(M,\tau_M)$ a t-motive ${M_\goth C'}$ --- the
$\goth C$-dual of $M$
--- is defined as follows. We put
$M_\goth C'=\Hom_{\w_C}(M,\goth C)$.
Since for any locally free $\hbox{\bf A}_C$-modules $M_1$, $M_2$ we have
$$\Hom_{\w_C}(M_1,M_2)^{(1)}=\Hom_{\w_C}(M_1^{(1)},M_2^{(1)})$$
we can define $\tau(M_\goth C')$ by the following formula:

$$\hbox{ For } \vf\in \Hom_{\w_C}(M,\goth C)^{(1)} \hbox{ we have }
\tau(M_\goth C')(\vf)=\tau_{\goth C}\circ \vf \circ \tau_M^{-1}$$
\medskip
{\bf 1.14. Duality for abelian $\tau$-sheaves.} We use notations of [BH],
Definition 2.1 if they do not differ from the notations of the present
paper; otherwise we continue to use notations of the present paper
(for example, $d$ (resp. $\sigma^*(\goth X)$ for any object $\goth X$)
of [BH] is
$n$ (resp. $\goth X^{(1)}$) of the present paper). For any abelian
$\tau$-sheaf $\underline{\Cal F}$ we denote its $\Pi_i$, $\tau_i$ by
$\Pi_i(\underline{\Cal F})$, $\tau_i(\underline{\Cal F})$
respectively. If $M$, $N$ are invertible sheaves on $X$ and $\rho: M
\to N$ a rational map then we denote by $\rho^{inv}: N \to M$ the
rational map which is inverse to $\rho$ with respect to the composition.
We define
$\tau_{\goth r,i-1}(\underline{\Cal F})$ (the rational $\tau_i$) as
the composition map $\tau_{i-1}(\underline{\Cal F})
\circ {\Pi_{i-1}^{(1)}}^{inv}(\underline{\Cal F})$, it is a rational
map from $\Cal F_i^{(1)}$ to $\Cal F_i$.

Let $\underline{\Cal O}$ be a fixed abelian
$\tau$-sheaf having $r=n=1$. The $\underline{\Cal O}$-dual abelian
$\tau$-sheaf $\underline{\Cal F}'= \underline{\Cal
F}'_{\underline{\Cal O}}$ is defined by the formulas
$$\Cal F'_0=\Hom_{X}(\Cal F_0, \Cal O_0)$$ where Hom is the sheaf's
one, and the map $\tau_{\goth r,-1}(\underline{\Cal F}'): {\Cal
F'_0}^{(1)}\to \Cal F'_0$ is defined as follows. We have ${\Cal
F'_0}^{(1)}= \Hom_{X}({\Cal F_0}^{(1)}, {\Cal O_0}^{(1)})$. Let
$\gamma\in\Hom_{X}({\Cal F_0}^{(1)}, {\Cal O_0}^{(1)})(U)$ where $U$
is a sufficiently small affine subset of $X_{\p}$, such that $\gamma:
{\Cal F_0}^{(1)}(U) \rightarrow {\Cal O_0}^{(1)}(U)$.
\medskip
{\bf 1.14.1.} We define: $[[\tau_{\goth r,-1}(\underline{\Cal
F}')](U)](\gamma)$ is the following composition map:
$$\Cal F_0(U)\overset{[\tau_{\goth r,-1}^{inv}(\underline{\Cal
F})](U)}\to{\longrightarrow}
\Cal F_0^{(1)}(U)\overset{\gamma}\to{\to}\Cal O_0^{(1)}(U)
\overset{[\tau_{\goth r,-1}(\underline{\Cal
O})](U)}\to{\longrightarrow}\Cal O_0(U)\in \Hom_{X}({\Cal F_0}, {\Cal
O_0})(U)$$

Clearly that this definition and the definitions 1.8, 1.13 are
compatible with the forgetting functor $\underline{M}(\underline{\Cal
F})$ from abelian $\tau$-sheaves to pure Anderson t-motives of [BH],
Section 3, page 8.
\medskip
{\bf 1.15. Duality over fields.} Let $L\supset \n F_q(\theta)$ be a field extension of $\n F_q(\theta)$, and $M$ a t-motive over $L$ (i.e. a pair ($M$, an $L$-structure on $M$)). Obviously we have
\medskip
{\bf Proposition 1.15.1.} The notion of duality for $M$ over $L$ is well-defined. $\square$
\medskip
Similarly, we have a proposition for Galois action:
\medskip
{\bf Proposition 1.15.2.} Let $M$ be defined over $\overline{ \n F_q(\theta)}$ and $\gamma\in \Gal( \n F_q(\theta))$. Then $(\gamma(M))'=\gamma(M')$. $\square$
\medskip
\medskip
{\bf 1.16. Case of finite characteristic.} Let $\iota: \w\to \bar \n F_q$ be a map of finite characteristic, we denote $\Ker \iota$ by $\Cal P$. The definition of t-motive for this case is similar to 1.3, see [G] for the details. The definition of duality also is similar to the one of the case of generic characteristic. Duality commutes with reduction. Namely, let $M$ be from 1.15, $\goth P$ a prime of $L$ not over the infinity of $\n F_q(\theta)$, $\Cal P\subset \w$ is $\iota^{-1}(\goth P\cap \n F_q[\theta]$) --- the finite characteristic. We consider the case of good reduction of $M$ at $\goth P$, we denote it by $\tilde M$. It is a t-motive in characteristic $\Cal P$.  Let $M$ have dual $M'$.
\medskip
{\bf Proposition 1.16.1.} $\tilde M$ has dual iff $M'$ has good reduction at $\goth P$; in this case they coincide. $\square$
\medskip
{\bf Remark 1.16.2.} Apparently if $M$ has good reduction and dual, then $M'$ also has good reduction (in this case 1.16.1 means that $M'$ exists implies $(\tilde M)'$ exists). For standard-3 t-motives (this is a simple tipe of t-motives, see 11.8.1) apparently this can be shown by explicit calculations.
\medskip
{\bf Remark 1.16.3.} Clearly 1.16.1 is true for the case of bad reductions. I do not give exact definitions for this case.
\medskip
{\bf 1.16.4. Ordinarity.} Let $M$ be of finite characteristic. By analogy with the number field case, $M$ is called ordinary if its Newton polygon consists of 2 segments. If $N=0$ then the Newton polygon of $M'$ is the dual of the one of $M$ (the notion of duality of polygons is clear; apparently the condition $N=0$ can be omitted). So, we have
\medskip
{\bf Proposition 1.16.5.} $M$ is ordinary $\iff M'$ is ordinary. $\square$
\medskip
See 13.4.1 for a more exact result.
\medskip
{\bf 2. Analytic duality.}
\medskip
We consider in the present section the case of arbitrary $\w\supset \n F_q[T]$ (and $N=0$ as usually).
\medskip
Condition $N=0$ implies that an element $a\in \w$ acts on $\Lie(M)$ by multiplication by $\iota(a)$. Hence, we have a
\medskip
{\bf Definition 2.1.} Let $V$ be the space $\p^n$. A locally free $r$-dimensional
$\n Z_\infty$-submodule
$L$ of $V$ is called a lattice if

(a) $L$ generates $V$ as a $\p$-module and

(b) The $\n R_\infty$-linear span of $L$ has dimension $r$ over
$\n R_\infty$.
\medskip
Numbers $n$, $r$ are called the dimension and the rank of $L$
respectively. Attached to $(L,V)$ is the tautological inclusion $\vf=\vf(L,V):
L \to V$. We shall consider the category of triples $(\vf, L, V)$; a map $\psi:
(\vf, L,V)\to(\vf_1, L_1,V_1)$ is a pair $(\psi_L, \psi_V)$ where $\psi_L: L \to L_1$ is a $\n Z_\infty$-linear map, $\psi_V: V\to V_1$ is a $\p$-linear map such that $\vf_1\circ \psi_L=\psi_V \circ \vf$.

Inclusion $\vf$ can be extended to a map
$L\underset{\n Z_\infty}\to{\otimes}\p \to V$ (which is surjective
by 2.1a), we denote it by $\vf=\vf(L,V)$ as well.
We can also attach to $(L,V)$ an exact sequence
$$0\to \Ker \vf \to L\underset{\n Z_\infty} \to{\otimes}\p \overset{\vf}\to{\to} V \to
0\eqno{(2.2)}$$

Let $\Cal I\in \Cl(\w)$ be a class of ideals; we shall use the same notation $\Cal I$ to
denote a representative in the $\iota$-image of this class. Let $(\vf', L', V')$ be another lattice and $D$ a structure of a perfect $\Cal I$-pairing $<* , * >_D$ between $L$
and $L'$. Let us fix an isomorphism $$\alpha: \Cal I \underset{\n Z_\infty} \to{\otimes} \p \to \p\eqno{(2.2')}$$ $D$ extends via $\alpha$ to a perfect $\p$-pairing between $L\underset{\n Z_\infty}\to{\otimes}\p$
and $L'\underset{\n Z_\infty}\to{\otimes}\p$, we denote this pairing by $D_{\alpha, \infty}$.
\medskip
{\bf Definition 2.3.} Two lattices $(\vf, L, V)$ and $(\vf', L', V')$ are called $(\alpha, \Cal I)$-dual if there exists a perfect $\Cal I$-pairing $D$ between $L$
and $L'$ such that $\Ker \vf \subset L\underset{\n Z_\infty}\to{\otimes}\p$, $\Ker \vf' \subset L'\underset{\n Z_\infty}\to{\otimes}\p$ are mutually orthogonal with respect to $D_{\alpha, \infty}$.
\medskip
Let $(n,r)$, $(n',r')$ be the dimension and rank of $(\vf, L, V)$ and $(\vf', L', V')$ respectively. If they are $(\alpha, \Cal I)$-dual then $r'=r$, $n'=r-n$. There exists the following reformulation of the definition of duality. $D_{\alpha, \infty}$ induces an isomorphism
$\gamma_{\alpha, D}: (L\underset{\n Z_\infty}\to{\otimes}\p)^* \to
L'\underset{\n Z_\infty}\to{\otimes}\p$ (here and below for any object $W$ we
denote $W^*=\Hom_{\p}(W,\p)$ ).
\medskip
{\bf Property 2.4.} $(\vf, L, V)$ and $(\vf', L', V')$ are $(\alpha, \Cal I)$-dual iff there exists an isomorphism from $(\Ker \vf)^*$
to $V'$ making the following diagram commutative: $$\matrix 0 & \to &
V^* & \overset{\vf^*}\to{\to} &  (L\underset{\n Z_\infty}
\to{\otimes}\p)^* & \to & (\Ker \vf)^* & \to & 0\\
&&&&&&&& \\
& & \downarrow & & \gamma_{\alpha, D} \downarrow && \downarrow \\ &&&&&&&& \\ 0
& \to & \Ker \vf' & \rightarrow  & L'\underset{\n Z_\infty}\to{\otimes}\p &
\overset{\vf'}\to{\to} & V' & \to &
0\endmatrix \eqno{(2.5)}$$

Further, this property is equivalent to the following two conditions:
\medskip
{\bf 2.6.} $\dim V'=r-n$;
\medskip
{\bf 2.7.} The composition map $\vf'\circ \gamma_D \circ \vf^*:
V^* \to V'$ is 0.
\medskip
Both 2.4 and (2.6, 2.7) are obvious.
\medskip
{\bf Remark 2.8.} It is
easy to see that the functor $(\vf, L, V) \mapsto (\vf', L', V')$ is well-defined
on a subcategory
(not all lattices have duals, see
below) of the category of the triples $(\vf, L, V)$, it is contravariant and
involutive.
\medskip
{\bf 3. Explicit formulas for analytic duality.}
\medskip
Here we consider the case $\w=\n
F_q[T]$. In this case $\Cl(\w)=0$, and $(\alpha, \Cal I)$-dual is called simply
dual. The coordinate
description of the dual lattice is the following. Let
$e_1, ..., e_r$ be a $\n Z_\infty$-basis of
$L$ such that $\vf(e_1), ..., \vf(e_n)$ form a $\p$-basis of $V$. Like in the
theory of abelian
varieties, we denote by $Z=(z_{ij})$ the Siegel matrix whose lines are
coordinates of $\vf(e_{n+1}), ..., \vf(e_r)$ in the basis $\vf(e_1), ..., \vf(e_n)$, more
exactly, the size of $Z$ is $(r-n)\times n$ and
$$\forall i =1,..., r-n \ \ \ \ \vf(e_{n+i})=\sum_{j=1}^n z_{ij}\vf(e_j)\eqno{(3.1)}$$ $Z$
defines $L$, we denote $L$ by $\goth L(Z)$.
\medskip
{\bf Proposition 3.2.} $[\goth L(Z)]'=\goth L(-Z^t)$, i.e. a Siegel matrix of the dual
lattice is the minus transposed Siegel matrix.
\medskip
{\bf Proof.} Follows immediately from the definitions. Really, let $f_1, ..., f_r$ be a
basis of $L'$, we define the pairing by the formula $$<e_i,
f_j>=\delta_i^j\eqno{(3.3)}$$ and the map $\vf'$ by the formula
$$\forall i =1,... ,n \ \ \ \ \vf'(f_{i})=\sum_{j=1}^{r-n} -z_{ji}\vf'(f_{n+j})$$
(minus transposed Siegel matrix).
$\Ker \vf$ is generated by elements $$v_i=e_{n+i}-\sum_{j=1}^n z_{ij}e_j, \ \ \ \ i
=1,..., r-n$$ and $\Ker \vf'$ is generated by elements $$w_i=f_{i}+\sum_{j=1}^{r-n}
z_{ji}f_{n+j}, \ \ \ \ i =1,..., n\eqno{(3.4)}$$ It is sufficient to check that $\forall i,j$ we have
$<v_i, w_j>=0$; this follows immediately from 3.3. $\square$
\medskip
{\bf Remark 3.5.} $L'$ exists not for all $L$.
Trivial counterexample: case $n=r=1$. To get another counterexamples, we use that for
$n=1$ (lattices of Drinfeld modules) a Siegel matrix is a column matrix $Z=\left(\matrix
z_1&...&z_{r-1} \endmatrix \right)^t$ and
$$ \goth L(Z) \hbox{ is not a lattice } \iff 1,z_1, ... ,z_{r-1} \hbox{ are linearly
dependent over } \n R_\infty \eqno{(3.6)}$$ while for $n=r-1$ a Siegel matrix is a
row matrix $Z=\left(\matrix -z_1&...&-z_{r-1} \endmatrix \right)$ and
$$ \goth L(Z) \hbox{ is not a lattice } \iff \forall i \ \ z_i\in \n R_\infty
\eqno{(3.7)}$$ Since condition (3.7) is strictly stronger than (3.6) we see
that all lattices having $n=1$, $r>1$ have duals while not all lattices having $n=r-1$,
$r>2$ have duals.

It is clear that almost all matrices have duals. Here "almost all" has the same meaning
that as "Almost all matrices $Z$ are a Siegel matrice of a lattice", i.e. if we choose an
(infinite) basis of $\p/\n R_\infty$, then coordinates of the entries of $Z$ in this
basis must satisfy some polynomial relations in order that $Z$ is not a Siegel matrice of
a lattice.
\medskip
{\bf Remark 3.8.} The coordinate proof of the theorem that the notion
of the dual lattice is well-defined, is the following. Two Siegel matrices
$Z$, $Z_1$ are called equivalent iff there exists an isomorphism of their pairs $(\goth
L(Z), V)$, $(\goth L(Z_1),V_1)$. Like in the classical theory of modular forms, $Z$,
$Z_1$ are equivalent iff there exists a matrix $\gamma \in GL_r(\n Z_\infty)=\left(\matrix A&B\\ C&D \endmatrix \right)$ ($A,B,C,D$ are the ($n\times
n$), ($n\times r-n$), ($r-n\times n$), ($r-n\times r-n$)-blocks of $\gamma$ respectively;
we shall call this block structire by the $(n, r-n)$-block structure) such that
$$C+DZ=Z_1(A+BZ)\eqno{(3.8.1)}$$

Let $A_1,B_1, C_1, D_1$ be the $(n, r-n)$-block structure of the matrix $\gamma^{-1}$.
The equality
$$-C_1^t+A_1^tZ^t={Z_1}^t(D_1^t-B_1^tZ^t)\eqno{(3.8.2)}$$
shows that if $Z$, $Z_1$ are equivalent
then $-Z^t$, $-Z_1^t$ are equivalent. [Proof of (3.8.2): (3.8.1) implies $Z_1=(C+DZ)(A+BZ)^{-1}$; substituting this value of $Z_1$ to the transposed
(3.8.2), we get $-C_1+ZA_1=(D_1-ZB_1)(C+DZ)(A+BZ)^{-1}$, or $(-C_1+ZA_1)(A+BZ)=(D_1-ZB_1)(C+DZ)$. This formula follows immediately from
$\left(\matrix A_1&B_1\\C_1&D_1\endmatrix \right)\left(\matrix A&B\\C&D\endmatrix \right)=\left(\matrix I_n&0\\0&I_{r-n}\endmatrix \right)$].

Further, let $\alpha: (L_1\subset \p^n) \to
(L_2\subset \p^n)$ be a map of lattices. If $L'_1$, $L'_2$
exist, then the map
$\alpha': (L'_2\subset \p^{r-n}) \to (L'_1\subset \p^{r-n})$ is defined by the
following formulas. Let $Z_i$ be the Siegel matrices of $L_i$ in the
bases
$e_{i1}, ... e_{ir}$ of $L_i$ ($i=1,2$). Let us consider the matrix
$\goth M=(m_{ij})\in
M_{r}(\n Z_\infty)$ of $\alpha$ in the bases $e_{i1}, ...,
e_{ir}$ (i.e.
$\alpha(e_{1i})=\sum_j m_{ij} e_{2j}$). Let $f_{i1}, ..., f_{ir}$ be the dual
base of $L'_i$ (see 3.3) and $e'_{i1}, ... e'_{ir}$ another base of $L'_i$ defined by $$e'_{ij}=f_{i,j+n}, \ \ \ \ j+n \mod r\eqno{(3.8.3)}$$ Formulas (3.8.3), (3.4) show that an analog of 3.1 is satisfied for both bases $e'_{i1}, ..., e'_{ir}$, their Siegel matrices are $-Z_i^t$.

Let
$$\goth M=\left(\matrix \goth M_{11} & \goth M_{12} \\ \goth M_{21} &
\goth M_{22}
\endmatrix \right)$$
be the $(n, r-n)$-block structure of $\goth M$. The matrix of $\alpha'$ in the bases $f_{i1}, ..., f_{ir}$ is $\goth M^t$, and using the matrix 3.8.3 of change of base, we get that $\goth M'$ --- the matrix of $\alpha'$ in
the bases $e'_{i1}, ..., e'_{ir}$ --- has the following $(r-n,
n)$-block structure:
$$\goth M'=\left(\matrix \goth M_{22}^t & \goth M_{12}^t \\
\goth M_{21}^t & \goth
M_{11}^t \endmatrix \right)\eqno{(3.8.4)} $$
The property that $\goth M$ comes from a $\p$-linear map $\p^n \to \p^n$ implies
that $\goth M'$ comes from a $\p$-linear map $\p^{r-n} \to
\p^{r-n}$. This follows immediately from the definition of dual lattice, or can be easily checked algebraically.
\medskip
{\bf Remark 3.9.} Taking $\gamma =\left(\matrix 1&0\\ 0&-1 \endmatrix \right)$ we get
that $Z$ is equivalent to $-Z$, hence $Z'$ is also a Siegel matrix of the dual lattice.
\medskip
{\bf 4. Main conjecture for arbitrary $\w$}.
\medskip
The main result of the paper is the following Theorem 5 on coincidence of algebraic and
analytic duality. We formulate it as a conjecture 4.1 for any $\w$, but we prove it only for the case $\w=\n
F_q[T]$.  Let $M$ be a uniformizable t-motive. Its lattice $L(M)$ is really a lattice in
the meaning of Definition 2.1, because [A],
Corollary 3.3.6 (resp. [G], Lemma 5.9.12) means that it satisfies 2.1a (resp. 2.1b);
recall that we consider the case $N=0$,
i.e. the action of $T$ on $\Lie(M)$ is simply multiplication by $\theta$.
Let us fix (like in 1.13) $\goth C=(\goth C, \tau_{\goth C})$ --- a t-motive of rank 1 over
{\bf A}, and let $L(\goth C)$ be its lattice. It is a $\n Z_\infty$-module. $\Omega=\Omega(\w)$ is an $\w$-module, we consider a $\n Z_\infty$-module $\iota^{-1}(\Omega)$. There exists the   notion of the $L(\goth C)\otimes \iota^{-1}(\Omega)$-duality.
\medskip
{\bf Conjecture 4.1.} Let $M$ be a uniformizable t-motive having $N=0$ such
that its $\goth C$-dual $M'$ exists. Then $M'$ is uniformizable, it has $N':=N(M')=0$, and $(L(M), \Lie(M))$ and $(L(M'), \Lie(M'))$ are $\alpha, L(\goth C)\otimes \iota^{-1}(\Omega)$-dual for some $\alpha$ from $2.2'$ (it can be explicitly described).
\medskip
We give in Section 5 the first step of the proof of this conjecture.
\medskip
{\bf 5. Main theorem.}
\medskip
Recall that the word "canonical" means "canonical up to multiplication by elements of $\n F_q^*$".
\medskip
{\bf Theorem 5.}\footnotemark \footnotetext{The proof of this theorem was inspired by a result of Anderson, see Section 6 for details.} Let $M$ be a uniformizable t-motive over $\w=\n F_q[T]$ having $N=0$ such
that its dual $M'$ exists and has $N':=N(M')=0$. Then $M'$ is uniformizable, and $(L(M), \Lie(M))$ and $(L(M'), \Lie(M'))$ are dual.
\medskip
{\bf Remark 5A.} Condition $N'=0$ holds for pure $M$ (Theorem 10.3) and for a large class of non-pure $M$ (Theorem 11.5). Most likely, a modification of the end of the proof of the present theorem will permit us to prove that $N'=0$ holds for all $M$ having $N=0$ and having dual.
\medskip
{\bf Remark 5B.} A reformulation of the theorem in terms of Hodge-Pink structures is given in Section 9.
Proof of the theorem for the case $N\ne0$ is given in [GL18].
\medskip
{\bf Corollary 5.1.1.} If $\w=\n F_q[T]$ then a Siegel matrix of $M'$ is the minus transposed
of a Siegel matrix of $M$.
\medskip
In the section 8 below we give a corollary of this theorem and some conjectures related to the
problem of 1 -- 1 correspondence between t-motives and lattices.
\medskip
{\bf 5.1.2. Some definitions.} Recall that $E=E(M)$ is isomorphic to $\p^n$. There is a structure of $\w$-module on $E$; multiplication by $T$ is denoted by $m_T$, and this operator $m_T$ is defined in coordinates by the formula
$$m_T(x)=\sum_{i=0}^l\goth A_ix^{(i)}$$ where $x\in E=\p^n$ is a vector column, $\goth A_i$ are from 1.9.1. There is a map $\exp: \Lie(M) \to E$ making the following diagram commutative:
$$\matrix \Lie(M) &
\overset{\Exp}\to{\to} & E \\ \\ \theta \downarrow & & m_T
\downarrow \\ \\ \Lie(M) & \overset{\Exp}\to{\to} & E \endmatrix \eqno{(5.1.3)}$$
By definition, $L(M)=\Ker \Exp$.

We need another space $\Lie_T(M)$ together with an isomorphism $\goth a:\Lie_T(M) \to \Lie(M)$ and a structure of $\w$-module on $\Lie_T(M)$ such that the multiplication by $T$ on $\Lie_T(M)$ is simply the multiplication by $\theta$ on $\Lie(M)$, i.e. $$\goth a(Tx)=\theta\cdot(\goth a(x))\eqno{(5.1.4)}$$ where $x\in\Lie_T(M)$. Commutativity of 5.1.3 means that $\Exp\circ\goth a:\Lie_T(M)\to E$ is a map of $\w$-modules.
\medskip
{\bf 5.1.5.} We shall work merely with $L_T(M):=\Ker (\Exp\circ\goth a)\subset \Lie_T(M)$ rather than $L(M)$. Clearly $L_T(M)$ is an $\w$-module, $\goth a:L_T(M)\to L(M)$ is an isomorphism satisfying 5.1.4 for $x\in L_T(M)$.
\medskip
The proof of Theorem 5 consists of two steps. We formulate and prove Step 1 for the case of arbitrary $\w$.
\medskip
{\bf Step 1.} For the above $M$, $M'$ we have:
\medskip
(A) Uniformizability of $M$ implies uniformizability of $M'$.
\medskip
(B) There exists a canonical $\w$-linear $L_T(\goth C)\otimes
\Omega$-valued perfect pairing $<* , * >_M$
between $L_T(M)$ and $L_T(M')$ (by 5.1.5, this is the same as the $\n Z_\infty$-linear pairing between $L(M)$ and $L(M')$, which, in its turn, is $D$ of Definition 2.3). It is functorial.
\medskip
{\bf Remark 5.1.6.} Practically, (B) comes from [T], Theorem 4.3 (case
$\w=\n F_q[T]$). Really, to define a pairing between $L(M)$ and
$L(M')$ it is sufficient to define (concordant) pairings between $L(M)
/aL(M)$ and $L(M') /aL(M')$ for any $a\in \w$. Since $M_a:=E(M)_a=L(M)/aL(M)$ and because of Proposition 1.12.3 which affirms that
$M_a$ and $M'_a$ are Taguchi-dual, we see that [T], Theorem 4.3 gives
exactly the desired pairing.
\medskip
We give two versions of the proof of Step 1: the first one --- for the general case of arbitrary $\w$ and the second one --- for the case $\w=\n F_q[T]$ --- is based on explicit calculations, it is used for the proof of Step 2.
\medskip
{\bf 5.2. Proof: Step 1, Version 1.} Here we consider the general case of arbitrary $\w$.
Let $\Omega=\Omega(\w/\n
F_q)$ be the module of differential forms; we can consider it as an
element of $\Cl(\w)$. We use formulas and notations of [G], Section 5.9
modifying them to the case of arbitrary $\w$. For example, {\bf
A} (resp. {\bf K}) of [G], 5.9.16 is {\bf A} (resp. {\bf K})
of the present paper (recall that $\bar K$ (resp. $\bar K[T,\tau]$) of [G] is
$\p$ (resp. $\w_C[\tau]$, see 1.13) of the present paper). Hence, we denote $\bar K\{T\}$ of
[G], Definition 5.9.10 by $\p\{T\}$. For the general case it must be
replaced by a ring $Z_0$ defined by the formula
$$Z_0:=\w\underset{\n F_q[T]}\to{\otimes}\p\{T\}\eqno{(5.2.1)}$$
$Z_0$ is a $\w_C[\tau]$-module, i.e. $\tau$ acts on $Z_0$, and $Z_0^\tau=\w$.

$Z_1$ for the present case is defined by the same formula [G], 5.9.22. Explicitly,
$$Z_1:=\Hom^{cont}_{\w}(\x/\w,\p)\eqno{(5.2.1a)}$$
It is a locally free $Z_0$-module of dimension 1 (the module structure
is compatible with the action of $\tau$; see [G], p. 168, lines 3 - 4
for the case $\w=\n F_q[T]$). We have: $Z_1^\tau$ is a
$Z_0^\tau$-module ( = $\w$-module) which is isomorphic to $\Omega(\w)$
(see the last lines of the proof of [G], Corollary 5.9.35 for the case
$\w=\n F_q[T]$), and $Z_1$ is isomorphic to $Z_0\otimes_{\w}\Omega(\w)$.

We shall consider $M$ as a $\w_C[\tau]$-module, like in 1.13.3. We denote
$M\{T\}:=M\underset{\w_C}\to{\otimes}Z_0$ ( = [G],
Definition 5.9.11.1 for the case $\w=\n F_q[T]$) and
$H^1(M):=M\{T\}^\tau$ like in [G], Definition 5.9.11.2.
Analogous to [G], Corollary 5.9.25 we get that for the present case
$$H_1(M):=\Hom_{\w_C[\tau]}(M,Z_1)=L_T(M)$$
($H_1(M)=H_1(E)$ of [G], 5.9). Particularly, for $M=\goth C$ we have
$$L_T(\goth C)=\Hom_{\w_C[\tau]}(\goth C,Z_1)$$

{\bf Lemma 5.2.2. } $H_1(M')=H^1(M)\underset{\w}\to{\otimes}L_T(\goth C)$.
\medskip
{\bf Proof.} By definition,
$\Hom_{\w_C} (M',Z_1)=\Hom_{\w_C} (\Hom_{\w_C} (M, \goth C), Z_1)$. Further,
$$\Hom_{\w_C} (\Hom_{\w_C} (M, \goth C), Z_1)=(M\underset{\w_C}
\to{\otimes}Z_0)\underset{Z_0}\to{\otimes} (\Hom_{\w_C}(\goth C,
Z_1))\eqno{(5.2.3)}$$ (an equality of  linear algebra). In order to show that we can
consider $\tau$-invariant subspaces, we need the following objects.
Let $I$ be an ideal of $\w$, $\Cal M_0=IZ_0$. It is clear that $\Cal M_0^\tau=I$.
Further, let $\Cal M_1$ be a locally free $Z_0$-module. We have a formula: $$(\Cal
M_0\underset{Z_0}\to{\otimes} \Cal M_1)^\tau=\Cal M_0^\tau \underset{\w}\to{\otimes} \Cal
M_1^\tau \eqno{(5.2.4)}$$
Really, $\Cal M_0\underset{Z_0}\to{\otimes} \Cal M_1=I\Cal M_1$, and $$(I\Cal
M_1)^\tau=I\Cal M_1^\tau\eqno{(5.2.5)}$$ where this formula is true by the following
reason. Obviously $(I\Cal M_1)^\tau\supset I\Cal M_1^\tau$. Let $J$ be an ideal of $\w$
such that $IJ$ is a principal ideal. We have $(IJ(J^{-1}\Cal M_1))^\tau=IJ(J^{-1}\Cal
M_1)^\tau$ and $(IJ(J^{-1}\Cal M_1))^\tau\supset I(J(J^{-1}\Cal M_1))^\tau\supset
IJ(J^{-1}\Cal M_1)^\tau$, hence all these objects are equal and we get 5.2.5 and hence
5.2.4.

The action of $\tau$ on both sides of 5.2.3 coincide.  Considering $\tau$-invariant
elements of both sides of 5.2.3 and taking into consideration 5.2.4 ($\Cal
M_0=\Hom_{\w_C}(\goth C,
Z_1)$ and $\Cal M_1=M\underset{\w_C}
\to{\otimes}Z_0$) we get the lemma. $\square$
\medskip
This lemma proves (A) of Step 1.
\medskip
{\bf Lemma 5.2.6.} Let $\Cal M_i$ ($i=0,1$) be two locally free $Z_0$-modules with $\tau$-action satisfying $\tau(cm)=\tau(c) \tau(m)$ ($c\in Z_0$, $m\in \Cal M_i$), and $\psi: \Cal M_0\otimes_{Z_0}\Cal M_1\to Z_1$ a perfect pairing of $Z_0$-modules with $\tau$-action. Let, further, both $\Cal M_i$ satisfy $\Cal M_i^\tau\otimes_{\w} Z_0=\Cal M_i$. Then the restriction of $\psi$ to $\Cal M_0^\tau\otimes_{\w}\Cal M_1^\tau\to \Omega$ is a perfect pairing as well.
\medskip
{\bf Proof.} Let $\alpha: \Cal M_0^\tau \to \Omega$ be an $\w$-linear map. We prolonge it to a map $\bar \alpha: \Cal M_0 \to Z_1$ by $Z_0$-$\tau$-linearity. By perfectness of $\psi$, there exists $m_1\in \Cal M_1$ such that $\bar \alpha(m_0)= \psi(m_0 \otimes m_1)$. It is easy to see that $m_1$ is $\tau$-invariant (we use the fact that $\tau: Z_0 \to Z_0$ is surjective). $\square$
\medskip
{\bf Lemma 5.2.7.} There is a natural perfect $\w$-linear
$\Omega$-valued pairing
between $H_1(M)$ and $H^1(M)$:
$H_1(M)\underset{\w}\to{\otimes}H^1(M)\to \Omega$.
\medskip
{\bf Proof.} For the case $\w=\n F_q[T]$ this is [G], Corollary
5.9.35. General case: we have a perfect $Z_0$-pairing
$$\Hom_{\w_C}(M,Z_1)\underset{Z_0}\to{\otimes} (M\underset{\w_C}
\to{\otimes} Z_0) \to Z_1$$
Now we take $\Cal M_0=\Hom_{\w_C}(M,Z_1)$, $\Cal M_1=M\underset{\w_C} \to{\otimes}Z_0$ and we apply Lemma 5.2.6. $\square$
\medskip
Step 1 of the theorem follows from these lemmas. 
\medskip
{\bf Remark 5.2.8.} The pairing can be defined also as the composition of
$$\matrix H_1(M)\underset{\w} \to{\otimes} H_1(M')=\Hom_{\w_C[\tau]}(M,Z_1)
\underset{\w}
\to{\otimes} \Hom_{\w_C[\tau]}(M',Z_1) \\
\to \Hom_{\w_C[\tau]}(M\underset{\w_C} \to{\otimes}M'
,Z_1\underset{Z_0} \to{\otimes}Z_1) \to
\Hom_{\w_C[\tau]}(\goth C,Z_1\underset{Z_0} \to{\otimes}Z_1) = L_T(\goth C) \underset{\w} \to{\otimes}\Omega\endmatrix \eqno{(5.2.9)}$$
where the second map comes from a canonical map $\delta: \goth C \to M\underset{\w_C}
\to{\otimes}M'$ of Remark 1.11.1 (more exactly, of its analog for arbitrary $\w$).
\medskip
{\bf Remark 5.2.10.} Recall that the explicit formula for functoriality is
the following. Let $\alpha: M_1 \to M_2$ be a map of t-motives,
$\alpha': M'_2 \to M'_1$
the dual map and $L_T(\alpha): L_T(M_2) \to L_T(M_1)$, $L_T(\alpha'): L_T(M'_1)
\to L_T(M'_2)$ the
corresponding maps on lattices. For any $l_1'\in L_T(M_1')$, $l_2\in L_T(M_2)$ we
have:
$$<L_T(\alpha)(l_2), l_1'>_{M_1}=<l_2, L_T(\alpha')(l_1')>_{M_2}\eqno{(5.2.11)}$$

{\bf 5.3. Proof: Step 1, Version 2.} Case $\w=\n F_q[T]$. We identify $Z_1$ of [G], p.168, lines 3 -- 4 with $\p\{T\}$ (see [G], Definition 5.9.10) and $\w$ with $\Omega$. Like above, we have an isomorphism of $\w$-modules (recall that $\w$ is the center of $\p[T,\tau]$):
$$L_T(M)=\Hom_{\p[T,\tau]}(M,Z_1)\eqno{(5.3.1)}$$ ([G], first terms of 5.9.25, 5.9.19). Let $\vf: M \to Z_1$, $\vf':
M' \to Z_1$ be elements of $L_T(M)$, $L_T(M')$ respectively, and let
$f_*$, $f'_*$, $Q$, $Q'$ be from
1.9.3, 1.10. We denote $$\vf(f_*)=v_*\eqno{(5.3.2)}$$ where $v_*\in
(Z_1)^r$ is a vector column (it is a column of the scattering matrix ([A], p. 486) of
$M$, see 5.4.1 below). The same notation for
the dual:
$\vf'(f'_*)=v'_*$. Condition
that $\vf$, $\vf'$ are
$\tau$-homomorphisms is equivalent to $$Qv_*=v_*^{(1)}, \ \
Q'v'_*={v'}_*^{(1)}\eqno{(5.3.3)}$$ (analog of the formula for scattering matrices [A],
(3.2.2)). Let us consider $\Xi=\sum_{i=0}^\infty
a_iT^i\in\p\{T\}\subset\p[[T]]$ of [G], p. 172, line 1; recall that it is the only
element (up to
multiplication by $\n F_q^*$) satisfying
$$\Xi=(T-\theta)\Xi^{(1)}, \ \ \ \lim_{i\to\infty}a_i=0, \ \ \ |a_0|>|a_i| \ \ \forall
i>0\eqno{(5.3.4)}$$ (see [G], p. 171,
(*); there is a formula $\Xi=a_0\prod_{i\ge0}(1-T/\theta^{q^i})$ where $a_0$ satisfies
$a_0^{q-1}=-1/\theta$). Finally, we define
$$<\vf, \vf'>=\Xi v_{*}^tv'_{*}\eqno{(5.3.5)}$$
Obviously $<\vf, \vf'>$ does not depend on a choice of a basis $f_*$.
\medskip
{\bf Lemma 5.3.6.} $<\vf, \vf'>\in \w$.
\medskip
{\bf Proof.} Firstly, this element belongs to $\n F_q[[T]]$, because $$\Xi
v_{*}^tv'_{*} -(\Xi v_{*}^tv'_{*})^{(1)}=\Xi (v_{*}^tv'_{*}-(T-\theta)^{-1}
v_{*}^{(1)t}{v'}_{*}^{(1)})= \Xi v_{*}^t(I_r-(T-\theta)^{-1}
Q^tQ')v'_{*}$$ because of (5.3.3). But we have (see (1.10.1) --- the
definition of $Q'$)
$$I_r-(T-\theta)^{-1} Q^tQ'=0$$
Secondly, let $<\vf, \vf'>=\sum_{i=0}^\infty c_iT^i$. Since
coefficients of all factors of (5.3.5): $\Xi$, $v_*$ and $v'_*$ ---
tend to 0, we get
that $c_i$ also tend to 0. But $c_i\in\n F_q$, i.e. they are almost
all 0. $\square$
\medskip
{\bf Lemma 5.3.7.} The above pairing is perfect.
\medskip
{\bf Proof.} We have an isomorphism (here $M\{T\}=M\otimes_{\p[T]}\p\{T\}$ with
the natural action of $\tau$ (see [G], Definition 5.9.11))
$$\alpha: \Hom_{\p[T,\tau]}(M,Z_1)\to\Hom_{\w}(M\{T\}^\tau,\w)\eqno{(5.3.8)}$$ defined as the composition of the maps
$$\Hom_{\p[T,\tau]}(M,Z_1)=\Hom_{\p[T]}(M,Z_1)^\tau\overset{\beta'}\to{\to}\Hom_{\p\{T\}}(M\{T\},\p\{T\})^\tau$$
$$\overset{\gamma}\to{\to} \Hom_{\w}(M\{T\}^\tau,\w)$$ where
$\beta: \Hom_{\p[T]}(M,Z_1)\to \Hom_{\p\{T\}}(M\{T\},\p\{T\} )$ is the natural map and
$\beta'$
is the restriction of $\beta$ to $\tau$-invariant elements. Using the Anderson's
criterion of uniformizability
of $M$ (see, for example, [G], 5.9.14.3 and 5.9.13) we get immediately that both
$\gamma$, $\beta$, and hence
$\beta'$, and hence $\alpha$ are isomorphisms. Further, let us consider a homomorphism
$$i: \Hom_{\p[T,\tau]}(M',Z_1)\to M\{T\}^\tau\eqno{(5.3.9)}$$ defined as follows. Let
$\vf'$, $f'_*$, $v'_*$ be as above. We set
$$i(\vf')=\Xi{v'}^t_*f_*\in M\underset{\p[T]}\to{\otimes}\p[[T]]$$ Since $\Xi\in \p
\{T\}$, we get that $\Xi{v'}^t_*f_*\in M\{T\}$. A simple calculation
(like in the Lemma 5.3.6, but simpler) shows that $i(\vf')$ is $\tau$-invariant, hence
$i$ really defines a map from $\Hom_{\p[T,\tau]}(M',Z_1)$ to $M\{T\}^\tau$. Obviously it
is an inclusion. Let us prove that $i$ is surjective. Really, let $c_*\in (Z_1)^r$ be a
column
vector such that $c_*^tf_*\in M\{T\}^\tau$. An analog of the above calculation shows that
if we define $\vf'$ by the formula $\vf'(f'_*)=\Xi^{-1}c_*$ then $\vf'\in
\Hom_{\p[T,\tau]}(M',Z_1)$, and $i(\vf')=c_*^tf_*\in M\{T\}^\tau$. Finally, the
combination of isomorphisms (5.3.8) and (5.3.9) corresponds to the
pairing (5.3.5). $\square$
\medskip
{\bf 5.4. Step 2 -- End of the proof of Theorem 5.} It is easy to see that the converse of the Corollary 5.1.1 (taking into consideration
Proposition 3.2) is also true, i.e. in order to prove Theorem 5 it is sufficient
to prove that a Siegel matrix of $M'$ is $-Z^t$ where $Z$ is a Siegel matrix of $M$. Let
us consider a basis $l_1,...,l_r$ of $L_T(M)$ and for each $l_i$ we consider the
corresponding (under identification 5.3.1) $\vf_i\in \Hom_{\p[T,\tau]}(M,Z_1)$. Let
$\Psi$ be the scattering matrix of $M$ ([A], p. 486) with respect to the bases $l_1,...,l_r$,
$f_1,...,f_r$, and we denote $\vf_i(f_*)$ by $v_{i*}$ (notations of 5.3.2).
\medskip
{\bf Lemma 5.4.1.} $v_{i*}$ is the $i$-th column of $\Psi$ ($Z_1$ is identified with $\p\{T\}$, see the proof).
\medskip
{\bf Proof.} Follows from the definitions. Recall that $\hbox{\bf K}=\n F_q((1/T))$. The isomorphism 5.3.1 is the composition of
2 isomorphisms $i_1: L_T(M) \to \Hom_{\w}^c(\hbox{\bf K}/\w, E)$ ([G], 5.9.19) and $i_2:
\Hom_{\w}^c(\hbox{\bf K}/\w, E) \to \Hom_{\p[T,\tau]}(M,\Hom^c(\hbox{\bf K}/\w, \p)$
([G], 5.9.24; recall that $Z_1=\Hom^c(\hbox{\bf K}/\w, \p)$). For $l_i\in L_T(M)$ we have
$(i_1(l_i))(T^{-k})=\exp(\theta^{-k}l_i)$ ([G], line above the lemma 5.9.18) and
$$((i_2\circ i_1(l_i))(f_j))(T^{-k})=<f_j,\exp(\theta^{-k}l_i)>$$ ([G], two lines above the
lemma 5.9.24). Using the identification of $Z_1$ and $\p\{T\}$ ([G], p. 168, lines 3 - 4)
and the definition of $\Psi$ ([A], p. 486, first formula of 3.2) we get immediately the
lemma. $\square$
\medskip
Let $l'_1,...,l'_r$ be a basis of $L_T(M')$ which is dual to a basis $l_1,...,l_r$ of
$L_T(M)$ with respect to the pairing 5.3.5.
\medskip
{\bf Lemma 5.4.2.} The scattering matrix of $M'$ with respect to the bases
$l'_1,...,l'_r$, $f'_1,...,f'_r$ (denoted by $\Psi'$) is $\Xi^{-1}\Psi^{t-1}$.
\medskip
{\bf Proof.} Follows immediately from 5.4.1 applied to both $M$, $M'$, and formula 5.3.5.
$\square$
\medskip
{\bf Remark 5.4.3.} An alternative proof for the case of pure $M$ (for $some$ basis of $L_T(M')$)
is the following. We denote $\Xi^{-1}\Psi^{t-1}$ by $\Psi_1$. It satisfies
$\Psi_1^{(1)}=(T-\theta) Q^{t-1}\Psi_1$ and other conditions of [A], 3.1. According [A],
Theorem 5, p. 488, there exists a pure uniformizable t-motive $M_1$ with
$\sigma$-structure such that its scattering matrix is $\Psi_1$. Since $\Psi_1$ satisfies
$$\Psi_1^{(1)}=Q'\Psi_1$$
we get that $Q(M_1)=Q'$, i.e. $M_1=M'$. $\square$
\medskip
Let us recall the statement of the crucial proposition 3.3.2 of [A]. Here we consider the
case of those $M$ whose $N$ is not necessarily 0. Let $\Psi$ be a scattering matrix of
$M$. We consider the $(T-\theta)$-Laurent series for $\Psi$ (here $k(M)<0$ is a number,
and $D_{-i}\in M_{r}(\p)$): $$\Psi=\sum_{i=k(M)}^\infty D_{-i}(T-\theta)^i$$ We consider its
negative part $$\Psi^-:=\sum_{i=k(M)}^{-1} D_{-i}(T-\theta)^i$$ as an element of $M_{r}(\p)((T-\theta))/M_{r}(\p)[[T-\theta]]$.
\medskip
We consider the space $(T-\theta)^{k(M)}\p[[T-\theta]]/\p[[T-\theta]]$ as a $\p$-vector space endowed by the action of $\w$, and we denote by $\goth V$ its $r$-th direct sum written as vector columns of length $r$. Obviously $$k(M)=-1 \iff \hbox{ the action of $T$ on $\goth V$ coincides with multiplication by $\theta$} \eqno{(5.4.3a)}$$
We denote the $i$-th column of $\Psi^-$ by $\Psi^-_{i*}$, it belongs to $\goth V$. Following [A], we denote by $\Prin(M)$ (resp. by $\Prin_0(M)$)
the $\p[T]$-linear span (resp. the $\w$-linear span) of all $\Psi^-_{i*}$ in $\goth V$. Finally, we obviously extend the definition of $\Lie_T(M)$, $L_T(M)$ to the case $N\ne 0$; formula 5.1.4 becomes
$$\goth a(Tx)=(\theta+N)(\goth a(x))\eqno{(5.4.3b)}$$
\medskip
{\bf Proposition 3.3.2, [A]} (see also Remark 5.5 below). There exists a $\p[T]$-linear isomorphism $\psi_E: \Lie_T(M)
\to \Prin(M)$ such that its restriction to $L_T(M)\subset \Lie_T(M)$ defines an isomorphism
$L_T(M) \to \Prin_0(M)$ (denoted by $\psi_E$ as well). $\square$
\medskip
{\bf Corollary 5.4.4.} $N=0 \iff k(M)=-1$ (because $N=0 \iff $ the action of $T$ on both $\Lie_T(M)$, $\goth V$ coincides with multiplication by $\theta$, by 5.4.3a). $\square$
\medskip
We return to the case $N=0$.
\medskip
Let us consider the $(T-\theta)$-Laurent series for $\Psi'$ and $\Xi^{-1}$:
$$\Psi'=\sum_{i=k(M')}^\infty D'_{-i}(T-\theta)^i,  \ \ \ \Xi^{-1}=\sum_{i=k(\xi)}^\infty
a_i(T-\theta)^i$$
Since for both $M$, $M'$ we have $N=N'=0$, we get $k(M)=k(M')=-1$. An
elementary calculation shows that $k(\xi)$ is also $-1$. Hence, equality
$\Psi'\Psi^t=\Xi^{-1}$ (Lemma 5.4.2) implies that $D'_{1}D^t_{1}=0$.

Further, there exist $n$ columns of $D_{1}$ which are
$\p$-linerly independent (they are $\psi_E$-images of elements of $L_T(M)$ which form a
$\p$-basis of $\Lie_T(M)$) and all other columns of $D_{1}$ are their linear combinations.
Interchanging columns of $D_{1}$ if necessary we can assume that these columns are the
last $n$
columns. We denote by $D_{12}$ (resp. $D_{11}$ ) the $r\times n$ (resp. $r\times
(r-n)$ ) matrix formed by the last $n$ (resp. the first $r-n$) columns of $D_{1}$. There
exists a matrix $S$ such that $D_{11}=D_{12}S^t$. Again according Proposition 3.3.2,
[A], we have:
$$S \hbox{ is a Siegel matrix of $L(M)$ }\eqno{(5.4.5)}$$
(see also Remark 5.5 below).

Analogous objects are defined for $D'_{1}$. We denote by $D'_{12}$ (resp.
$D'_{11}$) the $r\times n$- (resp. $r\times (r-n)$)-matrix formed by the last $n$
(resp. the first $r-n$) columns of $D'_{1}$. Since
$D'_{1}D^t_{1}=D'_{11}D_{11}^t+D'_{12}D_{12}^t$ we get that
$D'_{12}D_{12}^t+D'_{11}SD_{12}^t=0$. Since $D_{12}^t$ is a $n\times
r$-matrix of rank $n$, it is not a zero-divisor from the right, so $$D'_{12}=
-D'_{11}S\eqno{(5.4.6)}$$ Since the rank of $D'_{1}$ is $r-n$ and $D'_{11}$ is a
$r\times (r-n)$ matrix, (5.4.6) implies that columns of $D'_{11}$ are linearly independent,
and by (5.4.6) and Proposition 3.3.2, [A] we get that $-S$ is a Siegel matrix of
$M'$. $\square$
\medskip
{\bf Remark 5.5.} Since the notations of [A] differ from the ones of the present
paper, for the reader's convenience we give here a sketch of the proof for the case $N=0$
of two  crucial facts: Corollary 5.4.4 and 5.4.5 ([A], Theorem 3.3.2).

Let $\alpha: \Lie(M)\to E(M)$ be a linear isomorphism which is the first term of the
series for $\exp: \Lie(M)\to E(M)$, and let $l\in \Lie(M)$, $f \in M$ be arbitrary. We
consider the $(T-\theta)$-Laurent series $\sum_{i=k}^\infty b_i(T-\theta)^i$ of
$\sum_{j=0}^\infty <\exp(\frac{1}{\theta^{j+1}}l),f>T^j$.
\medskip
{\bf Lemma 5.6.} If $N=0$ then $k=-1$, and $b_{-1}=-<\alpha(l),f>$ (this is [A],
3.3.4).
\medskip
{\bf Sketch of the proof.} For $z\in \Lie(M)$ we denote $\exp(z)- \alpha(z)$ by $\ve(z)$,
hence
$\sum_{j=0}^\infty <\exp(\frac{1}{\theta^{j+1}}l),f>T^j=\underline{A} +\underline{E}$, where
$$\underline{A}=\sum_{j=0}^\infty <\alpha(\frac{1}{\theta^{j+1}}l),f>T^j; \ \ \ \ \underline{E}=\sum_{j=0}^\infty <\ve(\frac{1}{\theta^{j+1}}l),f>T^j$$
We consider their $(T-\theta)$-Laurent series:
$$\underline{A}=\sum_{i=k(\underline{A})}^\infty \underline{a}_i(T-\theta)^i; \ \ \ \ \underline{E}=\sum_{i=k(\underline{E})}^\infty \underline{e}_i(T-\theta)^i$$
Since we have $\exp(z)=\sum_{i=0}^{\infty}C_iz^{(i)}$ where $C_0=I_n$ we get that
$\ve(z)=\sum_{i=1}^{\infty}C_iz^{(i)}$. This means that for large $j$ the element
$\ve(\frac{1}{\theta^{j+1}}l)$ is small, and hence $k(\underline{E})=0$, because finitely many
terms having small $j$ do not contribute to the pole of the $(T-\theta)$-Laurent series
of $\underline{E}$ (the reader can prove easily the exact estimations himself, or to look [A],
p. 491). Since $\alpha$ is $\p$-linear, equality $\sum_{j=0}^\infty
\frac{1}{\theta^{j+1}}T^j= -(T-\theta)^{-1}$ implies that $k(\underline{A})=-1$ and $\underline{a}_{\ -1}=-<\alpha(l),f>$ (and other $\underline{a}_i=0$), hence the lemma. $\square$
\medskip
This lemma obviously implies Corollary 5.4.4. Further, elements $f_1,...,f_r$ generate
the $\p$-space $M/\tau M$, because multiplication by $T$ on $M/\tau M$ coincides with
multiplication by $\theta$, hence the fact that $f_1,...,f_r$ \ \ $\p[T]$-generate $M/\tau
M$ implies that they $\p$-generate $M/\tau M$.

Let $l_1,...,l_n$ form a $\p$-basis of $\Lie(M)$ (here we identify $\Lie_T(M)$ and $\Lie(M)$ via $\goth a$). Since the pairing $<*,*>$ between
$E(M)$ and $M/\tau M$ is non-degenerate and $\alpha$ is an isomorphism, we get that
columns $<\alpha(l_1),f_*>,..., <\alpha(l_n),f_*>$ are linearly independent. Again since
$\alpha$ is an isomorphism and the pairing with $f_*$ is linear, we get that
$$(<\alpha(l_{n+1}),f_*> \ ... \ <\alpha(l_r),f_*>)=(<\alpha(l_{1}),f_*> \ ... \
<\alpha(l_n),f_*>)Z^t $$ Applying the lemma 5.6 to this formula we get immediately 5.4.5.
\medskip
{\bf 6. Tensor products.}
\medskip
There exists an analog of the Theorem 5 for the case of tensor
products of
t-motives. It describes the lattice $L(M_1\otimes M_2)$ in terms of $L(M_1)$,
$L(M_2)$. This is a theorem of Anderson; it is formulated in [P], end of page 3, but its
proof was not published. We recall its statement for the case of arbitrary $N\ne 0$, and
we give its proof for the case $N=0$ (case of arbitrary $N$ can be obtained easily using
the same ideas).

Let $M$ be an uniformizable t-motive whose $N$ is not necessarily 0. Since $N$ is nilpotent, formula 5.4.3b shows that $\Lie_T(M)$ is a $\p[[T-\theta]]$-module. There exists an epimorphism of
$\p[[T-\theta]]$-modules
$$L_T(M)\underset{\w}\to{\otimes}\p[[T-\theta]]\to \Lie_T(M)$$
whose kernel $\goth q=\goth q(M)$ carries information on the pair $(L(M), \Lie(M))$.
\medskip
{\bf Theorem 6} (Anderson). Let $M$, $\bar M$ be any two uniformizable abelian
t-motives. Then $$\goth q(M\otimes \bar M)=\goth q(M)\underset{\p[[T-\theta]]}
\to{\otimes}\goth q(\bar M)\eqno{(6.1)}$$
\medskip
{\bf Remark 6A.} $M\otimes \bar M$ is a uniformizable t-motive ([G], Corollary 5.9.38).
\medskip
{\bf Proof of Theorem 6 (case $N=0$).} We define notations for $M$, and all notations for $\bar M$
will carry bar. Let $e_i$ and $Z$ be from the beginning of Section 3. We denote $\goth a^{-1}(e_i)\in \Lie_T(M)$ by $e_i$ (there is no possibility of confusion). So, $\{e_i\}$ is a
$\p[[T-\theta]]$-basis of $L_T(M)\underset{\w}\to{\otimes}\p[[T-\theta]]$. Elements
$b_i:=(T-\theta)e_i$, $i=1,...,n$ and $b_{n+i}:=e_{n+i}-\sum_{j=1}^n z_{ij}e_j$,
$i=1,...,r-n$ form a $\p[[T-\theta]]$-basis of $\goth q$. We need a
\medskip
{\bf Lemma 6.2.} $\Psi(M\otimes \bar M)=\Psi(M)\otimes\Psi(\bar M)$ where $\Psi(M)$
(resp. $\Psi(\bar M)$; $\Psi(M\otimes \bar M)$) is taken with respect to bases $e_*$ of
$L_T(M)$, $f_*$ of $M_{\p[T]}$ (resp. $\bar e_*$ of $L_T(\bar M)$, $\bar f_*$ of $\bar
M_{\p[T]}$; $e_*\otimes \bar e_*$ of $L_T(M\otimes \bar M)$, $f_*\otimes \bar f_*$ of
$(M\otimes \bar M)_{\p[T]}$) (see the proof for the notations).
\medskip
{\bf Proof.} We consider a map $$\alpha:\Hom_{\p[T]}(M,Z_1)^\tau
\underset{\w}\to{\otimes} \Hom_{\p[T]}(\bar M,Z_1)^\tau\to \Hom_{\p[T]}(M\otimes \bar
M,Z_1)^\tau$$ defined as follows: for $\vf\in \Hom_{\p[T]}(M,Z_1)^\tau$, $\bar \vf\in
\Hom_{\p[T]}(\bar M,Z_1)^\tau$ we let $[\alpha(\vf \otimes \bar \vf)](f\otimes \bar
f)=\vf(f)\cdot\bar \vf(\bar f)$ (it is obvious that $\alpha(\vf \otimes \bar \vf)$ is
$\tau$-stable). Since $e_1, ..., e_r$ (resp. $\bar e_1, ..., \bar e_{\bar r}$) is a basis
of $\Hom_{\p[T]}(M,Z_1)^\tau$ (resp. $\Hom_{\p[T]}(\bar M,Z_1)^\tau$; we identify $L_T(M)$,
resp. $L_T(\bar M)$ with $\Hom_{\p[T]}(M,Z_1)^\tau$ (resp. $\Hom_{\p[T]}(\bar M,Z_1)^\tau$)
we get (using Lemma 5.4.1) that $\Psi(M)$, $\Psi(\bar M)$ are non-degenerate. Since their
product is also non-degenerate, we get $\alpha(e_i\otimes \bar e_{\bar i})$ are linearly
independent and hence a basis of $\Hom_{\p[T]}(M\otimes \bar M,Z_1)^\tau$. Applying once
again Lemma 5.4.1 we get the lemma. $\square$
\medskip
If $A$, $B$ are two matrices then columns of $A\otimes B$ are indexed by pairs $(k,l)$
where $k$ (resp. $l$) is the number of a column of $A$ (resp. $B$). We denote by $A_k$,
$B_l$, $A\otimes B_{(k,l)}$ the respective columns. Obviosly we have: $A\otimes
B_{(k,l)}=A_k\otimes B_l$ (tensor product of column matrices).
\medskip
Let us prove that for $i=1,...,r-n$, $\bar i=1,...,\bar r-\bar n$ the element
$b_{n+i}\otimes \bar b_{\bar n+\bar i}\in \goth q(M\otimes \bar M)$. According [A],
Proposition 3.3.2, it is sufficient to prove that the corresponding linear combination
(see 6.3 below) of the columns of the matrix $\Psi^-_{M\otimes \bar M}$ is 0. Since
$$b_{n+i}\otimes \bar b_{\bar n+\bar i}=\sum_{j,\bar j}z_{ij}\bar z_{\bar i\bar
j}e_j\otimes \bar e_{\bar j}- \sum_{j}z_{ij}e_j\otimes \bar e_{\bar n+\bar i}
-\sum_{\bar j}\bar z_{\bar i\bar j}e_{n+i}\otimes \bar e_{\bar j}+ e_{n+i}\otimes \bar
e_{\bar n+\bar i}$$ we get the explicit form of this linear combination: it is sufficient
to prove that for all $i$, $\bar i$ we have
$$\sum_{j,\bar j}z_{ij}\bar z_{\bar i\bar j} (\Psi^-_{M\otimes \bar M})_{(j,\bar j)} -
\sum_{j}z_{ij}(\Psi^-_{M\otimes \bar M})_{(j, \bar n+\bar i)} $$ $$ -\sum_{\bar j}\bar
z_{\bar i\bar j}(\Psi^-_{M\otimes \bar M})_{(n+i,\bar j)} +(\Psi^-_{M\otimes \bar
M})_{(n+i,\bar n+\bar i)} =0\eqno{(6.3)}$$
Further, 6.2 implies that $$(\Psi^-_{M\otimes \bar M})_{(k,\bar
k)}=\frac{A_{-1,k}\otimes \bar A_{-1,\bar k}}{(T- \theta)^2}+\frac{A_{-1,k}\otimes \bar
A_{0,\bar k}+A_{0,k}\otimes \bar A_{-1,\bar k}}{T- \theta}$$ hence 6.3 becomes
$$\sum_{j,\bar j}z_{ij}\bar z_{\bar i\bar j} (\frac{A_{-1,j}\otimes \bar A_{-1,\bar
j}}{(T- \theta)^2}+\frac{A_{-1,j}\otimes \bar A_{0,\bar j}+A_{0,j}\otimes \bar A_{-1,\bar
j}}{T- \theta}) $$ $$- \sum_{j}z_{ij}(\frac{A_{-1,j}\otimes \bar A_{-1,\bar n+\bar i}}{(T-
\theta)^2}+\frac{A_{-1,j}\otimes \bar A_{0,\bar n+\bar i}+A_{0,j}\otimes \bar A_{-1,\bar
n+\bar i}}{T- \theta}) $$ $$-\sum_{\bar j}\bar z_{\bar i\bar j}(\frac{A_{-1,n+i}\otimes
\bar A_{-1,\bar j}}{(T- \theta)^2}+\frac{A_{-1,n+i}\otimes \bar A_{0,\bar
j}+A_{0,n+i}\otimes \bar A_{-1,\bar j}}{T- \theta}) $$ $$+\frac{A_{-1,n+i}\otimes \bar
A_{-1,\bar n+\bar i}}{(T- \theta)^2}+\frac{A_{-1,n+i}\otimes \bar A_{0,\bar n+\bar
i}+A_{0,n+i}\otimes \bar A_{-1,\bar n+\bar i}}{T- \theta}=0\eqno{(6.4)}$$
It is easy to see that 6.4 follows immediately from the equalities
$$A_{-1,n+i}=\sum_j z_{ij}A_{-1,j}\eqno{(6.5)}$$ $$ \bar A_{-1,\bar n+\bar
i}=\sum_{\bar j} \bar z_{\bar i\bar j}\bar A_{-1,\bar j}$$ For example, the left hand
side of (6.4) has 2 terms containing $\bar A_{0,\bar j}$ (in the middle of the first
and the third lines of (6.4)). Multiplying (6.5) by $\bar z_{\bar i\bar j}\bar
A_{0,\bar j}$ we get that the sum of these 2 terms of (6.4) is 0. For other pairs of
terms of (6.4) the situation is the same.
\medskip
The proof that for $i=1,...,r-n$, $\bar i=1,...,\bar n$ the element $b_{n+i}\otimes \bar
b_{\bar i}\in \goth q(M\otimes \bar M)$ is analogous but simpler. We have
$$b_{n+i}\otimes \bar b_{\bar i}=(T-\theta) (-\sum_{j} z_{ij} e_j\otimes \bar e_{\bar i}
+ e_{n+i}\otimes \bar e_{\bar i})$$ The analog of (6.3)) is $$(T-\theta) (-
\sum_{j}z_{ij}(\Psi^-_{M\otimes \bar M})_{(j, \bar i)}  +(\Psi^-_{M\otimes \bar
M})_{(n+i,\bar i)}) =0$$ and the analog of (6.4)) is $$-
\sum_{j}z_{ij}\frac{A_{-1,j}\otimes \bar A_{-1,\bar i}}{T- \theta}
+\frac{A_{-1,n+i}\otimes \bar A_{-1,\bar i}}{T- \theta}=0$$ This equality follows
immediately from (6.5).
\medskip
Finally, elements $b_{i}\otimes \bar b_{\bar i}$ ($i=1,...,n$, $\bar i=1,...,\bar n$)
obviously belong to $ \goth q(M\otimes \bar M)$.
\medskip
So, we proved that $\goth q(M)\underset{\p[[T-\theta]]}\to{\otimes}\goth q(\bar M)
\subset \goth q(M\otimes \bar M)$. Since the $\p$-codimension of both subspaces in
$L_T(M)\underset{\w}\to{\otimes} L_T(\bar M)\underset{\w}\to{\otimes}
\p[[T-\theta]]$ is $n\bar n$, they are equal. $\square$
\medskip
{\bf 7. Self-dual t-motives.}
\medskip
{\bf Case $\w=\n F_q[T]$.} A uniformizable t-motive $M$ is called self-dual if there exists an isogeny $\alpha: M \to M'$. It defines an $\w$-valued, $\w$-bilinear form $<*,*>_\alpha$ on $L_T(M')$ as follows:
$$<\vf_1, \vf_2>_\alpha=<L_T(\alpha)(\vf_1), \vf_2>_M$$
5.2.11 implies that if $\alpha'=-\alpha$ (resp. $\alpha'=\alpha$) then
$<*,*>_\alpha$ is skew symmetric (resp. symmetric). $M$ is called positively (resp. negatively) self-dual if $\alpha$ satisfies $\alpha'=\alpha$ (resp. $\alpha'=-\alpha$). Hence, we have an
\medskip
{\bf Analogy 7a.} The number field case analog of a pair: $\{$negatively self-dual t-motive of rank $2n$, dimension $n$; negative $\alpha: M \to M'\}$ is a (generic) abelian variety of dimension $n$ with a fixed polarization form.
\medskip
For example, like in the number field case, we can define the Rosati involution $I_\alpha$ on  $\End_0(M):=\End(M)\otimes \n F_q(T)$ by the same formula $I_\alpha(\vf)=\alpha^{-1}\circ \vf'\circ \alpha$.
\medskip
Further, we have a
\medskip
{\bf Conjecture 7b.} The dimension of the moduly variety of negatively self-dual t-motives (if it exists) is $n(n+1)/2$.
\medskip
{\bf Examples.} Let $e_*$ be from 1.9, and let $M=M(A)$ given by the equation (here $A\in M_n(\p)$ is $\goth A_1$ of 1.9.1)
$$Te_*= \theta e_*+ A\tau e_* + \tau^2 e_*\eqno{(7.1)}$$
be a t-motive of dimension $n$ and rank $2n$. Elements $f_i=e_i$,
$f_{n+i}=\tau e_i$ $(i=1,...,n)$ form a $\p[T]$-basis of $M$. We have
(see, for example, Section 11): $M'$ is given by the equation $$Te'_*=
\theta e'_*- A^t\tau e'_* + \tau^2 e'_*$$ and if we define $$f'_i=\tau
e'_i, \ \ f'_{n+i}= e'_i\eqno{(7.2)}$$ then bases $f_*$, $f'_*$ are dual in the
meaning of Lemma 1.10.

Let $\alpha: M \to M'$ be given by the formula $\alpha(e_*)=De'_*$
where $D\in M_n(\p)$ (we impose this essential restriction only in order to simplify exposition. In the general case $D\in M_n(\p[\tau])$, $D_f\in M_{2n}(\p[T])$, $D_f$ from 7.4). Condition that $\alpha$ is a $\p[T,\tau]$-map is
equivalent to
$$D^{(2)}=D, \ \ AD^{(1)}=-DA^t\eqno{(7.3)}$$
Further, we have $$\alpha(f_*)=D_ff'_*\eqno{(7.4)}$$ where
$D_f=\left(\matrix 0 & D\\ D^{(1)}&0 \endmatrix \right)$, hence
$$\alpha'=\pm \alpha\iff
D_f^t=\pm D_f\iff D^{(1)}=\pm D^t\eqno{(7.5)}$$
Let us fix $\varepsilon_0\in \n F_{q^2}$ satisfying $\varepsilon_0^{q-1}=-1$. Then $D=\varepsilon_0I_n$ satisfies 7.5 with the sign minus, and the set of $A$ satisfying 7.3 with this $D$ is the set of symmetric matrices. This justifies 7b, because the set of $A_1\in  M_n(\p)$ such that $M(A)=M(A_1)$ is conjecturally discrete.

For $D=I_n$ the sign in 7.5 is plus and hence a skew symmetric $A$ defines a positively self-dual $M(A)$.
\medskip
{\bf Remark 7.6.} The below statements are conjectures based on arguments similar to the ones which justify the below Conjecture 11.8.3. Since they are of secondary importance, we do not give any details of justification here.
\medskip
{\bf 7.6.1. Conjecture.} If $n\ge 3$ then for a
generic skew symmetric $A$ we have: $\End(M(A))=\w$.
\medskip
{\bf 7.6.2. Corollary.} Conjecture 7.6.1 implies that the "minimal" $\alpha: M \to M'$ is defined uniquely up to an element of
$\n F_q^*$, and hence the symmetric pairing $<*,*>_\alpha$ is also defined uniquely up to an element of $\n F_q^*$.
\medskip
{\bf 7.6.3. Conjecture.} If $n=2$, $\alpha'=\alpha$ then $\End (M)$ is strictly larger than $\w$.
\medskip
Other examples of a self-dual t-motive are $M\oplus M'$ where $M$ is any t-motive, but they do not give interesting examples of pairings.
\medskip
{\bf 7.6.4. Conjecture.} There exist other (distinct from the ones defined by 7.1) self-dual t-motives $M$ having $\End (M)=\w$ (we can use a version of standard t-motives of Section 11).
\medskip
{\bf Example 7.7.} Case $A=0$, $D=I_n$.
\medskip
In this case we can find explicitly the matrix of the symmetric form
$<*,*>_\alpha$ in some basis of $L_T(M')$. Let $\goth C_2$ be the
Carlitz module over the field $\n F_{q^2}$ considered as a rank 2
Drinfeld module over $\n F_{q}$ given by the equation $$Te=\theta e
+ \tau^2e$$ We have $M=\goth C_2^{\oplus n}$. Let $\goth T_T(\goth
C_2)$ be the convergent $T$-Tate module of $\goth C_2$, i.e. the set
of elements $\{z_i\}\in E(\goth C_2)=\p$ $(i\ge -1, \ \ z_{-1}=0)$ such that
$$\hbox{ $T z_i=z_{i-1}$ for $i \ge 0$ (i.e. $z_i^{q^2}+\theta
z_i=z_{i-1}$) and $z_i\to 0$}$$ It is a free 1-dimensional module over
$\n F_{q^2}[T]$. We choose and fix its generator; its $\{z_i\}$ satisfy (like in 5.3.4)
$|z_0|>|z_i| \ \ \forall i>0$. We denote $\sum_{k=0}^\infty z_kT^k$ by $\goth Z$.

Let $c$ be a fixed element of $\n F_{q^2}-\n F_{q}$. Formulas (5.3.3)
show that the following elements $\vf_i$, $\vf'_i$ ($i=1, ... , 2n$)
form bases of $L(M)$, $L(M')$ respectively ($j=1, ..., n$; clearly that thanks to 7.2
we have $\vf'_i(f'_j)=\vf_i(f_{n+j})$, $n+j$ mod $2n$):
$$i\le n: \ \ \vf_i(f_j)=\goth Z\delta_i^j, \ \
\vf_i(f_{n+j})=\goth Z^{(1)}\delta_i^j $$

$$i>n: \ \ \vf_i(f_j)=c\goth Z\delta_{i-n}^j, \ \
\vf_i(f_{n+j})=c^q\goth Z^{(1)}\delta_{i-n}^j $$

$$i\le n: \ \ \vf'_i(f'_j)=\goth Z^{(1)}\delta_i^j, \ \
\vf'_i(f'_{n+j})=\goth Z\delta_i^j $$

$$i>n: \ \ \vf'_i(f'_j)=c^q\goth Z^{(1)}\delta_{i-n}^j, \
\ \vf'_i(f'_{n+j})=c\goth Z\delta_{i-n}^j $$
(by the way, it is clear that the same relation between elements of $\goth T_T(M)$ and
$\Hom_{\p[T,\tau]}(M,Z_1)$ holds for all $M$). Formula 7.4 shows
that $\alpha'(\vf'_i)=\vf_{i+n}$, where $i+n \mod 2n$. Let us denote
$\Xi\cdot \goth Z\cdot \goth Z^{(1)}\in\n F_q^*$
by $\gamma$. The above definitions and formulas show that the matrix
of $<*,*>_\alpha$ in the basis $\vf_{1}, \vf_{n+1}, ... ,\vf_{n},
\vf_{2n}$ consists of $n$ \ \ $(2\times 2)$-blocks (trace and norm of $\n
F_{q^2}/\n F_q$)
$$\gamma\left(\matrix \tr(1) & \tr(c)\\ \tr(c)&\tr(N(c))  \endmatrix
\right)=\gamma\left(\matrix 2 & c+c^q \\ c+c^q & 2c^{q+1} \endmatrix
\right)$$ The determinant of
this block is $-(c-c^q)^2\gamma^2$; it belongs to $\n F_q^{*2}\iff
q\equiv 3 \mod 4$ or $q$ is even. Since we have $n$ blocks, we have:

$$\hbox{det $<*,*>_\alpha\not\in\n F_q^{*2}\iff q\equiv 1 \mod 4$ and
$n$ is odd.}$$
\medskip
{\bf Remark 7.8 (Jorge Morales).}
There is a theorem of Harder (see
e.g. W. Scharlau, "Quadratic and Hermitian forms", Springer-Verlag,
Berlin, 1985, Chapter 6, Theorem 3.3) that states that a unimodular form over
$k[X]$ \ \ --- \ \ $k$ being any field of characteristic not 2 --- is the
extension of a form over $k$, i.e. there is a basis in which all the
entries of the associated symmetric matrix are constant. This means
that the classification of the above quadratic forms over $\n F_q[T]$
($q$ odd) is very simple.
\medskip
{\bf Remark 7.9.} Let $M$ be a t-motive which is both negatively and positively self-dual. There is a natural idea 7.9.2 to define an analog of
Hodge structure on $M$. Nevertheless, this idea fails. Namely, the exact sequence
$$0\to \Ker \vf \to L(M)\otimes \p \overset{\vf}\to{\to} \Lie(M)\to 0$$
is the functional field analog of an exact sequence for an abelian variety $A$:
$$ 0 \to H^{0,-1}(A) \to H^{-1}(A) \to (H^{1,0})^*(A) \to 0 $$
Hence, we can define $H^{0,-1}(M):= \Ker \vf$, and the problem is to define an analog of $H^{-1,0}(M)$.

Let us fix a negative isogeny $\alpha: M\to M'$, and let us extend the skew form $<*,*>_\alpha$ to $L(M)\otimes \p$ by $\p$-linearity. It is easy to check that $\Ker \vf$ is isotropic with respect to this form (there is an analogy with the number field case). Let us consider the following elementary lemma of linear algebra:
\medskip
{\bf Lemma 7.9.1.} Let $W$ be a vector space of dimension $2n$ over a field of characteristic $\ne 2$, $B^+$ (resp. $B^-$) a symmetric (resp. skew symmetric) non-degenerate bilinear form on $W$, and $W_0\subset W$ a subspace of dimension $n$ which is isotropic with respect to both $B^+$, $B^-$. Then almost always there exists the only $W_1\subset W$ of dimension $n$ having properties:
$$W_0\cap W_1=0; \ \ \ W_1 \hbox{  is isotropic with respect to both } B^+, \  B^-$$
where almost always means that entries of the matrices of $B^+$, $B^-$ in a basis of $W$ must not satisfy (at least one of) polynomial relations. $\square$
\medskip
If $\End_0(M)\ne \n F_q(T)$ and the action of $I_\alpha$ on $\End_0(M)$ is not identical, then there exists a positive isogeny $\beta: M\to M'$ and hence the symmetric form $<*,*>_\beta$ on $L(M)\otimes \p$. $\Ker \vf$ is isotropic with respect to  $<*,*>_\beta$. Let us fix $\beta$.
\medskip
{\bf Idea 7.9.2.} To apply Lemma 7.9.1 to this situation ($W=L(M)\otimes \p$, $W_0=\Ker \vf$, $B^+=<*,*>_\beta$, $B^-=<*,*>_\alpha$) in order to get a canonical subspace of $L(M)\otimes \p$ which is complementary to $\Ker \vf$ and hence can be considered as an analog of $H^{-1,0}(M)$.
\medskip
Clearly there is no complete analogy with the number field case. But the situation is even worse:
\medskip
{\bf Proposition 7.9.3.} For all $M$, $\alpha$, $\beta$ the "almost always" condition of Lemma 7.9.1 is not satisfied. $\square$
\medskip
{\bf 8. Relations between lattices and t-motives.}
\medskip
We have\footnotemark \footnotetext{I am grateful to Urs
Hartl
who indicated me this reference.}
\medskip
{\bf Theorem 8.1.} ([H],  Theorem 3.2). The dimension of the moduli set of pure t-motives of
dimension $n$ and rank $r$ is $n(r-n)$. $\square$
\medskip
{\bf Remark.} A tuple
$(e_1,...,e_r)$ of
integers entering in the statement of this theorem in [H] is
(0,...,0,1,...,1) with 0 repeated
$r-n$ times and 1 repeated $n$ times for the case under
consideration.
\medskip
Since this number $n(r-n)$ is equal to the dimension of the set of lattices of rank $r$
and dimension $n$, we can state an
\medskip
{\bf Open question 8.2.} Let $r$, $n$ be given. Let us consider the lattice map from the set of the pure uniformizable
t-motives of rank $r$ and dimension $n$ to the set
of lattices of rank $r$ and dimension $n$. Is it true that its image is open and the fibre at a generic point is discrete? If yes, what is the fibre?
\medskip
{\bf Remark.} Results of [GL17] give some evidence that for the case $r=2n$ in a "neighborhood" of the $n$-th power of the rank 2 Carlitz module the fibre consists of 1 point.
\medskip
Theorem 5 implies that for $n=r-1$ the answer to 8.2 is yes (the below Proposition 11.8.5 shows that most likely the condition of purity is essential):
\medskip
{\bf Corollary 8.4.} All pure t-motives of
dimension $r-1$ and rank $r$ having $N=0$ are uniformizable. There is a 1 -- 1 functorial
correspondence between pure t-motives of dimension $r-1$ and rank $r$ having $N=0$ ($r\ge
2$), and lattices of rank $r$ in $\p^{r-1}$ having dual.
\medskip
{\bf Proof.} Let $L$ be a lattice of rank $r$ in $\p^{r-1}$ having dual $L'$. There
exists the only Drinfeld module $M'$ such that $L(M')=L'$, and let $M$ be its dual.
Theorem 5 implies that $L(M)=L$. If there exists another pure t-motive $M_1$ of
dimension $r-1$ and rank $r$ having $N=0$ such that $L(M_1)=L$ then by Corollary 10.4 (its proof is logically independent: there is no vicious circle) the
dual $M'_1$ is a Drinfeld module, according Theorem 5 it satisfies $L(M'_1)=L'$, hence
$M'_1=M'$ and hence $M_1=M$. $\square$
\medskip
{\bf Remark 8.5.} Recall that lattices of rank $r$ in $\p^{r-1}$ having dual are
described in 3.5 (formulas 3.6, 3.7).
We see that for the case $n=r-1$, $N=0$ purity implies uniformizability. We have
\medskip
{\bf Question 8.5a.} Do exist non-uniformizable t-motives having $n=r-1$, $N=0$?
\medskip
{\bf Question 8.5b.} Do exist uniformizable t-motives having $n=r-1$, $N=0$ such that its lattice has no dual? (Clearly this is a subquestion of 8.2).
\medskip
{\bf Remark 8.6.} Clearly for any $r$, $n$ we have: if a lattice $L$ of rank $r$ and dimension $n$ has no dual then $L\ne L(M)$ for any pure uniformizable $M$. I do not know whether Theorem 6 (which is an analog of Theorem 5 for another tensor operation) imposes a more strong similar restriction on the property of $L$ to be the $L(M)$ of some pure uniformizable $M$, or not.
\medskip
Further, for any uniformizable t-motive $M$ we have a
\medskip
{\bf Corollary 8.7.} If the dual of $(L(M), \Lie(M))$ does not exist then the
dual of $M$ does not exist. Example: the Carlitz module.
\medskip

{\bf 9. Main theorem in terms of Hodge-Pink structure.}
\medskip
Let us consider a version of a special case of the general definition of  Hodge-Pink structure ([P], 0.2; 9.1).
\medskip
{\bf Definition.} A Hodge-Pink structure of constant weight and complete dimension is a pair $\underline{H}=(H,\goth q_H)$ where $H$ is a free finite dimensional $\w$-module and $\goth q_H$ is a $\p[[T-\theta]]$-lattice in $H\underset{\w}\to{\otimes}\p[[T-\theta]]$ such that the dimension of $\goth q_H$ over $\p[[T-\theta]]$ is equal to the dimension of $H$ over $\w$ (condition of complete dimension).
\medskip
Let $\vf: L \hookrightarrow \p^n$ be a lattice. It defines a Hodge-Pink structure $\underline{H}=\underline{H}(L)$ of constant weight and complete dimension. Firstly, instead of a $\n F_q[\theta]$-module $L$ we consider an isomorphic $\w$-module $H$ formally defined by the property $H\underset{\w}\to{\otimes}\n F_q[\theta]=L$ where the map $\w \to \n F_q[\theta]$ is $\iota$. We denote the isomorphism $H \to L$ by $\iota$ as well; the composition $\vf \circ \iota: H \to \p^n$ is a map of $\w$-modules where $T\in \w$ acts on $ \p^n$ by multipication by $\theta$. Further, $\vf \circ \iota$ extends to a surjection of $\p[[T-\theta]]$-modules $H\underset{\w}\to{\otimes}\p[[T-\theta]] \to \p^n$ denoted by $\vf \circ \iota$ as well.  Finally, $\goth q_H$ is defined as $\Ker \vf \circ \iota$.

If $M$ is a pure uniformizable t-motive then we associate it a Hodge-Pink structure of constant weight and complete dimension $\underline{H}(M)=\underline{H}(L(M))$.

Let $m=m(\underline{H})$ be the minimal number such that $\goth q_H \supset (T-\theta)^m H\underset{\w}\to{\otimes}\p[[T-\theta]]$. For $\mu\ge m$ we define the $\mu$-dual structure ${\underline{H}'}^{\mu}=({H'}^{\mu}, \goth q_{{H'}^{\mu}})$ as follows:

$${H'}^{\mu}=H^*, \ \  \goth q_{{H'}^{\mu}}=\{\chi\in H^*\underset{\w}\to{\otimes}\p[[T-\theta]]  $$ $$\hbox{ such that }  \forall y\in \goth q_H  \hbox{ we have } \chi(y)\in (T-\theta)^\mu\p[[T-\theta]]\} $$ It is obvious that it is really a Hodge-Pink structure of constant weight and complete dimension.
\medskip
If $\underline{H}=\underline{H}(L)$ for a lattice $L$ then $m=1$ and if $L$ has dual then $${\underline{H}'}^1=\underline{H}(L')\eqno{(9.1)}$$ this is easy to prove.
\medskip
{\bf Remark 9.2.} And if $L$ has no dual? Really, $\underline{H}(L)$ exists even if $L$ does not satisfy Definition 2.1 (b). If $L$ is a lattice having no dual this means that $L'$  does not satisfy Definition 2.1 (b). Nevertheless, equality ${\underline{H}'}^1=\underline{H}(L')$ is meaningful and holds. We are not interested in these lattices because they cannot be lattices of uniformizable t-motives having dual.
\medskip
Proof of the duality theorem for $M$ having $N\ne0$ is given in [GL18].
\medskip
{\bf 10. Duals of pures, and other elementary results. }

\nopagebreak
\medskip
We consider in this section the case of arbitrary $N$ (i.e. not necessarily $N=0$), and $\w=\n F_q[T]$. The definition 1.8 extends to the case of pr\'e-t-motives, and
remarks 1.11 hold for this case.
\medskip
{\bf Lemma 10.2. } Let $M$ be a pr\'e-t-motive, $m=m(M)$ from its
(1.3.1), and $\mu\ge m$.
Then $M'$ --- the $\mu$-dual of $M$ --- exists as a pr\'e-t-motive,
and $m(M')\le\mu$. If
$M'$ is a t-motive then $\dim M'=r\mu - \dim M$ ($r$ is the rank of $M$).
\medskip
{\bf Proof.} We must check that $Q'$ has no denominators, and the
condition (1.3.1). The
module $\tau M$ is a $\p[T]$-submodule of $M$ (because $a \tau x =
\tau a^{1/q} x$ for
$x\in M$),
hence there are $\p[T]$-bases $f_*=(f_1, ... f_r)^t$, $g_*=(g_1, ...
g_r)^t$ of $M$,
$\tau M$ respectively such that $g_i=P_if_i$, where $P_1 | P_2 | ... |
P_r$, $P_i \in
\p[T]$. Condition (1.3.1) means that $\forall i$ \ $(T-\theta)^m f_i
\in \tau M$, i.e.
$P_i|(T-\theta)^m$, i.e. $\forall i$ \ $P_i=(T-\theta)^{m_i}$ where $0
\le m_i \le
m_{i+1} \le m$. There exists a matrix $\goth Q=\{\goth q_{ij}\}\in
M_r(\p[T])$ such that
$$\tau f_i = \sum_{j=1}^r \goth q_{ij}g_j= \sum_{j=1}^r \goth q_{ij} P_j
f_j\eqno{(10.2.1)}$$
Although $\tau$ is not a linear operator, it is easy to see that
$\goth Q\in GL_r(\p[T])$
(really, there exists $C=\{c_{ij}\}\in M_r(\p[T])$ such that $g_i=P_i
f_i=\tau (\sum_{j=1}^r
c_{ij}f_j)$, we have $C^{(1)}\goth Q =I_r$).

We denote the matrix $\diag(P_1, P_2, ... ,P_r)$ by $\goth P$, so
(10.2.1) means that
$$Q= \goth Q \goth P \eqno{(10.2.2)}$$

{\bf Remark 10.2.3.} Since
$\goth Q \goth P \in
GL_r(\p(T))$, we get that the action of $\tau$ on $i_2(M)$ is invertible.
\medskip
It is clear that if $M$ is a t-motive then
$$\dim M = \sum_{j=1}^r m_j\eqno{(10.2.4)}$$
(because $\dim M = \dim_{\p}(M/\tau M)$. Further, (10.2.2) implies that
for $Q'=Q({M'})$
we have
$$Q'=\goth Q^{t-1}\diag((T-\theta)^{\mu-m_1}, ...,
(T-\theta)^{\mu-m_r})\eqno{(10.2.5)}$$
This means that elements of $Q'$ have no denominators. The condition
(1.3.1) for $M'$
follows easily from (10.2.5) (because $\goth Q^{t-1}\in GL_r(\p[T])$),
and the dimension
formula (for the case $M'$ is a t-motive) follows immediately from
(10.2.4) applied to
$M'$. $\square$
\medskip
A definition of a pure t-motive can be found in [G] ((5.5.2),
(5.5.6) of [G] + formula (1.3.1) of the present paper).
\medskip
{\bf Theorem 10.3.} Let $M$ be a pure t-motive and $m=m(M)$ from
(1.3.1). Then (if
$rm-n>0$) its $m$-dual $M'$ exists, and it is pure.
\medskip
{\bf Proof.} The definition of pure ([G], (5.5.2)) is valid for
pr\'e-t-motives. We use
its following
matrix form. We denote $T^{-1}$ by $S$ and for any $C$ we let
$$C^{[i]}=C^{(i-1)}\cdot C^{(i-2)}\cdot ... \cdot C^{(1)}\cdot C$$

{\bf Lemma 10.3.1.} Let $Q\in M_r(\p[T])$ be a matrix such that formula
(1.9.3) defines an
t-motive $M$. Then it is pure iff there exists $C\in
GL_r(\p((S)) \ )$ such that
for some $\goth q$, $s>0$
$$S^\goth q C^{(s)}Q^{[s]}C^{-1}\in GL_r(\p[[S]])$$
i.e. iff $S^\goth q C^{(s)}Q^{[s]}C^{-1}$ is $S$-integer and its inicial
coefficient is
invertible.
\medskip
{\bf Proof.} Elementary matrix calculations. We take $C$ as a matrix
of base change of
$f_*$ to a $\p[[S]]$-basis of $W$ of (5.5.2) of [G]. $\square$
\medskip
{\bf Lemma 10.3.2.} Let $\mu = m$. We have: $M'={M'}^\mu$ of Lemma 10.2 is a pure
pr\'e-t-motive.
\medskip
{\bf Proof.} Let $\goth q$, $s$ and $C$ be from Lemma 10.3.1. We have
$${Q'}^{[s]}=((T-\theta)^{[s]})^\mu Q^{[s]\ t-1}$$
(we use (1.2)). We take $C'=C^{t-1}$. We have
$$S^{s\mu-\goth q} {C'}^{(s)}{Q'}^{[s]}{C'}^{-1}=$$
$$=S^{s\mu-\goth q}C^{(s)\ t-1}Q^{[s]\ t-1}((\frac1S-\theta)^{[s]})^\mu C^t$$
$$=((1-S\theta)^{[s]})^\mu S^{-\goth q}C^{(s)\ t-1}Q^{[s]\ t-1}C^t$$
$$=((1-S\theta)^{[s]})^\mu (S^\goth q C^{(s)}Q^{[s]}C^{-1})^{t-1}$$
We have: $\goth q/s=n/r$ ([G], (5.5.6)), hence $(s\mu-\goth q)/s=(r\mu-n)/r$ and
$s\mu-\goth q>0$. Further,
$((1-S\theta)^{[s]})^\mu\in GL_r(\p[[S]])$, and the result follows from
Lemma 10.3.1. $\square$
\medskip
{\bf Remark.} This result holds also for $\mu>m$.
\medskip
The theorem 10.3 follows from Lemma 10.2, the above lemmas and the
proposition that a pure
pr\'e-t-motive satisfying (1.3.1) is a t-motive ([G], (5.5.6),
(5.5.7)). $\square$
\medskip
{\bf Corollary 10.4.} Let $M$ be a t-motive such that $m=1$,
$n=r-1$. Then $M$
has dual $\iff$ $M$ is pure $\iff$ $M$ is dual to a Drinfeld module.
\medskip
{\bf Proof.} Dimension formula shows that $M'$ (if it exists) is a
Drinfeld module, and they are all pure. $\square$
\medskip
{\bf Example 10.5.} Let $M$ be given by (notations of 1.9.1)
$$\goth A_0=\theta I_2, \ \ \goth A_{1}= \left(\matrix a_{111} & 0
\\ a_{121} & 1 \endmatrix \right), \ \ \goth A_2= \left(\matrix 1 & 0
\\0 & 0 \endmatrix \right)$$
This $M$ has $m=1$, $n=2$, $r=3$, and it is easy to see that it has no dual. Really, for
this $M$ we have (notations of 1.9) $f_1=e_1$, $f_2=\tau e_1$, $f_3=e_2$,
$$Q=\left(\matrix 0&1&0\\T-\theta&-a_{111} & 0
\\ 0&-a_{121} & t-\theta \endmatrix \right), \ \ Q'=\left(\matrix a_{111} & t-\theta
& a_{121} \\ 1 &0&0\\0&0&1\endmatrix \right)$$ The last line of $Q'$ means that $\tau
f_3'=f_3'$. This is a contradiction to the property that $M'_{\p[\tau]}$ is free. It is
possible also to show (Proposition 11.3.4) that $M$ is not pure, and to use 10.4 in order to
prove that it has no dual.
\medskip
Later (Section 11) we shall construct examples of non-pure abelian
t-motives which have dual. Considerations of 11.8 predict that there is enough such t-motives.
\medskip
{\bf Theorem 10.6. } For any t-motive $M$ there exists $\mu_0$ such that for all
$\mu\ge\mu_0$ the object
${M'}^{\mu}$ exists as a t-motive. For these $\mu$ we have
$${M'}^{\mu+1}={M'}^{\mu}\otimes \goth
C\eqno{(10.6.1)}$$

{\bf Proof.} (10.6.1) holds at the level of pr\'e-t-motives, because $Q(\goth
C)=(T-\theta)I_1$. According [G], Lemma 5.4.10 it is sufficient to prove that
${M'}^{\mu}$ is finitely generated as a $\p [\tau]$-module. We shall
use notations of
Lemma 10.2. We take $$\mu_0=1+\hbox{ \{the maximum of the degrees of
entries of $\goth
Q(M)$ as polynomials in $T$\} } $$ $$+ \max(m_k)$$
Let $f'_1, ...f'_r$ be the basis of ${M'}^{\mu}$ over $\p [T]$ dual to
$f_1, ...f_r$. It
is sufficient to prove the
\medskip
{\bf Lemma 10.6.2.} Let $i_0=\mu-\min(m_k)$. Then elements $T^if'_j$,
$i<i_0$, $j=1,
....,r$, generate ${M'}^{\mu}$ as a $\p [\tau]$-module.
\medskip
{\bf Proof of the lemma.} By induction, it is sufficient to show that
for all $\alpha\ge
i_0$ the equation
$$\tau x = (T-\theta)^\alpha f'_j\eqno{(10.6.3)}$$
(equality in ${M'}^{\mu}$) has a solution
$$x=\sum_{k=1}^r C_kf'_k$$
where $C_k \in \p [T]$, $\deg(C_k)< \alpha$. According (10.2.5), the
solution to (10.6.3)
is given by
$$(C_1^{(1)}, ... , C_r^{(1)})=(0,... 0,
(T-\theta)^{\alpha-\mu+m_j},0,... 0)\goth Q^t$$
(the non-0 element of the row matrix is at the $j$-th place).
Unequalities satisfied by
$\mu$ and $\alpha$ show that all $C_k^{(1)}$ are polynomials of degree
$<\alpha$. Since
$c \mapsto c^q$ is surjective on $\p$, we get the desired. $\square$
\medskip
{\bf 10.7. Virtual t-motives.} \footnotemark \footnotetext{This
notion was
indicated me by
Taguchi.} We need two elementary lemmas.

\nopagebreak
\medskip
{\bf Lemma 10.7.0.}\footnotemark \footnotetext{Anderson proved (not published) that the tensor product of any t-motives is also a t-motive.} If $M$ is a t-motive then $M\otimes \goth C$ is also a t-motive.
\medskip
{\bf Proof.} Let $f_{j}$ ($j=1,...,r$) be a $\p[T]$-basis of
$M_{\p[T]}$ and $\goth f$
from 1.10.2, so $f_{j}\otimes \goth f$ is a $\p[T]$-basis of $(M\otimes \goth
C)_{\p[T]}$. It is sufficient to prove that $(M\otimes \goth
C)_{\p[\tau]}$ is finitely
generated. Since $M_{\p[\tau]}$ is finitely generated, it is easy to
see that there
exists $a$ such that elements
$$(T-\theta)^i f_j, \ \ i=0, ... , a, \ \ j=1, ... ,r$$
generate $M_{\p[\tau]}$. This means that $\forall j=1, ... ,r$ there
exist $c_{ijkl}\in
\p$ such that
$$(T-\theta)^{a+1} f_j=\sum_{i=0}^a \sum_{k=0}^\gamma\sum_{l=1}^r
c_{ijkl}(T-\theta)^i
\tau^k f_l\eqno{(10.7.0.1)}$$
where $\gamma$ is a number.

Let us multiply (10.7.0.1) by $(T-\theta)^\gamma$. Taking into
consideration the formula
of the action of $\tau$ on $M\otimes \goth C$ we get that the result
gives us the
following formula in $M\otimes \goth C$:
$$(T-\theta)^{a+\gamma+1} f_j\otimes \goth f=\sum_{i=0}^a
\sum_{k=0}^\gamma\sum_{l=1}^r
c_{ijkl} (T-\theta)^{i+\gamma-k} \tau^k \cdot (f_l\otimes \goth f)\eqno{(10.7.0.2)}$$
This proves that for all $j$ the element $(T-\theta)^{a+\gamma+1} f_j\otimes \goth f$ is
a linear combination of
$$(T-\theta)^i f_l\otimes \goth f, \ \ i=0, ... , a+\gamma, \ \ l=1, ...
,r\eqno{(10.7.0.3)}$$
in $(M\otimes \goth C)_{\p[\tau]}$.
Multiplying (10.7.0.2) by consecutive powers of $T-\theta$ we get by
induction that
elements of 10.7.0.3 generate $(M\otimes \goth C)_{\p[\tau]}$. $\square$
\medskip
{\bf Lemma 10.7.1.} If $M_1\otimes \goth C$ is isomorphic to
$M_2\otimes \goth C$ then
$M_1$ is isomorphic to $M_2$.
\medskip
{\bf Proof.} Let $f_{i*}$ ($i=1,2$) be a $\p[T]$-basis of
$(M_i)_{\p[T]}$, $Q_i$ from
1.9.3, $\alpha: M_1\otimes \goth C\to M_2\otimes \goth C$ an
isomorphism and $C\in
GL_r(\p[T])$ the matrix of $\alpha$ in $f_{1*} \otimes \goth f $,
$f_{2*}\otimes \goth f
$. The matrix of the action of $\tau$ on $M_i\otimes \goth C$ in the
base $f_{i*} \otimes
\goth f $ is $(T-\theta)Q_i$, and the condition that $\alpha$ commutes with
multiplication by $\tau$ is
$$(T-\theta)Q_1 C = C^{(1)}(T-\theta)Q_2$$
Dividing this equality by $T-\theta$ we get that the map $\alpha_0$
from $M_1$ to $M_2$
having the same matrix $C$ in the bases $f_{i*}$, commutes with
$\tau$, i.e. defines an
isomorphism from $M_1$ to $M_2$. $\square$
\medskip
Using Lemma 10.7.1 we can state the following
\medskip
{\bf Definition.} A virtual t-motive is an object $M\otimes \goth C^{\otimes
\mu}$ where $M$ is
a t-motive
and $\mu\in \n Z$, with the standard equivalence relation (here
$\mu_1\ge\mu_2$):
$$M_1\otimes \goth C^{\otimes \mu_1}=M_2\otimes \goth C^{\otimes \mu_2}\iff
M_2=M_1\otimes \goth C^{\otimes (\mu_1-\mu_2)}$$
$$\iff \exists \mu \hbox{ such that $\mu+\mu_1\ge 0$, $\mu+\mu_2\ge 0$
and } M_1\otimes
\goth C^{\otimes (\mu+\mu_1)}=M_2\otimes \goth C^{\otimes (\mu+\mu_2)}$$
Lemma 10.7.1 shows that these conditions are really equivalent.
\medskip
{\bf Corollary 10.7.2.} The $\mu$-dual of a virtual t-motive is
well-defined and
always exists as a virtual t-motive. $\square$
\medskip
{\bf Proposition 10.8.} The following formula is valid at the level of
pr\'e-t-motives: for
any $\mu_1$, $\mu_2$, if ${M_1'}^{\mu_1}$, ${M_2'}^{\mu_2}$ exist then
${(M_1\otimes
M_2)'}^{(\mu_1+\mu_2)}$ exists and
$${(M_1\otimes M_2)'}^{(\mu_1+\mu_2)}={M_1'}^{\mu_1} \otimes {M_2'}^{\mu_2}$$
{\bf Proof.} This is a functorial equality; also we can check it by
means of elementary
matrix calculations. $\square$
\medskip
{\bf Proposition 10.9.} Let $P\in \hbox{\bf A}$ be an irreducible
element. The Tate
module
$T_P({M'}^{\mu})$ is equal to $$T_P(\goth C)^{\otimes \mu}\otimes
\widehat{T_P(M)}$$
(equality of Galois modules) where $\widehat{T_P(M)}$ is the dual
Galois module.
\medskip
{\bf Proof.} It is completely analogous to the proof of the corresponding theorem for
tensor products
([G], Proposition 5.7.3, p. 157).
All modules in the below proof will be the Galois modules, and equalities of modules will
be equalities of Galois
modules. Recall that $E=E(M)$. Since $T_P(M)=\invlim_n E_{P^n}$, it is sufficient to prove that for any
$a\in \w$ we have $E({M'}^{\mu})_a=E(\goth C^{\otimes \mu})_a\otimes \hat E_a$, where $\hat E_a$ is the dual of $E_a$ in the meaning of [T], Definition 4.1.
We have the following sequence of equalities of modules:
$${M'}^{\mu}/a{M'}^{\mu}=\Hom_{\p[T]}(M/aM, \goth C^{\otimes \mu}/a\goth C^{\otimes
\mu})\eqno{(10.9.2)}$$
such that the action of $\tau$ on both sides of this equality coincide (to define the
action of $\tau$ on the
right and side of (10.9.2) we need the action of $\tau^{-1}$ on $M/aM$; it is
well-defined, because the determinant
of the action of $\tau$ on $M$ is a power of $T-\theta$, hence its image in
$\p[T]/a\p[T]$ is invertible). 10.9.2
follows immediately from the definition of ${M'}^{\mu}$;
$$({M'}^{\mu}/a{M'}^{\mu})^\tau=\Hom_{\n F_q[T]}((M/aM)^\tau,
(\goth C^{\otimes \mu}/a\goth C^{\otimes \mu})^\tau)\eqno{(10.9.3)}$$
This follows from 10.9.2 and the Lang's theorem
$$\goth M/a\goth M=(\goth M/a\goth M)^\tau\underset{\n F_q[T]/a\n
F_q[T]}\to{\otimes}\p[T]/a\p[T]$$
applied to both $\goth M=M$, $\goth M={M'}^{\mu}$ (we use that both $M$, ${M'}^{\mu}$ are
free $\p[T]$-modules).
Finally, we have a formula
$$E(\goth M)_a=\Hom_{\n F_q}((\goth M/a\goth M)^\tau, \n F_q)$$
([G], p. 152, last line of the proof of Proposition 5.6.3). Applying this formula to 10.9.3 we get the desired. $\square$
\medskip
{\bf 11. An explicit formula.}

\nopagebreak
\medskip
We return to the case $N=0$. Let $e_*$, $\goth A$, $\goth A_i$, $l$, $n$ be from (1.9). We consider in the
present section two simple types of t-motives (called standard-1 and standard-2
t-motives respectively) whose $\goth A_i$ have a row echelon form, and we give an
explicit formula for the dual of some standard-1
t-motives. Analogous formula can be easily obtained for more general types of
t-motives. These results are the
first step of the problem of description of all t-motives
having duals.
\medskip
{\bf 11.1.} For the reader's convenience, we give here the definition of standard-1
t-motives
for the case $n=2$ (here $\lambda_1$ and $\lambda_2 $ satisfying
$\lambda_1=l$,
$l>\lambda_2 \ge 2$ are parameters):

$$\goth A_0=\theta I_2, \hbox{ for } 0 < i < \lambda_2 \ \ \goth A_i \hbox{ is
arbitrary, } $$
$$\goth A_{\lambda_2}= \left(\matrix * & 0 \\ * & 1 \endmatrix \right),
\hbox{ for } \lambda_2
< i < l \ \ \goth A_i= \left(\matrix * & 0 \\ * & 0 \endmatrix \right), \ \
\goth A_l= \left(\matrix
1 & 0
\\0 & 0 \endmatrix \right) $$

{\bf 11.2.} To define standard-2 t-motives of dimension $n$, we need to fix
\medskip
1. A permutation $\vf\in S_n$, i.e. a 1 -- 1 map $\vf: (1, ..., n) \to (1, ..., n)$;
\medskip
2. A function $k: (1, ..., n) \to \n Z^+$ where $\n Z^+$ is the set of integers $\ge 1$.
\medskip
{\bf Definition.} A standard-2 t-motive of the type $(\vf, k)$ is an abelian
t-motive of dimension $n$ given by the formulas ($i=1,...,n$):
$$Te_{\vf(i)}= \theta e_{\vf(i)}
+\sum_{\alpha=1}^n\sum_{j=1}^{k(\alpha)-1}a_{j,\vf(i),\alpha}\
\tau^je_\alpha + \tau^{k(i)}e_i\eqno{(11.2.1)}$$
where $a_{j,\vf(i),\alpha}\in\p$ is the $(\vf(i),\alpha)$-th entry of the matrix $\goth A_j$.
\medskip
{\bf Proposition 11.2.2.} Formula 11.2.1 really defines a t-motive denoted by
$M=M(\vf,k)=M(\vf,k,a_{***})$. Its rank is $\sum_{\alpha=1}^n k(\alpha)$ and elements
$X_{\alpha j}:=\tau^je_\alpha$, $\alpha=1, ..., n$, $j=0, ..., k(\alpha)-1$, form its
$\p[T]$-basis.
$\square$
\medskip
The group $S_n$ acts on the set of types $(\vf, k)$ and on the set of the above $M$;
clearly for any $\psi\in S_n$ we have $\psi(M)$ is isomorphic to $M$. Particularly, we
can consider only $\vf$ of the following form of the product of $i$ cycles ($\alpha_0=0,
\ \alpha_i=n$):
$$\vf=(\alpha_0+1, ..., \alpha_1)(\alpha_1+1,...,\alpha_2)
...(\alpha_{i-1}+1,...,\alpha_i)\eqno{(11.2.3)}$$
(standard notation of the theory of permutations, for $\gamma\ne \alpha_j$ we have
$\vf(\gamma)=\gamma+1$, for $\gamma= \alpha_j$ we have  $\vf(\alpha_j)=\alpha_{j-1}+1$).
\medskip
{\bf Example 11.2.4.} Let $\vf$ be defined by 11.2.3, the quantity of cycles $i$ is equal
to $1$ and all $a_{***}=0$. Then the corresponding $M$ is of complete multiplication by a
CM-field $\n F_{q^r}(T)$ and its CM-type $\Phi$ is $\{\Id, \fr^{k(1)}, \fr^{k(1)+k(2)},
..., \fr^{k(1)+k(2)+...+k(n-1)}\}$ where $\fr$ is the Frobenuis homomorphism $\n
F_{q^r}\to \bar \n F_q$ (see 13.3, first case: formulas 13.3.1, 13.3.2 coinside with 11.2.1
for the given $\vf$ and $a_{***}=0$; $i_j$ of 13.3.0 is $k(1)+k(2)+...+k(j-1)$ of the
present notations).
\medskip
{\bf Definition 11.3.} A standard-1 t-motive is a standard-2 t-motive
whose $\vf$ is the identical permutation $Id$.
\medskip
{\bf 11.3.0.} Let $M=M(Id,k)$ be a standard-1 t-motive. Acting
by $\psi\in S_n$ we can consider only the case of non-increasing $k(j)$. We introduce a
number $\goth m\ge1$ --- the quantity of jumps of $k(j)$, and two sequences
$$0=\gamma_0< \gamma_1<...< \gamma_\goth m=n$$
(sequence of arguments of points of jumps of the function $k$) and
$$0=\lambda_{\goth m+1}< \lambda_{\goth m}<...< \lambda_2 <\lambda_1=l$$
(sequence of values of $k$ on segments $[\gamma_{i-1}+1, ...,\gamma_i]$) by the formulas
$$\matrix k(1)=...=k(\gamma_1)=\lambda_1 \\ \\
k(\gamma_1+1)=...=k(\gamma_2)=\lambda_2\\ ... \\
k(\gamma_{\goth m-1}+1)=...=k(\gamma_\goth m)=\lambda_\goth m\endmatrix \eqno{(11.3.1)}$$
\medskip
{\bf Example 11.3.2.} The t-motive $M$ of 11.1 is a standard-1 having $\goth m=2$,
$\gamma_1=1$, $\gamma_2=2$ and $\lambda_1$, $\lambda_2$ as in 11.1. Its rank $r=\lambda_1+\lambda_2$.
\medskip
{\bf Conjecture 11.3.3.} A standard-2 t-motive of the type $(\vf, k)$ (notations of
11.2.3) is pure iff $\forall j=1, ..., i$ we have:
$$\frac{\alpha_j-\alpha_{j-1}}{\sum_{\gamma=\alpha_{j-1}+1}^{\alpha_{j}}k(\gamma)} = \frac{n}{r}$$
This conjecture is obviously true if all $a_{***}$ are 0.
\medskip
To simplify exposition, we prove here only the following particular case of this
conjecture.
\medskip
{\bf Proposition 11.3.4.} Let $M$ be a standard-1 t-motive having $\goth m>1$,
defined over $\n F_q(\theta)$, having a good reduction at a point of degree 1 of $\n
F_q(\theta)$ (i.e. a point $\theta+c$, $c\in \n F_q$). Then $M$ is not pure.
\medskip
{\bf Proof.} Let $M$ be defined by 11.2.1, we use notations of 11.3.1. We consider the
action of Frobenius on $\tilde M$ --- the reduction of
$M$ at $\theta+c$. According [G], Theorem 5.6.10, it is sufficient to prove that orders
of the roots of the characteristic polynomial of Frobenius over $\w$ are not equal. More
exactly, we consider the valuation infinity on $\w$ (defined by the condition
$\ord(T)=-1$); the order corresponds to a continuation of this valuation to $\End(\tilde
M)$. The
action of Frobenius on $\tilde M$ coincides with multiplication by $\tau$, because the
degree of the reduction point is 1.

A basis $f_*$ of $M_{\p[T]}$ is the set of $X_{\alpha j}:=
\tau^je_{\alpha}$ of 11.2.2. The matrix
$Q(M)$ is defined by the following formulas for the action of
$\tau$ on $X_{\alpha j}$:
$$\tau(X_{\alpha j})=X_{\alpha,j+1} \hbox{ if } j <
k(\alpha)-1\eqno{(11.3.4.1)}$$
$$\tau(X_{\alpha,k(\alpha)-1})=TX_{\alpha,0} -
\sum_{\delta=1}^\goth m
\sum_{d=\lambda_{\delta+1}}^{\lambda_{\delta}-1}
\sum_{c=1}^{\gamma_{\delta}} a_{d\alpha c}X_{cd} \eqno{(11.3.4.2)}$$
This means that if we arrange $X_{\alpha j}$ in lexicographic order ($X_{\alpha_1 j_1}$ precedes to $X_{\alpha_2 j_2}$ if $\alpha_1 <
\alpha_2$) then the matrix
$Q(M)$ has the block form: $$Q(M)=(C_{ij})\ \ \ (i,j=1,...,n)$$ where $C_{ij}$ is a
$k(i)\times k(j)$-matrix of the form
$$C_{ii}=\left(\matrix 0&1&0&...&0\\0&0&1&...&0\\ ...&...&...&...&... \\ 0&0&0&...&1 \\
T-\theta & *&*&...&*\endmatrix \right), \ \ C_{ij}=\left(\matrix 0&0&...&0\\
...&...&...&...\\ 0&0&...&0 \\ 0& *&...&*\endmatrix \right) (i\ne j)$$
where asterisks mean elements $a_{***}$ (in some order). We consider the characteristic polynomial $P(X)\in (\p[T])[X]$ of $Q(M)$. We
have $$C_{ii}-XI_{k(i)}=\left(\matrix -X&1&0&...&0\\0&-X&1&...&0\\ ...&...&...&...&... \\
0&0&0&...&1 \\ t-\theta & *&*&...&*-X\endmatrix \right)$$

A subset of the set of entries of a matrix is called (following N.N.Luzin) a lightning if each row and each column of the matrix contains exactly one element of this subset. The product of elements of a lightning is called the value of this lightning (i.e. the determinant is the alternating sum of the values of all lightnings).
\medskip
{\bf Lemma 11.3.4.3.} If a non-zero lightning of $C_{ii}-XI_{k(i)}$ contains the term
$T-\theta$, then it does not contain any term containing $X$. $\square$
\medskip
Let $J$ be a subset of the set $1,...,n$ and $J'$ its complement.
\medskip
{\bf Corollary 11.3.4.4.} If a non-zero lightning of $Q(M)-XI_{r}$ contains terms
$T-\theta$ of blocks $C_{\al}$, $j\in J$, then its value is a polynomial in $X$ of degree
$\le \sum_{j'\in J'}k(j')$, and there exists exactly one such lightning (called the
principal $J$-lightning) whose value is a polynomial in $X$ of degree $\sum_{j'\in
J'}k(j')$. $\square$

Since the characteristic polynomial of Frobenius of $\tilde M$ is $\tilde P$
(respectively the valuation infinity of $\p[T]$), it is sufficient to prove that the
Newton polygon of $P(X)$ is not reduced to the segment $((0,-n); (r,0))$ defined by its
extreme terms $(T-\theta)^n$ and $X^r$. To do it, it is sufficient to find a point on its
Newton polygon which is below this segment. We consider $J_{min}=$ the set of all
$\gamma_\goth m - \gamma_{\goth m-1}$ diagonal blocks $C_{ii}$ ($i=\gamma_{\goth m-1}+1,
..., \gamma_\goth m$) of $Q(M)$ of minimal size $\lambda_\goth m$. The value of the
principal $J_{min}$-lightning is $(T-\theta)^{\gamma_\goth m - \gamma_{\goth m-1}}$ times
polynomial in $X$ of degree $d:=r-(\gamma_\goth m - \gamma_{\goth m-1})\lambda_\goth m$.
Corollary 11.3.4.4 implies that if the value of any other lightning of $Q(M)-XI_r$ contains a
term whose $X$-degree is equal to $d$, then the $T$-degree of this term is strictly less
than $\gamma_\goth m - \gamma_{\goth m-
 1}$. This means that if we write $P(X)=\sum_{i=0}^r C_iX^i$, $C_i\in \p[T]$, then
$\ord_\infty(C_d)= -(\gamma_\goth m - \gamma_{\goth m-1})$, i.e. the point with
coordinates $[-(\gamma_\goth m - \gamma_{\goth m-1}), d]$ belongs to the Newton diagram
of $P(X)$, i.e. it is above (really, at) the Newton polygon of $P(X)$. This point is
below the segment $((0,-n); (r,0))$. $\square$
\medskip
{\bf Remark 11.3.4.5.} It is easy to see that the Newton polygon of $P(X)$ coincides with
the Newton polygon of the direct sum
of trivial Drinfeld modules of ranks $\lambda_*$, i.e. with the Newton polygon of the
polynomial
$$\prod_{i=1}^{\goth m} (X^{\lambda_i}-T)^{\gamma_i - \gamma_{i-1}}$$
\medskip
{\bf 11.4.} To formulate the below theorem 11.5 we need some notations. Let $M$ be a
standard-1 t-motive defined by formulas 11.2.1, 11.3.1. We impose the condition
$\lambda_\goth m \ge 3$. Theorem 11.5 affirms that it has dual. To find explicitly the
dual of $M$, we need to choose an arbitrary function
$\nu: (i,j) \to \nu(i,j)$ which is a 1 - 1 map from the set of pairs
$(i,j)$ such that

$$1 \le i \le n; \ \ 1 \le j \le k(i)-2 \eqno{(11.4.1)}$$
to the set $[n+1, ..., r-n]$ (recall that $r=\sum_{i=1}^n k(i)= \sum_{i=1}^\goth m
(\gamma_{i}-\gamma_{i-1})\lambda_i$).
\medskip
Let the $(r-n)\times(r-n)$-matrices $B_1$, $B_2$ be defined by the
following formulas
(here and until the end of the proof of 11.5 we have $i ,\alpha = 1, ... , n$; \ \
$b_{\beta\gamma\delta}$
is the $(\gamma\delta)$-th entry of $B_\beta$, all
entries of $B_1$,
$B_2$ that are not in the below list are 0):
\medskip
{\bf 11.4.2.} $b_{1i\alpha}= - a_{k(i)-1,\alpha,i}$;
\medskip
$b_{1,\nu(i,j),\alpha} = - a_{j,\alpha,i}$ for $1 \le j \le k(i)-2$;
\medskip
$b_{1,\nu(i,j+1),\nu(i,j)}=1$ for $1 \le j \le k(i)-3$;
\medskip
$b_{1,i,\nu(i,k(i)-2)} =1$;
\medskip
$b_{2,\nu(i,1),i}=1$.
\medskip
We let $B=\theta I_{r-n}+B_1\tau+B_2\tau^2$ and consider a t-motive $M(B)$ (see 11.5.1
below). Formulas
11.4.2 mean that $M(B)$ is standard-2, its $\vf=\vf_B$ is a product of $n$ cycles
$$i\overset{\vf_B}\to{\to}\nu(i,1)\overset{\vf_B}\to{\to}\nu(i,2)\overset{\vf_B}\to{\to}...
\overset{\vf_B}\to{\to}\nu(i,k(i)-2)\overset{\vf_B}\to{\to}i$$ and its $k=k_B$ is defined
by the
formulas $k_B(\gamma)=2$ for $\gamma \in [1, ...,n]$, $k_B(\gamma)=1$ for $\gamma \in
[n+1, ...,r-n]$.
\medskip
{\bf Theorem 11.5.} Let $M$ be from 11.4 (i.e. a standard-1 t-motive having
$\lambda_\goth m \ge 3$). Then $M'=M(B)$.
\medskip
{\bf Proof.}\footnotemark \footnotetext{This proof is a generalization
of the corresponding proof of
Taguchi; we keep his notations.} Let $e'_*=(e'_1, ... e'_{r-n})^t$ be the vector column
of elements of a basis
of $M(B)$ over $\p[\tau]$ satisfying
$$T e'_*=Be'_*\eqno{(11.5.1)}$$
Let us consider the set of pairs $(j,\goth k)$ such that either $j=1,...,n$,
$\goth k=0,1$ or
$j=n+1,...,r-n$, $\goth k=0$. For each pair $(j,\goth k)$ of this set we let (as
in [T], p. 580)
$Y_{j \goth k} = \tau^{\goth k}e'_j$. Formulas (11.4.2) show that these $Y_{**} $
form a basis of
$M(B)_{\p[T]}$, and the action of $\tau$ on this basis is given by the
following formulas
(here $j=1, ... , k(i)-2$):
$$\tau(Y_{i,0})=Y_{i,1} \eqno{(11.5.2.1)}$$
$$\tau(Y_{i,1})= (T-\theta) Y_{\nu(i,1),0} + \sum_{\gamma=1}^n
a_{1\gamma i}Y_{\gamma,1} \eqno{(11.5.2.2)}$$
$$\tau(Y_{\nu(i,j),0})= (T-\theta) Y_{\nu(i,j+1),0} + \sum_{\gamma=1}^n
a_{j+1,\gamma,i}Y_{\gamma,1} \hbox{ if }j<k(i)-2 \eqno{(11.5.2.3)} $$
$$\tau(Y_{\nu(i,k(i)-2),0})= (T-\theta) Y_{i,0} + \sum_{\gamma=1}^n
a_{k(i)-1,\gamma,i}Y_{\gamma,1} \eqno{(11.5.2.4)}$$
Let $X'_{**}$ be the dual basis to the basis $X_{**}$ of 11.2.2.
\medskip
{\bf 11.5.3.} Let us consider the following correspondence between
$X'_{**}$ and $Y_{**}$:
\medskip
$X'_{ij}$ corresponds to $Y_{\nu(i,j),0}$ for
the pair $(i,j)$ like in (11.4.1),
\medskip
$X'_{i0}$ corresponds to $Y_{i1}$ for $1 \le i \le n$;
\medskip
$X'_{i, k(i)-1}$ corresponds to $Y_{i0}$ for $1 \le i \le n$. 
\medskip
Therefore, in order to prove the Theorem 11.5 we must check that
matrices defined by the dual to (11.3.4.*) and by (11.5.2.*) satisfy (1.10.1) under
identification (11.5.3). This is an elementary exercise. $\square $
\medskip
{\bf Remark 11.6.} Clearly it is possible to generalize the Theorem 11.5 to a larger class
of t-motives --- some subclass of standard-3 t-motives, see Definition 11.8.1.
The below example of the proof of Proposition 11.8.7 shows that probably the condition
$\lambda_\goth m \ge 3$ of
the Theorem 11.5 can be changed by
$\lambda_\goth m \ge 2$: it is necessary
to modify slightly formulas 11.4.2. From another side, a standard-1
t-motive of the Example 2.5 shows that this condition cannot
be changed to $\lambda_\goth m\ge1$.
\medskip
{\bf 11.7. An elementary transformation.} To formulate the proposition
11.7.3, we change slightly notations in 1.9.1, namely, instead of $\goth A =
\sum_{i=0}^l \goth A_i \tau^i$ we consider polynomials $P_k(M)$ of $x_1, ...
,x_n$ ($k=1, ...,n$) defined by the formula
$$P_k(M)= \sum_{i=0}^l\sum_{j=1}^n a_{ikj} x^{q^i}_j \eqno{(11.7.1)} $$
Particularly, if $E$ is the t-module associated to $M$ (see [G], 5.4.5),
$x_*=(x_1, ..., x_n)^t$ an element of $E$ then 11.7.1 is equivalent to
$Tx_*=P_*(x_*)$ where $P_*=(P_1(M),...,P_n(M))^t$ is the vector column.
For a standard-1 t-motive $M$ (we use notations of 11.3.0) having $\goth m\ge 2$ we denote vector columns $\goth
P_1(M)=(P_1(M),...,P_{\gamma_1}(M))^t$, $\goth
P_2(M)=(P_{\gamma_1+1}(M),...,P_{\gamma_2}(M))^t$. We use similar
notations for $M'$.
\medskip
{\bf 11.7.2.} Let $M$ be as above, we consider the case
$\lambda_2=\lambda_1-1$. Let $C$ be a fixed $\gamma_1\times (\gamma_
2-\gamma_1)$-matrix. We define a transformed t-motive $M_1$ by the
formulas

$$\goth P_1(M_1)= \goth P_1(M)+C\goth P_2(M)^q$$

$$P_i(M_1)=P_i(M) \hbox{ for } i>\gamma_1$$
\medskip
{\bf Proposition 11.7.3.} For $M$, $C$, $M_1$ of 11.7.2 the dual $M'_1$
of $M_1$ is described by the following formulas:
$$\goth P_2(M'_1)= \goth P_2(M')-C^t\goth P_1(M')^q$$
$$P_i(M'_1)=P_i(M')\hbox{ for }i\not\in[\gamma_1+1, ... , \gamma_2]$$

{\bf Proof} is similar to the proof of the Theorem 11.5, it is omitted.
$\square$
\medskip
\medskip
{\bf 11.8. Non-pure t-motives.} Most results of this subsection are conditional. We shall show that under some natural conjecture the condition of purity in 8.2 and 8.4 is essential, and that for non-pure t-motives the notion of algebraic duality is richer than the notion of analytic duality.

We generalize slightly the definition 11.2.1 as follows. Let
$\succ$ be a linear ordering on the set $[1,...,n]$, and let $\vf$, $k$ be as in 11.2.
\medskip
{\bf Definition 11.8.1.} A standard-3 t-motive of the type $(\vf, k,\succ)$ is
a t-motive of
dimension
$n$ given by the formulas
$$Te_{\vf(i)}= \theta e_{\vf(i)} +\sum_{j=1}^n\sum_{l=1}^{k(j)-1}a_{l,\vf(i),j}\
\tau^le_j +\sum_{j\succ i} a_{k(j),\vf(i),j}\ \tau^{k(j)}e_j +
\tau^{k(i)}e_i\eqno{(11.8.2)}$$
where $a_{***}\in\p$ are coefficients (the only difference with 11.2.1 is the term
$\sum_{j\succ i} a_{k(j),\vf(i),j}\ \tau^{k(j)}e_j$). We denote it by $M(a_{***})$.

Let $M_1=M(a_{1***})$, $M_2=M(a_{2***})$ be two isomorphic standard-3 t-motives of the same type $(\vf,
k,\succ)$ with
$\p[\tau]$-bases $e_{1*}$, $e_{2*}$ respectively (we use notations of 11.8.2 for both
$M_1$, $M_2$). There exists $C\in M_n(\p[\tau])$ such that the formula defining an isomorphism
between $M_1$ and $M_2$ is the following: $e_{2*}=Ce_{1*}$.
\medskip
{\bf Conjecture 11.8.3.} For a generic set of $a_{1***}$ there exists only a countable set of $a_{2***}$ such that $M_2$ is isomorphic to $M_1$.
\medskip
This conjecture is based on calculations in some explicit cases. Particularly, it is proved if $M_1$, $M_2$ are given by the below formula 11.8.5.1 and entries of $C$ are polynomials in $\tau$ of degree $\le 1$.

We denote by $\Cal M_{u}(r,n)$ the moduli space of uniformizable t-motives of the rank
$r$ and dimension $n$,
by $\Cal L(r,n)$ the moduli space of lattices of the rank $r$ and dimension $n$ and by
$\goth L:\Cal M_{u}(r,n) \to \Cal L(r,n)$ the functor of lattice associated to an uniformizable
t-motive.
\medskip
{\bf Proposition 11.8.5.} Conjecture 11.8.3 implies that the dimension of the fibers of
$\goth L$ is $> 0$ for $r=3$, $n=2$. Particularly, we cannot omit condition of purity in the
statement of 8.2.
\medskip
{\bf Proof.} We consider standard-3 t-motives of the type $n=2$,
$\vf=Id$, $k(1)=2$,
$k(2)=1$, $2\succ 1$. Such $M_1=M_1(a_{111}, a_{112}, a_{121})$ is given by
$$\goth A_0=\theta I_2, \ \ \goth A_{1}= \left(\matrix a_{111} & a_{112}
\\ a_{121} & 1 \endmatrix \right), \ \ \goth A_2= \left(\matrix 1 & 0
\\0 & 0 \endmatrix \right)\eqno{(11.8.5.1)}$$ (notations of Example 10.5).
It has $r=3$, it is not pure, hence it has no dual.
Conjecture 11.8.3 implies
that the dimension of the moduli space of these t-motives is 3 (because there are 3
coefficients
$a_{111}, a_{112}, a_{121}$). Uniformizable t-motives form an open subset of this moduli
space, while
the moduli space of lattices of $n=2$ and $r=3$ has dimension 2. $\square$
\medskip
{\bf Remark.} Similar calculations are valid for any sufficiently large $r$, $n$.
\medskip
Standard-3 t-motives of the above type have not dual. The following proposition
shows that the same
phenomenon holds for t-motives having dual. We denote by $\Cal M_{u,d}(r,n)$ the
moduli space of uniformizable t-motives of the rank $r$ and dimension $n$ having dual, by
$\Cal L_d(r,n)$ the moduli space of lattices of the rank $r$ and dimension $n$ having
dual, by $\goth L_d:\Cal M_{u,d}(r,n) \to \Cal L_d(r,n)$ the functor of lattice and by $D_M: \Cal
M_{u,d}(r,n)\to \Cal M_{u,d}(r,r-n)$, $D_L: \Cal L_d(r,n)\to \Cal L_d(r,r-n)$ the
functors of duality on t-motives and lattices respectively. Practically, Theorem 5
means that the following diagram is commutative:
$$\matrix \Cal M_{u,d}(r,n)&\overset{D_M}\to{\to}&\Cal M_{u,d}(r,r-n)
\\ \\ \goth L_d\downarrow&&\goth L_d\downarrow\\  \\ \Cal L_d(r,n)&\overset{D_L}\to{\to}&\Cal
L_d(r,r-n)\endmatrix\eqno{(11.8.6)}$$

{\bf Proposition 11.8.7.} Conjecture 11.8.3 implies that the dimension of the fibers of
$\goth L_d$ in the diagram (11.8.6) is $> 0$ for $r=5$, $n=2$.
\medskip
Practically, this means that the notion of algebraic duality is "richer" than the notion
of analytic duality.
\medskip
{\bf Proof.} We consider standard-3 t-motives of the type $n=2$,
$\vf=Id$, $k(1)=3$,
$k(2)=2$, $2\succ 1$, $r=5$. Such $M$ is given by
$$\goth A_0=\theta I_2, \ \ \goth A_{1}= \left(\matrix a_{111} & a_{112}
\\ a_{121} & a_{122}  \endmatrix \right), \ \ \goth A_{2}= \left(\matrix a_{211} & a_{212}
\\ a_{221} & 1 \endmatrix \right), \ \ \goth A_3= \left(\matrix 1 & 0
\\0 & 0 \endmatrix \right)$$ (notations of Example 10.5). It has dual. Really, we denote by
$A_{i*j}$ the $j$-th column of $\goth A_i$, and we denote by $ (C_1 | C_2 )$ the matrix formed
by union of columns
$C_1$, $C_2$. Then $M'=M(B)$ is also a standard-3 t-motive, where
\medskip
$B_1= \left(\matrix - \det \goth A_2 & -a_{221} & 1
\\ - \det ( A_{1*2} | A_{2*2} ) & -a_{122} & 0
\\ - \det ( A_{1*1} | A_{2*2} ) & -a_{121} & 0
\endmatrix \right)$,
$B_2 = \left(\matrix 0 & 0 & 0
\\ -a_{212}^q & 1 & 0
\\ 1 & 0 & 0 \endmatrix \right)$

The same arguments as in the proof of Proposition 11.8.5 show that the conjecture 11.8.3
implies that the dimension of the moduli space of these t-motives is 7, while
the moduli space of lattices of $n=2$ and $r=5$ has dimension 6. $\square$
\medskip
As above, similar calculations are valid for any sufficiently large $r$, $n$; clearly the dimension of fibers of $\goth L_d$ becomes larger as $r$, $n$ grow.
\medskip
Let us mention two open questions related to the functor $\goth L$. Firstly, let $L$ be a self-dual lattice such that $L\in \goth L(\Cal M_{u,d}(2n,n))$. This means that $D_M: \goth L_d^{-1}(L) \to \goth L_d^{-1}(L)$ is defined.
\medskip
{\bf Open question 11.8.8.} What can we tell on this functor, for example, what is the dimension of its stable elements?
\medskip
Secondly, let us consider $M_1$, $M_2$ of CM-type with CM-field $\n F_{q^r}(T)$, see 13.3.
\medskip
{\bf Open question 11.8.9.} Let the CM-types $\Phi_1$, $\Phi_2$ of the above $M_1$, $M_2$ satisfy $\Phi_1\ne \alpha \Phi_2$, where $\alpha\in \Gal(\n F_{q^r}(T)/\n F_{q}(T))$. Are lattices $L(M_1)$, $L(M_2)$ non-isomorphic?
\medskip
Clearly the negative answer to this question implies the negative answer to the Question 8.2.
\medskip
For any given $M_1$, $M_2$ the answer can be easily found by computer calculation.
Really, let $M$ be one of $M_1$, $M_2$, $c_1,...,c_r$ a basis of $\n F_{q^r}/\n F_{q}$ and
$\alpha_{1},...,\alpha_{n} \subset \Gal(\n F_{q^r}(\theta)/\n F_{q}(\theta))$ the CM-type of $M$.
We define matrices $\Cal M$, $\Cal N$ as follows: $(\Cal M)_{ij}=\alpha_j(c_i)$ $(i,j=1,...,n$),
$(\Cal N)_{ij}=\alpha_{j}(c_{n+i})$, $j=1,...,n$, $i=1,...,r-n$. The Siegel matrix $Z(M)$ is obviously $\Cal N\Cal
M^{-1}$. So, we can find explicitly $Z(M_1)$, $Z(M_2)$ for both $M_1$, $M_2$. To check
whether $Z(M_1)$, $Z(M_2)$ are equivalent or not, it is sufficient to find a solution to
3.8.1 such that the entries of $A$, $B$, $C$, $D$ are in $M_{*,*}(\n F_q)$ (this is
obvious: the condition $\exists \gamma \in GL_r(\n Z_\infty)$ is equivalent to the
condition $\exists \gamma \in GL_r(\n F_q)$, because entries of $Z(M_1)$, $Z(M_2)$ are in $\n F_{q^r}$). The equation 3.8.1 is linear with respect
to $A$, $B$, $C$, $D$, and we can check whether its solution satisfying $\det
\gamma\ne 0$ exists or not.
\medskip
For the case $q=2$, $r=4$, $n=2$, CM-types of $M_1$, $M_2$ are $(Id, Fr)$, $(Id, Fr^2)$
respectively, a calculation shows that the answer is positive: lattices $L(M_1)$, $L(M_2)$ are not isomorphic.
\medskip
{\bf 12. t-motives having multiplications.}
\medskip
Let $\goth K$ be a separable extension of $\n F_q(T)$ such that $\goth K_C:=\goth K\underset{\n F_q}\to{\otimes} \p$ is also a field, $\pi: X\to P^1(\p)$ the projection of curves over $\p$ corresponding to $\p(T)\subset \goth K_C$. Let $\goth K$, $X$ satisfy the condition: $\infty\in X$ is the only point
on $X$ over $\infty\in P^1(\p)$. Let $\w_{\goth K}$ be the subring of
$\goth K$ consisting of elements regular outside of infinity. We
denote $g=\dim \goth K/\n F_q(T)$ and
$ \alpha_1, ... , \alpha_g: \goth K\to\p$ --- inclusions over $\iota: \n F_q(T)\to\p$
(recall that $\iota(T)=\theta$). Let $\Cal W$ be a central simple algebra
over $\goth K$ of dimension $\goth q^2$. Each $ \alpha_i: \goth K \to
\p$ can be extended to a representation $\chi_i: \Cal W \to M_\goth
q(\p)$.
\medskip
{\bf 12.1. Analytic CM-type.} Let $(L, V)$ be as in Section 2 (recall that
$\w=\n F_q[T]$) such that there exists an inclusion $i: \Cal W \to
\End^0(L, V)$, where $\End^0(L, V)=\End(L, V)\underset{\w}\to{\otimes} \n F_q(T)$.
It defines a representation of $\Cal W$
on $V$ denoted by $\Psi$ which is isomorphic to $\sum_{i=1}^g\goth r_i\chi_i$ where
$\{\goth r_i\}$ are some multiplicities (the CM-type of the action of
$\Cal W$ on $(L, V)$). [Proof: restriction of $\Psi$ on $\goth K$ is a sum of one-dimensional representations, i.e. $V=\oplus_{i=1}^g V_i$ where $k\in \goth K$ acts on $V_i$ by multiplication by $\alpha_i(k)$. Spaces $V_i$ are $\Psi$-invariant. We consider an isomorphism $\Cal W\otimes_\goth K\p=M_\goth q(\p)$ where the inclusion of $\goth K$ in $\p$ is $\alpha_i$. We extend $\Psi|_{V_i}$ to $\Cal W\otimes_\goth K\p$ by $\p$-linearity using the inclusion $\alpha_i$ of $\goth K$ in $\p$. It remains to show that a representation of $M_\goth q(\p)$ is a direct sum of its $\goth q$-dimensional standard representations. We consider the corresponding representation of Lie algebra $\goth s\goth l_\goth q(\p)$. It is a sum of irreducible representations. Let $\omega$ be the highest weight of any of these irreducible representations. $\omega$ is extended uniquely to the set of diagonal matrices of $M_\goth q(\p)$, because $\omega$ is identical on scalars. Since our representation is not only of Lie algebra but of algebra $M_\goth q(\p)$, we get that $\omega$ is a ring homomorphism $\Diag(M_\goth q(\p))\to \p$. There exists the only such $\omega$ corresponding to the $\goth q$-dimensional standard representation].

Further, we
denote $m=\dim_\Cal W L \otimes \n F_q(T)$ ($g$, $\goth q$, $\Psi$, $\goth
r_i$, $m$ are analogs of $g$, $q$, $\Phi$, $r_i$, $m$ of [Sh63] respectively).
Clearly we have
$$n=\goth q\sum_{i=1}^g\goth r_i, \ \ r=mg\goth q^2\eqno{(12.2)}$$
By functoriality, we have the dual inclusion $i': \Cal W^{op} \to
\End^0(L',V')$ where $\Cal W^{op}$ is the opposite algebra.
\medskip
{\bf Remark.} A construction of Hilbert-Blumental modules ([A], 4.3, p. 498) practically is a particular case
of the present
construction: for Hilbert-Blumental modules we have $\goth q=1$, i.e. $\goth K=\Cal W$, and all $\goth r_i=1$. Anderson considers the case when $\infty$ splits completely; this difference with the present case is not essential.
\medskip
{\bf Proposition 12.3.} If the dual pair $(L',V')$ exists then the
CM-type of the dual inclusion is $\{m\goth
q - \goth r_i\}$, $i=1, ... ,g$.
\medskip
{\bf Proof.} We have $L\underset{\n Z_\infty}\to{\otimes}\p$ is
isomorphic to $(\Cal W\underset{\n F_q(\theta)}\to{\otimes}\p)^m$ as a
$\Cal W$-module. Since the natural representation of $\Cal W$ on $\Cal
W\underset{\n F_q(\theta)}\to{\otimes}\p$ is isomorphic to $\goth
q\sum_{i=1}^g\chi_i$ we get that $L\underset{\n Z_\infty}\to{\otimes}\p$ is isomorphic to $m\goth q\sum_{i=1}^g\chi_i$ as a $\Cal W$-module. Consideration of the exact
sequence $0\to {V'}^* \to L\underset{\n Z_\infty}\to{\otimes}\p \to V\to0$ gives us the desired. $\square$
\medskip
{\bf Remark 12.4.1.} This result is an analog of the corresponding
theorem in the number field case. We use notations of [Sh63], Section 2. Let
$A$ be an abelian variety having endomorphism algebra of type IV, and
$(r_\nu, s_\nu)=(r_\nu(A), s_\nu(A))$ are from [Sh63], Section 2, (8).
Then
$$r_\nu(A')=mq-r_\nu(A)=s_\nu(A), \ s_\nu(A')=mq-s_\nu(A)=r_\nu(A)$$
By the way, Shimura writes that the CM-types of $A$ and $A'$ coincide
([Sh98], 6.3, second line below (5), case $A$ of CM-type). We see that
his affirmation is not natural: he considers the complex conjugate
action of the endomorphism ring on $A'$. It is necessary to take into
consideration this difference of notations comparing formulas of 12.3 and 13.2
with the corresponding formulas of Shimura.
\medskip
{\bf Remark 12.4.2.} According [L09], a t-motive $M$ is an analog of an abelian variety $A$ with multiplication by an imaginary quadratic field $K$. The above consideration shows that this analogy holds for $M$ and $A$ having more multiplications. Really, if $A$ has more multiplications then (we use notations of [Sh63], Section 2) $F_0=FK$, and numbers $(r_\nu(A), s_\nu(A))$ satisfy $n(A)=q\sum_{i=1}^g r_\nu(A)$, where $(n(A), \dim(A)-n(A))$ is the signature of $A$ treated as an abelian variety with multiplication by $K$. This is an analog of 12.2.
\medskip
{\bf 12.5. Complete multiplication.} Here we consider the case $\goth
q=m=1$, i.e. $\goth K=\Cal W$ and $g=r$.
\medskip
{\bf Lemma 12.5.1.} In this case the condition $N=0$
implies that the CM-type $$\sum_{i=1}^r\goth r_i\alpha_i\eqno{(12.5.2)}$$ of the action of
$\goth K$ on
on $(L,V)$ has the property: all $\goth r_i$ are 0 or 1.
\medskip
{\bf Proof.} $N=0$ means that
the action of $T\in\w$ on $V$ is simply multiplication by $\theta$. We write the
CM-type $\sum_{i=1}^r\goth r_i\chi_i$ in the form $\sum_{i=1}^n\chi_{\alpha_i}$ where
$\alpha_1,...,\alpha_n\in[1,...,r]$ are not necessarily distinct. Let $l_1$ be an (only)
element of a basis of $L\otimes_{\w_{\goth K}}\goth K$ over $\goth K$ and $e_1,...,e_n$ a
basis of $V$
over $\p$ such that the action of $\goth K$ on $V$ is given by the formulas
$$k(e_i)=\chi_{\alpha_i}(k)e_i, \ \ \ k\in \goth K$$
Multiplying $e_i$ by scalars if necessary, we can assume that $l_1=\sum e_i$. Therefore,
if $\alpha_i=\alpha_j$ (i.e. not all $\goth r_*$ in (12.5.2) are 0, \ 1) then the
$e_{\alpha_i}$-th coordinate of any element of $L$ coincide with its $e_{\alpha_j}$-th
coordinate, hence $L$ does not $\p$-generate $V$ --- a contradiction. $\square$
\medskip
{\bf 12.5.3.} Let $M$ be a t-motive of
rank $r$ and dimension $n$ having multiplication by $\w_{\goth K}$. Recall that we
consider only
the case $N=0$. This means that the character of the action of $\goth K$ on $M/\tau M$
is isomorphic to $\sum_{i=1}^r\goth r_i\alpha_i$. Since $E(M)=(M/\tau M)^*$ we get that
the character of the action of $\goth K$ on $E(M)$ is the same. If
$$\hbox{all $\goth r_i$ are 0 or 1}\eqno{(12.5.4)}$$
we shall use the terminology that $M$ has the CM-type $\Phi\subset \{\alpha_1, ... ,
\alpha_r \}$
where $\Phi$ is defined by the condition $\alpha_i\in \Phi \iff \goth r_i=1$.

It is easy to see that this case occurs for uniformizable $M$. Really, if $M$ is
uniformizable then
the action of $\goth K$ can be prolonged on $(L(M),V(M))$, and the character
of the action of $\goth K$ on $V(M)$ coincides with the one on $E(M)$. The result follows
from Lemma 12.5.1.
\medskip
{\bf Lemma 12.5.5.} There exists a canonical isomorphism $\gamma$ from the set of
inclusions $\alpha_1, ... , \alpha_r$ to the set of points $ \theta_{\alpha_1}, ... ,
\theta_{\alpha_r}$ of $X$ over $\theta\in P^1(\p)$.
\medskip
{\bf Proof.} A point $t\in X$ over $\theta\in P^1(\p)$ defines a function
$\vf_t: \goth K_{C} \to P^1(\p)$ --- the value of
an element $f\in \goth K_{C}$ treated as a function on $X$ at the point $t$. This
function
must satisfy the standard axioms of valuation and the condition $\vf_t(T)=\theta$. Let
$\alpha_i$
be an inclusion of $\goth K$ to $\p$ over $\iota$. It defines a valuation
$\vf_{\alpha_i}: \goth K_{C} \to P^1(\p)$ by the formula $\vf_{\alpha_i}(k\otimes
f)=\alpha_i(k) f(\theta)$,
where $k\in\goth K$, $f\in \p(T)$. We define $\gamma(\alpha_i)$ by the condition
$\vf_{\gamma(\alpha_i)}=\vf_{\alpha_i}$; it is easy to see that $\gamma$ is an
isomorphism. $\square$
\medskip
{\bf Theorem 12.6.} For any above \{$\goth K$, $\Phi$\} there exists an
t-motive $(M,\tau)$ with
complete multiplication by $\goth K$ having CM-type $\Phi$.
\medskip
{\bf Proof (Drinfeld).} We denote the divisor $\sum_{\alpha_i\in\Phi}\gamma({\alpha_i})$
by
$\theta_\Phi$. We construct a $\Cal F$-sheaf $F$ of dimension 1 over
$\goth K$ which will give us $M$.
Let $\fr$ be the Frobenius map on $ \Pic_0(X)$. It is an algebraic
map, and the $\fr-\Id: \Pic_0(X) \to \Pic_0(X)$ is an algebraic map as
well. Since the action of $\fr$ on the tangent space of $ \Pic_0(X)$
at 0 is the zero map, the action of $\fr-\Id$ on the tangent space of
$ \Pic_0(X)$ at 0 is the minus identical map and hence $\fr-\Id$ is an
isogeny of $ \Pic_0(X)$. Particularly, there exists a divisor $D$ of
degree 0 on $X$ such that we have the following equality in $\Pic_0(X)$:
$$\fr(D)-D=-\theta_\Phi+n\infty \eqno{(12.6.0)}$$
This means that if we let $F=F_\Phi=O(D)$ then there exists a rational map
$\tau_X=\tau_{X,\Phi}: F^{(1)}\to F$ such
that
$$\Div(\tau_X)=\theta_\Phi-n\infty\eqno{(12.6.1)}$$
The pair $(F_\Phi, \tau_{X,\Phi})$ is the desired $\Cal F$-sheaf.
\medskip
{\bf Remark.} It is easy to see that if the genus of $X$ is $> 0$ then different CM-types $\Phi_1$, $\Phi_2$ give us different sheaves $F_{\Phi_1}$, $F_{\Phi_2}$, while if the genus of $X$ is 0 then $F_{\Phi_1}=F_{\Phi_2}=\Cal O$, but the maps $\tau_{X,\Phi_1}$, $\tau_{X,\Phi_2}$ are clearly different.
\medskip
Let $U_0=X-\{\infty\}$ be an open part of
$X$. We denote $F(U_0)$ by $\Cal M$, hence $F^{(1)}(U_0)=\Cal M^{(1)}$.
Since the support of the negative part of the right hand side of
12.6.1 is $\{\infty\}$,
we get that the (a priory rational) map $\tau_X(U_0): \Cal M^{(1)}\to \Cal M$ is
really a map of
$\w_{\goth K}$-modules.

Let $M$ be a $\p[T]$-module obtained from $\Cal M$ by restriction of scalars from
$\w_{\goth K}$ to $\p[T]$. Construction $F\mapsto M$ is functorial, and we denote this functor by $\delta$. Further, we denote by $\alpha$ the tautological isomorphism
$\Cal M\to M$. $M$ is a free $r$-dimensional $\p[T]$-module, and (because
$M^{(1)}$ is isomorphic to $M$) the
same restriction of scalars of $\tau_X(U_0)$ defines us a $\p[T]$-skew map from $M$ to
$M$ denoted by $\tau$ (skew means that $\tau(zm)=z^q\tau(m)$, $z\in\p$). $\tau$ is
defined by the formula $\tau(m)=\alpha\circ \tau_X((\alpha^{-1}(m))^{(1)})$.

It is easy to check that
$(M,\tau)$ is the required t-motive. Really, $M$ is a $\w_{\goth K}$-module, and
$\tau$
commutes with this multiplication. The fact that the positive part
of the right hand
side of 12.6.1 is $\theta_\Phi$ means that 1.13.2 holds for $M$ and that the CM-type of
the action of
$\w_{\goth K}$ is $\Phi$.
\medskip
{\bf Remark 12.6.2.} It is easy to prove for this case that $M$ is a free
$\p[\tau]$-module.
Really, it is
sufficient to prove (see [G], Lemma 12.4.10)
that $M$ is finitely generated as a $\p[\tau]$-module. We choose $D$ such that
$\infty\not\in \Supp (D)$. There exists $P\in \goth K_{C}^*$
such that $\tau_X(U_0): \Cal M^{(1)}\to \Cal M$ is multiplication by $P$ (recall that both
$\Cal M^{(1)}$, $\Cal M$ are $\w_{\goth K}$-submodules of $\goth K$). 12.6.0 implies that
$-\ord_\infty(P)=n$.
There
exists a number $n_1$ such that $$\hbox{(a) $h^0(X,\Cal O(D+n_1\infty))>0$; \ \ (b)
for any $k\ge 0$ we have}$$ $$h^0(X,\Cal O(D+(n_1+k)\infty))=h^0(X,\Cal
O(D+n_1\infty))+k\eqno{(12.6.3)}$$ $$h^0(X,\Cal O(D^{(1)}+(n_1+k)\infty))=h^0(X,\Cal
O(D^{(1)}+n_1\infty))+k\eqno{(12.6.4)}$$
It is sufficient to prove that if $g_1,...,g_k$ are elements of a basis of $H^0(X,\Cal
O(D+(n_1+n)\infty))$, then
for any $Q\in \Cal M$ the element $\alpha(Q)\in M$ is generated by
$\alpha(g_1),...,\alpha(g_k)$ over
$\p[\tau]$. We prove it by induction by $n_2:=-\ord_\infty(Q)$. If $n_2\le n_1+n$ the
result is trivial. If not
then 12.6.3, 12.6.4 imply that the multiplication by $P$ defines an isomorphism
$$H^0(X, \Cal O(D^{(1)}+(n_2-n)\infty))/H^0(X, \Cal O(D^{(1)}+(n_2-n-1)\infty))\to$$
$$\to H^0(X,\Cal O(D+n_2\infty))/H^0(X, \Cal O(D+(n_2-1)\infty))$$
This means that $\exists Q_1\in H^0(X,\Cal O(D^{(1)}+(n_2-n)\infty))$,
$-\ord_\infty(Q_1)=n_2-n$ such
that
$-\ord_\infty(Q-PQ_1)\le n_2-1$. An element $Q_1^{(-1)}\in \Cal M$ exists;
since $\alpha(Q)=\tau(\alpha(Q_1^{(-1)}))+\alpha(Q-PQ_1)$, the result follows by
induction.
$\square$
\medskip
If $\goth K$ and $\Phi$ are given then the construction of the Theorem 12.6 defines $F$ uniquely up to tensoring by $O(D)$ where $D\in \Div (X(\goth K))$. We denote the set of these $F$ by $F$(\{$\goth K$, $\Phi$\}), and we denote by $M$(\{$\goth K$, $\Phi$\}) the set $\delta(F$(\{$\goth K$, $\Phi$\})). Further, we denote by $\Phi'=\{\alpha_1, ... ,
\alpha_r\}-\Phi$ the complementary CM-type.
\medskip
{\bf Theorem 12.7.} Let $M\in M(\{\goth K,\Phi\})$. Then $M'$ exists, and $M'\in M(\{\goth K,\Phi'\})$. More exactly, if $F\in F$(\{$\goth K$, $\Phi$\}) then $F^{-1}\otimes \Cal D^{-1}\in F$(\{$\goth K$, $\Phi'$\}) where $\Cal D$ is the different
sheaf on $X$, and if $M=\delta(F)$ then $M'=\delta(F^{-1}\otimes \Cal D^{-1})$.
\medskip
{\bf Proof.} Let $G$ be any invertible sheaf on $X$. We have a
\medskip
{\bf Lemma 12.7.0.} There exists the canonical isomorphism $\varphi_G:
\pi_*(G^{-1}\otimes \Cal D^{-1}) \to \Hom_{P^1}(\pi_*(G), \Cal O)$.
\medskip
{\bf Proof.} At the level of affine open sets $\varphi_G$ comes from the trace bilinear form of field extension $\goth K/\n F_q(T)$. Concordance with glueing is obvious. $\square$
\medskip
We need the relative version of this lemma. Let $G_1$, $G_2$ be invertible sheaves on $X$,
$\rho: G_1 \to G_2$ any rational map. Obviously there exists a rational map
$\rho^{-1}: G_1^{-1} \to G_2^{-1}$. Recall that we denote by
$\rho^{inv}: G_2 \to G_1$ the rational map which is inverse to $\rho$
respectively the composition. The map $\pi_*(\rho^{-1}\otimes \Cal D^{-1}): \pi_*(G_1^{-1}\otimes \Cal D^{-1})\to \pi_*(G_2^{-1}\otimes \Cal D^{-1})$ is obviously defined. The map (denoted by $\beta(\rho)$) from $\Hom_{P^1}(\pi_*(G_1), \Cal O)$ to $\Hom_{P^1}(\pi_*(G_2),
\Cal O)$ is defined as follows at the level of affine open sets: let
$\gamma\in\Hom_{P^1}(\pi_*(G_1), \Cal O)(U)$ where $U$ is a sufficiently
small affine subset of $P^1$, such that we have a map $\gamma(U): \pi_*(G_1)(U)
\rightarrow \Cal O(U)$. Then $(\beta(\gamma))(U)$ is the composition
map $\gamma(U)\circ \pi_*(\rho^{inv})(U)$:
$$\pi_*(G_2)(U)\overset{\pi_*(\rho^{inv})(U)}\to{\longrightarrow}
\pi_*(G_1)(U)\overset{\gamma(U)}\to{\to}\Cal O(U)$$

{\bf Lemma 12.7.1.} The above maps form a commutative diagram:
$$\matrix
\pi_*(G_1^{-1}\otimes \Cal D^{-1})&\overset{\pi_*(\rho^{-1}\otimes \Cal
D^{-1})}\to{\longrightarrow}&\pi_*(G_2^{-1}\otimes \Cal D^{-1})&&\\&&&&\\
\varphi_{G_1}\downarrow&&\varphi_{G_2}\downarrow &&\\&&&&\\
\Hom_{P^1}(\pi_*(G_1), \Cal
O)&\overset{\beta(\rho)}\to{\longrightarrow}&\Hom_{P^1}(\pi_*(G_2),
\Cal O)&&\square
\endmatrix$$

We apply this lemma to the case $\{\rho: G_1 \to G_2\}=\{\tau_{X, \Phi}:
F^{(1)}\to F\}$. We have:
$$\Div(\tau_{X, \Phi}^{-1}\otimes \Cal D^{-1})=-\Div(\tau_{X, \Phi})=-\theta_{\Phi}+n\infty$$
Futher, we multiply $\tau_{X, \Phi}^{-1}\otimes \Cal D^{-1}$ by $T-\theta$. We have:
$$\Div((T-\theta)\tau_{X, \Phi}^{-1}\otimes \Cal
D^{-1})=\Div(T-\theta)+\Div(\tau_{X, \Phi}^{-1}\otimes \Cal D^{-1})=\theta_{\Phi'}-(r-n)\infty$$
i.e. $(T-\theta)\tau_{X, \Phi}^{-1}\otimes \Cal D$ is one of
$\tau_{X,\Phi'}$, i.e. $F^{-1}\otimes \Cal D^{-1}\in F$(\{$\goth K$, $\Phi'$\}). Further, $(T-\theta)\beta(\tau_{X, \Phi})$ is the map which is used in the definition of duality of
$M$. This means that the lemma 12.7.1 implies the theorem. $\square$
\medskip
{\bf Remark 12.8.} There exists a simple proof of the second part of the Theorem 5 for uniformizable abelian
t-motives $M$ with complete multiplication by $\w_\goth K\subset \goth K$. Recall that this second part is the proof of 2.7 for $M$. Really, let us
consider the diagram 2.5. The CM-types
of action of $\goth K$ on $\Lie(M)$ and on $E(M)$ coincide,
and the CM-types of action of $\goth K$ on a vector space and on its
dual space coincide. This means that the CM-type of $V^*$ is $\Phi$
and the CM-type of $V'$ is $\Phi'$. Further, $\gamma_D$ of 2.5 commutes with complete multiplication: this follows immediately for example from a description of $\gamma_D$ given in Remark 5.2.8. Really, all homomorphisms of 5.2.9 commute with complete multiplication. For example, this condition for $\delta$ of 1.11.1 is written as follows: if $k\in \goth K$, $\goth m_k(M)$, resp. $\goth m_k(M')$ is the map of complete multiplication by $k$ of $M$, resp. $M'$, then $(\goth m_k(M)\otimes Id) \circ \delta= (Id \otimes \goth m_k(M'))\circ \delta$ --- see any textbook on linear algebra.

Finally, since
$\Phi\cap \Phi'=\emptyset$ and the map $\vf'\circ \gamma_D \circ
\vf^*$ commutes with complete multiplication, we get that it must be
0.
\medskip
{\bf 13. Miscellaneous.}

\nopagebreak
\medskip
Let now $(L,V)$ be from 12.1, case $\goth q=m=1$, i.e. $\goth K=\Cal W$
and $r=g$, and let the ring of complete multiplication be the maximal
order $\w_\goth K$. We identify $\w$ and $\n Z_\infty$ via $\iota$, i.e. we consider $\goth K$ as an extension of $\n F_q(\theta)$. Let $\Phi$ be the CM-type of the action of $\goth
K$ on $V$.
This means that --- as an
$\w_\goth K$-module --- $L$ is isomorphic to $I$ where $I$ is an ideal of
$\w_\goth K$. The class of $I$ in $\Cl(\w_\goth K)$ is defined by $L$
and $\Phi$ uniquely; we denote it by $\Cl(L,\Phi)$.
\medskip
{\bf Remark.} $\Cl(L,\Phi)$ depends on $\Phi$, because the action of $\w_\goth K$ on $V$ depends on $\Phi$. Really, let $a\in L \subset V$, $a=(a_1,...,a_n)$ its coordinates, $\Phi=\{\alpha_{i_1},...,\alpha_{i_n}\}\subset \{\alpha_{1},...,\alpha_{r}\}$ and $k\in \w_\goth K$. Then $ka$ has coordinates $(\alpha_{i_1}(k)a_1,...,\alpha_{i_n}(k)a_n)$, i.e. depends de $\Phi$. Particularly, the $\w_\goth K$-module structure on $L$ depends on $\Phi$, and hence $\Cl(L,\Phi)$ depends on $\Phi$. For example, if $n=1$, $r=2$, $\Phi_1=\{\alpha_{1}\}$, $\Phi_2=\{\alpha_{2}\}$, then $\Cl(L,\Phi_2)$ is the conjugate of $\Cl(L,\Phi_1)$.
\medskip
{\bf Theorem 13.1.} $\Cl(L',\Phi')=(\Cl(\goth d))^{-1}
(\Cl(L,\Phi))^{-1}$ where
$\goth d$ is the different ideal of the ring extension $\w_\goth K/\w$.
\medskip
{\bf Proof.} This theorem follows from the above results;
nevertheless, I give here an explicit elementary proof. Let $a_*=(a_1, ..., a_r)^t$ be a
basis (considered
as a vector column) of $\goth K$ over $\n F_q(\theta)$ and $b_*=(b_1,
..., b_r)^t$ the dual basis. Recall that it satisfies 2 properties:
$$(1) \ \ \forall i\ne j \ \ \alpha_i(a_*)^t\alpha_j(b_*)=0 \ \ \
(\hbox{ i.e. } \sum_{k=1}^r \alpha_i(a_k) \alpha_j(b_k)=0
)\eqno{(13.1.1)} $$
(2) For $x\in \goth K$ let $\goth m_{x,a_*}$ (resp $\goth m_{x,b_*}$)
be the matrix of multiplication by $x$ in the basis $a_*$ (resp.
$b_*$). Then for all $x\in \goth K$ we have
$$\goth m_{x,a_*}=\goth m_{x,b_*}^t\eqno{(13.1.2)} $$

We define $\goth I_{n,r-n}$
as an $r\times r$ block matrix $\left(\matrix 0&I_{r-n}\\ -I_n&0
\endmatrix\right)$, and we define a new basis $\tilde b_*=(\tilde b_1,
..., \tilde b_r)^t$ by
$$\tilde b_*=\goth I_{n,r-n}b_*\eqno{(13.1.3)} $$
(explicit formula: $(\tilde b_1, ..., \tilde b_r)=(b_{n+1},
..., b_r, -b_1, ..., -b_n)$).

We can assume that $\Phi=\{\alpha_1, ... , \alpha_n \}$. Since $L$ has multiplication by $\w_\goth K$ and the CM-type of this
multiplication is $\Phi$, it is possible to choose $a_*$ such that $L\subset
\p^{n}$ is generated over $\n Z_\infty$ by $e_1, ... , e_r$ where
$$e_i=(\alpha_1(a_i), ... , \alpha_n(a_i))\eqno{(13.1.4)}$$ Let $\hat L\subset \p^{r-n}$ be
generated over $\n Z_\infty$ by $\hat e_1, ... , \hat e_r$ where
$$\hat e_i=(\alpha_{n+1}(\tilde b_i), ... , \alpha_r(\tilde b_i))\eqno{(13.1.5)}$$

{\bf Lemma 13.1.6. } $L'=\hat L$.
\medskip
{\bf Proof. } Let $A$ (resp. $B$) be a matrix whose lines are the lines of
coordinates of $e_1, ... , e_n$ (resp. $e_{n+1}, ... ,
e_r$) in 13.1.4, and $C$ (resp. $D$) a matrix whose lines are the lines of
coordinates of $\hat e_1, ... , \hat e_{r-n}$ (resp. $\hat e_{r-n+1}, ... ,
\hat e_r$) in 13.1.5. By definition of Siegel matrix, we have $L=\goth L(BA^{-1})$,
$\hat L=\goth L(DC^{-1})$ ($\goth L$ is defined in 3.1, 3.2). So, it is sufficient to prove that
$(BA^{-1})^t=DC^{-1}$, i.e. $A^tD=B^tC$. This follows immediately from
the definition of $A,B,C,D$ and (13.1.1). $\square$
\medskip
For $x\in \w_\goth K$ we denote by $\goth
M_x(L)$ the matrix of multiplication by $x$ in the basis $e_*$ (see
the notations of Remark 3.8). Obviously $\goth M_x(L)=\goth
m_{x,a_*}$.

Let now $\w_\goth K$ acts on $\p^{r-n}$ (the ambient space of $L'$) by CM-type
$\Phi'$. According (13.1.2) and
(13.1.3), the matrix of the action of $x\in \w_\goth K$ in the basis
$\tilde b_*$ is
$$\goth I_{n,r-n}\goth m_{x,a_*}^t\goth I_{n,r-n}^{-1}\eqno{(13.1.7)} $$
Let $\goth M$, $\goth M'$ be from Remark 3.8. Formula 3.8.4
shows that
$$\goth M'=\goth I_{n,r-n}\goth M^t\goth I_{n,r-n}^{-1}\eqno{(13.1.8)} $$
Formulas (13.1.7) and (13.1.8) --- because of Lemma 1.10.3 --- prove the theorem. $\square$
\medskip
{\bf 13.2. Compatibility with the weak form of the main theorem of
complete multiplication.}

\nopagebreak
\medskip
The reader can think that Theorem 13.1 is incompatible with the main theorem of complex
multiplication, because of the $-1$-th power in its statement. The reason is a bad choice
of notations of Shimura, he affirms that the CM-type of an abelian variety $A$ over a number field coincides with the CM-type
of $A'$, while we see that it is really the complement. Since an analog of even the weak form of the main theorem of complex multiplication --- Theorem 13.2.6 --- for the function field case is not proved yet, the main result of the present section --- Theorem 13.2.8 --- is conditional: it affirms that if this weak form of the main
theorem --- Conjecture 13.2.7 --- is true for a t-motive with complete multiplication $M$, then it is true for $M'$ as well. By the way, even if it will turn out that the statement of the Conjecture 13.2.7 is not correct, the proof of 13.2.8 will not be affected, because the main
ingredient of the proof is the formula 13.2.10 "neutralizing" the $-1$-th power of the
Theorem 13.1.

Let us recall some definitions of [Sh71], Section 5.5. We consider an abelian variety
$A=\n C^n/L$ with
complex multiplication by $K$. The set $\Hom(K, \bar \n Q)$ consists of $n$ pairs of
mutually
conjugate inclusions $\{\vf_1,\bar \vf_1, ..., \vf_n,\bar\vf_n\}$. $\Phi$ is a subset of
the set
$\Hom(K, \bar \n Q)$ such that $\forall i=1,...,n$ we have:
$$\hbox{$\Phi\cap\{\vf_i,\bar \vf_i\}$ consists of one element.}\eqno{(13.2.1)}$$
It is defined by the condition that the action of complex multiplication
of $K$ on $\n C^n$ is isomorphic to the direct sum of the elements of $\Phi$. Let $F$ be
the
Galois envelope of $K/\n Q$,
$$G:=\Gal (F/\n Q), \ \ H:=\Gal(F/K), \ \
S:=\bigcup_{\alpha\in\Phi}H\alpha\eqno{(13.2.2)}$$
(the elements of Galois group act on $x\in F$ from the right, i.e. by the formula
$x^{\alpha\beta}=(x^\alpha)^\beta$; for $\alpha\in\Phi$ we denote by $\alpha$ also a
representative in $G$ of the coset $\alpha$). We denote
$$H^{ref}:=\{\gamma\in G\vert S\gamma=S\}\eqno{(13.2.3)}$$
and let $K^{ref}$ be the subfield of $F$ corresponding to $H^{ref}$. We have:
$$H^{ref}S^{-1}=S^{-1}\eqno{(13.2.4)}$$
i.e. $S^{-1}$ is an union of cosets of $H^{ref}$ in $G$. We can identify these cosets
with elements of $\Hom(K^{ref}, \bar \n Q)$. $\Phi^{ref}\subset \Hom(K^{ref}, \bar \n Q)$
is, by definition, the set of these cosets. There is a map $\det
\Phi^{ref}:{K^{ref}}^\times \to K^\times$ defined as follows:
$$\det \Phi^{ref}(x):=\prod_{\alpha\in\Phi}\alpha(x)\eqno{(13.2.5)}$$
(it follows easily from the above formulas and definitions that $\det \Phi^{ref}(x)$
really belongs to $K^\times$). It can be extended to the group of ideles and factorized to
the group of classes of ideals, we denote this map
by $\det_{Cl} \Phi^{ref}: \Cl(K^{ref}) \to \Cl(K)$. Finally, let
$\theta^{ref}: \Gal(K^{ref\ Hilb}/K^{ref}) \to \Cl(K^{ref})$ be an isomorphism defined by
the Artin reciprocity law.

We consider the case $\End(A)=O_K$. In this case $L$ is isomorphic to an ideal of $O_K$,
its class
is well-defined by the class of isomorphism of $A$, we denote it by $\Cl(A)$.
\medskip
{\bf Theorem 13.2.6.} $A$ is defined over $K^{ref\ Hilb}$;

For any $\gamma\in\Gal(K^{ref\ Hilb}/K^{ref})$ we have

$\Cl(\gamma(A))=\det_{Cl} \Phi^{ref}\circ \theta^{ref} (\gamma)^{-1}(\Cl(A))$. $\square$
\medskip
This is a weak form of [SH71], Theorem 5.15 --- the main theorem of
complex multiplication.
\medskip
Now we define analogous objects for the function field case in order to formulate a
conjectural
analog of Theorem 13.2.6. Let $\goth K$, $\Phi$ be from 12.5.3. $\goth K^{ref}$,
$\Phi^{ref}$,
$\det \Phi^{ref}$ are defined by the same formulas 13.2.2 -- 13.2.5 like in the number
field case
($\n Q$ must be replaced by $\n F_q(T)$). The facts that 13.2.1 has no meaning in the
function field case and
that the order of $S$ is not necessarily the half of the order of $G$ do not affect the
definitions.

The $\infty$-Hilbert class field of $\goth K$ (denoted by $\goth K^{Hilb \ \infty}$) is an abelian
extension of
$\goth K$ corresponding to the subgroup
$$\goth K_\infty^*\cdot \prod_{v\ne\infty}O_{\goth K_v}^*\cdot \goth K^*$$
of the idele group of $\goth K$. We have an isomorphism $\theta: \Gal(\goth K^{Hilb \ \infty}/\goth K) \to
\Cl(\w_{\goth K})$.

We formulate the function field analog of Theorem 13.2.6 only for the case when
\medskip
{\bf (*)} There exists only one
point over $\infty\in P^1(\n F_q)$ in the extension $\goth K^{ref}/\n F_q(T)$.
\medskip
In this case the field
$\goth K^{ref \ Hilb \ \infty}$ and the ring $\w_{\goth K^{ref}}$ are naturally defined, and we have
an
isomorphism
$\theta^{ref}: \Gal(\goth K^{ref \ Hilb \ \infty}/\goth K^{ref}) \to \Cl(\w_{\goth K^{ref}})$.

Let $M$ be an uniformizable
t-motive of rank $r$ and dimension $n$ having complete multiplication
by $\w_\goth K$, and $\Phi$ its CM-type. $\Cl(M)$ is defined like $\Cl(A)$ in the number
field case,
it is $\Cl(L,\Phi)$ of 13.1.
\medskip
{\bf Conjecture 13.2.7.} If (*) holds, then $M$ is defined over $\goth K^{ref\ Hilb \
\infty}$, and for any
$\gamma\in\Gal(\goth K^{ref\ Hilb \ \infty}/\goth K^{ref})$
we have $\Cl(\gamma(M))=\det_{Cl} \Phi^{ref}\circ \theta^{ref} (\gamma)^{-1}\Cl(M)$.
\medskip

Now we can formulate the main theorem of this section.
\medskip
{\bf Theorem 13.2.8.} If conjecture 13.2.7 is true for $M$ then it is true for $M'$.
\medskip
{\bf Proof.} It follows immediately from the functional analogs of 13.2.2 -- 13.2.4 that
$$(\goth K,\Phi')^{ref}=(\goth K^{ref},(\Phi^{ref})')\eqno{(13.2.9)}$$
Further,
$$\hbox{det}_{Cl} {\Phi'}^{ref}=(\hbox{det}_{Cl} \Phi^{ref})^{-1}\eqno{(13.2.10)}$$
Really, $\det \Phi^{ref}(x) \cdot \det (\Phi^{ref})'(x) = N_{\goth
K^{ref}/\n F_q(T)}(x)\in\n F_q(T)^\times$, hence gives the trivial class of ideals (we use
here (13.2.9).
Finally, for $\gamma\in\Gal(\goth K^{ref})$ we have
$$(\gamma(M))'=\gamma(M')\eqno{(13.2.11)}$$
The theorem follows immediately from 13.1, 13.2.10, 13.2.11 (recall that $\Cl(M)$ is
$\Cl(L,\Phi)$ of 13.1). $\square$
\medskip
{\bf 13.3. Some explicit formulas.} We give
here an elementary explicit proof of the theorem 12.7 in two simple
cases: $\goth K=\n F_{q^r}(T)$ and $\n F_q(T^{1/r})$. By the way,
since the extension $\n F_{q^r}(T)/\n F_{q}(T)$ is not absolutely
irreducible, formally this case is not covered by the theorem 12.7.
\medskip
{\bf Case $\w_\goth K= \n F_{q^r}[T]$.} Let $\alpha_i$, where $i=0, ... ,r-1$, be inclusions $\goth K \to \p$. For $\omega \in\n F_{q^r}$
 we have $\alpha_i(\omega )=\omega^{q^{i}}$. Let $$0 \le
i_1 < i_2 < ... <
i_n \le r-1\eqno{(13.3.0)}$$ be numbers such that $\Phi=\{\alpha_{i_j}\}$, $j=1, ... , n$.
We consider the following t-motive $M=M(\goth K, \Phi)$. Let $e_1, ... , e_n$ be a basis of $M_{\p[\tau]}$ such that $\goth
m_\omega(e_j)=\omega^{q^{i_j}}e_j$ and the multiplication by $T$ is
defined by formulas
$$Te_1=\theta e_1 + \tau^{i_1-i_n+r}e_n\eqno{(13.3.1)} $$
$$Te_j=\theta e_j + \tau^{i_j-i_{j-1}}e_{j-1}, \ \ j=2, ... , n
\eqno{(13.3.2)} $$

It is easy to check that $M$ has complete multiplication by $\w_\goth K$, and its CM-type is $\Phi$.
\medskip
{\bf Remark.} It is possible to prove that $M(\goth K, \Phi)$ is the only t-motive having these properties; we omit the proof.
\medskip
{\bf Proposition 13.3.3.} For $\w_\goth K= \n F_{q^r}[T]$ we have: $M(\goth K, \Phi)'=M(\goth K, \Phi')$.
\medskip
{\bf Proof.} Elements $\tau^je_k$ for $k=1, ... ,n$, $j=0, ..., i_{k+1}-i_{k}-1$ for $k<n$
and $j=0,
..., i_1-i_n+r-1$ for $k=n$ form a basis of $M_{\p[T]}$. Let us
arrange these elements in
the lexicographic order ($\tau^{j_1}e_{k_1}$ precedes to
$\tau^{j_2}e_{k_2}$ if $k_1 <
k_2$) and make a cyclic shift of them by $i_1$ denoting $e_1$ by $f_{i_1+1}$,
$\tau^{i_2-i_1-1}e_1$ by $f_{i_2}$ etc. until
$\tau^{i_1-i_n+r-1}e_n=f_{i_1}$. Formulas
13.3.1, 13.3.2 become
$$\tau(f_i)=f_{i+1} \hbox{ if } i\not\in \{i_1, ..., i_n\} $$
$$\tau(f_i)=(T-\theta)f_{i+1} \hbox{ if } i\in \{i_1, ..., i_n\} $$
($i \mod r$, i.e. $f_{r+1}=f_1$). Formula 1.10.1 shows that in the
dual basis $f'_*$ we
have
$$\tau(f'_i)=f'_{i+1} \hbox{ if } i\in \{i_1, ..., i_n\} $$
$$\tau(f'_i)=(T-\theta)f'_{i+1} \hbox{ if } i\not\in \{i_1, ..., i_n\} $$
which proves the proposition. $\square$
\medskip
{\bf Case $\w_\goth K=\n F_q[T^{1/r}]$, $(r,q)=1$.} In order to define $M(\goth K, \Phi)$ we need more notations. We denote
$\theta^{1/r}$ and $T^{1/r}$ by $\goth s$ and $S$ respectively, and
let $\zeta_r$ be a primitive $r$-th root of 1. Let $\alpha_i$, $i_1 < i_2 < ... < i_n$ and $\Phi$ be the
same as in the case $\w_\goth K= \n F_{q^r}[T]$. We
have $\alpha_i(S)=\zeta_r^iS$. Further, we consider an overring $\p[S,\tau]$ of $\p[T,\tau]$ ($S$ is in the
center of this ring), and we consider the category of modules
over $\p[S,\tau]$ such that the condition 1.9.2 is changed by a
weakened condition 13.3.4 (here $A_{S,0}\in M_n(\p)$ is defined by the formula $Se_*=A_Se_*$, where $A_S\in M_n(\p)[\tau]$, $A_S=\sum_{i=0}^*A_{S,i}\tau^i$):
$$A_{S,0}^r=\theta I_n + N\eqno{(13.3.4)}$$
Let $\bar M$ be a $\p[S,\tau]$-module such that $\dim \bar
M_{\p[S]}=1$, $f_1$ the only element of a basis of $\bar M_{\p[S]}$ and
$$\tau f_1=(S-\zeta_r^{i_1}\goth s)\cdot ... \cdot
(S-\zeta_r^{i_n}\goth s)f_1$$
By definition, $M=M(\goth K, \Phi)$ is the restriction of scalars from
$\p[S,\tau]$ to $\p[T,\tau]$ of $\bar M$. Like in the case $\w_\goth K= \n F_{q^r}[T]$, it is easy to check that $M$ has complete multiplication by $\w_\goth K$ with CM-type $\Phi$, and it is possible to prove that it is the only t-motive having these properties.
\medskip
{\bf Proposition 13.3.5.} For $\w_\goth K= \n F_q[T^{1/r}]$, $(r,q)=1$ we have: $M(\goth K, \Phi)'=M(\goth K, \Phi')$.
\medskip
{\bf Proof.} For $i=1, ..., r$ we denote $f_i=S^{i-1}f_1$. These $f_*=f_*(\Phi)$
form a basis of $M_{\p[T]}$, and the matrix $Q=Q(f_*,\Phi)$ of
multiplication of $\tau$ in this basis has the following description.
We denote by $\sigma_k(\Phi)$ the elementary symmetric polynomial
$\sigma_k(\zeta_r^{i_1}, ... , \zeta_r^{i_n})$.

The first line of $Q$ is
$$\sigma_n(\Phi)\goth s^n \ \ \ \sigma_{n-1}(\Phi)\goth s^{n-1} \ \ \ ...
\ \ \  \sigma_1(\Phi)\goth s \ \ \ 1 \ \ \ 0 \ \ \ ... \ \ \ 0$$
and its $i$-th line is obtained from the first line by 2 operations:

1. Cyclic shift of elements of the first line by $i-1$ positions to the right;

2. Multiplication of the first $i+n-r$ elements of the obtained line by $T$.

We consider another basis $g_*=g_*(\Phi)$ of $M_{\p[T]}$ obtained by
inversion of order of $f_i$, i.e. $g_i=f_{r+1-i}$. The elements of
$Q(g_*)$ are obtained
by reflection of positions of elements of $Q(f_*)$ respectively the
center of the matrix.

The theorem for the present case follows from the formula
$$Q(f_*,\Phi)Q(g_*,\Phi')^t=(T-\theta)I_r$$
whose proof is an elementary exercise: let $\Phi'=\{j_1, ... ,
j_{r-n}\}$; we apply
equality
$$\sigma_k(x_1, ... x_r)=\sum_l \sigma_l(x_{i_1}, ... ,
x_{i_n})\sigma_{k-l}(x_{j_1}, ...
, x_{j_{r-n}})$$
to $1, \zeta_r, ... , \zeta_r^{r-1}$. $\square$
\medskip
{\bf 13.4. Reduction.} Recall notations of 1.16.  Let $L$ be a finite extension of $\n F_q(\theta)$, $\goth p$ a valuation of $L$ over a valuation $P\ne\infty$ of $\n F_q(\theta)$, and we denote $\iota^{-1}(P)\subset \w$ by $\Cal P$. Let $M$ be a t-motive defined over $L$ having a good ordinary reduction $\tilde M$ at $\goth p$ and such that the dual $M'$ exists. According 1.15.1, the $L$-structure on $M'$ is well-defined. We denote by $M_{\Cal P,0}$ the kernel of the reduction map
$M_\Cal P \to \tilde M_\Cal P$. Condition of ordinarity means that  $M_{\Cal P,0}=(\w/\Cal P)^n$.
\medskip
{\bf Conjecture 13.4.1.} For the above $M$, $M'$ we have:
\medskip
$M_{\Cal P,0}$ and $M'_{\Cal P,0}$ are mutually dual with respect to the pairing of Remarks 4.2, 5.1.6 (recall that conjecturally $M'$ also has good ordinary reduction at $\goth p$).
\medskip
{\bf Proof for a particular case:} $M$ is a Drinfeld module, $\Cal P=T$.
\medskip
(1.9.1) for $M$ has a form
$$Te=\theta e+a_1\tau e+... a_{r-1}\tau^{r-1}e+\tau^r e$$
Condition of good ordinary reduction means $a_i\in L$, $\ord_\goth p(a_i)\ge 0$,
$\ord_\goth p a_1=0$. Let $x\in
M_T$, $y\in M'_T$; we can consider $x$ (resp. $y$) as an element of
$\p$ (resp. $\p^{r-1}$) satisfying some polynomial equation(s). Considering
Newton polygon of these polynomials we get immediately (1) for both
$M$, $M'$. Let $y=(y_1, ... ,y_{r-1})$ be the coordinates of $y$; explicit
formula (5.3.5) for the present case has the form
$$<x,y>_M=\Xi(xy_{r-1}^q+x^qy_1+x^{q^2}y_2+ ... +x^{q^{r-1}}y_{r-1})$$
The same consideration of the Newton polygon of the above polynomials
shows that for $x\in M_{T,0}$, $y\in M'_{T,0}$ we have
$\ord_\goth p x$, $\ord_\goth p y_i \ge 1/(q-1)$. Since $\ord_\goth p \Xi = -1/(q-1)$ we get
that $\ord_\goth p (<x,y>_M)>0$ and hence (because $<x,y>_M\in \n F_q$) we have
$<x,y>_M=0$. Dimensions of $M_{T,0}$, $M'_{T,0}$ are
complementary, hence they are mutually dual. $\square$
\medskip
{\bf Remark 13.4.2.} Analogous explicit proof exists for any standard-3 $M$ of
Section 11.8.
\medskip
\medskip
{\bf References}

\nopagebreak
\medskip
[A] Anderson Greg W. $t$-motives. Duke Math. J. Volume 53, Number 2 (1986), 457 -- 502.
\medskip
[BH] Matthias Bornhofen, Urs Hartl. Pure Anderson Motives and Abelian $\tau$-Sheaves. arXiv:0709.2809
\medskip
[F] Faltings, Gerd. Group schemes with strict $\Cal O$-action. Mosc.
Math. J.  2
(2002),  no. 2, 249 -- 279.
\medskip
[G] Goss, David. Basic structures of function
field arithmetic.
\medskip
[GL17] Grishkov A., Logachev, D. Lattice map for Anderson t-motives: first approach.
J. Number Theory 180 (2017), 373 -- 402. 

http://arxiv.org/pdf/1109.0679.pdf
\medskip
[GL18] Grishkov A., Logachev, D. Duality of Anderson t-motives having $N\ne0$.
https://arxiv.org/pdf/1812.11576.pdf
\medskip
[H] Urs Hartl. Uniformizing the Stacks of Abelian Sheaves.

http://arxiv.org/abs/math.NT/0409341
\medskip
[L09] Logachev, D. Anderson t-motives are analogs of abelian varieties with multiplication by imaginary quadratic fields.
http://arxiv.org 0907.4712.
\medskip
[L] Logachev, D. Reductions of Hecke correspondences on Anderson varieties. In preparation.
\medskip
[P] Richard Pink, Hodge structures over function fields. Universit\"at Mannheim.
Preprint. September 17, 1997.
\medskip
[Sh63] Shimura, Goro. On analytic families of polarized abelian
varieties and automorphic functions. Annals of Math., 1 (1963), vol.
78, p. 149 -- 192
\medskip
[Sh71]  Shimura, Goro. Introduction to the
arithmetic theory of
automorphic functions.
\medskip
[Sh98]  Shimura, Goro. Abelian varieties with
complex multiplication
and modular functions. Princeton Mathematical Series, 46.
\medskip
[Tae] Taelman Lenny, Artin t-motifs. J. Number Theory, 129 (2009), 142 - 157
\medskip
[T] Taguchi, Yuichiro. A duality for finite $t$-modules.
J. Math. Sci. Univ. Tokyo 2 (1995), no. 3, 563--588.

\enddocument

and $\mu$ a number such that the $\mu$-dual ${M'}^{\mu}$ exists. Let $m$ be the minimal number such that $N^m=0$.
\medskip
{\bf Conjecture 6B1.} 1. $\mu \ge m$, i.e. $N^\mu=0$.
\medskip
2. ${M'}^{\mu}$ is uniformizable. There exists a canonical perfect $\w$-valued pairing between $L_T(M)$ and $L_T({M'}^{\mu})$.
\medskip
This pairing extends by $\p[[T-\theta]]$-linearity to the $\p[[T-\theta]]$-valued pairing between $L_T(M)\underset{\w}\to{\otimes}\p[[T-\theta]]$, $L_T({M'}^{\mu})\underset{\w}\to{\otimes}\p[[T-\theta]]$. We denote it by $<.,.>$.
\medskip
{\bf Conjecture 6B2.} We have:
$$\goth q({M'}^{\mu})=\{x\in L_T({M'}^{\mu})\underset{\w}\to{\otimes}\p[[T-\theta]]  \hbox{ such that }  \forall y\in \goth q(M) $$ $$\hbox{ we have } <x,y>\in (T-\theta)^\mu\p[[T-\theta]]\} $$

Let us consider uniformizable $M$ having $N=0$ (i.e. $m=1$) such that $M'$ --- the dual of $M$ --- exists and has $N'=0$.
\medskip
{\bf Corollary 6B3.} For this case Conjecture 6B2 is true, i.e.
$$\goth q({M'})=\{x\in L_T({M'})\underset{\w}\to{\otimes}\p[[T-\theta]]  \hbox{ such that }  \forall y\in \goth q(M) $$ $$\hbox{ we have } <x,y>\in (T-\theta)\p[[T-\theta]]\} $$
\medskip
This follows immediately from Theorem 5 and equivalence of Definition 2.3 and Property 2.4 (we consider the reduction of the above pairing $<.,.>$ modulo the maximal ideal $ (T-\theta) \p[[T-\theta]]$ of $\p[[T-\theta]]$; this reduction coincides with the pairing of Section 5).

\newpage
{\bf 1.16.4. $\Cal P$-rank of $M$ and of $M'$.} By analogy with the number field case, the $\Cal P$-rank of $M$ is the dimension of the $\Cal P$-torsion points of $E(M)$ over $\w/\Cal P$; it varies from $r-n$ (ordinary $M$) to 0 (completely suresingular $M$).
\medskip
{\bf Conjecture 1.16.5.} $M$ is ordinary $\iff M'$ is ordinary.
\medskip
For standard-3 t-motives  apparently this can be shown by explicit calculations. See 13.4.1 for the proof for the case of Drinfeld modules.
\medskip
{\bf Question 1.16.6.} What are possible values of $\Cal P$-rank of $M$, $\Cal P$-rank of $M'$? Is it true that the pair of numbers ($\Cal P$-rank of $M$, $\Cal P$-rank of $M'$) characterizes completely (in some meaning) the type of $M$? Particularly, whether the $\Cal P$-ranks of other tensor operations of $M$ (exterior powers etc.) are defined completely by ($\Cal P$-rank of $M$, $\Cal P$-rank of $M'$), or not?
\medskip
{\bf Example 1.16.7.} Let us consider $M$ defined by (8.2.2), entries of $A$ belong to $\bar \n F_q$, and let $\Cal P=T$. According Lang's Theorem, $T$-rank of $M$ is $n \iff \det A\ne 0$. An explicit calculation for the case $n=2$ shows that if the entries of $A$ belong to $\n F_q$ then the $T$-rank of $M$ is equal to the quantity of the non-zero eigenvalues of $A$ (if the entries of $A$ do not belong to $\n F_q$ then the formula for the $T$-rank of $M$ is more complicated). The dual $M'$ is defined by the same (8.2.2) with $A$ replaced by $-A^t$; this means that  if the entries of $A$ belong to $\n F_{q}$ then the $T$-rank of $M$ is equal to the $T$-rank of $M'$. It is easy to check that the same is true if the entries of $A$ belong to $\n F_{q^2}$ but not to a larger field:  for example, if $\gamma\in \bar \n F_q-\n F_{q^2} $ then for $A=\left(\matrix \gamma & 1 \\ -\gamma^{q+1}& -\gamma^{q}\endmatrix \right)$ we have $T$-rank of $M$ is 0, $T$-rank of $M'$ is 1.
\medskip
{\bf Example 1.16.8.} For pure non-ordinary standard-3 $M$ having $r=5$, $n=2$ we get easily that the case $T$-rank of $M$ is 2, $T$-rank of $M'$ is 0 cannot be realized; all other possible cases can be realised.
\medskip
{\bf Remark 13.4.3: Case $N\ne0$.} I do not know a definition of ordinarity for $M$ in finite characteristic for the case $N\ne0$; we can expect that these are those $M$ whose $\Cal P$-rank is the maximal possible in some family of these $M$. A version of Conjecture 13.4.1 should hold for this case; instead of the dual $M'$ we should consider the $m$-th dual ${M'}^m$ where $m$ satisfies $M^m=0$ (or $m$ is the minimal number with this property?)
\medskip
{\bf Example 13.4.3.1.} Let $M$ be given by the following modified formula (8.2.2):
$$Te_* =( \theta+N) e_* + A \tau e_* + \tau^2 e_*\eqno{(13.4.3.2)}$$
where $n=2$, $N=\left(\matrix 0&1\\0&0 \endmatrix\right)$, $A= \left(\matrix a_{11}&a_{12} \\a_{21}&a_{22} \endmatrix\right)$. For the case of finite characteristic (i.e. $a_{ij}\in \bar\n F_q$) and $\Cal P=T$ we have $M$ is ordinary iff $a_{21}\ne0$, in this case the $\Cal P$-rank of $M$ is 3. We have $m=2$, ${M'}^2$ is given by the formula
$$Te'_*=( \theta+N') e'_* + A' \tau e'_*$$
where $N'=\left(\matrix 0&-1&0&0&0&0\\0&0&0&0&0&0\\0&0&0&1&0&0\\0&0&0&0&0&0
\\0&0&0&0&0&1\\0&0&0&0&0&0\endmatrix\right)$,  $A'=\left(\matrix 0&0&0&0&1&0 \\0&0&1&0&0&0\\0&0&0&0&0&0\\0&1&-a_{11}&0&-a_{21}&0\\0&0&0&0&0&0
\\1&0&-a_{12}&0&-a_{22}&0 \endmatrix\right)$.

An explicit calculation shows that the $\Cal P$-rank of $M$ is $\le 1$ and it is 1 iff $a_{21}\ne0$. So, an analog of 13.4.1 for $M$ is the following:
\medskip
{\bf Conjecture 13.4.3.3.} Let $M$ be given by 13.4.3.2 in the generic characteristic, and let $M$ have good ordinary reduction at $\Cal P=T$.  We have: $M_{T,0}$ and $({M'}^2)_{T,0}$ are mutually dual with respect to the pairing.
\medskip
This conjecture can be easily checked by explicit calculation.

\enddocument

More exactly, if $\vf: M \to N$ is a map of abelian
t-motives and both ${M_{\goth C}'}^{\mu}$, ${N_{\goth C}'}^{\mu}$ exist then there exists
the map of rational pr\'e-t-motives ${\vf_{\goth C}'}^{\mu}: {N_{\goth C}'}^{\mu} \to
{M_{\goth C}'}^{\mu}$.


$\sigma: \p \to \p$ the Frobenius automorphism of $\p$, i.e. $\sigma(z)=z^q$.

We extend $\sigma$ to $K_C$, resp. $\w_C$ by the formula
$\sigma_{K}(k\otimes z)=k\otimes \sigma(z)$, $z\in\p$, $k\in K$, resp. $k\in
\w$.

: we consider $M\mapsto M^{(1)}$ as a
functor; if $M$ is free $\p[T]$-module and $\alpha: M \to M$ is a map whose matrix in
some basis $\{m_*\}$ is $C$, then the matrix of $\alpha^{(1)}: M^{(1)} \to M^{(1)}$ in the
base $\{m_*\otimes1\}$ is $C^{(1)}$

We have an inclusion $\w \hookrightarrow K_\infty$,
we prolonge $\iota$ to the inclusion $K_\infty \hookrightarrow \p$, and we denote
the image $\iota(K_\infty)\subset \p$ by $K_\infty$ as well.
zarhin@math.psu.edu
Dorogoj Yura,

mozhet byt', Vy znaete chto-nibud' po takim voprosam:

1. Rassmotrim ob'ekty: $C$-prostranstvo $V$, $Z$-reshyotka $L$ v nyom i ermitova forma $H$ na $V$, kotoraya ne >>0 i takaya, chto $im H$ celochislenna na $L$. Kakov kriterij togo, chto eti ob'ekty yavlyayutsya realizaciej motiva?

My znaem, chto esli $H >> 0$, to vsegda. Esli net, to po krajnej mere inogda, naprimer esli $X$ - trifold takoj, chto ego $J^3(X)$ - ne abelevo mnogoobrazie, to $H^3(X)$ - kak raz takogo tipa.

Otsyuda vytekaet takoj vopros

2. My znaem, chto abelevo mnogoobrazie - bolee obshchij ob'ekt, chem jacobian krivoj. Rassmotrim trifold $X$ takoj, chto ego $J^3(X)$ - ne abelevo mnogoobrazie, rassmotrim ego motiv $M=H^3(X)$ i sprosim: motivy kakogo tipa obobshchayut $M$ tak zhe, kak abelevo mnogoobrazie obobshchaet jacobian krivoj?

Mozhet byt', lyudi dumali nad etim 30 - 40 let nazad i uvideli, chto na eti voprosy net razumnyh otvetov?

Prois|hozhdenie etih voprosov takoe. God nazad ya poslal Vam stat'yu ob analogii mezhdu t-motivami s $N=0$ i abelevymi mnogoobraziyami s umnozheniem na mnimoe kvadratichnoe pole. Vopros: a kakov analog t-motivov s $n \ne 0$? Mozhet byt', eto motiv? My legko mozhem opredelit' vysheupomyanutye $V$, $L$ i $H$ dlya vneshnej stepeni abelevyh mnogoobrazij s umnozheniem na mnimoe kvadratichnoe pole. Esli by my znali, chto oni - realizaciya motiva, to kak raz eto i byl by motiv, sootvetstvuyushchij t-motivu s $n \ne 0$.

Yours, Dima

 We consider the case $P\ne \infty$, i.e. we can identify $P$ with an irreducible polynomial in $\n Z_\infty$ (up to a scalar factor).

C-MATRICY                 POZOR !

{\bf 8.2. Infinitesimal lattice map and C-lattices.} For $r=2n$, \ $n>1$ the degree of the lattice map is infinite, hence in order to get a 1 - 1 correspondence it is necessary to introduce a new object called a C-lattice. There is a C-lattice map from a neighbourhood of a t-motive $M_0$ to a neighbourhood of a C-lattice $L_{C,0}$ (see below for the definition of $M_0$, $L_{C,0}$) which is 1 - 1.

For the case $(r,n)=1$ I do not see an inicial t-motive $M_0$ such that in its neighbourhood the degree of the lattice map is not 1 (see Remark 8.2.14 for more details), hence we can ask whether in this case the answer to 8.1.2 is yes. I think that this is few likely.
\medskip
{\bf Definition 8.2.1.} A C-lattice is a pair $(L, e_*)$ where
\medskip
$L \subset V=\p^n$ is a lattice and
\medskip
$e_1, ... , e_r$ is a $\n Z_\infty$-basis of $L$.
\medskip
Further, two such pairs $(L; e_1, ... , e_r)$ and $(L'; e'_1, ... ,
e'_r)$ are called equivalent if there exists a $\p$-linear map $\psi:
V\to V$ such that $\psi(L; e_1, ... , e_r) = (L'; e'_1, ... , e'_r)$,
or $L=L'$ and the matrix of the change of basis from $e_1, ... , e_r$
to $e'_1, ... , e'_r$ (which a priori belongs to $GL_r(\n Z_\infty)$ ) belongs to $GL_r(\n F_q)$.
\medskip
The functor of forgetting the basis from C-lattices to lattices is
denoted by $\goth i$. The notion of a Siegel matrix for a C-lattice is
the same as for a lattice.

We consider infinitesimal degree of the lattice map in a neighbourhood
of some distinguished (having many endomorphisms) t-motive $M_0$. We
consider the case $M_0=\goth C_2^{\oplus n}$ (recall that $\goth C_2$ is the Carlitz module over $\n F_{q^2}[T]$; in 7.7 $M_0$ is denoted by $M$).
\medskip
{\bf Remark.} It is easy to see that $M_0$ is the (only) t-motive with complete multiplication by $\n F_{q^{2n}}[\theta]$ with CM-type $Id, \fr^2, \fr^4, ..., \fr^{2n-2}$, see 12.5.3 and 13.3, first case.
\medskip
Let $\omega \in \n F_{q^2} - \n F_q$ be a
fixed element. A Siegel matrix of $M_0$ is $\omega I_n$. We denote the lattice (resp. the C-lattice) corresponding to $\omega I_n$ by $L_0$ (resp. $L_{C,0}$). $L_0$ is the lattice of $M_0$. We consider 5 sets $S_1, ... , S_5$:
\medskip
$S_1$. The set of $n \times n$ matrices $A$.
\medskip
$S_2$. The set of t-motives $M$ given by the equation (see 1.9.1)

$$Te_* = \theta e_* + A \tau e_* + \tau^2 e_*\eqno{(8.2.2)}$$ where $A
\in S_1$, $e_* = (e_1,...,e_n)^t$ (notations of 1.9.1: $\goth A_1=A$, $\goth A_2=I_n$).
\medskip
$S_3$. The set of $n \times n$ Siegel matrices $Z$.
\medskip
$S_4$. The set of C-lattices of rank $r=2n$ in $C_\infty ^n$.
\medskip
$S_5$. The set of lattices of rank $r=2n$ in $C_\infty ^n$.
\medskip
We consider initial elements 0, $M_0$, $\omega I_n$, $L_{C,0}$, $L_0$ of $S_1, ...,S_5$ respectively and open neighbourhoods $U_i\subset S_i$ of these initial elements.
\medskip
{\bf Proposition 8.2.3.} There exist neighbourhoods $U_2$, $U_4$, $U_5$ such that

(a) The restriction of $\goth i$ to $U_4$ gives us an epimorphism $U_4\to U_5$.

(b) there exists a 1 -- 1 map $\mu_{24}$ from $U_2$ to $U_4$ such that $\mu_{25}:=\goth i\circ\mu_{24}$ (see the below diagram 8.2.4) is the lattice map from uniformizable t-motives to lattices. Particularly, for $n>1$ the fibre of $\mu_{25}$ is discrete infinite.
\medskip
{\bf Proof.} It is sufficient to prove that there exists a commutative diagram $$\matrix U_1 &
\overset{\mu_{12}}\to{\to} & U_2 && \\ \\ \mu_{13} \downarrow &
&\mu_{24}\downarrow &
\overset{\mu_{25}}\to{\searrow} &\\ \\ U_3 &
\overset{\mu_{34}}\to{\to} & U_4 & \overset{\goth i}\to{\to} & U_5
\endmatrix \eqno{(8.2.4)}$$ where $\mu_{12}(A)$ is the t-motive defined by
8.2.2, $\mu_{34}(Z)$ is the C-lattice corresponding to a Siegel matrix
$Z$ and $\mu_{13}$ is defined as
follows. We identify $\Lie(\goth C_2)$ with $\p$ and hence $\Lie(M_0)$ with $\p^n$.
We consider the following basis
$l_{0,1}, ... , l_{0,2n}$ of $L_0\subset \Lie(M_0)=\p^n$: $l_{0,i}=(0,...,0,1,0,...,0)$ (1 at the $i$-th place), $l_{0,n+i}=\omega l_{0,i}$, $i=1,...,n$. The Siegel
matrix of this basis is $\omega I_n$. Let $A\in U_1$, $M$ its $\mu_{12}$-image and $L$ its
lattice, i.e. its $\mu_{25}\circ\mu_{12}$-image. For any $l_0\in L_0$
there exists a well-defined $l\in L$ which is close to $l_0$ (because
entries of $A$ are near 0). So, we consider a basis $l_{1}, ... ,
l_{2n}$ of $L$ where any $l_i$ is near the corresponding $l_{0,i}$.
The Siegel matrix corresponding to $l_{1}, ... , l_{2n}$ is exactly
the $\mu_{13}$-image of $M$. By definitions, the outer quadrangle is
commutative.

We denote by $d_{\alpha\beta}$ the degree of
$\mu_{\alpha\beta}$ at a
generic point near the initial element of $U_\alpha$. We must prove that there exists a map $\mu_{24}$ preserving
commutativity, and that $d_{24}=1$.
\medskip
{\bf Lemma 8.2.5.} $d_{13}=1$.
\medskip
{\bf Proof.} We consider maps of the diagram 5.1.3 for $M$ given by 8.2.2. For $Z\in E=\p^n$ we have $m_T(Z)=\theta Z + A Z^{(1)} +
Z^{(2)}$, and for $Z\in \Lie(M)=\p^n$ we have $\Exp(Z)=\Exp_A(Z)=\sum_{i=0}^\infty C_iZ^{(i)}$ where
$C_i=C_i(A)$, $C_0=1$.

For the reader's convenience, we consider only the case $n=1$ (the
general case does not require any new ideas), hence $Z$, $A$ will be
denoted by $z$, $a$ respectively. We denote
$\theta_{ij}=\theta^{q^i}-\theta^{q^j}$. Recall that the exponent for $\goth C_2$ has the form
$$\Exp_0(z)=z+\frac{1}{\theta_{20}}z^{q^2}+\frac{1}{\theta_{42}\theta_{40}}z^{q^4}+...\eqno{(8.2.6)}$$
($C_{2i}(0)=\frac{1}{\prod_{j=0}^{i-1}\theta_{2i,2j}}$). We denote by
$y_0\in \p$ a nearest-to-zero root to $\Exp_0(z)=0$ (this is $\xi$ of $\goth C_2$ in notations of [G]). It is defined up to
multiplication by elements of $\n F_{q^2}^*$, and it generates over $\n F_{q^2}[\theta]$ the
lattice of the Carlitz module $\goth C_2$. We fix one such $y_0$.
If $a$ is sufficiently small then there is the only root to
$\Exp_a(z)=0$ near $y_{0}$, and there is the only root to
$\Exp_a(z)=0$ near $\omega y_{0}$. We denote these roots
by $z=z(a)$, $z'=z'(a)$ respectively, and we denote $z=y_{0}+\delta$,
$z'=\omega y_{0}+\delta'$. $\delta$ (resp. $\delta'$) is a root to the
power series
$$\sum_{i,j=0}^\infty d_{ij}a^i\delta^j=0 \eqno{(8.2.7)}$$
$$\hbox{resp. }\sum_{i,j=0}^\infty d'_{ij}a^i{\delta'}^j=0 \eqno{(8.2.7')}$$
where
$$d_{00}=0, \ \
d_{10}=\frac{y_{0}^{q}}{\theta_{10}}+\frac{y_{0}^{q^3}}{\theta_{31}\theta_{30}}
+ \frac{y_{0}^{q^5}}{\theta_{53}\theta_{51}\theta_{50}}+
\frac{y_{0}^{q^7}}{\theta_{75}\theta_{73}\theta_{71}\theta_{70}}+..., \ \ d_{01}=1$$
$$d'_{00}=0, \ \ d'_{10}=\frac{(\omega
y_{0})^{q}}{\theta_{10}}+\frac{(\omega
y_{0})^{q^3}}{\theta_{31}\theta_{30}} + \frac{(\omega
y_{0})^{q^5}}{\theta_{53}\theta_{51}\theta_{50}}+
\frac{(\omega y_{0})^{q^7}}{\theta_{75}\theta_{73}\theta_{71}\theta_{70}}+..., \ \ d'_{01}=1$$
Moreover, $\delta$, $\delta'$ are the nearest-to-0 roots to (8.2.7),
$(8.2.7')$ respectively. It is clear that $d_{10}, d'_{10}\ne 0$. This
means that the approximate value of $\delta$, $\delta'$ is
$$-d_{10}a, \ \ -d'_{10}a\eqno{(8.2.8)}$$
respectively. Exactly, both $\delta$, $\delta'$ are power series in
$a$ whose first term is given by (8.2.8). This means that
$$z=y_0-d_{10}a+\sum_{i=2}^\infty k_i
a^i$$
$$z'=\omega y_0-d'_{10}a+\sum_{i=2}^\infty k'_i
a^i$$ It is easy to see that $d'_{10}\ne \omega d_{10}$,
hence the Siegel matrix $\goth Z=z^{-1}z'$ (which is a number because $n=1$) is given by the formula $$\goth Z=\omega
+ \sum_{i=1}^\infty l_i a^i\eqno{(8.2.9)}$$ and $l_1\ne 0$. Since (for $n=1$)
$d_{13}$ is the minimal $i$ such that $l_i\ne 0$ we get that
$d_{13}=1$.
\medskip
For $n > 1$ the calculation is the same (this is the main part of the proof, because for the case $n=1$ the result is known). Analog of (8.2.9) is $$\goth Z=\omega I_n + l_1 A + P_{\ge 2}(A)\eqno{(8.2.10)}$$ where $P_{\ge 2}(A)$ is a power series of entries of $A$ such that all its terms have degree $\ge 2$.
Condition $l_1\ne 0$ implies $d_{13}=1$. $\square$
\medskip
In order to simplify notations, we identify until the end of the proof $L_T(M)$ and $L(M)$ via $\goth a$, and $\n Z_\infty$ and $\w$ via $\iota$. We denote the monodromy group of $\mu_{12}$, resp. $\mu_{35}:=\goth i \circ \mu_{34}$ by $\Cal M_{12}$, resp. $\Cal M_{35}$. We have $\Cal M_{12}=GL_n(\n F_{q^2})/\n F_q^*$. Really, the automorphism group of $M_0$ is $GL_n(\n F_{q^2}[T])$: an element $g$ of this group acts on the basis $e_*$ of (8.2.2) if $A=0$. But if $A$ is a generic matrix, then for $g\in GL_n(\n F_{q^2}[T]) - GL_n(\n F_{q^2})$ the result of the action of $g$ on (8.2.2) becomes a more complicated equation: terms having higher powers of $\tau$ appear, so these $g$ do not belong to $\Cal M_{12}$. Factorization by $\n F_q^*$ is obvious. Further, obviously $\Cal M_{35}=\{\gamma\in PGL_{2n}(\n Z_\infty)|\gamma(\omega I_n)=\omega I_n\}$. The outer quadrangle of 8.2.4 defines a map from $\Cal M_{12}$ to $\Cal M_{35}$ which we denote by $\alpha$.
\medskip
{\bf Lemma 8.2.11.} $\alpha$ is injective, and $\im \alpha= \{\gamma\in PGL_{2n}(\n F_{q})|\gamma(\omega I_n)=\omega I_n\}=\Cal M_{35} \cap PGL_{2n}(\n F_{q}) \subset PGL_{2n}(\n Z_\infty)$.
\medskip
{\bf Proof.} $\alpha$ is defined by the condition: for any $A\in S_1$, $\gamma \in \Cal M_{12}$ we have
$$\mu_{13}(\gamma(A))=(\alpha(\gamma))(\mu_{13}(A))\eqno{(8.2.12)}$$ The explicit formula for $\alpha$ is the following. For odd $q$ we fix $\omega$ satisfying $\omega^2\in \n F_q^*$ (it is an easy exercise to find analog of the below formula for even $q$). Let $\gamma=U+\omega V$, \ \ $U,V\in GL_n(\n F_q)$ (we consider cleary a representative of $\gamma\in GL_n(\n F_{q^2})/\n F_q^*$ in $GL_n(\n F_{q^2})$). Then $$\alpha(\gamma)=\left(\matrix U&-\omega^2V\\ -V&U \endmatrix\right)^{-1}\eqno{(8.2.13)}$$ ($\alpha$ is an antihomomorphism, because the functor of lattice is contravariant). It is checked immediately that (8.2.12) holds. We see that $\im \alpha= \{\gamma\in PGL_{2n}(\n F_{q})|\gamma(\omega I_n)=\omega I_n\}$. $\square$ 
\medskip
Proposition 8.2.3 follows immediately from these lemmas. $\square$
\medskip
{\bf Remark 8.2.14.} We see that if $r=2n$, $n>1$ then the lattice map $\mu_{25}$ is not a local isomorphism near $M_0$. The origin of this phenomenon is reducibility of $M_0$ which implies that the monodromy group of $\mu_{35}$ is much bigger then the one of $\mu_{12}$. For other values of $r$, $n$ a natural analog of $M_0$ is a t-motive with complete multiplication. Apparently if $(r,n)=1$ then all pure t-motives with complete multiplication are irreducible (example: CM-field is $\n F_{q^r}[\theta]$), and analog of Proposition 8.2.3 for this case shows that $\mu_{25}$ is a local isomorphism near this t-motive.
\medskip
{\bf Remark 8.2.15.} From the first sight, for $n=1$ Proposition 8.2.3 contradicts to a result of Drinfeld about 1 -- 1 correspondence between Drinfeld modules and lattices. Really, there is no contradiction: if $n=1$ then $\im \alpha=\Cal M_{35}$, and --- although $\goth i: S_4 \to S_5$ clearly is not an isomorphism --- its restriction to $U_4$ is an isomorphism $U_4 \to U_5$.
\medskip
{\bf 8.3. Duality of C-lattices.} We have no analog of 2.2 for C-lattices, so we use an analog of 3.2 as a definition of duality. Namely, if $Z$ is a Siegel matrix of a C-lattice $(L, e_*)$ then its dual $(L, e_*)'$ is a C-lattice whose Siegel matrix is $-Z^t$ (or $Z^t$ which is the same, but more convenient for further calculations, because $(-\omega I_n)^t\ne \omega I_n$). Equality 3.8.2 shows that this notion is well-defined (entries of $A$, $B$, $C$, $D$ belong to $\n F_q$).

An analog of Theorem 5 holds for C-lattices:
\medskip
{\bf Theorem 8.3.1.} Let $M\in U_2$. Then $\mu_{24}(M')=\mu_{24}(M)'$.
\medskip
{\bf Proof} is completely analogous to the proof of Theorem 5, so we omit it. Alternatively, we can show that the exact form of (8.2.10) is $$\goth Z=\omega I_n + \sum_{k=1}^\infty \sum_{d_1,...,d_k}l_{d_1,...,d_k}A^{(d_1)}\cdot...\cdot A^{(d_k)}$$ where coefficients $l_{d_1,...,d_k}$ satisfy $$l_{d_1,...,d_k}=l_{d_k,...,d_1}$$ This obviously implies the theorem. $\square$
\medskip

\newpage
{\bf Theorem.} ${M'}^{m}$ is uniformizable, and $\underline{H}({M'}^{m})= {\underline{H}(M)'}^{m}$.
\medskip
We shall need several elementary lemmas.
Let $w$ be a number such that $N^w=0$. We define numbers $k_i=k_i(M)$, $i=2, \dots,w+1$, as follows:
$$k_{i}:=\dim \Ker N^{i-1}/ \Ker N^{i-2}- \dim \Ker N^{i}/ \Ker N^{i-1} $$
$$=\dim \im N^{i-2} / \im N^{i-1} - \dim \im N^{i-1} / \im N^{i}\eqno{(9.3)}$$
Equivalently, let $n=d_1+...+d_\al$, where $d_1\ge d_2\ge...\ge d_\al>0$, be a partition of $n$ corresponding to the Jordan form of $N$, i.e. the Jordan form of $N$ consists of $\al$ 0-Jordan blocks of sizes $d_1, d_2,...,d_\al$. We have $w\ge d_1$, $\al\le r$. We shall call a 0-partition of length $\g r$ of a number $\g n$ a representation of $\g n$ as a sum
$$\g n = \g d_1+\g d_2+...+\g d_{\g r}$$ where $\g d_i\in \n Z$ and $\g d_1\ge\g d_2\ge...\ge\g d_{\g r}\ge0$, i.e. a 0-partition is a partition plus several zeroes at its end. We extend the partition $n=d_1+...+d_\al$ to a 0-partition of length $r$ denoted by $\g p=\g p(M)$:
$n=d_1+...+d_\al+d_{\al+1}+...+d_r$ where $d_{\al+1}=...=d_r=0$. Let $n=c_1+...+c_{d_1}+c_{d_1+1}+...+c_w$ be the 0-partition of length $w$ dual to $\g p$ (the definition of the dual 0-partition of a given length is clear). We have $\al=c_1\ge c_2\ge...\ge c_{d_1}>0$, $c_{d_1+1}=...=c_w=0$. We have $\dim \Ker N^{i}=c_1+...+c_i$, hence $k_i=c_{i-1}-c_i\ge0$ (for $i=w+1$ we let $c_{w+1}=0$), and
$$n=\sum_{i=1}^{w}ik_{i+1}\eqno{(9.4)}$$ We have $\al=c_1=\sum_{i=2}^{w+1}k_{i}$, hence $r\ge \sum_{i=2}^{w+1}k_{i}$. We let $k_1:=r- \sum_{i=2}^{w+1}k_{i}$.
\medskip
According ??, elements $T^il_j$, $i=0,1,\dots$, $j=1,\dots,r$, generate $\Lie(M)$ as a $\p$-vector space. Hence, elements $N^il_j$, $i=0,\dots,w-1$, $j=1,\dots,r$, also generate $\Lie(M)$ as a $\p$-vector space. Using this fact we arrange elements $l_1,\dots,l_r$ in $w+1$ segments as follows. First, elements $N^{w-1}l_j$, $j=1,\dots,r$, generate $N^{w-1}\Lie(M)$ as a $\p$-vector space. Its dimension is $k_{w+1}$, hence (first step of a process) we can choose $k_{w+1}$ elements from $l_1,\dots,l_r$ (we denote them by $l_{w+1,1}, \dots, l_{w+1,k_{w+1}}$ respectively) such that
\medskip
(9.6) $N^{w-1}(l_{w+1,i})$, $i=1, \dots, k_{w+1}$, form a $\p$-basis of $N^{w-1}\Lie(M)$.
\medskip
Further (second step), elements $N^{w-2}l_j$, $N^{w-1}l_j$, $j=1,\dots,r$, generate $N^{w-2}\Lie(M)$ as a $\p$-vector space. Elements $N^{w-2}(l_{w+1,i})$, $N^{w-1}(l_{w+1,i})$, $i=1, \dots, k_{w+1}$, are linearly independent over $\p$. Really, let $$\sum_{\al_2=1}^{k_{w+1}}c_{\al_2}N^{w-2}(l_{w+1,\al_2})+ \sum_{\al_1=1}^{k_{w+1}}c_{\al_1}N^{w-1}(l_{w+1,\al_1})=0\eqno{(9.7)}$$
be a non-trivial dependence relation. Applying $N$ to (9.7) we get $$\sum_{\al_2=1}^{k_{w+1}}c_{\al_2}N^{w-1}(l_{w+1,\al_2})=0$$ that contradicts (9.6). Hence, all $c_{\al_2}$ are 0 that again contradicts (9.6).
\medskip
We have $\dim N^{w-1}\Lie(M)=2k_{w+1}+k_{w}$, this follows immediately from (9.3). Hence, we get:
\medskip
Among $l_1,\dots,l_r$ there exist $k_w$ elements (we denote them by $l_{w,1}, \dots, l_{w,k_{w}}$ respectively) such that their intersection with $l_{w+1,1}, \dots, l_{w+1,k_{w+1}}$ is empty and such that
\medskip
(9.8) $N^{w-2}(l_{w,i})$, $i=1, \dots, k_{w}$, $N^{w-2}(l_{w+1,i})$, $i=1, \dots, k_{w+1}$, $N^{w-1}(l_{w+1,i})$, $i=1, \dots, k_{w+1}$, form a $\p$-basis of $N^{w-2}\Lie(M)$.
\medskip
Third step of the process: elements $N^{w-3}l_j$, $N^{w-2}l_j$, $N^{w-1}l_j$, $j=1,\dots,r$, generate $N^{w-3}\Lie(M)$ as a $\p$-vector space. Elements $N^{w-3}(l_{w,i})$, $N^{w-2}(l_{w,i})$, $i=1, \dots, k_{w}$, and $N^{w-3}(l_{w+1,i})$, $N^{w-2}(l_{w+1,i})$, $N^{w-1}(l_{w+1,i})$, $i=1, \dots, k_{w+1}$, are linearly independent over $\p$ (the proof is as the one above for the second step). Hence, we get:
\medskip
Among $l_1,\dots,l_r$ there exist $k_{w-1}$ elements (we denote them by $l_{w-1,1}, \dots,$ $ l_{w-1,k_{w-1}}$ respectively) such that their intersection with $l_{w,1}, \dots, l_{w,k_{w}}$, $l_{w+1,1}, \dots, l_{w+1,k_{w+1}}$ is empty and such that
\medskip
(9.9) $N^{w-3}(l_{w-1,i})$, $i=1, \dots, k_{w-1}$, $N^{w-3}(l_{w,i})$, $N^{w-2}(l_{w,i})$, $i=1, \dots, k_{w}$, and $N^{w-3}(l_{w+1,i})$, $N^{w-2}(l_{w+1,i})$, $N^{w-1}(l_{w+1,i})$, $i=1, \dots, k_{w+1}$, form a $\p$-basis of $N^{w-3}\Lie(M)$.
\medskip
Continuing this process, we represent $r$ as an ordered partition
$$r=k_1+...+k_{w+1}\eqno{(9.10)}$$
(recall that some $k_*$ can be 0) and we represent the set $\{l_1,\dots,l_{r}\}$ as a union of segments $$\{l_1,\dots,l_{r}\}=\bigcup_{u=1}^{w+1}\{l_{u1}, \dots, l_{uk_u}\}\eqno{(9.11)}$$ (the union is ordered and disjoint) such that $\forall \ u=0,\dots,w-1$ we have:
\medskip
$(9.12)$ A $\p$-basis of $N^u\Lie(M)$ is formed by elements $N^\al(l_{\be\ga})$, where $\al\in[u,\dots, w-1]$, $\be\in[\al+2,\dots, w+1]$, $\ga\in [1,\dots, k_\be]$.
\medskip
This implies that for any
$$u\in [1,w], \ \ z\in [u-1, w-1], \ \  y\in [z+2,w+1], \ \ \eqno{(9.13)}$$
$$v=u-1\eqno{(9.14)}$$
there exist matrices $S_{uvyz}$ of size $k_u\times k_y$ with entries in $\p$ (analogues of the Siegel matrix for the case $w=1$) such that $\forall \ u=1,\dots, w, \forall \ i=1,\dots, k_u$ the following holds:
$$N^{u-1}l_{ui}=-\sum_{z=u-1}^{w-1}\sum_{y=z+2}^{w+1} \sum_{j=1}^{k_y}(S_{uvyz})_{ij}N^zl_{yj}\eqno{(9.15)}$$
(if some $k_*$ are 0 then the corresponding $S_{****}$ do not exist).
\medskip
To simplify formulas, below for any $\al$ we consider $\hat l_\al:=l_{\al*}$ as matrix columns. (9.15) becomes a matrix equality
$$N^{u-1}\hat l_{u}=-\sum_{z=u-1}^{w-1}\sum_{y=z+2}^{w+1}S_{uvyz}N^z\hat l_{y}\eqno{(9.16)}$$

{\bf Remark 9.17.} Since always $v=u-1$, really matrices $S_{uvyz}$ depend on 3 parameters $u,y,z$. Number $v$ indicates the exponent of $N$ in the left hand side of (9.15), by analogy with $z$ which indicates the exponent of $N$ in the right hand side of (9.15). This notation is convenient to define a symmetry between $A_*$ and $P_*$, see below.
\medskip
{\bf 9.18.} Example for $w=3$, $u=1$:
$$\matrix l_{1i}=-(\sum_{j=1}^{k_2}(S_{1020})_{ij}l_{2j} &+& \sum_{j=1}^{k_3}(S_{1030})_{ij}l_{3j} &+& \sum_{j=1}^{k_4}(S_{1040})_{ij}l_{4j} \\ \\
&+& \sum_{j=1}^{k_3}(S_{1031})_{ij}N(l_{3j}) &+& \sum_{j=1}^{k_4}(S_{1041})_{ij} N(l_{4j}) \\ \\ &&&+& \sum_{j=1}^{k_4}(S_{1042})_{ij} N^2(l_{4j}) \ )\endmatrix $$
(terms of a fixed column of this formula correspond to a fixed $y$ and different $z$ of (9.15), and terms of a fixed row of this formula correspond to a fixed $z$ and different $y$ of (9.15) ).
\medskip
Applying powers of $N$ to (9.16), for any $v\in [0,\dots, w-1]$, $u\in [1,\dots, w+1]$ we can represent $N^{v}(\hat l_{u})$ as a linear combination of $N^z(\hat l_{y})$ where for a fixed $v$ the numbers $z, \ y$ satisfy $$z\in [v,\dots,w-1], \ \ y\in [z+2,\dots,w+1], \ \ \eqno{(9.19)}$$ Namely, there exist polynomials in $S_{****}$ denoted by $P_{uvyz}$ such that (matrix notations)

$$N^{v}\hat l_{u}=-\sum_{z=v}^{w-1}\sum_{y=z+2}^{w+1}P_{uvyz}N^z\hat l_{y}\eqno{(9.20)}$$
Clearly for $v=u-1$ we have $P_{uvyz}=S_{uvyz}$.

{\bf 9.21.} The domain $v\ge u-1 \ \ \wedge \ \ \{z, \ y$ satisfy (9.19)\} is called the non-trivial domain of the definition of $P_{****}$.

For $v< u-1$ (trivial domain) we have:

$$ P_{u,v,y,z}=-1, \hbox{ resp. } P_{u,v,y,z}=0  \eqno{(9.22)}$$ for $y$, $z$ satisfying (9.19), $(y,z)=(u,v)$, resp. $(y,z)\ne(u,v)$.

\medskip
{\bf 9.23.} Example for $w=3$:

\medskip

$N^2\hat l_2=(S_{2131}S_{3242}-S_{2141})N^2\hat l_4$, i.e. $P_{2242}=-S_{2131}S_{3242}+S_{2141}$;

\medskip

$N^2\hat l_1=(-S_{1020}S_{2131}S_{3242}+S_{1020}S_{2141} +S_{1030}S_{3242}-S_{1040})N^2\hat l_4$, i.e.

\medskip

$P_{1242}=S_{1020}S_{2131}S_{3242}-S_{1020}S_{2141} -S_{1030}S_{3242}+S_{1040}$;

\medskip

$N\hat l_1=(S_{1020}S_{2131}-S_{1030})N\hat l_3+(S_{1020}S_{2141}-S_{1040})N\hat l_4+$

\medskip

$+(S_{1020}S_{2142}+S_{1031}S_{3242}-S_{1041})N^2\hat l_4$, i.e.

\medskip

$P_{1131}=-S_{1020}S_{2131}+S_{1030}$, $P_{1141}=-S_{1020}S_{2141}+S_{1040}$,

\medskip

$P_{1142}=-S_{1020}S_{2142}-S_{1031}S_{3242}+S_{1041}$.
\medskip
{\bf Remark 9.24.} Although we do not need this fact, let us give a formula for some $P_{****}$. Let us define a block unitriangular matrix $\g S$ whose $(i,j)$-th block is $S_{i,i-1,j,i-1}$ for $j>i$, $I_{k_i}$ for $i=j$ and 0 for $j<i$. Further, we define a block unitriangular matrix $\g P$ whose $(i,j)$-th block is $-P_{i,j-2,j,j-2}$ for $j>i$, $I_{k_i}$ for $i=j$ and 0 for $j<i$. We have $\g P=\g S^{-1}$ (see the the below propositions).
\medskip
Some $P_{****}$ that enter in the below formula for $\bar \g D$ are not of the form of the elements of the inverse unitriangular matrix, for example $P_{1142}$, $w=3$.
\medskip
For the proof of Lemmas 9.32, 23, we need
\medskip
{\bf Lemma 9.25.} For all $i,\ j,\ \psi,\ \xi$ (domain?)
$$(\sum_{\be=0}^{j+\xi-w} \sum_{\al=i+\be}^{w+1-j+\be} S_{i-1,i-2,\al,i-2+\be} P_{\al,w-j+\be,\psi,\xi})-$$ $$-S_{i-1,i-2,\psi,i-2+\xi+j-w}+P_{i-1,w-j,\psi,\xi}=0\eqno{(9.25.1)}$$
(A recurrent formula for $P_{****}$).
\medskip
{\bf Proof.} First, we rewrite (9.16) for $u=i-1$:

$$N^{i-2}\hat l_{i-1}=-\sum_{z=i-2}^{w-1}\sum_{y=z+2}^{w+1}S_{i-1,i-2,y,z}N^z\hat l_{y}\eqno{(9.25.2)}$$

Now, for any $$j=1,\dots,w-i+2\eqno{(9.25.3)}$$ we apply $N^{w-i+2-j}$ to (9.25.2):

$$N^{w-j}\hat l_{i-1}=-\sum_{z=i-2}^{i+j-3}\sum_{y=z+2}^{w+1}S_{i-1,i-2,y,z}N^{z+w-i+2-j}\hat l_{y}\eqno{(9.25.4)}$$

(since $N^w=0$, we get that $z\le i+j-3$ ).

We change a summation variable: $z=i-2+\be$, and $y \to \al$, we get

$$N^{w-j}\hat l_{i-1}=-\sum_{\be=0}^{j-1}\sum_{\al=i+\be}^{w+1}S_{i-1,i-2,\al,i-2+\be}N^{w+\be-j}\hat l_{\al}\eqno{(9.25.5)}$$

Now we use (9.20), we make the following variable change:

$$u \to \al\ \ \ \ \ y\to \psi$$
$$v\to w-j+\be\ \ \ \ \ z\to \xi$$
we get
$$N^{w-j+\be}\hat l_\al=-\sum_{\xi=w-j+\be}^{w-1} \sum_{\psi=\xi+2}^{w+1} P_{\al,w-j+\be,\psi,\xi}N^{\xi}\hat l_\psi\eqno{(9.20a)}$$

We substitute (9.20a) in (9.25.5):

$$N^{w-j}\hat l_{i-1}=\sum_{\be=0}^{j-1}\sum_{\al=i+\be}^{w+1} \sum_{\xi=w-j+\be}^{w-1}\sum_{\psi=\xi+2}^{w+1} S_{i-1,i-2,\al,i-2+\be} P_{\al,w-j+\be,\psi,\xi}N^{\xi}\hat l_\psi\eqno{(9.25.6)}$$

We change the order of summation in (9.25.6):

$$N^{w-j}\hat l_{i-1}=\sum_{\xi=w-j}^{w-1} \sum_{\psi=\xi+2}^{w+1} (\sum_{\be=0}^{j+\xi-w} \sum_{\al=i+\be}^{w+1} S_{i-1,i-2,\al,i-2+\be} P_{\al,w-j+\be,\psi,\xi}) N^{\xi}\hat l_\psi\eqno{(9.25.7)}$$

We rewrite (9.20) making changes:

$$u\to i-1\ \ \ \ \ \ y \to \psi$$
$$v\to w-j\ \ \ \ \ \  z\to \xi$$
we get
$$N^{w-j}\hat l_{i-1}=-\sum_{\xi=w-j}^{w-1}\sum_{\psi=\xi+2}^{w+1} P_{i-1,w-j,\psi,\xi}N^\xi \hat l_{\psi}\eqno{(9.25.8)}$$

For $\psi\ge \xi+2$ elements $N^\xi l_{\psi i}$, $i=1,\dots,k_\psi$, are linearly independent over $\p$. Hence, (9.25.7), (9.25.8) imply

$$P_{i-1,w-j,\psi,\xi}=-\sum_{\be=0}^{j+\xi-w} \sum_{\al=i+\be}^{w+1} S_{i-1,i-2,\al,i-2+\be} P_{\al,w-j+\be,\psi,\xi}\eqno{(9.25.9)}$$

Here the domain of $\xi$, $\psi$ is:

$$\xi\in [w-j,\dots,w-1], \ \ \ \psi\in[\xi+2,\dots, w+1]$$
Taking into consideration (9.22) we can rewrite (9.25.9) as follows:

$$P_{i-1,w-j,\psi,\xi}=-(\sum_{\be=0}^{j+\xi-w} \sum_{\al=i+\be}^{w+1-j+\be} S_{i-1,i-2,\al,i-2+\be} P_{\al,w-j+\be,\psi,\xi})+$$ $$+S_{i-1,i-2,\psi,i-2+\xi+j-w}\eqno{(9.25.10)}$$ with the same domain of $\xi$, $\psi$. This is (9.25.1). Because of (9.25.3), this formula is valid for $w-j\ge i-2$ (the non-trivial case of the definition of $P_{****}$). $\square$
\medskip
Let us consider a symmetry $\g s: \n Z^4\to\n Z^4$ defined as follows: $\g s(\al,\be,\ga,\de)=(w+2-\ga,w-1-\de,w+2-\al,w-1-\be)$.
\medskip
{\bf Remark 9.26.} $\g s$ has the following geometric interpretation. Let us consider a matrix $NL$ whose $(i,j)$-th entry is a symbol $N^{i-1}\hat l_{j}$. We interpret a quadruple $(\al,\be,\ga,\de)$ as a vector from $N^\be \hat l_\al$ to $N^\de \hat l_\ga$ in $NL$. $\g s$ is the reflection of this vector with respect to the center of $NL$ and the inversion of its direction.
\medskip
{\bf Definition 9.27.} $\bar S_{uvyz}:=-P_{\g s(uvyz)}^t$ (defined if $P_{\g s(uvyz)}$ has meaning).
\medskip
We shall consider block $r\times r$-matrices having the following block structure: their block size is $(w+1)\times(w+1)$, quantities of columns in blocks are $k_{w+1},k_w,\dots,k_1$ (counting from the left to the right), and quantities of lines in blocks are $k_{1},k_2,\dots,k_{w+1}$ (counting from up to down). Hence, the $(\al,\be)$-th block of this matrix is a $k_{\al}\times k_{w+2-\be}$-matrix. These matrices will be called skew $k_*$-block matrices.

\medskip
$\forall \ i = 0,\dots, w$ we define skew $k_*$-block matrices $C_i=C_i(S_{****})$ as follows:
\medskip
The $(\al, \be)$-th block of $C_i$ is $S^t_{w+2-\be,w+1-\be,\al,i}$ if the quadruple $(w+2-\be,w+1-\be,\al,i)$ satisfies (9.13, 9.14)\footnotemark \footnotetext{Really, it always satisfies (9.14).} (i.e. if $S_{w+2-\be,w+1-\be,\al,i}$ exists);
\medskip
The $(i+1,w+1-i)$-th block of $C_i$ is $I_{k_{i+1}}$, all other blocks of $C_i$ are 0:
$$(C_i)_{\al\be}=S^t_{w+2-\be,w+1-\be,\al,i}\eqno{(9.28.1)}$$
$$(C_i)_{i+1,w+1-i}=I_{k_{i+1}}\eqno{(9.28.2)}$$
Particularly, the $(\al,\be)$-th block of $C_i$ is a $(k_\al\times k_{w+2-\be})$-th matrix.
\medskip
$\forall \ i = 0,\dots, w$ we define skew $k_*$-block matrices $\bar C_i=\bar C_i(S_{****})$ as follows:
\medskip
The $(\al,\be)$-th block of $\bar C_i$ is given by the formula

$$(\bar C_i)_{\al,\be}=-P_{\al,w-1-i,w+2-\be,w-\be}=\bar S^{t}_{\be,\be-1,w+2-\al,i}\eqno{(9.29)}$$
if the quadruple $(\al,w-1-i,w+2-\be,w-\be)$ belongs to the non-trivial domain of $P_{****}$;
\medskip
For $i=0,\dots,w$
$$(\bar C_i)_{w+1-i,i+1}=I_{k_{w+1-i}}\eqno{(9.30)}$$
other block entries of $\bar C_i$ are 0.
\medskip
{\bf Remark 9.31.} Formula (9.30) is concordant with (9.29), if we consider $P_{****}$ from (9.22). Nevertheless, some 0-blocks of $\bar C_i$ correspond to $P_{**yz}$ where $(y,z)$ do not satisfy (9.19), and hence this $P_{****}$ is not defined.
\medskip
Finally, we define elements $B(S_{****}):=\sum_{i=0}^w C_iN^i\in M_r(\p)[N]$ and $\bar B(S_{****}):=\sum_{i=0}^w \bar C_iN^i\in M_r(\p)[N]$.
\medskip
Example for $w=3$:
$$B(S_{****})=\left(\matrix 0&0&0&I_{k_1}\\0&0&0& S^t_{1020}\\0&0&0& S^t_{1030}\\0&0&0& S^t_{1040} \endmatrix \right)+
\left(\matrix 0&0&0&0\\0&0&I_{k_2}& 0\\0&0&S^t_{2131}& S^t_{1031}\\0&0&S^t_{2141}& S^t_{1041} \endmatrix \right)N+$$ $$+ \left(\matrix 0&0&0&0\\0&0&0& 0\\0&I_{k_3}&0& 0\\0& S^t_{3242}&S^t_{2142}& S^t_{1042} \endmatrix \right)N^2 + \left(\matrix 0&0&0&0\\0&0&0&0\\0&0&0& 0\\I_{k_4}&0&0&0 \endmatrix \right)N^3$$
\medskip
$$\bar B(S_{****})=\left(\matrix \bar S^{t}_{1040}&0&0&0&\\ \bar S^{t}_{1030}&0&0&0\\ \bar S^{t}_{1020}&0&0&0\\I_{k_4}&0&0&0  \endmatrix \right)+
\left(\matrix \bar S^{t}_{1041}&\bar S^{t}_{2141}&0&0\\ \bar S^{t}_{1031}&\bar S^{t}_{2131}&0& 0\\0&I_{k_3}&0&0 \\0&0&0&0  \endmatrix \right)N+$$ $$+ \left(\matrix \bar S^{t}_{1042}&\bar S^{t}_{2142}& \bar S^{t}_{3242}&0\\0&0&I_{k_2}& 0\\0&0&0& 0\\0&0&0& 0  \endmatrix \right)N^2 + \left(\matrix 0&0&0&I_{k_1}\\0&0&0&0\\0&0&0& 0\\0&0&0&0 \endmatrix \right)N^3=$$
\medskip
$$=\left(\matrix -P_{1242}&0&0&0&\\-P_{2242}&0&0&0\\-P_{3242}&0&0&0\\I_{k_4}&0&0&0  \endmatrix \right)+
\left(\matrix -P_{1142}&-P_{1131}&0&0\\-P_{2142}&-P_{2131}&0& 0\\0&I_{k_3}&0&0 \\0&0&0&0  \endmatrix \right)N+$$ $$+ \left(\matrix -P_{1042}&-P_{1031}& -P_{1020}&0\\0&0&I_{k_2}& 0\\0&0&0& 0\\0&0&0& 0  \endmatrix \right)N^2 + \left(\matrix 0&0&0&I_{k_1}\\0&0&0&0\\0&0&0& 0\\0&0&0&0 \endmatrix \right)N^3$$
\medskip
{\bf Lemma 9.32.} $B(S_{****})^t \cdot \bar B(S_{****})=I_r N^w\in M_r(\p)[N]$.
\medskip
{\bf Proof.} We denote $B(S_{****})^t \cdot \bar B(S_{****})$ by $\sum_\mu \Cal C_\mu N^\mu$. The fact that $\Cal C_w=I_{r}$ is obvious: the only non-zero factors that enter in the sum $\sum_{\ga=0}^{w} C_\ga^t \bar C_{w-\ga}$ are products of blocks of $C_*$, $\bar C_*$ containing $I_*$, and they form $I_r$. Also it is obvious that for $\mu>w$ we have $\Cal C_\mu=0$, because all products whose sum is $\Cal C_\mu$, have at least one factor 0. We need to consider $\Cal C_\mu$ for $\mu<w$. (9.28.1), (9.28.2), (9.29), (9.30) give us (here and below $(C^t_\ga)_{\nu\de}$ is the $(\nu\de)$-th block of $C^t_\ga$, i.e. $(C^t_\ga)_{\nu\de}=((C_\ga)_{\de\nu})^t$ )
$$(\Cal C_\mu)_{\nu\pi}=\sum_{\ga=0}^\mu \sum_{\de=1}^{w+1} (C^t_\ga)_{\nu\de}(\bar C_{\mu-\ga})_{\de\pi}\eqno{(9.32.1.1)}$$ $$=-\sum_{\ga,\de} S_{w+2-\nu, w+1-\nu, \de, \ga}P_{\de,w-1-\mu+\ga,w+2-\pi,w-\pi}+ \eqno{(9.32.1.2)}$$
$$-P_{w+2-\nu,2w-\mu-\nu,w+2-\pi,w-\pi}+S_{w+2-\nu,w+1-\nu,w+2-\pi,\mu+1-\pi}\eqno{(9.32.1.3)}$$
where (9.32.1.2) corresponds to the products $(C^t_\ga)_{\nu\de}(\bar C_{\mu-\ga})_{\de\pi}$ where both terms $\ne 0, \ I_*$, and (9.32.1.3) corresponds to the products where one of the terms is $I_*$.

Let us find the relations satisfied by $\mu,\ \nu, \ \pi$ and the domain of summation by $\ga, \ \de$ in (9.32.1.2). We have
$$\matrix (C^t_\ga)_{\nu\de}\ne0, I_{k_*} \ \iff \ \nu\ge w+1-\ga \ \ \wedge \ \ \de\ge\ga+2 \\
(\bar C_{\mu-\ga})_{\de\pi}\ne0, I_{k_*} \ \iff \ \de\le w-(\mu-\ga) \ \ \wedge \ \ \pi \le \mu-\ga+1\endmatrix \eqno{(9.32.2)}$$
The set of $\ga, \ \de$ is non-empty $\iff \ \mu\le w-2$ and $\mu +\nu -\pi\ge w$. In this case the conditions (9.32.2) on $\ga, \ \de$ become
$$\mu+1-\pi\ge\ga\ge w+1-\nu\eqno{(9.32.3.1)}$$
$$w+\ga-\mu\ge\de\ge\ga+2\eqno{(9.32.3.2)}$$
Now we use Proposition 9.25 for
$$\matrix i=w+3-\nu && \mu=j+i-3 \\
j=\mu+\nu-w &\iff & \nu=w-i+3 \\
\psi=w-\pi+2 && \pi=w-\xi=w-\psi+2\\
\xi=w-\pi \endmatrix\eqno{(9.32.4)}$$
and summation variables $\al, \ \be$ in (9.25.1) are
$$\matrix \al=\de \\ \be=\ga-i+2 \endmatrix\eqno{(9.32.5)}$$
Under this variable change, (9.32.1.2) becomes the double sum in (9.25.10), and (9.32.1.3) becomes
$$-P_{i-1,w-j,\psi,\xi}+S_{i-1,i-2,\psi,i-2+\xi+j-w}$$
hence the desired. $\square$
\medskip
Let for $i=0,\dots,w-1$ $X_i$ be skew $k_*$-matrices having the following property:
\medskip
If $(\al, \ \be)$ are such that the $(\al,\ \be)$-block of $\bar C_i$ is 0 or $I_*$ then $(X_i)_{\al\be}=(\bar C_i)_{\al\be}$;

If $(\al, \ \be)$ are such that the $(\al,\ \be)$-block of $\bar C_i$ is $\ne 0, \ I_*$ then $(X_i)_{\al\be}$ is arbitrary.
\medskip
We denote $X:=\sum_{i=0}^{w-1}X_iN^i$.
\medskip
{\bf Lemma 9.33.} If $B(S_{****})^t \cdot X\in N^wM_r(\p[N])$ then $X=\bar B(S_{****})$.
\medskip
{\bf Proof.} For any fixed $\mu, \ \nu, \ \pi$ (9.32.1.2), (9.32.1.3) become
$$\sum_{\ga,\de} S_{w+2-\nu, w+1-\nu, \de, \ga} (X_{\ga-\mu})_{\de,\pi}\eqno{(9.33.1a)}$$ $$
+ (X_{\mu+\nu-w-1})_{w+2-\nu,\pi}+S_{w+2-\nu,w+1-\nu,w+2-\pi,\mu+1-\pi}=0\eqno{(9.33.1b)}$$
This is system of linear equations with unknowns $(X_i)_{\al\be}$ where $$\matrix 0\le i \le w-1 \\ \\ 1\le \be\le i+1\\ \\ 1\le\al \le w-i \endmatrix \eqno{(9.33.2)}$$ (for other values of $i,\ \al, \ \be$ we have $(X_i)_{\al\be}=0$). We arrange $(X_i)_{\al\be}$ in decreasing order of $i+\al$ (for $(X_i)_{\al\be}$ having equal $i+\al$ their ordering is arbitrary), and we arrange equations (9.33.1) in decreasing order of $\mu$ (for equations having equal $\mu$ the order of equations corresponds to the order of $(X_i)_{\al\be}$ having equal $i+\al$). Under this arrangement of unknowns and equations, the matrix of the system (9.33.1) becomes unitriangular. Really, for any $i,\ \al, \ \be$ satisfying (9.33.2) there is exactly one values of $\mu, \ \nu, \ \pi$ --- namely, $$\matrix \mu=i+\al-1 \\ \\ \nu=w+2-\al\\ \\ \pi=\be \endmatrix $$ such that the first term of (9.33.1b) is $(X_i)_{\al\be}$. Other terms of (9.33.1) for these $\mu, \ \nu, \ \pi$ --- namely, the terms that enter in (9.33.1a) --- contain $(X_{\ga-\mu})_{\de,\pi}$ such that $\ga-\mu+\de>i+\al$ hence unitriangularity.
\medskip
Lemma 9.32 affirms that
$$(X_i)_{\de,\pi}=-P_{\de,w-1-i,w+2-\pi,w-\pi}$$
is a solution to this system. Unitriangularity implies that this solution is unique. $\square$
\medskip
Let us consider $\g q_H$ from above (9.1).
\medskip
{\bf Lemma 9.34.} $\forall \ u=1,\dots, w+1$, for $v=u-1$, $\forall \ i=1,\dots, k_u$ the elements

$$\om_{ui}:=N^{v}l_{ui}+\sum_{y=u+1}^{w+1}\sum_{z=u-1}^{y-2}\sum_{j=1}^{k_y}(S_{uvyz})_{ij}N^zl_{yj}\eqno{(9.34.1)}$$

form a basis of $\g q_H$. $\square$
\medskip
We denote the set of elements $\om_{ui}$ ($u$ is fixed, $i$ varies) by $\hat \om_{u}$ (matrix columns). So, (9.34.1) becomes ($v=u-1$)
$$\hat \om_{u}=N^{v}\hat l_{u}+\sum_{z=u-1}^{w-1}\sum_{y=z+2}^{w+1}S_{uvyz}N^z\hat l_{y}\eqno{(9.35)}$$

Let $\vf_i$, $i=1,...,r$ be the basis of $L'$ dual to $l_i$, i.e. $\vf_i(l_j)=\delta_{ij}$. We shall need the dual numbers $k'_i:=k_{w+2-i}$ (inverse order of $k_*$). We consider the analogous two-subscript notation of $\vf_i$, but the order of segments of the partition of $\vf_i$ is opposite, namely:
$$(\vf_{w+1,1},\dots,\vf_{w+1,k'_{w+1}},\ \ \vf_{w,1},\dots,\vf_{w,k'_{w}},\ \ \ \dots \ \ \ \vf_{11},\dots,\vf_{1,k'_{1}}):= (\vf_1,\dots, \vf_r)\eqno{(9.35.1)}$$
(order of elements $\vf_*$ is the same in both sides of this equality).
\medskip
{\bf Lemma 9.36.} $\forall \ u=1,\dots, w+1$, for $v=u-1$, $ \forall \ i=1,\dots, k'_u$ the elements

$$\chi_{ui}:=N^{v}\vf_{ui}+\sum_{y=u+1}^{w+1}\sum_{z=u-1}^{y-2}\sum_{j=1}^{k'_y}(\bar S_{uvyz})_{ij}N^z\vf_{yj}\eqno{(9.36.1)}$$

form a basis of ${\g q'_H}^w$.

{\bf Proof.} As above we denote the set of elements $\vf_{ui}$, resp. $\chi_{ui}$ ($u$ is fixed, $i$ varies) by $\hat \vf_{u}$, resp. $\hat \chi_{u}$ (matrix columns). (9.35), (9.36.1) can be written in terms of blocks of $C_i$, $\bar C_i$:

$$\hat \om_{u}=\sum_{z=0}^{w}\sum_{y=1}^{w+1} (C_z^t)_{w+2-u,y}N^z\hat l_{y}$$
$$\hat \chi_{u}=\sum_{z=0}^{w}\sum_{y=1}^{w+1} (\bar C_z^t)_{uy}N^z\hat \vf_{w+2-y}$$
We must prove that $\forall \ u_1, \ u_2$ we have $\hat \om_{u_1}\hat \chi_{u_2}^t=\de_{u_1}^{u_2}I_{k_{u_1}}N^w$ (product is pairing). This is immediate:

$$\hat \om_{u_1}\hat \chi_{u_2}^t=\sum_{z_1=0}^{w}\sum_{y_1=1}^{w+1}\sum_{z_2=0}^{w}\sum_{y_2=1}^{w+1} (C_{z_1}^t)_{w+2-u_1,y_1}N^{z_1}\hat l_{y_1} N^{z_2}\hat \vf_{w+2-y_2}^t (\bar C_{z_2})_{y_2,u_2}$$
We have $\hat l_{y_1} \hat \vf_{w+2-y_2}^t=\de_{y_1}^{y_2}I_{k_{y_1}}$, hence
$$\hat \om_{u_1}\hat \chi_{u_2}^t=\sum_{z_1=0}^{w}\sum_{y=1}^{w+1}\sum_{z_2=0}^{w} (C_{z_1}^t)_{w+2-u_1,y}N^{z_1+z_2} (\bar C_{z_2})_{y,u_2} = \sum_{z=0}^{2w}(\Cal C_z)_{w+2-u_1,u_2}N^z$$
Lemma 9.32 implies the desired. $\square$
\medskip
{\bf Corollary 9.37.} Matrices $S_{uvyz}(L')$ for the dual lattice $L'$ are $\bar S_{uvyz}(L)$ (order of segments of $\vf_*$, and hence of numbers $k_*$, is inverse).
\medskip
We define $D_i$ ($i=1,\dots, w$ ) as $-<N^{i-1}(l_*),f_*>$ where the ordering of $l_*$ is $l_{11},\dots, l_{w+1, k_{w+1}}$. Hence, $D_i$ is a union of $w+1$ matrices $D_{ij}$, where $D_{ij}$ is a $r\times k_j$-matrix, and $(D_{ij})_{\al\be}=-<N^{i-1}(l_{j\be}),f_\al>$. This means that there are relations between $D_{ij}$ coming from (9.20), namely:

$$D_{v+1,u}= -\sum_{z=v}^{w-1}\sum_{y=z+2}^{w+1}D_{z+1,y} P_{uvyz}^t\eqno{(9.38)}$$

Like in (9.21), for $(z,y)$ satisfying $y\ge z+1$ (resp. $y< z+1$) we shall call the corresponding $D_{zy}$ as belonging to the trivial (resp. non-trivial) domain.
\medskip
{\bf Lemma 9.39.} $\Psi_N B(S_{****})\in M_r(\p[[N]])$.
\medskip
{\bf Proof.} For any $\mu=0,\dots, w-1$ we must prove that $\g C_\mu:=\sum_{\de =0}^\mu D_{w-\de}C_{\mu-\de}$ is 0. The $\nu$-th block ($\nu=1,\dots,w+1$) of this matrix is
$$(\g C_\mu)_\nu:=\sum_{\de =0}^\mu \sum_{\ga=1}^{w+1} D_{w-\de,\ga}(C_{\mu-\de})_{\ga \nu}\eqno{(9.39.1)}$$
We have
$$(C_{\mu-\de})_{\ga \nu}\ne0, I_{k_*} \ \ \iff \de\le \nu+\mu-w-1 \ \ \wedge \ \ \ga\ge \mu-\de+2$$
$$(C_{\mu-\de})_{\ga \nu}=I_{k_*} \ \ \iff \de=\nu+\mu-w-1 \ \ \wedge \ \ \ga=w+2-\nu$$
hence (9.39.1) becomes
$$(\g C_\mu)_\nu=\sum_{\de =0}^{\nu+\mu-w-1}\sum_{\ga=\mu-\de+2}^{w+1}D_{w-\de,\ga}S^t_{w+2-\nu,w+1-\nu,\ga,\mu-\de}+\eqno{(9.39.2.1)}$$
$$+D_{2w+1-\nu-\mu,w+2-\nu}\eqno{(9.39.2.2)}$$
where (9.39.2.1) is non-empty if $\nu+\mu\ge w+1$.

Terms $D_{2w+1-\nu-\mu,w+2-\nu}$ always belong to the non-trivial domain. We separate the terms of (9.39.2.1) in terms of trivial and non-trivial domain:
$$(\g C_\mu)_\nu=\sum_{\de =0}^{\nu+\mu-w-1}\sum_{\ga=\mu-\de+2}^{w-\de} D_{w-\de,\ga}  S^t_{w+2-\nu,w+1-\nu,\ga,\mu-\de}+\eqno{(9.39.3.1)}$$
$$+\sum_{\de =0}^{\nu+\mu-w-1}\sum_{\ga=w-\de+1}^{w+1}  D_{w-\de,\ga}  S^t_{w+2-\nu,w+1-\nu,\ga,\mu-\de}+\eqno{(9.39.3.2)}$$
$$+D_{2w+1-\nu-\mu,w+2-\nu}\eqno{(9.39.3.3)}$$

Now we substitute non-trivial $D_{**}$ by linear combinations of the trivial ones, using (9.38):
$$(\g C_\mu)_\nu=-\sum_{\de =0}^{\nu+\mu-w-1}\sum_{\ga=\mu-\de+2}^{w-\de}   \sum_{z=w-\de-1}^{w-1}\sum_{y=z+2}^{w+1}D_{z+1,y} P_{\ga,w-\de-1,y,z}^t       S^t_{w+2-\nu,w+1-\nu,\ga,\mu-\de}+\eqno{(9.39.4.1)}$$
$$+\sum_{\de =0}^{\nu+\mu-w-1}\sum_{\ga=w-\de+1}^{w+1}D_{w-\de,\ga}S^t_{w+2-\nu,w+1-\nu,\ga,\mu-\de}-\eqno{(9.39.4.2)}$$
$$-\sum_{z=2w-\nu-\mu}^{w-1}\sum_{y=z+2}^{w+1}D_{z+1,y} P_{w+2-\nu,2w-\nu-\mu,y,z}^t \eqno{(9.39.4.3)}$$
Now we change variables in (9.39.4.2):
$$\de=w-z-1$$
$$\ga=y$$
interchange the order of summation and transpose:
$$(\g C_\mu)_\nu=\sum_{z=2w-\nu-\mu}^{w-1}\sum_{y=z+2}^{w+1} \g K(\mu,\nu,z,y) D^t_{z+1,y}$$

where $$\g K(\mu,\nu,z,y)=-( \sum_{\de =w-z-1}^{\nu+\mu-w-1}\sum_{\ga=\mu-\de+2}^{w-\de} S_{w+2-\nu,w+1-\nu,\ga,\mu-\de}    P_{\ga,w-\de-1,y,z}) +$$
$$+S_{w+2-\nu,w+1-\nu,y,\mu-w+z+1}-P_{w+2-\nu,2w-\nu-\mu,y,z} \eqno{(9.39.5)}$$
Change of variables in (9.39.5):
$$\matrix y=\psi & \nu=w+3-i & \ga=\al \\ z=\xi & \mu=i+\al-3 & \de=\al-\be-1 \endmatrix $$
transforms (9.39.5) to (9.25.1), hence all $\g K(\mu,\nu,z,y)$ are 0. $\square$
\medskip
Let $M'=M^{\prime w}$ be the $w$-dual of $M$.
\medskip
{\bf Lemma 9.40.} For $i=1,\dots, w+1$ we have $k_i(M')=k'_{i}(M)$.
\medskip
{\bf Proof.} We use Lemma 10.2. $\mu$ of 10.2 is $w$, $m_i$ of 10.2 are $d_{w+2-i}$. (10.2.5) means that $d_i(M')=w-d_i$. The result follows immediately from the properties of dual 0-partitions. $\square$

We shall another basis $\hat \eta$ of $L(M')$ obtained by a permutation matrix from the basis (9.35.1). Namely, we let $\eta_{ij}:=\vf_{ij}$ ($i=1, \dots, w+1, \ j=1,\dots, k'_i$), but the order of elements $\eta_{ij}$ is the following ($\hat \eta$ is a matrix column):

$$\hat \eta=(\eta_{11},\dots, \eta_{1,k'_1}, \eta_{21},\dots, \eta_{2,k'_2}, \dots, \eta_{w+1,1},\dots, \eta_{w+1,k'_{w+1}})^t$$

We shall consider only $M$ satisfying the following (compare with 9.12)
\medskip
{\bf Condition 9.41.} For any $u$ elements $N^\al(\eta_{\be\ga})$, where $\al\in[u,\dots, w-1]$, $\be\in[\al+2,\dots, w+1]$, $\ga\in [1,\dots, k'_\be]$, are linearly independent over $\p$.
\medskip
According the general principle that almost all $n$-uples $(v_1,\dots,v_n)$ of vectors in $n$-dimensional vector space form a basis of this space, we can guess that almost all $M$ satisfy 9.41. Really, Lemma 9.40 affirms that the dimension of $N^u \Lie(M')$ is exactly the quantity of elements $N^\al(\eta_{\be\ga})$ mentioned in 9.41. Again by Lemma 9.40, we see that Condition 9.41 implies that these elements are a basis of $N^u \Lie(M')$.
\medskip
We can use ideas of [GL] in order to show that a large class of Anderson t-motives $M$ satisfies Condition 9.41. Namely,
\medskip
Finally, we have
\medskip
{\bf Conjecture 9.42.} All $M$ satisfy Condition 9.41.
\medskip
{\bf Theorem 9.43.} Let $M$ satisfy Condition 9.41. Then
\medskip
{\bf Proof.} We denote the basis $\vf_*$ from (9.35.1) by $\hat \vf$ (matrix column). We have $\hat \eta= \g I\cdot \hat \vf$ where $\g I$ is a matrix of the change of bases. It is a skew $k'_*$-block anti-identity matrix, i.e. its antidiagonal block entries are identity matrices: the $(w+2-i,i)$-block is $I_{k_i}$, and other block entries are 0.

We denote by $\Psi'_N$, resp. $\Psi'_{N,\eta}$ the $N$-$\Psi$-series of $M'$ in bases $\hat \vf$, resp. $\hat \eta$ (as a base of $M'$ over $\p[T]$ we use $\hat f'$ in both cases). We have $\Psi'_N=\Psi_N^{t-1}\Xi_N^{-w}$, $\Psi'_{N,\eta}=\Psi'_N\g I^t$. Since $M$ satisfies Condition 9.41, there exists a set of Siegel matrices for $M'$ with respect to the basis $\hat \eta$. We denote it by $U_{****}$. It defines $B(U_{****})$ --- the corresponding $B$ in $M_r(\p)[N]$.
\medskip
We denote $\Psi_NB(A(M)_{****})$, resp. $\Psi'_{N,\eta}B(U_{****})$ by $\g Z(M)$, resp. $\g Z(M^d)$. We have $\g Z(M)\in M_r(\p[[N]])$, $\g Z(M^d)\in M_r(\p[[N]])$ (Lemma 9.39), hence
$$\g Z(M)^t\g Z(M^d)=B(A(M)_{****})^t\Psi_N^t \Psi'_N \g I^t B(U_{****})=$$ $$B(A(M)_{****})^t \g I^t B(U_{****}) \Xi_N^{-w}\in M_r(\p[[N]])$$ and hence
$$  B(A(M)_{****})^t \g I^t B(U_{****}) \g I^t  \in \Xi_N^{w}M_r(\p[[N]])$$
We have  $\g I^t B(U_{****}) \g I^t$ is of the form $X$ of Lemma 9.33. Further, $\Xi_N^{w}\in N^wM_r(\p[[N]])$, hence $\g I^t B(U_{****}) \g I^t=\bar B(A(M)_{****})$ (Lemma 9.33). This means that $U_{uvyz}=\bar A_{uvyz}$. The theorem follows from Lemma 9.36. $\square$
\medskip

\newpage
{\bf Case of linearly dependent columns.}

We shall consider for simplicity the case $w=2$. We use notations: $\Psi:=\Psi_N(M)$, $\Psi^d:=\Psi_N(M^d)$, $\Psi^d=D^d_2N^{-2}+D^d_1N^{-1}+...$ in the basis $\vf_1,...,\vf_r$ dual to $l_1,...,l_r$. We denote $D^d_2=D_{21}^d | D_{22}^d | D_{23}^d $ where $D_{2i}^d$ are $r\times k_i$ matrices.  Let some of the columns of $D_{21}^d$ are linearly dependent. Interchanging the order of columns we can assume that $\exists \ u>0$ such that the first $k_1-u$ columns of $D_{21}^d$ are linearly independent and the last $u$ columns of $D_{21}^d$ are their linear combinations. Further, we assume that the remaining $u$ columns of $D_{21}^d$ forming (together with the first $k_1-u$ columns of $D_{21}^d$) a basis of $NV$ are the first $u$ columns of $D_{22}^d$. This means that we can represent $D_{21}^d$, $D_{22}^d$ as unions:

$D_{21}^d=(D_{211}^d | D_{212}^d)$, $D_{22}^d=(D_{221}^d | D_{222}^d)$, where the quantities of rows in all these matrices are $r$ and the quantities of columns in $D_{211}^d$, resp. $ D_{212}^d$, $D_{221}^d, \ D_{222}^d$ are $k_1-u$, resp. $u, \ u, \ k_2-u$. We have: $\exists$ a $(k_1-u)\times u$-matrix $Y$ such that $D_{212}^d=D_{211}^dY$.

We denote by $\tilde \Psi^d$ the matrix obtained from $\Psi^d$ by permutation of (12)-block and (21)-block. We assume that (1,2,3)-blocks of $\tilde \Psi^d$ are "good", i.e. there exists $\tilde D^d$ such that

$$\tilde \Psi^d\cdot\tilde D^d\in M_r(\p[[N]])$$

We denote by $\hat \Psi^d$ the matrix obtained from $\Psi^d$ by elimination of (12)-block. It is a $5\times 4$-block matrix, and we have

$$\Psi^d=\hat \Psi^d U_3$$ where $U_3=\left(\matrix 1&Y&0&0&0\\ 0&0&1&0&0\\ 0&0&0&1&0\\ 0&0&0&0&1\endmatrix \right)$ (block structure; sizes of blocks: $k_1-u; \ u; \ k_2-u; \ k_3$ for columns; $k_1-u; \ u; \ u; \ k_2-u; \ k_3$ for lines).

$$\tilde \Psi^d=\hat \Psi^d U_2$$ where $U_2=\left(\matrix 1&0&Y&0&0\\ 0&1&0&0&0\\ 0&0&0&1&0\\ 0&0&0&0&1\endmatrix \right)$ (block structure; sizes of blocks: $k_1-u; \ u; \ k_2-u; \ k_3$ for columns; $k_1-u; \ u; \ u; \ k_2-u; \ k_3$ for lines).
\medskip
$U_3$ has a right inverse $U_3^{-1}=\left(\matrix 1&0&0&0\\ 0&0&0&0\\ 0&1&0&0\\ 0&0&1&0\\ 0&0&0&1\endmatrix \right)$ such that we have

$$\hat \Psi^d = \Psi^dU_3^{-1}$$

hence $B^t\Psi^t\tilde \Psi^d\cdot\tilde B^d\in M_r(\p[[N]])$

and $\tilde \Psi^d=\Psi^dU_3^{-1}U_2$

As earlier we have $$\Psi^t\Psi^d=\Xi^{-w}$$

hence $B^t U_3^{-1}U_2 B^d\in N^2 M_r(\p[[N]])$
\medskip
\medskip

Generalizations.
\medskip
1. To prove that $\g q_{M_1\otimes M_2}=\g q_{M_1} \otimes \g q_{M_2}$.
\medskip
2. It is known that if $A_1$, $A_2$ are abelian varieties then $A_1\otimes A_2$ is a mixed motive. Is it possible to get an analog of (1) for it?
\medskip
3. Analog of the main theorem for groups other than $GL_r$. We have a theorem that any full sublattice of $\n Z^r$ is described by its elementary divisors $d_1 | d_2| ... | d_r$. There is an analog of this theorem for any reductive group. We should find the corresponding generalization of the main theorem of the present paper.
\medskip
This is an analog of Theorem 5 for arbitrary $N$, and of Theorem 6 for the operator of duality. Its proof, as well as its generalization to the case of non-pure $M$, can be easily obtained using the same ideas of the proofs of Theorems 5, 6 (this explains the terminology: this is not a conjecture, but a theorem whose proof is not written yet).
\medskip
For $m=1$ Formula 9.1 and Theorem 5 imply immediately
\medskip
{\bf Proposition 9.4.} Result 9.3 is proved for $m=1$. $\square$
\medskip

\input amstex
\documentstyle{amsppt}
\magnification1200
\tolerance=10000
\overfullrule=0pt
\def\n#1{\Bbb #1}
\def\p{\Bbb C_{\infty}}
\def\fr{\hbox{fr}}

\def\im{\hbox{im }}
\def\invlim{\hbox{invlim}}
\def\tr{\hbox{tr }}

\def\Gal{\hbox{Gal }}

\def\Exp{\hbox{Exp}}

\def\Hom{\hbox{Hom}}

\def\End{\hbox{End}}
\def\Prin{\hbox{Prin}}
\def\Ker{\hbox{Ker }}
\def\Lie{\hbox{Lie}}

\def\Div{\hbox{Div}}
\def\Pic{\hbox{Pic}}

\def\ord{\hbox{ord}}

\def\Id{\hbox{Id}}
\def\Cl{\hbox{Cl}}
\def\Supp{\hbox{ Supp }}
\def\Spec{\hbox{ Spec }}

\def\diag{\hbox{ diag }}
\def\Diag{\hbox{ Diag }}

\def\e11{E_{11}}

\def\ga{\goth A}
\def\w{\hbox{\bf A}}
\def\x{\hbox{\bf K}}
\def\ve{\varepsilon}
\def\vf{\varphi}

\def\ve{\varepsilon}

\def\vf{\varphi}

\def\de{\delta}

\def\ga{\gamma}

\def\be{\beta}

\def\al{\alpha}

\def\om{\omega}
\def\g{\goth }

\topmatter
\title
Duality of Anderson $t$-motives
\endtitle
\author
A. Grishkov, D. Logachev\footnotemark \footnotetext{E-mails: shuragri{\@}gmail.com; logachev94{\@}gmail.com (corresponding author)\phantom{*******************}}
\endauthor
\thanks Thanks: The authors are grateful to FAPESP, S\~ao Paulo, Brazil for a financial support (process No. 2017/19777-6). The first author is grateful to SNPq, Brazil, to RFBR, Russia, grant 16-01-00577a (Secs. 1-4), and to Russian Science Foundation, project 16-11-10002 (Secs. 5-8) for a financial support. The second author is grateful to Gilles Lachaud for invitation to IML and to Laurent Lafforgue for invitation to
IHES where this paper was started, and to Vladimir
Drinfeld for invitation to the University of Chicago where
this paper was continued. Discussions with Greg Anderson, Vladimir
Drinfeld, Laurent Fargues, Alain Genestier, David Goss, Richard Pink, Yuichiro
Taguchi, Dinesh Thakur on
the subject of this paper were very important. Particularly, Alain
Genestier informed me about the paper of Taguchi where the notion of
the dual of a Drinfeld module is defined. Further, Richard Pink indicated me an important reference (see Section 6 for details); proof of the main theorem of the present paper grew from it. Finally, Vladimir
Drinfeld indicated me the proof of the Theorem 12.6 and Jorge Morales
gave me a reference on classification of quadratic forms over $\n
F_q[T]$ (Remark 7.8).
\endthanks
\NoRunningHeads
\address
First author: Departamento de Matem\'atica e estatistica
Universidade de S\~ao Paulo. Rua de Mat\~ao 1010, CEP 05508-090, S\~ao Paulo, Brasil, and Omsk State University n.a. F.M.Dostoevskii. Pr. Mira 55-A, Omsk 644077, Russia.
\medskip
Second author: Departamento de Matem\'atica, Universidade Federal do Amazonas, Manaus, Brasil
\endaddress

\abstract Let $M$ be a t-motive. We introduce the notion of duality for $M$. Main results of the paper (we consider uniformizable $M$ over $\n F_q[T]$ of rank $r$, dimension $n$, whose nilpotent operator $N$ is 0):

1. Algebraic duality implies analytic duality (Theorem 5). Explicitly, this means that the lattice of the dual of $M$ is the dual of the lattice of $M$, i.e. the transposed of a Siegel matrix of $M$ is a Siegel matrix of the dual of $M$.

2. Let $n=r-1$. There is a 1 -- 1 correspondence between pure t-motives (all they are uniformizable), and lattices of rank $r$ in $\p^{n}$ having dual (Corollary 8.4).

3. Pure t-motives have duals which are pure t-motives as well (Theorem 10.3).

4. For a self-dual uniformizable $M$ a polarization form on its lattice $L(M)$ is defined. For some $M$ this form is skew symmetric like in the number field case, but for some other $M$ it is symmetric. An example is given.

5. Some explicit results are proved for $M$ having complete
multiplication. The CM-type of the dual of $M$ is the complement of
the CM-type of $M$. Moreover, for $M$ having multiplication by a
division algebra there exists a simple formula for the
CM-type of the dual of $M$ (Section 12).

6. We construct a class of non-pure t-motives (t-motives having the
completely non-pure row echelon form) for which duals
are explicitly calculated (Theorem 11.5). This is the first step
of the problem of description of all t-motives having duals.

7. If $M$ has good ordinary reduction then the kernels of reduction maps on groups of  torsion points for $M$ and its dual are complementary with respect to a natural pairing (proof is given for a particular case, Conjecture 13.4.1).
\endabstract
\endtopmatter
\document
{\bf 0. Introduction.}
\nopagebreak
\medskip
t-motives are the function field analogs of abelian varieties (more exactly, of abelian varieties with multiplication by an imaginary quadratic field, see [L1]). Main references for t-motives are [A], [G]. Nevertheless, function field analogs of some basic results in the theory of abelian varieties are not known yet.

The present paper contains an analog of such result. Namely, we introduce the notion of duality for a t-motive $M$ (this is not the duality in a Tannakian category!), and we prove some properties of this notion, see the abstract. Particularly, if $M$ is uniformizable and has dual then the lattice of the dual of $M$ is the dual of the lattice of $M$ (Theorem 5). An immediate corollary of the above theorem and the result of Drinfeld on 1 -- 1 correspondence between Drinfeld modules and lattices in $\p$ (here $\p$ is the function field analog of $\n C$) is Corollary 8.4: there is a 1 -- 1 correspondence between pure t-motives of dimension $r-1$ and rank $r$, and lattices of rank $r$ in $\p^{r-1}$ having dual (not all such lattices have dual).

Let us give more details on the contents of the paper. For simplicity, most results are proved for t-motives over the ring $\hbox{\bf A}=\n F_q[T]$, and we consider, with few exceptions, only the case $N=0$. The main definition of duality of t-motives (definition 1.8 --- case $\hbox{\bf A}=\n F_q[T]$ and definition 1.13 --- general case) is given in Section 1.\footnotemark \footnotetext{A version of the definition of duality is obtained independently in [Tae], 2.2.} Section 2 contains the definition of the dual lattice. Section 3 contains explicit formulas for the dual lattice. Section 5 contains the statement and the proof of the main theorem 5 --- coincidence of algebraic and analytic duality for the case $\hbox{\bf A}=\n F_q[T]$ (section 4 contains the statement of the corresponding conjecture for the case of general $\hbox{\bf A}$). Section 6 contains the theorem 6 describing the lattice of the tensor product of two t-motives (case $N=0$; the proof for the general case was obtained, but not published, by Anderson). Section 7 contains the notion of self-dual t-motives and polarization form on them. Some examples are given. We discuss in Section 8 the problem of correspondence between uniformizable t-motives and lattices. Section 9 gives the statement of the main result for the case $N\ne 0$ without proof and a reformulation of the theorem 5 in terms of Pink-Hodge structures of constant weight.

Further, Lemma 1.10 gives the explicit matrix form of the definition of duality of t-motives. Since Taguchi in [T] gave a definition of dual to a Drinfeld module, we prove in Proposition 1.12.3 that the definition of the present paper is equivalent to the original definition of Taguchi. Section 1.14 contains a definition of duality for abelian $\tau$-sheaves ([BH], Definition 2.1), but we do not develop this subject. We prove in Section 10 that pure t-motives have duals which are pure
t-motives as well, and some related results (a proof that the dual of an abelian
$\tau$-sheaf is also an abelian $\tau$-sheaf can be obtained using ideas
of Section 10). In Section 11 we consider t-motives having the
completely non-pure row echelon form, and we give an explicit formula
for their duals. In Section 12 we consider t-motives with complete multiplication, and we give for them a very simple version of the proof of the first part of the main theorem.
Section 13 contains some explicit formulas for t-motives of complete multiplication. In 13.1 we describe the dual lattice, in 13.2 we show that the results of Section 12 are compatible with (the first form of) the main theorem of complete multiplication. Section 13.3 contains an explicit proof of the main theorem for t-motives with complete multiplication by two types of simplest fields.
Section 13.4 gives us an application of the notion of duality to the
reduction of t-motives (subject in development, see [L2]).
\medskip
{\bf Notations. }
\medskip
$q$ is a power of a prime $p$;
\medskip
{\it Case of $M$ over $\n F_q[T]$}:
\medskip
$\n Z_\infty:=\n F_q[\theta]$, $\n R_{\infty}:=\n F_q((1/\theta))$, $\p$ is the completion of its algebraic closure
($\n Z_\infty$, $\n R_{\infty}$, $\p$ are the function field analogs of $\n Z$, $\n R$, $\n C$ respectively);

$\w:=\n F_q[T]$, $\hbox{\bf K}:=\n F_q((1/T))$;

$\iota: \hbox{\bf A} \to \p$ ($\iota(T)=\theta$) is the standard map of generic characteristic (with one exception (1.16), we shall not consider the case of finite characteristic);
we extend $\iota$ to $\hbox{\bf K}$, and we have $\n Z_\infty=\iota(\w)\subset \p$, $\n R_\infty=\iota(\hbox{\bf K})\subset \p$.

$\goth C$ (resp. $\goth C_2$) is the Carlitz module over $\w=\n F_q[T]$ (resp. over $\n F_{q^2}[T]$).
\medskip
{\it Case of $M$ over an extension of $\n F_q[T]$}:
\medskip
$\n Q_\infty$ is a finite separable extension of $\n F_q(\theta)$;

$\infty$ is a fixed valuation of $\n Q_\infty$ over the infinity of $\n F_q(\theta)$;

${\n Z_\infty \subset \n Q_\infty}$ is the subring of elements which are regular outside $\infty$;

$\n R_\infty$ is the completion of $\n Q_\infty$ at infinity, and $\p$ --- the completion of its algebraic closure --- is the same as of the case of $M$ over $\n F_q[T]$.

$\w\supset \n F_q[T]$, $\hbox{\bf K}\supset \n F_q((1/T))$ are defined by the condition that $\iota: \w \to \n Z_\infty$, $\iota: \hbox{\bf K} \to \n R_\infty$ are isomorphisms.

$\w_C:=\w\underset{\n F_q}\to{\otimes} \p$ (i.e. $\w_C=\p[T]$ for the case of $M$ over $\n F_q[T]$).

$\goth C$ is a Drinfeld module of rank 1 over {\bf A}.
\medskip
If $P=\frac{\sum a_iT^i}{\sum b_iT^i}\in \p(T)$ then $P^{(k)}:=\frac{\sum a_i^{q^k}T^i}{\sum b_i^{q^k}T^i}$. For $x\in \w_C$, $x=a\otimes z$, $a\in \w$, $z\in \p$ we let $x^{(k)}:= a\otimes z^{q^k}$.

$M_r$ is the set of $r\times r$ matrices. If $C=\{c_{ij}\}$ is a matrix with entries $c_{ij}\in \p(T)$ then $C^{(k)}:=\{c_{ij}^{(k)}\}$, $C^t$ is the  transposed of $C$, $C^{(k)\ -1}=(C^{(k)})^{-1}$, $C^{t-1}=(C^t)^{-1}$.

If $M$ is an $\w_C$-module, we define $M^{(1)}$ as the tensor
product $M\otimes_{\w_C,*^{(k)}}\w_C$ with respect to the map $*^{(k)}: \w_C\to \w_C$ (this notation is concordant in the obvious sense with the above notation $C^{(1)}$).

For a t-motive $M$ we denote by $E=E(M)$ the corresponding t-module (see [G], Theorem 5.4.11; Goss uses the inverse functor $E\mapsto M=M(E)$).

$\Lie(M)$ is $\Lie(E(M))$ ([G], 5.4).

$E_k$ is the unit matrix of size $k$.

Throughout the whole paper the word "canonical" will mean "canonical up to multiplication by elements of $\n F_q^*$".
\medskip
{\bf 1. Definitions. }

\nopagebreak
\medskip
If otherwise is not explicitly stated, throughout the whole paper we consider the case of t-motives $M$ over the ring $\hbox{\bf A}=\n F_q[T]$ such that $N=N(M)=0$. Exceptions: case of arbitrary $\w$ is treated in Sections 1.13, 1.14, 2, 4, 5.2. Case of arbitrary $N$ is treated in Sections 1, 10 and in statements of some results of Anderson in Sections 5, 6.
\medskip
In the present section we consider $M$ such that $N(M)$ is arbitary.
\medskip
Let $\p[T,\tau]$ be the Anderson ring, i.e. the ring of non-commutative polynomials satisfying the following relations (here
$a \in \p$):
$$Ta=aT, \ T\tau = \tau T, \ \tau a = a^q \tau \eqno{(1.1)}$$
We need also an extension of $\p[T,\tau]$ --- the ring $\p(T)[\tau]$
which is the ring of
non-commutative polynomials in $\tau$ over the field of rational
functions $\p(T)$ with
the same relations (1.1). For a left $\p[T,\tau]$-module $M$ we denote by $M_{\p[T]}$ the
same
$M$ treated as a
$\p[T]$-module with respect to the natural inclusion
$\p[T]\hookrightarrow\p[T,\tau]$.
Analogously, we define $M_{\p[\tau]}$; we shall use similar notations
also for the left
$\p(T)[\tau]$-modules.
\medskip
Obviously we have:
\medskip
{\bf (1.2)} For $C\in M_{r}(\p(T))$ operations $C^t$, $C^{-1}$
and $C^{(i)}$ commute.
\medskip
{\bf Definition 1.3.} ([G], 5.4.2, 5.4.12, 5.4.10). A t-motive $M$ is a left
$\p[T, \tau]$-module which is free and finitely generated as both $\p[T]$-,
$\p[\tau]$-module and such that
$$ \exists m=m(M) \ \hbox{such that}\ (T-\theta)^m M/\tau M=0\eqno{(1.3.1)}$$
\medskip
{\bf Remark.} The above object is called "abelian t-motive" (resp. "t-motive") in [G] (resp. [A]), while the name "t-motive" is used in [G] for a more general object ([G], Definition 5.4.2). Since we shall not use objects defined in [G], 5.4.2, I prefer to use a shorter name for the above $M$.
\medskip
t-motives are main objects of the present paper. If we affirm
that an object exists this means that it exists as a t-motive if otherwise
is not stated. We denote dimension of $M$ over $\p[\tau]$ (resp. $\p[T]$) by $n$ (resp.
$r$), these numbers are called dimension and rank of $M$. Morphisms of abelian
t-motives are morphisms of left $\p[T, \tau]$-modules.

To define a left $\p[T, \tau]$-module $M$ is the same as to define a left $\p[T]$-module
$M_{\p[T]}$
endowed by an action of $\tau$ satisfying $\tau(Pm)=P^{(1)}\tau(m)$, $P\in \p[T]$. In
this situation
we can also treate $\tau$ as a $\p[T]$-linear map $M^{(1)}\to M$. This interpretation is
necessary if
we consider the general case $\w\supset \n F_q[T]$.

We need two categories which are larger than the category of abelian
t-motives.
\medskip
{\bf Definition 1.4.} A pr\'e-t-motive is a left $\p[T,
\tau]$-module which is free
and finitely generated as $\p[T]$-module, and satisfies (1.3.1).
\medskip
{\bf Definition 1.5.} A rational pr\'e-t-motive is a left
$\p(T)[\tau]$-module which is
free and finitely generated as $\p(T)$-module.
\medskip
{\bf Remark 1.6.} An analog of (1.3.1) does not exist for them.
\medskip
There is an obvious functor from the category of t-motives to
the category of
pr\'e-t-motives which is fully faithful, and an obvious functor from the category of
pr\'e-t-motives to the
category of rational pr\'e-t-motives. We denote these functors by
$i_1$, $i_2$
respectively. It is easy to see (Remark 10.2.3) that if $M$ is a
pr\'e-t-motive then the
action of $\tau$ on $i_2(M)$ is invertible.

Let $M_1$, $M_2$ be rational pr\'e-t-motives such
that the action of
$\tau$ on $(M_1)_{\p(T)}$ is invertible.
\medskip
{\bf Definition 1.7.} $\Hom(M_1,M_2)$ is a rational pr\'e-t-motive such that

$$\Hom(M_1,M_2)_{\p(T)}=\Hom_{\p(T)}((M_1)_{\p(T)}, (M_2)_{\p(T)})$$
and the action of $\tau$ is defined by the usual manner: for
$\varphi:M_1 \to M_2$, $m\in M_1$
$$(\tau\varphi)(m)=\tau(\varphi(\tau^{-1}(m)))$$
\medskip
{\bf Definition 1.8.} Let $M$ be a t-motive and
$\mu$ a positive
number. A t-motive
$M'={M'}^{\mu}$ is called the $\mu$-dual of $M$ (dual if $\mu=1$) if
$M'=\Hom(M,\goth
C^{\otimes
\mu})$ as a rational pr\'e-t-motive, i.e. $$i_2\circ
i_1(M')=\Hom(i_2\circ i_1(M),\goth
C^{\otimes
\mu})\eqno{(1.8.1)}$$

{\bf Remark.} This definition generalizes the original one of Taguchi
([T], Section 5), see 1.12 below. A similar definition is in [F].
\medskip
{\bf 1.9.} We shall need the explicit matrix description of the above
objects. Let
$e_*=(e_1, ..., e_n)^t$ be the vector column of elements of a basis
of $M$ over $\p[\tau]$. There exists a matrix $\goth A\in M_n(\p[\tau])$ such that

$$T e_* = \goth A e_*, \ \ \goth A = \sum_{i=0}^l \goth A_i \tau^i \hbox{ where } \goth A_i
\in M_n(\p)\eqno{(1.9.1)}$$
Condition (1.3.1) is equivalent to the condition
$$\goth A_0=\theta E_n + N\eqno{(1.9.2)}$$
where $N$ is a
nilpotent matrix, and the
condition
$m(M)=1$ is equivalent to the condition $N=0$.
\medskip
Let $f_*=(f_1, ..., f_r)^t$ be the vector column of elements of a basis
of $M$ over
$\p[T]$. There exists a matrix $Q=Q(f_*)\in M_r(\p[T])$ such that
$$\tau f_* = Q f_*\eqno{(1.9.3)}$$

{\bf Lemma 1.10.} Let $M$ be as above. A t-motive
$M'$ is the $\mu$-dual of $M$ iff there exists a basis
$f'_*=(f'_1, ..., f'_r)^t$ of $M'$ over $\p[T]$ such that its matrix
$Q'=Q(f'_*)$ satisfies
$$Q'=(T-\theta)^{\mu}Q^{t-1} \eqno{(1.10.1)}$$
$\square$
\medskip
{\bf 1.10.2.} For further applications we shall need the following
lemma. The above
$f_*$, $f'_*$ are the dual bases (i.e. if we consider $f'_i$ as
elements of $\Hom(M,\goth
C)$ then $f'_i(f_j)=\delta^i_j\goth f$, where $\goth f$ is canonically defined by the condition that it generates $\goth C_{\p[T]}$ and satisfies $\tau \goth
f=(T-\theta)\goth f$). Let $\gamma$ be an endomorphism of $M$ and $D$ its
matrix in the basis
$f_*$ (i.e. $\gamma(f_*)=Df_*$). Let $\gamma'$ be the dual endomorphism.
\medskip
{\bf Lemma 1.10.3.} The matrix of $\gamma'$ in the basis $f'_*$ is $D^t$.
$\square$
\medskip
{\bf Remark 1.11.1.} For any $M$ having dual there exists a canonical homomorphism $\delta: \goth
C \to M\otimes M'$.
This is a well-known theorem of linear algebra. Really, in the above notations we have
$\goth f \mapsto \sum_i f_i \times f'_i$. It is obvious that $\delta$ is well-defined,
canonical and compatible
with the action of $\tau$.
\medskip
{\bf Remark 1.11.2.} The $\mu$-dual of $M$ --- if it exists --- is
unique, i.e. does not
depend on base change. This follows immediately from Definition 1.8,
but can be deduced
easily from 1.10.1. Really, let $g_*=(g_1, ..., g_r)^t$ be another basis
of $M$ over
$\p[T]$ and $C\in GL_r(\p[T])$ the matrix of base change (i.e. $g_*= C
f_*$). Then
$Q(g_*)=C^{(1)}QC^{-1}$. Let $g'_*=(g'_1, ..., g'_r)^t$ be a basis of
$M'$ over $\p[T]$
satisfying $g'_*= C^{t-1} f'_*$. Elementary calculation shows that
matrices $Q(g_*)$,
$Q(g'_*)$ satisfy (1.10.1).
\medskip
{\bf Remark 1.11.3.} The operation $M \mapsto {M'}^{\mu}$ is
obviously contravariant functorial. I leave as an exercise to the reader to give an exact definition of the corresponding category such that the functor of duality is defined on it, and is involutive (recall that not all t-motives have duals, and the dual of a map of t-motives is a priori a map of rational pr\'e-t-motives).
\medskip
{\bf 1.12.} The original definition of duality ([T], Definition 4.1; Theorem 5.1) from
the first
sight seems to be more
restrictive than the definition 1.8 of the present paper, but really
they are
equivalent. We recall some notations and definitions of [T] in a slightly less
general setting (rough statements; see [T] for the exact statements). Let $G$ be a finite affine group scheme over $\p$, i.e. $G=\Spec R$
where $R$ is a finite-dimensional $\p$-algebra. Let $\mu:R \to R\otimes R$ be the
comultiplication of $R$. Such group $G$ is called a finite $v$-module ([T], Definition
3.1) if there is a homomorphism $\psi: \w \to \End_{gr. \ sch.}(G)$ satisfying some natural conditions (for
example, an analog of 1.3.1). Further, let $\Cal E_G$ be a $\p$-subspace of $R$ defined as
follows: $$\Cal E_G=\{x\in R \ \vert \ \mu(x)=x\otimes 1+1\otimes x\}$$
The map $x\mapsto x^q$ is a $\p$-linear map $\fr: \Cal E_G^{(1)}\to \Cal E_G$. Further,
the map $\psi(T): G \to G$ can be defined on $\Cal E_G$. Let $v: \Cal E_G\to\Cal
E_G^{(1)}$ be a map satisfying $\fr \circ v= \psi(T)-\theta$.

We consider two finite $v$-modules $G$, $H$, the above objects $\fr$, $v$ etc. will carry
the respective subscript. Let * be the dual in the meaning of linear algebra.
\medskip
{\bf Definition 1.12.1} ([T], 4.1). Two finite $v$-modules $G$, $H$
are called dual if there exists an
isomorphism $\alpha: \Cal E^*_H\to\Cal E_G$ such that if we denote by
$\goth v: \Cal E_G\to \Cal E^{(1)}_G$ a map which enters in the commutative diagram
$$ \matrix \Cal E^*_H & \overset{\fr^*_H}\to{\longrightarrow} & \Cal E^{*(1)}_H \\ & & \\
\alpha \downarrow & & \alpha ^{(1)}\downarrow \\ & & \\ \Cal E_G & \overset{\goth
v}\to{\longrightarrow} & \Cal E^{(1)}_G \endmatrix $$
then we have:
$$\fr_G\circ \goth v= \psi_G(T)-\theta \eqno{(1.12.2)}$$
i.e. $\goth v=v_G$.
\medskip
Let $M$ be a t-motive having $m(M)=1$, $E=E(M)$ the corresponding t-module and
$a\in \hbox{\bf
A}$. We denote $E_a$ --- the set of $a$-torsion elements of $E$ --- by $M_a$. It is a
finite
$v$-module.
\medskip
{\bf Proposition 1.12.3.} Let $M$, $M'$ be t-motives which are
dual in the meaning
of Definition 1.8. Then $\forall a\in \hbox{\bf A}$, $a\ne 0$ we have: $M_a$,
$M'_a$ are dual in the meaning of 1.12.1 = [T], Definition 4.1.
\medskip
{\bf Proof.} Condition $a\in \n F_q[T]$ implies that multiplication by $\tau$ is
well-defined on $M/aM$.
\medskip
{\bf Lemma 1.12.3.1.} We have canonical isomorphisms $i: M/aM \to
\Cal E_{M_a}$,
$i^{(1)}: M/aM \to \Cal E^{(1)}_{M_a}$ such that the following
diagrams are commutative:
$$ \matrix M/aM & \overset{\tau}\to{\longrightarrow} & M/aM &&& M/aM &
\overset{T}\to{\longrightarrow} & M/aM
\\ & & &&&&&
\\ i^{(1)}\downarrow & & i\downarrow &&& i\downarrow & & i\downarrow \\ & & &&&&&
\\ \Cal E^{(1)}_{M_a} & \overset{fr}\to{\longrightarrow} & \Cal E_{M_a} &&&
\Cal E_{M_a} & \overset{\psi_T}\to{\longrightarrow} & \Cal E_{M_a}
\endmatrix $$
\medskip
{\bf Proof.} Let $R$ be a ring such that $\Spec R=M_a$. The pairing between $M$ and $E$
shows that there exists a map $M\to R$ which is obviously factorized via an inclusion
$M/aM\to R$. It is easy to see that the image of this inclusion is contained in $\Cal
E_{M_a}$, i.e. we get $i$. Since $\dim_{\p}(M/aM)=\deg a \cdot r(M)$ and
$\dim_{\p}(R)=q^{\deg a \cdot r(M)}$ we get from [T], Definition 1.3 that $i$ is an
isomorphism. Other statements of the lemma are obvious. $\square$
\medskip
This lemma means that we can rewrite Definition 1.12.1 for the case $G=M_a$, $H=N_a$ by
the following way:\footnotemark \footnotetext{Here and below a t-motive $N$ should not be confused with $N$ of 1.9.2.}
\medskip
{\bf 1.12.3.2.} Two finite $v$-modules $M_a$, $N_a$ are dual if there
exists an
isomorphism $\alpha: (N/aN)^* \to M/aM $ such that after identification
via $\alpha$ of
$\tau^*: (N/aN)^* \to (N/aN)^*$ with a map $\goth v: M/aM \to M/aM$ we
have on $M/aM$:
$$\tau\circ \goth v= t-\theta \eqno{(1.12.3.3)}$$
We need a
\medskip
{\bf Lemma 1.12.3.4.} For $i=1,2$ let $N_i$ be a free $\p[T]$-module of
dimension $r$ with a
base $f_{i*} =(f_{i1}, ... , f_{ir} )$, let $\varphi_i:N_i \to N_i$ be
$\p[T]$-linear
maps having matrices $\goth Q_i$ in $f_{i*} $ such that $\goth
Q_2=\goth Q_1^t$, and let
$a$ be as above. Let, further, $\varphi_{i,a}: N_i/aN_i \to N_i/aN_i $
be the natural
quotient of $\varphi_i$. Then there exist $\p$-bases $\tilde f_{i*} $
of $N_i/aN_i $ such
that the matrix of $\varphi_{1,a}$ in the base $\tilde f_{1*} $ is
transposed to the
matrix of $\varphi_{2,a}$ in the base $\tilde f_{2*} $.
\medskip
{\bf Proof.} We can identify elements of $N_2$ with $\p[T]$-linear
forms on $N_1$
(notation: for $x\in N_2$ the corresponding form is denoted by
$\chi_x$) such that
$\chi_{\varphi_2(x)}=\chi_x\circ\varphi_1$. Any $\p[T]$-linear form
$\chi$ on $N_i$
defines a $\p[T]/a\p[T]$-linear form on $N_i/aN_i$ which is denoted by
$\chi_a$. Let now
$x\in N_2/aN_2$, $\bar x$ its lift on $N_2$, then
$\chi_{x,a}=(\chi_{\bar x})_a$ is a
well-defined $\p[T]/a\p[T]$-linear form on $N_1/aN_1$. For $x\in
N_2/aN_2$ we have
$$\chi_{\vf_{2,a}(x),a}= \chi_{x,a}\circ \vf_{1,a}$$
Further, let $\lambda: \p[T] \to \p$ be a $\p$-linear map such that
\medskip
{\bf 1.12.3.5.} Its kernel does not contain any non-zero ideal of
$\p[T]/a\p[T]$.
\medskip
(such $\lambda$ obviously exist.) For $x\in N_2/aN_2$ we denote $\lambda \circ
\chi_{x,a}$ by $\psi_x$, it is a $\p$-linear form on $\p$-vector space
$N_1/aN_1$.
Obviously condition (1.12.3.5) implies that the map $x\mapsto \psi_x$
is an isomorphism
from $N_2/aN_2$ to the space of $\p$-linear forms on $\p$-vector space
$N_1/aN_1$, and we
have
$$\psi_{\vf_{2,a}(x)}=\psi_x\circ\vf_{1,a}$$
which is equivalent to the statement of the lemma. $\square$
\medskip
Finally, the proposition follows immediately from this lemma
multiplied by $T-\theta$,
formula 1.10.1 and 1.12.3.2. $\square$
\medskip
{\bf Remark.} Let $a=\sum_{i=0}^k g_iT^i$, $g_i\in \n F_q$, $g_k=1$.
Taguchi ([T], proof
of 5.1 (iv)) uses the following $\lambda$: $\lambda(T^j)=0$ for $j<k-1$,
$\lambda(T^{k-1})=1$. It is easy to check that for $x=(T^i+
T^{i-1}g_{k-1}+T^{i-2}g_{k-2}+ ...
+g_{k-i})f_{2j}$ for this $\lambda$ we have: $\psi_x(T^i f_{1j})=1$,
$\psi_x(T^{i'}
f_{1j'})=0$ for other $i'$, $j'$.
\medskip
{\bf 1.13.} We consider in Sections 1.13, 1.14 the case of arbitrary $\w\supset \n F_q[T]$.
\medskip
A t-motive over {\bf A} is defined for example in [BH], p.1. Let us reproduce this
definition for the case of characteristic 0. Let $J$ be an ideal of
$\w_C$ generated by the elements $a\otimes 1 - 1 \otimes \iota(a)$ for all $a\in \w$. The ring $\w_C[\tau]$ is defined by the formula $\tau\cdot (a\otimes z)=(a\otimes z^q)\cdot \tau$, $a\in \w$, $z\in \p$.
\medskip
{\bf Definition 1.13.1.} A t-motive $M$ over {\bf A} is a pair $(M, \tau)$
where $M$ is
a locally free $\hbox{\bf A}_C$-module and $\tau$ is an $\hbox{\bf A}_C$-linear map
$M^{(1)}\to M$ satisfying the following analog of 1.3.1, 1.9.2:
$$\exists m \hbox{ such that } J^m(M/\tau (M^{(1)}))=0\eqno{(1.13.2)}$$

{\bf Remark 1.13.3.} We can consider $M$ as an $\w_C[\tau]$-module using the following formula for the product $\tau \cdot m$:
$$\tau \cdot m = \tau(m\otimes 1)$$
where $m\in M$, $m\otimes 1\in M^{(1)}$.
\medskip
The rank of $M$ as a locally free $\hbox{\bf A}_C$-module is called the rank of the
corresponding t-motive $(M, \tau)$. If $\w=\n F_q[T]$ then $M^{(1)}$ is isomorphic to
$M$, we can
consider $M$ as a $\p[T,\tau]$-module, and it is possible to show that in this case
1.13.2 implies
that $M_{\p[\tau]}$ is a free $\p[\tau]$-module. In the general case, the dimension $n$
of $(M, \tau)$
is defined as $\dim_{\p}(M/\tau (M^{(1)}))$.

Let us fix $\goth C=(\goth C, \tau_{\goth C})$ --- a t-motive of rank 1 over
{\bf A}.
For a t-motive $M=(M,\tau_M)$ a t-motive ${M_\goth C'}$ --- the
$\goth C$-dual of $M$
--- is defined as follows. We put
$M_\goth C'=\Hom_{\w_C}(M,\goth C)$.
Since for any locally free $\hbox{\bf A}_C$-modules $M_1$, $M_2$ we have
$$\Hom_{\w_C}(M_1,M_2)^{(1)}=\Hom_{\w_C}(M_1^{(1)},M_2^{(1)})$$
we can define $\tau(M_\goth C')$ by the following formula:

$$\hbox{ For } \vf\in \Hom_{\w_C}(M,\goth C)^{(1)} \hbox{ we have }
\tau(M_\goth C')(\vf)=\tau_{\goth C}\circ \vf \circ \tau_M^{-1}$$
\medskip
{\bf 1.14. Duality for abelian $\tau$-sheaves.} We use notations of [BH],
Definition 2.1 if they do not differ from the notations of the present
paper; otherwise we continue to use notations of the present paper
(for example, $d$ (resp. $\sigma^*(\goth X)$ for any object $\goth X$)
of [BH] is
$n$ (resp. $\goth X^{(1)}$) of the present paper). For any abelian
$\tau$-sheaf $\underline{\Cal F}$ we denote its $\Pi_i$, $\tau_i$ by
$\Pi_i(\underline{\Cal F})$, $\tau_i(\underline{\Cal F})$
respectively. If $M$, $N$ are invertible sheaves on $X$ and $\rho: M
\to N$ a rational map then we denote by $\rho^{inv}: N \to M$ the
rational map which is inverse to $\rho$ with respect to the composition.
We define
$\tau_{\goth r,i-1}(\underline{\Cal F})$ (the rational $\tau_i$) as
the composition map $\tau_{i-1}(\underline{\Cal F})
\circ {\Pi_{i-1}^{(1)}}^{inv}(\underline{\Cal F})$, it is a rational
map from $\Cal F_i^{(1)}$ to $\Cal F_i$.

Let $\underline{\Cal O}$ be a fixed abelian
$\tau$-sheaf having $r=n=1$. The $\underline{\Cal O}$-dual abelian
$\tau$-sheaf $\underline{\Cal F}'= \underline{\Cal
F}'_{\underline{\Cal O}}$ is defined by the formulas
$$\Cal F'_0=\Hom_{X}(\Cal F_0, \Cal O_0)$$ where Hom is the sheaf's
one, and the map $\tau_{\goth r,-1}(\underline{\Cal F}'): {\Cal
F'_0}^{(1)}\to \Cal F'_0$ is defined as follows. We have ${\Cal
F'_0}^{(1)}= \Hom_{X}({\Cal F_0}^{(1)}, {\Cal O_0}^{(1)})$. Let
$\gamma\in\Hom_{X}({\Cal F_0}^{(1)}, {\Cal O_0}^{(1)})(U)$ where $U$
is a sufficiently small affine subset of $X_{\p}$, such that $\gamma:
{\Cal F_0}^{(1)}(U) \rightarrow {\Cal O_0}^{(1)}(U)$.
\medskip
{\bf 1.14.1.} We define: $[[\tau_{\goth r,-1}(\underline{\Cal
F}')](U)](\gamma)$ is the following composition map:
$$\Cal F_0(U)\overset{[\tau_{\goth r,-1}^{inv}(\underline{\Cal
F})](U)}\to{\longrightarrow}
\Cal F_0^{(1)}(U)\overset{\gamma}\to{\to}\Cal O_0^{(1)}(U)
\overset{[\tau_{\goth r,-1}(\underline{\Cal
O})](U)}\to{\longrightarrow}\Cal O_0(U)\in \Hom_{X}({\Cal F_0}, {\Cal
O_0})(U)$$

Clearly that this definition and the definitions 1.8, 1.13 are
compatible with the forgetting functor $\underline{M}(\underline{\Cal
F})$ from abelian $\tau$-sheaves to pure Anderson t-motives of [BH],
Section 3, page 8.
\medskip
{\bf 1.15. Duality over fields.} Let $L\supset \n F_q(\theta)$ be a field extension of $\n F_q(\theta)$, and $M$ a t-motive over $L$ (i.e. a pair ($M$, an $L$-structure on $M$)). Obviously we have
\medskip
{\bf Proposition 1.15.1.} The notion of duality for $M$ over $L$ is well-defined. $\square$
\medskip
Similarly, we have a proposition for Galois action:
\medskip
{\bf Proposition 1.15.2.} Let $M$ be defined over $\overline{ \n F_q(\theta)}$ and $\gamma\in \Gal( \n F_q(\theta))$. Then $(\gamma(M))'=\gamma(M')$. $\square$
\medskip
\medskip
{\bf 1.16. Case of finite characteristic.} Let $\iota: \w\to \bar \n F_q$ be a map of finite characteristic, we denote $\Ker \iota$ by $\Cal P$. The definition of t-motive for this case is similar to 1.3, see [G] for the details. The definition of duality also is similar to the one of the case of generic characteristic. Duality commutes with reduction. Namely, let $M$ be from 1.15, $\goth P$ a prime of $L$ not over the infinity of $\n F_q(\theta)$, $\Cal P\subset \w$ is $\iota^{-1}(\goth P\cap \n F_q[\theta]$) --- the finite characteristic. We consider the case of good reduction of $M$ at $\goth P$, we denote it by $\tilde M$. It is a t-motive in characteristic $\Cal P$.  Let $M$ have dual $M'$.
\medskip
{\bf Proposition 1.16.1.} $\tilde M$ has dual iff $M'$ has good reduction at $\goth P$; in this case they coincide. $\square$
\medskip
{\bf Remark 1.16.2.} Apparently if $M$ has good reduction and dual, then $M'$ also has good reduction (in this case 1.16.1 means that $M'$ exists implies $(\tilde M)'$ exists). For standard-3 t-motives (this is a simple tipe of t-motives, see 11.8.1) apparently this can be shown by explicit calculations.
\medskip
{\bf Remark 1.16.3.} Clearly 1.16.1 is true for the case of bad reductions. I do not give exact definitions for this case.
\medskip
{\bf 1.16.4. Ordinarity.} Let $M$ be of finite characteristic. By analogy with the number field case, $M$ is called ordinary if its Newton polygon consists of 2 segments. If $N=0$ then the Newton polygon of $M'$ is the dual of the one of $M$ (the notion of duality of polygons is clear; apparently the condition $N=0$ can be omitted). So, we have
\medskip
{\bf Proposition 1.16.5.} $M$ is ordinary $\iff M'$ is ordinary. $\square$
\medskip
See 13.4.1 for a more exact result.
\medskip
{\bf 2. Analytic duality.}
\medskip
We consider in the present section the case of arbitrary $\w\supset \n F_q[T]$ (and $N=0$ as usually).
\medskip
Condition $N=0$ implies that an element $a\in \w$ acts on $\Lie(M)$ by multiplication by $\iota(a)$. Hence, we have a
\medskip
{\bf Definition 2.1.} Let $V$ be the space $\p^n$. A locally free $r$-dimensional
$\n Z_\infty$-submodule
$L$ of $V$ is called a lattice if

(a) $L$ generates $V$ as a $\p$-module and

(b) The $\n R_\infty$-linear envelope of $L$ has dimension $r$ over
$\n R_\infty$.
\medskip
Numbers $n$, $r$ are called the dimension and the rank of $L$
respectively. Attached to $(L,V)$ is the tautological inclusion $\vf=\vf(L,V):
L \to V$. We shall consider the category of triples $(\vf, L, V)$; a map $\psi:
(\vf, L,V)\to(\vf_1, L_1,V_1)$ is a pair $(\psi_L, \psi_V)$ where $\psi_L: L \to L_1$ is a $\n Z_\infty$-linear map, $\psi_V: V\to V_1$ is a $\p$-linear map such that $\vf_1\circ \psi_L=\psi_V \circ \vf$.

Inclusion $\vf$ can be extended to a map
$L\underset{\n Z_\infty}\to{\otimes}\p \to V$ (which is surjective
by 2.1 a), we denote it by $\vf=\vf(L,V)$ as well.
We can also attach to $(L,V)$ an exact sequence
$$0\to \Ker \vf \to L\underset{\n Z_\infty} \to{\otimes}\p \overset{\vf}\to{\to} V \to
0\eqno{(2.2)}$$

Let $\Cal I\in \Cl(\w)$ be a class of ideals; we shall use the same notation $\Cal I$ to
denote a representative in the $\iota$-image of this class. Let $(\vf', L', V')$ be another lattice and $D$ a structure of a perfect $\Cal I$-pairing $<* , * >_D$ between $L$
and $L'$. Let us fix an isomorphism $$\alpha: \Cal I \underset{\n Z_\infty} \to{\otimes} \p \to \p\eqno{(2.2')}$$ $D$ extends via $\alpha$ to a perfect $\p$-pairing between $L\underset{\n Z_\infty}\to{\otimes}\p$
and $L'\underset{\n Z_\infty}\to{\otimes}\p$, we denote this pairing by $D_{\alpha, \infty}$.
\medskip
{\bf Definition 2.3.} Two lattices $(\vf, L, V)$ and $(\vf', L', V')$ are called $(\alpha, \Cal I)$-dual if there exists a perfect $\Cal I$-pairing $D$ between $L$
and $L'$ such that $\Ker \vf \subset L\underset{\n Z_\infty}\to{\otimes}\p$, $\Ker \vf' \subset L'\underset{\n Z_\infty}\to{\otimes}\p$ are mutually orthogonal with respect to $D_{\alpha, \infty}$.
\medskip
Let $(n,r)$, $(n',r')$ be the dimension and rank of $(\vf, L, V)$ and $(\vf', L', V')$ respectively. If they are $(\alpha, \Cal I)$-dual then $r'=r$, $n'=r-n$. There exists the following reformulation of the definition of duality. $D_{\alpha, \infty}$ induces an isomorphism
$\gamma_{\alpha, D}: (L\underset{\n Z_\infty}\to{\otimes}\p)^* \to
L'\underset{\n Z_\infty}\to{\otimes}\p$ (here and below for any object $W$ we
denote $W^*=\Hom_{\p}(W,\p)$ ).
\medskip
{\bf Property 2.4.} $(\vf, L, V)$ and $(\vf', L', V')$ are $(\alpha, \Cal I)$-dual iff there exists an isomorphism from $(\Ker \vf)^*$
to $V'$ making the following diagram commutative: $$\matrix 0 & \to &
V^* & \overset{\vf^*}\to{\to} &  (L\underset{\n Z_\infty}
\to{\otimes}\p)^* & \to & (\Ker \vf)^* & \to & 0\\
&&&&&&&& \\
& & \downarrow & & \gamma_{\alpha, D} \downarrow && \downarrow \\ &&&&&&&& \\ 0
& \to & \Ker \vf' & \rightarrow  & L'\underset{\n Z_\infty}\to{\otimes}\p &
\overset{\vf'}\to{\to} & V' & \to &
0\endmatrix \eqno{(2.5)}$$

Further, this property is equivalent to the following two conditions:
\medskip
{\bf 2.6.} $\dim V'=r-n$;
\medskip
{\bf 2.7.} The composition map $\vf'\circ \gamma_D \circ \vf^*:
V^* \to V'$ is 0.
\medskip
Both 2.4 and (2.6, 2.7) are obvious.
\medskip
{\bf Remark 2.8.} It is
easy to see that the functor $(\vf, L, V) \mapsto (\vf', L', V')$ is well-defined
on a subcategory
(not all lattices have duals, see
below) of the category of the triples $(\vf, L, V)$, it is contravariant and
involutive.
\medskip
{\bf 3. Explicit formulas for analytic duality.}
\medskip
Here we consider the case $\w=\n
F_q[T]$. In this case $\Cl(\w)=0$, and $(\alpha, \Cal I)$-dual is called simply
dual. The coordinate
description of the dual lattice is the following. Let
$e_1, ..., e_r$ be a $\n Z_\infty$-basis of
$L$ such that $\vf(e_1), ..., \vf(e_n)$ form a $\p$-basis of $V$. Like in the
theory of abelian
varieties, we denote by $Z=(z_{ij})$ the Siegel matrix whose lines are
coordinates of $\vf(e_{n+1}), ..., \vf(e_r)$ in the basis $\vf(e_1), ..., \vf(e_n)$, more
exactly, the size of $Z$ is $(r-n)\times n$ and
$$\forall i =1,..., r-n \ \ \ \ \vf(e_{n+i})=\sum_{j=1}^n z_{ij}\vf(e_j)\eqno{(3.1)}$$ $Z$
defines $L$, we denote $L$ by $\goth L(Z)$.
\medskip
{\bf Proposition 3.2.} $[\goth L(Z)]'=\goth L(-Z^t)$, i.e. a Siegel matrix of the dual
lattice is the minus transposed Siegel matrix.
\medskip
{\bf Proof.} Follows immediately from the definitions. Really, let $f_1, ..., f_r$ be a
basis of $L'$, we define the pairing by the formula $$<e_i,
f_j>=\delta_i^j\eqno{(3.3)}$$ and the map $\vf'$ by the formula
$$\forall i =1,... ,n \ \ \ \ \vf'(f_{i})=\sum_{j=1}^{r-n} -z_{ji}\vf'(f_{n+j})$$
(minus transposed Siegel matrix).
$\Ker \vf$ is generated by elements $$v_i=e_{n+i}-\sum_{j=1}^n z_{ij}e_j, \ \ \ \ i
=1,..., r-n$$ and $\Ker \vf'$ is generated by elements $$w_i=f_{i}+\sum_{j=1}^{r-n}
z_{ji}f_{n+j}, \ \ \ \ i =1,..., n\eqno{(3.4)}$$ It is sufficient to check that $\forall i,j$ we have
$<v_i, w_j>=0$; this follows immediately from 3.3. $\square$
\medskip
{\bf Remark 3.5.} $L'$ exists not for all $L$.
Trivial counterexample: case $n=r=1$. To get another counterexamples, we use that for
$n=1$ (lattices of Drinfeld modules) a Siegel matrix is a column matrix $Z=\left(\matrix
z_1&...&z_{r-1} \endmatrix \right)^t$ and
$$ \goth L(Z) \hbox{ is not a lattice } \iff 1,z_1, ... ,z_{r-1} \hbox{ are linearly
dependent over } \n R_\infty \eqno{(3.6)}$$ while for $n=r-1$ a Siegel matrix is a
row matrix $Z=\left(\matrix -z_1&...&-z_{r-1} \endmatrix \right)$ and
$$ \goth L(Z) \hbox{ is not a lattice } \iff \forall i \ \ z_i\in \n R_\infty
\eqno{(3.7)}$$ Since condition (3.7) is strictly stronger than (3.6) we see
that all lattices having $n=1$, $r>1$ have duals while not all lattices having $n=r-1$,
$r>2$ have duals.

It is clear that almost all matrices have duals. Here "almost all" has the same meaning
that as "Almost all matrices $Z$ are a Siegel matrice of a lattice", i.e. if we choose an
(infinite) basis of $\p/\n R_\infty$, then coordinates of the entries of $Z$ in this
basis must satisfy some polynomial relations in order that $Z$ is not a Siegel matrice of
a lattice.
\medskip
{\bf Remark 3.8.} The coordinate proof of the theorem that the notion
of the dual lattice is well-defined, is the following. Two Siegel matrices
$Z$, $Z_1$ are called equivalent iff there exists an isomorphism of their pairs $(\goth
L(Z), V)$, $(\goth L(Z_1),V_1)$. Like in the classical theory of modular forms, $Z$,
$Z_1$ are equivalent iff there exists a matrix $\gamma \in GL_r(\n Z_\infty)=\left(\matrix A&B\\ C&D \endmatrix \right)$ ($A,B,C,D$ are the ($n\times
n$), ($n\times r-n$), ($r-n\times n$), ($r-n\times r-n$)-blocks of $\gamma$ respectively;
we shall call this block structire by the $(n, r-n)$-block structure) such that
$$C+DZ=Z_1(A+BZ)\eqno{(3.8.1)}$$

Let $A_1,B_1, C_1, D_1$ be the $(n, r-n)$-block structure of the matrix $\gamma^{-1}$.
The equality
$$-C_1^t+A_1^tZ^t={Z_1}^t(D_1^t-B_1^tZ^t)\eqno{(3.8.2)}$$
shows that if $Z$, $Z_1$ are equivalent
then $-Z^t$, $-Z_1^t$ are equivalent. [Proof of (3.8.2): (3.8.1) implies $Z_1=(C+DZ)(A+BZ)^{-1}$; substituting this value of $Z_1$ to the transposed (3.8.2), we get $-C_1+ZA_1=(D_1-ZB_1)(C+DZ)(A+BZ)^{-1}$, or $(-C_1+ZA_1)(A+BZ)=(D_1-ZB_1)(C+DZ)$. This formula follows immediately from $\left(\matrix A_1&B_1\\C_1&D_1\endmatrix \right)\left(\matrix A&B\\C&D\endmatrix \right)=\left(\matrix E_n&0\\0&E_{r-n}\endmatrix \right)$].  

Further, let $\alpha: (L_1\subset \p^n) \to
(L_2\subset \p^n)$ be a map of lattices. If $L'_1$, $L'_2$
exist, then the map
$\alpha': (L'_2\subset \p^{r-n}) \to (L'_1\subset \p^{r-n})$ is defined by the
following formulas. Let $Z_i$ be the Siegel matrices of $L_i$ in the
bases
$e_{i1}, ... e_{ir}$ of $L_i$ ($i=1,2$). Let us consider the matrix
$\goth M=(m_{ij})\in
M_{r}(\n Z_\infty)$ of $\alpha$ in the bases $e_{i1}, ...,
e_{ir}$ (i.e.
$\alpha(e_{1i})=\sum_j m_{ij} e_{2j}$). Let $f_{i1}, ..., f_{ir}$ be the dual
base of $L'_i$ (see 3.3) and $e'_{i1}, ... e'_{ir}$ another base of $L'_i$ defined by $$e'_{ij}=f_{i,j+n}, \ \ \ \ j+n \mod r\eqno{(3.8.3)}$$ Formulas (3.8.3), (3.4) show that an analog of 3.1 is satisfied for both bases $e'_{i1}, ..., e'_{ir}$, their Siegel matrices are $-Z_i^t$.

Let
$$\goth M=\left(\matrix \goth M_{11} & \goth M_{12} \\ \goth M_{21} &
\goth M_{22}
\endmatrix \right)$$
be the $(n, r-n)$-block structure of $\goth M$. The matrix of $\alpha'$ in the bases $f_{i1}, ..., f_{ir}$ is $\goth M^t$, and using the matrix 3.8.3 of change of base, we get that $\goth M'$ --- the matrix of $\alpha'$ in
the bases $e'_{i1}, ..., e'_{ir}$ --- has the following $(r-n,
n)$-block structure:
$$\goth M'=\left(\matrix \goth M_{22}^t & \goth M_{12}^t \\
\goth M_{21}^t & \goth
M_{11}^t \endmatrix \right)\eqno{(3.8.4)} $$
The property that $\goth M$ comes from a $\p$-linear map $\p^n \to \p^n$ implies
that $\goth M'$ comes from a $\p$-linear map $\p^{r-n} \to
\p^{r-n}$. This follows immediately from the definition of dual lattice, or can be easily checked algebraically.
\medskip
{\bf Remark 3.9.} Taking $\gamma =\left(\matrix 1&0\\ 0&-1 \endmatrix \right)$ we get
that $Z$ is equivalent to $-Z$, hence $Z'$ is also a Siegel matrix of the dual lattice.
\medskip
{\bf 4. Main conjecture for arbitrary $\w$}.
\medskip
The main result of the paper is the following Theorem 5 on coincidence of algebraic and
analytic duality. We formulate it as a conjecture 4.1 for any $\w$, but we prove it only for the case $\w=\n
F_q[T]$.  Let $M$ be a uniformizable t-motive. Its lattice $L(M)$ is really a lattice in
the meaning of Definition 2.1, because [A],
Corollary 3.3.6 (resp. [G], Lemma 5.9.12) means that it satisfies 2.1a (resp. 2.1b);
recall that we consider the case $N=0$,
i.e. the action of $T$ on $\Lie(M)$ is simply multiplication by $\theta$.
Let us fix (like in 1.13) $\goth C=(\goth C, \tau_{\goth C})$ --- a t-motive of rank 1 over
{\bf A}, and let $L(\goth C)$ be its lattice. It is a $\n Z_\infty$-module. $\Omega=\Omega(\w)$ is an $\w$-module, we consider a $\n Z_\infty$-module $\iota^{-1}(\Omega)$. There exists the   notion of the $L(\goth C)\otimes \iota^{-1}(\Omega)$-duality.
\medskip
{\bf Conjecture 4.1.} Let $M$ be a uniformizable t-motive having $N=0$ such
that its $\goth C$-dual $M'$ exists. Then $M'$ is uniformizable, it has $N':=N(M')=0$, and $(L(M), \Lie(M))$ and $(L(M'), \Lie(M'))$ are $\alpha, L(\goth C)\otimes \iota^{-1}(\Omega)$-dual for some $\alpha$ from $2.2'$ (it can be explicitly described).
\medskip
We prove in Section 5 the first step of the proof of this conjecture.
\medskip
{\bf Remark 4.2.} It is possible to generalize the above conjecture to the case of non-uniformizable $M$, $M'$. The pairing is defined between $\Hom_{\w_C[\tau]}(M,Z_1)$ and $\Hom_{\w_C[\tau]}(M',Z_1)$ (see (5.2.1a) for the definition of $Z_1$), or, the same, between $M_a$ and $M'_a$ for any $a\in \w$ (see 5.1.6).
\medskip
{\bf 5. Main theorem.}
\medskip
Recall that the word "canonical" means "canonical up to multiplication by elements of $\n F_q^*$".
\medskip
{\bf Theorem 5.}\footnotemark \footnotetext{The proof of this theorem was inspired by a result of Anderson, see Section 6 for details.} Let $M$ be a uniformizable t-motive over $\w=\n F_q[T]$ having $N=0$ such
that its dual $M'$ exists and has $N':=N(M')=0$. Then $M'$ is uniformizable, and $(L(M), \Lie(M))$ and $(L(M'), \Lie(M'))$ are dual.
\medskip
{\bf Remark 5A.} Condition $N'=0$ holds for pure $M$ (Theorem 10.3) and for a large class of non-pure $M$ (Theorem 11.5). Most likely, a modification of the end of the proof of the present theorem will permit us to prove that $N'=0$ holds for all $M$ having $N=0$ and having dual.
\medskip
{\bf Remark 5B.} A reformulation of this theorem in terms of Pink-Hodge structures, as well as its statement for $N\ne0$ (without proof), are given in Section 9 (Result 9.3 and Proposition 9.4 respectively).
\medskip
{\bf Corollary 5.1.1.} If $\w=\n F_q[T]$ then a Siegel matrix of $M'$ is the minus transposed
of a Siegel matrix of $M$.
\medskip
In the section 8 below we give a corollary of this theorem and some conjectures related to the
problem of 1 -- 1 correspondence between t-motives and lattices.
\medskip
{\bf 5.1.2. Some definitions.} Recall that $E=E(M)$ is isomorphic to $\p^n$. There is a structure of $\w$-module on $E$; multiplication by $T$ is denoted by $m_T$, and this operator $m_T$ is defined in coordinates by the formula
$$m_T(x)=\sum_{i=0}^l\goth A_ix^{(i)}$$ where $x\in E=\p^n$ is a vector column, $\goth A_i$ are from 1.9.1. There is a map $\exp: \Lie(M) \to E$ making the following diagram commutative:
$$\matrix \Lie(M) &
\overset{\Exp}\to{\to} & E \\ \\ \theta \downarrow & & m_T
\downarrow \\ \\ \Lie(M) & \overset{\Exp}\to{\to} & E \endmatrix \eqno{(5.1.3)}$$
By definition, $L(M)=\Ker \Exp$.

We need another space $\Lie_T(M)$ together with an isomorphism $\goth a:\Lie_T(M) \to \Lie(M)$ and a structure of $\w$-module on $\Lie_T(M)$ such that the multiplication by $T$ on $\Lie_T(M)$ is simply the multiplication by $\theta$ on $\Lie(M)$, i.e. $$\goth a(Tx)=\theta\cdot(\goth a(x))\eqno{(5.1.4)}$$ where $x\in\Lie_T(M)$. Commutativity of 5.1.3 means that $\Exp\circ\goth a:\Lie_T(M)\to E$ is a map of $\w$-modules.
\medskip
{\bf 5.1.5.} We shall work merely with $L_T(M):=\Ker (\Exp\circ\goth a)\subset \Lie_T(M)$ rather than $L(M)$. Clearly $L_T(M)$ is an $\w$-module, $\goth a:L_T(M)\to L(M)$ is an isomorphism satisfying 5.1.4 for $x\in L_T(M)$.
\medskip
The proof of Theorem 5 consists of two steps. We formulate and prove Step 1 for the case of arbitrary $\w$.
\medskip
{\bf Step 1.} For the above $M$, $M'$ we have:
\medskip
(A) Uniformizability of $M$ implies uniformizability of $M'$.
\medskip
(B) There exists a canonical $\w$-linear $L_T(\goth C)\otimes
\Omega$-valued perfect pairing $<* , * >_M$
between $L_T(M)$ and $L_T(M')$ (by 5.1.5, this is the same as the $\n Z_\infty$-linear pairing between $L(M)$ and $L(M')$, which, in its turn, is $D$ of Definition 2.3). It is functorial.
\medskip
{\bf Remark 5.1.6.} Practically, (B) comes from [T], Theorem 4.3 (case
$\w=\n F_q[T]$). Really, to define a pairing between $L(M)$ and
$L(M')$ it is sufficient to define (concordant) pairings between $L(M)
/aL(M)$ and $L(M') /aL(M')$ for any $a\in \w$. Since $M_a:=E(M)_a=L(M)/aL(M)$ and because of Proposition 1.12.3 which affirms that
$M_a$ and $M'_a$ are Taguchi-dual, we see that [T], Theorem 4.3 gives
exactly the desired pairing.
\medskip
We give two versions of the proof of Step 1: the first one --- for the general case of arbitrary $\w$ and the second one --- for the case $\w=\n F_q[T]$ --- is based on explicit calculations, it is used for the proof of Step 2.
\medskip
{\bf 5.2. Proof: Step 1, Version 1.} Here we consider the general case of arbitrary $\w$.
Let $\Omega=\Omega(\w/\n
F_q)$ be the module of differential forms; we can consider it as an
element of $\Cl(\w)$. We use formulas and notations of [G], Section 5.9
modifying them to the case of arbitrary $\w$. For example, {\bf
A} (resp. {\bf K}) of [G], 5.9.16 is {\bf A} (resp. {\bf K})
of the present paper (recall that $\bar K$ (resp. $\bar K[T,\tau]$) of [G] is
$\p$ (resp. $\w_C[\tau]$, see 1.13) of the present paper). Hence, we denote $\bar K\{T\}$ of
[G], Definition 5.9.10 by $\p\{T\}$. For the general case it must be
replaced by a ring $Z_0$ defined by the formula
$$Z_0:=\w\underset{\n F_q[T]}\to{\otimes}\p\{T\}\eqno{(5.2.1)}$$
$Z_0$ is a $\w_C[\tau]$-module, i.e. $\tau$ acts on $Z_0$, and $Z_0^\tau=\w$.

$Z_1$ for the present case is defined by the same formula [G], 5.9.22. Explicitly,
$$Z_1:=\Hom^{cont}_{\w}(\x/\w,\p)\eqno{(5.2.1a)}$$
It is a locally free $Z_0$-module of dimension 1 (the module structure
is compatible with the action of $\tau$; see [G], p. 168, lines 3 - 4
for the case $\w=\n F_q[T]$). We have: $Z_1^\tau$ is a
$Z_0^\tau$-module ( = $\w$-module) which is isomorphic to $\Omega(\w)$
(see the last lines of the proof of [G], Corollary 5.9.35 for the case
$\w=\n F_q[T]$), and $Z_1$ is isomorphic to $Z_0\otimes_{\w}\Omega(\w)$.

We shall consider $M$ as a $\w_C[\tau]$-module, like in 1.13.3. We denote
$M\{T\}:=M\underset{\w_C}\to{\otimes}Z_0$ ( = [G],
Definition 5.9.11.1 for the case $\w=\n F_q[T]$) and
$H^1(M):=M\{T\}^\tau$ like in [G], Definition 5.9.11.2.
Analogous to [G], Corollary 5.9.25 we get that for the present case
$$H_1(M):=\Hom_{\w_C[\tau]}(M,Z_1)=L_T(M)$$
($H_1(M)=H_1(E)$ of [G], 5.9). Particularly, for $M=\goth C$ we have
$$L_T(\goth C)=\Hom_{\w_C[\tau]}(\goth C,Z_1)$$

{\bf Lemma 5.2.2. } $H_1(M')=H^1(M)\underset{\w}\to{\otimes}L_T(\goth C)$.
\medskip
{\bf Proof.} By definition,
$\Hom_{\w_C} (M',Z_1)=\Hom_{\w_C} (\Hom_{\w_C} (M, \goth C), Z_1)$. Further,
$$\Hom_{\w_C} (\Hom_{\w_C} (M, \goth C), Z_1)=(M\underset{\w_C}
\to{\otimes}Z_0)\underset{Z_0}\to{\otimes} (\Hom_{\w_C}(\goth C,
Z_1))\eqno{(5.2.3)}$$ (an equality of  linear algebra). In order to show that we can
consider $\tau$-invariant subspaces, we need the following objects.
Let $I$ be an ideal of $\w$, $\Cal M_0=IZ_0$. It is clear that $\Cal M_0^\tau=I$.
Further, let $\Cal M_1$ be a locally free $Z_0$-module. We have a formula: $$(\Cal
M_0\underset{Z_0}\to{\otimes} \Cal M_1)^\tau=\Cal M_0^\tau \underset{\w}\to{\otimes} \Cal
M_1^\tau \eqno{(5.2.4)}$$
Really, $\Cal M_0\underset{Z_0}\to{\otimes} \Cal M_1=I\Cal M_1$, and $$(I\Cal
M_1)^\tau=I\Cal M_1^\tau\eqno{(5.2.5)}$$ where this formula is true by the following
reason. Obviously $(I\Cal M_1)^\tau\supset I\Cal M_1^\tau$. Let $J$ be an ideal of $\w$
such that $IJ$ is a principal ideal. We have $(IJ(J^{-1}\Cal M_1))^\tau=IJ(J^{-1}\Cal
M_1)^\tau$ and $(IJ(J^{-1}\Cal M_1))^\tau\supset I(J(J^{-1}\Cal M_1))^\tau\supset
IJ(J^{-1}\Cal M_1)^\tau$, hence all these objects are equal and we get 5.2.5 and hence
5.2.4.

The action of $\tau$ on both sides of 5.2.3 coincide.  Considering $\tau$-invariant
elements of both sides of 5.2.3 and taking into consideration 5.2.4 ($\Cal
M_0=\Hom_{\w_C}(\goth C,
Z_1)$ and $\Cal M_1=M\underset{\w_C}
\to{\otimes}Z_0$) we get the lemma. $\square$
\medskip
This lemma proves (A) of Step 1.
\medskip
{\bf Lemma 5.2.6.} Let $\Cal M_i$ ($i=0,1$) be two locally free $Z_0$-modules with $\tau$-action satisfying $\tau(cm)=\tau(c) \tau(m)$ ($c\in Z_0$, $m\in \Cal M_i$), and $\psi: \Cal M_0\otimes_{Z_0}\Cal M_1\to Z_1$ a perfect pairing of $Z_0$-modules with $\tau$-action. Let, further, both $\Cal M_i$ satisfy $\Cal M_i^\tau\otimes_{\w} Z_0=\Cal M_i$. Then the restriction of $\psi$ to $\Cal M_0^\tau\otimes_{\w}\Cal M_1^\tau\to \Omega$ is a perfect pairing as well.
\medskip
{\bf Proof.} Let $\alpha: \Cal M_0^\tau \to \Omega$ be an $\w$-linear map. We prolonge it to a map $\bar \alpha: \Cal M_0 \to Z_1$ by $Z_0$-$\tau$-linearity. By perfectness of $\psi$, there exists $m_1\in \Cal M_1$ such that $\bar \alpha(m_0)= \psi(m_0 \otimes m_1)$. It is easy to see that $m_1$ is $\tau$-invariant (we use the fact that $\tau: Z_0 \to Z_0$ is surjective). $\square$
\medskip
{\bf Lemma 5.2.7.} There is a natural perfect $\w$-linear
$\Omega$-valued pairing
between $H_1(M)$ and $H^1(M)$:
$H_1(M)\underset{\w}\to{\otimes}H^1(M)\to \Omega$.
\medskip
{\bf Proof.} For the case $\w=\n F_q[T]$ this is [G], Corollary
5.9.35. General case: we have a perfect $Z_0$-pairing
$$\Hom_{\w_C}(M,Z_1)\underset{Z_0}\to{\otimes} (M\underset{\w_C}
\to{\otimes} Z_0) \to Z_1$$
Now we take $\Cal M_0=\Hom_{\w_C}(M,Z_1)$, $\Cal M_1=M\underset{\w_C} \to{\otimes}Z_0$ and we apply Lemma 5.2.6. $\square$
\medskip
Step 1 of the theorem follows from these lemmas. 
\medskip
{\bf Remark 5.2.8.} The pairing can be defined also as the composition of
$$\matrix H_1(M)\underset{\w} \to{\otimes} H_1(M')=\Hom_{\w_C[\tau]}(M,Z_1)
\underset{\w}
\to{\otimes} \Hom_{\w_C[\tau]}(M',Z_1) \\
\to \Hom_{\w_C[\tau]}(M\underset{\w_C} \to{\otimes}M'
,Z_1\underset{Z_0} \to{\otimes}Z_1) \to
\Hom_{\w_C[\tau]}(\goth C,Z_1\underset{Z_0} \to{\otimes}Z_1) = L_T(\goth C) \underset{\w} \to{\otimes}\Omega\endmatrix \eqno{(5.2.9)}$$
where the second map comes from a canonical map $\delta: \goth C \to M\underset{\w_C}
\to{\otimes}M'$ of Remark 1.11.1 (more exactly, of its analog for arbitrary $\w$).
\medskip
{\bf Remark 5.2.10.} Recall that the explicit formula for functoriality is
the following. Let $\alpha: M_1 \to M_2$ be a map of t-motives,
$\alpha': M'_2 \to M'_1$
the dual map and $L_T(\alpha): L_T(M_2) \to L_T(M_1)$, $L_T(\alpha'): L_T(M'_1)
\to L_T(M'_2)$ the
corresponding maps on lattices. For any $l_1'\in L_T(M_1')$, $l_2\in L_T(M_2)$ we
have:
$$<L_T(\alpha)(l_2), l_1'>_{M_1}=<l_2, L_T(\alpha')(l_1')>_{M_2}\eqno{(5.2.11)}$$

{\bf 5.3. Proof: Step 1, Version 2.} Case $\w=\n F_q[T]$. We identify $Z_1$ of [G], p.168, lines 3 -- 4 with $\p\{T\}$ (see [G], Definition 5.9.10) and $\w$ with $\Omega$. Like above, we have an isomorphism of $\w$-modules (recall that $\w$ is the center of $\p[T,\tau]$):
$$L_T(M)=\Hom_{\p[T,\tau]}(M,Z_1)\eqno{(5.3.1)}$$ ([G], first terms of 5.9.25, 5.9.19). Let $\vf: M \to Z_1$, $\vf':
M' \to Z_1$ be elements of $L_T(M)$, $L_T(M')$ respectively, and let
$f_*$, $f'_*$, $Q$, $Q'$ be from
1.9.3, 1.10. We denote $$\vf(f_*)=v_*\eqno{(5.3.2)}$$ where $v_*\in
(Z_1)^r$ is a vector column (it is a column of the scattering matrix ([A], p. 486) of
$M$, see 5.4.1 below). The same notation for
the dual:
$\vf'(f'_*)=v'_*$. Condition
that $\vf$, $\vf'$ are
$\tau$-homomorphisms is equivalent to $$Qv_*=v_*^{(1)}, \ \
Q'v'_*={v'}_*^{(1)}\eqno{(5.3.3)}$$ (analog of the formula for scattering matrices [A],
(3.2.2)). Let us consider $\Xi=\sum_{i=0}^\infty
a_iT^i\in\p\{T\}\subset\p[[T]]$ of [G], p. 172, line 1; recall that it is the only
element (up to
multiplication by $\n F_q^*$) satisfying
$$\Xi=(T-\theta)\Xi^{(1)}, \ \ \ \lim_{i\to\infty}a_i=0, \ \ \ |a_0|>|a_i| \ \ \forall
i>0\eqno{(5.3.4)}$$ (see [G], p. 171,
(*); there is a formula $\Xi=a_0\prod_{i\ge0}(1-T/\theta^{q^i})$ where $a_0$ satisfies
$a_0^{q-1}=-1/\theta$). Finally, we define
$$<\vf, \vf'>=\Xi v_{*}^tv'_{*}\eqno{(5.3.5)}$$
Obviously $<\vf, \vf'>$ does not depend on a choice of a basis $f_*$.
\medskip
{\bf Lemma 5.3.6.} $<\vf, \vf'>\in \w$.
\medskip
{\bf Proof.} Firstly, this element belongs to $\n F_q[[T]]$, because $$\Xi
v_{*}^tv'_{*} -(\Xi v_{*}^tv'_{*})^{(1)}=\Xi (v_{*}^tv'_{*}-(T-\theta)^{-1}
v_{*}^{(1)t}{v'}_{*}^{(1)})= \Xi v_{*}^t(E_r-(T-\theta)^{-1}
Q^tQ')v'_{*}$$ because of (5.3.3). But we have (see (1.10.1) --- the
definition of $Q'$)
$$E_r-(T-\theta)^{-1} Q^tQ'=0$$
Secondly, let $<\vf, \vf'>=\sum_{i=0}^\infty c_iT^i$. Since
coefficients of all factors of (5.3.5): $\Xi$, $v_*$ and $v'_*$ ---
tend to 0, we get
that $c_i$ also tend to 0. But $c_i\in\n F_q$, i.e. they are almost
all 0. $\square$
\medskip
{\bf Lemma 5.3.7.} The above pairing is perfect.
\medskip
{\bf Proof.} We have an isomorphism (here $M\{T\}=M\otimes_{\p[T]}\p\{T\}$ with
the natural action of $\tau$ (see [G], Definition 5.9.11))
$$\alpha: \Hom_{\p[T,\tau]}(M,Z_1)\to\Hom_{\w}(M\{T\}^\tau,\w)\eqno{(5.3.8)}$$ defined as the composition of the maps
$$\Hom_{\p[T,\tau]}(M,Z_1)=\Hom_{\p[T]}(M,Z_1)^\tau\overset{\beta'}\to{\to}\Hom_{\p\{T\}}(M\{T\},\p\{T\})^\tau$$
$$\overset{\gamma}\to{\to} \Hom_{\w}(M\{T\}^\tau,\w)$$ where
$\beta: \Hom_{\p[T]}(M,Z_1)\to \Hom_{\p\{T\}}(M\{T\},\p\{T\} )$ is the natural map and
$\beta'$
is the restriction of $\beta$ to $\tau$-invariant elements. Using the Anderson's
criterion of uniformizability
of $M$ (see, for example, [G], 5.9.14.3 and 5.9.13) we get immediately that both
$\gamma$, $\beta$, and hence
$\beta'$, and hence $\alpha$ are isomorphisms. Further, let us consider a homomorphism
$$i: \Hom_{\p[T,\tau]}(M',Z_1)\to M\{T\}^\tau\eqno{(5.3.9)}$$ defined as follows. Let
$\vf'$, $f'_*$, $v'_*$ be as above. We set
$$i(\vf')=\Xi{v'}^t_*f_*\in M\underset{\p[T]}\to{\otimes}\p[[T]]$$ Since $\Xi\in \p
\{T\}$, we get that $\Xi{v'}^t_*f_*\in M\{T\}$. A simple calculation
(like in the Lemma 5.3.6, but simpler) shows that $i(\vf')$ is $\tau$-invariant, hence
$i$ really defines a map from $\Hom_{\p[T,\tau]}(M',Z_1)$ to $M\{T\}^\tau$. Obviously it
is an inclusion. Let us prove that $i$ is surjective. Really, let $c_*\in (Z_1)^r$ be a
column
vector such that $c_*^tf_*\in M\{T\}^\tau$. An analog of the above calculation shows that
if we define $\vf'$ by the formula $\vf'(f'_*)=\Xi^{-1}c_*$ then $\vf'\in
\Hom_{\p[T,\tau]}(M',Z_1)$, and $i(\vf')=c_*^tf_*\in M\{T\}^\tau$. Finally, the
combination of isomorphisms (5.3.8) and (5.3.9) corresponds to the
pairing (5.3.5). $\square$
\medskip
{\bf 5.4. Step 2 -- End of the proof of Theorem 5.} It is easy to see that the converse of the Corollary 5.1.1 (taking into consideration
Proposition 3.2) is also true, i.e. in order to prove Theorem 5 it is sufficient
to prove that a Siegel matrix of $M'$ is $-Z^t$ where $Z$ is a Siegel matrix of $M$. Let
us consider a basis $l_1,...,l_r$ of $L_T(M)$ and for each $l_i$ we consider the
corresponding (under identification 5.3.1) $\vf_i\in \Hom_{\p[T,\tau]}(M,Z_1)$. Let
$\Psi$ be the scattering matrix of $M$ ([A], p. 486) with respect to the bases $l_1,...,l_r$,
$f_1,...,f_r$, and we denote $\vf_i(f_*)$ by $v_{i*}$ (notations of 5.3.2).
\medskip
{\bf Lemma 5.4.1.} $v_{i*}$ is the $i$-th column of $\Psi$ ($Z_1$ is identified with $\p\{T\}$, see the proof).
\medskip
{\bf Proof.} Follows from the definitions. Recall that $\hbox{\bf K}=\n F_q((1/T))$. The isomorphism 5.3.1 is the composition of
2 isomorphisms $i_1: L_T(M) \to \Hom_{\w}^c(\hbox{\bf K}/\w, E)$ ([G], 5.9.19) and $i_2:
\Hom_{\w}^c(\hbox{\bf K}/\w, E) \to \Hom_{\p[T,\tau]}(M,\Hom^c(\hbox{\bf K}/\w, \p)$
([G], 5.9.24; recall that $Z_1=\Hom^c(\hbox{\bf K}/\w, \p)$). For $l_i\in L_T(M)$ we have
$(i_1(l_i))(T^{-k})=\exp(\theta^{-k}l_i)$ ([G], line above the lemma 5.9.18) and
$$((i_2\circ i_1(l_i))(f_j))(T^{-k})=<f_j,\exp(\theta^{-k}l_i)>$$ ([G], two lines above the
lemma 5.9.24). Using the identification of $Z_1$ and $\p\{T\}$ ([G], p. 168, lines 3 - 4)
and the definition of $\Psi$ ([A], p. 486, first formula of 3.2) we get immediately the
lemma. $\square$
\medskip
Let $l'_1,...,l'_r$ be a basis of $L_T(M')$ which is dual to a basis $l_1,...,l_r$ of
$L_T(M)$ with respect to the pairing 5.3.5.
\medskip
{\bf Lemma 5.4.2.} The scattering matrix of $M'$ with respect to the bases
$l'_1,...,l'_r$, $f'_1,...,f'_r$ (denoted by $\Psi'$) is $\Xi^{-1}\Psi^{t-1}$.
\medskip
{\bf Proof.} Follows immediately from 5.4.1 applied to both $M$, $M'$, and formula 5.3.5.
$\square$
\medskip
{\bf Remark 5.4.3.} An alternative proof for the case of pure $M$ (for $some$ basis of $L_T(M')$)
is the following. We denote $\Xi^{-1}\Psi^{t-1}$ by $\Psi_1$. It satisfies
$\Psi_1^{(1)}=(T-\theta) Q^{t-1}\Psi_1$ and other conditions of [A], 3.1. According [A],
Theorem 5, p. 488, there exists a pure uniformizable t-motive $M_1$ with
$\sigma$-structure such that its scattering matrix is $\Psi_1$. Since $\Psi_1$ satisfies
$$\Psi_1^{(1)}=Q'\Psi_1$$
we get that $Q(M_1)=Q'$, i.e. $M_1=M'$. $\square$
\medskip
Let us recall the statement of the crucial proposition 3.3.2 of [A]. Here we consider the
case of those $M$ whose $N$ is not necessarily 0. Let $\Psi$ be a scattering matrix of
$M$. We consider the $(T-\theta)$-Laurent series for $\Psi$ (here $k(M)<0$ is a number,
and $D_{-i}\in M_{r}(\p)$): $$\Psi=\sum_{i=k(M)}^\infty D_{-i}(T-\theta)^i$$ We consider its
negative part $$\Psi^-:=\sum_{i=k(M)}^{-1} D_{-i}(T-\theta)^i$$ as an element of $M_{r}(\p)((T-\theta))/M_{r}(\p)[[T-\theta]]$.
\medskip
We consider the space $(T-\theta)^{k(M)}\p[[T-\theta]]/\p[[T-\theta]]$ as a $\p$-vector space endowed by the action of $\w$, and we denote by $\goth V$ its $r$-th direct sum written as vector columns of length $r$. Obviously $$k(M)=-1 \iff \hbox{ the action of $T$ on $\goth V$ coincides with multiplication by $\theta$} \eqno{(5.4.3a)}$$
We denote the $i$-th column of $\Psi^-$ by $\Psi^-_{i*}$, it belongs to $\goth V$. Further, we denote by $\Prin(M)$ (resp. by $\Prin_0(M)$) the $\p[T]$-linear envelope (resp. the $\w$-linear envelope) of all $\Psi^-_{i*}$ in $\goth V$. Finally, we obviously extend the definition of $\Lie_T(M)$, $L_T(M)$ to the case $N\ne 0$; formula 5.1.4 becomes
$$\goth a(Tx)=(\theta+N)(\goth a(x))\eqno{(5.4.3b)}$$
\medskip
{\bf Proposition 3.3.2, [A]} (see also Remark 5.5 below). There exists a $\p[T]$-linear isomorphism $\psi_E: \Lie_T(M)
\to \Prin(M)$ such that its restriction to $L_T(M)\subset \Lie_T(M)$ defines an isomorphism
$L_T(M) \to \Prin_0(M)$ (denoted by $\psi_E$ as well). $\square$
\medskip
{\bf Corollary 5.4.4.} $N=0 \iff k(M)=-1$ (because $N=0 \iff $ the action of $T$ on both $\Lie_T(M)$, $\goth V$ coincides with multiplication by $\theta$, by 5.4.3a). $\square$
\medskip
We return to the case $N=0$.
\medskip
Let us consider the $(T-\theta)$-Laurent series for $\Psi'$ and $\Xi^{-1}$:
$$\Psi'=\sum_{i=k(M')}^\infty D'_{-i}(T-\theta)^i,  \ \ \ \Xi^{-1}=\sum_{i=k(\xi)}^\infty
a_i(T-\theta)^i$$
Since for both $M$, $M'$ we have $N=N'=0$, we get $k(M)=k(M')=-1$. An
elementary calculation shows that $k(\xi)$ is also $-1$. Hence, equality
$\Psi'\Psi^t=\Xi^{-1}$ (Lemma 5.4.2) implies that $D'_{1}D^t_{1}=0$.

Further, there exist $n$ columns of $D_{1}$ which are
$\p$-linerly independent (they are $\psi_E$-images of elements of $L_T(M)$ which form a
$\p$-basis of $\Lie_T(M)$) and all other columns of $D_{1}$ are their linear combinations.
Interchanging columns of $D_{1}$ if necessary we can assume that these columns are the
last $n$
columns. We denote by $D_{12}$ (resp. $D_{11}$ ) the $r\times n$ (resp. $r\times
(r-n)$ ) matrix formed by the last $n$ (resp. the first $r-n$) columns of $D_{1}$. There
exists a matrix $S$ such that $D_{11}=D_{12}S^t$. Again according Proposition 3.3.2,
[A], we have:
$$S \hbox{ is a Siegel matrix of $L(M)$ }\eqno{(5.4.5)}$$
(see also Remark 5.5 below).

Analogous objects are defined for $D'_{1}$. We denote by $D'_{12}$ (resp.
$D'_{11}$) the $r\times n$- (resp. $r\times (r-n)$)-matrix formed by the last $n$
(resp. the first $r-n$) columns of $D'_{1}$. Since
$D'_{1}D^t_{1}=D'_{11}D_{11}^t+D'_{12}D_{12}^t$ we get that
$D'_{12}D_{12}^t+D'_{11}SD_{12}^t=0$. Since $D_{12}^t$ is a $n\times
r$-matrix of rank $n$, it is not a zero-divisor from the right, so $$D'_{12}=
-D'_{11}S\eqno{(5.4.6)}$$ Since the rank of $D'_{1}$ is $r-n$ and $D'_{11}$ is a
$r\times (r-n)$ matrix, (5.4.6) implies that columns of $D'_{11}$ are linearly independent,
and by (5.4.6) and Proposition 3.3.2, [A] we get that $-S$ is a Siegel matrix of
$M'$. $\square$
\medskip
{\bf Remark 5.5.} Since the notations of [A] differ from the ones of the present
paper, for the reader's convenience we give here a sketch of the proof for the case $N=0$
of two  crucial facts: Corollary 5.4.4 and 5.4.5 ([A], Theorem 3.3.2).

Let $\alpha: \Lie(M)\to E(M)$ be a linear isomorphism which is the first term of the
series for $\exp: \Lie(M)\to E(M)$, and let $l\in \Lie(M)$, $f \in M$ be arbitrary. We
consider the $(T-\theta)$-Laurent series $\sum_{i=k}^\infty b_i(T-\theta)^i$ of
$\sum_{j=0}^\infty <\exp(\frac{1}{\theta^{j+1}}l),f>T^j$.
\medskip
{\bf Lemma 5.6.} If $N=0$ then $k=-1$, and $b_{-1}=-<\alpha(l),f>$ (this is [A],
3.3.4).
\medskip
{\bf Sketch of the proof.} For $z\in \Lie(M)$ we denote $\exp(z)- \alpha(z)$ by $\ve(z)$,
hence
$\sum_{j=0}^\infty <\exp(\frac{1}{\theta^{j+1}}l),f>T^j=\underline{A} +\underline{E}$, where
$$\underline{A}=\sum_{j=0}^\infty <\alpha(\frac{1}{\theta^{j+1}}l),f>T^j; \ \ \ \ \underline{E}=\sum_{j=0}^\infty <\ve(\frac{1}{\theta^{j+1}}l),f>T^j$$
We consider their $(T-\theta)$-Laurent series:
$$\underline{A}=\sum_{i=k(\underline{A})}^\infty \underline{a}_i(T-\theta)^i; \ \ \ \ \underline{E}=\sum_{i=k(\underline{E})}^\infty \underline{e}_i(T-\theta)^i$$
Since we have $\exp(z)=\sum_{i=0}^{\infty}C_iz^{(i)}$ where $C_0=E_n$ we get that
$\ve(z)=\sum_{i=1}^{\infty}C_iz^{(i)}$. This means that for large $j$ the element
$\ve(\frac{1}{\theta^{j+1}}l)$ is small, and hence $k(\underline{E})=0$, because finitely many
terms having small $j$ do not contribute to the pole of the $(T-\theta)$-Laurent series
of $\underline{E}$ (the reader can prove easily the exact estimations himself, or to look [A],
p. 491). Since $\alpha$ is $\p$-linear, equality $\sum_{j=0}^\infty
\frac{1}{\theta^{j+1}}T^j= -(T-\theta)^{-1}$ implies that $k(\underline{A})=-1$ and $\underline{a}_{\ -1}=-<\alpha(l),f>$ (and other $\underline{a}_i=0$), hence the lemma. $\square$
\medskip
This lemma obviously implies Corollary 5.4.4. Further, elements $f_1,...,f_r$ generate
the $\p$-space $M/\tau M$, because multiplication by $T$ on $M/\tau M$ coincides with
multiplication by $\theta$, hence the fact that $f_1,...,f_r$ \ \ $\p[T]$-generate $M/\tau
M$ implies that they $\p$-generate $M/\tau M$.

Let $l_1,...,l_n$ form a $\p$-basis of $\Lie(M)$ (here we identify $\Lie_T(M)$ and $\Lie(M)$ via $\goth a$). Since the pairing $<*,*>$ between
$E(M)$ and $M/\tau M$ is non-degenerate and $\alpha$ is an isomorphism, we get that
columns $<\alpha(l_1),f_*>,..., <\alpha(l_n),f_*>$ are linearly independent. Again since
$\alpha$ is an isomorphism and the pairing with $f_*$ is linear, we get that
$$(<\alpha(l_{n+1}),f_*> \ ... \ <\alpha(l_r),f_*>)=(<\alpha(l_{1}),f_*> \ ... \
<\alpha(l_n),f_*>)Z^t $$ Applying the lemma 5.6 to this formula we get immediately 5.4.5.
\medskip
{\bf 6. Tensor products.}
\medskip
There exists an analog of the Theorem 5 for the case of tensor
products of
t-motives. It describes the lattice $L(M_1\otimes M_2)$ in terms of $L(M_1)$,
$L(M_2)$. This is a theorem of Anderson; it is formulated in [P], end of page 3, but its
proof was not published. We recall its statement for the case of arbitrary $N\ne 0$, and
we give its proof for the case $N=0$ (case of arbitrary $N$ can be obtained easily using
the same ideas).

Let $M$ be an uniformizable t-motive whose $N$ is not necessarily 0. Since $N$ is nilpotent, formula 5.4.3b shows that $\Lie_T(M)$ is a $\p[[T-\theta]]$-module. There exists an epimorphism of
$\p[[T-\theta]]$-modules
$$L_T(M)\underset{\w}\to{\otimes}\p[[T-\theta]]\to \Lie_T(M)$$
whose kernel $\goth q=\goth q(M)$ carries information on the pair $(L(M), \Lie(M))$.
\medskip
{\bf Theorem 6} (Anderson). Let $M$, $\bar M$ be any two uniformizable abelian
t-motives. Then $$\goth q(M\otimes \bar M)=\goth q(M)\underset{\p[[T-\theta]]}
\to{\otimes}\goth q(\bar M)\eqno{(6.1)}$$
\medskip
{\bf Remark 6A.} $M\otimes \bar M$ is a uniformizable t-motive ([G], Corollary 5.9.38).
\medskip
{\bf Proof of Theorem 6 (case $N=0$).} We define notations for $M$, and all notations for $\bar M$
will carry bar. Let $e_i$ and $Z$ be from the beginning of Section 3. We denote $\goth a^{-1}(e_i)\in \Lie_T(M)$ by $e_i$ (there is no possibility of confusion). So, $\{e_i\}$ is a
$\p[[T-\theta]]$-basis of $L_T(M)\underset{\w}\to{\otimes}\p[[T-\theta]]$. Elements
$b_i:=(T-\theta)e_i$, $i=1,...,n$ and $b_{n+i}:=e_{n+i}-\sum_{j=1}^n z_{ij}e_j$,
$i=1,...,r-n$ form a $\p[[T-\theta]]$-basis of $\goth q$. We need a
\medskip
{\bf Lemma 6.2.} $\Psi(M\otimes \bar M)=\Psi(M)\otimes\Psi(\bar M)$ where $\Psi(M)$
(resp. $\Psi(\bar M)$; $\Psi(M\otimes \bar M)$) is taken with respect to bases $e_*$ of
$L_T(M)$, $f_*$ of $M_{\p[T]}$ (resp. $\bar e_*$ of $L_T(\bar M)$, $\bar f_*$ of $\bar
M_{\p[T]}$; $e_*\otimes \bar e_*$ of $L_T(M\otimes \bar M)$, $f_*\otimes \bar f_*$ of
$(M\otimes \bar M)_{\p[T]}$) (see the proof for the notations).
\medskip
{\bf Proof.} We consider a map $$\alpha:\Hom_{\p[T]}(M,Z_1)^\tau
\underset{\w}\to{\otimes} \Hom_{\p[T]}(\bar M,Z_1)^\tau\to \Hom_{\p[T]}(M\otimes \bar
M,Z_1)^\tau$$ defined as follows: for $\vf\in \Hom_{\p[T]}(M,Z_1)^\tau$, $\bar \vf\in
\Hom_{\p[T]}(\bar M,Z_1)^\tau$ we let $[\alpha(\vf \otimes \bar \vf)](f\otimes \bar
f)=\vf(f)\cdot\bar \vf(\bar f)$ (it is obvious that $\alpha(\vf \otimes \bar \vf)$ is
$\tau$-stable). Since $e_1, ..., e_r$ (resp. $\bar e_1, ..., \bar e_{\bar r}$) is a basis
of $\Hom_{\p[T]}(M,Z_1)^\tau$ (resp. $\Hom_{\p[T]}(\bar M,Z_1)^\tau$; we identify $L_T(M)$,
resp. $L_T(\bar M)$ with $\Hom_{\p[T]}(M,Z_1)^\tau$ (resp. $\Hom_{\p[T]}(\bar M,Z_1)^\tau$)
we get (using Lemma 5.4.1) that $\Psi(M)$, $\Psi(\bar M)$ are non-degenerate. Since their
product is also non-degenerate, we get $\alpha(e_i\otimes \bar e_{\bar i})$ are linearly
independent and hence a basis of $\Hom_{\p[T]}(M\otimes \bar M,Z_1)^\tau$. Applying once
again Lemma 5.4.1 we get the lemma. $\square$
\medskip
If $A$, $B$ are two matrices then columns of $A\otimes B$ are indexed by pairs $(k,l)$
where $k$ (resp. $l$) is the number of a column of $A$ (resp. $B$). We denote by $A_k$,
$B_l$, $A\otimes B_{(k,l)}$ the respective columns. Obviosly we have: $A\otimes
B_{(k,l)}=A_k\otimes B_l$ (tensor product of column matrices).
\medskip
Let us prove that for $i=1,...,r-n$, $\bar i=1,...,\bar r-\bar n$ the element
$b_{n+i}\otimes \bar b_{\bar n+\bar i}\in \goth q(M\otimes \bar M)$. According [A],
Proposition 3.3.2, it is sufficient to prove that the corresponding linear combination
(see 6.3 below) of the columns of the matrix $\Psi^-_{M\otimes \bar M}$ is 0. Since
$$b_{n+i}\otimes \bar b_{\bar n+\bar i}=\sum_{j,\bar j}z_{ij}\bar z_{\bar i\bar
j}e_j\otimes \bar e_{\bar j}- \sum_{j}z_{ij}e_j\otimes \bar e_{\bar n+\bar i}
-\sum_{\bar j}\bar z_{\bar i\bar j}e_{n+i}\otimes \bar e_{\bar j}+ e_{n+i}\otimes \bar
e_{\bar n+\bar i}$$ we get the explicit form of this linear combination: it is sufficient
to prove that for all $i$, $\bar i$ we have
$$\sum_{j,\bar j}z_{ij}\bar z_{\bar i\bar j} (\Psi^-_{M\otimes \bar M})_{(j,\bar j)} -
\sum_{j}z_{ij}(\Psi^-_{M\otimes \bar M})_{(j, \bar n+\bar i)} $$ $$ -\sum_{\bar j}\bar
z_{\bar i\bar j}(\Psi^-_{M\otimes \bar M})_{(n+i,\bar j)} +(\Psi^-_{M\otimes \bar
M})_{(n+i,\bar n+\bar i)} =0\eqno{(6.3)}$$
Further, 6.2 implies that $$(\Psi^-_{M\otimes \bar M})_{(k,\bar
k)}=\frac{A_{-1,k}\otimes \bar A_{-1,\bar k}}{(T- \theta)^2}+\frac{A_{-1,k}\otimes \bar
A_{0,\bar k}+A_{0,k}\otimes \bar A_{-1,\bar k}}{T- \theta}$$ hence 6.3 becomes
$$\sum_{j,\bar j}z_{ij}\bar z_{\bar i\bar j} (\frac{A_{-1,j}\otimes \bar A_{-1,\bar
j}}{(T- \theta)^2}+\frac{A_{-1,j}\otimes \bar A_{0,\bar j}+A_{0,j}\otimes \bar A_{-1,\bar
j}}{T- \theta}) $$ $$- \sum_{j}z_{ij}(\frac{A_{-1,j}\otimes \bar A_{-1,\bar n+\bar i}}{(T-
\theta)^2}+\frac{A_{-1,j}\otimes \bar A_{0,\bar n+\bar i}+A_{0,j}\otimes \bar A_{-1,\bar
n+\bar i}}{T- \theta}) $$ $$-\sum_{\bar j}\bar z_{\bar i\bar j}(\frac{A_{-1,n+i}\otimes
\bar A_{-1,\bar j}}{(T- \theta)^2}+\frac{A_{-1,n+i}\otimes \bar A_{0,\bar
j}+A_{0,n+i}\otimes \bar A_{-1,\bar j}}{T- \theta}) $$ $$+\frac{A_{-1,n+i}\otimes \bar
A_{-1,\bar n+\bar i}}{(T- \theta)^2}+\frac{A_{-1,n+i}\otimes \bar A_{0,\bar n+\bar
i}+A_{0,n+i}\otimes \bar A_{-1,\bar n+\bar i}}{T- \theta}=0\eqno{(6.4)}$$
It is easy to see that 6.4 follows immediately from the equalities
$$A_{-1,n+i}=\sum_j z_{ij}A_{-1,j}\eqno{(6.5)}$$ $$ \bar A_{-1,\bar n+\bar
i}=\sum_{\bar j} \bar z_{\bar i\bar j}\bar A_{-1,\bar j}$$ For example, the left hand
side of (6.4) has 2 terms containing $\bar A_{0,\bar j}$ (in the middle of the first
and the third lines of (6.4)). Multiplying (6.5) by $\bar z_{\bar i\bar j}\bar
A_{0,\bar j}$ we get that the sum of these 2 terms of (6.4) is 0. For other pairs of
terms of (6.4) the situation is the same.
\medskip
The proof that for $i=1,...,r-n$, $\bar i=1,...,\bar n$ the element $b_{n+i}\otimes \bar
b_{\bar i}\in \goth q(M\otimes \bar M)$ is analogous but simpler. We have
$$b_{n+i}\otimes \bar b_{\bar i}=(T-\theta) (-\sum_{j} z_{ij} e_j\otimes \bar e_{\bar i}
+ e_{n+i}\otimes \bar e_{\bar i})$$ The analog of (6.3)) is $$(T-\theta) (-
\sum_{j}z_{ij}(\Psi^-_{M\otimes \bar M})_{(j, \bar i)}  +(\Psi^-_{M\otimes \bar
M})_{(n+i,\bar i)}) =0$$ and the analog of (6.4)) is $$-
\sum_{j}z_{ij}\frac{A_{-1,j}\otimes \bar A_{-1,\bar i}}{T- \theta}
+\frac{A_{-1,n+i}\otimes \bar A_{-1,\bar i}}{T- \theta}=0$$ This equality follows
immediately from (6.5).
\medskip
Finally, elements $b_{i}\otimes \bar b_{\bar i}$ ($i=1,...,n$, $\bar i=1,...,\bar n$)
obviously belong to $ \goth q(M\otimes \bar M)$.
\medskip
So, we proved that $\goth q(M)\underset{\p[[T-\theta]]}\to{\otimes}\goth q(\bar M)
\subset \goth q(M\otimes \bar M)$. Since the $\p$-codimension of both subspaces in
$L_T(M)\underset{\w}\to{\otimes} L_T(\bar M)\underset{\w}\to{\otimes}
\p[[T-\theta]]$ is $n\bar n$, they are equal. $\square$
\medskip
{\bf 7. Self-dual t-motives.}
\medskip
{\bf Case $\w=\n F_q[T]$.} A uniformizable t-motive $M$ is called self-dual if there exists an isogeny $\alpha: M \to M'$. It defines an $\w$-valued, $\w$-bilinear form $<*,*>_\alpha$ on $L_T(M')$ as follows:
$$<\vf_1, \vf_2>_\alpha=<L_T(\alpha)(\vf_1), \vf_2>_M$$
5.2.11 implies that if $\alpha'=-\alpha$ (resp. $\alpha'=\alpha$) then
$<*,*>_\alpha$ is skew symmetric (resp. symmetric). $M$ is called positively (resp. negatively) self-dual if $\alpha$ satisfies $\alpha'=\alpha$ (resp. $\alpha'=-\alpha$). Hence, we have an
\medskip
{\bf Analogy 7a.} The number field case analog of a pair: $\{$negatively self-dual t-motive of rank $2n$, dimension $n$; negative $\alpha: M \to M'\}$ is a (generic) abelian variety of dimension $n$ with a fixed polarization form.
\medskip
For example, like in the number field case, we can define the Rosati involution $I_\alpha$ on  $\End_0(M):=\End(M)\otimes \n F_q(T)$ by the same formula $I_\alpha(\vf)=\alpha^{-1}\circ \vf'\circ \alpha$.
\medskip
Further, we have a
\medskip
{\bf Conjecture 7b.} The dimension of the moduly variety of negatively self-dual t-motives (if it exists) is $n(n+1)/2$.
\medskip
{\bf Examples.} Let $e_*$ be from 1.9, and let $M=M(A)$ given by the equation (here $A\in M_n(\p)$ is $\goth A_1$ of 1.9.1)
$$Te_*= \theta e_*+ A\tau e_* + \tau^2 e_*\eqno{(7.1)}$$
be a t-motive of dimension $n$ and rank $2n$. Elements $f_i=e_i$,
$f_{n+i}=\tau e_i$ $(i=1,...,n)$ form a $\p[T]$-basis of $M$. We have
(see, for example, Section 11): $M'$ is given by the equation $$Te'_*=
\theta e'_*- A^t\tau e'_* + \tau^2 e'_*$$ and if we define $$f'_i=\tau
e'_i, \ \ f'_{n+i}= e'_i\eqno{(7.2)}$$ then bases $f_*$, $f'_*$ are dual in the
meaning of Lemma 1.10.

Let $\alpha: M \to M'$ be given by the formula $\alpha(e_*)=De'_*$
where $D\in M_n(\p)$ (we impose this essential restriction only in order to simplify exposition. In the general case $D\in M_n(\p[\tau])$, $D_f\in M_{2n}(\p[T])$, $D_f$ from 7.4). Condition that $\alpha$ is a $\p[T,\tau]$-map is
equivalent to
$$D^{(2)}=D, \ \ AD^{(1)}=-DA^t\eqno{(7.3)}$$
Further, we have $$\alpha(f_*)=D_ff'_*\eqno{(7.4)}$$ where
$D_f=\left(\matrix 0 & D\\ D^{(1)}&0 \endmatrix \right)$, hence
$$\alpha'=\pm \alpha\iff
D_f^t=\pm D_f\iff D^{(1)}=\pm D^t\eqno{(7.5)}$$
Let us fix $\varepsilon_0\in \n F_{q^2}$ satisfying $\varepsilon_0^{q-1}=-1$. Then $D=\varepsilon_0E_n$ satisfies 7.5 with the sign minus, and the set of $A$ satisfying 7.3 with this $D$ is the set of symmetric matrices. This justifies 7b, because the set of $A_1\in  M_n(\p)$ such that $M(A)=M(A_1)$ is conjecturally discrete.

For $D=E_n$ the sign in 7.5 is plus and hence a skew symmetric $A$ defines a positively self-dual $M(A)$.
\medskip
{\bf Remark 7.6.} The below statements are conjectures based on arguments similar to the ones which justify the below Conjecture 11.8.3. Since they are of secondary importance, we do not give any details of justification here.
\medskip
{\bf 7.6.1. Conjecture.} If $n\ge 3$ then for a
generic skew symmetric $A$ we have: $\End(M(A))=\w$.
\medskip
{\bf 7.6.2. Corollary.} Conjecture 7.6.1 implies that the "minimal" $\alpha: M \to M'$ is defined uniquely up to an element of
$\n F_q^*$, and hence the symmetric pairing $<*,*>_\alpha$ is also defined uniquely up to an element of $\n F_q^*$.
\medskip
{\bf 7.6.3. Conjecture.} If $n=2$, $\alpha'=\alpha$ then $\End (M)$ is strictly larger than $\w$.
\medskip
Other examples of a self-dual t-motive are $M\oplus M'$ where $M$ is any t-motive, but they do not give interesting examples of pairings.
\medskip
{\bf 7.6.4. Conjecture.} There exist other (distinct from the ones defined by 7.1) self-dual t-motives $M$ having $\End (M)=\w$ (we can use a version of standard t-motives of Section 11).
\medskip
{\bf Example 7.7.} Case $A=0$, $D=E_n$.
\medskip
In this case we can find explicitly the matrix of the symmetric form
$<*,*>_\alpha$ in some basis of $L_T(M')$. Let $\goth C_2$ be the
Carlitz module over the field $\n F_{q^2}$ considered as a rank 2
Drinfeld module over $\n F_{q}$ given by the equation $$Te=\theta e
+ \tau^2e$$ We have $M=\goth C_2^{\oplus n}$. Let $\goth T_T(\goth
C_2)$ be the convergent $T$-Tate module of $\goth C_2$, i.e. the set
of elements $\{z_i\}\in E(\goth C_2)=\p$ $(i\ge -1, \ \ z_{-1}=0)$ such that
$$\hbox{ $T z_i=z_{i-1}$ for $i \ge 0$ (i.e. $z_i^{q^2}+\theta
z_i=z_{i-1}$) and $z_i\to 0$}$$ It is a free 1-dimensional module over
$\n F_{q^2}[T]$. We choose and fix its generator; its $\{z_i\}$ satisfy (like in 5.3.4)
$|z_0|>|z_i| \ \ \forall i>0$. We denote $\sum_{k=0}^\infty z_kT^k$ by $\goth Z$.

Let $c$ be a fixed element of $\n F_{q^2}-\n F_{q}$. Formulas (5.3.3)
show that the following elements $\vf_i$, $\vf'_i$ ($i=1, ... , 2n$)
form bases of $L(M)$, $L(M')$ respectively ($j=1, ..., n$; clearly that thanks to 7.2
we have $\vf'_i(f'_j)=\vf_i(f_{n+j})$, $n+j$ mod $2n$):
$$i\le n: \ \ \vf_i(f_j)=\goth Z\delta_i^j, \ \
\vf_i(f_{n+j})=\goth Z^{(1)}\delta_i^j $$

$$i>n: \ \ \vf_i(f_j)=c\goth Z\delta_{i-n}^j, \ \
\vf_i(f_{n+j})=c^q\goth Z^{(1)}\delta_{i-n}^j $$

$$i\le n: \ \ \vf'_i(f'_j)=\goth Z^{(1)}\delta_i^j, \ \
\vf'_i(f'_{n+j})=\goth Z\delta_i^j $$

$$i>n: \ \ \vf'_i(f'_j)=c^q\goth Z^{(1)}\delta_{i-n}^j, \
\ \vf'_i(f'_{n+j})=c\goth Z\delta_{i-n}^j $$
(by the way, it is clear that the same relation between elements of $\goth T_T(M)$ and
$\Hom_{\p[T,\tau]}(M,Z_1)$ holds for all $M$). Formula 7.4 shows
that $\alpha'(\vf'_i)=\vf_{i+n}$, where $i+n \mod 2n$. Let us denote
$\Xi\cdot \goth Z\cdot \goth Z^{(1)}\in\n F_q^*$
by $\gamma$. The above definitions and formulas show that the matrix
of $<*,*>_\alpha$ in the basis $\vf_{1}, \vf_{n+1}, ... ,\vf_{n},
\vf_{2n}$ consists of $n$ \ \ $(2\times 2)$-blocks (trace and norm of $\n
F_{q^2}/\n F_q$)
$$\gamma\left(\matrix \tr(1) & \tr(c)\\ \tr(c)&\tr(N(c))  \endmatrix
\right)=\gamma\left(\matrix 2 & c+c^q \\ c+c^q & 2c^{q+1} \endmatrix
\right)$$ The determinant of
this block is $-(c-c^q)^2\gamma^2$; it belongs to $\n F_q^{*2}\iff
q\equiv 3 \mod 4$ or $q$ is even. Since we have $n$ blocks, we have:

$$\hbox{det $<*,*>_\alpha\not\in\n F_q^{*2}\iff q\equiv 1 \mod 4$ and
$n$ is odd.}$$
\medskip
{\bf Remark 7.8 (Jorge Morales).}
There is a theorem of Harder (see
e.g. W. Scharlau, "Quadratic and Hermitian forms", Springer-Verlag,
Berlin, 1985, Chapter 6, Theorem 3.3) that states that a unimodular form over
$k[X]$ \ \ --- \ \ $k$ being any field of characteristic not 2 --- is the
extension of a form over $k$, i.e. there is a basis in which all the
entries of the associated symmetric matrix are constant. This means
that the classification of the above quadratic forms over $\n F_q[T]$
($q$ odd) is very simple.
\medskip
{\bf Remark 7.9.} Let $M$ be a t-motive which is both negatively and positively self-dual. There is a natural idea 7.9.2 to define an analog of Hodge structure on $M$. Nevertheless, this idea fails. Namely, the exact sequence
$$0\to \Ker \vf \to L(M)\otimes \p \overset{\vf}\to{\to} \Lie(M)\to 0$$
is the functional field analog of an exact sequence for an abelian variety $A$:
$$ 0 \to H^{0,-1}(A) \to H^{-1}(A) \to (H^{1,0})^*(A) \to 0 $$
Hence, we can define $H^{0,-1}(M):= \Ker \vf$, and the problem is to define an analog of $H^{-1,0}(M)$.

Let us fix a negative isogeny $\alpha: M\to M'$, and let us extend the skew form $<*,*>_\alpha$ to $L(M)\otimes \p$ by $\p$-linearity. It is easy to check that $\Ker \vf$ is isotropic with respect to this form (there is an analogy with the number field case). Let us consider the following elementary lemma of linear algebra:
\medskip
{\bf Lemma 7.9.1.} Let $W$ be a vector space of dimension $2n$ over a field of characteristic $\ne 2$, $B^+$ (resp. $B^-$) a symmetric (resp. skew symmetric) non-degenerate bilinear form on $W$, and $W_0\subset W$ a subspace of dimension $n$ which is isotropic with respect to both $B^+$, $B^-$. Then almost always there exists the only $W_1\subset W$ of dimension $n$ having properties:
$$W_0\cap W_1=0; \ \ \ W_1 \hbox{  is isotropic with respect to both } B^+, \  B^-$$
where almost always means that entries of the matrices of $B^+$, $B^-$ in a basis of $W$ must not satisfy (at least one of) polynomial relations. $\square$
\medskip
If $\End_0(M)\ne \n F_q(T)$ and the action of $I_\alpha$ on $\End_0(M)$ is not identical, then there exists a positive isogeny $\beta: M\to M'$ and hence the symmetric form $<*,*>_\beta$ on $L(M)\otimes \p$. $\Ker \vf$ is isotropic with respect to  $<*,*>_\beta$. Let us fix $\beta$.
\medskip
{\bf Idea 7.9.2.} To apply Lemma 7.9.1 to this situation ($W=L(M)\otimes \p$, $W_0=\Ker \vf$, $B^+=<*,*>_\beta$, $B^-=<*,*>_\alpha$) in order to get a canonical subspace of $L(M)\otimes \p$ which is complementary to $\Ker \vf$ and hence can be considered as an analog of $H^{-1,0}(M)$.
\medskip
Clearly there is no complete analogy with the number field case. But the situation is even worse:
\medskip
{\bf Proposition 7.9.3.} For all $M$, $\alpha$, $\beta$ the "almost always" condition of Lemma 7.9.1 is not satisfied. $\square$
\medskip
{\bf 8. Relations between lattices and t-motives.}
\medskip
We have\footnotemark \footnotetext{I am grateful to Urs
Hartl
who indicated me this reference.}
\medskip
{\bf Theorem 8.1.} ([H],  Theorem 3.2). The dimension of the moduli set of pure t-motives of
dimension $n$ and rank $r$ is $n(r-n)$. $\square$
\medskip
{\bf Remark.} A tuple
$(e_1,...,e_r)$ of
integers entering in the statement of this theorem in [H] is
(0,...,0,1,...,1) with 0 repeated
$r-n$ times and 1 repeated $n$ times for the case under
consideration.
\medskip
Since this number $n(r-n)$ is equal to the dimension of the set of lattices of rank $r$
and dimension $n$, we can state an
\medskip
{\bf Open question 8.2.} Let $r$, $n$ be given. Let us consider the lattice map from the set of the pure uniformizable
t-motives of rank $r$ and dimension $n$ to the set
of lattices of rank $r$ and dimension $n$. Is it true that its image is open and the fibre at a generic point is discrete? If yes, what is the fibre?
\medskip
{\bf Remark.} Results of [L3] give some evidence that for the case $r=2n$ in a "neighborhood" of the $n$-th power of the rank 2 Carlitz module the fibre consists of 1 point.
\medskip
Theorem 5 implies that for $n=r-1$ the answer to 8.2 is yes (the below Proposition 11.8.5 shows that most likely the condition of purity is essential):
\medskip
{\bf Corollary 8.4.} All pure t-motives of
dimension $r-1$ and rank $r$ having $N=0$ are uniformizable. There is a 1 -- 1 functorial
correspondence between pure t-motives of dimension $r-1$ and rank $r$ having $N=0$ ($r\ge
2$), and lattices of rank $r$ in $\p^{r-1}$ having dual.
\medskip
{\bf Proof.} Let $L$ be a lattice of rank $r$ in $\p^{r-1}$ having dual $L'$. There
exists the only Drinfeld module $M'$ such that $L(M')=L'$, and let $M$ be its dual.
Theorem 5 implies that $L(M)=L$. If there exists another pure t-motive $M_1$ of
dimension $r-1$ and rank $r$ having $N=0$ such that $L(M_1)=L$ then by Corollary 10.4 (its proof is logically independent: there is no vicious circle) the
dual $M'_1$ is a Drinfeld module, according Theorem 5 it satisfies $L(M'_1)=L'$, hence
$M'_1=M'$ and hence $M_1=M$. $\square$
\medskip
{\bf Remark 8.5.} Recall that lattices of rank $r$ in $\p^{r-1}$ having dual are
described in 3.5 (formulas 3.6, 3.7).
We see that for the case $n=r-1$, $N=0$ purity implies uniformizability. We have
\medskip
{\bf Question 8.5a.} Do exist non-uniformizable t-motives having $n=r-1$, $N=0$?
\medskip
{\bf Question 8.5b.} Do exist uniformizable t-motives having $n=r-1$, $N=0$ such that its lattice has no dual? (Clearly this is a subquestion of 8.2).
\medskip
{\bf Remark 8.6.} Clearly for any $r$, $n$ we have: if a lattice $L$ of rank $r$ and dimension $n$ has no dual then $L\ne L(M)$ for any pure uniformizable $M$. I do not know whether Theorem 6 (which is an analog of Theorem 5 for another tensor operation) imposes a more strong similar restriction on the property of $L$ to be the $L(M)$ of some pure uniformizable $M$, or not.
\medskip
Further, for any uniformizable t-motive $M$ we have a
\medskip
{\bf Corollary 8.7.} If the dual of $(L(M), \Lie(M))$ does not exist then the
dual of $M$ does not exist. Example: the Carlitz module.
\medskip

{\bf 9. Main theorem in terms of Pink-Hodge structure.}
\medskip
Let us consider a version of a special case of the general definition of  Pink-Hodge structure ([P], 0.2; 9.1).
\medskip
{\bf Definition.} A Pink-Hodge structure of constant weight and complete dimension is a pair $\underline{H}=(H,\goth q_H)$ where $H$ is a free finite dimensional $\w$-module and $\goth q_H$ is a $\p[[T-\theta]]$-lattice in $H\underset{\w}\to{\otimes}\p[[T-\theta]]$ such that the dimension of $\goth q_H$ over $\p[[T-\theta]]$ is equal to the dimension of $H$ over $\w$ (condition of complete dimension).
\medskip
Let $\vf: L \hookrightarrow \p^n$ be a lattice. It defines a Pink-Hodge structure $\underline{H}=\underline{H}(L)$ of constant weight and complete dimension. Firstly, instead of a $\n F_q[\theta]$-module $L$ we consider an isomorphic $\w$-module $H$ formally defined by the property $H\underset{\w}\to{\otimes}\n F_q[\theta]=L$ where the map $\w \to \n F_q[\theta]$ is $\iota$. We denote the isomorphism $H \to L$ by $\iota$ as well; the composition $\vf \circ \iota: H \to \p^n$ is a map of $\w$-modules where $T\in \w$ acts on $ \p^n$ by multipication by $\theta$. Further, $\vf \circ \iota$ extends to a surjection of $\p[[T-\theta]]$-modules $H\underset{\w}\to{\otimes}\p[[T-\theta]] \to \p^n$ denoted by $\vf \circ \iota$ as well.  Finally, $\goth q_H$ is defined as $\Ker \vf \circ \iota$.

If $M$ is a pure uniformizable t-motive then we associate it a Pink-Hodge structure of constant weight and complete dimension $\underline{H}(M)=\underline{H}(L(M))$.

Let $m=m(\underline{H})$ be the minimal number such that $\goth q_H \supset (T-\theta)^m H\underset{\w}\to{\otimes}\p[[T-\theta]]$. For $\mu\ge m$ we define the $\mu$-dual structure ${\underline{H}'}^{\mu}=({H'}^{\mu}, \goth q_{{H'}^{\mu}})$ as follows:

$${H'}^{\mu}=H^*, \ \  \goth q_{{H'}^{\mu}}=\{\chi\in H^*\underset{\w}\to{\otimes}\p[[T-\theta]]  $$ $$\hbox{ such that }  \forall y\in \goth q_H  \hbox{ we have } \chi(y)\in (T-\theta)^\mu\p[[T-\theta]]\} $$ It is obvious that it is really a Pink-Hodge structure of constant weight and complete dimension.
\medskip
If $\underline{H}=\underline{H}(L)$ for a lattice $L$ then $m=1$ and if $L$ has dual then $${\underline{H}'}^1=\underline{H}(L')\eqno{(9.1)}$$ this is easy to prove.
\medskip
{\bf Remark 9.2.} And if $L$ has no dual? Really, $\underline{H}(L)$ exists even if $L$ does not satisfy Definition 2.1 (b). If $L$ is a lattice having no dual this means that $L'$  does not satisfy Definition 2.1 (b). Nevertheless, equality ${\underline{H}'}^1=\underline{H}(L')$ is meaningful and holds. We are not interested in these lattices because they cannot be lattices of uniformizable t-motives having dual.
\medskip
An analog of Theorem 6 for dual t-motives is the following. Let $M$ be as above. Obviously $m=m(\underline{H}(M))$ is the minimal number such that $N^m=0$. According Theorem 10.3, ${M'}^{m}$ exists, and it is pure.
\medskip

{\bf 10. Duals of pures, and other elementary results. }

\nopagebreak
\medskip
We consider in this section the case of arbitrary $N$ (i.e. not necessarily $N=0$), and $\w=\n F_q[T]$. The definition 1.8 extends to the case of pr\'e-t-motives, and
remarks 1.11 hold for this case.
\medskip
{\bf Lemma 10.2. } Let $M$ be a pr\'e-t-motive, $m=m(M)$ from its
(1.3.1), and $\mu\ge m$.
Then $M'$ --- the $\mu$-dual of $M$ --- exists as a pr\'e-t-motive,
and $m(M')\le\mu$. If
$M'$ is a t-motive then $\dim M'=r\mu - \dim M$ ($r$ is the rank of $M$).
\medskip
{\bf Proof.} We must check that $Q'$ has no denominators, and the
condition (1.3.1). The
module $\tau M$ is a $\p[T]$-submodule of $M$ (because $a \tau x =
\tau a^{1/q} x$ for
$x\in M$),
hence there are $\p[T]$-bases $f_*=(f_1, ... f_r)^t$, $g_*=(g_1, ...
g_r)^t$ of $M$,
$\tau M$ respectively such that $g_i=P_if_i$, where $P_1 | P_2 | ... |
P_r$, $P_i \in
\p[T]$. Condition (1.3.1) means that $\forall i$ \ $(T-\theta)^m f_i
\in \tau M$, i.e.
$P_i|(T-\theta)^m$, i.e. $\forall i$ \ $P_i=(T-\theta)^{m_i}$ where $0
\le m_i \le
m_{i+1} \le m$. There exists a matrix $\goth Q=\{\goth q_{ij}\}\in
M_r(\p[T])$ such that
$$\tau f_i = \sum_{j=1}^r \goth q_{ij}g_j= \sum_{j=1}^r \goth q_{ij} P_j
f_j\eqno{(10.2.1)}$$
Although $\tau$ is not a linear operator, it is easy to see that
$\goth Q\in GL_r(\p[T])$
(really, there exists $C=\{c_{ij}\}\in M_r(\p[T])$ such that $g_i=P_i
f_i=\tau (\sum_{j=1}^r
c_{ij}f_j)$, we have $C^{(1)}\goth Q =E_r$).

We denote the matrix $\diag(P_1, P_2, ... ,P_r)$ by $\goth P$, so
(10.2.1) means that
$$Q= \goth Q \goth P \eqno{(10.2.2)}$$

{\bf Remark 10.2.3.} Since
$\goth Q \goth P \in
GL_r(\p(T))$, we get that the action of $\tau$ on $i_2(M)$ is invertible.
\medskip
It is clear that if $M$ is a t-motive then
$$\dim M = \sum_{j=1}^r m_j\eqno{(10.2.4)}$$
(because $\dim M = \dim_{\p}(M/\tau M)$. Further, (10.2.2) implies that
for $Q'=Q({M'})$
we have
$$Q'=\goth Q^{t-1}\diag((T-\theta)^{\mu-m_1}, ...,
(T-\theta)^{\mu-m_r})\eqno{(10.2.5)}$$
This means that elements of $Q'$ have no denominators. The condition
(1.3.1) for $M'$
follows easily from (10.2.5) (because $\goth Q^{t-1}\in GL_r(\p[T])$),
and the dimension
formula (for the case $M'$ is a t-motive) follows immediately from
(10.2.4) applied to
$M'$. $\square$
\medskip
A definition of a pure t-motive can be found in [G] ((5.5.2),
(5.5.6) of [G] + formula (1.3.1) of the present paper).
\medskip
{\bf Theorem 10.3.} Let $M$ be a pure t-motive and $m=m(M)$ from
(1.3.1). Then (if
$rm-n>0$) its $m$-dual $M'$ exists, and it is pure.
\medskip
{\bf Proof.} The definition of pure ([G], (5.5.2)) is valid for
pr\'e-t-motives. We use
its following
matrix form. We denote $T^{-1}$ by $S$ and for any $C$ we let
$$C^{[i]}=C^{(i-1)}\cdot C^{(i-2)}\cdot ... \cdot C^{(1)}\cdot C$$

{\bf Lemma 10.3.1.} Let $Q\in M_r(\p[T])$ be a matrix such that formula
(1.9.3) defines an
t-motive $M$. Then it is pure iff there exists $C\in
GL_r(\p((S)) \ )$ such that
for some $\goth q$, $s>0$
$$S^\goth q C^{(s)}Q^{[s]}C^{-1}\in GL_r(\p[[S]])$$
i.e. iff $S^\goth q C^{(s)}Q^{[s]}C^{-1}$ is $S$-integer and its inicial
coefficient is
invertible.
\medskip
{\bf Proof.} Elementary matrix calculations. We take $C$ as a matrix
of base change of
$f_*$ to a $\p[[S]]$-basis of $W$ of (5.5.2) of [G]. $\square$
\medskip
{\bf Lemma 10.3.2.} Let $\mu = m$. We have: $M'={M'}^\mu$ of Lemma 10.2 is a pure
pr\'e-t-motive.
\medskip
{\bf Proof.} Let $\goth q$, $s$ and $C$ be from Lemma 10.3.1. We have
$${Q'}^{[s]}=((T-\theta)^{[s]})^\mu Q^{[s]\ t-1}$$
(we use (1.2)). We take $C'=C^{t-1}$. We have
$$S^{s\mu-\goth q} {C'}^{(s)}{Q'}^{[s]}{C'}^{-1}=$$
$$=S^{s\mu-\goth q}C^{(s)\ t-1}Q^{[s]\ t-1}((\frac1S-\theta)^{[s]})^\mu C^t$$
$$=((1-S\theta)^{[s]})^\mu S^{-\goth q}C^{(s)\ t-1}Q^{[s]\ t-1}C^t$$
$$=((1-S\theta)^{[s]})^\mu (S^\goth q C^{(s)}Q^{[s]}C^{-1})^{t-1}$$
We have: $\goth q/s=n/r$ ([G], (5.5.6)), hence $(s\mu-\goth q)/s=(r\mu-n)/r$ and
$s\mu-\goth q>0$. Further,
$((1-S\theta)^{[s]})^\mu\in GL_r(\p[[S]])$, and the result follows from
Lemma 10.3.1. $\square$
\medskip
{\bf Remark.} This result holds also for $\mu>m$.
\medskip
The theorem 10.3 follows from Lemma 10.2, the above lemmas and the
proposition that a pure
pr\'e-t-motive satisfying (1.3.1) is a t-motive ([G], (5.5.6),
(5.5.7)). $\square$
\medskip
{\bf Corollary 10.4.} Let $M$ be a t-motive such that $m=1$,
$n=r-1$. Then $M$
has dual $\iff$ $M$ is pure $\iff$ $M$ is dual to a Drinfeld module.
\medskip
{\bf Proof.} Dimension formula shows that $M'$ (if it exists) is a
Drinfeld module, and they are all pure. $\square$
\medskip
{\bf Example 10.5.} Let $M$ be given by (notations of 1.9.1)
$$\goth A_0=\theta E_2, \ \ \goth A_{1}= \left(\matrix a_{111} & 0
\\ a_{121} & 1 \endmatrix \right), \ \ \goth A_2= \left(\matrix 1 & 0
\\0 & 0 \endmatrix \right)$$
This $M$ has $m=1$, $n=2$, $r=3$, and it is easy to see that it has no dual. Really, for
this $M$ we have (notations of 1.9) $f_1=e_1$, $f_2=\tau e_1$, $f_3=e_2$,
$$Q=\left(\matrix 0&1&0\\T-\theta&-a_{111} & 0
\\ 0&-a_{121} & t-\theta \endmatrix \right), \ \ Q'=\left(\matrix a_{111} & t-\theta
& a_{121} \\ 1 &0&0\\0&0&1\endmatrix \right)$$ The last line of $Q'$ means that $\tau
f_3'=f_3'$. This is a contradiction to the property that $M'_{\p[\tau]}$ is free. It is
possible also to show (Proposition 11.3.4) that $M$ is not pure, and to use 10.4 in order to
prove that it has no dual.
\medskip
Later (Section 11) we shall construct examples of non-pure abelian
t-motives which have dual. Considerations of 11.8 predict that there is enough such t-motives.
\medskip
{\bf Theorem 10.6. } For any t-motive $M$ there exists $\mu_0$ such that for all
$\mu\ge\mu_0$ the object
${M'}^{\mu}$ exists as a t-motive. For these $\mu$ we have
$${M'}^{\mu+1}={M'}^{\mu}\otimes \goth
C\eqno{(10.6.1)}$$

{\bf Proof.} (10.6.1) holds at the level of pr\'e-t-motives, because $Q(\goth
C)=(T-\theta)E_1$. According [G], Lemma 5.4.10 it is sufficient to prove that
${M'}^{\mu}$ is finitely generated as a $\p [\tau]$-module. We shall
use notations of
Lemma 10.2. We take $$\mu_0=1+\hbox{ \{the maximum of the degrees of
entries of $\goth
Q(M)$ as polynomials in $T$\} } $$ $$+ \max(m_k)$$
Let $f'_1, ...f'_r$ be the basis of ${M'}^{\mu}$ over $\p [T]$ dual to
$f_1, ...f_r$. It
is sufficient to prove the
\medskip
{\bf Lemma 10.6.2.} Let $i_0=\mu-\min(m_k)$. Then elements $T^if'_j$,
$i<i_0$, $j=1,
....,r$, generate ${M'}^{\mu}$ as a $\p [\tau]$-module.
\medskip
{\bf Proof of the lemma.} By induction, it is sufficient to show that
for all $\alpha\ge
i_0$ the equation
$$\tau x = (T-\theta)^\alpha f'_j\eqno{(10.6.3)}$$
(equality in ${M'}^{\mu}$) has a solution
$$x=\sum_{k=1}^r C_kf'_k$$
where $C_k \in \p [T]$, $\deg(C_k)< \alpha$. According (10.2.5), the
solution to (10.6.3)
is given by
$$(C_1^{(1)}, ... , C_r^{(1)})=(0,... 0,
(T-\theta)^{\alpha-\mu+m_j},0,... 0)\goth Q^t$$
(the non-0 element of the row matrix is at the $j$-th place).
Unequalities satisfied by
$\mu$ and $\alpha$ show that all $C_k^{(1)}$ are polynomials of degree
$<\alpha$. Since
$c \mapsto c^q$ is surjective on $\p$, we get the desired. $\square$
\medskip
{\bf 10.7. Virtual t-motives.} \footnotemark \footnotetext{This
notion was
indicated me by
Taguchi.} We need two elementary lemmas.

\nopagebreak
\medskip
{\bf Lemma 10.7.0.}\footnotemark \footnotetext{Anderson proved (not published) that the tensor product of any t-motives is also a t-motive.} If $M$ is a t-motive then $M\otimes \goth C$ is also a t-motive.
\medskip
{\bf Proof.} Let $f_{j}$ ($j=1,...,r$) be a $\p[T]$-basis of
$M_{\p[T]}$ and $\goth f$
from 1.10.2, so $f_{j}\otimes \goth f$ is a $\p[T]$-basis of $(M\otimes \goth
C)_{\p[T]}$. It is sufficient to prove that $(M\otimes \goth
C)_{\p[\tau]}$ is finitely
generated. Since $M_{\p[\tau]}$ is finitely generated, it is easy to
see that there
exists $a$ such that elements
$$(T-\theta)^i f_j, \ \ i=0, ... , a, \ \ j=1, ... ,r$$
generate $M_{\p[\tau]}$. This means that $\forall j=1, ... ,r$ there
exist $c_{ijkl}\in
\p$ such that
$$(T-\theta)^{a+1} f_j=\sum_{i=0}^a \sum_{k=0}^\gamma\sum_{l=1}^r
c_{ijkl}(T-\theta)^i
\tau^k f_l\eqno{(10.7.0.1)}$$
where $\gamma$ is a number.

Let us multiply (10.7.0.1) by $(T-\theta)^\gamma$. Taking into
consideration the formula
of the action of $\tau$ on $M\otimes \goth C$ we get that the result
gives us the
following formula in $M\otimes \goth C$:
$$(T-\theta)^{a+\gamma+1} f_j\otimes \goth f=\sum_{i=0}^a
\sum_{k=0}^\gamma\sum_{l=1}^r
c_{ijkl} (T-\theta)^{i+\gamma-k} \tau^k \cdot (f_l\otimes \goth f)\eqno{(10.7.0.2)}$$
This proves that for all $j$ the element $(T-\theta)^{a+\gamma+1} f_j\otimes \goth f$ is
a linear combination of
$$(T-\theta)^i f_l\otimes \goth f, \ \ i=0, ... , a+\gamma, \ \ l=1, ...
,r\eqno{(10.7.0.3)}$$
in $(M\otimes \goth C)_{\p[\tau]}$.
Multiplying (10.7.0.2) by consecutive powers of $T-\theta$ we get by
induction that
elements of 10.7.0.3 generate $(M\otimes \goth C)_{\p[\tau]}$. $\square$
\medskip
{\bf Lemma 10.7.1.} If $M_1\otimes \goth C$ is isomorphic to
$M_2\otimes \goth C$ then
$M_1$ is isomorphic to $M_2$.
\medskip
{\bf Proof.} Let $f_{i*}$ ($i=1,2$) be a $\p[T]$-basis of
$(M_i)_{\p[T]}$, $Q_i$ from
1.9.3, $\alpha: M_1\otimes \goth C\to M_2\otimes \goth C$ an
isomorphism and $C\in
GL_r(\p[T])$ the matrix of $\alpha$ in $f_{1*} \otimes \goth f $,
$f_{2*}\otimes \goth f
$. The matrix of the action of $\tau$ on $M_i\otimes \goth C$ in the
base $f_{i*} \otimes
\goth f $ is $(T-\theta)Q_i$, and the condition that $\alpha$ commutes with
multiplication by $\tau$ is
$$(T-\theta)Q_1 C = C^{(1)}(T-\theta)Q_2$$
Dividing this equality by $T-\theta$ we get that the map $\alpha_0$
from $M_1$ to $M_2$
having the same matrix $C$ in the bases $f_{i*}$, commutes with
$\tau$, i.e. defines an
isomorphism from $M_1$ to $M_2$. $\square$
\medskip
Using Lemma 10.7.1 we can state the following
\medskip
{\bf Definition.} A virtual t-motive is an object $M\otimes \goth C^{\otimes
\mu}$ where $M$ is
a t-motive
and $\mu\in \n Z$, with the standard equivalence relation (here
$\mu_1\ge\mu_2$):
$$M_1\otimes \goth C^{\otimes \mu_1}=M_2\otimes \goth C^{\otimes \mu_2}\iff
M_2=M_1\otimes \goth C^{\otimes (\mu_1-\mu_2)}$$
$$\iff \exists \mu \hbox{ such that $\mu+\mu_1\ge 0$, $\mu+\mu_2\ge 0$
and } M_1\otimes
\goth C^{\otimes (\mu+\mu_1)}=M_2\otimes \goth C^{\otimes (\mu+\mu_2)}$$
Lemma 10.7.1 shows that these conditions are really equivalent.
\medskip
{\bf Corollary 10.7.2.} The $\mu$-dual of a virtual t-motive is
well-defined and
always exists as a virtual t-motive. $\square$
\medskip
{\bf Proposition 10.8.} The following formula is valid at the level of
pr\'e-t-motives: for
any $\mu_1$, $\mu_2$, if ${M_1'}^{\mu_1}$, ${M_2'}^{\mu_2}$ exist then
${(M_1\otimes
M_2)'}^{(\mu_1+\mu_2)}$ exists and
$${(M_1\otimes M_2)'}^{(\mu_1+\mu_2)}={M_1'}^{\mu_1} \otimes {M_2'}^{\mu_2}$$
{\bf Proof.} This is a functorial equality; also we can check it by
means of elementary
matrix calculations. $\square$
\medskip
{\bf Proposition 10.9.} Let $P\in \hbox{\bf A}$ be an irreducible
element. The Tate
module
$T_P({M'}^{\mu})$ is equal to $$T_P(\goth C)^{\otimes \mu}\otimes
\widehat{T_P(M)}$$
(equality of Galois modules) where $\widehat{T_P(M)}$ is the dual
Galois module.
\medskip
{\bf Proof.} It is completely analogous to the proof of the corresponding theorem for
tensor products
([G], Proposition 5.7.3, p. 157).
All modules in the below proof will be the Galois modules, and equalities of modules will
be equalities of Galois
modules. Recall that $E=E(M)$. Since $T_P(M)=\invlim_n E_{P^n}$, it is sufficient to prove that for any
$a\in \w$ we have $E({M'}^{\mu})_a=E(\goth C^{\otimes \mu})_a\otimes \hat E_a$, where $\hat E_a$ is the dual of $E_a$ in the meaning of [T], Definition 4.1.
We have the following sequence of equalities of modules:
$${M'}^{\mu}/a{M'}^{\mu}=\Hom_{\p[T]}(M/aM, \goth C^{\otimes \mu}/a\goth C^{\otimes
\mu})\eqno{(10.9.2)}$$
such that the action of $\tau$ on both sides of this equality coincide (to define the
action of $\tau$ on the
right and side of (10.9.2) we need the action of $\tau^{-1}$ on $M/aM$; it is
well-defined, because the determinant
of the action of $\tau$ on $M$ is a power of $T-\theta$, hence its image in
$\p[T]/a\p[T]$ is invertible). 10.9.2
follows immediately from the definition of ${M'}^{\mu}$;
$$({M'}^{\mu}/a{M'}^{\mu})^\tau=\Hom_{\n F_q[T]}((M/aM)^\tau,
(\goth C^{\otimes \mu}/a\goth C^{\otimes \mu})^\tau)\eqno{(10.9.3)}$$
This follows from 10.9.2 and the Lang's theorem
$$\goth M/a\goth M=(\goth M/a\goth M)^\tau\underset{\n F_q[T]/a\n
F_q[T]}\to{\otimes}\p[T]/a\p[T]$$
applied to both $\goth M=M$, $\goth M={M'}^{\mu}$ (we use that both $M$, ${M'}^{\mu}$ are
free $\p[T]$-modules).
Finally, we have a formula
$$E(\goth M)_a=\Hom_{\n F_q}((\goth M/a\goth M)^\tau, \n F_q)$$
([G], p. 152, last line of the proof of Proposition 5.6.3). Applying this formula to 10.9.3 we get the desired. $\square$
\medskip
{\bf 11. An explicit formula.}

\nopagebreak
\medskip
We return to the case $N=0$. Let $e_*$, $\goth A$, $\goth A_i$, $l$, $n$ be from (1.9). We consider in the
present section two simple types of t-motives (called standard-1 and standard-2
t-motives respectively) whose $\goth A_i$ have a row echelon form, and we give an
explicit formula for the dual of some standard-1
t-motives. Analogous formula can be easily obtained for more general types of
t-motives. These results are the
first step of the problem of description of all t-motives
having duals.
\medskip
{\bf 11.1.} For the reader's convenience, we give here the definition of standard-1
t-motives
for the case $n=2$ (here $\lambda_1$ and $\lambda_2 $ satisfying
$\lambda_1=l$,
$l>\lambda_2 \ge 2$ are parameters):

$$\goth A_0=\theta E_2, \hbox{ for } 0 < i < \lambda_2 \ \ \goth A_i \hbox{ is
arbitrary, } $$
$$\goth A_{\lambda_2}= \left(\matrix * & 0 \\ * & 1 \endmatrix \right),
\hbox{ for } \lambda_2
< i < l \ \ \goth A_i= \left(\matrix * & 0 \\ * & 0 \endmatrix \right), \ \
\goth A_l= \left(\matrix
1 & 0
\\0 & 0 \endmatrix \right) $$

{\bf 11.2.} To define standard-2 t-motives of dimension $n$, we need to fix
\medskip
1. A permutation $\vf\in S_n$, i.e. a 1 -- 1 map $\vf: (1, ..., n) \to (1, ..., n)$;
\medskip
2. A function $k: (1, ..., n) \to \n Z^+$ where $\n Z^+$ is the set of integers $\ge 1$.
\medskip
{\bf Definition.} A standard-2 t-motive of the type $(\vf, k)$ is an abelian
t-motive of dimension $n$ given by the formulas ($i=1,...,n$):
$$Te_{\vf(i)}= \theta e_{\vf(i)}
+\sum_{\alpha=1}^n\sum_{j=1}^{k(\alpha)-1}a_{j,\vf(i),\alpha}\
\tau^je_\alpha + \tau^{k(i)}e_i\eqno{(11.2.1)}$$
where $a_{j,\vf(i),\alpha}\in\p$ is the $(\vf(i),\alpha)$-th entry of the matrix $\goth A_j$.
\medskip
{\bf Proposition 11.2.2.} Formula 11.2.1 really defines a t-motive denoted by
$M=M(\vf,k)=M(\vf,k,a_{***})$. Its rank is $\sum_{\alpha=1}^n k(\alpha)$ and elements
$X_{\alpha j}:=\tau^je_\alpha$, $\alpha=1, ..., n$, $j=0, ..., k(\alpha)-1$, form its
$\p[T]$-basis.
$\square$
\medskip
The group $S_n$ acts on the set of types $(\vf, k)$ and on the set of the above $M$;
clearly for any $\psi\in S_n$ we have $\psi(M)$ is isomorphic to $M$. Particularly, we
can consider only $\vf$ of the following form of the product of $i$ cycles ($\alpha_0=0,
\ \alpha_i=n$):
$$\vf=(\alpha_0+1, ..., \alpha_1)(\alpha_1+1,...,\alpha_2)
...(\alpha_{i-1}+1,...,\alpha_i)\eqno{(11.2.3)}$$
(standard notation of the theory of permutations, for $\gamma\ne \alpha_j$ we have
$\vf(\gamma)=\gamma+1$, for $\gamma= \alpha_j$ we have  $\vf(\alpha_j)=\alpha_{j-1}+1$).
\medskip
{\bf Example 11.2.4.} Let $\vf$ be defined by 11.2.3, the quantity of cycles $i$ is equal
to $1$ and all $a_{***}=0$. Then the corresponding $M$ is of complete multiplication by a
CM-field $\n F_{q^r}(T)$ and its CM-type $\Phi$ is $\{\Id, \fr^{k(1)}, \fr^{k(1)+k(2)},
..., \fr^{k(1)+k(2)+...+k(n-1)}\}$ where $\fr$ is the Frobenuis homomorphism $\n
F_{q^r}\to \bar \n F_q$ (see 13.3, first case: formulas 13.3.1, 13.3.2 coinside with 11.2.1
for the given $\vf$ and $a_{***}=0$; $i_j$ of 13.3.0 is $k(1)+k(2)+...+k(j-1)$ of the
present notations).
\medskip
{\bf Definition 11.3.} A standard-1 t-motive is a standard-2 t-motive
whose $\vf$ is the identical permutation $Id$.
\medskip
{\bf 11.3.0.} Let $M=M(Id,k)$ be a standard-1 t-motive. Acting
by $\psi\in S_n$ we can consider only the case of non-increasing $k(j)$. We introduce a
number $\goth m\ge1$ --- the quantity of jumps of $k(j)$, and two sequences
$$0=\gamma_0< \gamma_1<...< \gamma_\goth m=n$$
(sequence of arguments of points of jumps of the function $k$) and
$$0=\lambda_{\goth m+1}< \lambda_{\goth m}<...< \lambda_2 <\lambda_1=l$$
(sequence of values of $k$ on segments $[\gamma_{i-1}+1, ...,\gamma_i]$) by the formulas
$$\matrix k(1)=...=k(\gamma_1)=\lambda_1 \\ \\
k(\gamma_1+1)=...=k(\gamma_2)=\lambda_2\\ ... \\
k(\gamma_{\goth m-1}+1)=...=k(\gamma_\goth m)=\lambda_\goth m\endmatrix \eqno{(11.3.1)}$$
\medskip
{\bf Example 11.3.2.} The t-motive $M$ of 11.1 is a standard-1 having $\goth m=2$,
$\gamma_1=1$, $\gamma_2=2$ and $\lambda_1$, $\lambda_2$ as in 11.1. Its rank $r=\lambda_1+\lambda_2$.
\medskip
{\bf Conjecture 11.3.3.} A standard-2 t-motive of the type $(\vf, k)$ (notations of
11.2.3) is pure iff $\forall j=1, ..., i$ we have:
$$\frac{\alpha_j-\alpha_{j-1}}{\sum_{\gamma=\alpha_{j-1}+1}^{\alpha_{j}}k(\gamma)} = \frac{n}{r}$$
This conjecture is obviously true if all $a_{***}$ are 0.
\medskip
To simplify exposition, we prove here only the following particular case of this
conjecture.
\medskip
{\bf Proposition 11.3.4.} Let $M$ be a standard-1 t-motive having $\goth m>1$,
defined over $\n F_q(\theta)$, having a good reduction at a point of degree 1 of $\n
F_q(\theta)$ (i.e. a point $\theta+c$, $c\in \n F_q$). Then $M$ is not pure.
\medskip
{\bf Proof.} Let $M$ be defined by 11.2.1, we use notations of 11.3.1. We consider the
action of Frobenius on $\tilde M$ --- the reduction of
$M$ at $\theta+c$. According [G], Theorem 5.6.10, it is sufficient to prove that orders
of the roots of the characteristic polynomial of Frobenius over $\w$ are not equal. More
exactly, we consider the valuation infinity on $\w$ (defined by the condition
$\ord(T)=-1$); the order corresponds to a continuation of this valuation to $\End(\tilde
M)$. The
action of Frobenius on $\tilde M$ coincides with multiplication by $\tau$, because the
degree of the reduction point is 1.

A basis $f_*$ of $M_{\p[T]}$ is the set of $X_{\alpha j}:=
\tau^je_{\alpha}$ of 11.2.2. The matrix
$Q(M)$ is defined by the following formulas for the action of
$\tau$ on $X_{\alpha j}$:
$$\tau(X_{\alpha j})=X_{\alpha,j+1} \hbox{ if } j <
k(\alpha)-1\eqno{(11.3.4.1)}$$
$$\tau(X_{\alpha,k(\alpha)-1})=TX_{\alpha,0} -
\sum_{\delta=1}^\goth m
\sum_{d=\lambda_{\delta+1}}^{\lambda_{\delta}-1}
\sum_{c=1}^{\gamma_{\delta}} a_{d\alpha c}X_{cd} \eqno{(11.3.4.2)}$$
This means that if we arrange $X_{\alpha j}$ in lexicographic order ($X_{\alpha_1 j_1}$ precedes to $X_{\alpha_2 j_2}$ if $\alpha_1 <
\alpha_2$) then the matrix
$Q(M)$ has the block form: $$Q(M)=(C_{ij})\ \ \ (i,j=1,...,n)$$ where $C_{ij}$ is a
$k(i)\times k(j)$-matrix of the form
$$C_{ii}=\left(\matrix 0&1&0&...&0\\0&0&1&...&0\\ ...&...&...&...&... \\ 0&0&0&...&1 \\
T-\theta & *&*&...&*\endmatrix \right), \ \ C_{ij}=\left(\matrix 0&0&...&0\\
...&...&...&...\\ 0&0&...&0 \\ 0& *&...&*\endmatrix \right) (i\ne j)$$
where asterisks mean elements $a_{***}$ (in some order). We consider the characteristic polynomial $P(X)\in (\p[T])[X]$ of $Q(M)$. We
have $$C_{ii}-XE_{k(i)}=\left(\matrix -X&1&0&...&0\\0&-X&1&...&0\\ ...&...&...&...&... \\
0&0&0&...&1 \\ t-\theta & *&*&...&*-X\endmatrix \right)$$

A subset of the set of entries of a matrix is called (following N.N.Luzin) a lightning if each row and each column of the matrix contains exactly one element of this subset. The product of elements of a lightning is called the value of this lightning (i.e. the determinant is the alternating sum of the values of all lightnings).
\medskip
{\bf Lemma 11.3.4.3.} If a non-zero lightning of $C_{ii}-XE_{k(i)}$ contains the term
$T-\theta$, then it does not contain any term containing $X$. $\square$
\medskip
Let $J$ be a subset of the set $1,...,n$ and $J'$ its complement.
\medskip
{\bf Corollary 11.3.4.4.} If a non-zero lightning of $Q(M)-XE_{r}$ contains terms
$T-\theta$ of blocks $C_{\al}$, $j\in J$, then its value is a polynomial in $X$ of degree
$\le \sum_{j'\in J'}k(j')$, and there exists exactly one such lightning (called the
principal $J$-lightning) whose value is a polynomial in $X$ of degree $\sum_{j'\in
J'}k(j')$. $\square$

Since the characteristic polynomial of Frobenius of $\tilde M$ is $\tilde P$
(respectively the valuation infinity of $\p[T]$), it is sufficient to prove that the
Newton polygon of $P(X)$ is not reduced to the segment $((0,-n); (r,0))$ defined by its
extreme terms $(T-\theta)^n$ and $X^r$. To do it, it is sufficient to find a point on its
Newton polygon which is below this segment. We consider $J_{min}=$ the set of all
$\gamma_\goth m - \gamma_{\goth m-1}$ diagonal blocks $C_{ii}$ ($i=\gamma_{\goth m-1}+1,
..., \gamma_\goth m$) of $Q(M)$ of minimal size $\lambda_\goth m$. The value of the
principal $J_{min}$-lightning is $(T-\theta)^{\gamma_\goth m - \gamma_{\goth m-1}}$ times
polynomial in $X$ of degree $d:=r-(\gamma_\goth m - \gamma_{\goth m-1})\lambda_\goth m$.
Corollary 11.3.4.4 implies that if the value of any other lightning of $Q(M)-XE_r$ contains a
term whose $X$-degree is equal to $d$, then the $T$-degree of this term is strictly less
than $\gamma_\goth m - \gamma_{\goth m-
 1}$. This means that if we write $P(X)=\sum_{i=0}^r C_iX^i$, $C_i\in \p[T]$, then
$\ord_\infty(C_d)= -(\gamma_\goth m - \gamma_{\goth m-1})$, i.e. the point with
coordinates $[-(\gamma_\goth m - \gamma_{\goth m-1}), d]$ belongs to the Newton diagram
of $P(X)$, i.e. it is above (really, at) the Newton polygon of $P(X)$. This point is
below the segment $((0,-n); (r,0))$. $\square$
\medskip
{\bf Remark 11.3.4.5.} It is easy to see that the Newton polygon of $P(X)$ coincides with
the Newton polygon of the direct sum
of trivial Drinfeld modules of ranks $\lambda_*$, i.e. with the Newton polygon of the
polynomial
$$\prod_{i=1}^{\goth m} (X^{\lambda_i}-T)^{\gamma_i - \gamma_{i-1}}$$
\medskip
{\bf 11.4.} To formulate the below theorem 11.5 we need some notations. Let $M$ be a
standard-1 t-motive defined by formulas 11.2.1, 11.3.1. We impose the condition
$\lambda_\goth m \ge 3$. Theorem 11.5 affirms that it has dual. To find explicitly the
dual of $M$, we need to choose an arbitrary function
$\nu: (i,j) \to \nu(i,j)$ which is a 1 - 1 map from the set of pairs
$(i,j)$ such that

$$1 \le i \le n; \ \ 1 \le j \le k(i)-2 \eqno{(11.4.1)}$$
to the set $[n+1, ..., r-n]$ (recall that $r=\sum_{i=1}^n k(i)= \sum_{i=1}^\goth m
(\gamma_{i}-\gamma_{i-1})\lambda_i$).
\medskip
Let the $(r-n)\times(r-n)$-matrices $B_1$, $B_2$ be defined by the
following formulas
(here and until the end of the proof of 11.5 we have $i ,\alpha = 1, ... , n$; \ \
$b_{\beta\gamma\delta}$
is the $(\gamma\delta)$-th entry of $B_\beta$, all
entries of $B_1$,
$B_2$ that are not in the below list are 0):
\medskip
{\bf 11.4.2.} $b_{1i\alpha}= - a_{k(i)-1,\alpha,i}$;
\medskip
$b_{1,\nu(i,j),\alpha} = - a_{j,\alpha,i}$ for $1 \le j \le k(i)-2$;
\medskip
$b_{1,\nu(i,j+1),\nu(i,j)}=1$ for $1 \le j \le k(i)-3$;
\medskip
$b_{1,i,\nu(i,k(i)-2)} =1$;
\medskip
$b_{2,\nu(i,1),i}=1$.
\medskip
We let $B=\theta E_{r-n}+B_1\tau+B_2\tau^2$ and consider a t-motive $M(B)$ (see 11.5.1
below). Formulas
11.4.2 mean that $M(B)$ is standard-2, its $\vf=\vf_B$ is a product of $n$ cycles
$$i\overset{\vf_B}\to{\to}\nu(i,1)\overset{\vf_B}\to{\to}\nu(i,2)\overset{\vf_B}\to{\to}...
\overset{\vf_B}\to{\to}\nu(i,k(i)-2)\overset{\vf_B}\to{\to}i$$ and its $k=k_B$ is defined
by the
formulas $k_B(\gamma)=2$ for $\gamma \in [1, ...,n]$, $k_B(\gamma)=1$ for $\gamma \in
[n+1, ...,r-n]$.
\medskip
{\bf Theorem 11.5.} Let $M$ be from 11.4 (i.e. a standard-1 t-motive having
$\lambda_\goth m \ge 3$). Then $M'=M(B)$.
\medskip
{\bf Proof.}\footnotemark \footnotetext{This proof is a generalization
of the corresponding proof of
Taguchi; we keep his notations.} Let $e'_*=(e'_1, ... e'_{r-n})^t$ be the vector column
of elements of a basis
of $M(B)$ over $\p[\tau]$ satisfying
$$T e'_*=Be'_*\eqno{(11.5.1)}$$
Let us consider the set of pairs $(j,\goth k)$ such that either $j=1,...,n$,
$\goth k=0,1$ or
$j=n+1,...,r-n$, $\goth k=0$. For each pair $(j,\goth k)$ of this set we let (as
in [T], p. 580)
$Y_{j \goth k} = \tau^{\goth k}e'_j$. Formulas (11.4.2) show that these $Y_{**} $
form a basis of
$M(B)_{\p[T]}$, and the action of $\tau$ on this basis is given by the
following formulas
(here $j=1, ... , k(i)-2$):
$$\tau(Y_{i,0})=Y_{i,1} \eqno{(11.5.2.1)}$$
$$\tau(Y_{i,1})= (T-\theta) Y_{\nu(i,1),0} + \sum_{\gamma=1}^n
a_{1\gamma i}Y_{\gamma,1} \eqno{(11.5.2.2)}$$
$$\tau(Y_{\nu(i,j),0})= (T-\theta) Y_{\nu(i,j+1),0} + \sum_{\gamma=1}^n
a_{j+1,\gamma,i}Y_{\gamma,1} \hbox{ if }j<k(i)-2 \eqno{(11.5.2.3)} $$
$$\tau(Y_{\nu(i,k(i)-2),0})= (T-\theta) Y_{i,0} + \sum_{\gamma=1}^n
a_{k(i)-1,\gamma,i}Y_{\gamma,1} \eqno{(11.5.2.4)}$$
Let $X'_{**}$ be the dual basis to the basis $X_{**}$ of 11.2.2.
\medskip
{\bf 11.5.3.} Let us consider the following correspondence between
$X'_{**}$ and $Y_{**}$:
\medskip
$X'_{ij}$ corresponds to $Y_{\nu(i,j),0}$ for
the pair $(i,j)$ like in (11.4.1),
\medskip
$X'_{i0}$ corresponds to $Y_{i1}$ for $1 \le i \le n$;
\medskip
$X'_{i, k(i)-1}$ corresponds to $Y_{i0}$ for $1 \le i \le n$. 
\medskip
Therefore, in order to prove the Theorem 11.5 we must check that
matrices defined by the dual to (11.3.4.*) and by (11.5.2.*) satisfy (1.10.1) under
identification (11.5.3). This is an elementary exercise. $\square $
\medskip
{\bf Remark 11.6.} Clearly it is possible to generalize the Theorem 11.5 to a larger class
of t-motives --- some subclass of standard-3 t-motives, see Definition 11.8.1.
The below example of the proof of Proposition 11.8.7 shows that probably the condition
$\lambda_\goth m \ge 3$ of
the Theorem 11.5 can be changed by
$\lambda_\goth m \ge 2$: it is necessary
to modify slightly formulas 11.4.2. From another side, a standard-1
t-motive of the Example 2.5 shows that this condition cannot
be changed to $\lambda_\goth m\ge1$.
\medskip
{\bf 11.7. An elementary transformation.} To formulate the proposition
11.7.3, we change slightly notations in 1.9.1, namely, instead of $\goth A =
\sum_{i=0}^l \goth A_i \tau^i$ we consider polynomials $P_k(M)$ of $x_1, ...
,x_n$ ($k=1, ...,n$) defined by the formula
$$P_k(M)= \sum_{i=0}^l\sum_{j=1}^n a_{ikj} x^{q^i}_j \eqno{(11.7.1)} $$
Particularly, if $E$ is the t-module associated to $M$ (see [G], 5.4.5),
$x_*=(x_1, ..., x_n)^t$ an element of $E$ then 11.7.1 is equivalent to
$Tx_*=P_*(x_*)$ where $P_*=(P_1(M),...,P_n(M))^t$ is the vector column.
For a standard-1 t-motive $M$ (we use notations of 11.3.0) having $\goth m\ge 2$ we denote vector columns $\goth
P_1(M)=(P_1(M),...,P_{\gamma_1}(M))^t$, $\goth
P_2(M)=(P_{\gamma_1+1}(M),...,P_{\gamma_2}(M))^t$. We use similar
notations for $M'$.
\medskip
{\bf 11.7.2.} Let $M$ be as above, we consider the case
$\lambda_2=\lambda_1-1$. Let $C$ be a fixed $\gamma_1\times (\gamma_
2-\gamma_1)$-matrix. We define a transformed t-motive $M_1$ by the
formulas

$$\goth P_1(M_1)= \goth P_1(M)+C\goth P_2(M)^q$$

$$P_i(M_1)=P_i(M) \hbox{ for } i>\gamma_1$$
\medskip
{\bf Proposition 11.7.3.} For $M$, $C$, $M_1$ of 11.7.2 the dual $M'_1$
of $M_1$ is described by the following formulas:
$$\goth P_2(M'_1)= \goth P_2(M')-C^t\goth P_1(M')^q$$
$$P_i(M'_1)=P_i(M')\hbox{ for }i\not\in[\gamma_1+1, ... , \gamma_2]$$

{\bf Proof} is similar to the proof of the Theorem 11.5, it is omitted.
$\square$
\medskip
\medskip
{\bf 11.8. Non-pure t-motives.} Most results of this subsection are conditional. We shall show that under some natural conjecture the condition of purity in 8.2 and 8.4 is essential, and that for non-pure t-motives the notion of algebraic duality is richer than the notion of analytic duality.

We generalize slightly the definition 11.2.1 as follows. Let
$\succ$ be a linear ordering on the set $[1,...,n]$, and let $\vf$, $k$ be as in 11.2.
\medskip
{\bf Definition 11.8.1.} A standard-3 t-motive of the type $(\vf, k,\succ)$ is
a t-motive of
dimension
$n$ given by the formulas
$$Te_{\vf(i)}= \theta e_{\vf(i)} +\sum_{j=1}^n\sum_{l=1}^{k(j)-1}a_{l,\vf(i),j}\
\tau^le_j +\sum_{j\succ i} a_{k(j),\vf(i),j}\ \tau^{k(j)}e_j +
\tau^{k(i)}e_i\eqno{(11.8.2)}$$
where $a_{***}\in\p$ are coefficients (the only difference with 11.2.1 is the term
$\sum_{j\succ i} a_{k(j),\vf(i),j}\ \tau^{k(j)}e_j$). We denote it by $M(a_{***})$.

Let $M_1=M(a_{1***})$, $M_2=M(a_{2***})$ be two isomorphic standard-3 t-motives of the same type $(\vf,
k,\succ)$ with
$\p[\tau]$-bases $e_{1*}$, $e_{2*}$ respectively (we use notations of 11.8.2 for both
$M_1$, $M_2$). There exists $C\in M_n(\p[\tau])$ such that the formula defining an isomorphism
between $M_1$ and $M_2$ is the following: $e_{2*}=Ce_{1*}$.
\medskip
{\bf Conjecture 11.8.3.} For a generic set of $a_{1***}$ there exists only a countable set of $a_{2***}$ such that $M_2$ is isomorphic to $M_1$.
\medskip
This conjecture is based on calculations in some explicit cases. Particularly, it is proved if $M_1$, $M_2$ are given by the below formula 11.8.5.1 and entries of $C$ are polynomials in $\tau$ of degree $\le 1$.

We denote by $\Cal M_{u}(r,n)$ the moduli space of uniformizable t-motives of the rank
$r$ and dimension $n$,
by $\Cal L(r,n)$ the moduli space of lattices of the rank $r$ and dimension $n$ and by
$\goth L:\Cal M_{u}(r,n) \to \Cal L(r,n)$ the functor of lattice associated to an uniformizable
t-motive.
\medskip
{\bf Proposition 11.8.5.} Conjecture 11.8.3 implies that the dimension of the fibers of
$\goth L$ is $> 0$ for $r=3$, $n=2$. Particularly, we cannot omit condition of purity in the
statement of 8.2.
\medskip
{\bf Proof.} We consider standard-3 t-motives of the type $n=2$,
$\vf=Id$, $k(1)=2$,
$k(2)=1$, $2\succ 1$. Such $M_1=M_1(a_{111}, a_{112}, a_{121})$ is given by
$$\goth A_0=\theta E_2, \ \ \goth A_{1}= \left(\matrix a_{111} & a_{112}
\\ a_{121} & 1 \endmatrix \right), \ \ \goth A_2= \left(\matrix 1 & 0
\\0 & 0 \endmatrix \right)\eqno{(11.8.5.1)}$$ (notations of Example 10.5).
It has $r=3$, it is not pure, hence it has no dual.
Conjecture 11.8.3 implies
that the dimension of the moduli space of these t-motives is 3 (because there are 3
coefficients
$a_{111}, a_{112}, a_{121}$). Uniformizable t-motives form an open subset of this moduli
space, while
the moduli space of lattices of $n=2$ and $r=3$ has dimension 2. $\square$
\medskip
{\bf Remark.} Similar calculations are valid for any sufficiently large $r$, $n$.
\medskip
Standard-3 t-motives of the above type have not dual. The following proposition
shows that the same
phenomenon holds for t-motives having dual. We denote by $\Cal M_{u,d}(r,n)$ the
moduli space of uniformizable t-motives of the rank $r$ and dimension $n$ having dual, by
$\Cal L_d(r,n)$ the moduli space of lattices of the rank $r$ and dimension $n$ having
dual, by $\goth L_d:\Cal M_{u,d}(r,n) \to \Cal L_d(r,n)$ the functor of lattice and by $D_M: \Cal
M_{u,d}(r,n)\to \Cal M_{u,d}(r,r-n)$, $D_L: \Cal L_d(r,n)\to \Cal L_d(r,r-n)$ the
functors of duality on t-motives and lattices respectively. Practically, Theorem 5
means that the following diagram is commutative:
$$\matrix \Cal M_{u,d}(r,n)&\overset{D_M}\to{\to}&\Cal M_{u,d}(r,r-n)
\\ \\ \goth L_d\downarrow&&\goth L_d\downarrow\\  \\ \Cal L_d(r,n)&\overset{D_L}\to{\to}&\Cal
L_d(r,r-n)\endmatrix\eqno{(11.8.6)}$$

{\bf Proposition 11.8.7.} Conjecture 11.8.3 implies that the dimension of the fibers of
$\goth L_d$ in the diagram (11.8.6) is $> 0$ for $r=5$, $n=2$.
\medskip
Practically, this means that the notion of algebraic duality is "richer" than the notion
of analytic duality.
\medskip
{\bf Proof.} We consider standard-3 t-motives of the type $n=2$,
$\vf=Id$, $k(1)=3$,
$k(2)=2$, $2\succ 1$, $r=5$. Such $M$ is given by
$$\goth A_0=\theta E_2, \ \ \goth A_{1}= \left(\matrix a_{111} & a_{112}
\\ a_{121} & a_{122}  \endmatrix \right), \ \ \goth A_{2}= \left(\matrix a_{211} & a_{212}
\\ a_{221} & 1 \endmatrix \right), \ \ \goth A_3= \left(\matrix 1 & 0
\\0 & 0 \endmatrix \right)$$ (notations of Example 10.5). It has dual. Really, we denote by
$A_{i*j}$ the $j$-th column of $\goth A_i$, and we denote by $ (C_1 | C_2 )$ the matrix formed
by union of columns
$C_1$, $C_2$. Then $M'=M(B)$ is also a standard-3 t-motive, where
\medskip
$B_1= \left(\matrix - \det \goth A_2 & -a_{221} & 1
\\ - \det ( A_{1*2} | A_{2*2} ) & -a_{122} & 0
\\ - \det ( A_{1*1} | A_{2*2} ) & -a_{121} & 0
\endmatrix \right)$,
$B_2 = \left(\matrix 0 & 0 & 0
\\ -a_{212}^q & 1 & 0
\\ 1 & 0 & 0 \endmatrix \right)$

The same arguments as in the proof of Proposition 11.8.5 show that the conjecture 11.8.3
implies that the dimension of the moduli space of these t-motives is 7, while
the moduli space of lattices of $n=2$ and $r=5$ has dimension 6. $\square$
\medskip
As above, similar calculations are valid for any sufficiently large $r$, $n$; clearly the dimension of fibers of $\goth L_d$ becomes larger as $r$, $n$ grow.
\medskip
Let us mention two open questions related to the functor $\goth L$. Firstly, let $L$ be a self-dual lattice such that $L\in \goth L(\Cal M_{u,d}(2n,n))$. This means that $D_M: \goth L_d^{-1}(L) \to \goth L_d^{-1}(L)$ is defined.
\medskip
{\bf Open question 11.8.8.} What can we tell on this functor, for example, what is the dimension of its stable elements?
\medskip
Secondly, let us consider $M_1$, $M_2$ of CM-type with CM-field $\n F_{q^r}(T)$, see 13.3.
\medskip
{\bf Open question 11.8.9.} Let the CM-types $\Phi_1$, $\Phi_2$ of the above $M_1$, $M_2$ satisfy $\Phi_1\ne \alpha \Phi_2$, where $\alpha\in \Gal(\n F_{q^r}(T)/\n F_{q}(T))$. Are lattices $L(M_1)$, $L(M_2)$ non-isomorphic?
\medskip
Clearly the negative answer to this question implies the negative answer to the Question 8.2.
\medskip
For any given $M_1$, $M_2$ the answer can be easily found by computer calculation.
Really, let $M$ be one of $M_1$, $M_2$, $c_1,...,c_r$ a basis of $\n F_{q^r}/\n F_{q}$ and
$\alpha_{1},...,\alpha_{n} \subset \Gal(\n F_{q^r}(\theta)/\n F_{q}(\theta))$ the CM-type of $M$.
We define matrices $\Cal M$, $\Cal N$ as follows: $(\Cal M)_{ij}=\alpha_j(c_i)$ $(i,j=1,...,n$),
$(\Cal N)_{ij}=\alpha_{j}(c_{n+i})$, $j=1,...,n$, $i=1,...,r-n$. The Siegel matrix $Z(M)$ is obviously $\Cal N\Cal
M^{-1}$. So, we can find explicitly $Z(M_1)$, $Z(M_2)$ for both $M_1$, $M_2$. To check
whether $Z(M_1)$, $Z(M_2)$ are equivalent or not, it is sufficient to find a solution to
3.8.1 such that the entries of $A$, $B$, $C$, $D$ are in $M_{*,*}(\n F_q)$ (this is
obvious: the condition $\exists \gamma \in GL_r(\n Z_\infty)$ is equivalent to the
condition $\exists \gamma \in GL_r(\n F_q)$, because entries of $Z(M_1)$, $Z(M_2)$ are in $\n F_{q^r}$). The equation 3.8.1 is linear with respect
to $A$, $B$, $C$, $D$, and we can check whether its solution satisfying $\det
\gamma\ne 0$ exists or not.
\medskip
For the case $q=2$, $r=4$, $n=2$, CM-types of $M_1$, $M_2$ are $(Id, Fr)$, $(Id, Fr^2)$
respectively, a calculation shows that the answer is positive: lattices $L(M_1)$, $L(M_2)$ are not isomorphic.
\medskip
{\bf 12. t-motives having multiplications.}
\medskip
Let $\goth K$ be a separable extension of $\n F_q(T)$ such that $\goth K_C:=\goth K\underset{\n F_q}\to{\otimes} \p$ is also a field, $\pi: X\to P^1(\p)$ the projection of curves over $\p$ corresponding to $\p(T)\subset \goth K_C$. Let $\goth K$, $X$ satisfy the condition: $\infty\in X$ is the only point
on $X$ over $\infty\in P^1(\p)$. Let $\w_{\goth K}$ be the subring of
$\goth K$ consisting of elements regular outside of infinity. We
denote $g=\dim \goth K/\n F_q(T)$ and
$ \alpha_1, ... , \alpha_g: \goth K\to\p$ --- inclusions over $\iota: \n F_q(T)\to\p$
(recall that $\iota(T)=\theta$). Let $\Cal W$ be a central simple algebra
over $\goth K$ of dimension $\goth q^2$. Each $ \alpha_i: \goth K \to
\p$ can be extended to a representation $\chi_i: \Cal W \to M_\goth
q(\p)$.
\medskip
{\bf 12.1. Analytic CM-type.} Let $(L, V)$ be as in Section 2 (recall that
$\w=\n F_q[T]$) such that there exists an inclusion $i: \Cal W \to
\End^0(L, V)$, where $\End^0(L, V)=\End(L, V)\underset{\w}\to{\otimes} \n F_q(T)$.
It defines a representation of $\Cal W$
on $V$ denoted by $\Psi$ which is isomorphic to $\sum_{i=1}^g\goth r_i\chi_i$ where
$\{\goth r_i\}$ are some multiplicities (the CM-type of the action of
$\Cal W$ on $(L, V)$). [Proof: restriction of $\Psi$ on $\goth K$ is a sum of one-dimensional representations, i.e. $V=\oplus_{i=1}^g V_i$ where $k\in \goth K$ acts on $V_i$ by multiplication by $\alpha_i(k)$. Spaces $V_i$ are $\Psi$-invariant. We consider an isomorphism $\Cal W\otimes_\goth K\p=M_\goth q(\p)$ where the inclusion of $\goth K$ in $\p$ is $\alpha_i$. We extend $\Psi|_{V_i}$ to $\Cal W\otimes_\goth K\p$ by $\p$-linearity using the inclusion $\alpha_i$ of $\goth K$ in $\p$. It remains to show that a representation of $M_\goth q(\p)$ is a direct sum of its $\goth q$-dimensional standard representations. We consider the corresponding representation of Lie algebra $\goth s\goth l_\goth q(\p)$. It is a sum of irreducible representations. Let $\omega$ be the highest weight of any of these irreducible representations. $\omega$ is extended uniquely to the set of diagonal matrices of $M_\goth q(\p)$, because $\omega$ is identical on scalars. Since our representation is not only of Lie algebra but of algebra $M_\goth q(\p)$, we get that $\omega$ is a ring homomorphism $\Diag(M_\goth q(\p))\to \p$. There exists the only such $\omega$ corresponding to the $\goth q$-dimensional standard representation].

Further, we
denote $m=\dim_\Cal W L \otimes \n F_q(T)$ ($g$, $\goth q$, $\Psi$, $\goth
r_i$, $m$ are analogs of $g$, $q$, $\Phi$, $r_i$, $m$ of [Sh63] respectively).
Clearly we have
$$n=\goth q\sum_{i=1}^g\goth r_i, \ \ r=mg\goth q^2\eqno{(12.2)}$$
By functoriality, we have the dual inclusion $i': \Cal W^{op} \to
\End^0(L',V')$ where $\Cal W^{op}$ is the opposite algebra.
\medskip
{\bf Remark.} A construction of Hilbert-Blumental modules ([A], 4.3, p. 498) practically is a particular case
of the present
construction: for Hilbert-Blumental modules we have $\goth q=1$, i.e. $\goth K=\Cal W$, and all $\goth r_i=1$. Anderson considers the case when $\infty$ splits completely; this difference with the present case is not essential.
\medskip
{\bf Proposition 12.3.} If the dual pair $(L',V')$ exists then the
CM-type of the dual inclusion is $\{m\goth
q - \goth r_i\}$, $i=1, ... ,g$.
\medskip
{\bf Proof.} We have $L\underset{\n Z_\infty}\to{\otimes}\p$ is
isomorphic to $(\Cal W\underset{\n F_q(\theta)}\to{\otimes}\p)^m$ as a
$\Cal W$-module. Since the natural representation of $\Cal W$ on $\Cal
W\underset{\n F_q(\theta)}\to{\otimes}\p$ is isomorphic to $\goth
q\sum_{i=1}^g\chi_i$ we get that $L\underset{\n Z_\infty}\to{\otimes}\p$ is isomorphic to $m\goth q\sum_{i=1}^g\chi_i$ as a $\Cal W$-module. Consideration of the exact
sequence $0\to {V'}^* \to L\underset{\n Z_\infty}\to{\otimes}\p \to V\to0$ gives us the desired. $\square$
\medskip
{\bf Remark 12.4.1.} This result is an analog of the corresponding
theorem in the number field case. We use notations of [Sh63], Section 2. Let
$A$ be an abelian variety having endomorphism algebra of type IV, and
$(r_\nu, s_\nu)=(r_\nu(A), s_\nu(A))$ are from [Sh63], Section 2, (8).
Then
$$r_\nu(A')=mq-r_\nu(A)=s_\nu(A), \ s_\nu(A')=mq-s_\nu(A)=r_\nu(A)$$
By the way, Shimura writes that the CM-types of $A$ and $A'$ coincide
([Sh98], 6.3, second line below (5), case $A$ of CM-type). We see that
his affirmation is not natural: he considers the complex conjugate
action of the endomorphism ring on $A'$. It is necessary to take into
consideration this difference of notations comparing formulas of 12.3 and 13.2
with the corresponding formulas of Shimura.
\medskip
{\bf Remark 12.4.2.} According [L1], a t-motive $M$ is an analog of an abelian variety $A$ with multiplication by an imaginary quadratic field $K$. The above consideration shows that this analogy holds for $M$ and $A$ having more multiplications. Really, if $A$ has more multiplications then (we use notations of [Sh63], Section 2) $F_0=FK$, and numbers $(r_\nu(A), s_\nu(A))$ satisfy $n(A)=q\sum_{i=1}^g r_\nu(A)$, where $(n(A), \dim(A)-n(A))$ is the signature of $A$ treated as an abelian variety with multiplication by $K$. This is an analog of 12.2.
\medskip
{\bf 12.5. Complete multiplication.} Here we consider the case $\goth
q=m=1$, i.e. $\goth K=\Cal W$ and $g=r$.
\medskip
{\bf Lemma 12.5.1.} In this case the condition $N=0$
implies that the CM-type $$\sum_{i=1}^r\goth r_i\alpha_i\eqno{(12.5.2)}$$ of the action of
$\goth K$ on
on $(L,V)$ has the property: all $\goth r_i$ are 0 or 1.
\medskip
{\bf Proof.} $N=0$ means that
the action of $T\in\w$ on $V$ is simply multiplication by $\theta$. We write the
CM-type $\sum_{i=1}^r\goth r_i\chi_i$ in the form $\sum_{i=1}^n\chi_{\alpha_i}$ where
$\alpha_1,...,\alpha_n\in[1,...,r]$ are not necessarily distinct. Let $l_1$ be an (only)
element of a basis of $L\otimes_{\w_{\goth K}}\goth K$ over $\goth K$ and $e_1,...,e_n$ a
basis of $V$
over $\p$ such that the action of $\goth K$ on $V$ is given by the formulas
$$k(e_i)=\chi_{\alpha_i}(k)e_i, \ \ \ k\in \goth K$$
Multiplying $e_i$ by scalars if necessary, we can assume that $l_1=\sum e_i$. Therefore,
if $\alpha_i=\alpha_j$ (i.e. not all $\goth r_*$ in (12.5.2) are 0, \ 1) then the
$e_{\alpha_i}$-th coordinate of any element of $L$ coincide with its $e_{\alpha_j}$-th
coordinate, hence $L$ does not $\p$-generate $V$ --- a contradiction. $\square$
\medskip
{\bf 12.5.3.} Let $M$ be a t-motive of
rank $r$ and dimension $n$ having multiplication by $\w_{\goth K}$. Recall that we
consider only
the case $N=0$. This means that the character of the action of $\goth K$ on $M/\tau M$
is isomorphic to $\sum_{i=1}^r\goth r_i\alpha_i$. Since $E(M)=(M/\tau M)^*$ we get that
the character of the action of $\goth K$ on $E(M)$ is the same. If
$$\hbox{all $\goth r_i$ are 0 or 1}\eqno{(12.5.4)}$$
we shall use the terminology that $M$ has the CM-type $\Phi\subset \{\alpha_1, ... ,
\alpha_r \}$
where $\Phi$ is defined by the condition $\alpha_i\in \Phi \iff \goth r_i=1$.

It is easy to see that this case occurs for uniformizable $M$. Really, if $M$ is
uniformizable then
the action of $\goth K$ can be prolonged on $(L(M),V(M))$, and the character
of the action of $\goth K$ on $V(M)$ coincides with the one on $E(M)$. The result follows
from Lemma 12.5.1.
\medskip
{\bf Lemma 12.5.5.} There exists a canonical isomorphism $\gamma$ from the set of
inclusions $\alpha_1, ... , \alpha_r$ to the set of points $ \theta_{\alpha_1}, ... ,
\theta_{\alpha_r}$ of $X$ over $\theta\in P^1(\p)$.
\medskip
{\bf Proof.} A point $t\in X$ over $\theta\in P^1(\p)$ defines a function
$\vf_t: \goth K_{C} \to P^1(\p)$ --- the value of
an element $f\in \goth K_{C}$ treated as a function on $X$ at the point $t$. This
function
must satisfy the standard axioms of valuation and the condition $\vf_t(T)=\theta$. Let
$\alpha_i$
be an inclusion of $\goth K$ to $\p$ over $\iota$. It defines a valuation
$\vf_{\alpha_i}: \goth K_{C} \to P^1(\p)$ by the formula $\vf_{\alpha_i}(k\otimes
f)=\alpha_i(k) f(\theta)$,
where $k\in\goth K$, $f\in \p(T)$. We define $\gamma(\alpha_i)$ by the condition
$\vf_{\gamma(\alpha_i)}=\vf_{\alpha_i}$; it is easy to see that $\gamma$ is an
isomorphism. $\square$
\medskip
{\bf Theorem 12.6.} For any above \{$\goth K$, $\Phi$\} there exists an
t-motive $(M,\tau)$ with
complete multiplication by $\goth K$ having CM-type $\Phi$.
\medskip
{\bf Proof (Drinfeld).} We denote the divisor $\sum_{\alpha_i\in\Phi}\gamma({\alpha_i})$
by
$\theta_\Phi$. We construct a $\Cal F$-sheaf $F$ of dimension 1 over
$\goth K$ which will give us $M$.
Let $\fr$ be the Frobenius map on $ \Pic_0(X)$. It is an algebraic
map, and the $\fr-\Id: \Pic_0(X) \to \Pic_0(X)$ is an algebraic map as
well. Since the action of $\fr$ on the tangent space of $ \Pic_0(X)$
at 0 is the zero map, the action of $\fr-\Id$ on the tangent space of
$ \Pic_0(X)$ at 0 is the minus identical map and hence $\fr-\Id$ is an
isogeny of $ \Pic_0(X)$. Particularly, there exists a divisor $D$ of
degree 0 on $X$ such that we have the following equality in $\Pic_0(X)$:
$$\fr(D)-D=-\theta_\Phi+n\infty \eqno{(12.6.0)}$$
This means that if we let $F=F_\Phi=O(D)$ then there exists a rational map
$\tau_X=\tau_{X,\Phi}: F^{(1)}\to F$ such
that
$$\Div(\tau_X)=\theta_\Phi-n\infty\eqno{(12.6.1)}$$
The pair $(F_\Phi, \tau_{X,\Phi})$ is the desired $\Cal F$-sheaf.
\medskip
{\bf Remark.} It is easy to see that if the genus of $X$ is $> 0$ then different CM-types $\Phi_1$, $\Phi_2$ give us different sheaves $F_{\Phi_1}$, $F_{\Phi_2}$, while if the genus of $X$ is 0 then $F_{\Phi_1}=F_{\Phi_2}=\Cal O$, but the maps $\tau_{X,\Phi_1}$, $\tau_{X,\Phi_2}$ are clearly different.
\medskip
Let $U_0=X-\{\infty\}$ be an open part of
$X$. We denote $F(U_0)$ by $\Cal M$, hence $F^{(1)}(U_0)=\Cal M^{(1)}$.
Since the support of the negative part of the right hand side of
12.6.1 is $\{\infty\}$,
we get that the (a priory rational) map $\tau_X(U_0): \Cal M^{(1)}\to \Cal M$ is
really a map of
$\w_{\goth K}$-modules.

Let $M$ be a $\p[T]$-module obtained from $\Cal M$ by restriction of scalars from
$\w_{\goth K}$ to $\p[T]$. Construction $F\mapsto M$ is functorial, and we denote this functor by $\delta$. Further, we denote by $\alpha$ the tautological isomorphism
$\Cal M\to M$. $M$ is a free $r$-dimensional $\p[T]$-module, and (because
$M^{(1)}$ is isomorphic to $M$) the
same restriction of scalars of $\tau_X(U_0)$ defines us a $\p[T]$-skew map from $M$ to
$M$ denoted by $\tau$ (skew means that $\tau(zm)=z^q\tau(m)$, $z\in\p$). $\tau$ is
defined by the formula $\tau(m)=\alpha\circ \tau_X((\alpha^{-1}(m))^{(1)})$.

It is easy to check that
$(M,\tau)$ is the required t-motive. Really, $M$ is a $\w_{\goth K}$-module, and
$\tau$
commutes with this multiplication. The fact that the positive part
of the right hand
side of 12.6.1 is $\theta_\Phi$ means that 1.13.2 holds for $M$ and that the CM-type of
the action of
$\w_{\goth K}$ is $\Phi$.
\medskip
{\bf Remark 12.6.2.} It is easy to prove for this case that $M$ is a free
$\p[\tau]$-module.
Really, it is
sufficient to prove (see [G], Lemma 12.4.10)
that $M$ is finitely generated as a $\p[\tau]$-module. We choose $D$ such that
$\infty\not\in \Supp (D)$. There exists $P\in \goth K_{C}^*$
such that $\tau_X(U_0): \Cal M^{(1)}\to \Cal M$ is multiplication by $P$ (recall that both
$\Cal M^{(1)}$, $\Cal M$ are $\w_{\goth K}$-submodules of $\goth K$). 12.6.0 implies that
$-\ord_\infty(P)=n$.
There
exists a number $n_1$ such that $$\hbox{(a) $h^0(X,\Cal O(D+n_1\infty))>0$; \ \ (b)
for any $k\ge 0$ we have}$$ $$h^0(X,\Cal O(D+(n_1+k)\infty))=h^0(X,\Cal
O(D+n_1\infty))+k\eqno{(12.6.3)}$$ $$h^0(X,\Cal O(D^{(1)}+(n_1+k)\infty))=h^0(X,\Cal
O(D^{(1)}+n_1\infty))+k\eqno{(12.6.4)}$$
It is sufficient to prove that if $g_1,...,g_k$ are elements of a basis of $H^0(X,\Cal
O(D+(n_1+n)\infty))$, then
for any $Q\in \Cal M$ the element $\alpha(Q)\in M$ is generated by
$\alpha(g_1),...,\alpha(g_k)$ over
$\p[\tau]$. We prove it by induction by $n_2:=-\ord_\infty(Q)$. If $n_2\le n_1+n$ the
result is trivial. If not
then 12.6.3, 12.6.4 imply that the multiplication by $P$ defines an isomorphism
$$H^0(X, \Cal O(D^{(1)}+(n_2-n)\infty))/H^0(X, \Cal O(D^{(1)}+(n_2-n-1)\infty))\to$$
$$\to H^0(X,\Cal O(D+n_2\infty))/H^0(X, \Cal O(D+(n_2-1)\infty))$$
This means that $\exists Q_1\in H^0(X,\Cal O(D^{(1)}+(n_2-n)\infty))$,
$-\ord_\infty(Q_1)=n_2-n$ such
that
$-\ord_\infty(Q-PQ_1)\le n_2-1$. An element $Q_1^{(-1)}\in \Cal M$ exists;
since $\alpha(Q)=\tau(\alpha(Q_1^{(-1)}))+\alpha(Q-PQ_1)$, the result follows by
induction.
$\square$
\medskip
If $\goth K$ and $\Phi$ are given then the construction of the Theorem 12.6 defines $F$ uniquely up to tensoring by $O(D)$ where $D\in \Div (X(\goth K))$. We denote the set of these $F$ by $F$(\{$\goth K$, $\Phi$\}), and we denote by $M$(\{$\goth K$, $\Phi$\}) the set $\delta(F$(\{$\goth K$, $\Phi$\})). Further, we denote by $\Phi'=\{\alpha_1, ... ,
\alpha_r\}-\Phi$ the complementary CM-type.
\medskip
{\bf Theorem 12.7.} Let $M\in M(\{\goth K,\Phi\})$. Then $M'$ exists, and $M'\in M(\{\goth K,\Phi'\})$. More exactly, if $F\in F$(\{$\goth K$, $\Phi$\}) then $F^{-1}\otimes \Cal D^{-1}\in F$(\{$\goth K$, $\Phi'$\}) where $\Cal D$ is the different
sheaf on $X$, and if $M=\delta(F)$ then $M'=\delta(F^{-1}\otimes \Cal D^{-1})$.
\medskip
{\bf Proof.} Let $G$ be any invertible sheaf on $X$. We have a
\medskip
{\bf Lemma 12.7.0.} There exists the canonical isomorphism $\varphi_G:
\pi_*(G^{-1}\otimes \Cal D^{-1}) \to \Hom_{P^1}(\pi_*(G), \Cal O)$.
\medskip
{\bf Proof.} At the level of affine open sets $\varphi_G$ comes from the trace bilinear form of field extension $\goth K/\n F_q(T)$. Concordance with glueing is obvious. $\square$
\medskip
We need the relative version of this lemma. Let $G_1$, $G_2$ be invertible sheaves on $X$,
$\rho: G_1 \to G_2$ any rational map. Obviously there exists a rational map
$\rho^{-1}: G_1^{-1} \to G_2^{-1}$. Recall that we denote by
$\rho^{inv}: G_2 \to G_1$ the rational map which is inverse to $\rho$
respectively the composition. The map $\pi_*(\rho^{-1}\otimes \Cal D^{-1}): \pi_*(G_1^{-1}\otimes \Cal D^{-1})\to \pi_*(G_2^{-1}\otimes \Cal D^{-1})$ is obviously defined. The map (denoted by $\beta(\rho)$) from $\Hom_{P^1}(\pi_*(G_1), \Cal O)$ to $\Hom_{P^1}(\pi_*(G_2),
\Cal O)$ is defined as follows at the level of affine open sets: let
$\gamma\in\Hom_{P^1}(\pi_*(G_1), \Cal O)(U)$ where $U$ is a sufficiently
small affine subset of $P^1$, such that we have a map $\gamma(U): \pi_*(G_1)(U)
\rightarrow \Cal O(U)$. Then $(\beta(\gamma))(U)$ is the composition
map $\gamma(U)\circ \pi_*(\rho^{inv})(U)$:
$$\pi_*(G_2)(U)\overset{\pi_*(\rho^{inv})(U)}\to{\longrightarrow}
\pi_*(G_1)(U)\overset{\gamma(U)}\to{\to}\Cal O(U)$$

{\bf Lemma 12.7.1.} The above maps form a commutative diagram:
$$\matrix
\pi_*(G_1^{-1}\otimes \Cal D^{-1})&\overset{\pi_*(\rho^{-1}\otimes \Cal
D^{-1})}\to{\longrightarrow}&\pi_*(G_2^{-1}\otimes \Cal D^{-1})&&\\&&&&\\
\varphi_{G_1}\downarrow&&\varphi_{G_2}\downarrow &&\\&&&&\\
\Hom_{P^1}(\pi_*(G_1), \Cal
O)&\overset{\beta(\rho)}\to{\longrightarrow}&\Hom_{P^1}(\pi_*(G_2),
\Cal O)&&\square
\endmatrix$$

We apply this lemma to the case $\{\rho: G_1 \to G_2\}=\{\tau_{X, \Phi}:
F^{(1)}\to F\}$. We have:
$$\Div(\tau_{X, \Phi}^{-1}\otimes \Cal D^{-1})=-\Div(\tau_{X, \Phi})=-\theta_{\Phi}+n\infty$$
Futher, we multiply $\tau_{X, \Phi}^{-1}\otimes \Cal D^{-1}$ by $T-\theta$. We have:
$$\Div((T-\theta)\tau_{X, \Phi}^{-1}\otimes \Cal
D^{-1})=\Div(T-\theta)+\Div(\tau_{X, \Phi}^{-1}\otimes \Cal D^{-1})=\theta_{\Phi'}-(r-n)\infty$$
i.e. $(T-\theta)\tau_{X, \Phi}^{-1}\otimes \Cal D$ is one of
$\tau_{X,\Phi'}$, i.e. $F^{-1}\otimes \Cal D^{-1}\in F$(\{$\goth K$, $\Phi'$\}). Further, $(T-\theta)\beta(\tau_{X, \Phi})$ is the map which is used in the definition of duality of
$M$. This means that the lemma 12.7.1 implies the theorem. $\square$
\medskip
{\bf Remark 12.8.} There exists a simple proof of the second part of the Theorem 5 for uniformizable abelian
t-motives $M$ with complete multiplication by $\w_\goth K\subset \goth K$. Recall that this second part is the proof of 2.7 for $M$. Really, let us
consider the diagram 2.5. The CM-types
of action of $\goth K$ on $\Lie(M)$ and on $E(M)$ coincide,
and the CM-types of action of $\goth K$ on a vector space and on its
dual space coincide. This means that the CM-type of $V^*$ is $\Phi$
and the CM-type of $V'$ is $\Phi'$. Further, $\gamma_D$ of 2.5 commutes with complete multiplication: this follows immediately for example from a description of $\gamma_D$ given in Remark 5.2.8. Really, all homomorphisms of 5.2.9 commute with complete multiplication. For example, this condition for $\delta$ of 1.11.1 is written as follows: if $k\in \goth K$, $\goth m_k(M)$, resp. $\goth m_k(M')$ is the map of complete multiplication by $k$ of $M$, resp. $M'$, then $(\goth m_k(M)\otimes Id) \circ \delta= (Id \otimes \goth m_k(M'))\circ \delta$ --- see any textbook on linear algebra.

Finally, since
$\Phi\cap \Phi'=\emptyset$ and the map $\vf'\circ \gamma_D \circ
\vf^*$ commutes with complete multiplication, we get that it must be
0.
\medskip
{\bf 13. Miscellaneous.}

\nopagebreak
\medskip
Let now $(L,V)$ be from 12.1, case $\goth q=m=1$, i.e. $\goth K=\Cal W$
and $r=g$, and let the ring of complete multiplication be the maximal
order $\w_\goth K$. We identify $\w$ and $\n Z_\infty$ via $\iota$, i.e. we consider $\goth K$ as an extension of $\n F_q(\theta)$. Let $\Phi$ be the CM-type of the action of $\goth
K$ on $V$.
This means that --- as an
$\w_\goth K$-module --- $L$ is isomorphic to $I$ where $I$ is an ideal of
$\w_\goth K$. The class of $I$ in $\Cl(\w_\goth K)$ is defined by $L$
and $\Phi$ uniquely; we denote it by $\Cl(L,\Phi)$.
\medskip
{\bf Remark.} $\Cl(L,\Phi)$ depends on $\Phi$, because the action of $\w_\goth K$ on $V$ depends on $\Phi$. Really, let $a\in L \subset V$, $a=(a_1,...,a_n)$ its coordinates, $\Phi=\{\alpha_{i_1},...,\alpha_{i_n}\}\subset \{\alpha_{1},...,\alpha_{r}\}$ and $k\in \w_\goth K$. Then $ka$ has coordinates $(\alpha_{i_1}(k)a_1,...,\alpha_{i_n}(k)a_n)$, i.e. depends de $\Phi$. Particularly, the $\w_\goth K$-module structure on $L$ depends on $\Phi$, and hence $\Cl(L,\Phi)$ depends on $\Phi$. For example, if $n=1$, $r=2$, $\Phi_1=\{\alpha_{1}\}$, $\Phi_2=\{\alpha_{2}\}$, then $\Cl(L,\Phi_2)$ is the conjugate of $\Cl(L,\Phi_1)$.
\medskip
{\bf Theorem 13.1.} $\Cl(L',\Phi')=(\Cl(\goth d))^{-1}
(\Cl(L,\Phi))^{-1}$ where
$\goth d$ is the different ideal of the ring extension $\w_\goth K/\w$.
\medskip
{\bf Proof.} This theorem follows from the above results;
nevertheless, I give here an explicit elementary proof. Let $a_*=(a_1, ..., a_r)^t$ be a
basis (considered
as a vector column) of $\goth K$ over $\n F_q(\theta)$ and $b_*=(b_1,
..., b_r)^t$ the dual basis. Recall that it satisfies 2 properties:
$$(1) \ \ \forall i\ne j \ \ \alpha_i(a_*)^t\alpha_j(b_*)=0 \ \ \
(\hbox{ i.e. } \sum_{k=1}^r \alpha_i(a_k) \alpha_j(b_k)=0
)\eqno{(13.1.1)} $$
(2) For $x\in \goth K$ let $\goth m_{x,a_*}$ (resp $\goth m_{x,b_*}$)
be the matrix of multiplication by $x$ in the basis $a_*$ (resp.
$b_*$). Then for all $x\in \goth K$ we have
$$\goth m_{x,a_*}=\goth m_{x,b_*}^t\eqno{(13.1.2)} $$

We define $\goth E_{n,r-n}$
as an $r\times r$ block matrix $\left(\matrix 0&E_{r-n}\\ -E_n&0
\endmatrix\right)$, and we define a new basis $\tilde b_*=(\tilde b_1,
..., \tilde b_r)^t$ by
$$\tilde b_*=\goth E_{n,r-n}b_*\eqno{(13.1.3)} $$
(explicit formula: $(\tilde b_1, ..., \tilde b_r)=(b_{n+1},
..., b_r, -b_1, ..., -b_n)$).

We can assume that $\Phi=\{\alpha_1, ... , \alpha_n \}$. Since $L$ has multiplication by $\w_\goth K$ and the CM-type of this
multiplication is $\Phi$, it is possible to choose $a_*$ such that $L\subset
\p^{n}$ is generated over $\n Z_\infty$ by $e_1, ... , e_r$ where
$$e_i=(\alpha_1(a_i), ... , \alpha_n(a_i))\eqno{(13.1.4)}$$ Let $\hat L\subset \p^{r-n}$ be
generated over $\n Z_\infty$ by $\hat e_1, ... , \hat e_r$ where
$$\hat e_i=(\alpha_{n+1}(\tilde b_i), ... , \alpha_r(\tilde b_i))\eqno{(13.1.5)}$$

{\bf Lemma 13.1.6. } $L'=\hat L$.
\medskip
{\bf Proof. } Let $A$ (resp. $B$) be a matrix whose lines are the lines of
coordinates of $e_1, ... , e_n$ (resp. $e_{n+1}, ... ,
e_r$) in 13.1.4, and $C$ (resp. $D$) a matrix whose lines are the lines of
coordinates of $\hat e_1, ... , \hat e_{r-n}$ (resp. $\hat e_{r-n+1}, ... ,
\hat e_r$) in 13.1.5. By definition of Siegel matrix, we have $L=\goth L(BA^{-1})$,
$\hat L=\goth L(DC^{-1})$ ($\goth L$ is defined in 3.1, 3.2). So, it is sufficient to prove that
$(BA^{-1})^t=DC^{-1}$, i.e. $A^tD=B^tC$. This follows immediately from
the definition of $A,B,C,D$ and (13.1.1). $\square$
\medskip
For $x\in \w_\goth K$ we denote by $\goth
M_x(L)$ the matrix of multiplication by $x$ in the basis $e_*$ (see
the notations of Remark 3.8). Obviously $\goth M_x(L)=\goth
m_{x,a_*}$.

Let now $\w_\goth K$ acts on $\p^{r-n}$ (the ambient space of $L'$) by CM-type
$\Phi'$. According (13.1.2) and
(13.1.3), the matrix of the action of $x\in \w_\goth K$ in the basis
$\tilde b_*$ is
$$\goth E_{n,r-n}\goth m_{x,a_*}^t\goth E_{n,r-n}^{-1}\eqno{(13.1.7)} $$
Let $\goth M$, $\goth M'$ be from Remark 3.8. Formula 3.8.4
shows that
$$\goth M'=\goth E_{n,r-n}\goth M^t\goth E_{n,r-n}^{-1}\eqno{(13.1.8)} $$
Formulas (13.1.7) and (13.1.8) --- because of Lemma 1.10.3 --- prove the theorem. $\square$
\medskip
{\bf 13.2. Compatibility with the weak form of the main theorem of
complete multiplication.}

\nopagebreak
\medskip
The reader can think that Theorem 13.1 is incompatible with the main theorem of complex
multiplication, because of the $-1$-th power in its statement. The reason is a bad choice
of notations of Shimura, he affirms that the CM-type of an abelian variety $A$ over a number field coincides with the CM-type
of $A'$, while we see that it is really the complement. Since an analog of even the weak form of the main theorem of complex multiplication --- Theorem 13.2.6 --- for the function field case is not proved yet, the main result of the present section --- Theorem 13.2.8 --- is conditional: it affirms that if this weak form of the main
theorem --- Conjecture 13.2.7 --- is true for a t-motive with complete multiplication $M$, then it is true for $M'$ as well. By the way, even if it will turn out that the statement of the Conjecture 13.2.7 is not correct, the proof of 13.2.8 will not be affected, because the main
ingredient of the proof is the formula 13.2.10 "neutralizing" the $-1$-th power of the
Theorem 13.1.

Let us recall some definitions of [Sh71], Section 5.5. We consider an abelian variety
$A=\n C^n/L$ with
complex multiplication by $K$. The set $\Hom(K, \bar \n Q)$ consists of $n$ pairs of
mutually
conjugate inclusions $\{\vf_1,\bar \vf_1, ..., \vf_n,\bar\vf_n\}$. $\Phi$ is a subset of
the set
$\Hom(K, \bar \n Q)$ such that $\forall i=1,...,n$ we have:
$$\hbox{$\Phi\cap\{\vf_i,\bar \vf_i\}$ consists of one element.}\eqno{(13.2.1)}$$
It is defined by the condition that the action of complex multiplication
of $K$ on $\n C^n$ is isomorphic to the direct sum of the elements of $\Phi$. Let $F$ be
the
Galois envelope of $K/\n Q$,
$$G:=\Gal (F/\n Q), \ \ H:=\Gal(F/K), \ \
S:=\bigcup_{\alpha\in\Phi}H\alpha\eqno{(13.2.2)}$$
(the elements of Galois group act on $x\in F$ from the right, i.e. by the formula
$x^{\alpha\beta}=(x^\alpha)^\beta$; for $\alpha\in\Phi$ we denote by $\alpha$ also a
representative in $G$ of the coset $\alpha$). We denote
$$H^{ref}:=\{\gamma\in G\vert S\gamma=S\}\eqno{(13.2.3)}$$
and let $K^{ref}$ be the subfield of $F$ corresponding to $H^{ref}$. We have:
$$H^{ref}S^{-1}=S^{-1}\eqno{(13.2.4)}$$
i.e. $S^{-1}$ is an union of cosets of $H^{ref}$ in $G$. We can identify these cosets
with elements of $\Hom(K^{ref}, \bar \n Q)$. $\Phi^{ref}\subset \Hom(K^{ref}, \bar \n Q)$
is, by definition, the set of these cosets. There is a map $\det
\Phi^{ref}:{K^{ref}}^\times \to K^\times$ defined as follows:
$$\det \Phi^{ref}(x):=\prod_{\alpha\in\Phi}\alpha(x)\eqno{(13.2.5)}$$
(it follows easily from the above formulas and definitions that $\det \Phi^{ref}(x)$
really belongs to $K^\times$). It can be extended to the group of ideles and factorized to
the group of classes of ideals, we denote this map
by $\det_{Cl} \Phi^{ref}: \Cl(K^{ref}) \to \Cl(K)$. Finally, let
$\theta^{ref}: \Gal(K^{ref\ Hilb}/K^{ref}) \to \Cl(K^{ref})$ be an isomorphism defined by
the Artin reciprocity law.

We consider the case $\End(A)=O_K$. In this case $L$ is isomorphic to an ideal of $O_K$,
its class
is well-defined by the class of isomorphism of $A$, we denote it by $\Cl(A)$.
\medskip
{\bf Theorem 13.2.6.} $A$ is defined over $K^{ref\ Hilb}$;

For any $\gamma\in\Gal(K^{ref\ Hilb}/K^{ref})$ we have

$\Cl(\gamma(A))=\det_{Cl} \Phi^{ref}\circ \theta^{ref} (\gamma)^{-1}(\Cl(A))$. $\square$
\medskip
This is a weak form of [SH71], Theorem 5.15 --- the main theorem of
complex multiplication.
\medskip
Now we define analogous objects for the function field case in order to formulate a
conjectural
analog of Theorem 13.2.6. Let $\goth K$, $\Phi$ be from 12.5.3. $\goth K^{ref}$,
$\Phi^{ref}$,
$\det \Phi^{ref}$ are defined by the same formulas 13.2.2 -- 13.2.5 like in the number
field case
($\n Q$ must be replaced by $\n F_q(T)$). The facts that 13.2.1 has no meaning in the
function field case and
that the order of $S$ is not necessarily the half of the order of $G$ do not affect the
definitions.

The $\infty$-Hilbert class field of $\goth K$ (denoted by $\goth K^{Hilb \ \infty}$) is an abelian
extension of
$\goth K$ corresponding to the subgroup
$$\goth K_\infty^*\cdot \prod_{v\ne\infty}O_{\goth K_v}^*\cdot \goth K^*$$
of the idele group of $\goth K$. We have an isomorphism $\theta: \Gal(\goth K^{Hilb \ \infty}/\goth K) \to
\Cl(\w_{\goth K})$.

We formulate the function field analog of Theorem 13.2.6 only for the case when
\medskip
{\bf (*)} There exists only one
point over $\infty\in P^1(\n F_q)$ in the extension $\goth K^{ref}/\n F_q(T)$.
\medskip
In this case the field
$\goth K^{ref \ Hilb \ \infty}$ and the ring $\w_{\goth K^{ref}}$ are naturally defined, and we have
an
isomorphism
$\theta^{ref}: \Gal(\goth K^{ref \ Hilb \ \infty}/\goth K^{ref}) \to \Cl(\w_{\goth K^{ref}})$.

Let $M$ be an uniformizable
t-motive of rank $r$ and dimension $n$ having complete multiplication
by $\w_\goth K$, and $\Phi$ its CM-type. $\Cl(M)$ is defined like $\Cl(A)$ in the number
field case,
it is $\Cl(L,\Phi)$ of 13.1.
\medskip
{\bf Conjecture 13.2.7.} If (*) holds, then $M$ is defined over $\goth K^{ref\ Hilb \
\infty}$, and for any
$\gamma\in\Gal(\goth K^{ref\ Hilb \ \infty}/\goth K^{ref})$
we have $\Cl(\gamma(M))=\det_{Cl} \Phi^{ref}\circ \theta^{ref} (\gamma)^{-1}\Cl(M)$.
\medskip

Now we can formulate the main theorem of this section.
\medskip
{\bf Theorem 13.2.8.} If conjecture 13.2.7 is true for $M$ then it is true for $M'$.
\medskip
{\bf Proof.} It follows immediately from the functional analogs of 13.2.2 -- 13.2.4 that
$$(\goth K,\Phi')^{ref}=(\goth K^{ref},(\Phi^{ref})')\eqno{(13.2.9)}$$
Further,
$$\hbox{det}_{Cl} {\Phi'}^{ref}=(\hbox{det}_{Cl} \Phi^{ref})^{-1}\eqno{(13.2.10)}$$
Really, $\det \Phi^{ref}(x) \cdot \det (\Phi^{ref})'(x) = N_{\goth
K^{ref}/\n F_q(T)}(x)\in\n F_q(T)^\times$, hence gives the trivial class of ideals (we use
here (13.2.9).
Finally, for $\gamma\in\Gal(\goth K^{ref})$ we have
$$(\gamma(M))'=\gamma(M')\eqno{(13.2.11)}$$
The theorem follows immediately from 13.1, 13.2.10, 13.2.11 (recall that $\Cl(M)$ is
$\Cl(L,\Phi)$ of 13.1). $\square$
\medskip
{\bf 13.3. Some explicit formulas.} We give
here an elementary explicit proof of the theorem 12.7 in two simple
cases: $\goth K=\n F_{q^r}(T)$ and $\n F_q(T^{1/r})$. By the way,
since the extension $\n F_{q^r}(T)/\n F_{q}(T)$ is not absolutely
irreducible, formally this case is not covered by the theorem 12.7.
\medskip
{\bf Case $\w_\goth K= \n F_{q^r}[T]$.} Let $\alpha_i$, where $i=0, ... ,r-1$, be inclusions $\goth K \to \p$. For $\omega \in\n F_{q^r}$
 we have $\alpha_i(\omega )=\omega^{q^{i}}$. Let $$0 \le
i_1 < i_2 < ... <
i_n \le r-1\eqno{(13.3.0)}$$ be numbers such that $\Phi=\{\alpha_{i_j}\}$, $j=1, ... , n$.
We consider the following t-motive $M=M(\goth K, \Phi)$. Let $e_1, ... , e_n$ be a basis of $M_{\p[\tau]}$ such that $\goth
m_\omega(e_j)=\omega^{q^{i_j}}e_j$ and the multiplication by $T$ is
defined by formulas
$$Te_1=\theta e_1 + \tau^{i_1-i_n+r}e_n\eqno{(13.3.1)} $$
$$Te_j=\theta e_j + \tau^{i_j-i_{j-1}}e_{j-1}, \ \ j=2, ... , n
\eqno{(13.3.2)} $$

It is easy to check that $M$ has complete multiplication by $\w_\goth K$, and its CM-type is $\Phi$.
\medskip
{\bf Remark.} It is possible to prove that $M(\goth K, \Phi)$ is the only t-motive having these properties; we omit the proof.
\medskip
{\bf Proposition 13.3.3.} For $\w_\goth K= \n F_{q^r}[T]$ we have: $M(\goth K, \Phi)'=M(\goth K, \Phi')$.
\medskip
{\bf Proof.} Elements $\tau^je_k$ for $k=1, ... ,n$, $j=0, ..., i_{k+1}-i_{k}-1$ for $k<n$
and $j=0,
..., i_1-i_n+r-1$ for $k=n$ form a basis of $M_{\p[T]}$. Let us
arrange these elements in
the lexicographic order ($\tau^{j_1}e_{k_1}$ precedes to
$\tau^{j_2}e_{k_2}$ if $k_1 <
k_2$) and make a cyclic shift of them by $i_1$ denoting $e_1$ by $f_{i_1+1}$,
$\tau^{i_2-i_1-1}e_1$ by $f_{i_2}$ etc. until
$\tau^{i_1-i_n+r-1}e_n=f_{i_1}$. Formulas
13.3.1, 13.3.2 become
$$\tau(f_i)=f_{i+1} \hbox{ if } i\not\in \{i_1, ..., i_n\} $$
$$\tau(f_i)=(T-\theta)f_{i+1} \hbox{ if } i\in \{i_1, ..., i_n\} $$
($i \mod r$, i.e. $f_{r+1}=f_1$). Formula 1.10.1 shows that in the
dual basis $f'_*$ we
have
$$\tau(f'_i)=f'_{i+1} \hbox{ if } i\in \{i_1, ..., i_n\} $$
$$\tau(f'_i)=(T-\theta)f'_{i+1} \hbox{ if } i\not\in \{i_1, ..., i_n\} $$
which proves the proposition. $\square$
\medskip
{\bf Case $\w_\goth K=\n F_q[T^{1/r}]$, $(r,q)=1$.} In order to define $M(\goth K, \Phi)$ we need more notations. We denote
$\theta^{1/r}$ and $T^{1/r}$ by $\goth s$ and $S$ respectively, and
let $\zeta_r$ be a primitive $r$-th root of 1. Let $\alpha_i$, $i_1 < i_2 < ... < i_n$ and $\Phi$ be the
same as in the case $\w_\goth K= \n F_{q^r}[T]$. We
have $\alpha_i(S)=\zeta_r^iS$. Further, we consider an overring $\p[S,\tau]$ of $\p[T,\tau]$ ($S$ is in the
center of this ring), and we consider the category of modules
over $\p[S,\tau]$ such that the condition 1.9.2 is changed by a
weakened condition 13.3.4 (here $A_{S,0}\in M_n(\p)$ is defined by the formula $Se_*=A_Se_*$, where $A_S\in M_n(\p)[\tau]$, $A_S=\sum_{i=0}^*A_{S,i}\tau^i$):
$$A_{S,0}^r=\theta E_n + N\eqno{(13.3.4)}$$
Let $\bar M$ be a $\p[S,\tau]$-module such that $\dim \bar
M_{\p[S]}=1$, $f_1$ the only element of a basis of $\bar M_{\p[S]}$ and
$$\tau f_1=(S-\zeta_r^{i_1}\goth s)\cdot ... \cdot
(S-\zeta_r^{i_n}\goth s)f_1$$
By definition, $M=M(\goth K, \Phi)$ is the restriction of scalars from
$\p[S,\tau]$ to $\p[T,\tau]$ of $\bar M$. Like in the case $\w_\goth K= \n F_{q^r}[T]$, it is easy to check that $M$ has complete multiplication by $\w_\goth K$ with CM-type $\Phi$, and it is possible to prove that it is the only t-motive having these properties.
\medskip
{\bf Proposition 13.3.5.} For $\w_\goth K= \n F_q[T^{1/r}]$, $(r,q)=1$ we have: $M(\goth K, \Phi)'=M(\goth K, \Phi')$.
\medskip
{\bf Proof.} For $i=1, ..., r$ we denote $f_i=S^{i-1}f_1$. These $f_*=f_*(\Phi)$
form a basis of $M_{\p[T]}$, and the matrix $Q=Q(f_*,\Phi)$ of
multiplication of $\tau$ in this basis has the following description.
We denote by $\sigma_k(\Phi)$ the elementary symmetric polynomial
$\sigma_k(\zeta_r^{i_1}, ... , \zeta_r^{i_n})$.

The first line of $Q$ is
$$\sigma_n(\Phi)\goth s^n \ \ \ \sigma_{n-1}(\Phi)\goth s^{n-1} \ \ \ ...
\ \ \  \sigma_1(\Phi)\goth s \ \ \ 1 \ \ \ 0 \ \ \ ... \ \ \ 0$$
and its $i$-th line is obtained from the first line by 2 operations:

1. Cyclic shift of elements of the first line by $i-1$ positions to the right;

2. Multiplication of the first $i+n-r$ elements of the obtained line by $T$.

We consider another basis $g_*=g_*(\Phi)$ of $M_{\p[T]}$ obtained by
inversion of order of $f_i$, i.e. $g_i=f_{r+1-i}$. The elements of
$Q(g_*)$ are obtained
by reflection of positions of elements of $Q(f_*)$ respectively the
center of the matrix.

The theorem for the present case follows from the formula
$$Q(f_*,\Phi)Q(g_*,\Phi')^t=(T-\theta)E_r$$
whose proof is an elementary exercise: let $\Phi'=\{j_1, ... ,
j_{r-n}\}$; we apply
equality
$$\sigma_k(x_1, ... x_r)=\sum_l \sigma_l(x_{i_1}, ... ,
x_{i_n})\sigma_{k-l}(x_{j_1}, ...
, x_{j_{r-n}})$$
to $1, \zeta_r, ... , \zeta_r^{r-1}$. $\square$
\medskip
{\bf 13.4. Reduction.} Recall notations of 1.16.  Let $L$ be a finite extension of $\n F_q(\theta)$, $\goth p$ a valuation of $L$ over a valuation $P\ne\infty$ of $\n F_q(\theta)$, and we denote $\iota^{-1}(P)\subset \w$ by $\Cal P$. Let $M$ be a t-motive defined over $L$ having a good ordinary reduction $\tilde M$ at $\goth p$ and such that the dual $M'$ exists. According 1.15.1, the $L$-structure on $M'$ is well-defined. We denote by $M_{\Cal P,0}$ the kernel of the reduction map
$M_\Cal P \to \tilde M_\Cal P$. Condition of ordinarity means that  $M_{\Cal P,0}=(\w/\Cal P)^n$.
\medskip
{\bf Conjecture 13.4.1.} For the above $M$, $M'$ we have:
\medskip
$M_{\Cal P,0}$ and $M'_{\Cal P,0}$ are mutually dual with respect to the pairing of Remarks 4.2, 5.1.6 (recall that conjecturally $M'$ also has good ordinary reduction at $\goth p$).
\medskip
{\bf Proof for a particular case:} $M$ is a Drinfeld module, $\Cal P=T$.
\medskip
(1.9.1) for $M$ has a form
$$Te=\theta e+a_1\tau e+... a_{r-1}\tau^{r-1}e+\tau^r e$$
Condition of good ordinary reduction means $a_i\in L$, $\ord_\goth p(a_i)\ge 0$,
$\ord_\goth p a_1=0$. Let $x\in
M_T$, $y\in M'_T$; we can consider $x$ (resp. $y$) as an element of
$\p$ (resp. $\p^{r-1}$) satisfying some polynomial equation(s). Considering
Newton polygon of these polynomials we get immediately (1) for both
$M$, $M'$. Let $y=(y_1, ... ,y_{r-1})$ be the coordinates of $y$; explicit
formula (5.3.5) for the present case has the form
$$<x,y>_M=\Xi(xy_{r-1}^q+x^qy_1+x^{q^2}y_2+ ... +x^{q^{r-1}}y_{r-1})$$
The same consideration of the Newton polygon of the above polynomials
shows that for $x\in M_{T,0}$, $y\in M'_{T,0}$ we have
$\ord_\goth p x$, $\ord_\goth p y_i \ge 1/(q-1)$. Since $\ord_\goth p \Xi = -1/(q-1)$ we get
that $\ord_\goth p (<x,y>_M)>0$ and hence (because $<x,y>_M\in \n F_q$) we have
$<x,y>_M=0$. Dimensions of $M_{T,0}$, $M'_{T,0}$ are
complementary, hence they are mutually dual. $\square$
\medskip
{\bf Remark 13.4.2.} Analogous explicit proof exists for any standard-3 $M$ of
Section 11.8.
\medskip
\medskip
{\bf References}

\nopagebreak
\medskip
[A] Anderson Greg W. $t$-motives. Duke Math. J. Volume 53, Number 2 (1986), 457-502.
\medskip
[BH] Matthias Bornhofen, Urs Hartl, Pure Anderson Motives and Abelian $\tau$-Sheaves. arXiv:0709.2809
\medskip
[F] Faltings, Gerd, Group schemes with strict $\Cal O$-action. Mosc.
Math. J.  2
(2002),  no. 2, 249--279.
\medskip
[G] Goss, David Basic structures of function
field arithmetic.
\medskip
[H] Urs Hartl, Uniformizing the Stacks of Abelian Sheaves.

http://arxiv.org/abs/math.NT/0409341
\medskip
[L1] Logachev, Anderson t-motives are analogs of abelian varieties with multiplication by imaginary quadratic fields.
http://arxiv.org 0907.4712.
\medskip
[L2] Logachev, Reductions of Hecke correspondences on Anderson varieties. In preparation.
\medskip
[L3] Logachev, Lattice map for Anderson t-motives: first approach. http://arxiv.org/pdf/1109.0679.pdf
\medskip
[P] Richard Pink, Hodge structures over function fields. Universit\"at Mannheim.
Preprint. September 17, 1997.
\medskip
[Sh63] Shimura, Goro On analytic families of polarized abelian
varieties and automorphic functions. Annals of Math., 1 (1963), vol.
78, p. 149 -- 192
\medskip
[Sh71]  Shimura, Goro Introduction to the
arithmetic theory of
automorphic functions.
\medskip
[Sh98]  Shimura, Goro Abelian varieties with
complex multiplication
and modular functions. Princeton Mathematical Series, 46.
\medskip
[Tae] Taelman Lenny, Artin t-motifs. J. Number Theory, 129 (2009), 142 - 157
\medskip
[T] Taguchi, Yuichiro A duality for finite $t$-modules.
J. Math. Sci. Univ. Tokyo 2 (1995), no. 3, 563--588.
\medskip
E-mail: logachev94{\@}gmail.com

\enddocument

and $\mu$ a number such that the $\mu$-dual ${M'}^{\mu}$ exists. Let $m$ be the minimal number such that $N^m=0$.
\medskip
{\bf Conjecture 6B1.} 1. $\mu \ge m$, i.e. $N^\mu=0$.
\medskip
2. ${M'}^{\mu}$ is uniformizable. There exists a canonical perfect $\w$-valued pairing between $L_T(M)$ and $L_T({M'}^{\mu})$.
\medskip
This pairing extends by $\p[[T-\theta]]$-linearity to the $\p[[T-\theta]]$-valued pairing between $L_T(M)\underset{\w}\to{\otimes}\p[[T-\theta]]$, $L_T({M'}^{\mu})\underset{\w}\to{\otimes}\p[[T-\theta]]$. We denote it by $<.,.>$.
\medskip
{\bf Conjecture 6B2.} We have:
$$\goth q({M'}^{\mu})=\{x\in L_T({M'}^{\mu})\underset{\w}\to{\otimes}\p[[T-\theta]]  \hbox{ such that }  \forall y\in \goth q(M) $$ $$\hbox{ we have } <x,y>\in (T-\theta)^\mu\p[[T-\theta]]\} $$

Let us consider uniformizable $M$ having $N=0$ (i.e. $m=1$) such that $M'$ --- the dual of $M$ --- exists and has $N'=0$.
\medskip
{\bf Corollary 6B3.} For this case Conjecture 6B2 is true, i.e.
$$\goth q({M'})=\{x\in L_T({M'})\underset{\w}\to{\otimes}\p[[T-\theta]]  \hbox{ such that }  \forall y\in \goth q(M) $$ $$\hbox{ we have } <x,y>\in (T-\theta)\p[[T-\theta]]\} $$
\medskip
This follows immediately from Theorem 5 and equivalence of Definition 2.3 and Property 2.4 (we consider the reduction of the above pairing $<.,.>$ modulo the maximal ideal $ (T-\theta) \p[[T-\theta]]$ of $\p[[T-\theta]]$; this reduction coincides with the pairing of Section 5).

\newpage
{\bf 1.16.4. $\Cal P$-rank of $M$ and of $M'$.} By analogy with the number field case, the $\Cal P$-rank of $M$ is the dimension of the $\Cal P$-torsion points of $E(M)$ over $\w/\Cal P$; it varies from $r-n$ (ordinary $M$) to 0 (completely suresingular $M$).
\medskip
{\bf Conjecture 1.16.5.} $M$ is ordinary $\iff M'$ is ordinary.
\medskip
For standard-3 t-motives  apparently this can be shown by explicit calculations. See 13.4.1 for the proof for the case of Drinfeld modules.
\medskip
{\bf Question 1.16.6.} What are possible values of $\Cal P$-rank of $M$, $\Cal P$-rank of $M'$? Is it true that the pair of numbers ($\Cal P$-rank of $M$, $\Cal P$-rank of $M'$) characterizes completely (in some meaning) the type of $M$? Particularly, whether the $\Cal P$-ranks of other tensor operations of $M$ (exterior powers etc.) are defined completely by ($\Cal P$-rank of $M$, $\Cal P$-rank of $M'$), or not?
\medskip
{\bf Example 1.16.7.} Let us consider $M$ defined by (8.2.2), entries of $A$ belong to $\bar \n F_q$, and let $\Cal P=T$. According Lang's Theorem, $T$-rank of $M$ is $n \iff \det A\ne 0$. An explicit calculation for the case $n=2$ shows that if the entries of $A$ belong to $\n F_q$ then the $T$-rank of $M$ is equal to the quantity of the non-zero eigenvalues of $A$ (if the entries of $A$ do not belong to $\n F_q$ then the formula for the $T$-rank of $M$ is more complicated). The dual $M'$ is defined by the same (8.2.2) with $A$ replaced by $-A^t$; this means that  if the entries of $A$ belong to $\n F_{q}$ then the $T$-rank of $M$ is equal to the $T$-rank of $M'$. It is easy to check that the same is true if the entries of $A$ belong to $\n F_{q^2}$ but not to a larger field:  for example, if $\gamma\in \bar \n F_q-\n F_{q^2} $ then for $A=\left(\matrix \gamma & 1 \\ -\gamma^{q+1}& -\gamma^{q}\endmatrix \right)$ we have $T$-rank of $M$ is 0, $T$-rank of $M'$ is 1.
\medskip
{\bf Example 1.16.8.} For pure non-ordinary standard-3 $M$ having $r=5$, $n=2$ we get easily that the case $T$-rank of $M$ is 2, $T$-rank of $M'$ is 0 cannot be realized; all other possible cases can be realised.
\medskip
{\bf Remark 13.4.3: Case $N\ne0$.} I do not know a definition of ordinarity for $M$ in finite characteristic for the case $N\ne0$; we can expect that these are those $M$ whose $\Cal P$-rank is the maximal possible in some family of these $M$. A version of Conjecture 13.4.1 should hold for this case; instead of the dual $M'$ we should consider the $m$-th dual ${M'}^m$ where $m$ satisfies $M^m=0$ (or $m$ is the minimal number with this property?)
\medskip
{\bf Example 13.4.3.1.} Let $M$ be given by the following modified formula (8.2.2):
$$Te_* =( \theta+N) e_* + A \tau e_* + \tau^2 e_*\eqno{(13.4.3.2)}$$
where $n=2$, $N=\left(\matrix 0&1\\0&0 \endmatrix\right)$, $A= \left(\matrix a_{11}&a_{12} \\a_{21}&a_{22} \endmatrix\right)$. For the case of finite characteristic (i.e. $a_{ij}\in \bar\n F_q$) and $\Cal P=T$ we have $M$ is ordinary iff $a_{21}\ne0$, in this case the $\Cal P$-rank of $M$ is 3. We have $m=2$, ${M'}^2$ is given by the formula
$$Te'_*=( \theta+N') e'_* + A' \tau e'_*$$
where $N'=\left(\matrix 0&-1&0&0&0&0\\0&0&0&0&0&0\\0&0&0&1&0&0\\0&0&0&0&0&0
\\0&0&0&0&0&1\\0&0&0&0&0&0\endmatrix\right)$,  $A'=\left(\matrix 0&0&0&0&1&0 \\0&0&1&0&0&0\\0&0&0&0&0&0\\0&1&-a_{11}&0&-a_{21}&0\\0&0&0&0&0&0
\\1&0&-a_{12}&0&-a_{22}&0 \endmatrix\right)$.

An explicit calculation shows that the $\Cal P$-rank of $M$ is $\le 1$ and it is 1 iff $a_{21}\ne0$. So, an analog of 13.4.1 for $M$ is the following:
\medskip
{\bf Conjecture 13.4.3.3.} Let $M$ be given by 13.4.3.2 in the generic characteristic, and let $M$ have good ordinary reduction at $\Cal P=T$.  We have: $M_{T,0}$ and $({M'}^2)_{T,0}$ are mutually dual with respect to the pairing.
\medskip
This conjecture can be easily checked by explicit calculation.

\enddocument

More exactly, if $\vf: M \to N$ is a map of abelian
t-motives and both ${M_{\goth C}'}^{\mu}$, ${N_{\goth C}'}^{\mu}$ exist then there exists
the map of rational pr\'e-t-motives ${\vf_{\goth C}'}^{\mu}: {N_{\goth C}'}^{\mu} \to
{M_{\goth C}'}^{\mu}$.


$\sigma: \p \to \p$ the Frobenius automorphism of $\p$, i.e. $\sigma(z)=z^q$.

We extend $\sigma$ to $K_C$, resp. $\w_C$ by the formula
$\sigma_{K}(k\otimes z)=k\otimes \sigma(z)$, $z\in\p$, $k\in K$, resp. $k\in
\w$.

: we consider $M\mapsto M^{(1)}$ as a
functor; if $M$ is free $\p[T]$-module and $\alpha: M \to M$ is a map whose matrix in
some basis $\{m_*\}$ is $C$, then the matrix of $\alpha^{(1)}: M^{(1)} \to M^{(1)}$ in the
base $\{m_*\otimes1\}$ is $C^{(1)}$

We have an inclusion $\w \hookrightarrow K_\infty$,
we prolonge $\iota$ to the inclusion $K_\infty \hookrightarrow \p$, and we denote
the image $\iota(K_\infty)\subset \p$ by $K_\infty$ as well.
zarhin@math.psu.edu
Dorogoj Yura,

mozhet byt', Vy znaete chto-nibud' po takim voprosam:

1. Rassmotrim ob'ekty: $C$-prostranstvo $V$, $Z$-reshyotka $L$ v nyom i ermitova forma $H$ na $V$, kotoraya ne >>0 i takaya, chto $im H$ celochislenna na $L$. Kakov kriterij togo, chto eti ob'ekty yavlyayutsya realizaciej motiva?

My znaem, chto esli $H >> 0$, to vsegda. Esli net, to po krajnej mere inogda, naprimer esli $X$ - trifold takoj, chto ego $J^3(X)$ - ne abelevo mnogoobrazie, to $H^3(X)$ - kak raz takogo tipa.

Otsyuda vytekaet takoj vopros

2. My znaem, chto abelevo mnogoobrazie - bolee obshchij ob'ekt, chem jacobian krivoj. Rassmotrim trifold $X$ takoj, chto ego $J^3(X)$ - ne abelevo mnogoobrazie, rassmotrim ego motiv $M=H^3(X)$ i sprosim: motivy kakogo tipa obobshchayut $M$ tak zhe, kak abelevo mnogoobrazie obobshchaet jacobian krivoj?

Mozhet byt', lyudi dumali nad etim 30 - 40 let nazad i uvideli, chto na eti voprosy net razumnyh otvetov?

Prois|hozhdenie etih voprosov takoe. God nazad ya poslal Vam stat'yu ob analogii mezhdu t-motivami s $N=0$ i abelevymi mnogoobraziyami s umnozheniem na mnimoe kvadratichnoe pole. Vopros: a kakov analog t-motivov s $n \ne 0$? Mozhet byt', eto motiv? My legko mozhem opredelit' vysheupomyanutye $V$, $L$ i $H$ dlya vneshnej stepeni abelevyh mnogoobrazij s umnozheniem na mnimoe kvadratichnoe pole. Esli by my znali, chto oni - realizaciya motiva, to kak raz eto i byl by motiv, sootvetstvuyushchij t-motivu s $n \ne 0$.

Yours, Dima

 We consider the case $P\ne \infty$, i.e. we can identify $P$ with an irreducible polynomial in $\n Z_\infty$ (up to a scalar factor).

C-MATRICY                 POZOR !

{\bf 8.2. Infinitesimal lattice map and C-lattices.} For $r=2n$, \ $n>1$ the degree of the lattice map is infinite, hence in order to get a 1 - 1 correspondence it is necessary to introduce a new object called a C-lattice. There is a C-lattice map from a neighbourhood of a t-motive $M_0$ to a neighbourhood of a C-lattice $L_{C,0}$ (see below for the definition of $M_0$, $L_{C,0}$) which is 1 - 1.

For the case $(r,n)=1$ I do not see an inicial t-motive $M_0$ such that in its neighbourhood the degree of the lattice map is not 1 (see Remark 8.2.14 for more details), hence we can ask whether in this case the answer to 8.1.2 is yes. I think that this is few likely.
\medskip
{\bf Definition 8.2.1.} A C-lattice is a pair $(L, e_*)$ where
\medskip
$L \subset V=\p^n$ is a lattice and
\medskip
$e_1, ... , e_r$ is a $\n Z_\infty$-basis of $L$.
\medskip
Further, two such pairs $(L; e_1, ... , e_r)$ and $(L'; e'_1, ... ,
e'_r)$ are called equivalent if there exists a $\p$-linear map $\psi:
V\to V$ such that $\psi(L; e_1, ... , e_r) = (L'; e'_1, ... , e'_r)$,
or $L=L'$ and the matrix of the change of basis from $e_1, ... , e_r$
to $e'_1, ... , e'_r$ (which a priori belongs to $GL_r(\n Z_\infty)$ ) belongs to $GL_r(\n F_q)$.
\medskip
The functor of forgetting the basis from C-lattices to lattices is
denoted by $\goth i$. The notion of a Siegel matrix for a C-lattice is
the same as for a lattice.

We consider infinitesimal degree of the lattice map in a neighbourhood
of some distinguished (having many endomorphisms) t-motive $M_0$. We
consider the case $M_0=\goth C_2^{\oplus n}$ (recall that $\goth C_2$ is the Carlitz module over $\n F_{q^2}[T]$; in 7.7 $M_0$ is denoted by $M$).
\medskip
{\bf Remark.} It is easy to see that $M_0$ is the (only) t-motive with complete multiplication by $\n F_{q^{2n}}[\theta]$ with CM-type $Id, \fr^2, \fr^4, ..., \fr^{2n-2}$, see 12.5.3 and 13.3, first case.
\medskip
Let $\omega \in \n F_{q^2} - \n F_q$ be a
fixed element. A Siegel matrix of $M_0$ is $\omega E_n$. We denote the lattice (resp. the C-lattice) corresponding to $\omega E_n$ by $L_0$ (resp. $L_{C,0}$). $L_0$ is the lattice of $M_0$. We consider 5 sets $S_1, ... , S_5$:
\medskip
$S_1$. The set of $n \times n$ matrices $A$.
\medskip
$S_2$. The set of t-motives $M$ given by the equation (see 1.9.1)

$$Te_* = \theta e_* + A \tau e_* + \tau^2 e_*\eqno{(8.2.2)}$$ where $A
\in S_1$, $e_* = (e_1,...,e_n)^t$ (notations of 1.9.1: $\goth A_1=A$, $\goth A_2=E_n$).
\medskip
$S_3$. The set of $n \times n$ Siegel matrices $Z$.
\medskip
$S_4$. The set of C-lattices of rank $r=2n$ in $C_\infty ^n$.
\medskip
$S_5$. The set of lattices of rank $r=2n$ in $C_\infty ^n$.
\medskip
We consider initial elements 0, $M_0$, $\omega E_n$, $L_{C,0}$, $L_0$ of $S_1, ...,S_5$ respectively and open neighbourhoods $U_i\subset S_i$ of these initial elements.
\medskip
{\bf Proposition 8.2.3.} There exist neighbourhoods $U_2$, $U_4$, $U_5$ such that

(a) The restriction of $\goth i$ to $U_4$ gives us an epimorphism $U_4\to U_5$.

(b) there exists a 1 -- 1 map $\mu_{24}$ from $U_2$ to $U_4$ such that $\mu_{25}:=\goth i\circ\mu_{24}$ (see the below diagram 8.2.4) is the lattice map from uniformizable t-motives to lattices. Particularly, for $n>1$ the fibre of $\mu_{25}$ is discrete infinite.
\medskip
{\bf Proof.} It is sufficient to prove that there exists a commutative diagram $$\matrix U_1 &
\overset{\mu_{12}}\to{\to} & U_2 && \\ \\ \mu_{13} \downarrow &
&\mu_{24}\downarrow &
\overset{\mu_{25}}\to{\searrow} &\\ \\ U_3 &
\overset{\mu_{34}}\to{\to} & U_4 & \overset{\goth i}\to{\to} & U_5
\endmatrix \eqno{(8.2.4)}$$ where $\mu_{12}(A)$ is the t-motive defined by
8.2.2, $\mu_{34}(Z)$ is the C-lattice corresponding to a Siegel matrix
$Z$ and $\mu_{13}$ is defined as
follows. We identify $\Lie(\goth C_2)$ with $\p$ and hence $\Lie(M_0)$ with $\p^n$.
We consider the following basis
$l_{0,1}, ... , l_{0,2n}$ of $L_0\subset \Lie(M_0)=\p^n$: $l_{0,i}=(0,...,0,1,0,...,0)$ (1 at the $i$-th place), $l_{0,n+i}=\omega l_{0,i}$, $i=1,...,n$. The Siegel
matrix of this basis is $\omega E_n$. Let $A\in U_1$, $M$ its $\mu_{12}$-image and $L$ its
lattice, i.e. its $\mu_{25}\circ\mu_{12}$-image. For any $l_0\in L_0$
there exists a well-defined $l\in L$ which is close to $l_0$ (because
entries of $A$ are near 0). So, we consider a basis $l_{1}, ... ,
l_{2n}$ of $L$ where any $l_i$ is near the corresponding $l_{0,i}$.
The Siegel matrix corresponding to $l_{1}, ... , l_{2n}$ is exactly
the $\mu_{13}$-image of $M$. By definitions, the outer quadrangle is
commutative.

We denote by $d_{\alpha\beta}$ the degree of
$\mu_{\alpha\beta}$ at a
generic point near the initial element of $U_\alpha$. We must prove that there exists a map $\mu_{24}$ preserving
commutativity, and that $d_{24}=1$.
\medskip
{\bf Lemma 8.2.5.} $d_{13}=1$.
\medskip
{\bf Proof.} We consider maps of the diagram 5.1.3 for $M$ given by 8.2.2. For $Z\in E=\p^n$ we have $m_T(Z)=\theta Z + A Z^{(1)} +
Z^{(2)}$, and for $Z\in \Lie(M)=\p^n$ we have $\Exp(Z)=\Exp_A(Z)=\sum_{i=0}^\infty C_iZ^{(i)}$ where
$C_i=C_i(A)$, $C_0=1$.

For the reader's convenience, we consider only the case $n=1$ (the
general case does not require any new ideas), hence $Z$, $A$ will be
denoted by $z$, $a$ respectively. We denote
$\theta_{ij}=\theta^{q^i}-\theta^{q^j}$. Recall that the exponent for $\goth C_2$ has the form
$$\Exp_0(z)=z+\frac{1}{\theta_{20}}z^{q^2}+\frac{1}{\theta_{42}\theta_{40}}z^{q^4}+...\eqno{(8.2.6)}$$
($C_{2i}(0)=\frac{1}{\prod_{j=0}^{i-1}\theta_{2i,2j}}$). We denote by
$y_0\in \p$ a nearest-to-zero root to $\Exp_0(z)=0$ (this is $\xi$ of $\goth C_2$ in notations of [G]). It is defined up to
multiplication by elements of $\n F_{q^2}^*$, and it generates over $\n F_{q^2}[\theta]$ the
lattice of the Carlitz module $\goth C_2$. We fix one such $y_0$.
If $a$ is sufficiently small then there is the only root to
$\Exp_a(z)=0$ near $y_{0}$, and there is the only root to
$\Exp_a(z)=0$ near $\omega y_{0}$. We denote these roots
by $z=z(a)$, $z'=z'(a)$ respectively, and we denote $z=y_{0}+\delta$,
$z'=\omega y_{0}+\delta'$. $\delta$ (resp. $\delta'$) is a root to the
power series
$$\sum_{i,j=0}^\infty d_{ij}a^i\delta^j=0 \eqno{(8.2.7)}$$
$$\hbox{resp. }\sum_{i,j=0}^\infty d'_{ij}a^i{\delta'}^j=0 \eqno{(8.2.7')}$$
where
$$d_{00}=0, \ \
d_{10}=\frac{y_{0}^{q}}{\theta_{10}}+\frac{y_{0}^{q^3}}{\theta_{31}\theta_{30}}
+ \frac{y_{0}^{q^5}}{\theta_{53}\theta_{51}\theta_{50}}+
\frac{y_{0}^{q^7}}{\theta_{75}\theta_{73}\theta_{71}\theta_{70}}+..., \ \ d_{01}=1$$
$$d'_{00}=0, \ \ d'_{10}=\frac{(\omega
y_{0})^{q}}{\theta_{10}}+\frac{(\omega
y_{0})^{q^3}}{\theta_{31}\theta_{30}} + \frac{(\omega
y_{0})^{q^5}}{\theta_{53}\theta_{51}\theta_{50}}+
\frac{(\omega y_{0})^{q^7}}{\theta_{75}\theta_{73}\theta_{71}\theta_{70}}+..., \ \ d'_{01}=1$$
Moreover, $\delta$, $\delta'$ are the nearest-to-0 roots to (8.2.7),
$(8.2.7')$ respectively. It is clear that $d_{10}, d'_{10}\ne 0$. This
means that the approximate value of $\delta$, $\delta'$ is
$$-d_{10}a, \ \ -d'_{10}a\eqno{(8.2.8)}$$
respectively. Exactly, both $\delta$, $\delta'$ are power series in
$a$ whose first term is given by (8.2.8). This means that
$$z=y_0-d_{10}a+\sum_{i=2}^\infty k_i
a^i$$
$$z'=\omega y_0-d'_{10}a+\sum_{i=2}^\infty k'_i
a^i$$ It is easy to see that $d'_{10}\ne \omega d_{10}$,
hence the Siegel matrix $\goth Z=z^{-1}z'$ (which is a number because $n=1$) is given by the formula $$\goth Z=\omega
+ \sum_{i=1}^\infty l_i a^i\eqno{(8.2.9)}$$ and $l_1\ne 0$. Since (for $n=1$)
$d_{13}$ is the minimal $i$ such that $l_i\ne 0$ we get that
$d_{13}=1$.
\medskip
For $n > 1$ the calculation is the same (this is the main part of the proof, because for the case $n=1$ the result is known). Analog of (8.2.9) is $$\goth Z=\omega E_n + l_1 A + P_{\ge 2}(A)\eqno{(8.2.10)}$$ where $P_{\ge 2}(A)$ is a power series of entries of $A$ such that all its terms have degree $\ge 2$.
Condition $l_1\ne 0$ implies $d_{13}=1$. $\square$
\medskip
In order to simplify notations, we identify until the end of the proof $L_T(M)$ and $L(M)$ via $\goth a$, and $\n Z_\infty$ and $\w$ via $\iota$. We denote the monodromy group of $\mu_{12}$, resp. $\mu_{35}:=\goth i \circ \mu_{34}$ by $\Cal M_{12}$, resp. $\Cal M_{35}$. We have $\Cal M_{12}=GL_n(\n F_{q^2})/\n F_q^*$. Really, the automorphism group of $M_0$ is $GL_n(\n F_{q^2}[T])$: an element $g$ of this group acts on the basis $e_*$ of (8.2.2) if $A=0$. But if $A$ is a generic matrix, then for $g\in GL_n(\n F_{q^2}[T]) - GL_n(\n F_{q^2})$ the result of the action of $g$ on (8.2.2) becomes a more complicated equation: terms having higher powers of $\tau$ appear, so these $g$ do not belong to $\Cal M_{12}$. Factorization by $\n F_q^*$ is obvious. Further, obviously $\Cal M_{35}=\{\gamma\in PGL_{2n}(\n Z_\infty)|\gamma(\omega E_n)=\omega E_n\}$. The outer quadrangle of 8.2.4 defines a map from $\Cal M_{12}$ to $\Cal M_{35}$ which we denote by $\alpha$.
\medskip
{\bf Lemma 8.2.11.} $\alpha$ is injective, and $\im \alpha= \{\gamma\in PGL_{2n}(\n F_{q})|\gamma(\omega E_n)=\omega E_n\}=\Cal M_{35} \cap PGL_{2n}(\n F_{q}) \subset PGL_{2n}(\n Z_\infty)$.
\medskip
{\bf Proof.} $\alpha$ is defined by the condition: for any $A\in S_1$, $\gamma \in \Cal M_{12}$ we have
$$\mu_{13}(\gamma(A))=(\alpha(\gamma))(\mu_{13}(A))\eqno{(8.2.12)}$$ The explicit formula for $\alpha$ is the following. For odd $q$ we fix $\omega$ satisfying $\omega^2\in \n F_q^*$ (it is an easy exercise to find analog of the below formula for even $q$). Let $\gamma=U+\omega V$, \ \ $U,V\in GL_n(\n F_q)$ (we consider cleary a representative of $\gamma\in GL_n(\n F_{q^2})/\n F_q^*$ in $GL_n(\n F_{q^2})$). Then $$\alpha(\gamma)=\left(\matrix U&-\omega^2V\\ -V&U \endmatrix\right)^{-1}\eqno{(8.2.13)}$$ ($\alpha$ is an antihomomorphism, because the functor of lattice is contravariant). It is checked immediately that (8.2.12) holds. We see that $\im \alpha= \{\gamma\in PGL_{2n}(\n F_{q})|\gamma(\omega E_n)=\omega E_n\}$. $\square$ 
\medskip
Proposition 8.2.3 follows immediately from these lemmas. $\square$
\medskip
{\bf Remark 8.2.14.} We see that if $r=2n$, $n>1$ then the lattice map $\mu_{25}$ is not a local isomorphism near $M_0$. The origin of this phenomenon is reducibility of $M_0$ which implies that the monodromy group of $\mu_{35}$ is much bigger then the one of $\mu_{12}$. For other values of $r$, $n$ a natural analog of $M_0$ is a t-motive with complete multiplication. Apparently if $(r,n)=1$ then all pure t-motives with complete multiplication are irreducible (example: CM-field is $\n F_{q^r}[\theta]$), and analog of Proposition 8.2.3 for this case shows that $\mu_{25}$ is a local isomorphism near this t-motive.
\medskip
{\bf Remark 8.2.15.} From the first sight, for $n=1$ Proposition 8.2.3 contradicts to a result of Drinfeld about 1 -- 1 correspondence between Drinfeld modules and lattices. Really, there is no contradiction: if $n=1$ then $\im \alpha=\Cal M_{35}$, and --- although $\goth i: S_4 \to S_5$ clearly is not an isomorphism --- its restriction to $U_4$ is an isomorphism $U_4 \to U_5$.
\medskip
{\bf 8.3. Duality of C-lattices.} We have no analog of 2.2 for C-lattices, so we use an analog of 3.2 as a definition of duality. Namely, if $Z$ is a Siegel matrix of a C-lattice $(L, e_*)$ then its dual $(L, e_*)'$ is a C-lattice whose Siegel matrix is $-Z^t$ (or $Z^t$ which is the same, but more convenient for further calculations, because $(-\omega E_n)^t\ne \omega E_n$). Equality 3.8.2 shows that this notion is well-defined (entries of $A$, $B$, $C$, $D$ belong to $\n F_q$).

An analog of Theorem 5 holds for C-lattices:
\medskip
{\bf Theorem 8.3.1.} Let $M\in U_2$. Then $\mu_{24}(M')=\mu_{24}(M)'$.
\medskip
{\bf Proof} is completely analogous to the proof of Theorem 5, so we omit it. Alternatively, we can show that the exact form of (8.2.10) is $$\goth Z=\omega E_n + \sum_{k=1}^\infty \sum_{d_1,...,d_k}l_{d_1,...,d_k}A^{(d_1)}\cdot...\cdot A^{(d_k)}$$ where coefficients $l_{d_1,...,d_k}$ satisfy $$l_{d_1,...,d_k}=l_{d_k,...,d_1}$$ This obviously implies the theorem. $\square$
\medskip

\newpage
{\bf Theorem.} ${M'}^{m}$ is uniformizable, and $\underline{H}({M'}^{m})= {\underline{H}(M)'}^{m}$.
\medskip
We shall need several elementary lemmas.
Let $w$ be a number such that $N^w=0$. We define numbers $k_i=k_i(M)$, $i=2, \dots,w+1$, as follows:
$$k_{i}:=\dim \Ker N^{i-1}/ \Ker N^{i-2}- \dim \Ker N^{i}/ \Ker N^{i-1} $$
$$=\dim \im N^{i-2} / \im N^{i-1} - \dim \im N^{i-1} / \im N^{i}\eqno{(9.3)}$$
Equivalently, let $n=d_1+...+d_\al$, where $d_1\ge d_2\ge...\ge d_\al>0$, be a partition of $n$ corresponding to the Jordan form of $N$, i.e. the Jordan form of $N$ consists of $\al$ 0-Jordan blocks of sizes $d_1, d_2,...,d_\al$. We have $w\ge d_1$, $\al\le r$. We shall call a 0-partition of length $\g r$ of a number $\g n$ a representation of $\g n$ as a sum
$$\g n = \g d_1+\g d_2+...+\g d_{\g r}$$ where $\g d_i\in \n Z$ and $\g d_1\ge\g d_2\ge...\ge\g d_{\g r}\ge0$, i.e. a 0-partition is a partition plus several zeroes at its end. We extend the partition $n=d_1+...+d_\al$ to a 0-partition of length $r$ denoted by $\g p=\g p(M)$:
$n=d_1+...+d_\al+d_{\al+1}+...+d_r$ where $d_{\al+1}=...=d_r=0$. Let $n=c_1+...+c_{d_1}+c_{d_1+1}+...+c_w$ be the 0-partition of length $w$ dual to $\g p$ (the definition of the dual 0-partition of a given length is clear). We have $\al=c_1\ge c_2\ge...\ge c_{d_1}>0$, $c_{d_1+1}=...=c_w=0$. We have $\dim \Ker N^{i}=c_1+...+c_i$, hence $k_i=c_{i-1}-c_i\ge0$ (for $i=w+1$ we let $c_{w+1}=0$), and
$$n=\sum_{i=1}^{w}ik_{i+1}\eqno{(9.4)}$$ We have $\al=c_1=\sum_{i=2}^{w+1}k_{i}$, hence $r\ge \sum_{i=2}^{w+1}k_{i}$. We let $k_1:=r- \sum_{i=2}^{w+1}k_{i}$.
\medskip
According ??, elements $T^il_j$, $i=0,1,\dots$, $j=1,\dots,r$, generate $\Lie(M)$ as a $\p$-vector space. Hence, elements $N^il_j$, $i=0,\dots,w-1$, $j=1,\dots,r$, also generate $\Lie(M)$ as a $\p$-vector space. Using this fact we arrange elements $l_1,\dots,l_r$ in $w+1$ segments as follows. First, elements $N^{w-1}l_j$, $j=1,\dots,r$, generate $N^{w-1}\Lie(M)$ as a $\p$-vector space. Its dimension is $k_{w+1}$, hence (first step of a process) we can choose $k_{w+1}$ elements from $l_1,\dots,l_r$ (we denote them by $l_{w+1,1}, \dots, l_{w+1,k_{w+1}}$ respectively) such that
\medskip
(9.6) $N^{w-1}(l_{w+1,i})$, $i=1, \dots, k_{w+1}$, form a $\p$-basis of $N^{w-1}\Lie(M)$.
\medskip
Further (second step), elements $N^{w-2}l_j$, $N^{w-1}l_j$, $j=1,\dots,r$, generate $N^{w-2}\Lie(M)$ as a $\p$-vector space. Elements $N^{w-2}(l_{w+1,i})$, $N^{w-1}(l_{w+1,i})$, $i=1, \dots, k_{w+1}$, are linearly independent over $\p$. Really, let $$\sum_{\al_2=1}^{k_{w+1}}c_{\al_2}N^{w-2}(l_{w+1,\al_2})+ \sum_{\al_1=1}^{k_{w+1}}c_{\al_1}N^{w-1}(l_{w+1,\al_1})=0\eqno{(9.7)}$$
be a non-trivial dependence relation. Applying $N$ to (9.7) we get $$\sum_{\al_2=1}^{k_{w+1}}c_{\al_2}N^{w-1}(l_{w+1,\al_2})=0$$ that contradicts (9.6). Hence, all $c_{\al_2}$ are 0 that again contradicts (9.6).
\medskip
We have $\dim N^{w-1}\Lie(M)=2k_{w+1}+k_{w}$, this follows immediately from (9.3). Hence, we get:
\medskip
Among $l_1,\dots,l_r$ there exist $k_w$ elements (we denote them by $l_{w,1}, \dots, l_{w,k_{w}}$ respectively) such that their intersection with $l_{w+1,1}, \dots, l_{w+1,k_{w+1}}$ is empty and such that
\medskip
(9.8) $N^{w-2}(l_{w,i})$, $i=1, \dots, k_{w}$, $N^{w-2}(l_{w+1,i})$, $i=1, \dots, k_{w+1}$, $N^{w-1}(l_{w+1,i})$, $i=1, \dots, k_{w+1}$, form a $\p$-basis of $N^{w-2}\Lie(M)$.
\medskip
Third step of the process: elements $N^{w-3}l_j$, $N^{w-2}l_j$, $N^{w-1}l_j$, $j=1,\dots,r$, generate $N^{w-3}\Lie(M)$ as a $\p$-vector space. Elements $N^{w-3}(l_{w,i})$, $N^{w-2}(l_{w,i})$, $i=1, \dots, k_{w}$, and $N^{w-3}(l_{w+1,i})$, $N^{w-2}(l_{w+1,i})$, $N^{w-1}(l_{w+1,i})$, $i=1, \dots, k_{w+1}$, are linearly independent over $\p$ (the proof is as the one above for the second step). Hence, we get:
\medskip
Among $l_1,\dots,l_r$ there exist $k_{w-1}$ elements (we denote them by $l_{w-1,1}, \dots,$ $ l_{w-1,k_{w-1}}$ respectively) such that their intersection with $l_{w,1}, \dots, l_{w,k_{w}}$, $l_{w+1,1}, \dots, l_{w+1,k_{w+1}}$ is empty and such that
\medskip
(9.9) $N^{w-3}(l_{w-1,i})$, $i=1, \dots, k_{w-1}$, $N^{w-3}(l_{w,i})$, $N^{w-2}(l_{w,i})$, $i=1, \dots, k_{w}$, and $N^{w-3}(l_{w+1,i})$, $N^{w-2}(l_{w+1,i})$, $N^{w-1}(l_{w+1,i})$, $i=1, \dots, k_{w+1}$, form a $\p$-basis of $N^{w-3}\Lie(M)$.
\medskip
Continuing this process, we represent $r$ as an ordered partition
$$r=k_1+...+k_{w+1}\eqno{(9.10)}$$
(recall that some $k_*$ can be 0) and we represent the set $\{l_1,\dots,l_{r}\}$ as a union of segments $$\{l_1,\dots,l_{r}\}=\bigcup_{u=1}^{w+1}\{l_{u1}, \dots, l_{uk_u}\}\eqno{(9.11)}$$ (the union is ordered and disjoint) such that $\forall \ u=0,\dots,w-1$ we have:
\medskip
$(9.12)$ A $\p$-basis of $N^u\Lie(M)$ is formed by elements $N^\al(l_{\be\ga})$, where $\al\in[u,\dots, w-1]$, $\be\in[\al+2,\dots, w+1]$, $\ga\in [1,\dots, k_\be]$.
\medskip
This implies that for any
$$u\in [1,w], \ \ z\in [u-1, w-1], \ \  y\in [z+2,w+1], \ \ \eqno{(9.13)}$$
$$v=u-1\eqno{(9.14)}$$
there exist matrices $S_{uvyz}$ of size $k_u\times k_y$ with entries in $\p$ (analogues of the Siegel matrix for the case $w=1$) such that $\forall \ u=1,\dots, w, \forall \ i=1,\dots, k_u$ the following holds:
$$N^{u-1}l_{ui}=-\sum_{z=u-1}^{w-1}\sum_{y=z+2}^{w+1} \sum_{j=1}^{k_y}(S_{uvyz})_{ij}N^zl_{yj}\eqno{(9.15)}$$
(if some $k_*$ are 0 then the corresponding $S_{****}$ do not exist).
\medskip
To simplify formulas, below for any $\al$ we consider $\hat l_\al:=l_{\al*}$ as matrix columns. (9.15) becomes a matrix equality
$$N^{u-1}\hat l_{u}=-\sum_{z=u-1}^{w-1}\sum_{y=z+2}^{w+1}S_{uvyz}N^z\hat l_{y}\eqno{(9.16)}$$

{\bf Remark 9.17.} Since always $v=u-1$, really matrices $S_{uvyz}$ depend on 3 parameters $u,y,z$. Number $v$ indicates the exponent of $N$ in the left hand side of (9.15), by analogy with $z$ which indicates the exponent of $N$ in the right hand side of (9.15). This notation is convenient to define a symmetry between $A_*$ and $P_*$, see below.
\medskip
{\bf 9.18.} Example for $w=3$, $u=1$:
$$\matrix l_{1i}=-(\sum_{j=1}^{k_2}(S_{1020})_{ij}l_{2j} &+& \sum_{j=1}^{k_3}(S_{1030})_{ij}l_{3j} &+& \sum_{j=1}^{k_4}(S_{1040})_{ij}l_{4j} \\ \\
&+& \sum_{j=1}^{k_3}(S_{1031})_{ij}N(l_{3j}) &+& \sum_{j=1}^{k_4}(S_{1041})_{ij} N(l_{4j}) \\ \\ &&&+& \sum_{j=1}^{k_4}(S_{1042})_{ij} N^2(l_{4j}) \ )\endmatrix $$
(terms of a fixed column of this formula correspond to a fixed $y$ and different $z$ of (9.15), and terms of a fixed row of this formula correspond to a fixed $z$ and different $y$ of (9.15) ).
\medskip
Applying powers of $N$ to (9.16), for any $v\in [0,\dots, w-1]$, $u\in [1,\dots, w+1]$ we can represent $N^{v}(\hat l_{u})$ as a linear combination of $N^z(\hat l_{y})$ where for a fixed $v$ the numbers $z, \ y$ satisfy $$z\in [v,\dots,w-1], \ \ y\in [z+2,\dots,w+1], \ \ \eqno{(9.19)}$$ Namely, there exist polynomials in $S_{****}$ denoted by $P_{uvyz}$ such that (matrix notations)

$$N^{v}\hat l_{u}=-\sum_{z=v}^{w-1}\sum_{y=z+2}^{w+1}P_{uvyz}N^z\hat l_{y}\eqno{(9.20)}$$
Clearly for $v=u-1$ we have $P_{uvyz}=S_{uvyz}$.

{\bf 9.21.} The domain $v\ge u-1 \ \ \wedge \ \ \{z, \ y$ satisfy (9.19)\} is called the non-trivial domain of the definition of $P_{****}$.

For $v< u-1$ (trivial domain) we have:

$$ P_{u,v,y,z}=-1, \hbox{ resp. } P_{u,v,y,z}=0  \eqno{(9.22)}$$ for $y$, $z$ satisfying (9.19), $(y,z)=(u,v)$, resp. $(y,z)\ne(u,v)$.

\medskip
{\bf 9.23.} Example for $w=3$:

\medskip

$N^2\hat l_2=(S_{2131}S_{3242}-S_{2141})N^2\hat l_4$, i.e. $P_{2242}=-S_{2131}S_{3242}+S_{2141}$;

\medskip

$N^2\hat l_1=(-S_{1020}S_{2131}S_{3242}+S_{1020}S_{2141} +S_{1030}S_{3242}-S_{1040})N^2\hat l_4$, i.e.

\medskip

$P_{1242}=S_{1020}S_{2131}S_{3242}-S_{1020}S_{2141} -S_{1030}S_{3242}+S_{1040}$;

\medskip

$N\hat l_1=(S_{1020}S_{2131}-S_{1030})N\hat l_3+(S_{1020}S_{2141}-S_{1040})N\hat l_4+$

\medskip

$+(S_{1020}S_{2142}+S_{1031}S_{3242}-S_{1041})N^2\hat l_4$, i.e.

\medskip

$P_{1131}=-S_{1020}S_{2131}+S_{1030}$, $P_{1141}=-S_{1020}S_{2141}+S_{1040}$,

\medskip

$P_{1142}=-S_{1020}S_{2142}-S_{1031}S_{3242}+S_{1041}$.
\medskip
{\bf Remark 9.24.} Although we do not need this fact, let us give a formula for some $P_{****}$. Let us define a block unitriangular matrix $\g S$ whose $(i,j)$-th block is $S_{i,i-1,j,i-1}$ for $j>i$, $I_{k_i}$ for $i=j$ and 0 for $j<i$. Further, we define a block unitriangular matrix $\g P$ whose $(i,j)$-th block is $-P_{i,j-2,j,j-2}$ for $j>i$, $I_{k_i}$ for $i=j$ and 0 for $j<i$. We have $\g P=\g S^{-1}$ (see the the below propositions).
\medskip
Some $P_{****}$ that enter in the below formula for $\bar \g D$ are not of the form of the elements of the inverse unitriangular matrix, for example $P_{1142}$, $w=3$.
\medskip
For the proof of Lemmas 9.32, 23, we need
\medskip
{\bf Lemma 9.25.} For all $i,\ j,\ \psi,\ \xi$ (domain?)
$$(\sum_{\be=0}^{j+\xi-w} \sum_{\al=i+\be}^{w+1-j+\be} S_{i-1,i-2,\al,i-2+\be} P_{\al,w-j+\be,\psi,\xi})-$$ $$-S_{i-1,i-2,\psi,i-2+\xi+j-w}+P_{i-1,w-j,\psi,\xi}=0\eqno{(9.25.1)}$$
(A recurrent formula for $P_{****}$).
\medskip
{\bf Proof.} First, we rewrite (9.16) for $u=i-1$:

$$N^{i-2}\hat l_{i-1}=-\sum_{z=i-2}^{w-1}\sum_{y=z+2}^{w+1}S_{i-1,i-2,y,z}N^z\hat l_{y}\eqno{(9.25.2)}$$

Now, for any $$j=1,\dots,w-i+2\eqno{(9.25.3)}$$ we apply $N^{w-i+2-j}$ to (9.25.2):

$$N^{w-j}\hat l_{i-1}=-\sum_{z=i-2}^{i+j-3}\sum_{y=z+2}^{w+1}S_{i-1,i-2,y,z}N^{z+w-i+2-j}\hat l_{y}\eqno{(9.25.4)}$$

(since $N^w=0$, we get that $z\le i+j-3$ ).

We change a summation variable: $z=i-2+\be$, and $y \to \al$, we get

$$N^{w-j}\hat l_{i-1}=-\sum_{\be=0}^{j-1}\sum_{\al=i+\be}^{w+1}S_{i-1,i-2,\al,i-2+\be}N^{w+\be-j}\hat l_{\al}\eqno{(9.25.5)}$$

Now we use (9.20), we make the following variable change:

$$u \to \al\ \ \ \ \ y\to \psi$$
$$v\to w-j+\be\ \ \ \ \ z\to \xi$$
we get
$$N^{w-j+\be}\hat l_\al=-\sum_{\xi=w-j+\be}^{w-1} \sum_{\psi=\xi+2}^{w+1} P_{\al,w-j+\be,\psi,\xi}N^{\xi}\hat l_\psi\eqno{(9.20a)}$$

We substitute (9.20a) in (9.25.5):

$$N^{w-j}\hat l_{i-1}=\sum_{\be=0}^{j-1}\sum_{\al=i+\be}^{w+1} \sum_{\xi=w-j+\be}^{w-1}\sum_{\psi=\xi+2}^{w+1} S_{i-1,i-2,\al,i-2+\be} P_{\al,w-j+\be,\psi,\xi}N^{\xi}\hat l_\psi\eqno{(9.25.6)}$$

We change the order of summation in (9.25.6):

$$N^{w-j}\hat l_{i-1}=\sum_{\xi=w-j}^{w-1} \sum_{\psi=\xi+2}^{w+1} (\sum_{\be=0}^{j+\xi-w} \sum_{\al=i+\be}^{w+1} S_{i-1,i-2,\al,i-2+\be} P_{\al,w-j+\be,\psi,\xi}) N^{\xi}\hat l_\psi\eqno{(9.25.7)}$$

We rewrite (9.20) making changes:

$$u\to i-1\ \ \ \ \ \ y \to \psi$$
$$v\to w-j\ \ \ \ \ \  z\to \xi$$
we get
$$N^{w-j}\hat l_{i-1}=-\sum_{\xi=w-j}^{w-1}\sum_{\psi=\xi+2}^{w+1} P_{i-1,w-j,\psi,\xi}N^\xi \hat l_{\psi}\eqno{(9.25.8)}$$

For $\psi\ge \xi+2$ elements $N^\xi l_{\psi i}$, $i=1,\dots,k_\psi$, are linearly independent over $\p$. Hence, (9.25.7), (9.25.8) imply

$$P_{i-1,w-j,\psi,\xi}=-\sum_{\be=0}^{j+\xi-w} \sum_{\al=i+\be}^{w+1} S_{i-1,i-2,\al,i-2+\be} P_{\al,w-j+\be,\psi,\xi}\eqno{(9.25.9)}$$

Here the domain of $\xi$, $\psi$ is:

$$\xi\in [w-j,\dots,w-1], \ \ \ \psi\in[\xi+2,\dots, w+1]$$
Taking into consideration (9.22) we can rewrite (9.25.9) as follows:

$$P_{i-1,w-j,\psi,\xi}=-(\sum_{\be=0}^{j+\xi-w} \sum_{\al=i+\be}^{w+1-j+\be} S_{i-1,i-2,\al,i-2+\be} P_{\al,w-j+\be,\psi,\xi})+$$ $$+S_{i-1,i-2,\psi,i-2+\xi+j-w}\eqno{(9.25.10)}$$ with the same domain of $\xi$, $\psi$. This is (9.25.1). Because of (9.25.3), this formula is valid for $w-j\ge i-2$ (the non-trivial case of the definition of $P_{****}$). $\square$
\medskip
Let us consider a symmetry $\g s: \n Z^4\to\n Z^4$ defined as follows: $\g s(\al,\be,\ga,\de)=(w+2-\ga,w-1-\de,w+2-\al,w-1-\be)$.
\medskip
{\bf Remark 9.26.} $\g s$ has the following geometric interpretation. Let us consider a matrix $NL$ whose $(i,j)$-th entry is a symbol $N^{i-1}\hat l_{j}$. We interpret a quadruple $(\al,\be,\ga,\de)$ as a vector from $N^\be \hat l_\al$ to $N^\de \hat l_\ga$ in $NL$. $\g s$ is the reflection of this vector with respect to the center of $NL$ and the inversion of its direction.
\medskip
{\bf Definition 9.27.} $\bar S_{uvyz}:=-P_{\g s(uvyz)}^t$ (defined if $P_{\g s(uvyz)}$ has meaning).
\medskip
We shall consider block $r\times r$-matrices having the following block structure: their block size is $(w+1)\times(w+1)$, quantities of columns in blocks are $k_{w+1},k_w,\dots,k_1$ (counting from the left to the right), and quantities of lines in blocks are $k_{1},k_2,\dots,k_{w+1}$ (counting from up to down). Hence, the $(\al,\be)$-th block of this matrix is a $k_{\al}\times k_{w+2-\be}$-matrix. These matrices will be called skew $k_*$-block matrices.

\medskip
$\forall \ i = 0,\dots, w$ we define skew $k_*$-block matrices $C_i=C_i(S_{****})$ as follows:
\medskip
The $(\al, \be)$-th block of $C_i$ is $S^t_{w+2-\be,w+1-\be,\al,i}$ if the quadruple $(w+2-\be,w+1-\be,\al,i)$ satisfies (9.13, 9.14)\footnotemark \footnotetext{Really, it always satisfies (9.14).} (i.e. if $S_{w+2-\be,w+1-\be,\al,i}$ exists);
\medskip
The $(i+1,w+1-i)$-th block of $C_i$ is $I_{k_{i+1}}$, all other blocks of $C_i$ are 0:
$$(C_i)_{\al\be}=S^t_{w+2-\be,w+1-\be,\al,i}\eqno{(9.28.1)}$$
$$(C_i)_{i+1,w+1-i}=I_{k_{i+1}}\eqno{(9.28.2)}$$
Particularly, the $(\al,\be)$-th block of $C_i$ is a $(k_\al\times k_{w+2-\be})$-th matrix.
\medskip
$\forall \ i = 0,\dots, w$ we define skew $k_*$-block matrices $\bar C_i=\bar C_i(S_{****})$ as follows:
\medskip
The $(\al,\be)$-th block of $\bar C_i$ is given by the formula

$$(\bar C_i)_{\al,\be}=-P_{\al,w-1-i,w+2-\be,w-\be}=\bar S^{t}_{\be,\be-1,w+2-\al,i}\eqno{(9.29)}$$
if the quadruple $(\al,w-1-i,w+2-\be,w-\be)$ belongs to the non-trivial domain of $P_{****}$;
\medskip
For $i=0,\dots,w$
$$(\bar C_i)_{w+1-i,i+1}=I_{k_{w+1-i}}\eqno{(9.30)}$$
other block entries of $\bar C_i$ are 0.
\medskip
{\bf Remark 9.31.} Formula (9.30) is concordant with (9.29), if we consider $P_{****}$ from (9.22). Nevertheless, some 0-blocks of $\bar C_i$ correspond to $P_{**yz}$ where $(y,z)$ do not satisfy (9.19), and hence this $P_{****}$ is not defined.
\medskip
Finally, we define elements $B(S_{****}):=\sum_{i=0}^w C_iN^i\in M_r(\p)[N]$ and $\bar B(S_{****}):=\sum_{i=0}^w \bar C_iN^i\in M_r(\p)[N]$.
\medskip
Example for $w=3$:
$$B(S_{****})=\left(\matrix 0&0&0&I_{k_1}\\0&0&0& S^t_{1020}\\0&0&0& S^t_{1030}\\0&0&0& S^t_{1040} \endmatrix \right)+
\left(\matrix 0&0&0&0\\0&0&I_{k_2}& 0\\0&0&S^t_{2131}& S^t_{1031}\\0&0&S^t_{2141}& S^t_{1041} \endmatrix \right)N+$$ $$+ \left(\matrix 0&0&0&0\\0&0&0& 0\\0&I_{k_3}&0& 0\\0& S^t_{3242}&S^t_{2142}& S^t_{1042} \endmatrix \right)N^2 + \left(\matrix 0&0&0&0\\0&0&0&0\\0&0&0& 0\\I_{k_4}&0&0&0 \endmatrix \right)N^3$$
\medskip
$$\bar B(S_{****})=\left(\matrix \bar S^{t}_{1040}&0&0&0&\\ \bar S^{t}_{1030}&0&0&0\\ \bar S^{t}_{1020}&0&0&0\\I_{k_4}&0&0&0  \endmatrix \right)+
\left(\matrix \bar S^{t}_{1041}&\bar S^{t}_{2141}&0&0\\ \bar S^{t}_{1031}&\bar S^{t}_{2131}&0& 0\\0&I_{k_3}&0&0 \\0&0&0&0  \endmatrix \right)N+$$ $$+ \left(\matrix \bar S^{t}_{1042}&\bar S^{t}_{2142}& \bar S^{t}_{3242}&0\\0&0&I_{k_2}& 0\\0&0&0& 0\\0&0&0& 0  \endmatrix \right)N^2 + \left(\matrix 0&0&0&I_{k_1}\\0&0&0&0\\0&0&0& 0\\0&0&0&0 \endmatrix \right)N^3=$$
\medskip
$$=\left(\matrix -P_{1242}&0&0&0&\\-P_{2242}&0&0&0\\-P_{3242}&0&0&0\\I_{k_4}&0&0&0  \endmatrix \right)+
\left(\matrix -P_{1142}&-P_{1131}&0&0\\-P_{2142}&-P_{2131}&0& 0\\0&I_{k_3}&0&0 \\0&0&0&0  \endmatrix \right)N+$$ $$+ \left(\matrix -P_{1042}&-P_{1031}& -P_{1020}&0\\0&0&I_{k_2}& 0\\0&0&0& 0\\0&0&0& 0  \endmatrix \right)N^2 + \left(\matrix 0&0&0&I_{k_1}\\0&0&0&0\\0&0&0& 0\\0&0&0&0 \endmatrix \right)N^3$$
\medskip
{\bf Lemma 9.32.} $B(S_{****})^t \cdot \bar B(S_{****})=I_r N^w\in M_r(\p)[N]$.
\medskip
{\bf Proof.} We denote $B(S_{****})^t \cdot \bar B(S_{****})$ by $\sum_\mu \Cal C_\mu N^\mu$. The fact that $\Cal C_w=I_{r}$ is obvious: the only non-zero factors that enter in the sum $\sum_{\ga=0}^{w} C_\ga^t \bar C_{w-\ga}$ are products of blocks of $C_*$, $\bar C_*$ containing $I_*$, and they form $I_r$. Also it is obvious that for $\mu>w$ we have $\Cal C_\mu=0$, because all products whose sum is $\Cal C_\mu$, have at least one factor 0. We need to consider $\Cal C_\mu$ for $\mu<w$. (9.28.1), (9.28.2), (9.29), (9.30) give us (here and below $(C^t_\ga)_{\nu\de}$ is the $(\nu\de)$-th block of $C^t_\ga$, i.e. $(C^t_\ga)_{\nu\de}=((C_\ga)_{\de\nu})^t$ )
$$(\Cal C_\mu)_{\nu\pi}=\sum_{\ga=0}^\mu \sum_{\de=1}^{w+1} (C^t_\ga)_{\nu\de}(\bar C_{\mu-\ga})_{\de\pi}\eqno{(9.32.1.1)}$$ $$=-\sum_{\ga,\de} S_{w+2-\nu, w+1-\nu, \de, \ga}P_{\de,w-1-\mu+\ga,w+2-\pi,w-\pi}+ \eqno{(9.32.1.2)}$$
$$-P_{w+2-\nu,2w-\mu-\nu,w+2-\pi,w-\pi}+S_{w+2-\nu,w+1-\nu,w+2-\pi,\mu+1-\pi}\eqno{(9.32.1.3)}$$
where (9.32.1.2) corresponds to the products $(C^t_\ga)_{\nu\de}(\bar C_{\mu-\ga})_{\de\pi}$ where both terms $\ne 0, \ I_*$, and (9.32.1.3) corresponds to the products where one of the terms is $I_*$.

Let us find the relations satisfied by $\mu,\ \nu, \ \pi$ and the domain of summation by $\ga, \ \de$ in (9.32.1.2). We have
$$\matrix (C^t_\ga)_{\nu\de}\ne0, I_{k_*} \ \iff \ \nu\ge w+1-\ga \ \ \wedge \ \ \de\ge\ga+2 \\
(\bar C_{\mu-\ga})_{\de\pi}\ne0, I_{k_*} \ \iff \ \de\le w-(\mu-\ga) \ \ \wedge \ \ \pi \le \mu-\ga+1\endmatrix \eqno{(9.32.2)}$$
The set of $\ga, \ \de$ is non-empty $\iff \ \mu\le w-2$ and $\mu +\nu -\pi\ge w$. In this case the conditions (9.32.2) on $\ga, \ \de$ become
$$\mu+1-\pi\ge\ga\ge w+1-\nu\eqno{(9.32.3.1)}$$
$$w+\ga-\mu\ge\de\ge\ga+2\eqno{(9.32.3.2)}$$
Now we use Proposition 9.25 for
$$\matrix i=w+3-\nu && \mu=j+i-3 \\
j=\mu+\nu-w &\iff & \nu=w-i+3 \\
\psi=w-\pi+2 && \pi=w-\xi=w-\psi+2\\
\xi=w-\pi \endmatrix\eqno{(9.32.4)}$$
and summation variables $\al, \ \be$ in (9.25.1) are
$$\matrix \al=\de \\ \be=\ga-i+2 \endmatrix\eqno{(9.32.5)}$$
Under this variable change, (9.32.1.2) becomes the double sum in (9.25.10), and (9.32.1.3) becomes
$$-P_{i-1,w-j,\psi,\xi}+S_{i-1,i-2,\psi,i-2+\xi+j-w}$$
hence the desired. $\square$
\medskip
Let for $i=0,\dots,w-1$ $X_i$ be skew $k_*$-matrices having the following property:
\medskip
If $(\al, \ \be)$ are such that the $(\al,\ \be)$-block of $\bar C_i$ is 0 or $I_*$ then $(X_i)_{\al\be}=(\bar C_i)_{\al\be}$;

If $(\al, \ \be)$ are such that the $(\al,\ \be)$-block of $\bar C_i$ is $\ne 0, \ I_*$ then $(X_i)_{\al\be}$ is arbitrary.
\medskip
We denote $X:=\sum_{i=0}^{w-1}X_iN^i$.
\medskip
{\bf Lemma 9.33.} If $B(S_{****})^t \cdot X\in N^wM_r(\p[N])$ then $X=\bar B(S_{****})$.
\medskip
{\bf Proof.} For any fixed $\mu, \ \nu, \ \pi$ (9.32.1.2), (9.32.1.3) become
$$\sum_{\ga,\de} S_{w+2-\nu, w+1-\nu, \de, \ga} (X_{\ga-\mu})_{\de,\pi}\eqno{(9.33.1a)}$$ $$
+ (X_{\mu+\nu-w-1})_{w+2-\nu,\pi}+S_{w+2-\nu,w+1-\nu,w+2-\pi,\mu+1-\pi}=0\eqno{(9.33.1b)}$$
This is system of linear equations with unknowns $(X_i)_{\al\be}$ where $$\matrix 0\le i \le w-1 \\ \\ 1\le \be\le i+1\\ \\ 1\le\al \le w-i \endmatrix \eqno{(9.33.2)}$$ (for other values of $i,\ \al, \ \be$ we have $(X_i)_{\al\be}=0$). We arrange $(X_i)_{\al\be}$ in decreasing order of $i+\al$ (for $(X_i)_{\al\be}$ having equal $i+\al$ their ordering is arbitrary), and we arrange equations (9.33.1) in decreasing order of $\mu$ (for equations having equal $\mu$ the order of equations corresponds to the order of $(X_i)_{\al\be}$ having equal $i+\al$). Under this arrangement of unknowns and equations, the matrix of the system (9.33.1) becomes unitriangular. Really, for any $i,\ \al, \ \be$ satisfying (9.33.2) there is exactly one values of $\mu, \ \nu, \ \pi$ --- namely, $$\matrix \mu=i+\al-1 \\ \\ \nu=w+2-\al\\ \\ \pi=\be \endmatrix $$ such that the first term of (9.33.1b) is $(X_i)_{\al\be}$. Other terms of (9.33.1) for these $\mu, \ \nu, \ \pi$ --- namely, the terms that enter in (9.33.1a) --- contain $(X_{\ga-\mu})_{\de,\pi}$ such that $\ga-\mu+\de>i+\al$ hence unitriangularity.
\medskip
Lemma 9.32 affirms that
$$(X_i)_{\de,\pi}=-P_{\de,w-1-i,w+2-\pi,w-\pi}$$
is a solution to this system. Unitriangularity implies that this solution is unique. $\square$
\medskip
Let us consider $\g q_H$ from above (9.1).
\medskip
{\bf Lemma 9.34.} $\forall \ u=1,\dots, w+1$, for $v=u-1$, $\forall \ i=1,\dots, k_u$ the elements

$$\om_{ui}:=N^{v}l_{ui}+\sum_{y=u+1}^{w+1}\sum_{z=u-1}^{y-2}\sum_{j=1}^{k_y}(S_{uvyz})_{ij}N^zl_{yj}\eqno{(9.34.1)}$$

form a basis of $\g q_H$. $\square$
\medskip
We denote the set of elements $\om_{ui}$ ($u$ is fixed, $i$ varies) by $\hat \om_{u}$ (matrix columns). So, (9.34.1) becomes ($v=u-1$)
$$\hat \om_{u}=N^{v}\hat l_{u}+\sum_{z=u-1}^{w-1}\sum_{y=z+2}^{w+1}S_{uvyz}N^z\hat l_{y}\eqno{(9.35)}$$

Let $\vf_i$, $i=1,...,r$ be the basis of $L'$ dual to $l_i$, i.e. $\vf_i(l_j)=\delta_{ij}$. We shall need the dual numbers $k'_i:=k_{w+2-i}$ (inverse order of $k_*$). We consider the analogous two-subscript notation of $\vf_i$, but the order of segments of the partition of $\vf_i$ is opposite, namely:
$$(\vf_{w+1,1},\dots,\vf_{w+1,k'_{w+1}},\ \ \vf_{w,1},\dots,\vf_{w,k'_{w}},\ \ \ \dots \ \ \ \vf_{11},\dots,\vf_{1,k'_{1}}):= (\vf_1,\dots, \vf_r)\eqno{(9.35.1)}$$
(order of elements $\vf_*$ is the same in both sides of this equality).
\medskip
{\bf Lemma 9.36.} $\forall \ u=1,\dots, w+1$, for $v=u-1$, $ \forall \ i=1,\dots, k'_u$ the elements

$$\chi_{ui}:=N^{v}\vf_{ui}+\sum_{y=u+1}^{w+1}\sum_{z=u-1}^{y-2}\sum_{j=1}^{k'_y}(\bar S_{uvyz})_{ij}N^z\vf_{yj}\eqno{(9.36.1)}$$

form a basis of ${\g q'_H}^w$.

{\bf Proof.} As above we denote the set of elements $\vf_{ui}$, resp. $\chi_{ui}$ ($u$ is fixed, $i$ varies) by $\hat \vf_{u}$, resp. $\hat \chi_{u}$ (matrix columns). (9.35), (9.36.1) can be written in terms of blocks of $C_i$, $\bar C_i$:

$$\hat \om_{u}=\sum_{z=0}^{w}\sum_{y=1}^{w+1} (C_z^t)_{w+2-u,y}N^z\hat l_{y}$$
$$\hat \chi_{u}=\sum_{z=0}^{w}\sum_{y=1}^{w+1} (\bar C_z^t)_{uy}N^z\hat \vf_{w+2-y}$$
We must prove that $\forall \ u_1, \ u_2$ we have $\hat \om_{u_1}\hat \chi_{u_2}^t=\de_{u_1}^{u_2}I_{k_{u_1}}N^w$ (product is pairing). This is immediate:

$$\hat \om_{u_1}\hat \chi_{u_2}^t=\sum_{z_1=0}^{w}\sum_{y_1=1}^{w+1}\sum_{z_2=0}^{w}\sum_{y_2=1}^{w+1} (C_{z_1}^t)_{w+2-u_1,y_1}N^{z_1}\hat l_{y_1} N^{z_2}\hat \vf_{w+2-y_2}^t (\bar C_{z_2})_{y_2,u_2}$$
We have $\hat l_{y_1} \hat \vf_{w+2-y_2}^t=\de_{y_1}^{y_2}I_{k_{y_1}}$, hence
$$\hat \om_{u_1}\hat \chi_{u_2}^t=\sum_{z_1=0}^{w}\sum_{y=1}^{w+1}\sum_{z_2=0}^{w} (C_{z_1}^t)_{w+2-u_1,y}N^{z_1+z_2} (\bar C_{z_2})_{y,u_2} = \sum_{z=0}^{2w}(\Cal C_z)_{w+2-u_1,u_2}N^z$$
Lemma 9.32 implies the desired. $\square$
\medskip
{\bf Corollary 9.37.} Matrices $S_{uvyz}(L')$ for the dual lattice $L'$ are $\bar S_{uvyz}(L)$ (order of segments of $\vf_*$, and hence of numbers $k_*$, is inverse).
\medskip
We define $D_i$ ($i=1,\dots, w$ ) as $-<N^{i-1}(l_*),f_*>$ where the ordering of $l_*$ is $l_{11},\dots, l_{w+1, k_{w+1}}$. Hence, $D_i$ is a union of $w+1$ matrices $D_{ij}$, where $D_{ij}$ is a $r\times k_j$-matrix, and $(D_{ij})_{\al\be}=-<N^{i-1}(l_{j\be}),f_\al>$. This means that there are relations between $D_{ij}$ coming from (9.20), namely:

$$D_{v+1,u}= -\sum_{z=v}^{w-1}\sum_{y=z+2}^{w+1}D_{z+1,y} P_{uvyz}^t\eqno{(9.38)}$$

Like in (9.21), for $(z,y)$ satisfying $y\ge z+1$ (resp. $y< z+1$) we shall call the corresponding $D_{zy}$ as belonging to the trivial (resp. non-trivial) domain.
\medskip
{\bf Lemma 9.39.} $\Psi_N B(S_{****})\in M_r(\p[[N]])$.
\medskip
{\bf Proof.} For any $\mu=0,\dots, w-1$ we must prove that $\g C_\mu:=\sum_{\de =0}^\mu D_{w-\de}C_{\mu-\de}$ is 0. The $\nu$-th block ($\nu=1,\dots,w+1$) of this matrix is
$$(\g C_\mu)_\nu:=\sum_{\de =0}^\mu \sum_{\ga=1}^{w+1} D_{w-\de,\ga}(C_{\mu-\de})_{\ga \nu}\eqno{(9.39.1)}$$
We have
$$(C_{\mu-\de})_{\ga \nu}\ne0, I_{k_*} \ \ \iff \de\le \nu+\mu-w-1 \ \ \wedge \ \ \ga\ge \mu-\de+2$$
$$(C_{\mu-\de})_{\ga \nu}=I_{k_*} \ \ \iff \de=\nu+\mu-w-1 \ \ \wedge \ \ \ga=w+2-\nu$$
hence (9.39.1) becomes
$$(\g C_\mu)_\nu=\sum_{\de =0}^{\nu+\mu-w-1}\sum_{\ga=\mu-\de+2}^{w+1}D_{w-\de,\ga}S^t_{w+2-\nu,w+1-\nu,\ga,\mu-\de}+\eqno{(9.39.2.1)}$$
$$+D_{2w+1-\nu-\mu,w+2-\nu}\eqno{(9.39.2.2)}$$
where (9.39.2.1) is non-empty if $\nu+\mu\ge w+1$.

Terms $D_{2w+1-\nu-\mu,w+2-\nu}$ always belong to the non-trivial domain. We separate the terms of (9.39.2.1) in terms of trivial and non-trivial domain:
$$(\g C_\mu)_\nu=\sum_{\de =0}^{\nu+\mu-w-1}\sum_{\ga=\mu-\de+2}^{w-\de} D_{w-\de,\ga}  S^t_{w+2-\nu,w+1-\nu,\ga,\mu-\de}+\eqno{(9.39.3.1)}$$
$$+\sum_{\de =0}^{\nu+\mu-w-1}\sum_{\ga=w-\de+1}^{w+1}  D_{w-\de,\ga}  S^t_{w+2-\nu,w+1-\nu,\ga,\mu-\de}+\eqno{(9.39.3.2)}$$
$$+D_{2w+1-\nu-\mu,w+2-\nu}\eqno{(9.39.3.3)}$$

Now we substitute non-trivial $D_{**}$ by linear combinations of the trivial ones, using (9.38):
$$(\g C_\mu)_\nu=-\sum_{\de =0}^{\nu+\mu-w-1}\sum_{\ga=\mu-\de+2}^{w-\de}   \sum_{z=w-\de-1}^{w-1}\sum_{y=z+2}^{w+1}D_{z+1,y} P_{\ga,w-\de-1,y,z}^t       S^t_{w+2-\nu,w+1-\nu,\ga,\mu-\de}+\eqno{(9.39.4.1)}$$
$$+\sum_{\de =0}^{\nu+\mu-w-1}\sum_{\ga=w-\de+1}^{w+1}D_{w-\de,\ga}S^t_{w+2-\nu,w+1-\nu,\ga,\mu-\de}-\eqno{(9.39.4.2)}$$
$$-\sum_{z=2w-\nu-\mu}^{w-1}\sum_{y=z+2}^{w+1}D_{z+1,y} P_{w+2-\nu,2w-\nu-\mu,y,z}^t \eqno{(9.39.4.3)}$$
Now we change variables in (9.39.4.2):
$$\de=w-z-1$$
$$\ga=y$$
interchange the order of summation and transpose:
$$(\g C_\mu)_\nu=\sum_{z=2w-\nu-\mu}^{w-1}\sum_{y=z+2}^{w+1} \g K(\mu,\nu,z,y) D^t_{z+1,y}$$

where $$\g K(\mu,\nu,z,y)=-( \sum_{\de =w-z-1}^{\nu+\mu-w-1}\sum_{\ga=\mu-\de+2}^{w-\de} S_{w+2-\nu,w+1-\nu,\ga,\mu-\de}    P_{\ga,w-\de-1,y,z}) +$$
$$+S_{w+2-\nu,w+1-\nu,y,\mu-w+z+1}-P_{w+2-\nu,2w-\nu-\mu,y,z} \eqno{(9.39.5)}$$
Change of variables in (9.39.5):
$$\matrix y=\psi & \nu=w+3-i & \ga=\al \\ z=\xi & \mu=i+\al-3 & \de=\al-\be-1 \endmatrix $$
transforms (9.39.5) to (9.25.1), hence all $\g K(\mu,\nu,z,y)$ are 0. $\square$
\medskip
Let $M'=M^{\prime w}$ be the $w$-dual of $M$.
\medskip
{\bf Lemma 9.40.} For $i=1,\dots, w+1$ we have $k_i(M')=k'_{i}(M)$.
\medskip
{\bf Proof.} We use Lemma 10.2. $\mu$ of 10.2 is $w$, $m_i$ of 10.2 are $d_{w+2-i}$. (10.2.5) means that $d_i(M')=w-d_i$. The result follows immediately from the properties of dual 0-partitions. $\square$

We shall another basis $\hat \eta$ of $L(M')$ obtained by a permutation matrix from the basis (9.35.1). Namely, we let $\eta_{ij}:=\vf_{ij}$ ($i=1, \dots, w+1, \ j=1,\dots, k'_i$), but the order of elements $\eta_{ij}$ is the following ($\hat \eta$ is a matrix column):

$$\hat \eta=(\eta_{11},\dots, \eta_{1,k'_1}, \eta_{21},\dots, \eta_{2,k'_2}, \dots, \eta_{w+1,1},\dots, \eta_{w+1,k'_{w+1}})^t$$

We shall consider only $M$ satisfying the following (compare with 9.12)
\medskip
{\bf Condition 9.41.} For any $u$ elements $N^\al(\eta_{\be\ga})$, where $\al\in[u,\dots, w-1]$, $\be\in[\al+2,\dots, w+1]$, $\ga\in [1,\dots, k'_\be]$, are linearly independent over $\p$.
\medskip
According the general principle that almost all $n$-uples $(v_1,\dots,v_n)$ of vectors in $n$-dimensional vector space form a basis of this space, we can guess that almost all $M$ satisfy 9.41. Really, Lemma 9.40 affirms that the dimension of $N^u \Lie(M')$ is exactly the quantity of elements $N^\al(\eta_{\be\ga})$ mentioned in 9.41. Again by Lemma 9.40, we see that Condition 9.41 implies that these elements are a basis of $N^u \Lie(M')$.
\medskip
We can use ideas of [GL] in order to show that a large class of Anderson t-motives $M$ satisfies Condition 9.41. Namely,
\medskip
Finally, we have
\medskip
{\bf Conjecture 9.42.} All $M$ satisfy Condition 9.41.
\medskip
{\bf Theorem 9.43.} Let $M$ satisfy Condition 9.41. Then
\medskip
{\bf Proof.} We denote the basis $\vf_*$ from (9.35.1) by $\hat \vf$ (matrix column). We have $\hat \eta= \g I\cdot \hat \vf$ where $\g I$ is a matrix of the change of bases. It is a skew $k'_*$-block anti-identity matrix, i.e. its antidiagonal block entries are identity matrices: the $(w+2-i,i)$-block is $I_{k_i}$, and other block entries are 0.

We denote by $\Psi'_N$, resp. $\Psi'_{N,\eta}$ the $N$-$\Psi$-series of $M'$ in bases $\hat \vf$, resp. $\hat \eta$ (as a base of $M'$ over $\p[T]$ we use $\hat f'$ in both cases). We have $\Psi'_N=\Psi_N^{t-1}\Xi_N^{-w}$, $\Psi'_{N,\eta}=\Psi'_N\g I^t$. Since $M$ satisfies Condition 9.41, there exists a set of Siegel matrices for $M'$ with respect to the basis $\hat \eta$. We denote it by $U_{****}$. It defines $B(U_{****})$ --- the corresponding $B$ in $M_r(\p)[N]$.
\medskip
We denote $\Psi_NB(A(M)_{****})$, resp. $\Psi'_{N,\eta}B(U_{****})$ by $\g Z(M)$, resp. $\g Z(M^d)$. We have $\g Z(M)\in M_r(\p[[N]])$, $\g Z(M^d)\in M_r(\p[[N]])$ (Lemma 9.39), hence
$$\g Z(M)^t\g Z(M^d)=B(A(M)_{****})^t\Psi_N^t \Psi'_N \g I^t B(U_{****})=$$ $$B(A(M)_{****})^t \g I^t B(U_{****}) \Xi_N^{-w}\in M_r(\p[[N]])$$ and hence
$$  B(A(M)_{****})^t \g I^t B(U_{****}) \g I^t  \in \Xi_N^{w}M_r(\p[[N]])$$
We have  $\g I^t B(U_{****}) \g I^t$ is of the form $X$ of Lemma 9.33. Further, $\Xi_N^{w}\in N^wM_r(\p[[N]])$, hence $\g I^t B(U_{****}) \g I^t=\bar B(A(M)_{****})$ (Lemma 9.33). This means that $U_{uvyz}=\bar A_{uvyz}$. The theorem follows from Lemma 9.36. $\square$
\medskip

\newpage
{\bf Case of linearly dependent columns.}

We shall consider for simplicity the case $w=2$. We use notations: $\Psi:=\Psi_N(M)$, $\Psi^d:=\Psi_N(M^d)$, $\Psi^d=D^d_2N^{-2}+D^d_1N^{-1}+...$ in the basis $\vf_1,...,\vf_r$ dual to $l_1,...,l_r$. We denote $D^d_2=D_{21}^d | D_{22}^d | D_{23}^d $ where $D_{2i}^d$ are $r\times k_i$ matrices.  Let some of the columns of $D_{21}^d$ are linearly dependent. Interchanging the order of columns we can assume that $\exists \ u>0$ such that the first $k_1-u$ columns of $D_{21}^d$ are linearly independent and the last $u$ columns of $D_{21}^d$ are their linear combinations. Further, we assume that the remaining $u$ columns of $D_{21}^d$ forming (together with the first $k_1-u$ columns of $D_{21}^d$) a basis of $NV$ are the first $u$ columns of $D_{22}^d$. This means that we can represent $D_{21}^d$, $D_{22}^d$ as unions:

$D_{21}^d=(D_{211}^d | D_{212}^d)$, $D_{22}^d=(D_{221}^d | D_{222}^d)$, where the quantities of rows in all these matrices are $r$ and the quantities of columns in $D_{211}^d$, resp. $ D_{212}^d$, $D_{221}^d, \ D_{222}^d$ are $k_1-u$, resp. $u, \ u, \ k_2-u$. We have: $\exists$ a $(k_1-u)\times u$-matrix $Y$ such that $D_{212}^d=D_{211}^dY$.

We denote by $\tilde \Psi^d$ the matrix obtained from $\Psi^d$ by permutation of (12)-block and (21)-block. We assume that (1,2,3)-blocks of $\tilde \Psi^d$ are "good", i.e. there exists $\tilde D^d$ such that

$$\tilde \Psi^d\cdot\tilde D^d\in M_r(\p[[N]])$$

We denote by $\hat \Psi^d$ the matrix obtained from $\Psi^d$ by elimination of (12)-block. It is a $5\times 4$-block matrix, and we have

$$\Psi^d=\hat \Psi^d U_3$$ where $U_3=\left(\matrix 1&Y&0&0&0\\ 0&0&1&0&0\\ 0&0&0&1&0\\ 0&0&0&0&1\endmatrix \right)$ (block structure; sizes of blocks: $k_1-u; \ u; \ k_2-u; \ k_3$ for columns; $k_1-u; \ u; \ u; \ k_2-u; \ k_3$ for lines).

$$\tilde \Psi^d=\hat \Psi^d U_2$$ where $U_2=\left(\matrix 1&0&Y&0&0\\ 0&1&0&0&0\\ 0&0&0&1&0\\ 0&0&0&0&1\endmatrix \right)$ (block structure; sizes of blocks: $k_1-u; \ u; \ k_2-u; \ k_3$ for columns; $k_1-u; \ u; \ u; \ k_2-u; \ k_3$ for lines).
\medskip
$U_3$ has a right inverse $U_3^{-1}=\left(\matrix 1&0&0&0\\ 0&0&0&0\\ 0&1&0&0\\ 0&0&1&0\\ 0&0&0&1\endmatrix \right)$ such that we have

$$\hat \Psi^d = \Psi^dU_3^{-1}$$

hence $B^t\Psi^t\tilde \Psi^d\cdot\tilde B^d\in M_r(\p[[N]])$

and $\tilde \Psi^d=\Psi^dU_3^{-1}U_2$

As earlier we have $$\Psi^t\Psi^d=\Xi^{-w}$$

hence $B^t U_3^{-1}U_2 B^d\in N^2 M_r(\p[[N]])$
\medskip
\medskip

Generalizations.
\medskip
1. To prove that $\g q_{M_1\otimes M_2}=\g q_{M_1} \otimes \g q_{M_2}$.
\medskip
2. It is known that if $A_1$, $A_2$ are abelian varieties then $A_1\otimes A_2$ is a mixed motive. Is it possible to get an analog of (1) for it?
\medskip
3. Analog of the main theorem for groups other than $GL_r$. We have a theorem that any full sublattice of $\n Z^r$ is described by its elementary divisors $d_1 | d_2| ... | d_r$. There is an analog of this theorem for any reductive group. We should find the corresponding generalization of the main theorem of the present paper.
\medskip
This is an analog of Theorem 5 for arbitrary $N$, and of Theorem 6 for the operator of duality. Its proof, as well as its generalization to the case of non-pure $M$, can be easily obtained using the same ideas of the proofs of Theorems 5, 6 (this explains the terminology: this is not a conjecture, but a theorem whose proof is not written yet).
\medskip
For $m=1$ Formula 9.1 and Theorem 5 imply immediately
\medskip
{\bf Proposition 9.4.} Result 9.3 is proved for $m=1$. $\square$
\medskip

\input amstex
\documentstyle{amsppt}
\magnification1200
\tolerance=10000
\overfullrule=0pt
\def\n#1{\Bbb #1}
\def\p{\Bbb C_{\infty}}
\def\fr{\hbox{fr}}

\def\im{\hbox{im }}
\def\invlim{\hbox{invlim}}
\def\tr{\hbox{tr }}

\def\Gal{\hbox{Gal }}

\def\Exp{\hbox{Exp}}

\def\Hom{\hbox{Hom}}

\def\End{\hbox{End}}
\def\Prin{\hbox{Prin}}
\def\Ker{\hbox{Ker }}
\def\Lie{\hbox{Lie}}

\def\Div{\hbox{Div}}
\def\Pic{\hbox{Pic}}

\def\ord{\hbox{ord}}

\def\Id{\hbox{Id}}
\def\Cl{\hbox{Cl}}
\def\Supp{\hbox{ Supp }}
\def\Spec{\hbox{ Spec }}

\def\diag{\hbox{ diag }}
\def\Diag{\hbox{ Diag }}

\def\e11{I_{11}}

\def\ga{\goth A}
\def\w{\hbox{\bf A}}
\def\x{\hbox{\bf K}}
\def\ve{\varepsilon}
\def\vf{\varphi}

\def\ve{\varepsilon}

\def\vf{\varphi}

\def\de{\delta}

\def\ga{\gamma}

\def\be{\beta}

\def\al{\alpha}

\def\om{\omega}
\def\g{\goth }

\topmatter
\title
Duality of Anderson $t$-motives
\endtitle
\author
A. Grishkov, D. Logachev\footnotemark \footnotetext{E-mails: shuragri{\@}gmail.com; logachev94{\@}gmail.com (corresponding author)\phantom{*******************}}
\endauthor
\thanks Thanks: The authors are grateful to FAPESP, S\~ao Paulo, Brazil for a financial support (process No. 2017/19777-6). The first author is grateful to SNPq, Brazil, to RFBR, Russia, grant 16-01-00577a (Secs. 1-4), and to Russian Science Foundation, project 16-11-10002 (Secs. 5-8) for a financial support. The second author is grateful to Gilles Lachaud for invitation to IML and to Laurent Lafforgue for invitation to
IHES where this paper was started, and to Vladimir
Drinfeld for invitation to the University of Chicago where
this paper was continued. Discussions with Greg Anderson, Vladimir
Drinfeld, Laurent Fargues, Alain Genestier, David Goss, Richard Pink, Yuichiro
Taguchi, Dinesh Thakur on
the subject of this paper were very important. Particularly, Alain
Genestier informed me about the paper of Taguchi where the notion of
the dual of a Drinfeld module is defined. Further, Richard Pink indicated me an important reference (see Section 6 for details); proof of the main theorem of the present paper grew from it. Finally, Vladimir
Drinfeld indicated me the proof of the Theorem 12.6 and Jorge Morales
gave me a reference on classification of quadratic forms over $\n
F_q[T]$ (Remark 7.8).
\endthanks
\NoRunningHeads
\address
First author: Departamento de Matem\'atica e estatistica
Universidade de S\~ao Paulo. Rua de Mat\~ao 1010, CEP 05508-090, S\~ao Paulo, Brasil, and Omsk State University n.a. F.M.Dostoevskii. Pr. Mira 55-A, Omsk 644077, Russia.
\medskip
Second author: Departamento de Matem\'atica, Universidade Federal do Amazonas, Manaus, Brasil
\endaddress

\abstract Let $M$ be a t-motive. We introduce the notion of duality for $M$. Main results of the paper (we consider uniformizable $M$ over $\n F_q[T]$ of rank $r$, dimension $n$, whose nilpotent operator $N$ is 0):

1. Algebraic duality implies analytic duality (Theorem 5). Explicitly, this means that the lattice of the dual of $M$ is the dual of the lattice of $M$, i.e. the transposed of a Siegel matrix of $M$ is a Siegel matrix of the dual of $M$.

2. Let $n=r-1$. There is a 1 -- 1 correspondence between pure t-motives (all they are uniformizable), and lattices of rank $r$ in $\p^{n}$ having dual (Corollary 8.4).

3. Pure t-motives have duals which are pure t-motives as well (Theorem 10.3).

4. Some explicit results are proved for $M$ having complete
multiplication. The CM-type of the dual of $M$ is the complement of
the CM-type of $M$. Moreover, for $M$ having multiplication by a
division algebra there exists a simple formula for the
CM-type of the dual of $M$ (Section 12).

5. We construct a class of non-pure t-motives (t-motives having the
completely non-pure row echelon form) for which duals
are explicitly calculated (Theorem 11.5). This is the first step
of the problem of description of all t-motives having duals.

6. If $M$ has good ordinary reduction then the kernels of reduction maps on groups of  torsion points for $M$ and its dual are complementary
with respect to a natural pairing (proof is given for a particular case, Conjecture 13.4.1).
\endabstract
\endtopmatter
\document
{\bf 0. Introduction.}
\nopagebreak
\medskip
t-motives are the function field analogs of abelian varieties (more exactly, of abelian varieties with multiplication by an imaginary quadratic field, see [L09]). Main references for t-motives are [A], [G]. Nevertheless, function field analogs of some basic results in the theory of abelian varieties are not known yet.

The present paper contains an analog of such result. Namely, we introduce the notion of duality for a t-motive $M$ (this is not the duality in
a Tannakian category!), and we prove some properties of this notion, see the abstract. Particularly, if $M$ is uniformizable and has dual then
the lattice of the dual of $M$ is the dual of the lattice of $M$ (Theorem 5)\footnotemark \footnotetext{Here this result is proved for $M$ having
the associated nilpotent operator $N$ (see (1.9.2))
equal to 0. The same result for $M$ having $N\ne0$ is proved in [GL18].}. An immediate corollary of the above theorem and the result of
Drinfeld on 1 -- 1 correspondence between Drinfeld modules and lattices in $\p$ (here $\p$ is the function field analog of $\n C$) is Corollary 8.4: there is a 1 -- 1 correspondence between pure t-motives of dimension $r-1$ and rank $r$, and lattices of rank $r$ in $\p^{r-1}$ having dual (not all such lattices have dual).

Let us give more details on the contents of the paper. For simplicity, most results are proved for t-motives over the ring $\hbox{\bf A}=\n F_q[T]$,
and we consider, with few exceptions, only the case $N=0$. The main definition of duality of t-motives
(definition 1.8 --- case $\hbox{\bf A}=\n F_q[T]$ and definition 1.13 --- general case) is given in
Section 1.\footnotemark \footnotetext{A version of the definition of duality is obtained independently in [Tae], 2.2.}
Lemma 1.10 gives the explicit matrix form of the definition of duality of t-motives. Since Taguchi in [T] gave a definition of dual to
a Drinfeld module, we prove in Proposition 1.12.3 that the definition of the present paper is equivalent to the original definition of Taguchi.
Section 1.14 contains a definition of duality for abelian $\tau$-sheaves ([BH], Definition 2.1), but we do not develop this subject.

Section 2 contains the definition of the dual lattice. Section 3 contains explicit formulas for the dual lattice.
Section 5 contains the statement and the proof of the main theorem 5 --- coincidence of algebraic and analytic duality for the case
$\hbox{\bf A}=\n F_q[T]$ (section 4 contains the statement of the corresponding conjecture for the case of general $\hbox{\bf A}$).
Section 6 contains the theorem 6 describing the lattice of the tensor product of two t-motives (case $N=0$; the proof for the general
case was obtained, but not published, by Anderson). Section 7 contains the notion of self-dual t-motives and polarization form on them.
Some examples are given. We discuss in Section 8 the problem of correspondence between uniformizable t-motives and lattices.
Section 9 gives the statement of the main result for the case $N\ne 0$ without proof and a reformulation of the theorem 5 in terms
of Hodge-Pink structures of constant weight.

Further on, we prove in Section 10 that pure t-motives have duals which are pure
t-motives as well, and some related results (a proof that the dual of an abelian
$\tau$-sheaf is also an abelian $\tau$-sheaf can be obtained using ideas
of Section 10). In Section 11 we consider t-motives having the
completely non-pure row echelon form, and we give an explicit formula
for their duals. In Section 12 we consider t-motives with complete multiplication, and we give for them a very simple version of the proof of the first part of the main theorem.
Section 13 contains some explicit formulas for t-motives of complete multiplication. In 13.1 we describe the dual lattice, in 13.2 we show that the results of Section 12 are compatible with (the first form of) the main theorem of complete multiplication. Section 13.3 contains an explicit proof of the main theorem for t-motives with complete multiplication by two types of simplest fields.
Section 13.4 gives us an application of the notion of duality to the
reduction of t-motives (subject in development, see [L]).
\medskip
{\bf Notations. }
\medskip
$q$ is a power of a prime $p$;
\medskip
{\it Case of $M$ over $\n F_q[T]$}:
\medskip
$\n Z_\infty:=\n F_q[\theta]$, $\n R_{\infty}:=\n F_q((1/\theta))$, $\p$ is the completion of its algebraic closure
($\n Z_\infty$, $\n R_{\infty}$, $\p$ are the function field analogs of $\n Z$, $\n R$, $\n C$ respectively);

$\w:=\n F_q[T]$, $\hbox{\bf K}:=\n F_q((1/T))$;

$\iota: \hbox{\bf A} \to \p$ ($\iota(T)=\theta$) is the standard map of generic characteristic (with one exception (1.16), we shall not consider the case of finite characteristic);
we extend $\iota$ to $\hbox{\bf K}$, and we have $\n Z_\infty=\iota(\w)\subset \p$, $\n R_\infty=\iota(\hbox{\bf K})\subset \p$.

$\goth C$ (resp. $\goth C_2$) is the Carlitz module over $\w=\n F_q[T]$ (resp. over $\n F_{q^2}[T]$).
\medskip
{\it Case of $M$ over an extension of $\n F_q[T]$}:
\medskip
$\n Q_\infty$ is a finite separable extension of $\n F_q(\theta)$;

$\infty$ is a fixed valuation of $\n Q_\infty$ over the infinity of $\n F_q(\theta)$;

${\n Z_\infty \subset \n Q_\infty}$ is the subring of elements which are regular outside $\infty$;

$\n R_\infty$ is the completion of $\n Q_\infty$ at infinity, and $\p$ --- the completion of its algebraic closure --- is the same as of the case of $M$ over $\n F_q[T]$.

$\w\supset \n F_q[T]$, $\hbox{\bf K}\supset \n F_q((1/T))$ are defined by the condition that $\iota: \w \to \n Z_\infty$, $\iota: \hbox{\bf K} \to \n R_\infty$ are isomorphisms.

$\w_C:=\w\underset{\n F_q}\to{\otimes} \p$ (i.e. $\w_C=\p[T]$ for the case of $M$ over $\n F_q[T]$).

$\goth C$ is a Drinfeld module of rank 1 over {\bf A}.
\medskip
If $P=\frac{\sum a_iT^i}{\sum b_iT^i}\in \p(T)$ then $P^{(k)}:=\frac{\sum a_i^{q^k}T^i}{\sum b_i^{q^k}T^i}$. For $x\in \w_C$, $x=a\otimes z$, $a\in \w$, $z\in \p$ we let $x^{(k)}:= a\otimes z^{q^k}$.

$M_r$ is the set of $r\times r$ matrices. If $C=\{c_{ij}\}$ is a matrix with entries $c_{ij}\in \p(T)$ then $C^{(k)}:=\{c_{ij}^{(k)}\}$, $C^t$ is the  transposed of $C$, $C^{(k)\ -1}=(C^{(k)})^{-1}$, $C^{t-1}=(C^t)^{-1}$.

If $M$ is an $\w_C$-module, we define $M^{(1)}$ as the tensor
product $M\otimes_{\w_C,*^{(1)}}\w_C$ with respect to the map $*^{(1)}: \w_C\to \w_C$ (this notation is concordant in the obvious sense with
the above notation $C^{(1)}$).

For a t-motive $M$ we denote by $E=E(M)$ the corresponding t-module (see [G], Theorem 5.4.11; Goss uses the inverse functor $E\mapsto M=M(E)$).

$\Lie(M)$ is $\Lie(E(M))$ ([G], 5.4).

$I_k$ is the unit matrix of size $k$.

Throughout the whole paper the word "canonical" will mean "canonical up to multiplication by elements of $\n F_q^*$".
\medskip
{\bf 1. Definitions. }

\nopagebreak
\medskip
If otherwise is not explicitly stated, throughout the whole paper we consider the case of t-motives $M$ over the ring $\hbox{\bf A}=\n F_q[T]$ such that $N=N(M)=0$. Exceptions: case of arbitrary $\w$ is treated in Sections 1.13, 1.14, 2, 4, 5.2. Case of arbitrary $N$ is treated in Sections 1, 10 and in statements of some results of Anderson in Sections 5, 6.
\medskip
In the present section we consider $M$ such that $N(M)$ is arbitary.
\medskip
Let $\p[T,\tau]$ be the Anderson ring, i.e. the ring of non-commutative polynomials satisfying the following relations (here
$a \in \p$):
$$Ta=aT, \ T\tau = \tau T, \ \tau a = a^q \tau \eqno{(1.1)}$$
We need also an extension of $\p[T,\tau]$ --- the ring $\p(T)[\tau]$
which is the ring of
non-commutative polynomials in $\tau$ over the field of rational
functions $\p(T)$ with
the same relations (1.1). For a left $\p[T,\tau]$-module $M$ we denote by $M_{\p[T]}$ the
same
$M$ treated as a
$\p[T]$-module with respect to the natural inclusion
$\p[T]\hookrightarrow\p[T,\tau]$.
Analogously, we define $M_{\p[\tau]}$; we shall use similar notations
also for the left
$\p(T)[\tau]$-modules.
\medskip
Obviously we have:
\medskip
{\bf (1.2)} For $C\in M_{r}(\p(T))$ operations $C^t$, $C^{-1}$
and $C^{(i)}$ commute.
\medskip
{\bf Definition 1.3.} ([G], 5.4.2, 5.4.12, 5.4.10). A t-motive $M$ is a left
$\p[T, \tau]$-module which is free and finitely generated as both $\p[T]$-,
$\p[\tau]$-module and such that
$$ \exists m=m(M) \ \hbox{such that}\ (T-\theta)^m M/\tau M=0\eqno{(1.3.1)}$$
\medskip
{\bf Remark.} The above object is called "abelian t-motive" (resp. "t-motive") in [G] (resp. [A]), while the name "t-motive" is used in [G] for a more general object ([G], Definition 5.4.2). Since we shall not use objects defined in [G], 5.4.2, I prefer to use a shorter name for the above $M$.
\medskip
t-motives are main objects of the present paper. If we affirm
that an object exists this means that it exists as a t-motive if otherwise
is not stated. We denote dimension of $M$ over $\p[\tau]$ (resp. $\p[T]$) by $n$ (resp.
$r$), these numbers are called dimension and rank of $M$. Morphisms of abelian
t-motives are morphisms of left $\p[T, \tau]$-modules.

To define a left $\p[T, \tau]$-module $M$ is the same as to define a left $\p[T]$-module
$M_{\p[T]}$
endowed by an action of $\tau$ satisfying $\tau(Pm)=P^{(1)}\tau(m)$, $P\in \p[T]$. In
this situation
we can also treate $\tau$ as a $\p[T]$-linear map $M^{(1)}\to M$. This interpretation is
necessary if
we consider the general case $\w\supset \n F_q[T]$.

We need two categories which are larger than the category of abelian
t-motives.
\medskip
{\bf Definition 1.4.} A pr\'e-t-motive is a left $\p[T,
\tau]$-module which is free
and finitely generated as $\p[T]$-module, and satisfies (1.3.1).
\medskip
{\bf Definition 1.5.} A rational pr\'e-t-motive is a left
$\p(T)[\tau]$-module which is
free and finitely generated as $\p(T)$-module.
\medskip
{\bf Remark 1.6.} An analog of (1.3.1) does not exist for them.
\medskip
There is an obvious functor from the category of t-motives to
the category of
pr\'e-t-motives which is fully faithful, and an obvious functor from the category of
pr\'e-t-motives to the
category of rational pr\'e-t-motives. We denote these functors by
$i_1$, $i_2$
respectively. It is easy to see (Remark 10.2.3) that if $M$ is a
pr\'e-t-motive then the
action of $\tau$ on $i_2(M)$ is invertible.

Let $M_1$, $M_2$ be rational pr\'e-t-motives such
that the action of
$\tau$ on $(M_1)_{\p(T)}$ is invertible.
\medskip
{\bf Definition 1.7. (a)} $\Hom(M_1,M_2)$ is a rational pr\'e-t-motive such that

$$\Hom(M_1,M_2)_{\p(T)}=\Hom_{\p(T)}((M_1)_{\p(T)}, (M_2)_{\p(T)})$$
and the action of $\tau$ is defined by the usual manner: for
$\varphi:M_1 \to M_2$, $m\in M_1$
$$(\tau\varphi)(m)=\tau(\varphi(\tau^{-1}(m)))$$
\medskip
{\bf (b)} $M_1\otimes M_2$
\medskip
{\bf Definition 1.8.} Let $M$ be a t-motive and
$\mu$ a positive
number. A t-motive
$M'={M'}^{\mu}$ is called the $\mu$-dual of $M$ (dual if $\mu=1$) if
$M'=\Hom(M,\goth
C^{\otimes
\mu})$ as a rational pr\'e-t-motive, i.e. $$i_2\circ
i_1(M')=\Hom(i_2\circ i_1(M),\goth
C^{\otimes
\mu})\eqno{(1.8.1)}$$

{\bf Remark.} This definition generalizes the original one of Taguchi
([T], Section 5), see 1.12 below. A similar definition is in [F].
\medskip
{\bf 1.9.} We shall need the explicit matrix description of the above
objects. Let
$e_*=(e_1, ..., e_n)^t$ be the vector column of elements of a basis
of $M$ over $\p[\tau]$. There exists a matrix $\goth A\in M_n(\p[\tau])$ such that

$$T e_* = \goth A e_*, \ \ \goth A = \sum_{i=0}^l \goth A_i \tau^i \hbox{ where } \goth A_i
\in M_n(\p)\eqno{(1.9.1)}$$
Condition (1.3.1) is equivalent to the condition
$$\goth A_0=\theta I_n + N\eqno{(1.9.2)}$$
where $N$ is a
nilpotent matrix, and the
condition
$m(M)=1$ is equivalent to the condition $N=0$.
\medskip
Let $f_*=(f_1, ..., f_r)^t$ be the vector column of elements of a basis
of $M$ over
$\p[T]$. There exists a matrix $Q=Q(f_*)\in M_r(\p[T])$ such that
$$\tau f_* = Q f_*\eqno{(1.9.3)}$$

{\bf Lemma 1.10.} Let $M$ be as above. A t-motive
$M'$ is the $\mu$-dual of $M$ iff there exists a basis
$f'_*=(f'_1, ..., f'_r)^t$ of $M'$ over $\p[T]$ such that its matrix
$Q'=Q(f'_*)$ satisfies
$$Q'=(T-\theta)^{\mu}Q^{t-1} \eqno{(1.10.1)}$$
$\square$
\medskip
{\bf 1.10.2.} For further applications we shall need the following
lemma. The above
$f_*$, $f'_*$ are the dual bases (i.e. if we consider $f'_i$ as
elements of $\Hom(M,\goth
C)$ then $f'_i(f_j)=\delta^i_j\goth f$, where $\goth f$ is canonically defined by the condition that it generates $\goth C_{\p[T]}$ and satisfies $\tau \goth
f=(T-\theta)\goth f$). Let $\gamma$ be an endomorphism of $M$ and $D$ its
matrix in the basis
$f_*$ (i.e. $\gamma(f_*)=Df_*$). Let $\gamma'$ be the dual endomorphism.
\medskip
{\bf Lemma 1.10.3.} The matrix of $\gamma'$ in the basis $f'_*$ is $D^t$.
$\square$
\medskip
{\bf Remark 1.11.1.} For any $M$ having dual there exists a canonical homomorphism $\delta: \goth
C \to M\otimes M'$.
This is a well-known theorem of linear algebra. Really, in the above notations we have
$\goth f \mapsto \sum_i f_i \times f'_i$. It is obvious that $\delta$ is well-defined,
canonical and compatible
with the action of $\tau$.
\medskip
{\bf Remark 1.11.2.} The $\mu$-dual of $M$ --- if it exists --- is
unique, i.e. does not
depend on base change. This follows immediately from Definition 1.8,
but can be deduced
easily from 1.10.1. Really, let $g_*=(g_1, ..., g_r)^t$ be another basis
of $M$ over
$\p[T]$ and $C\in GL_r(\p[T])$ the matrix of base change (i.e. $g_*= C
f_*$). Then
$Q(g_*)=C^{(1)}QC^{-1}$. Let $g'_*=(g'_1, ..., g'_r)^t$ be a basis of
$M'$ over $\p[T]$
satisfying $g'_*= C^{t-1} f'_*$. Elementary calculation shows that
matrices $Q(g_*)$,
$Q(g'_*)$ satisfy (1.10.1).
\medskip
{\bf Remark 1.11.3.} The operation $M \mapsto {M'}^{\mu}$ is
obviously contravariant functorial. It is an exercise to the reader to give an exact definition of the corresponding category such that the functor of duality is defined on it, and is involutive (recall that not all t-motives have duals, and the dual of a map of t-motives is a priori a map of rational pr\'e-t-motives).
\medskip
{\bf 1.12.} The original definition of duality ([T], Definition 4.1; Theorem 5.1) from
the first
sight seems to be more
restrictive than the definition 1.8 of the present paper, but really
they are
equivalent. We recall some notations and definitions of [T] in a slightly less
general setting (rough statements; see [T] for the exact statements). Let $G$ be a finite affine group scheme over $\p$, i.e. $G=\Spec R$
where $R$ is a finite-dimensional $\p$-algebra. Let $\mu:R \to R\otimes R$ be the
comultiplication of $R$. Such group $G$ is called a finite $v$-module ([T], Definition
3.1) if there is a homomorphism $\psi: \w \to \End_{gr. \ sch.}(G)$ satisfying some natural conditions (for
example, an analog of 1.3.1). Further, let $\Cal E_G$ be a $\p$-subspace of $R$ defined as
follows: $$\Cal E_G=\{x\in R \ \vert \ \mu(x)=x\otimes 1+1\otimes x\}$$
The map $x\mapsto x^q$ is a $\p$-linear map $\fr: \Cal E_G^{(1)}\to \Cal E_G$. Further,
the map $\psi(T): G \to G$ can be defined on $\Cal E_G$. Let $v: \Cal E_G\to\Cal
E_G^{(1)}$ be a map satisfying $\fr \circ v= \psi(T)-\theta$.

We consider two finite $v$-modules $G$, $H$, the above objects $\fr$, $v$ etc. will carry
the respective subscript. Let * be the dual in the meaning of linear algebra.
\medskip
{\bf Definition 1.12.1} ([T], 4.1). Two finite $v$-modules $G$, $H$
are called dual if there exists an
isomorphism $\alpha: \Cal E^*_H\to\Cal E_G$ such that if we denote by
$\goth v: \Cal E_G\to \Cal E^{(1)}_G$ a map which enters in the commutative diagram
$$ \matrix \Cal E^*_H & \overset{\fr^*_H}\to{\longrightarrow} & \Cal E^{*(1)}_H \\ & & \\
\alpha \downarrow & & \alpha ^{(1)}\downarrow \\ & & \\ \Cal E_G & \overset{\goth
v}\to{\longrightarrow} & \Cal E^{(1)}_G \endmatrix $$
then we have:
$$\fr_G\circ \goth v= \psi_G(T)-\theta \eqno{(1.12.2)}$$
i.e. $\goth v=v_G$.
\medskip
Let $M$ be a t-motive having $m(M)=1$, $E=E(M)$ the corresponding t-module and
$a\in \hbox{\bf
A}$. We denote $E_a$ --- the set of $a$-torsion elements of $E$ --- by $M_a$. It is a
finite
$v$-module.
\medskip
{\bf Proposition 1.12.3.} Let $M$, $M'$ be t-motives which are
dual in the meaning
of Definition 1.8. Then $\forall a\in \hbox{\bf A}$, $a\ne 0$ we have: $M_a$,
$M'_a$ are dual in the meaning of 1.12.1 = [T], Definition 4.1.
\medskip
{\bf Proof.} Condition $a\in \n F_q[T]$ implies that multiplication by $\tau$ is
well-defined on $M/aM$.
\medskip
{\bf Lemma 1.12.3.1.} We have canonical isomorphisms $i: M/aM \to
\Cal E_{M_a}$,
$i^{(1)}: M/aM \to \Cal E^{(1)}_{M_a}$ such that the following
diagrams are commutative:
$$ \matrix M/aM & \overset{\tau}\to{\longrightarrow} & M/aM &&& M/aM &
\overset{T}\to{\longrightarrow} & M/aM
\\ & & &&&&&
\\ i^{(1)}\downarrow & & i\downarrow &&& i\downarrow & & i\downarrow \\ & & &&&&&
\\ \Cal E^{(1)}_{M_a} & \overset{fr}\to{\longrightarrow} & \Cal E_{M_a} &&&
\Cal E_{M_a} & \overset{\psi_T}\to{\longrightarrow} & \Cal E_{M_a}
\endmatrix $$
\medskip
{\bf Proof.} Let $R$ be a ring such that $\Spec R=M_a$. The pairing between $M$ and $E$
shows that there exists a map $M\to R$ which is obviously factorized via an inclusion
$M/aM\to R$. It is easy to see that the image of this inclusion is contained in $\Cal
E_{M_a}$, i.e. we get $i$. Since $\dim_{\p}(M/aM)=\deg a \cdot r(M)$ and
$\dim_{\p}(R)=q^{\deg a \cdot r(M)}$ we get from [T], Definition 1.3 that $i$ is an
isomorphism. Other statements of the lemma are obvious. $\square$
\medskip
This lemma means that we can rewrite Definition 1.12.1 for the case $G=M_a$, $H=N_a$ by
the following way:\footnotemark \footnotetext{Here and below a t-motive $N$ should not be confused with $N$ of 1.9.2.}
\medskip
{\bf 1.12.3.2.} Two finite $v$-modules $M_a$, $N_a$ are dual if there
exists an
isomorphism $\alpha: (N/aN)^* \to M/aM $ such that after identification
via $\alpha$ of
$\tau^*: (N/aN)^* \to (N/aN)^*$ with a map $\goth v: M/aM \to M/aM$ we
have on $M/aM$:
$$\tau\circ \goth v= t-\theta \eqno{(1.12.3.3)}$$
We need a
\medskip
{\bf Lemma 1.12.3.4.} For $i=1,2$ let $N_i$ be a free $\p[T]$-module of
dimension $r$ with a
base $f_{i*} =(f_{i1}, ... , f_{ir} )$, let $\varphi_i:N_i \to N_i$ be
$\p[T]$-linear
maps having matrices $\goth Q_i$ in $f_{i*} $ such that $\goth
Q_2=\goth Q_1^t$, and let
$a$ be as above. Let, further, $\varphi_{i,a}: N_i/aN_i \to N_i/aN_i $
be the natural
quotient of $\varphi_i$. Then there exist $\p$-bases $\tilde f_{i*} $
of $N_i/aN_i $ such
that the matrix of $\varphi_{1,a}$ in the base $\tilde f_{1*} $ is
transposed to the
matrix of $\varphi_{2,a}$ in the base $\tilde f_{2*} $.
\medskip
{\bf Proof.} We can identify elements of $N_2$ with $\p[T]$-linear
forms on $N_1$
(notation: for $x\in N_2$ the corresponding form is denoted by
$\chi_x$) such that
$\chi_{\varphi_2(x)}=\chi_x\circ\varphi_1$. Any $\p[T]$-linear form
$\chi$ on $N_i$
defines a $\p[T]/a\p[T]$-linear form on $N_i/aN_i$ which is denoted by
$\chi_a$. Let now
$x\in N_2/aN_2$, $\bar x$ its lift on $N_2$, then
$\chi_{x,a}=(\chi_{\bar x})_a$ is a
well-defined $\p[T]/a\p[T]$-linear form on $N_1/aN_1$. For $x\in
N_2/aN_2$ we have
$$\chi_{\vf_{2,a}(x),a}= \chi_{x,a}\circ \vf_{1,a}$$
Further, let $\lambda: \p[T] \to \p$ be a $\p$-linear map such that
\medskip
{\bf 1.12.3.5.} Its kernel does not contain any non-zero ideal of
$\p[T]/a\p[T]$.
\medskip
(such $\lambda$ obviously exist.) For $x\in N_2/aN_2$ we denote $\lambda \circ
\chi_{x,a}$ by $\psi_x$, it is a $\p$-linear form on $\p$-vector space
$N_1/aN_1$.
Obviously condition (1.12.3.5) implies that the map $x\mapsto \psi_x$
is an isomorphism
from $N_2/aN_2$ to the space of $\p$-linear forms on $\p$-vector space
$N_1/aN_1$, and we
have
$$\psi_{\vf_{2,a}(x)}=\psi_x\circ\vf_{1,a}$$
which is equivalent to the statement of the lemma. $\square$
\medskip
Finally, the proposition follows immediately from this lemma
multiplied by $T-\theta$,
formula 1.10.1 and 1.12.3.2. $\square$
\medskip
{\bf Remark.} Let $a=\sum_{i=0}^k g_iT^i$, $g_i\in \n F_q$, $g_k=1$.
Taguchi ([T], proof
of 5.1 (iv)) uses the following $\lambda$: $\lambda(T^j)=0$ for $j<k-1$,
$\lambda(T^{k-1})=1$. It is easy to check that for $x=(T^i+
T^{i-1}g_{k-1}+T^{i-2}g_{k-2}+ ...
+g_{k-i})f_{2j}$ for this $\lambda$ we have: $\psi_x(T^i f_{1j})=1$,
$\psi_x(T^{i'}
f_{1j'})=0$ for other $i'$, $j'$.
\medskip
{\bf 1.13.} We consider in Sections 1.13, 1.14 the case of arbitrary $\w\supset \n F_q[T]$.
\medskip
A t-motive over {\bf A} is defined for example in [BH], p.1. Let us reproduce this
definition for the case of characteristic 0. Let $J$ be an ideal of
$\w_C$ generated by the elements $a\otimes 1 - 1 \otimes \iota(a)$ for all $a\in \w$. The ring $\w_C[\tau]$ is defined by the formula $\tau\cdot (a\otimes z)=(a\otimes z^q)\cdot \tau$, $a\in \w$, $z\in \p$.
\medskip
{\bf Definition 1.13.1.} A t-motive $M$ over {\bf A} is a pair $(M, \tau)$
where $M$ is
a locally free $\hbox{\bf A}_C$-module and $\tau$ is an $\hbox{\bf A}_C$-linear map
$M^{(1)}\to M$ satisfying the following analog of 1.3.1, 1.9.2:
$$\exists m \hbox{ such that } J^m(M/\tau (M^{(1)}))=0\eqno{(1.13.2)}$$

{\bf Remark 1.13.3.} We can consider $M$ as an $\w_C[\tau]$-module using the following formula for the product $\tau \cdot m$:
$$\tau \cdot m = \tau(m\otimes 1)$$
where $m\in M$, $m\otimes 1\in M^{(1)}$.
\medskip
The rank of $M$ as a locally free $\hbox{\bf A}_C$-module is called the rank of the
corresponding t-motive $(M, \tau)$. If $\w=\n F_q[T]$ then $M^{(1)}$ is isomorphic to
$M$, we can
consider $M$ as a $\p[T,\tau]$-module, and it is possible to show that in this case
1.13.2 implies
that $M_{\p[\tau]}$ is a free $\p[\tau]$-module. In the general case, the dimension $n$
of $(M, \tau)$
is defined as $\dim_{\p}(M/\tau (M^{(1)}))$.

Let us fix $\goth C=(\goth C, \tau_{\goth C})$ --- a t-motive of rank 1 over
{\bf A}.
For a t-motive $M=(M,\tau_M)$ a t-motive ${M_\goth C'}$ --- the
$\goth C$-dual of $M$
--- is defined as follows. We put
$M_\goth C'=\Hom_{\w_C}(M,\goth C)$.
Since for any locally free $\hbox{\bf A}_C$-modules $M_1$, $M_2$ we have
$$\Hom_{\w_C}(M_1,M_2)^{(1)}=\Hom_{\w_C}(M_1^{(1)},M_2^{(1)})$$
we can define $\tau(M_\goth C')$ by the following formula:

$$\hbox{ For } \vf\in \Hom_{\w_C}(M,\goth C)^{(1)} \hbox{ we have }
\tau(M_\goth C')(\vf)=\tau_{\goth C}\circ \vf \circ \tau_M^{-1}$$
\medskip
{\bf 1.14. Duality for abelian $\tau$-sheaves.} We use notations of [BH],
Definition 2.1 if they do not differ from the notations of the present
paper; otherwise we continue to use notations of the present paper
(for example, $d$ (resp. $\sigma^*(\goth X)$ for any object $\goth X$)
of [BH] is
$n$ (resp. $\goth X^{(1)}$) of the present paper). For any abelian
$\tau$-sheaf $\underline{\Cal F}$ we denote its $\Pi_i$, $\tau_i$ by
$\Pi_i(\underline{\Cal F})$, $\tau_i(\underline{\Cal F})$
respectively. If $M$, $N$ are invertible sheaves on $X$ and $\rho: M
\to N$ a rational map then we denote by $\rho^{inv}: N \to M$ the
rational map which is inverse to $\rho$ with respect to the composition.
We define
$\tau_{\goth r,i-1}(\underline{\Cal F})$ (the rational $\tau_i$) as
the composition map $\tau_{i-1}(\underline{\Cal F})
\circ {\Pi_{i-1}^{(1)}}^{inv}(\underline{\Cal F})$, it is a rational
map from $\Cal F_i^{(1)}$ to $\Cal F_i$.

Let $\underline{\Cal O}$ be a fixed abelian
$\tau$-sheaf having $r=n=1$. The $\underline{\Cal O}$-dual abelian
$\tau$-sheaf $\underline{\Cal F}'= \underline{\Cal
F}'_{\underline{\Cal O}}$ is defined by the formulas
$$\Cal F'_0=\Hom_{X}(\Cal F_0, \Cal O_0)$$ where Hom is the sheaf's
one, and the map $\tau_{\goth r,-1}(\underline{\Cal F}'): {\Cal
F'_0}^{(1)}\to \Cal F'_0$ is defined as follows. We have ${\Cal
F'_0}^{(1)}= \Hom_{X}({\Cal F_0}^{(1)}, {\Cal O_0}^{(1)})$. Let
$\gamma\in\Hom_{X}({\Cal F_0}^{(1)}, {\Cal O_0}^{(1)})(U)$ where $U$
is a sufficiently small affine subset of $X_{\p}$, such that $\gamma:
{\Cal F_0}^{(1)}(U) \rightarrow {\Cal O_0}^{(1)}(U)$.
\medskip
{\bf 1.14.1.} We define: $[[\tau_{\goth r,-1}(\underline{\Cal
F}')](U)](\gamma)$ is the following composition map:
$$\Cal F_0(U)\overset{[\tau_{\goth r,-1}^{inv}(\underline{\Cal
F})](U)}\to{\longrightarrow}
\Cal F_0^{(1)}(U)\overset{\gamma}\to{\to}\Cal O_0^{(1)}(U)
\overset{[\tau_{\goth r,-1}(\underline{\Cal
O})](U)}\to{\longrightarrow}\Cal O_0(U)\in \Hom_{X}({\Cal F_0}, {\Cal
O_0})(U)$$

Clearly that this definition and the definitions 1.8, 1.13 are
compatible with the forgetting functor $\underline{M}(\underline{\Cal
F})$ from abelian $\tau$-sheaves to pure Anderson t-motives of [BH],
Section 3, page 8.
\medskip
{\bf 1.15. Duality over fields.} Let $L\supset \n F_q(\theta)$ be a field extension of $\n F_q(\theta)$, and $M$ a t-motive over $L$ (i.e. a pair ($M$, an $L$-structure on $M$)). Obviously we have
\medskip
{\bf Proposition 1.15.1.} The notion of duality for $M$ over $L$ is well-defined. $\square$
\medskip
Similarly, we have a proposition for Galois action:
\medskip
{\bf Proposition 1.15.2.} Let $M$ be defined over $\overline{ \n F_q(\theta)}$ and $\gamma\in \Gal( \n F_q(\theta))$. Then $(\gamma(M))'=\gamma(M')$. $\square$
\medskip
\medskip
{\bf 1.16. Case of finite characteristic.} Let $\iota: \w\to \bar \n F_q$ be a map of finite characteristic, we denote $\Ker \iota$ by $\Cal P$. The definition of t-motive for this case is similar to 1.3, see [G] for the details. The definition of duality also is similar to the one of the case of generic characteristic. Duality commutes with reduction. Namely, let $M$ be from 1.15, $\goth P$ a prime of $L$ not over the infinity of $\n F_q(\theta)$, $\Cal P\subset \w$ is $\iota^{-1}(\goth P\cap \n F_q[\theta]$) --- the finite characteristic. We consider the case of good reduction of $M$ at $\goth P$, we denote it by $\tilde M$. It is a t-motive in characteristic $\Cal P$.  Let $M$ have dual $M'$.
\medskip
{\bf Proposition 1.16.1.} $\tilde M$ has dual iff $M'$ has good reduction at $\goth P$; in this case they coincide. $\square$
\medskip
{\bf Remark 1.16.2.} Apparently if $M$ has good reduction and dual, then $M'$ also has good reduction (in this case 1.16.1 means that $M'$ exists implies $(\tilde M)'$ exists). For standard-3 t-motives (this is a simple tipe of t-motives, see 11.8.1) apparently this can be shown by explicit calculations.
\medskip
{\bf Remark 1.16.3.} Clearly 1.16.1 is true for the case of bad reductions. I do not give exact definitions for this case.
\medskip
{\bf 1.16.4. Ordinarity.} Let $M$ be of finite characteristic. By analogy with the number field case, $M$ is called ordinary if its Newton polygon consists of 2 segments. If $N=0$ then the Newton polygon of $M'$ is the dual of the one of $M$ (the notion of duality of polygons is clear; apparently the condition $N=0$ can be omitted). So, we have
\medskip
{\bf Proposition 1.16.5.} $M$ is ordinary $\iff M'$ is ordinary. $\square$
\medskip
See 13.4.1 for a more exact result.
\medskip
{\bf 2. Analytic duality.}
\medskip
We consider in the present section the case of arbitrary $\w\supset \n F_q[T]$ (and $N=0$ as usually).
\medskip
Condition $N=0$ implies that an element $a\in \w$ acts on $\Lie(M)$ by multiplication by $\iota(a)$. Hence, we have a
\medskip
{\bf Definition 2.1.} Let $V$ be the space $\p^n$. A locally free $r$-dimensional
$\n Z_\infty$-submodule
$L$ of $V$ is called a lattice if

(a) $L$ generates $V$ as a $\p$-module and

(b) The $\n R_\infty$-linear span of $L$ has dimension $r$ over
$\n R_\infty$.
\medskip
Numbers $n$, $r$ are called the dimension and the rank of $L$
respectively. Attached to $(L,V)$ is the tautological inclusion $\vf=\vf(L,V):
L \to V$. We shall consider the category of triples $(\vf, L, V)$; a map $\psi:
(\vf, L,V)\to(\vf_1, L_1,V_1)$ is a pair $(\psi_L, \psi_V)$ where $\psi_L: L \to L_1$ is a $\n Z_\infty$-linear map, $\psi_V: V\to V_1$ is a $\p$-linear map such that $\vf_1\circ \psi_L=\psi_V \circ \vf$.

Inclusion $\vf$ can be extended to a map
$L\underset{\n Z_\infty}\to{\otimes}\p \to V$ (which is surjective
by 2.1a), we denote it by $\vf=\vf(L,V)$ as well.
We can also attach to $(L,V)$ an exact sequence
$$0\to \Ker \vf \to L\underset{\n Z_\infty} \to{\otimes}\p \overset{\vf}\to{\to} V \to
0\eqno{(2.2)}$$

Let $\Cal I\in \Cl(\w)$ be a class of ideals; we shall use the same notation $\Cal I$ to
denote a representative in the $\iota$-image of this class. Let $(\vf', L', V')$ be another lattice and $D$ a structure of a perfect $\Cal I$-pairing $<* , * >_D$ between $L$
and $L'$. Let us fix an isomorphism $$\alpha: \Cal I \underset{\n Z_\infty} \to{\otimes} \p \to \p\eqno{(2.2')}$$ $D$ extends via $\alpha$ to a perfect $\p$-pairing between $L\underset{\n Z_\infty}\to{\otimes}\p$
and $L'\underset{\n Z_\infty}\to{\otimes}\p$, we denote this pairing by $D_{\alpha, \infty}$.
\medskip
{\bf Definition 2.3.} Two lattices $(\vf, L, V)$ and $(\vf', L', V')$ are called $(\alpha, \Cal I)$-dual if there exists a perfect $\Cal I$-pairing $D$ between $L$
and $L'$ such that $\Ker \vf \subset L\underset{\n Z_\infty}\to{\otimes}\p$, $\Ker \vf' \subset L'\underset{\n Z_\infty}\to{\otimes}\p$ are mutually orthogonal with respect to $D_{\alpha, \infty}$.
\medskip
Let $(n,r)$, $(n',r')$ be the dimension and rank of $(\vf, L, V)$ and $(\vf', L', V')$ respectively. If they are $(\alpha, \Cal I)$-dual then $r'=r$, $n'=r-n$. There exists the following reformulation of the definition of duality. $D_{\alpha, \infty}$ induces an isomorphism
$\gamma_{\alpha, D}: (L\underset{\n Z_\infty}\to{\otimes}\p)^* \to
L'\underset{\n Z_\infty}\to{\otimes}\p$ (here and below for any object $W$ we
denote $W^*=\Hom_{\p}(W,\p)$ ).
\medskip
{\bf Property 2.4.} $(\vf, L, V)$ and $(\vf', L', V')$ are $(\alpha, \Cal I)$-dual iff there exists an isomorphism from $(\Ker \vf)^*$
to $V'$ making the following diagram commutative: $$\matrix 0 & \to &
V^* & \overset{\vf^*}\to{\to} &  (L\underset{\n Z_\infty}
\to{\otimes}\p)^* & \to & (\Ker \vf)^* & \to & 0\\
&&&&&&&& \\
& & \downarrow & & \gamma_{\alpha, D} \downarrow && \downarrow \\ &&&&&&&& \\ 0
& \to & \Ker \vf' & \rightarrow  & L'\underset{\n Z_\infty}\to{\otimes}\p &
\overset{\vf'}\to{\to} & V' & \to &
0\endmatrix \eqno{(2.5)}$$

Further, this property is equivalent to the following two conditions:
\medskip
{\bf 2.6.} $\dim V'=r-n$;
\medskip
{\bf 2.7.} The composition map $\vf'\circ \gamma_D \circ \vf^*:
V^* \to V'$ is 0.
\medskip
Both 2.4 and (2.6, 2.7) are obvious.
\medskip
{\bf Remark 2.8.} It is
easy to see that the functor $(\vf, L, V) \mapsto (\vf', L', V')$ is well-defined
on a subcategory
(not all lattices have duals, see
below) of the category of the triples $(\vf, L, V)$, it is contravariant and
involutive.
\medskip
{\bf 3. Explicit formulas for analytic duality.}
\medskip
Here we consider the case $\w=\n
F_q[T]$. In this case $\Cl(\w)=0$, and $(\alpha, \Cal I)$-dual is called simply
dual. The coordinate
description of the dual lattice is the following. Let
$e_1, ..., e_r$ be a $\n Z_\infty$-basis of
$L$ such that $\vf(e_1), ..., \vf(e_n)$ form a $\p$-basis of $V$. Like in the
theory of abelian
varieties, we denote by $Z=(z_{ij})$ the Siegel matrix whose lines are
coordinates of $\vf(e_{n+1}), ..., \vf(e_r)$ in the basis $\vf(e_1), ..., \vf(e_n)$, more
exactly, the size of $Z$ is $(r-n)\times n$ and
$$\forall i =1,..., r-n \ \ \ \ \vf(e_{n+i})=\sum_{j=1}^n z_{ij}\vf(e_j)\eqno{(3.1)}$$ $Z$
defines $L$, we denote $L$ by $\goth L(Z)$.
\medskip
{\bf Proposition 3.2.} $[\goth L(Z)]'=\goth L(-Z^t)$, i.e. a Siegel matrix of the dual
lattice is the minus transposed Siegel matrix.
\medskip
{\bf Proof.} Follows immediately from the definitions. Really, let $f_1, ..., f_r$ be a
basis of $L'$, we define the pairing by the formula $$<e_i,
f_j>=\delta_i^j\eqno{(3.3)}$$ and the map $\vf'$ by the formula
$$\forall i =1,... ,n \ \ \ \ \vf'(f_{i})=\sum_{j=1}^{r-n} -z_{ji}\vf'(f_{n+j})$$
(minus transposed Siegel matrix).
$\Ker \vf$ is generated by elements $$v_i=e_{n+i}-\sum_{j=1}^n z_{ij}e_j, \ \ \ \ i
=1,..., r-n$$ and $\Ker \vf'$ is generated by elements $$w_i=f_{i}+\sum_{j=1}^{r-n}
z_{ji}f_{n+j}, \ \ \ \ i =1,..., n\eqno{(3.4)}$$ It is sufficient to check that $\forall i,j$ we have
$<v_i, w_j>=0$; this follows immediately from 3.3. $\square$
\medskip
{\bf Remark 3.5.} $L'$ exists not for all $L$.
Trivial counterexample: case $n=r=1$. To get another counterexamples, we use that for
$n=1$ (lattices of Drinfeld modules) a Siegel matrix is a column matrix $Z=\left(\matrix
z_1&...&z_{r-1} \endmatrix \right)^t$ and
$$ \goth L(Z) \hbox{ is not a lattice } \iff 1,z_1, ... ,z_{r-1} \hbox{ are linearly
dependent over } \n R_\infty \eqno{(3.6)}$$ while for $n=r-1$ a Siegel matrix is a
row matrix $Z=\left(\matrix -z_1&...&-z_{r-1} \endmatrix \right)$ and
$$ \goth L(Z) \hbox{ is not a lattice } \iff \forall i \ \ z_i\in \n R_\infty
\eqno{(3.7)}$$ Since condition (3.7) is strictly stronger than (3.6) we see
that all lattices having $n=1$, $r>1$ have duals while not all lattices having $n=r-1$,
$r>2$ have duals.

It is clear that almost all matrices have duals. Here "almost all" has the same meaning
that as "Almost all matrices $Z$ are a Siegel matrice of a lattice", i.e. if we choose an
(infinite) basis of $\p/\n R_\infty$, then coordinates of the entries of $Z$ in this
basis must satisfy some polynomial relations in order that $Z$ is not a Siegel matrice of
a lattice.
\medskip
{\bf Remark 3.8.} The coordinate proof of the theorem that the notion
of the dual lattice is well-defined, is the following. Two Siegel matrices
$Z$, $Z_1$ are called equivalent iff there exists an isomorphism of their pairs $(\goth
L(Z), V)$, $(\goth L(Z_1),V_1)$. Like in the classical theory of modular forms, $Z$,
$Z_1$ are equivalent iff there exists a matrix $\gamma \in GL_r(\n Z_\infty)=\left(\matrix A&B\\ C&D \endmatrix \right)$ ($A,B,C,D$ are the ($n\times
n$), ($n\times r-n$), ($r-n\times n$), ($r-n\times r-n$)-blocks of $\gamma$ respectively;
we shall call this block structire by the $(n, r-n)$-block structure) such that
$$C+DZ=Z_1(A+BZ)\eqno{(3.8.1)}$$

Let $A_1,B_1, C_1, D_1$ be the $(n, r-n)$-block structure of the matrix $\gamma^{-1}$.
The equality
$$-C_1^t+A_1^tZ^t={Z_1}^t(D_1^t-B_1^tZ^t)\eqno{(3.8.2)}$$
shows that if $Z$, $Z_1$ are equivalent
then $-Z^t$, $-Z_1^t$ are equivalent. [Proof of (3.8.2): (3.8.1) implies $Z_1=(C+DZ)(A+BZ)^{-1}$; substituting this value of $Z_1$ to the transposed
(3.8.2), we get $-C_1+ZA_1=(D_1-ZB_1)(C+DZ)(A+BZ)^{-1}$, or $(-C_1+ZA_1)(A+BZ)=(D_1-ZB_1)(C+DZ)$. This formula follows immediately from
$\left(\matrix A_1&B_1\\C_1&D_1\endmatrix \right)\left(\matrix A&B\\C&D\endmatrix \right)=\left(\matrix I_n&0\\0&I_{r-n}\endmatrix \right)$].

Further, let $\alpha: (L_1\subset \p^n) \to
(L_2\subset \p^n)$ be a map of lattices. If $L'_1$, $L'_2$
exist, then the map
$\alpha': (L'_2\subset \p^{r-n}) \to (L'_1\subset \p^{r-n})$ is defined by the
following formulas. Let $Z_i$ be the Siegel matrices of $L_i$ in the
bases
$e_{i1}, ... e_{ir}$ of $L_i$ ($i=1,2$). Let us consider the matrix
$\goth M=(m_{ij})\in
M_{r}(\n Z_\infty)$ of $\alpha$ in the bases $e_{i1}, ...,
e_{ir}$ (i.e.
$\alpha(e_{1i})=\sum_j m_{ij} e_{2j}$). Let $f_{i1}, ..., f_{ir}$ be the dual
base of $L'_i$ (see 3.3) and $e'_{i1}, ... e'_{ir}$ another base of $L'_i$ defined by $$e'_{ij}=f_{i,j+n}, \ \ \ \ j+n \mod r\eqno{(3.8.3)}$$ Formulas (3.8.3), (3.4) show that an analog of 3.1 is satisfied for both bases $e'_{i1}, ..., e'_{ir}$, their Siegel matrices are $-Z_i^t$.

Let
$$\goth M=\left(\matrix \goth M_{11} & \goth M_{12} \\ \goth M_{21} &
\goth M_{22}
\endmatrix \right)$$
be the $(n, r-n)$-block structure of $\goth M$. The matrix of $\alpha'$ in the bases $f_{i1}, ..., f_{ir}$ is $\goth M^t$, and using the matrix 3.8.3 of change of base, we get that $\goth M'$ --- the matrix of $\alpha'$ in
the bases $e'_{i1}, ..., e'_{ir}$ --- has the following $(r-n,
n)$-block structure:
$$\goth M'=\left(\matrix \goth M_{22}^t & \goth M_{12}^t \\
\goth M_{21}^t & \goth
M_{11}^t \endmatrix \right)\eqno{(3.8.4)} $$
The property that $\goth M$ comes from a $\p$-linear map $\p^n \to \p^n$ implies
that $\goth M'$ comes from a $\p$-linear map $\p^{r-n} \to
\p^{r-n}$. This follows immediately from the definition of dual lattice, or can be easily checked algebraically.
\medskip
{\bf Remark 3.9.} Taking $\gamma =\left(\matrix 1&0\\ 0&-1 \endmatrix \right)$ we get
that $Z$ is equivalent to $-Z$, hence $Z'$ is also a Siegel matrix of the dual lattice.
\medskip
{\bf 4. Main conjecture for arbitrary $\w$}.
\medskip
The main result of the paper is the following Theorem 5 on coincidence of algebraic and
analytic duality. We formulate it as a conjecture 4.1 for any $\w$, but we prove it only for the case $\w=\n
F_q[T]$.  Let $M$ be a uniformizable t-motive. Its lattice $L(M)$ is really a lattice in
the meaning of Definition 2.1, because [A],
Corollary 3.3.6 (resp. [G], Lemma 5.9.12) means that it satisfies 2.1a (resp. 2.1b);
recall that we consider the case $N=0$,
i.e. the action of $T$ on $\Lie(M)$ is simply multiplication by $\theta$.
Let us fix (like in 1.13) $\goth C=(\goth C, \tau_{\goth C})$ --- a t-motive of rank 1 over
{\bf A}, and let $L(\goth C)$ be its lattice. It is a $\n Z_\infty$-module. $\Omega=\Omega(\w)$ is an $\w$-module, we consider a $\n Z_\infty$-module $\iota^{-1}(\Omega)$. There exists the   notion of the $L(\goth C)\otimes \iota^{-1}(\Omega)$-duality.
\medskip
{\bf Conjecture 4.1.} Let $M$ be a uniformizable t-motive having $N=0$ such
that its $\goth C$-dual $M'$ exists. Then $M'$ is uniformizable, it has $N':=N(M')=0$, and $(L(M), \Lie(M))$ and $(L(M'), \Lie(M'))$ are $\alpha, L(\goth C)\otimes \iota^{-1}(\Omega)$-dual for some $\alpha$ from $2.2'$ (it can be explicitly described).
\medskip
We prove in Section 5 the first step of the proof of this conjecture.
\medskip
{\bf Remark 4.2.} It is possible to generalize the above conjecture to the case of non-uniformizable $M$, $M'$. The pairing is defined between $\Hom_{\w_C[\tau]}(M,Z_1)$ and $\Hom_{\w_C[\tau]}(M',Z_1)$ (see (5.2.1a) for the definition of $Z_1$), or, the same, between $M_a$ and $M'_a$ for any $a\in \w$ (see 5.1.6).
\medskip
{\bf 5. Main theorem.}
\medskip
Recall that the word "canonical" means "canonical up to multiplication by elements of $\n F_q^*$".
\medskip
{\bf Theorem 5.}\footnotemark \footnotetext{The proof of this theorem was inspired by a result of Anderson, see Section 6 for details.} Let $M$ be a uniformizable t-motive over $\w=\n F_q[T]$ having $N=0$ such
that its dual $M'$ exists and has $N':=N(M')=0$. Then $M'$ is uniformizable, and $(L(M), \Lie(M))$ and $(L(M'), \Lie(M'))$ are dual.
\medskip
{\bf Remark 5A.} Condition $N'=0$ holds for pure $M$ (Theorem 10.3) and for a large class of non-pure $M$ (Theorem 11.5). Most likely, a modification of the end of the proof of the present theorem will permit us to prove that $N'=0$ holds for all $M$ having $N=0$ and having dual.
\medskip
{\bf Remark 5B.} A reformulation of the theorem in terms of Hodge-Pink structures is given in Section 9.
Proof of the theorem for the case $N\ne0$ is given in [GL18].
\medskip
{\bf Corollary 5.1.1.} If $\w=\n F_q[T]$ then a Siegel matrix of $M'$ is the minus transposed
of a Siegel matrix of $M$.
\medskip
In the section 8 below we give a corollary of this theorem and some conjectures related to the
problem of 1 -- 1 correspondence between t-motives and lattices.
\medskip
{\bf 5.1.2. Some definitions.} Recall that $E=E(M)$ is isomorphic to $\p^n$. There is a structure of $\w$-module on $E$; multiplication by $T$ is denoted by $m_T$, and this operator $m_T$ is defined in coordinates by the formula
$$m_T(x)=\sum_{i=0}^l\goth A_ix^{(i)}$$ where $x\in E=\p^n$ is a vector column, $\goth A_i$ are from 1.9.1. There is a map $\exp: \Lie(M) \to E$ making the following diagram commutative:
$$\matrix \Lie(M) &
\overset{\Exp}\to{\to} & E \\ \\ \theta \downarrow & & m_T
\downarrow \\ \\ \Lie(M) & \overset{\Exp}\to{\to} & E \endmatrix \eqno{(5.1.3)}$$
By definition, $L(M)=\Ker \Exp$.

We need another space $\Lie_T(M)$ together with an isomorphism $\goth a:\Lie_T(M) \to \Lie(M)$ and a structure of $\w$-module on $\Lie_T(M)$ such that the multiplication by $T$ on $\Lie_T(M)$ is simply the multiplication by $\theta$ on $\Lie(M)$, i.e. $$\goth a(Tx)=\theta\cdot(\goth a(x))\eqno{(5.1.4)}$$ where $x\in\Lie_T(M)$. Commutativity of 5.1.3 means that $\Exp\circ\goth a:\Lie_T(M)\to E$ is a map of $\w$-modules.
\medskip
{\bf 5.1.5.} We shall work merely with $L_T(M):=\Ker (\Exp\circ\goth a)\subset \Lie_T(M)$ rather than $L(M)$. Clearly $L_T(M)$ is an $\w$-module, $\goth a:L_T(M)\to L(M)$ is an isomorphism satisfying 5.1.4 for $x\in L_T(M)$.
\medskip
The proof of Theorem 5 consists of two steps. We formulate and prove Step 1 for the case of arbitrary $\w$.
\medskip
{\bf Step 1.} For the above $M$, $M'$ we have:
\medskip
(A) Uniformizability of $M$ implies uniformizability of $M'$.
\medskip
(B) There exists a canonical $\w$-linear $L_T(\goth C)\otimes
\Omega$-valued perfect pairing $<* , * >_M$
between $L_T(M)$ and $L_T(M')$ (by 5.1.5, this is the same as the $\n Z_\infty$-linear pairing between $L(M)$ and $L(M')$, which, in its turn, is $D$ of Definition 2.3). It is functorial.
\medskip
{\bf Remark 5.1.6.} Practically, (B) comes from [T], Theorem 4.3 (case
$\w=\n F_q[T]$). Really, to define a pairing between $L(M)$ and
$L(M')$ it is sufficient to define (concordant) pairings between $L(M)
/aL(M)$ and $L(M') /aL(M')$ for any $a\in \w$. Since $M_a:=E(M)_a=L(M)/aL(M)$ and because of Proposition 1.12.3 which affirms that
$M_a$ and $M'_a$ are Taguchi-dual, we see that [T], Theorem 4.3 gives
exactly the desired pairing.
\medskip
We give two versions of the proof of Step 1: the first one --- for the general case of arbitrary $\w$ and the second one --- for the case $\w=\n F_q[T]$ --- is based on explicit calculations, it is used for the proof of Step 2.
\medskip
{\bf 5.2. Proof: Step 1, Version 1.} Here we consider the general case of arbitrary $\w$.
Let $\Omega=\Omega(\w/\n
F_q)$ be the module of differential forms; we can consider it as an
element of $\Cl(\w)$. We use formulas and notations of [G], Section 5.9
modifying them to the case of arbitrary $\w$. For example, {\bf
A} (resp. {\bf K}) of [G], 5.9.16 is {\bf A} (resp. {\bf K})
of the present paper (recall that $\bar K$ (resp. $\bar K[T,\tau]$) of [G] is
$\p$ (resp. $\w_C[\tau]$, see 1.13) of the present paper). Hence, we denote $\bar K\{T\}$ of
[G], Definition 5.9.10 by $\p\{T\}$. For the general case it must be
replaced by a ring $Z_0$ defined by the formula
$$Z_0:=\w\underset{\n F_q[T]}\to{\otimes}\p\{T\}\eqno{(5.2.1)}$$
$Z_0$ is a $\w_C[\tau]$-module, i.e. $\tau$ acts on $Z_0$, and $Z_0^\tau=\w$.

$Z_1$ for the present case is defined by the same formula [G], 5.9.22. Explicitly,
$$Z_1:=\Hom^{cont}_{\w}(\x/\w,\p)\eqno{(5.2.1a)}$$
It is a locally free $Z_0$-module of dimension 1 (the module structure
is compatible with the action of $\tau$; see [G], p. 168, lines 3 - 4
for the case $\w=\n F_q[T]$). We have: $Z_1^\tau$ is a
$Z_0^\tau$-module ( = $\w$-module) which is isomorphic to $\Omega(\w)$
(see the last lines of the proof of [G], Corollary 5.9.35 for the case
$\w=\n F_q[T]$), and $Z_1$ is isomorphic to $Z_0\otimes_{\w}\Omega(\w)$.

We shall consider $M$ as a $\w_C[\tau]$-module, like in 1.13.3. We denote
$M\{T\}:=M\underset{\w_C}\to{\otimes}Z_0$ ( = [G],
Definition 5.9.11.1 for the case $\w=\n F_q[T]$) and
$H^1(M):=M\{T\}^\tau$ like in [G], Definition 5.9.11.2.
Analogous to [G], Corollary 5.9.25 we get that for the present case
$$H_1(M):=\Hom_{\w_C[\tau]}(M,Z_1)=L_T(M)$$
($H_1(M)=H_1(E)$ of [G], 5.9). Particularly, for $M=\goth C$ we have
$$L_T(\goth C)=\Hom_{\w_C[\tau]}(\goth C,Z_1)$$

{\bf Lemma 5.2.2. } $H_1(M')=H^1(M)\underset{\w}\to{\otimes}L_T(\goth C)$.
\medskip
{\bf Proof.} By definition,
$\Hom_{\w_C} (M',Z_1)=\Hom_{\w_C} (\Hom_{\w_C} (M, \goth C), Z_1)$. Further,
$$\Hom_{\w_C} (\Hom_{\w_C} (M, \goth C), Z_1)=(M\underset{\w_C}
\to{\otimes}Z_0)\underset{Z_0}\to{\otimes} (\Hom_{\w_C}(\goth C,
Z_1))\eqno{(5.2.3)}$$ (an equality of  linear algebra). In order to show that we can
consider $\tau$-invariant subspaces, we need the following objects.
Let $I$ be an ideal of $\w$, $\Cal M_0=IZ_0$. It is clear that $\Cal M_0^\tau=I$.
Further, let $\Cal M_1$ be a locally free $Z_0$-module. We have a formula: $$(\Cal
M_0\underset{Z_0}\to{\otimes} \Cal M_1)^\tau=\Cal M_0^\tau \underset{\w}\to{\otimes} \Cal
M_1^\tau \eqno{(5.2.4)}$$
Really, $\Cal M_0\underset{Z_0}\to{\otimes} \Cal M_1=I\Cal M_1$, and $$(I\Cal
M_1)^\tau=I\Cal M_1^\tau\eqno{(5.2.5)}$$ where this formula is true by the following
reason. Obviously $(I\Cal M_1)^\tau\supset I\Cal M_1^\tau$. Let $J$ be an ideal of $\w$
such that $IJ$ is a principal ideal. We have $(IJ(J^{-1}\Cal M_1))^\tau=IJ(J^{-1}\Cal
M_1)^\tau$ and $(IJ(J^{-1}\Cal M_1))^\tau\supset I(J(J^{-1}\Cal M_1))^\tau\supset
IJ(J^{-1}\Cal M_1)^\tau$, hence all these objects are equal and we get 5.2.5 and hence
5.2.4.

The action of $\tau$ on both sides of 5.2.3 coincide.  Considering $\tau$-invariant
elements of both sides of 5.2.3 and taking into consideration 5.2.4 ($\Cal
M_0=\Hom_{\w_C}(\goth C,
Z_1)$ and $\Cal M_1=M\underset{\w_C}
\to{\otimes}Z_0$) we get the lemma. $\square$
\medskip
This lemma proves (A) of Step 1.
\medskip
{\bf Lemma 5.2.6.} Let $\Cal M_i$ ($i=0,1$) be two locally free $Z_0$-modules with $\tau$-action satisfying $\tau(cm)=\tau(c) \tau(m)$ ($c\in Z_0$, $m\in \Cal M_i$), and $\psi: \Cal M_0\otimes_{Z_0}\Cal M_1\to Z_1$ a perfect pairing of $Z_0$-modules with $\tau$-action. Let, further, both $\Cal M_i$ satisfy $\Cal M_i^\tau\otimes_{\w} Z_0=\Cal M_i$. Then the restriction of $\psi$ to $\Cal M_0^\tau\otimes_{\w}\Cal M_1^\tau\to \Omega$ is a perfect pairing as well.
\medskip
{\bf Proof.} Let $\alpha: \Cal M_0^\tau \to \Omega$ be an $\w$-linear map. We prolonge it to a map $\bar \alpha: \Cal M_0 \to Z_1$ by $Z_0$-$\tau$-linearity. By perfectness of $\psi$, there exists $m_1\in \Cal M_1$ such that $\bar \alpha(m_0)= \psi(m_0 \otimes m_1)$. It is easy to see that $m_1$ is $\tau$-invariant (we use the fact that $\tau: Z_0 \to Z_0$ is surjective). $\square$
\medskip
{\bf Lemma 5.2.7.} There is a natural perfect $\w$-linear
$\Omega$-valued pairing
between $H_1(M)$ and $H^1(M)$:
$H_1(M)\underset{\w}\to{\otimes}H^1(M)\to \Omega$.
\medskip
{\bf Proof.} For the case $\w=\n F_q[T]$ this is [G], Corollary
5.9.35. General case: we have a perfect $Z_0$-pairing
$$\Hom_{\w_C}(M,Z_1)\underset{Z_0}\to{\otimes} (M\underset{\w_C}
\to{\otimes} Z_0) \to Z_1$$
Now we take $\Cal M_0=\Hom_{\w_C}(M,Z_1)$, $\Cal M_1=M\underset{\w_C} \to{\otimes}Z_0$ and we apply Lemma 5.2.6. $\square$
\medskip
Step 1 of the theorem follows from these lemmas. 
\medskip
{\bf Remark 5.2.8.} The pairing can be defined also as the composition of
$$\matrix H_1(M)\underset{\w} \to{\otimes} H_1(M')=\Hom_{\w_C[\tau]}(M,Z_1)
\underset{\w}
\to{\otimes} \Hom_{\w_C[\tau]}(M',Z_1) \\
\to \Hom_{\w_C[\tau]}(M\underset{\w_C} \to{\otimes}M'
,Z_1\underset{Z_0} \to{\otimes}Z_1) \to
\Hom_{\w_C[\tau]}(\goth C,Z_1\underset{Z_0} \to{\otimes}Z_1) = L_T(\goth C) \underset{\w} \to{\otimes}\Omega\endmatrix \eqno{(5.2.9)}$$
where the second map comes from a canonical map $\delta: \goth C \to M\underset{\w_C}
\to{\otimes}M'$ of Remark 1.11.1 (more exactly, of its analog for arbitrary $\w$).
\medskip
{\bf Remark 5.2.10.} Recall that the explicit formula for functoriality is
the following. Let $\alpha: M_1 \to M_2$ be a map of t-motives,
$\alpha': M'_2 \to M'_1$
the dual map and $L_T(\alpha): L_T(M_2) \to L_T(M_1)$, $L_T(\alpha'): L_T(M'_1)
\to L_T(M'_2)$ the
corresponding maps on lattices. For any $l_1'\in L_T(M_1')$, $l_2\in L_T(M_2)$ we
have:
$$<L_T(\alpha)(l_2), l_1'>_{M_1}=<l_2, L_T(\alpha')(l_1')>_{M_2}\eqno{(5.2.11)}$$

{\bf 5.3. Proof: Step 1, Version 2.} Case $\w=\n F_q[T]$. We identify $Z_1$ of [G], p.168, lines 3 -- 4 with $\p\{T\}$ (see [G], Definition 5.9.10) and $\w$ with $\Omega$. Like above, we have an isomorphism of $\w$-modules (recall that $\w$ is the center of $\p[T,\tau]$):
$$L_T(M)=\Hom_{\p[T,\tau]}(M,Z_1)\eqno{(5.3.1)}$$ ([G], first terms of 5.9.25, 5.9.19). Let $\vf: M \to Z_1$, $\vf':
M' \to Z_1$ be elements of $L_T(M)$, $L_T(M')$ respectively, and let
$f_*$, $f'_*$, $Q$, $Q'$ be from
1.9.3, 1.10. We denote $$\vf(f_*)=v_*\eqno{(5.3.2)}$$ where $v_*\in
(Z_1)^r$ is a vector column (it is a column of the scattering matrix ([A], p. 486) of
$M$, see 5.4.1 below). The same notation for
the dual:
$\vf'(f'_*)=v'_*$. Condition
that $\vf$, $\vf'$ are
$\tau$-homomorphisms is equivalent to $$Qv_*=v_*^{(1)}, \ \
Q'v'_*={v'}_*^{(1)}\eqno{(5.3.3)}$$ (analog of the formula for scattering matrices [A],
(3.2.2)). Let us consider $\Xi=\sum_{i=0}^\infty
a_iT^i\in\p\{T\}\subset\p[[T]]$ of [G], p. 172, line 1; recall that it is the only
element (up to
multiplication by $\n F_q^*$) satisfying
$$\Xi=(T-\theta)\Xi^{(1)}, \ \ \ \lim_{i\to\infty}a_i=0, \ \ \ |a_0|>|a_i| \ \ \forall
i>0\eqno{(5.3.4)}$$ (see [G], p. 171,
(*); there is a formula $\Xi=a_0\prod_{i\ge0}(1-T/\theta^{q^i})$ where $a_0$ satisfies
$a_0^{q-1}=-1/\theta$). Finally, we define
$$<\vf, \vf'>=\Xi v_{*}^tv'_{*}\eqno{(5.3.5)}$$
Obviously $<\vf, \vf'>$ does not depend on a choice of a basis $f_*$.
\medskip
{\bf Lemma 5.3.6.} $<\vf, \vf'>\in \w$.
\medskip
{\bf Proof.} Firstly, this element belongs to $\n F_q[[T]]$, because $$\Xi
v_{*}^tv'_{*} -(\Xi v_{*}^tv'_{*})^{(1)}=\Xi (v_{*}^tv'_{*}-(T-\theta)^{-1}
v_{*}^{(1)t}{v'}_{*}^{(1)})= \Xi v_{*}^t(I_r-(T-\theta)^{-1}
Q^tQ')v'_{*}$$ because of (5.3.3). But we have (see (1.10.1) --- the
definition of $Q'$)
$$I_r-(T-\theta)^{-1} Q^tQ'=0$$
Secondly, let $<\vf, \vf'>=\sum_{i=0}^\infty c_iT^i$. Since
coefficients of all factors of (5.3.5): $\Xi$, $v_*$ and $v'_*$ ---
tend to 0, we get
that $c_i$ also tend to 0. But $c_i\in\n F_q$, i.e. they are almost
all 0. $\square$
\medskip
{\bf Lemma 5.3.7.} The above pairing is perfect.
\medskip
{\bf Proof.} We have an isomorphism (here $M\{T\}=M\otimes_{\p[T]}\p\{T\}$ with
the natural action of $\tau$ (see [G], Definition 5.9.11))
$$\alpha: \Hom_{\p[T,\tau]}(M,Z_1)\to\Hom_{\w}(M\{T\}^\tau,\w)\eqno{(5.3.8)}$$ defined as the composition of the maps
$$\Hom_{\p[T,\tau]}(M,Z_1)=\Hom_{\p[T]}(M,Z_1)^\tau\overset{\beta'}\to{\to}\Hom_{\p\{T\}}(M\{T\},\p\{T\})^\tau$$
$$\overset{\gamma}\to{\to} \Hom_{\w}(M\{T\}^\tau,\w)$$ where
$\beta: \Hom_{\p[T]}(M,Z_1)\to \Hom_{\p\{T\}}(M\{T\},\p\{T\} )$ is the natural map and
$\beta'$
is the restriction of $\beta$ to $\tau$-invariant elements. Using the Anderson's
criterion of uniformizability
of $M$ (see, for example, [G], 5.9.14.3 and 5.9.13) we get immediately that both
$\gamma$, $\beta$, and hence
$\beta'$, and hence $\alpha$ are isomorphisms. Further, let us consider a homomorphism
$$i: \Hom_{\p[T,\tau]}(M',Z_1)\to M\{T\}^\tau\eqno{(5.3.9)}$$ defined as follows. Let
$\vf'$, $f'_*$, $v'_*$ be as above. We set
$$i(\vf')=\Xi{v'}^t_*f_*\in M\underset{\p[T]}\to{\otimes}\p[[T]]$$ Since $\Xi\in \p
\{T\}$, we get that $\Xi{v'}^t_*f_*\in M\{T\}$. A simple calculation
(like in the Lemma 5.3.6, but simpler) shows that $i(\vf')$ is $\tau$-invariant, hence
$i$ really defines a map from $\Hom_{\p[T,\tau]}(M',Z_1)$ to $M\{T\}^\tau$. Obviously it
is an inclusion. Let us prove that $i$ is surjective. Really, let $c_*\in (Z_1)^r$ be a
column
vector such that $c_*^tf_*\in M\{T\}^\tau$. An analog of the above calculation shows that
if we define $\vf'$ by the formula $\vf'(f'_*)=\Xi^{-1}c_*$ then $\vf'\in
\Hom_{\p[T,\tau]}(M',Z_1)$, and $i(\vf')=c_*^tf_*\in M\{T\}^\tau$. Finally, the
combination of isomorphisms (5.3.8) and (5.3.9) corresponds to the
pairing (5.3.5). $\square$
\medskip
{\bf 5.4. Step 2 -- End of the proof of Theorem 5.} It is easy to see that the converse of the Corollary 5.1.1 (taking into consideration
Proposition 3.2) is also true, i.e. in order to prove Theorem 5 it is sufficient
to prove that a Siegel matrix of $M'$ is $-Z^t$ where $Z$ is a Siegel matrix of $M$. Let
us consider a basis $l_1,...,l_r$ of $L_T(M)$ and for each $l_i$ we consider the
corresponding (under identification 5.3.1) $\vf_i\in \Hom_{\p[T,\tau]}(M,Z_1)$. Let
$\Psi$ be the scattering matrix of $M$ ([A], p. 486) with respect to the bases $l_1,...,l_r$,
$f_1,...,f_r$, and we denote $\vf_i(f_*)$ by $v_{i*}$ (notations of 5.3.2).
\medskip
{\bf Lemma 5.4.1.} $v_{i*}$ is the $i$-th column of $\Psi$ ($Z_1$ is identified with $\p\{T\}$, see the proof).
\medskip
{\bf Proof.} Follows from the definitions. Recall that $\hbox{\bf K}=\n F_q((1/T))$. The isomorphism 5.3.1 is the composition of
2 isomorphisms $i_1: L_T(M) \to \Hom_{\w}^c(\hbox{\bf K}/\w, E)$ ([G], 5.9.19) and $i_2:
\Hom_{\w}^c(\hbox{\bf K}/\w, E) \to \Hom_{\p[T,\tau]}(M,\Hom^c(\hbox{\bf K}/\w, \p)$
([G], 5.9.24; recall that $Z_1=\Hom^c(\hbox{\bf K}/\w, \p)$). For $l_i\in L_T(M)$ we have
$(i_1(l_i))(T^{-k})=\exp(\theta^{-k}l_i)$ ([G], line above the lemma 5.9.18) and
$$((i_2\circ i_1(l_i))(f_j))(T^{-k})=<f_j,\exp(\theta^{-k}l_i)>$$ ([G], two lines above the
lemma 5.9.24). Using the identification of $Z_1$ and $\p\{T\}$ ([G], p. 168, lines 3 - 4)
and the definition of $\Psi$ ([A], p. 486, first formula of 3.2) we get immediately the
lemma. $\square$
\medskip
Let $l'_1,...,l'_r$ be a basis of $L_T(M')$ which is dual to a basis $l_1,...,l_r$ of
$L_T(M)$ with respect to the pairing 5.3.5.
\medskip
{\bf Lemma 5.4.2.} The scattering matrix of $M'$ with respect to the bases
$l'_1,...,l'_r$, $f'_1,...,f'_r$ (denoted by $\Psi'$) is $\Xi^{-1}\Psi^{t-1}$.
\medskip
{\bf Proof.} Follows immediately from 5.4.1 applied to both $M$, $M'$, and formula 5.3.5.
$\square$
\medskip
{\bf Remark 5.4.3.} An alternative proof for the case of pure $M$ (for $some$ basis of $L_T(M')$)
is the following. We denote $\Xi^{-1}\Psi^{t-1}$ by $\Psi_1$. It satisfies
$\Psi_1^{(1)}=(T-\theta) Q^{t-1}\Psi_1$ and other conditions of [A], 3.1. According [A],
Theorem 5, p. 488, there exists a pure uniformizable t-motive $M_1$ with
$\sigma$-structure such that its scattering matrix is $\Psi_1$. Since $\Psi_1$ satisfies
$$\Psi_1^{(1)}=Q'\Psi_1$$
we get that $Q(M_1)=Q'$, i.e. $M_1=M'$. $\square$
\medskip
Let us recall the statement of the crucial proposition 3.3.2 of [A]. Here we consider the
case of those $M$ whose $N$ is not necessarily 0. Let $\Psi$ be a scattering matrix of
$M$. We consider the $(T-\theta)$-Laurent series for $\Psi$ (here $k(M)<0$ is a number,
and $D_{-i}\in M_{r}(\p)$): $$\Psi=\sum_{i=k(M)}^\infty D_{-i}(T-\theta)^i$$ We consider its
negative part $$\Psi^-:=\sum_{i=k(M)}^{-1} D_{-i}(T-\theta)^i$$ as an element of $M_{r}(\p)((T-\theta))/M_{r}(\p)[[T-\theta]]$.
\medskip
We consider the space $(T-\theta)^{k(M)}\p[[T-\theta]]/\p[[T-\theta]]$ as a $\p$-vector space endowed by the action of $\w$, and we denote by $\goth V$ its $r$-th direct sum written as vector columns of length $r$. Obviously $$k(M)=-1 \iff \hbox{ the action of $T$ on $\goth V$ coincides with multiplication by $\theta$} \eqno{(5.4.3a)}$$
We denote the $i$-th column of $\Psi^-$ by $\Psi^-_{i*}$, it belongs to $\goth V$. Following [A], we denote by $\Prin(M)$ (resp. by $\Prin_0(M)$)
the $\p[T]$-linear span (resp. the $\w$-linear span) of all $\Psi^-_{i*}$ in $\goth V$. Finally, we obviously extend the definition of $\Lie_T(M)$, $L_T(M)$ to the case $N\ne 0$; formula 5.1.4 becomes
$$\goth a(Tx)=(\theta+N)(\goth a(x))\eqno{(5.4.3b)}$$
\medskip
{\bf Proposition 3.3.2, [A]} (see also Remark 5.5 below). There exists a $\p[T]$-linear isomorphism $\psi_E: \Lie_T(M)
\to \Prin(M)$ such that its restriction to $L_T(M)\subset \Lie_T(M)$ defines an isomorphism
$L_T(M) \to \Prin_0(M)$ (denoted by $\psi_E$ as well). $\square$
\medskip
{\bf Corollary 5.4.4.} $N=0 \iff k(M)=-1$ (because $N=0 \iff $ the action of $T$ on both $\Lie_T(M)$, $\goth V$ coincides with multiplication by $\theta$, by 5.4.3a). $\square$
\medskip
We return to the case $N=0$.
\medskip
Let us consider the $(T-\theta)$-Laurent series for $\Psi'$ and $\Xi^{-1}$:
$$\Psi'=\sum_{i=k(M')}^\infty D'_{-i}(T-\theta)^i,  \ \ \ \Xi^{-1}=\sum_{i=k(\xi)}^\infty
a_i(T-\theta)^i$$
Since for both $M$, $M'$ we have $N=N'=0$, we get $k(M)=k(M')=-1$. An
elementary calculation shows that $k(\xi)$ is also $-1$. Hence, equality
$\Psi'\Psi^t=\Xi^{-1}$ (Lemma 5.4.2) implies that $D'_{1}D^t_{1}=0$.

Further, there exist $n$ columns of $D_{1}$ which are
$\p$-linerly independent (they are $\psi_E$-images of elements of $L_T(M)$ which form a
$\p$-basis of $\Lie_T(M)$) and all other columns of $D_{1}$ are their linear combinations.
Interchanging columns of $D_{1}$ if necessary we can assume that these columns are the
last $n$
columns. We denote by $D_{12}$ (resp. $D_{11}$ ) the $r\times n$ (resp. $r\times
(r-n)$ ) matrix formed by the last $n$ (resp. the first $r-n$) columns of $D_{1}$. There
exists a matrix $S$ such that $D_{11}=D_{12}S^t$. Again according Proposition 3.3.2,
[A], we have:
$$S \hbox{ is a Siegel matrix of $L(M)$ }\eqno{(5.4.5)}$$
(see also Remark 5.5 below).

Analogous objects are defined for $D'_{1}$. We denote by $D'_{12}$ (resp.
$D'_{11}$) the $r\times n$- (resp. $r\times (r-n)$)-matrix formed by the last $n$
(resp. the first $r-n$) columns of $D'_{1}$. Since
$D'_{1}D^t_{1}=D'_{11}D_{11}^t+D'_{12}D_{12}^t$ we get that
$D'_{12}D_{12}^t+D'_{11}SD_{12}^t=0$. Since $D_{12}^t$ is a $n\times
r$-matrix of rank $n$, it is not a zero-divisor from the right, so $$D'_{12}=
-D'_{11}S\eqno{(5.4.6)}$$ Since the rank of $D'_{1}$ is $r-n$ and $D'_{11}$ is a
$r\times (r-n)$ matrix, (5.4.6) implies that columns of $D'_{11}$ are linearly independent,
and by (5.4.6) and Proposition 3.3.2, [A] we get that $-S$ is a Siegel matrix of
$M'$. $\square$
\medskip
{\bf Remark 5.5.} Since the notations of [A] differ from the ones of the present
paper, for the reader's convenience we give here a sketch of the proof for the case $N=0$
of two  crucial facts: Corollary 5.4.4 and 5.4.5 ([A], Theorem 3.3.2).

Let $\alpha: \Lie(M)\to E(M)$ be a linear isomorphism which is the first term of the
series for $\exp: \Lie(M)\to E(M)$, and let $l\in \Lie(M)$, $f \in M$ be arbitrary. We
consider the $(T-\theta)$-Laurent series $\sum_{i=k}^\infty b_i(T-\theta)^i$ of
$\sum_{j=0}^\infty <\exp(\frac{1}{\theta^{j+1}}l),f>T^j$.
\medskip
{\bf Lemma 5.6.} If $N=0$ then $k=-1$, and $b_{-1}=-<\alpha(l),f>$ (this is [A],
3.3.4).
\medskip
{\bf Sketch of the proof.} For $z\in \Lie(M)$ we denote $\exp(z)- \alpha(z)$ by $\ve(z)$,
hence
$\sum_{j=0}^\infty <\exp(\frac{1}{\theta^{j+1}}l),f>T^j=\underline{A} +\underline{E}$, where
$$\underline{A}=\sum_{j=0}^\infty <\alpha(\frac{1}{\theta^{j+1}}l),f>T^j; \ \ \ \ \underline{E}=\sum_{j=0}^\infty <\ve(\frac{1}{\theta^{j+1}}l),f>T^j$$
We consider their $(T-\theta)$-Laurent series:
$$\underline{A}=\sum_{i=k(\underline{A})}^\infty \underline{a}_i(T-\theta)^i; \ \ \ \ \underline{E}=\sum_{i=k(\underline{E})}^\infty \underline{e}_i(T-\theta)^i$$
Since we have $\exp(z)=\sum_{i=0}^{\infty}C_iz^{(i)}$ where $C_0=I_n$ we get that
$\ve(z)=\sum_{i=1}^{\infty}C_iz^{(i)}$. This means that for large $j$ the element
$\ve(\frac{1}{\theta^{j+1}}l)$ is small, and hence $k(\underline{E})=0$, because finitely many
terms having small $j$ do not contribute to the pole of the $(T-\theta)$-Laurent series
of $\underline{E}$ (the reader can prove easily the exact estimations himself, or to look [A],
p. 491). Since $\alpha$ is $\p$-linear, equality $\sum_{j=0}^\infty
\frac{1}{\theta^{j+1}}T^j= -(T-\theta)^{-1}$ implies that $k(\underline{A})=-1$ and $\underline{a}_{\ -1}=-<\alpha(l),f>$ (and other $\underline{a}_i=0$), hence the lemma. $\square$
\medskip
This lemma obviously implies Corollary 5.4.4. Further, elements $f_1,...,f_r$ generate
the $\p$-space $M/\tau M$, because multiplication by $T$ on $M/\tau M$ coincides with
multiplication by $\theta$, hence the fact that $f_1,...,f_r$ \ \ $\p[T]$-generate $M/\tau
M$ implies that they $\p$-generate $M/\tau M$.

Let $l_1,...,l_n$ form a $\p$-basis of $\Lie(M)$ (here we identify $\Lie_T(M)$ and $\Lie(M)$ via $\goth a$). Since the pairing $<*,*>$ between
$E(M)$ and $M/\tau M$ is non-degenerate and $\alpha$ is an isomorphism, we get that
columns $<\alpha(l_1),f_*>,..., <\alpha(l_n),f_*>$ are linearly independent. Again since
$\alpha$ is an isomorphism and the pairing with $f_*$ is linear, we get that
$$(<\alpha(l_{n+1}),f_*> \ ... \ <\alpha(l_r),f_*>)=(<\alpha(l_{1}),f_*> \ ... \
<\alpha(l_n),f_*>)Z^t $$ Applying the lemma 5.6 to this formula we get immediately 5.4.5.
\medskip
{\bf 6. Tensor products.}
\medskip
There exists an analog of the Theorem 5 for the case of tensor
products of
t-motives. It describes the lattice $L(M_1\otimes M_2)$ in terms of $L(M_1)$,
$L(M_2)$. This is a theorem of Anderson; it is formulated in [P], end of page 3, but its
proof was not published. We recall its statement for the case of arbitrary $N\ne 0$, and
we give its proof for the case $N=0$ (case of arbitrary $N$ can be obtained easily using
the same ideas).

Let $M$ be an uniformizable t-motive whose $N$ is not necessarily 0. Since $N$ is nilpotent, formula 5.4.3b shows that $\Lie_T(M)$ is a $\p[[T-\theta]]$-module. There exists an epimorphism of
$\p[[T-\theta]]$-modules
$$L_T(M)\underset{\w}\to{\otimes}\p[[T-\theta]]\to \Lie_T(M)$$
whose kernel $\goth q=\goth q(M)$ carries information on the pair $(L(M), \Lie(M))$.
\medskip
{\bf Theorem 6} (Anderson). Let $M$, $\bar M$ be any two uniformizable abelian
t-motives. Then $$\goth q(M\otimes \bar M)=\goth q(M)\underset{\p[[T-\theta]]}
\to{\otimes}\goth q(\bar M)\eqno{(6.1)}$$
\medskip
{\bf Remark 6A.} $M\otimes \bar M$ is a uniformizable t-motive ([G], Corollary 5.9.38).
\medskip
{\bf Proof of Theorem 6 (case $N=0$).} We define notations for $M$, and all notations for $\bar M$
will carry bar. Let $e_i$ and $Z$ be from the beginning of Section 3. We denote $\goth a^{-1}(e_i)\in \Lie_T(M)$ by $e_i$ (there is no possibility of confusion). So, $\{e_i\}$ is a
$\p[[T-\theta]]$-basis of $L_T(M)\underset{\w}\to{\otimes}\p[[T-\theta]]$. Elements
$b_i:=(T-\theta)e_i$, $i=1,...,n$ and $b_{n+i}:=e_{n+i}-\sum_{j=1}^n z_{ij}e_j$,
$i=1,...,r-n$ form a $\p[[T-\theta]]$-basis of $\goth q$. We need a
\medskip
{\bf Lemma 6.2.} $\Psi(M\otimes \bar M)=\Psi(M)\otimes\Psi(\bar M)$ where $\Psi(M)$
(resp. $\Psi(\bar M)$; $\Psi(M\otimes \bar M)$) is taken with respect to bases $e_*$ of
$L_T(M)$, $f_*$ of $M_{\p[T]}$ (resp. $\bar e_*$ of $L_T(\bar M)$, $\bar f_*$ of $\bar
M_{\p[T]}$; $e_*\otimes \bar e_*$ of $L_T(M\otimes \bar M)$, $f_*\otimes \bar f_*$ of
$(M\otimes \bar M)_{\p[T]}$) (see the proof for the notations).
\medskip
{\bf Proof.} We consider a map $$\alpha:\Hom_{\p[T]}(M,Z_1)^\tau
\underset{\w}\to{\otimes} \Hom_{\p[T]}(\bar M,Z_1)^\tau\to \Hom_{\p[T]}(M\otimes \bar
M,Z_1)^\tau$$ defined as follows: for $\vf\in \Hom_{\p[T]}(M,Z_1)^\tau$, $\bar \vf\in
\Hom_{\p[T]}(\bar M,Z_1)^\tau$ we let $[\alpha(\vf \otimes \bar \vf)](f\otimes \bar
f)=\vf(f)\cdot\bar \vf(\bar f)$ (it is obvious that $\alpha(\vf \otimes \bar \vf)$ is
$\tau$-stable). Since $e_1, ..., e_r$ (resp. $\bar e_1, ..., \bar e_{\bar r}$) is a basis
of $\Hom_{\p[T]}(M,Z_1)^\tau$ (resp. $\Hom_{\p[T]}(\bar M,Z_1)^\tau$; we identify $L_T(M)$,
resp. $L_T(\bar M)$ with $\Hom_{\p[T]}(M,Z_1)^\tau$ (resp. $\Hom_{\p[T]}(\bar M,Z_1)^\tau$)
we get (using Lemma 5.4.1) that $\Psi(M)$, $\Psi(\bar M)$ are non-degenerate. Since their
product is also non-degenerate, we get $\alpha(e_i\otimes \bar e_{\bar i})$ are linearly
independent and hence a basis of $\Hom_{\p[T]}(M\otimes \bar M,Z_1)^\tau$. Applying once
again Lemma 5.4.1 we get the lemma. $\square$
\medskip
If $A$, $B$ are two matrices then columns of $A\otimes B$ are indexed by pairs $(k,l)$
where $k$ (resp. $l$) is the number of a column of $A$ (resp. $B$). We denote by $A_k$,
$B_l$, $A\otimes B_{(k,l)}$ the respective columns. Obviosly we have: $A\otimes
B_{(k,l)}=A_k\otimes B_l$ (tensor product of column matrices).
\medskip
Let us prove that for $i=1,...,r-n$, $\bar i=1,...,\bar r-\bar n$ the element
$b_{n+i}\otimes \bar b_{\bar n+\bar i}\in \goth q(M\otimes \bar M)$. According [A],
Proposition 3.3.2, it is sufficient to prove that the corresponding linear combination
(see 6.3 below) of the columns of the matrix $\Psi^-_{M\otimes \bar M}$ is 0. Since
$$b_{n+i}\otimes \bar b_{\bar n+\bar i}=\sum_{j,\bar j}z_{ij}\bar z_{\bar i\bar
j}e_j\otimes \bar e_{\bar j}- \sum_{j}z_{ij}e_j\otimes \bar e_{\bar n+\bar i}
-\sum_{\bar j}\bar z_{\bar i\bar j}e_{n+i}\otimes \bar e_{\bar j}+ e_{n+i}\otimes \bar
e_{\bar n+\bar i}$$ we get the explicit form of this linear combination: it is sufficient
to prove that for all $i$, $\bar i$ we have
$$\sum_{j,\bar j}z_{ij}\bar z_{\bar i\bar j} (\Psi^-_{M\otimes \bar M})_{(j,\bar j)} -
\sum_{j}z_{ij}(\Psi^-_{M\otimes \bar M})_{(j, \bar n+\bar i)} $$ $$ -\sum_{\bar j}\bar
z_{\bar i\bar j}(\Psi^-_{M\otimes \bar M})_{(n+i,\bar j)} +(\Psi^-_{M\otimes \bar
M})_{(n+i,\bar n+\bar i)} =0\eqno{(6.3)}$$
Further, 6.2 implies that $$(\Psi^-_{M\otimes \bar M})_{(k,\bar
k)}=\frac{A_{-1,k}\otimes \bar A_{-1,\bar k}}{(T- \theta)^2}+\frac{A_{-1,k}\otimes \bar
A_{0,\bar k}+A_{0,k}\otimes \bar A_{-1,\bar k}}{T- \theta}$$ hence 6.3 becomes
$$\sum_{j,\bar j}z_{ij}\bar z_{\bar i\bar j} (\frac{A_{-1,j}\otimes \bar A_{-1,\bar
j}}{(T- \theta)^2}+\frac{A_{-1,j}\otimes \bar A_{0,\bar j}+A_{0,j}\otimes \bar A_{-1,\bar
j}}{T- \theta}) $$ $$- \sum_{j}z_{ij}(\frac{A_{-1,j}\otimes \bar A_{-1,\bar n+\bar i}}{(T-
\theta)^2}+\frac{A_{-1,j}\otimes \bar A_{0,\bar n+\bar i}+A_{0,j}\otimes \bar A_{-1,\bar
n+\bar i}}{T- \theta}) $$ $$-\sum_{\bar j}\bar z_{\bar i\bar j}(\frac{A_{-1,n+i}\otimes
\bar A_{-1,\bar j}}{(T- \theta)^2}+\frac{A_{-1,n+i}\otimes \bar A_{0,\bar
j}+A_{0,n+i}\otimes \bar A_{-1,\bar j}}{T- \theta}) $$ $$+\frac{A_{-1,n+i}\otimes \bar
A_{-1,\bar n+\bar i}}{(T- \theta)^2}+\frac{A_{-1,n+i}\otimes \bar A_{0,\bar n+\bar
i}+A_{0,n+i}\otimes \bar A_{-1,\bar n+\bar i}}{T- \theta}=0\eqno{(6.4)}$$
It is easy to see that 6.4 follows immediately from the equalities
$$A_{-1,n+i}=\sum_j z_{ij}A_{-1,j}\eqno{(6.5)}$$ $$ \bar A_{-1,\bar n+\bar
i}=\sum_{\bar j} \bar z_{\bar i\bar j}\bar A_{-1,\bar j}$$ For example, the left hand
side of (6.4) has 2 terms containing $\bar A_{0,\bar j}$ (in the middle of the first
and the third lines of (6.4)). Multiplying (6.5) by $\bar z_{\bar i\bar j}\bar
A_{0,\bar j}$ we get that the sum of these 2 terms of (6.4) is 0. For other pairs of
terms of (6.4) the situation is the same.
\medskip
The proof that for $i=1,...,r-n$, $\bar i=1,...,\bar n$ the element $b_{n+i}\otimes \bar
b_{\bar i}\in \goth q(M\otimes \bar M)$ is analogous but simpler. We have
$$b_{n+i}\otimes \bar b_{\bar i}=(T-\theta) (-\sum_{j} z_{ij} e_j\otimes \bar e_{\bar i}
+ e_{n+i}\otimes \bar e_{\bar i})$$ The analog of (6.3)) is $$(T-\theta) (-
\sum_{j}z_{ij}(\Psi^-_{M\otimes \bar M})_{(j, \bar i)}  +(\Psi^-_{M\otimes \bar
M})_{(n+i,\bar i)}) =0$$ and the analog of (6.4)) is $$-
\sum_{j}z_{ij}\frac{A_{-1,j}\otimes \bar A_{-1,\bar i}}{T- \theta}
+\frac{A_{-1,n+i}\otimes \bar A_{-1,\bar i}}{T- \theta}=0$$ This equality follows
immediately from (6.5).
\medskip
Finally, elements $b_{i}\otimes \bar b_{\bar i}$ ($i=1,...,n$, $\bar i=1,...,\bar n$)
obviously belong to $ \goth q(M\otimes \bar M)$.
\medskip
So, we proved that $\goth q(M)\underset{\p[[T-\theta]]}\to{\otimes}\goth q(\bar M)
\subset \goth q(M\otimes \bar M)$. Since the $\p$-codimension of both subspaces in
$L_T(M)\underset{\w}\to{\otimes} L_T(\bar M)\underset{\w}\to{\otimes}
\p[[T-\theta]]$ is $n\bar n$, they are equal. $\square$
\medskip
{\bf 7. Self-dual t-motives.}
\medskip
{\bf Case $\w=\n F_q[T]$.} A uniformizable t-motive $M$ is called self-dual if there exists an isogeny $\alpha: M \to M'$. It defines an $\w$-valued, $\w$-bilinear form $<*,*>_\alpha$ on $L_T(M')$ as follows:
$$<\vf_1, \vf_2>_\alpha=<L_T(\alpha)(\vf_1), \vf_2>_M$$
5.2.11 implies that if $\alpha'=-\alpha$ (resp. $\alpha'=\alpha$) then
$<*,*>_\alpha$ is skew symmetric (resp. symmetric). $M$ is called positively (resp. negatively) self-dual if $\alpha$ satisfies $\alpha'=\alpha$ (resp. $\alpha'=-\alpha$). Hence, we have an
\medskip
{\bf Analogy 7a.} The number field case analog of a pair: $\{$negatively self-dual t-motive of rank $2n$, dimension $n$; negative $\alpha: M \to M'\}$ is a (generic) abelian variety of dimension $n$ with a fixed polarization form.
\medskip
For example, like in the number field case, we can define the Rosati involution $I_\alpha$ on  $\End_0(M):=\End(M)\otimes \n F_q(T)$ by the same formula $I_\alpha(\vf)=\alpha^{-1}\circ \vf'\circ \alpha$.
\medskip
Further, we have a
\medskip
{\bf Conjecture 7b.} The dimension of the moduly variety of negatively self-dual t-motives (if it exists) is $n(n+1)/2$.
\medskip
{\bf Examples.} Let $e_*$ be from 1.9, and let $M=M(A)$ given by the equation (here $A\in M_n(\p)$ is $\goth A_1$ of 1.9.1)
$$Te_*= \theta e_*+ A\tau e_* + \tau^2 e_*\eqno{(7.1)}$$
be a t-motive of dimension $n$ and rank $2n$. Elements $f_i=e_i$,
$f_{n+i}=\tau e_i$ $(i=1,...,n)$ form a $\p[T]$-basis of $M$. We have
(see, for example, Section 11): $M'$ is given by the equation $$Te'_*=
\theta e'_*- A^t\tau e'_* + \tau^2 e'_*$$ and if we define $$f'_i=\tau
e'_i, \ \ f'_{n+i}= e'_i\eqno{(7.2)}$$ then bases $f_*$, $f'_*$ are dual in the
meaning of Lemma 1.10.

Let $\alpha: M \to M'$ be given by the formula $\alpha(e_*)=De'_*$
where $D\in M_n(\p)$ (we impose this essential restriction only in order to simplify exposition. In the general case $D\in M_n(\p[\tau])$, $D_f\in M_{2n}(\p[T])$, $D_f$ from 7.4). Condition that $\alpha$ is a $\p[T,\tau]$-map is
equivalent to
$$D^{(2)}=D, \ \ AD^{(1)}=-DA^t\eqno{(7.3)}$$
Further, we have $$\alpha(f_*)=D_ff'_*\eqno{(7.4)}$$ where
$D_f=\left(\matrix 0 & D\\ D^{(1)}&0 \endmatrix \right)$, hence
$$\alpha'=\pm \alpha\iff
D_f^t=\pm D_f\iff D^{(1)}=\pm D^t\eqno{(7.5)}$$
Let us fix $\varepsilon_0\in \n F_{q^2}$ satisfying $\varepsilon_0^{q-1}=-1$. Then $D=\varepsilon_0I_n$ satisfies 7.5 with the sign minus, and the set of $A$ satisfying 7.3 with this $D$ is the set of symmetric matrices. This justifies 7b, because the set of $A_1\in  M_n(\p)$ such that $M(A)=M(A_1)$ is conjecturally discrete.

For $D=I_n$ the sign in 7.5 is plus and hence a skew symmetric $A$ defines a positively self-dual $M(A)$.
\medskip
{\bf Remark 7.6.} The below statements are conjectures based on arguments similar to the ones which justify the below Conjecture 11.8.3. Since they are of secondary importance, we do not give any details of justification here.
\medskip
{\bf 7.6.1. Conjecture.} If $n\ge 3$ then for a
generic skew symmetric $A$ we have: $\End(M(A))=\w$.
\medskip
{\bf 7.6.2. Corollary.} Conjecture 7.6.1 implies that the "minimal" $\alpha: M \to M'$ is defined uniquely up to an element of
$\n F_q^*$, and hence the symmetric pairing $<*,*>_\alpha$ is also defined uniquely up to an element of $\n F_q^*$.
\medskip
{\bf 7.6.3. Conjecture.} If $n=2$, $\alpha'=\alpha$ then $\End (M)$ is strictly larger than $\w$.
\medskip
Other examples of a self-dual t-motive are $M\oplus M'$ where $M$ is any t-motive, but they do not give interesting examples of pairings.
\medskip
{\bf 7.6.4. Conjecture.} There exist other (distinct from the ones defined by 7.1) self-dual t-motives $M$ having $\End (M)=\w$ (we can use a version of standard t-motives of Section 11).
\medskip
{\bf Example 7.7.} Case $A=0$, $D=I_n$.
\medskip
In this case we can find explicitly the matrix of the symmetric form
$<*,*>_\alpha$ in some basis of $L_T(M')$. Let $\goth C_2$ be the
Carlitz module over the field $\n F_{q^2}$ considered as a rank 2
Drinfeld module over $\n F_{q}$ given by the equation $$Te=\theta e
+ \tau^2e$$ We have $M=\goth C_2^{\oplus n}$. Let $\goth T_T(\goth
C_2)$ be the convergent $T$-Tate module of $\goth C_2$, i.e. the set
of elements $\{z_i\}\in E(\goth C_2)=\p$ $(i\ge -1, \ \ z_{-1}=0)$ such that
$$\hbox{ $T z_i=z_{i-1}$ for $i \ge 0$ (i.e. $z_i^{q^2}+\theta
z_i=z_{i-1}$) and $z_i\to 0$}$$ It is a free 1-dimensional module over
$\n F_{q^2}[T]$. We choose and fix its generator; its $\{z_i\}$ satisfy (like in 5.3.4)
$|z_0|>|z_i| \ \ \forall i>0$. We denote $\sum_{k=0}^\infty z_kT^k$ by $\goth Z$.

Let $c$ be a fixed element of $\n F_{q^2}-\n F_{q}$. Formulas (5.3.3)
show that the following elements $\vf_i$, $\vf'_i$ ($i=1, ... , 2n$)
form bases of $L(M)$, $L(M')$ respectively ($j=1, ..., n$; clearly that thanks to 7.2
we have $\vf'_i(f'_j)=\vf_i(f_{n+j})$, $n+j$ mod $2n$):
$$i\le n: \ \ \vf_i(f_j)=\goth Z\delta_i^j, \ \
\vf_i(f_{n+j})=\goth Z^{(1)}\delta_i^j $$

$$i>n: \ \ \vf_i(f_j)=c\goth Z\delta_{i-n}^j, \ \
\vf_i(f_{n+j})=c^q\goth Z^{(1)}\delta_{i-n}^j $$

$$i\le n: \ \ \vf'_i(f'_j)=\goth Z^{(1)}\delta_i^j, \ \
\vf'_i(f'_{n+j})=\goth Z\delta_i^j $$

$$i>n: \ \ \vf'_i(f'_j)=c^q\goth Z^{(1)}\delta_{i-n}^j, \
\ \vf'_i(f'_{n+j})=c\goth Z\delta_{i-n}^j $$
(by the way, it is clear that the same relation between elements of $\goth T_T(M)$ and
$\Hom_{\p[T,\tau]}(M,Z_1)$ holds for all $M$). Formula 7.4 shows
that $\alpha'(\vf'_i)=\vf_{i+n}$, where $i+n \mod 2n$. Let us denote
$\Xi\cdot \goth Z\cdot \goth Z^{(1)}\in\n F_q^*$
by $\gamma$. The above definitions and formulas show that the matrix
of $<*,*>_\alpha$ in the basis $\vf_{1}, \vf_{n+1}, ... ,\vf_{n},
\vf_{2n}$ consists of $n$ \ \ $(2\times 2)$-blocks (trace and norm of $\n
F_{q^2}/\n F_q$)
$$\gamma\left(\matrix \tr(1) & \tr(c)\\ \tr(c)&\tr(N(c))  \endmatrix
\right)=\gamma\left(\matrix 2 & c+c^q \\ c+c^q & 2c^{q+1} \endmatrix
\right)$$ The determinant of
this block is $-(c-c^q)^2\gamma^2$; it belongs to $\n F_q^{*2}\iff
q\equiv 3 \mod 4$ or $q$ is even. Since we have $n$ blocks, we have:

$$\hbox{det $<*,*>_\alpha\not\in\n F_q^{*2}\iff q\equiv 1 \mod 4$ and
$n$ is odd.}$$
\medskip
{\bf Remark 7.8 (Jorge Morales).}
There is a theorem of Harder (see
e.g. W. Scharlau, "Quadratic and Hermitian forms", Springer-Verlag,
Berlin, 1985, Chapter 6, Theorem 3.3) that states that a unimodular form over
$k[X]$ \ \ --- \ \ $k$ being any field of characteristic not 2 --- is the
extension of a form over $k$, i.e. there is a basis in which all the
entries of the associated symmetric matrix are constant. This means
that the classification of the above quadratic forms over $\n F_q[T]$
($q$ odd) is very simple.
\medskip
{\bf Remark 7.9.} Let $M$ be a t-motive which is both negatively and positively self-dual. There is a natural idea 7.9.2 to define an analog of
Hodge structure on $M$. Nevertheless, this idea fails. Namely, the exact sequence
$$0\to \Ker \vf \to L(M)\otimes \p \overset{\vf}\to{\to} \Lie(M)\to 0$$
is the functional field analog of an exact sequence for an abelian variety $A$:
$$ 0 \to H^{0,-1}(A) \to H^{-1}(A) \to (H^{1,0})^*(A) \to 0 $$
Hence, we can define $H^{0,-1}(M):= \Ker \vf$, and the problem is to define an analog of $H^{-1,0}(M)$.

Let us fix a negative isogeny $\alpha: M\to M'$, and let us extend the skew form $<*,*>_\alpha$ to $L(M)\otimes \p$ by $\p$-linearity. It is easy to check that $\Ker \vf$ is isotropic with respect to this form (there is an analogy with the number field case). Let us consider the following elementary lemma of linear algebra:
\medskip
{\bf Lemma 7.9.1.} Let $W$ be a vector space of dimension $2n$ over a field of characteristic $\ne 2$, $B^+$ (resp. $B^-$) a symmetric (resp. skew symmetric) non-degenerate bilinear form on $W$, and $W_0\subset W$ a subspace of dimension $n$ which is isotropic with respect to both $B^+$, $B^-$. Then almost always there exists the only $W_1\subset W$ of dimension $n$ having properties:
$$W_0\cap W_1=0; \ \ \ W_1 \hbox{  is isotropic with respect to both } B^+, \  B^-$$
where almost always means that entries of the matrices of $B^+$, $B^-$ in a basis of $W$ must not satisfy (at least one of) polynomial relations. $\square$
\medskip
If $\End_0(M)\ne \n F_q(T)$ and the action of $I_\alpha$ on $\End_0(M)$ is not identical, then there exists a positive isogeny $\beta: M\to M'$ and hence the symmetric form $<*,*>_\beta$ on $L(M)\otimes \p$. $\Ker \vf$ is isotropic with respect to  $<*,*>_\beta$. Let us fix $\beta$.
\medskip
{\bf Idea 7.9.2.} To apply Lemma 7.9.1 to this situation ($W=L(M)\otimes \p$, $W_0=\Ker \vf$, $B^+=<*,*>_\beta$, $B^-=<*,*>_\alpha$) in order to get a canonical subspace of $L(M)\otimes \p$ which is complementary to $\Ker \vf$ and hence can be considered as an analog of $H^{-1,0}(M)$.
\medskip
Clearly there is no complete analogy with the number field case. But the situation is even worse:
\medskip
{\bf Proposition 7.9.3.} For all $M$, $\alpha$, $\beta$ the "almost always" condition of Lemma 7.9.1 is not satisfied. $\square$
\medskip
{\bf 8. Relations between lattices and t-motives.}
\medskip
We have\footnotemark \footnotetext{I am grateful to Urs
Hartl
who indicated me this reference.}
\medskip
{\bf Theorem 8.1.} ([H],  Theorem 3.2). The dimension of the moduli set of pure t-motives of
dimension $n$ and rank $r$ is $n(r-n)$. $\square$
\medskip
{\bf Remark.} A tuple
$(e_1,...,e_r)$ of
integers entering in the statement of this theorem in [H] is
(0,...,0,1,...,1) with 0 repeated
$r-n$ times and 1 repeated $n$ times for the case under
consideration.
\medskip
Since this number $n(r-n)$ is equal to the dimension of the set of lattices of rank $r$
and dimension $n$, we can state an
\medskip
{\bf Open question 8.2.} Let $r$, $n$ be given. Let us consider the lattice map from the set of the pure uniformizable
t-motives of rank $r$ and dimension $n$ to the set
of lattices of rank $r$ and dimension $n$. Is it true that its image is open and the fibre at a generic point is discrete? If yes, what is the fibre?
\medskip
{\bf Remark.} Results of [GL17] give some evidence that for the case $r=2n$ in a "neighborhood" of the $n$-th power of the rank 2 Carlitz module the fibre consists of 1 point.
\medskip
Theorem 5 implies that for $n=r-1$ the answer to 8.2 is yes (the below Proposition 11.8.5 shows that most likely the condition of purity is essential):
\medskip
{\bf Corollary 8.4.} All pure t-motives of
dimension $r-1$ and rank $r$ having $N=0$ are uniformizable. There is a 1 -- 1 functorial
correspondence between pure t-motives of dimension $r-1$ and rank $r$ having $N=0$ ($r\ge
2$), and lattices of rank $r$ in $\p^{r-1}$ having dual.
\medskip
{\bf Proof.} Let $L$ be a lattice of rank $r$ in $\p^{r-1}$ having dual $L'$. There
exists the only Drinfeld module $M'$ such that $L(M')=L'$, and let $M$ be its dual.
Theorem 5 implies that $L(M)=L$. If there exists another pure t-motive $M_1$ of
dimension $r-1$ and rank $r$ having $N=0$ such that $L(M_1)=L$ then by Corollary 10.4 (its proof is logically independent: there is no vicious circle) the
dual $M'_1$ is a Drinfeld module, according Theorem 5 it satisfies $L(M'_1)=L'$, hence
$M'_1=M'$ and hence $M_1=M$. $\square$
\medskip
{\bf Remark 8.5.} Recall that lattices of rank $r$ in $\p^{r-1}$ having dual are
described in 3.5 (formulas 3.6, 3.7).
We see that for the case $n=r-1$, $N=0$ purity implies uniformizability. We have
\medskip
{\bf Question 8.5a.} Do exist non-uniformizable t-motives having $n=r-1$, $N=0$?
\medskip
{\bf Question 8.5b.} Do exist uniformizable t-motives having $n=r-1$, $N=0$ such that its lattice has no dual? (Clearly this is a subquestion of 8.2).
\medskip
{\bf Remark 8.6.} Clearly for any $r$, $n$ we have: if a lattice $L$ of rank $r$ and dimension $n$ has no dual then $L\ne L(M)$ for any pure uniformizable $M$. I do not know whether Theorem 6 (which is an analog of Theorem 5 for another tensor operation) imposes a more strong similar restriction on the property of $L$ to be the $L(M)$ of some pure uniformizable $M$, or not.
\medskip
Further, for any uniformizable t-motive $M$ we have a
\medskip
{\bf Corollary 8.7.} If the dual of $(L(M), \Lie(M))$ does not exist then the
dual of $M$ does not exist. Example: the Carlitz module.
\medskip

{\bf 9. Main theorem in terms of Hodge-Pink structure.}
\medskip
Let us consider a version of a special case of the general definition of  Hodge-Pink structure ([P], 0.2; 9.1).
\medskip
{\bf Definition.} A Hodge-Pink structure of constant weight and complete dimension is a pair $\underline{H}=(H,\goth q_H)$ where $H$ is a free finite dimensional $\w$-module and $\goth q_H$ is a $\p[[T-\theta]]$-lattice in $H\underset{\w}\to{\otimes}\p[[T-\theta]]$ such that the dimension of $\goth q_H$ over $\p[[T-\theta]]$ is equal to the dimension of $H$ over $\w$ (condition of complete dimension).
\medskip
Let $\vf: L \hookrightarrow \p^n$ be a lattice. It defines a Hodge-Pink structure $\underline{H}=\underline{H}(L)$ of constant weight and complete dimension. Firstly, instead of a $\n F_q[\theta]$-module $L$ we consider an isomorphic $\w$-module $H$ formally defined by the property $H\underset{\w}\to{\otimes}\n F_q[\theta]=L$ where the map $\w \to \n F_q[\theta]$ is $\iota$. We denote the isomorphism $H \to L$ by $\iota$ as well; the composition $\vf \circ \iota: H \to \p^n$ is a map of $\w$-modules where $T\in \w$ acts on $ \p^n$ by multipication by $\theta$. Further, $\vf \circ \iota$ extends to a surjection of $\p[[T-\theta]]$-modules $H\underset{\w}\to{\otimes}\p[[T-\theta]] \to \p^n$ denoted by $\vf \circ \iota$ as well.  Finally, $\goth q_H$ is defined as $\Ker \vf \circ \iota$.

If $M$ is a pure uniformizable t-motive then we associate it a Hodge-Pink structure of constant weight and complete dimension $\underline{H}(M)=\underline{H}(L(M))$.

Let $m=m(\underline{H})$ be the minimal number such that $\goth q_H \supset (T-\theta)^m H\underset{\w}\to{\otimes}\p[[T-\theta]]$. For $\mu\ge m$ we define the $\mu$-dual structure ${\underline{H}'}^{\mu}=({H'}^{\mu}, \goth q_{{H'}^{\mu}})$ as follows:

$${H'}^{\mu}=H^*, \ \  \goth q_{{H'}^{\mu}}=\{\chi\in H^*\underset{\w}\to{\otimes}\p[[T-\theta]]  $$ $$\hbox{ such that }  \forall y\in \goth q_H  \hbox{ we have } \chi(y)\in (T-\theta)^\mu\p[[T-\theta]]\} $$ It is obvious that it is really a Hodge-Pink structure of constant weight and complete dimension.
\medskip
If $\underline{H}=\underline{H}(L)$ for a lattice $L$ then $m=1$ and if $L$ has dual then $${\underline{H}'}^1=\underline{H}(L')\eqno{(9.1)}$$ this is easy to prove.
\medskip
{\bf Remark 9.2.} And if $L$ has no dual? Really, $\underline{H}(L)$ exists even if $L$ does not satisfy Definition 2.1 (b). If $L$ is a lattice having no dual this means that $L'$  does not satisfy Definition 2.1 (b). Nevertheless, equality ${\underline{H}'}^1=\underline{H}(L')$ is meaningful and holds. We are not interested in these lattices because they cannot be lattices of uniformizable t-motives having dual.
\medskip
Proof of the duality theorem for $M$ having $n\ne0$ is given in [GL18].
\medskip

{\bf 10. Duals of pures, and other elementary results. }

\nopagebreak
\medskip
We consider in this section the case of arbitrary $N$ (i.e. not necessarily $N=0$), and $\w=\n F_q[T]$. The definition 1.8 extends to the case of pr\'e-t-motives, and
remarks 1.11 hold for this case.
\medskip
{\bf Lemma 10.2. } Let $M$ be a pr\'e-t-motive, $m=m(M)$ from its
(1.3.1), and $\mu\ge m$.
Then $M'$ --- the $\mu$-dual of $M$ --- exists as a pr\'e-t-motive,
and $m(M')\le\mu$. If
$M'$ is a t-motive then $\dim M'=r\mu - \dim M$ ($r$ is the rank of $M$).
\medskip
{\bf Proof.} We must check that $Q'$ has no denominators, and the
condition (1.3.1). The
module $\tau M$ is a $\p[T]$-submodule of $M$ (because $a \tau x =
\tau a^{1/q} x$ for
$x\in M$),
hence there are $\p[T]$-bases $f_*=(f_1, ... f_r)^t$, $g_*=(g_1, ...
g_r)^t$ of $M$,
$\tau M$ respectively such that $g_i=P_if_i$, where $P_1 | P_2 | ... |
P_r$, $P_i \in
\p[T]$. Condition (1.3.1) means that $\forall i$ \ $(T-\theta)^m f_i
\in \tau M$, i.e.
$P_i|(T-\theta)^m$, i.e. $\forall i$ \ $P_i=(T-\theta)^{m_i}$ where $0
\le m_i \le
m_{i+1} \le m$. There exists a matrix $\goth Q=\{\goth q_{ij}\}\in
M_r(\p[T])$ such that
$$\tau f_i = \sum_{j=1}^r \goth q_{ij}g_j= \sum_{j=1}^r \goth q_{ij} P_j
f_j\eqno{(10.2.1)}$$
Although $\tau$ is not a linear operator, it is easy to see that
$\goth Q\in GL_r(\p[T])$
(really, there exists $C=\{c_{ij}\}\in M_r(\p[T])$ such that $g_i=P_i
f_i=\tau (\sum_{j=1}^r
c_{ij}f_j)$, we have $C^{(1)}\goth Q =I_r$).

We denote the matrix $\diag(P_1, P_2, ... ,P_r)$ by $\goth P$, so
(10.2.1) means that
$$Q= \goth Q \goth P \eqno{(10.2.2)}$$

{\bf Remark 10.2.3.} Since
$\goth Q \goth P \in
GL_r(\p(T))$, we get that the action of $\tau$ on $i_2(M)$ is invertible.
\medskip
It is clear that if $M$ is a t-motive then
$$\dim M = \sum_{j=1}^r m_j\eqno{(10.2.4)}$$
(because $\dim M = \dim_{\p}(M/\tau M)$. Further, (10.2.2) implies that
for $Q'=Q({M'})$
we have
$$Q'=\goth Q^{t-1}\diag((T-\theta)^{\mu-m_1}, ...,
(T-\theta)^{\mu-m_r})\eqno{(10.2.5)}$$
This means that elements of $Q'$ have no denominators. The condition
(1.3.1) for $M'$
follows easily from (10.2.5) (because $\goth Q^{t-1}\in GL_r(\p[T])$),
and the dimension
formula (for the case $M'$ is a t-motive) follows immediately from
(10.2.4) applied to
$M'$. $\square$
\medskip
A definition of a pure t-motive can be found in [G] ((5.5.2),
(5.5.6) of [G] + formula (1.3.1) of the present paper).
\medskip
{\bf Theorem 10.3.} Let $M$ be a pure t-motive and $m=m(M)$ from
(1.3.1). Then (if
$rm-n>0$) its $m$-dual $M'$ exists, and it is pure.
\medskip
{\bf Proof.} The definition of pure ([G], (5.5.2)) is valid for
pr\'e-t-motives. We use
its following
matrix form. We denote $T^{-1}$ by $S$ and for any $C$ we let
$$C^{[i]}=C^{(i-1)}\cdot C^{(i-2)}\cdot ... \cdot C^{(1)}\cdot C$$

{\bf Lemma 10.3.1.} Let $Q\in M_r(\p[T])$ be a matrix such that formula
(1.9.3) defines an
t-motive $M$. Then it is pure iff there exists $C\in
GL_r(\p((S)) \ )$ such that
for some $\goth q$, $s>0$
$$S^\goth q C^{(s)}Q^{[s]}C^{-1}\in GL_r(\p[[S]])$$
i.e. iff $S^\goth q C^{(s)}Q^{[s]}C^{-1}$ is $S$-integer and its inicial
coefficient is
invertible.
\medskip
{\bf Proof.} Elementary matrix calculations. We take $C$ as a matrix
of base change of
$f_*$ to a $\p[[S]]$-basis of $W$ of (5.5.2) of [G]. $\square$
\medskip
{\bf Lemma 10.3.2.} Let $\mu = m$. We have: $M'={M'}^\mu$ of Lemma 10.2 is a pure
pr\'e-t-motive.
\medskip
{\bf Proof.} Let $\goth q$, $s$ and $C$ be from Lemma 10.3.1. We have
$${Q'}^{[s]}=((T-\theta)^{[s]})^\mu Q^{[s]\ t-1}$$
(we use (1.2)). We take $C'=C^{t-1}$. We have
$$S^{s\mu-\goth q} {C'}^{(s)}{Q'}^{[s]}{C'}^{-1}=$$
$$=S^{s\mu-\goth q}C^{(s)\ t-1}Q^{[s]\ t-1}((\frac1S-\theta)^{[s]})^\mu C^t$$
$$=((1-S\theta)^{[s]})^\mu S^{-\goth q}C^{(s)\ t-1}Q^{[s]\ t-1}C^t$$
$$=((1-S\theta)^{[s]})^\mu (S^\goth q C^{(s)}Q^{[s]}C^{-1})^{t-1}$$
We have: $\goth q/s=n/r$ ([G], (5.5.6)), hence $(s\mu-\goth q)/s=(r\mu-n)/r$ and
$s\mu-\goth q>0$. Further,
$((1-S\theta)^{[s]})^\mu\in GL_r(\p[[S]])$, and the result follows from
Lemma 10.3.1. $\square$
\medskip
{\bf Remark.} This result holds also for $\mu>m$.
\medskip
The theorem 10.3 follows from Lemma 10.2, the above lemmas and the
proposition that a pure
pr\'e-t-motive satisfying (1.3.1) is a t-motive ([G], (5.5.6),
(5.5.7)). $\square$
\medskip
{\bf Corollary 10.4.} Let $M$ be a t-motive such that $m=1$,
$n=r-1$. Then $M$
has dual $\iff$ $M$ is pure $\iff$ $M$ is dual to a Drinfeld module.
\medskip
{\bf Proof.} Dimension formula shows that $M'$ (if it exists) is a
Drinfeld module, and they are all pure. $\square$
\medskip
{\bf Example 10.5.} Let $M$ be given by (notations of 1.9.1)
$$\goth A_0=\theta I_2, \ \ \goth A_{1}= \left(\matrix a_{111} & 0
\\ a_{121} & 1 \endmatrix \right), \ \ \goth A_2= \left(\matrix 1 & 0
\\0 & 0 \endmatrix \right)$$
This $M$ has $m=1$, $n=2$, $r=3$, and it is easy to see that it has no dual. Really, for
this $M$ we have (notations of 1.9) $f_1=e_1$, $f_2=\tau e_1$, $f_3=e_2$,
$$Q=\left(\matrix 0&1&0\\T-\theta&-a_{111} & 0
\\ 0&-a_{121} & t-\theta \endmatrix \right), \ \ Q'=\left(\matrix a_{111} & t-\theta
& a_{121} \\ 1 &0&0\\0&0&1\endmatrix \right)$$ The last line of $Q'$ means that $\tau
f_3'=f_3'$. This is a contradiction to the property that $M'_{\p[\tau]}$ is free. It is
possible also to show (Proposition 11.3.4) that $M$ is not pure, and to use 10.4 in order to
prove that it has no dual.
\medskip
Later (Section 11) we shall construct examples of non-pure abelian
t-motives which have dual. Considerations of 11.8 predict that there is enough such t-motives.
\medskip
{\bf Theorem 10.6. } For any t-motive $M$ there exists $\mu_0$ such that for all
$\mu\ge\mu_0$ the object
${M'}^{\mu}$ exists as a t-motive. For these $\mu$ we have
$${M'}^{\mu+1}={M'}^{\mu}\otimes \goth
C\eqno{(10.6.1)}$$

{\bf Proof.} (10.6.1) holds at the level of pr\'e-t-motives, because $Q(\goth
C)=(T-\theta)I_1$. According [G], Lemma 5.4.10 it is sufficient to prove that
${M'}^{\mu}$ is finitely generated as a $\p [\tau]$-module. We shall
use notations of
Lemma 10.2. We take $$\mu_0=1+\hbox{ \{the maximum of the degrees of
entries of $\goth
Q(M)$ as polynomials in $T$\} } $$ $$+ \max(m_k)$$
Let $f'_1, ...f'_r$ be the basis of ${M'}^{\mu}$ over $\p [T]$ dual to
$f_1, ...f_r$. It
is sufficient to prove the
\medskip
{\bf Lemma 10.6.2.} Let $i_0=\mu-\min(m_k)$. Then elements $T^if'_j$,
$i<i_0$, $j=1,
....,r$, generate ${M'}^{\mu}$ as a $\p [\tau]$-module.
\medskip
{\bf Proof of the lemma.} By induction, it is sufficient to show that
for all $\alpha\ge
i_0$ the equation
$$\tau x = (T-\theta)^\alpha f'_j\eqno{(10.6.3)}$$
(equality in ${M'}^{\mu}$) has a solution
$$x=\sum_{k=1}^r C_kf'_k$$
where $C_k \in \p [T]$, $\deg(C_k)< \alpha$. According (10.2.5), the
solution to (10.6.3)
is given by
$$(C_1^{(1)}, ... , C_r^{(1)})=(0,... 0,
(T-\theta)^{\alpha-\mu+m_j},0,... 0)\goth Q^t$$
(the non-0 element of the row matrix is at the $j$-th place).
Unequalities satisfied by
$\mu$ and $\alpha$ show that all $C_k^{(1)}$ are polynomials of degree
$<\alpha$. Since
$c \mapsto c^q$ is surjective on $\p$, we get the desired. $\square$
\medskip
{\bf 10.7. Virtual t-motives.} \footnotemark \footnotetext{This
notion was
indicated me by
Taguchi.} We need two elementary lemmas.

\nopagebreak
\medskip
{\bf Lemma 10.7.0.}\footnotemark \footnotetext{Anderson proved (not published) that the tensor product of any t-motives is also a t-motive.} If $M$ is a t-motive then $M\otimes \goth C$ is also a t-motive.
\medskip
{\bf Proof.} Let $f_{j}$ ($j=1,...,r$) be a $\p[T]$-basis of
$M_{\p[T]}$ and $\goth f$
from 1.10.2, so $f_{j}\otimes \goth f$ is a $\p[T]$-basis of $(M\otimes \goth
C)_{\p[T]}$. It is sufficient to prove that $(M\otimes \goth
C)_{\p[\tau]}$ is finitely
generated. Since $M_{\p[\tau]}$ is finitely generated, it is easy to
see that there
exists $a$ such that elements
$$(T-\theta)^i f_j, \ \ i=0, ... , a, \ \ j=1, ... ,r$$
generate $M_{\p[\tau]}$. This means that $\forall j=1, ... ,r$ there
exist $c_{ijkl}\in
\p$ such that
$$(T-\theta)^{a+1} f_j=\sum_{i=0}^a \sum_{k=0}^\gamma\sum_{l=1}^r
c_{ijkl}(T-\theta)^i
\tau^k f_l\eqno{(10.7.0.1)}$$
where $\gamma$ is a number.

Let us multiply (10.7.0.1) by $(T-\theta)^\gamma$. Taking into
consideration the formula
of the action of $\tau$ on $M\otimes \goth C$ we get that the result
gives us the
following formula in $M\otimes \goth C$:
$$(T-\theta)^{a+\gamma+1} f_j\otimes \goth f=\sum_{i=0}^a
\sum_{k=0}^\gamma\sum_{l=1}^r
c_{ijkl} (T-\theta)^{i+\gamma-k} \tau^k \cdot (f_l\otimes \goth f)\eqno{(10.7.0.2)}$$
This proves that for all $j$ the element $(T-\theta)^{a+\gamma+1} f_j\otimes \goth f$ is
a linear combination of
$$(T-\theta)^i f_l\otimes \goth f, \ \ i=0, ... , a+\gamma, \ \ l=1, ...
,r\eqno{(10.7.0.3)}$$
in $(M\otimes \goth C)_{\p[\tau]}$.
Multiplying (10.7.0.2) by consecutive powers of $T-\theta$ we get by
induction that
elements of 10.7.0.3 generate $(M\otimes \goth C)_{\p[\tau]}$. $\square$
\medskip
{\bf Lemma 10.7.1.} If $M_1\otimes \goth C$ is isomorphic to
$M_2\otimes \goth C$ then
$M_1$ is isomorphic to $M_2$.
\medskip
{\bf Proof.} Let $f_{i*}$ ($i=1,2$) be a $\p[T]$-basis of
$(M_i)_{\p[T]}$, $Q_i$ from
1.9.3, $\alpha: M_1\otimes \goth C\to M_2\otimes \goth C$ an
isomorphism and $C\in
GL_r(\p[T])$ the matrix of $\alpha$ in $f_{1*} \otimes \goth f $,
$f_{2*}\otimes \goth f
$. The matrix of the action of $\tau$ on $M_i\otimes \goth C$ in the
base $f_{i*} \otimes
\goth f $ is $(T-\theta)Q_i$, and the condition that $\alpha$ commutes with
multiplication by $\tau$ is
$$(T-\theta)Q_1 C = C^{(1)}(T-\theta)Q_2$$
Dividing this equality by $T-\theta$ we get that the map $\alpha_0$
from $M_1$ to $M_2$
having the same matrix $C$ in the bases $f_{i*}$, commutes with
$\tau$, i.e. defines an
isomorphism from $M_1$ to $M_2$. $\square$
\medskip
Using Lemma 10.7.1 we can state the following
\medskip
{\bf Definition.} A virtual t-motive is an object $M\otimes \goth C^{\otimes
\mu}$ where $M$ is
a t-motive
and $\mu\in \n Z$, with the standard equivalence relation (here
$\mu_1\ge\mu_2$):
$$M_1\otimes \goth C^{\otimes \mu_1}=M_2\otimes \goth C^{\otimes \mu_2}\iff
M_2=M_1\otimes \goth C^{\otimes (\mu_1-\mu_2)}$$
$$\iff \exists \mu \hbox{ such that $\mu+\mu_1\ge 0$, $\mu+\mu_2\ge 0$
and } M_1\otimes
\goth C^{\otimes (\mu+\mu_1)}=M_2\otimes \goth C^{\otimes (\mu+\mu_2)}$$
Lemma 10.7.1 shows that these conditions are really equivalent.
\medskip
{\bf Corollary 10.7.2.} The $\mu$-dual of a virtual t-motive is
well-defined and
always exists as a virtual t-motive. $\square$
\medskip
{\bf Proposition 10.8.} The following formula is valid at the level of
pr\'e-t-motives: for
any $\mu_1$, $\mu_2$, if ${M_1'}^{\mu_1}$, ${M_2'}^{\mu_2}$ exist then
${(M_1\otimes
M_2)'}^{(\mu_1+\mu_2)}$ exists and
$${(M_1\otimes M_2)'}^{(\mu_1+\mu_2)}={M_1'}^{\mu_1} \otimes {M_2'}^{\mu_2}$$
{\bf Proof.} This is a functorial equality; also we can check it by
means of elementary
matrix calculations. $\square$
\medskip
{\bf Proposition 10.9.} Let $P\in \hbox{\bf A}$ be an irreducible
element. The Tate
module
$T_P({M'}^{\mu})$ is equal to $$T_P(\goth C)^{\otimes \mu}\otimes
\widehat{T_P(M)}$$
(equality of Galois modules) where $\widehat{T_P(M)}$ is the dual
Galois module.
\medskip
{\bf Proof.} It is completely analogous to the proof of the corresponding theorem for
tensor products
([G], Proposition 5.7.3, p. 157).
All modules in the below proof will be the Galois modules, and equalities of modules will
be equalities of Galois
modules. Recall that $E=E(M)$. Since $T_P(M)=\invlim_n E_{P^n}$, it is sufficient to prove that for any
$a\in \w$ we have $E({M'}^{\mu})_a=E(\goth C^{\otimes \mu})_a\otimes \hat E_a$, where $\hat E_a$ is the dual of $E_a$ in the meaning of [T], Definition 4.1.
We have the following sequence of equalities of modules:
$${M'}^{\mu}/a{M'}^{\mu}=\Hom_{\p[T]}(M/aM, \goth C^{\otimes \mu}/a\goth C^{\otimes
\mu})\eqno{(10.9.2)}$$
such that the action of $\tau$ on both sides of this equality coincide (to define the
action of $\tau$ on the
right and side of (10.9.2) we need the action of $\tau^{-1}$ on $M/aM$; it is
well-defined, because the determinant
of the action of $\tau$ on $M$ is a power of $T-\theta$, hence its image in
$\p[T]/a\p[T]$ is invertible). 10.9.2
follows immediately from the definition of ${M'}^{\mu}$;
$$({M'}^{\mu}/a{M'}^{\mu})^\tau=\Hom_{\n F_q[T]}((M/aM)^\tau,
(\goth C^{\otimes \mu}/a\goth C^{\otimes \mu})^\tau)\eqno{(10.9.3)}$$
This follows from 10.9.2 and the Lang's theorem
$$\goth M/a\goth M=(\goth M/a\goth M)^\tau\underset{\n F_q[T]/a\n
F_q[T]}\to{\otimes}\p[T]/a\p[T]$$
applied to both $\goth M=M$, $\goth M={M'}^{\mu}$ (we use that both $M$, ${M'}^{\mu}$ are
free $\p[T]$-modules).
Finally, we have a formula
$$E(\goth M)_a=\Hom_{\n F_q}((\goth M/a\goth M)^\tau, \n F_q)$$
([G], p. 152, last line of the proof of Proposition 5.6.3). Applying this formula to 10.9.3 we get the desired. $\square$
\medskip
{\bf 11. An explicit formula.}

\nopagebreak
\medskip
We return to the case $N=0$. Let $e_*$, $\goth A$, $\goth A_i$, $l$, $n$ be from (1.9). We consider in the
present section two simple types of t-motives (called standard-1 and standard-2
t-motives respectively) whose $\goth A_i$ have a row echelon form, and we give an
explicit formula for the dual of some standard-1
t-motives. Analogous formula can be easily obtained for more general types of
t-motives. These results are the
first step of the problem of description of all t-motives
having duals.
\medskip
{\bf 11.1.} For the reader's convenience, we give here the definition of standard-1
t-motives
for the case $n=2$ (here $\lambda_1$ and $\lambda_2 $ satisfying
$\lambda_1=l$,
$l>\lambda_2 \ge 2$ are parameters):

$$\goth A_0=\theta I_2, \hbox{ for } 0 < i < \lambda_2 \ \ \goth A_i \hbox{ is
arbitrary, } $$
$$\goth A_{\lambda_2}= \left(\matrix * & 0 \\ * & 1 \endmatrix \right),
\hbox{ for } \lambda_2
< i < l \ \ \goth A_i= \left(\matrix * & 0 \\ * & 0 \endmatrix \right), \ \
\goth A_l= \left(\matrix
1 & 0
\\0 & 0 \endmatrix \right) $$

{\bf 11.2.} To define standard-2 t-motives of dimension $n$, we need to fix
\medskip
1. A permutation $\vf\in S_n$, i.e. a 1 -- 1 map $\vf: (1, ..., n) \to (1, ..., n)$;
\medskip
2. A function $k: (1, ..., n) \to \n Z^+$ where $\n Z^+$ is the set of integers $\ge 1$.
\medskip
{\bf Definition.} A standard-2 t-motive of the type $(\vf, k)$ is an abelian
t-motive of dimension $n$ given by the formulas ($i=1,...,n$):
$$Te_{\vf(i)}= \theta e_{\vf(i)}
+\sum_{\alpha=1}^n\sum_{j=1}^{k(\alpha)-1}a_{j,\vf(i),\alpha}\
\tau^je_\alpha + \tau^{k(i)}e_i\eqno{(11.2.1)}$$
where $a_{j,\vf(i),\alpha}\in\p$ is the $(\vf(i),\alpha)$-th entry of the matrix $\goth A_j$.
\medskip
{\bf Proposition 11.2.2.} Formula 11.2.1 really defines a t-motive denoted by
$M=M(\vf,k)=M(\vf,k,a_{***})$. Its rank is $\sum_{\alpha=1}^n k(\alpha)$ and elements
$X_{\alpha j}:=\tau^je_\alpha$, $\alpha=1, ..., n$, $j=0, ..., k(\alpha)-1$, form its
$\p[T]$-basis.
$\square$
\medskip
The group $S_n$ acts on the set of types $(\vf, k)$ and on the set of the above $M$;
clearly for any $\psi\in S_n$ we have $\psi(M)$ is isomorphic to $M$. Particularly, we
can consider only $\vf$ of the following form of the product of $i$ cycles ($\alpha_0=0,
\ \alpha_i=n$):
$$\vf=(\alpha_0+1, ..., \alpha_1)(\alpha_1+1,...,\alpha_2)
...(\alpha_{i-1}+1,...,\alpha_i)\eqno{(11.2.3)}$$
(standard notation of the theory of permutations, for $\gamma\ne \alpha_j$ we have
$\vf(\gamma)=\gamma+1$, for $\gamma= \alpha_j$ we have  $\vf(\alpha_j)=\alpha_{j-1}+1$).
\medskip
{\bf Example 11.2.4.} Let $\vf$ be defined by 11.2.3, the quantity of cycles $i$ is equal
to $1$ and all $a_{***}=0$. Then the corresponding $M$ is of complete multiplication by a
CM-field $\n F_{q^r}(T)$ and its CM-type $\Phi$ is $\{\Id, \fr^{k(1)}, \fr^{k(1)+k(2)},
..., \fr^{k(1)+k(2)+...+k(n-1)}\}$ where $\fr$ is the Frobenuis homomorphism $\n
F_{q^r}\to \bar \n F_q$ (see 13.3, first case: formulas 13.3.1, 13.3.2 coinside with 11.2.1
for the given $\vf$ and $a_{***}=0$; $i_j$ of 13.3.0 is $k(1)+k(2)+...+k(j-1)$ of the
present notations).
\medskip
{\bf Definition 11.3.} A standard-1 t-motive is a standard-2 t-motive
whose $\vf$ is the identical permutation $Id$.
\medskip
{\bf 11.3.0.} Let $M=M(Id,k)$ be a standard-1 t-motive. Acting
by $\psi\in S_n$ we can consider only the case of non-increasing $k(j)$. We introduce a
number $\goth m\ge1$ --- the quantity of jumps of $k(j)$, and two sequences
$$0=\gamma_0< \gamma_1<...< \gamma_\goth m=n$$
(sequence of arguments of points of jumps of the function $k$) and
$$0=\lambda_{\goth m+1}< \lambda_{\goth m}<...< \lambda_2 <\lambda_1=l$$
(sequence of values of $k$ on segments $[\gamma_{i-1}+1, ...,\gamma_i]$) by the formulas
$$\matrix k(1)=...=k(\gamma_1)=\lambda_1 \\ \\
k(\gamma_1+1)=...=k(\gamma_2)=\lambda_2\\ ... \\
k(\gamma_{\goth m-1}+1)=...=k(\gamma_\goth m)=\lambda_\goth m\endmatrix \eqno{(11.3.1)}$$
\medskip
{\bf Example 11.3.2.} The t-motive $M$ of 11.1 is a standard-1 having $\goth m=2$,
$\gamma_1=1$, $\gamma_2=2$ and $\lambda_1$, $\lambda_2$ as in 11.1. Its rank $r=\lambda_1+\lambda_2$.
\medskip
{\bf Conjecture 11.3.3.} A standard-2 t-motive of the type $(\vf, k)$ (notations of
11.2.3) is pure iff $\forall j=1, ..., i$ we have:
$$\frac{\alpha_j-\alpha_{j-1}}{\sum_{\gamma=\alpha_{j-1}+1}^{\alpha_{j}}k(\gamma)} = \frac{n}{r}$$
This conjecture is obviously true if all $a_{***}$ are 0.
\medskip
To simplify exposition, we prove here only the following particular case of this
conjecture.
\medskip
{\bf Proposition 11.3.4.} Let $M$ be a standard-1 t-motive having $\goth m>1$,
defined over $\n F_q(\theta)$, having a good reduction at a point of degree 1 of $\n
F_q(\theta)$ (i.e. a point $\theta+c$, $c\in \n F_q$). Then $M$ is not pure.
\medskip
{\bf Proof.} Let $M$ be defined by 11.2.1, we use notations of 11.3.1. We consider the
action of Frobenius on $\tilde M$ --- the reduction of
$M$ at $\theta+c$. According [G], Theorem 5.6.10, it is sufficient to prove that orders
of the roots of the characteristic polynomial of Frobenius over $\w$ are not equal. More
exactly, we consider the valuation infinity on $\w$ (defined by the condition
$\ord(T)=-1$); the order corresponds to a continuation of this valuation to $\End(\tilde
M)$. The
action of Frobenius on $\tilde M$ coincides with multiplication by $\tau$, because the
degree of the reduction point is 1.

A basis $f_*$ of $M_{\p[T]}$ is the set of $X_{\alpha j}:=
\tau^je_{\alpha}$ of 11.2.2. The matrix
$Q(M)$ is defined by the following formulas for the action of
$\tau$ on $X_{\alpha j}$:
$$\tau(X_{\alpha j})=X_{\alpha,j+1} \hbox{ if } j <
k(\alpha)-1\eqno{(11.3.4.1)}$$
$$\tau(X_{\alpha,k(\alpha)-1})=TX_{\alpha,0} -
\sum_{\delta=1}^\goth m
\sum_{d=\lambda_{\delta+1}}^{\lambda_{\delta}-1}
\sum_{c=1}^{\gamma_{\delta}} a_{d\alpha c}X_{cd} \eqno{(11.3.4.2)}$$
This means that if we arrange $X_{\alpha j}$ in lexicographic order ($X_{\alpha_1 j_1}$ precedes to $X_{\alpha_2 j_2}$ if $\alpha_1 <
\alpha_2$) then the matrix
$Q(M)$ has the block form: $$Q(M)=(C_{ij})\ \ \ (i,j=1,...,n)$$ where $C_{ij}$ is a
$k(i)\times k(j)$-matrix of the form
$$C_{ii}=\left(\matrix 0&1&0&...&0\\0&0&1&...&0\\ ...&...&...&...&... \\ 0&0&0&...&1 \\
T-\theta & *&*&...&*\endmatrix \right), \ \ C_{ij}=\left(\matrix 0&0&...&0\\
...&...&...&...\\ 0&0&...&0 \\ 0& *&...&*\endmatrix \right) (i\ne j)$$
where asterisks mean elements $a_{***}$ (in some order). We consider the characteristic polynomial $P(X)\in (\p[T])[X]$ of $Q(M)$. We
have $$C_{ii}-XI_{k(i)}=\left(\matrix -X&1&0&...&0\\0&-X&1&...&0\\ ...&...&...&...&... \\
0&0&0&...&1 \\ t-\theta & *&*&...&*-X\endmatrix \right)$$

A subset of the set of entries of a matrix is called (following N.N.Luzin) a lightning if each row and each column of the matrix contains exactly one element of this subset. The product of elements of a lightning is called the value of this lightning (i.e. the determinant is the alternating sum of the values of all lightnings).
\medskip
{\bf Lemma 11.3.4.3.} If a non-zero lightning of $C_{ii}-XI_{k(i)}$ contains the term
$T-\theta$, then it does not contain any term containing $X$. $\square$
\medskip
Let $J$ be a subset of the set $1,...,n$ and $J'$ its complement.
\medskip
{\bf Corollary 11.3.4.4.} If a non-zero lightning of $Q(M)-XI_{r}$ contains terms
$T-\theta$ of blocks $C_{\al}$, $j\in J$, then its value is a polynomial in $X$ of degree
$\le \sum_{j'\in J'}k(j')$, and there exists exactly one such lightning (called the
principal $J$-lightning) whose value is a polynomial in $X$ of degree $\sum_{j'\in
J'}k(j')$. $\square$

Since the characteristic polynomial of Frobenius of $\tilde M$ is $\tilde P$
(respectively the valuation infinity of $\p[T]$), it is sufficient to prove that the
Newton polygon of $P(X)$ is not reduced to the segment $((0,-n); (r,0))$ defined by its
extreme terms $(T-\theta)^n$ and $X^r$. To do it, it is sufficient to find a point on its
Newton polygon which is below this segment. We consider $J_{min}=$ the set of all
$\gamma_\goth m - \gamma_{\goth m-1}$ diagonal blocks $C_{ii}$ ($i=\gamma_{\goth m-1}+1,
..., \gamma_\goth m$) of $Q(M)$ of minimal size $\lambda_\goth m$. The value of the
principal $J_{min}$-lightning is $(T-\theta)^{\gamma_\goth m - \gamma_{\goth m-1}}$ times
polynomial in $X$ of degree $d:=r-(\gamma_\goth m - \gamma_{\goth m-1})\lambda_\goth m$.
Corollary 11.3.4.4 implies that if the value of any other lightning of $Q(M)-XI_r$ contains a
term whose $X$-degree is equal to $d$, then the $T$-degree of this term is strictly less
than $\gamma_\goth m - \gamma_{\goth m-
 1}$. This means that if we write $P(X)=\sum_{i=0}^r C_iX^i$, $C_i\in \p[T]$, then
$\ord_\infty(C_d)= -(\gamma_\goth m - \gamma_{\goth m-1})$, i.e. the point with
coordinates $[-(\gamma_\goth m - \gamma_{\goth m-1}), d]$ belongs to the Newton diagram
of $P(X)$, i.e. it is above (really, at) the Newton polygon of $P(X)$. This point is
below the segment $((0,-n); (r,0))$. $\square$
\medskip
{\bf Remark 11.3.4.5.} It is easy to see that the Newton polygon of $P(X)$ coincides with
the Newton polygon of the direct sum
of trivial Drinfeld modules of ranks $\lambda_*$, i.e. with the Newton polygon of the
polynomial
$$\prod_{i=1}^{\goth m} (X^{\lambda_i}-T)^{\gamma_i - \gamma_{i-1}}$$
\medskip
{\bf 11.4.} To formulate the below theorem 11.5 we need some notations. Let $M$ be a
standard-1 t-motive defined by formulas 11.2.1, 11.3.1. We impose the condition
$\lambda_\goth m \ge 3$. Theorem 11.5 affirms that it has dual. To find explicitly the
dual of $M$, we need to choose an arbitrary function
$\nu: (i,j) \to \nu(i,j)$ which is a 1 - 1 map from the set of pairs
$(i,j)$ such that

$$1 \le i \le n; \ \ 1 \le j \le k(i)-2 \eqno{(11.4.1)}$$
to the set $[n+1, ..., r-n]$ (recall that $r=\sum_{i=1}^n k(i)= \sum_{i=1}^\goth m
(\gamma_{i}-\gamma_{i-1})\lambda_i$).
\medskip
Let the $(r-n)\times(r-n)$-matrices $B_1$, $B_2$ be defined by the
following formulas
(here and until the end of the proof of 11.5 we have $i ,\alpha = 1, ... , n$; \ \
$b_{\beta\gamma\delta}$
is the $(\gamma\delta)$-th entry of $B_\beta$, all
entries of $B_1$,
$B_2$ that are not in the below list are 0):
\medskip
{\bf 11.4.2.} $b_{1i\alpha}= - a_{k(i)-1,\alpha,i}$;
\medskip
$b_{1,\nu(i,j),\alpha} = - a_{j,\alpha,i}$ for $1 \le j \le k(i)-2$;
\medskip
$b_{1,\nu(i,j+1),\nu(i,j)}=1$ for $1 \le j \le k(i)-3$;
\medskip
$b_{1,i,\nu(i,k(i)-2)} =1$;
\medskip
$b_{2,\nu(i,1),i}=1$.
\medskip
We let $B=\theta I_{r-n}+B_1\tau+B_2\tau^2$ and consider a t-motive $M(B)$ (see 11.5.1
below). Formulas
11.4.2 mean that $M(B)$ is standard-2, its $\vf=\vf_B$ is a product of $n$ cycles
$$i\overset{\vf_B}\to{\to}\nu(i,1)\overset{\vf_B}\to{\to}\nu(i,2)\overset{\vf_B}\to{\to}...
\overset{\vf_B}\to{\to}\nu(i,k(i)-2)\overset{\vf_B}\to{\to}i$$ and its $k=k_B$ is defined
by the
formulas $k_B(\gamma)=2$ for $\gamma \in [1, ...,n]$, $k_B(\gamma)=1$ for $\gamma \in
[n+1, ...,r-n]$.
\medskip
{\bf Theorem 11.5.} Let $M$ be from 11.4 (i.e. a standard-1 t-motive having
$\lambda_\goth m \ge 3$). Then $M'=M(B)$.
\medskip
{\bf Proof.}\footnotemark \footnotetext{This proof is a generalization
of the corresponding proof of
Taguchi; we keep his notations.} Let $e'_*=(e'_1, ... e'_{r-n})^t$ be the vector column
of elements of a basis
of $M(B)$ over $\p[\tau]$ satisfying
$$T e'_*=Be'_*\eqno{(11.5.1)}$$
Let us consider the set of pairs $(j,\goth k)$ such that either $j=1,...,n$,
$\goth k=0,1$ or
$j=n+1,...,r-n$, $\goth k=0$. For each pair $(j,\goth k)$ of this set we let (as
in [T], p. 580)
$Y_{j \goth k} = \tau^{\goth k}e'_j$. Formulas (11.4.2) show that these $Y_{**} $
form a basis of
$M(B)_{\p[T]}$, and the action of $\tau$ on this basis is given by the
following formulas
(here $j=1, ... , k(i)-2$):
$$\tau(Y_{i,0})=Y_{i,1} \eqno{(11.5.2.1)}$$
$$\tau(Y_{i,1})= (T-\theta) Y_{\nu(i,1),0} + \sum_{\gamma=1}^n
a_{1\gamma i}Y_{\gamma,1} \eqno{(11.5.2.2)}$$
$$\tau(Y_{\nu(i,j),0})= (T-\theta) Y_{\nu(i,j+1),0} + \sum_{\gamma=1}^n
a_{j+1,\gamma,i}Y_{\gamma,1} \hbox{ if }j<k(i)-2 \eqno{(11.5.2.3)} $$
$$\tau(Y_{\nu(i,k(i)-2),0})= (T-\theta) Y_{i,0} + \sum_{\gamma=1}^n
a_{k(i)-1,\gamma,i}Y_{\gamma,1} \eqno{(11.5.2.4)}$$
Let $X'_{**}$ be the dual basis to the basis $X_{**}$ of 11.2.2.
\medskip
{\bf 11.5.3.} Let us consider the following correspondence between
$X'_{**}$ and $Y_{**}$:
\medskip
$X'_{ij}$ corresponds to $Y_{\nu(i,j),0}$ for
the pair $(i,j)$ like in (11.4.1),
\medskip
$X'_{i0}$ corresponds to $Y_{i1}$ for $1 \le i \le n$;
\medskip
$X'_{i, k(i)-1}$ corresponds to $Y_{i0}$ for $1 \le i \le n$. 
\medskip
Therefore, in order to prove the Theorem 11.5 we must check that
matrices defined by the dual to (11.3.4.*) and by (11.5.2.*) satisfy (1.10.1) under
identification (11.5.3). This is an elementary exercise. $\square $
\medskip
{\bf Remark 11.6.} Clearly it is possible to generalize the Theorem 11.5 to a larger class
of t-motives --- some subclass of standard-3 t-motives, see Definition 11.8.1.
The below example of the proof of Proposition 11.8.7 shows that probably the condition
$\lambda_\goth m \ge 3$ of
the Theorem 11.5 can be changed by
$\lambda_\goth m \ge 2$: it is necessary
to modify slightly formulas 11.4.2. From another side, a standard-1
t-motive of the Example 2.5 shows that this condition cannot
be changed to $\lambda_\goth m\ge1$.
\medskip
{\bf 11.7. An elementary transformation.} To formulate the proposition
11.7.3, we change slightly notations in 1.9.1, namely, instead of $\goth A =
\sum_{i=0}^l \goth A_i \tau^i$ we consider polynomials $P_k(M)$ of $x_1, ...
,x_n$ ($k=1, ...,n$) defined by the formula
$$P_k(M)= \sum_{i=0}^l\sum_{j=1}^n a_{ikj} x^{q^i}_j \eqno{(11.7.1)} $$
Particularly, if $E$ is the t-module associated to $M$ (see [G], 5.4.5),
$x_*=(x_1, ..., x_n)^t$ an element of $E$ then 11.7.1 is equivalent to
$Tx_*=P_*(x_*)$ where $P_*=(P_1(M),...,P_n(M))^t$ is the vector column.
For a standard-1 t-motive $M$ (we use notations of 11.3.0) having $\goth m\ge 2$ we denote vector columns $\goth
P_1(M)=(P_1(M),...,P_{\gamma_1}(M))^t$, $\goth
P_2(M)=(P_{\gamma_1+1}(M),...,P_{\gamma_2}(M))^t$. We use similar
notations for $M'$.
\medskip
{\bf 11.7.2.} Let $M$ be as above, we consider the case
$\lambda_2=\lambda_1-1$. Let $C$ be a fixed $\gamma_1\times (\gamma_
2-\gamma_1)$-matrix. We define a transformed t-motive $M_1$ by the
formulas

$$\goth P_1(M_1)= \goth P_1(M)+C\goth P_2(M)^q$$

$$P_i(M_1)=P_i(M) \hbox{ for } i>\gamma_1$$
\medskip
{\bf Proposition 11.7.3.} For $M$, $C$, $M_1$ of 11.7.2 the dual $M'_1$
of $M_1$ is described by the following formulas:
$$\goth P_2(M'_1)= \goth P_2(M')-C^t\goth P_1(M')^q$$
$$P_i(M'_1)=P_i(M')\hbox{ for }i\not\in[\gamma_1+1, ... , \gamma_2]$$

{\bf Proof} is similar to the proof of the Theorem 11.5, it is omitted.
$\square$
\medskip
\medskip
{\bf 11.8. Non-pure t-motives.} Most results of this subsection are conditional. We shall show that under some natural conjecture the condition of purity in 8.2 and 8.4 is essential, and that for non-pure t-motives the notion of algebraic duality is richer than the notion of analytic duality.

We generalize slightly the definition 11.2.1 as follows. Let
$\succ$ be a linear ordering on the set $[1,...,n]$, and let $\vf$, $k$ be as in 11.2.
\medskip
{\bf Definition 11.8.1.} A standard-3 t-motive of the type $(\vf, k,\succ)$ is
a t-motive of
dimension
$n$ given by the formulas
$$Te_{\vf(i)}= \theta e_{\vf(i)} +\sum_{j=1}^n\sum_{l=1}^{k(j)-1}a_{l,\vf(i),j}\
\tau^le_j +\sum_{j\succ i} a_{k(j),\vf(i),j}\ \tau^{k(j)}e_j +
\tau^{k(i)}e_i\eqno{(11.8.2)}$$
where $a_{***}\in\p$ are coefficients (the only difference with 11.2.1 is the term
$\sum_{j\succ i} a_{k(j),\vf(i),j}\ \tau^{k(j)}e_j$). We denote it by $M(a_{***})$.

Let $M_1=M(a_{1***})$, $M_2=M(a_{2***})$ be two isomorphic standard-3 t-motives of the same type $(\vf,
k,\succ)$ with
$\p[\tau]$-bases $e_{1*}$, $e_{2*}$ respectively (we use notations of 11.8.2 for both
$M_1$, $M_2$). There exists $C\in M_n(\p[\tau])$ such that the formula defining an isomorphism
between $M_1$ and $M_2$ is the following: $e_{2*}=Ce_{1*}$.
\medskip
{\bf Conjecture 11.8.3.} For a generic set of $a_{1***}$ there exists only a countable set of $a_{2***}$ such that $M_2$ is isomorphic to $M_1$.
\medskip
This conjecture is based on calculations in some explicit cases. Particularly, it is proved if $M_1$, $M_2$ are given by the below formula 11.8.5.1 and entries of $C$ are polynomials in $\tau$ of degree $\le 1$.

We denote by $\Cal M_{u}(r,n)$ the moduli space of uniformizable t-motives of the rank
$r$ and dimension $n$,
by $\Cal L(r,n)$ the moduli space of lattices of the rank $r$ and dimension $n$ and by
$\goth L:\Cal M_{u}(r,n) \to \Cal L(r,n)$ the functor of lattice associated to an uniformizable
t-motive.
\medskip
{\bf Proposition 11.8.5.} Conjecture 11.8.3 implies that the dimension of the fibers of
$\goth L$ is $> 0$ for $r=3$, $n=2$. Particularly, we cannot omit condition of purity in the
statement of 8.2.
\medskip
{\bf Proof.} We consider standard-3 t-motives of the type $n=2$,
$\vf=Id$, $k(1)=2$,
$k(2)=1$, $2\succ 1$. Such $M_1=M_1(a_{111}, a_{112}, a_{121})$ is given by
$$\goth A_0=\theta I_2, \ \ \goth A_{1}= \left(\matrix a_{111} & a_{112}
\\ a_{121} & 1 \endmatrix \right), \ \ \goth A_2= \left(\matrix 1 & 0
\\0 & 0 \endmatrix \right)\eqno{(11.8.5.1)}$$ (notations of Example 10.5).
It has $r=3$, it is not pure, hence it has no dual.
Conjecture 11.8.3 implies
that the dimension of the moduli space of these t-motives is 3 (because there are 3
coefficients
$a_{111}, a_{112}, a_{121}$). Uniformizable t-motives form an open subset of this moduli
space, while
the moduli space of lattices of $n=2$ and $r=3$ has dimension 2. $\square$
\medskip
{\bf Remark.} Similar calculations are valid for any sufficiently large $r$, $n$.
\medskip
Standard-3 t-motives of the above type have not dual. The following proposition
shows that the same
phenomenon holds for t-motives having dual. We denote by $\Cal M_{u,d}(r,n)$ the
moduli space of uniformizable t-motives of the rank $r$ and dimension $n$ having dual, by
$\Cal L_d(r,n)$ the moduli space of lattices of the rank $r$ and dimension $n$ having
dual, by $\goth L_d:\Cal M_{u,d}(r,n) \to \Cal L_d(r,n)$ the functor of lattice and by $D_M: \Cal
M_{u,d}(r,n)\to \Cal M_{u,d}(r,r-n)$, $D_L: \Cal L_d(r,n)\to \Cal L_d(r,r-n)$ the
functors of duality on t-motives and lattices respectively. Practically, Theorem 5
means that the following diagram is commutative:
$$\matrix \Cal M_{u,d}(r,n)&\overset{D_M}\to{\to}&\Cal M_{u,d}(r,r-n)
\\ \\ \goth L_d\downarrow&&\goth L_d\downarrow\\  \\ \Cal L_d(r,n)&\overset{D_L}\to{\to}&\Cal
L_d(r,r-n)\endmatrix\eqno{(11.8.6)}$$

{\bf Proposition 11.8.7.} Conjecture 11.8.3 implies that the dimension of the fibers of
$\goth L_d$ in the diagram (11.8.6) is $> 0$ for $r=5$, $n=2$.
\medskip
Practically, this means that the notion of algebraic duality is "richer" than the notion
of analytic duality.
\medskip
{\bf Proof.} We consider standard-3 t-motives of the type $n=2$,
$\vf=Id$, $k(1)=3$,
$k(2)=2$, $2\succ 1$, $r=5$. Such $M$ is given by
$$\goth A_0=\theta I_2, \ \ \goth A_{1}= \left(\matrix a_{111} & a_{112}
\\ a_{121} & a_{122}  \endmatrix \right), \ \ \goth A_{2}= \left(\matrix a_{211} & a_{212}
\\ a_{221} & 1 \endmatrix \right), \ \ \goth A_3= \left(\matrix 1 & 0
\\0 & 0 \endmatrix \right)$$ (notations of Example 10.5). It has dual. Really, we denote by
$A_{i*j}$ the $j$-th column of $\goth A_i$, and we denote by $ (C_1 | C_2 )$ the matrix formed
by union of columns
$C_1$, $C_2$. Then $M'=M(B)$ is also a standard-3 t-motive, where
\medskip
$B_1= \left(\matrix - \det \goth A_2 & -a_{221} & 1
\\ - \det ( A_{1*2} | A_{2*2} ) & -a_{122} & 0
\\ - \det ( A_{1*1} | A_{2*2} ) & -a_{121} & 0
\endmatrix \right)$,
$B_2 = \left(\matrix 0 & 0 & 0
\\ -a_{212}^q & 1 & 0
\\ 1 & 0 & 0 \endmatrix \right)$

The same arguments as in the proof of Proposition 11.8.5 show that the conjecture 11.8.3
implies that the dimension of the moduli space of these t-motives is 7, while
the moduli space of lattices of $n=2$ and $r=5$ has dimension 6. $\square$
\medskip
As above, similar calculations are valid for any sufficiently large $r$, $n$; clearly the dimension of fibers of $\goth L_d$ becomes larger as $r$, $n$ grow.
\medskip
Let us mention two open questions related to the functor $\goth L$. Firstly, let $L$ be a self-dual lattice such that $L\in \goth L(\Cal M_{u,d}(2n,n))$. This means that $D_M: \goth L_d^{-1}(L) \to \goth L_d^{-1}(L)$ is defined.
\medskip
{\bf Open question 11.8.8.} What can we tell on this functor, for example, what is the dimension of its stable elements?
\medskip
Secondly, let us consider $M_1$, $M_2$ of CM-type with CM-field $\n F_{q^r}(T)$, see 13.3.
\medskip
{\bf Open question 11.8.9.} Let the CM-types $\Phi_1$, $\Phi_2$ of the above $M_1$, $M_2$ satisfy $\Phi_1\ne \alpha \Phi_2$, where $\alpha\in \Gal(\n F_{q^r}(T)/\n F_{q}(T))$. Are lattices $L(M_1)$, $L(M_2)$ non-isomorphic?
\medskip
Clearly the negative answer to this question implies the negative answer to the Question 8.2.
\medskip
For any given $M_1$, $M_2$ the answer can be easily found by computer calculation.
Really, let $M$ be one of $M_1$, $M_2$, $c_1,...,c_r$ a basis of $\n F_{q^r}/\n F_{q}$ and
$\alpha_{1},...,\alpha_{n} \subset \Gal(\n F_{q^r}(\theta)/\n F_{q}(\theta))$ the CM-type of $M$.
We define matrices $\Cal M$, $\Cal N$ as follows: $(\Cal M)_{ij}=\alpha_j(c_i)$ $(i,j=1,...,n$),
$(\Cal N)_{ij}=\alpha_{j}(c_{n+i})$, $j=1,...,n$, $i=1,...,r-n$. The Siegel matrix $Z(M)$ is obviously $\Cal N\Cal
M^{-1}$. So, we can find explicitly $Z(M_1)$, $Z(M_2)$ for both $M_1$, $M_2$. To check
whether $Z(M_1)$, $Z(M_2)$ are equivalent or not, it is sufficient to find a solution to
3.8.1 such that the entries of $A$, $B$, $C$, $D$ are in $M_{*,*}(\n F_q)$ (this is
obvious: the condition $\exists \gamma \in GL_r(\n Z_\infty)$ is equivalent to the
condition $\exists \gamma \in GL_r(\n F_q)$, because entries of $Z(M_1)$, $Z(M_2)$ are in $\n F_{q^r}$). The equation 3.8.1 is linear with respect
to $A$, $B$, $C$, $D$, and we can check whether its solution satisfying $\det
\gamma\ne 0$ exists or not.
\medskip
For the case $q=2$, $r=4$, $n=2$, CM-types of $M_1$, $M_2$ are $(Id, Fr)$, $(Id, Fr^2)$
respectively, a calculation shows that the answer is positive: lattices $L(M_1)$, $L(M_2)$ are not isomorphic.
\medskip
{\bf 12. t-motives having multiplications.}
\medskip
Let $\goth K$ be a separable extension of $\n F_q(T)$ such that $\goth K_C:=\goth K\underset{\n F_q}\to{\otimes} \p$ is also a field, $\pi: X\to P^1(\p)$ the projection of curves over $\p$ corresponding to $\p(T)\subset \goth K_C$. Let $\goth K$, $X$ satisfy the condition: $\infty\in X$ is the only point
on $X$ over $\infty\in P^1(\p)$. Let $\w_{\goth K}$ be the subring of
$\goth K$ consisting of elements regular outside of infinity. We
denote $g=\dim \goth K/\n F_q(T)$ and
$ \alpha_1, ... , \alpha_g: \goth K\to\p$ --- inclusions over $\iota: \n F_q(T)\to\p$
(recall that $\iota(T)=\theta$). Let $\Cal W$ be a central simple algebra
over $\goth K$ of dimension $\goth q^2$. Each $ \alpha_i: \goth K \to
\p$ can be extended to a representation $\chi_i: \Cal W \to M_\goth
q(\p)$.
\medskip
{\bf 12.1. Analytic CM-type.} Let $(L, V)$ be as in Section 2 (recall that
$\w=\n F_q[T]$) such that there exists an inclusion $i: \Cal W \to
\End^0(L, V)$, where $\End^0(L, V)=\End(L, V)\underset{\w}\to{\otimes} \n F_q(T)$.
It defines a representation of $\Cal W$
on $V$ denoted by $\Psi$ which is isomorphic to $\sum_{i=1}^g\goth r_i\chi_i$ where
$\{\goth r_i\}$ are some multiplicities (the CM-type of the action of
$\Cal W$ on $(L, V)$). [Proof: restriction of $\Psi$ on $\goth K$ is a sum of one-dimensional representations, i.e. $V=\oplus_{i=1}^g V_i$ where $k\in \goth K$ acts on $V_i$ by multiplication by $\alpha_i(k)$. Spaces $V_i$ are $\Psi$-invariant. We consider an isomorphism $\Cal W\otimes_\goth K\p=M_\goth q(\p)$ where the inclusion of $\goth K$ in $\p$ is $\alpha_i$. We extend $\Psi|_{V_i}$ to $\Cal W\otimes_\goth K\p$ by $\p$-linearity using the inclusion $\alpha_i$ of $\goth K$ in $\p$. It remains to show that a representation of $M_\goth q(\p)$ is a direct sum of its $\goth q$-dimensional standard representations. We consider the corresponding representation of Lie algebra $\goth s\goth l_\goth q(\p)$. It is a sum of irreducible representations. Let $\omega$ be the highest weight of any of these irreducible representations. $\omega$ is extended uniquely to the set of diagonal matrices of $M_\goth q(\p)$, because $\omega$ is identical on scalars. Since our representation is not only of Lie algebra but of algebra $M_\goth q(\p)$, we get that $\omega$ is a ring homomorphism $\Diag(M_\goth q(\p))\to \p$. There exists the only such $\omega$ corresponding to the $\goth q$-dimensional standard representation].

Further, we
denote $m=\dim_\Cal W L \otimes \n F_q(T)$ ($g$, $\goth q$, $\Psi$, $\goth
r_i$, $m$ are analogs of $g$, $q$, $\Phi$, $r_i$, $m$ of [Sh63] respectively).
Clearly we have
$$n=\goth q\sum_{i=1}^g\goth r_i, \ \ r=mg\goth q^2\eqno{(12.2)}$$
By functoriality, we have the dual inclusion $i': \Cal W^{op} \to
\End^0(L',V')$ where $\Cal W^{op}$ is the opposite algebra.
\medskip
{\bf Remark.} A construction of Hilbert-Blumental modules ([A], 4.3, p. 498) practically is a particular case
of the present
construction: for Hilbert-Blumental modules we have $\goth q=1$, i.e. $\goth K=\Cal W$, and all $\goth r_i=1$. Anderson considers the case when $\infty$ splits completely; this difference with the present case is not essential.
\medskip
{\bf Proposition 12.3.} If the dual pair $(L',V')$ exists then the
CM-type of the dual inclusion is $\{m\goth
q - \goth r_i\}$, $i=1, ... ,g$.
\medskip
{\bf Proof.} We have $L\underset{\n Z_\infty}\to{\otimes}\p$ is
isomorphic to $(\Cal W\underset{\n F_q(\theta)}\to{\otimes}\p)^m$ as a
$\Cal W$-module. Since the natural representation of $\Cal W$ on $\Cal
W\underset{\n F_q(\theta)}\to{\otimes}\p$ is isomorphic to $\goth
q\sum_{i=1}^g\chi_i$ we get that $L\underset{\n Z_\infty}\to{\otimes}\p$ is isomorphic to $m\goth q\sum_{i=1}^g\chi_i$ as a $\Cal W$-module. Consideration of the exact
sequence $0\to {V'}^* \to L\underset{\n Z_\infty}\to{\otimes}\p \to V\to0$ gives us the desired. $\square$
\medskip
{\bf Remark 12.4.1.} This result is an analog of the corresponding
theorem in the number field case. We use notations of [Sh63], Section 2. Let
$A$ be an abelian variety having endomorphism algebra of type IV, and
$(r_\nu, s_\nu)=(r_\nu(A), s_\nu(A))$ are from [Sh63], Section 2, (8).
Then
$$r_\nu(A')=mq-r_\nu(A)=s_\nu(A), \ s_\nu(A')=mq-s_\nu(A)=r_\nu(A)$$
By the way, Shimura writes that the CM-types of $A$ and $A'$ coincide
([Sh98], 6.3, second line below (5), case $A$ of CM-type). We see that
his affirmation is not natural: he considers the complex conjugate
action of the endomorphism ring on $A'$. It is necessary to take into
consideration this difference of notations comparing formulas of 12.3 and 13.2
with the corresponding formulas of Shimura.
\medskip
{\bf Remark 12.4.2.} According [L09], a t-motive $M$ is an analog of an abelian variety $A$ with multiplication by an imaginary quadratic field $K$. The above consideration shows that this analogy holds for $M$ and $A$ having more multiplications. Really, if $A$ has more multiplications then (we use notations of [Sh63], Section 2) $F_0=FK$, and numbers $(r_\nu(A), s_\nu(A))$ satisfy $n(A)=q\sum_{i=1}^g r_\nu(A)$, where $(n(A), \dim(A)-n(A))$ is the signature of $A$ treated as an abelian variety with multiplication by $K$. This is an analog of 12.2.
\medskip
{\bf 12.5. Complete multiplication.} Here we consider the case $\goth
q=m=1$, i.e. $\goth K=\Cal W$ and $g=r$.
\medskip
{\bf Lemma 12.5.1.} In this case the condition $N=0$
implies that the CM-type $$\sum_{i=1}^r\goth r_i\alpha_i\eqno{(12.5.2)}$$ of the action of
$\goth K$ on
on $(L,V)$ has the property: all $\goth r_i$ are 0 or 1.
\medskip
{\bf Proof.} $N=0$ means that
the action of $T\in\w$ on $V$ is simply multiplication by $\theta$. We write the
CM-type $\sum_{i=1}^r\goth r_i\chi_i$ in the form $\sum_{i=1}^n\chi_{\alpha_i}$ where
$\alpha_1,...,\alpha_n\in[1,...,r]$ are not necessarily distinct. Let $l_1$ be an (only)
element of a basis of $L\otimes_{\w_{\goth K}}\goth K$ over $\goth K$ and $e_1,...,e_n$ a
basis of $V$
over $\p$ such that the action of $\goth K$ on $V$ is given by the formulas
$$k(e_i)=\chi_{\alpha_i}(k)e_i, \ \ \ k\in \goth K$$
Multiplying $e_i$ by scalars if necessary, we can assume that $l_1=\sum e_i$. Therefore,
if $\alpha_i=\alpha_j$ (i.e. not all $\goth r_*$ in (12.5.2) are 0, \ 1) then the
$e_{\alpha_i}$-th coordinate of any element of $L$ coincide with its $e_{\alpha_j}$-th
coordinate, hence $L$ does not $\p$-generate $V$ --- a contradiction. $\square$
\medskip
{\bf 12.5.3.} Let $M$ be a t-motive of
rank $r$ and dimension $n$ having multiplication by $\w_{\goth K}$. Recall that we
consider only
the case $N=0$. This means that the character of the action of $\goth K$ on $M/\tau M$
is isomorphic to $\sum_{i=1}^r\goth r_i\alpha_i$. Since $E(M)=(M/\tau M)^*$ we get that
the character of the action of $\goth K$ on $E(M)$ is the same. If
$$\hbox{all $\goth r_i$ are 0 or 1}\eqno{(12.5.4)}$$
we shall use the terminology that $M$ has the CM-type $\Phi\subset \{\alpha_1, ... ,
\alpha_r \}$
where $\Phi$ is defined by the condition $\alpha_i\in \Phi \iff \goth r_i=1$.

It is easy to see that this case occurs for uniformizable $M$. Really, if $M$ is
uniformizable then
the action of $\goth K$ can be prolonged on $(L(M),V(M))$, and the character
of the action of $\goth K$ on $V(M)$ coincides with the one on $E(M)$. The result follows
from Lemma 12.5.1.
\medskip
{\bf Lemma 12.5.5.} There exists a canonical isomorphism $\gamma$ from the set of
inclusions $\alpha_1, ... , \alpha_r$ to the set of points $ \theta_{\alpha_1}, ... ,
\theta_{\alpha_r}$ of $X$ over $\theta\in P^1(\p)$.
\medskip
{\bf Proof.} A point $t\in X$ over $\theta\in P^1(\p)$ defines a function
$\vf_t: \goth K_{C} \to P^1(\p)$ --- the value of
an element $f\in \goth K_{C}$ treated as a function on $X$ at the point $t$. This
function
must satisfy the standard axioms of valuation and the condition $\vf_t(T)=\theta$. Let
$\alpha_i$
be an inclusion of $\goth K$ to $\p$ over $\iota$. It defines a valuation
$\vf_{\alpha_i}: \goth K_{C} \to P^1(\p)$ by the formula $\vf_{\alpha_i}(k\otimes
f)=\alpha_i(k) f(\theta)$,
where $k\in\goth K$, $f\in \p(T)$. We define $\gamma(\alpha_i)$ by the condition
$\vf_{\gamma(\alpha_i)}=\vf_{\alpha_i}$; it is easy to see that $\gamma$ is an
isomorphism. $\square$
\medskip
{\bf Theorem 12.6.} For any above \{$\goth K$, $\Phi$\} there exists an
t-motive $(M,\tau)$ with
complete multiplication by $\goth K$ having CM-type $\Phi$.
\medskip
{\bf Proof (Drinfeld).} We denote the divisor $\sum_{\alpha_i\in\Phi}\gamma({\alpha_i})$
by
$\theta_\Phi$. We construct a $\Cal F$-sheaf $F$ of dimension 1 over
$\goth K$ which will give us $M$.
Let $\fr$ be the Frobenius map on $ \Pic_0(X)$. It is an algebraic
map, and the $\fr-\Id: \Pic_0(X) \to \Pic_0(X)$ is an algebraic map as
well. Since the action of $\fr$ on the tangent space of $ \Pic_0(X)$
at 0 is the zero map, the action of $\fr-\Id$ on the tangent space of
$ \Pic_0(X)$ at 0 is the minus identical map and hence $\fr-\Id$ is an
isogeny of $ \Pic_0(X)$. Particularly, there exists a divisor $D$ of
degree 0 on $X$ such that we have the following equality in $\Pic_0(X)$:
$$\fr(D)-D=-\theta_\Phi+n\infty \eqno{(12.6.0)}$$
This means that if we let $F=F_\Phi=O(D)$ then there exists a rational map
$\tau_X=\tau_{X,\Phi}: F^{(1)}\to F$ such
that
$$\Div(\tau_X)=\theta_\Phi-n\infty\eqno{(12.6.1)}$$
The pair $(F_\Phi, \tau_{X,\Phi})$ is the desired $\Cal F$-sheaf.
\medskip
{\bf Remark.} It is easy to see that if the genus of $X$ is $> 0$ then different CM-types $\Phi_1$, $\Phi_2$ give us different sheaves $F_{\Phi_1}$, $F_{\Phi_2}$, while if the genus of $X$ is 0 then $F_{\Phi_1}=F_{\Phi_2}=\Cal O$, but the maps $\tau_{X,\Phi_1}$, $\tau_{X,\Phi_2}$ are clearly different.
\medskip
Let $U_0=X-\{\infty\}$ be an open part of
$X$. We denote $F(U_0)$ by $\Cal M$, hence $F^{(1)}(U_0)=\Cal M^{(1)}$.
Since the support of the negative part of the right hand side of
12.6.1 is $\{\infty\}$,
we get that the (a priory rational) map $\tau_X(U_0): \Cal M^{(1)}\to \Cal M$ is
really a map of
$\w_{\goth K}$-modules.

Let $M$ be a $\p[T]$-module obtained from $\Cal M$ by restriction of scalars from
$\w_{\goth K}$ to $\p[T]$. Construction $F\mapsto M$ is functorial, and we denote this functor by $\delta$. Further, we denote by $\alpha$ the tautological isomorphism
$\Cal M\to M$. $M$ is a free $r$-dimensional $\p[T]$-module, and (because
$M^{(1)}$ is isomorphic to $M$) the
same restriction of scalars of $\tau_X(U_0)$ defines us a $\p[T]$-skew map from $M$ to
$M$ denoted by $\tau$ (skew means that $\tau(zm)=z^q\tau(m)$, $z\in\p$). $\tau$ is
defined by the formula $\tau(m)=\alpha\circ \tau_X((\alpha^{-1}(m))^{(1)})$.

It is easy to check that
$(M,\tau)$ is the required t-motive. Really, $M$ is a $\w_{\goth K}$-module, and
$\tau$
commutes with this multiplication. The fact that the positive part
of the right hand
side of 12.6.1 is $\theta_\Phi$ means that 1.13.2 holds for $M$ and that the CM-type of
the action of
$\w_{\goth K}$ is $\Phi$.
\medskip
{\bf Remark 12.6.2.} It is easy to prove for this case that $M$ is a free
$\p[\tau]$-module.
Really, it is
sufficient to prove (see [G], Lemma 12.4.10)
that $M$ is finitely generated as a $\p[\tau]$-module. We choose $D$ such that
$\infty\not\in \Supp (D)$. There exists $P\in \goth K_{C}^*$
such that $\tau_X(U_0): \Cal M^{(1)}\to \Cal M$ is multiplication by $P$ (recall that both
$\Cal M^{(1)}$, $\Cal M$ are $\w_{\goth K}$-submodules of $\goth K$). 12.6.0 implies that
$-\ord_\infty(P)=n$.
There
exists a number $n_1$ such that $$\hbox{(a) $h^0(X,\Cal O(D+n_1\infty))>0$; \ \ (b)
for any $k\ge 0$ we have}$$ $$h^0(X,\Cal O(D+(n_1+k)\infty))=h^0(X,\Cal
O(D+n_1\infty))+k\eqno{(12.6.3)}$$ $$h^0(X,\Cal O(D^{(1)}+(n_1+k)\infty))=h^0(X,\Cal
O(D^{(1)}+n_1\infty))+k\eqno{(12.6.4)}$$
It is sufficient to prove that if $g_1,...,g_k$ are elements of a basis of $H^0(X,\Cal
O(D+(n_1+n)\infty))$, then
for any $Q\in \Cal M$ the element $\alpha(Q)\in M$ is generated by
$\alpha(g_1),...,\alpha(g_k)$ over
$\p[\tau]$. We prove it by induction by $n_2:=-\ord_\infty(Q)$. If $n_2\le n_1+n$ the
result is trivial. If not
then 12.6.3, 12.6.4 imply that the multiplication by $P$ defines an isomorphism
$$H^0(X, \Cal O(D^{(1)}+(n_2-n)\infty))/H^0(X, \Cal O(D^{(1)}+(n_2-n-1)\infty))\to$$
$$\to H^0(X,\Cal O(D+n_2\infty))/H^0(X, \Cal O(D+(n_2-1)\infty))$$
This means that $\exists Q_1\in H^0(X,\Cal O(D^{(1)}+(n_2-n)\infty))$,
$-\ord_\infty(Q_1)=n_2-n$ such
that
$-\ord_\infty(Q-PQ_1)\le n_2-1$. An element $Q_1^{(-1)}\in \Cal M$ exists;
since $\alpha(Q)=\tau(\alpha(Q_1^{(-1)}))+\alpha(Q-PQ_1)$, the result follows by
induction.
$\square$
\medskip
If $\goth K$ and $\Phi$ are given then the construction of the Theorem 12.6 defines $F$ uniquely up to tensoring by $O(D)$ where $D\in \Div (X(\goth K))$. We denote the set of these $F$ by $F$(\{$\goth K$, $\Phi$\}), and we denote by $M$(\{$\goth K$, $\Phi$\}) the set $\delta(F$(\{$\goth K$, $\Phi$\})). Further, we denote by $\Phi'=\{\alpha_1, ... ,
\alpha_r\}-\Phi$ the complementary CM-type.
\medskip
{\bf Theorem 12.7.} Let $M\in M(\{\goth K,\Phi\})$. Then $M'$ exists, and $M'\in M(\{\goth K,\Phi'\})$. More exactly, if $F\in F$(\{$\goth K$, $\Phi$\}) then $F^{-1}\otimes \Cal D^{-1}\in F$(\{$\goth K$, $\Phi'$\}) where $\Cal D$ is the different
sheaf on $X$, and if $M=\delta(F)$ then $M'=\delta(F^{-1}\otimes \Cal D^{-1})$.
\medskip
{\bf Proof.} Let $G$ be any invertible sheaf on $X$. We have a
\medskip
{\bf Lemma 12.7.0.} There exists the canonical isomorphism $\varphi_G:
\pi_*(G^{-1}\otimes \Cal D^{-1}) \to \Hom_{P^1}(\pi_*(G), \Cal O)$.
\medskip
{\bf Proof.} At the level of affine open sets $\varphi_G$ comes from the trace bilinear form of field extension $\goth K/\n F_q(T)$. Concordance with glueing is obvious. $\square$
\medskip
We need the relative version of this lemma. Let $G_1$, $G_2$ be invertible sheaves on $X$,
$\rho: G_1 \to G_2$ any rational map. Obviously there exists a rational map
$\rho^{-1}: G_1^{-1} \to G_2^{-1}$. Recall that we denote by
$\rho^{inv}: G_2 \to G_1$ the rational map which is inverse to $\rho$
respectively the composition. The map $\pi_*(\rho^{-1}\otimes \Cal D^{-1}): \pi_*(G_1^{-1}\otimes \Cal D^{-1})\to \pi_*(G_2^{-1}\otimes \Cal D^{-1})$ is obviously defined. The map (denoted by $\beta(\rho)$) from $\Hom_{P^1}(\pi_*(G_1), \Cal O)$ to $\Hom_{P^1}(\pi_*(G_2),
\Cal O)$ is defined as follows at the level of affine open sets: let
$\gamma\in\Hom_{P^1}(\pi_*(G_1), \Cal O)(U)$ where $U$ is a sufficiently
small affine subset of $P^1$, such that we have a map $\gamma(U): \pi_*(G_1)(U)
\rightarrow \Cal O(U)$. Then $(\beta(\gamma))(U)$ is the composition
map $\gamma(U)\circ \pi_*(\rho^{inv})(U)$:
$$\pi_*(G_2)(U)\overset{\pi_*(\rho^{inv})(U)}\to{\longrightarrow}
\pi_*(G_1)(U)\overset{\gamma(U)}\to{\to}\Cal O(U)$$

{\bf Lemma 12.7.1.} The above maps form a commutative diagram:
$$\matrix
\pi_*(G_1^{-1}\otimes \Cal D^{-1})&\overset{\pi_*(\rho^{-1}\otimes \Cal
D^{-1})}\to{\longrightarrow}&\pi_*(G_2^{-1}\otimes \Cal D^{-1})&&\\&&&&\\
\varphi_{G_1}\downarrow&&\varphi_{G_2}\downarrow &&\\&&&&\\
\Hom_{P^1}(\pi_*(G_1), \Cal
O)&\overset{\beta(\rho)}\to{\longrightarrow}&\Hom_{P^1}(\pi_*(G_2),
\Cal O)&&\square
\endmatrix$$

We apply this lemma to the case $\{\rho: G_1 \to G_2\}=\{\tau_{X, \Phi}:
F^{(1)}\to F\}$. We have:
$$\Div(\tau_{X, \Phi}^{-1}\otimes \Cal D^{-1})=-\Div(\tau_{X, \Phi})=-\theta_{\Phi}+n\infty$$
Futher, we multiply $\tau_{X, \Phi}^{-1}\otimes \Cal D^{-1}$ by $T-\theta$. We have:
$$\Div((T-\theta)\tau_{X, \Phi}^{-1}\otimes \Cal
D^{-1})=\Div(T-\theta)+\Div(\tau_{X, \Phi}^{-1}\otimes \Cal D^{-1})=\theta_{\Phi'}-(r-n)\infty$$
i.e. $(T-\theta)\tau_{X, \Phi}^{-1}\otimes \Cal D$ is one of
$\tau_{X,\Phi'}$, i.e. $F^{-1}\otimes \Cal D^{-1}\in F$(\{$\goth K$, $\Phi'$\}). Further, $(T-\theta)\beta(\tau_{X, \Phi})$ is the map which is used in the definition of duality of
$M$. This means that the lemma 12.7.1 implies the theorem. $\square$
\medskip
{\bf Remark 12.8.} There exists a simple proof of the second part of the Theorem 5 for uniformizable abelian
t-motives $M$ with complete multiplication by $\w_\goth K\subset \goth K$. Recall that this second part is the proof of 2.7 for $M$. Really, let us
consider the diagram 2.5. The CM-types
of action of $\goth K$ on $\Lie(M)$ and on $E(M)$ coincide,
and the CM-types of action of $\goth K$ on a vector space and on its
dual space coincide. This means that the CM-type of $V^*$ is $\Phi$
and the CM-type of $V'$ is $\Phi'$. Further, $\gamma_D$ of 2.5 commutes with complete multiplication: this follows immediately for example from a description of $\gamma_D$ given in Remark 5.2.8. Really, all homomorphisms of 5.2.9 commute with complete multiplication. For example, this condition for $\delta$ of 1.11.1 is written as follows: if $k\in \goth K$, $\goth m_k(M)$, resp. $\goth m_k(M')$ is the map of complete multiplication by $k$ of $M$, resp. $M'$, then $(\goth m_k(M)\otimes Id) \circ \delta= (Id \otimes \goth m_k(M'))\circ \delta$ --- see any textbook on linear algebra.

Finally, since
$\Phi\cap \Phi'=\emptyset$ and the map $\vf'\circ \gamma_D \circ
\vf^*$ commutes with complete multiplication, we get that it must be
0.
\medskip
{\bf 13. Miscellaneous.}

\nopagebreak
\medskip
Let now $(L,V)$ be from 12.1, case $\goth q=m=1$, i.e. $\goth K=\Cal W$
and $r=g$, and let the ring of complete multiplication be the maximal
order $\w_\goth K$. We identify $\w$ and $\n Z_\infty$ via $\iota$, i.e. we consider $\goth K$ as an extension of $\n F_q(\theta)$. Let $\Phi$ be the CM-type of the action of $\goth
K$ on $V$.
This means that --- as an
$\w_\goth K$-module --- $L$ is isomorphic to $I$ where $I$ is an ideal of
$\w_\goth K$. The class of $I$ in $\Cl(\w_\goth K)$ is defined by $L$
and $\Phi$ uniquely; we denote it by $\Cl(L,\Phi)$.
\medskip
{\bf Remark.} $\Cl(L,\Phi)$ depends on $\Phi$, because the action of $\w_\goth K$ on $V$ depends on $\Phi$. Really, let $a\in L \subset V$, $a=(a_1,...,a_n)$ its coordinates, $\Phi=\{\alpha_{i_1},...,\alpha_{i_n}\}\subset \{\alpha_{1},...,\alpha_{r}\}$ and $k\in \w_\goth K$. Then $ka$ has coordinates $(\alpha_{i_1}(k)a_1,...,\alpha_{i_n}(k)a_n)$, i.e. depends de $\Phi$. Particularly, the $\w_\goth K$-module structure on $L$ depends on $\Phi$, and hence $\Cl(L,\Phi)$ depends on $\Phi$. For example, if $n=1$, $r=2$, $\Phi_1=\{\alpha_{1}\}$, $\Phi_2=\{\alpha_{2}\}$, then $\Cl(L,\Phi_2)$ is the conjugate of $\Cl(L,\Phi_1)$.
\medskip
{\bf Theorem 13.1.} $\Cl(L',\Phi')=(\Cl(\goth d))^{-1}
(\Cl(L,\Phi))^{-1}$ where
$\goth d$ is the different ideal of the ring extension $\w_\goth K/\w$.
\medskip
{\bf Proof.} This theorem follows from the above results;
nevertheless, I give here an explicit elementary proof. Let $a_*=(a_1, ..., a_r)^t$ be a
basis (considered
as a vector column) of $\goth K$ over $\n F_q(\theta)$ and $b_*=(b_1,
..., b_r)^t$ the dual basis. Recall that it satisfies 2 properties:
$$(1) \ \ \forall i\ne j \ \ \alpha_i(a_*)^t\alpha_j(b_*)=0 \ \ \
(\hbox{ i.e. } \sum_{k=1}^r \alpha_i(a_k) \alpha_j(b_k)=0
)\eqno{(13.1.1)} $$
(2) For $x\in \goth K$ let $\goth m_{x,a_*}$ (resp $\goth m_{x,b_*}$)
be the matrix of multiplication by $x$ in the basis $a_*$ (resp.
$b_*$). Then for all $x\in \goth K$ we have
$$\goth m_{x,a_*}=\goth m_{x,b_*}^t\eqno{(13.1.2)} $$

We define $\goth I_{n,r-n}$
as an $r\times r$ block matrix $\left(\matrix 0&I_{r-n}\\ -I_n&0
\endmatrix\right)$, and we define a new basis $\tilde b_*=(\tilde b_1,
..., \tilde b_r)^t$ by
$$\tilde b_*=\goth I_{n,r-n}b_*\eqno{(13.1.3)} $$
(explicit formula: $(\tilde b_1, ..., \tilde b_r)=(b_{n+1},
..., b_r, -b_1, ..., -b_n)$).

We can assume that $\Phi=\{\alpha_1, ... , \alpha_n \}$. Since $L$ has multiplication by $\w_\goth K$ and the CM-type of this
multiplication is $\Phi$, it is possible to choose $a_*$ such that $L\subset
\p^{n}$ is generated over $\n Z_\infty$ by $e_1, ... , e_r$ where
$$e_i=(\alpha_1(a_i), ... , \alpha_n(a_i))\eqno{(13.1.4)}$$ Let $\hat L\subset \p^{r-n}$ be
generated over $\n Z_\infty$ by $\hat e_1, ... , \hat e_r$ where
$$\hat e_i=(\alpha_{n+1}(\tilde b_i), ... , \alpha_r(\tilde b_i))\eqno{(13.1.5)}$$

{\bf Lemma 13.1.6. } $L'=\hat L$.
\medskip
{\bf Proof. } Let $A$ (resp. $B$) be a matrix whose lines are the lines of
coordinates of $e_1, ... , e_n$ (resp. $e_{n+1}, ... ,
e_r$) in 13.1.4, and $C$ (resp. $D$) a matrix whose lines are the lines of
coordinates of $\hat e_1, ... , \hat e_{r-n}$ (resp. $\hat e_{r-n+1}, ... ,
\hat e_r$) in 13.1.5. By definition of Siegel matrix, we have $L=\goth L(BA^{-1})$,
$\hat L=\goth L(DC^{-1})$ ($\goth L$ is defined in 3.1, 3.2). So, it is sufficient to prove that
$(BA^{-1})^t=DC^{-1}$, i.e. $A^tD=B^tC$. This follows immediately from
the definition of $A,B,C,D$ and (13.1.1). $\square$
\medskip
For $x\in \w_\goth K$ we denote by $\goth
M_x(L)$ the matrix of multiplication by $x$ in the basis $e_*$ (see
the notations of Remark 3.8). Obviously $\goth M_x(L)=\goth
m_{x,a_*}$.

Let now $\w_\goth K$ acts on $\p^{r-n}$ (the ambient space of $L'$) by CM-type
$\Phi'$. According (13.1.2) and
(13.1.3), the matrix of the action of $x\in \w_\goth K$ in the basis
$\tilde b_*$ is
$$\goth I_{n,r-n}\goth m_{x,a_*}^t\goth I_{n,r-n}^{-1}\eqno{(13.1.7)} $$
Let $\goth M$, $\goth M'$ be from Remark 3.8. Formula 3.8.4
shows that
$$\goth M'=\goth I_{n,r-n}\goth M^t\goth I_{n,r-n}^{-1}\eqno{(13.1.8)} $$
Formulas (13.1.7) and (13.1.8) --- because of Lemma 1.10.3 --- prove the theorem. $\square$
\medskip
{\bf 13.2. Compatibility with the weak form of the main theorem of
complete multiplication.}

\nopagebreak
\medskip
The reader can think that Theorem 13.1 is incompatible with the main theorem of complex
multiplication, because of the $-1$-th power in its statement. The reason is a bad choice
of notations of Shimura, he affirms that the CM-type of an abelian variety $A$ over a number field coincides with the CM-type
of $A'$, while we see that it is really the complement. Since an analog of even the weak form of the main theorem of complex multiplication --- Theorem 13.2.6 --- for the function field case is not proved yet, the main result of the present section --- Theorem 13.2.8 --- is conditional: it affirms that if this weak form of the main
theorem --- Conjecture 13.2.7 --- is true for a t-motive with complete multiplication $M$, then it is true for $M'$ as well. By the way, even if it will turn out that the statement of the Conjecture 13.2.7 is not correct, the proof of 13.2.8 will not be affected, because the main
ingredient of the proof is the formula 13.2.10 "neutralizing" the $-1$-th power of the
Theorem 13.1.

Let us recall some definitions of [Sh71], Section 5.5. We consider an abelian variety
$A=\n C^n/L$ with
complex multiplication by $K$. The set $\Hom(K, \bar \n Q)$ consists of $n$ pairs of
mutually
conjugate inclusions $\{\vf_1,\bar \vf_1, ..., \vf_n,\bar\vf_n\}$. $\Phi$ is a subset of
the set
$\Hom(K, \bar \n Q)$ such that $\forall i=1,...,n$ we have:
$$\hbox{$\Phi\cap\{\vf_i,\bar \vf_i\}$ consists of one element.}\eqno{(13.2.1)}$$
It is defined by the condition that the action of complex multiplication
of $K$ on $\n C^n$ is isomorphic to the direct sum of the elements of $\Phi$. Let $F$ be
the
Galois envelope of $K/\n Q$,
$$G:=\Gal (F/\n Q), \ \ H:=\Gal(F/K), \ \
S:=\bigcup_{\alpha\in\Phi}H\alpha\eqno{(13.2.2)}$$
(the elements of Galois group act on $x\in F$ from the right, i.e. by the formula
$x^{\alpha\beta}=(x^\alpha)^\beta$; for $\alpha\in\Phi$ we denote by $\alpha$ also a
representative in $G$ of the coset $\alpha$). We denote
$$H^{ref}:=\{\gamma\in G\vert S\gamma=S\}\eqno{(13.2.3)}$$
and let $K^{ref}$ be the subfield of $F$ corresponding to $H^{ref}$. We have:
$$H^{ref}S^{-1}=S^{-1}\eqno{(13.2.4)}$$
i.e. $S^{-1}$ is an union of cosets of $H^{ref}$ in $G$. We can identify these cosets
with elements of $\Hom(K^{ref}, \bar \n Q)$. $\Phi^{ref}\subset \Hom(K^{ref}, \bar \n Q)$
is, by definition, the set of these cosets. There is a map $\det
\Phi^{ref}:{K^{ref}}^\times \to K^\times$ defined as follows:
$$\det \Phi^{ref}(x):=\prod_{\alpha\in\Phi}\alpha(x)\eqno{(13.2.5)}$$
(it follows easily from the above formulas and definitions that $\det \Phi^{ref}(x)$
really belongs to $K^\times$). It can be extended to the group of ideles and factorized to
the group of classes of ideals, we denote this map
by $\det_{Cl} \Phi^{ref}: \Cl(K^{ref}) \to \Cl(K)$. Finally, let
$\theta^{ref}: \Gal(K^{ref\ Hilb}/K^{ref}) \to \Cl(K^{ref})$ be an isomorphism defined by
the Artin reciprocity law.

We consider the case $\End(A)=O_K$. In this case $L$ is isomorphic to an ideal of $O_K$,
its class
is well-defined by the class of isomorphism of $A$, we denote it by $\Cl(A)$.
\medskip
{\bf Theorem 13.2.6.} $A$ is defined over $K^{ref\ Hilb}$;

For any $\gamma\in\Gal(K^{ref\ Hilb}/K^{ref})$ we have

$\Cl(\gamma(A))=\det_{Cl} \Phi^{ref}\circ \theta^{ref} (\gamma)^{-1}(\Cl(A))$. $\square$
\medskip
This is a weak form of [SH71], Theorem 5.15 --- the main theorem of
complex multiplication.
\medskip
Now we define analogous objects for the function field case in order to formulate a
conjectural
analog of Theorem 13.2.6. Let $\goth K$, $\Phi$ be from 12.5.3. $\goth K^{ref}$,
$\Phi^{ref}$,
$\det \Phi^{ref}$ are defined by the same formulas 13.2.2 -- 13.2.5 like in the number
field case
($\n Q$ must be replaced by $\n F_q(T)$). The facts that 13.2.1 has no meaning in the
function field case and
that the order of $S$ is not necessarily the half of the order of $G$ do not affect the
definitions.

The $\infty$-Hilbert class field of $\goth K$ (denoted by $\goth K^{Hilb \ \infty}$) is an abelian
extension of
$\goth K$ corresponding to the subgroup
$$\goth K_\infty^*\cdot \prod_{v\ne\infty}O_{\goth K_v}^*\cdot \goth K^*$$
of the idele group of $\goth K$. We have an isomorphism $\theta: \Gal(\goth K^{Hilb \ \infty}/\goth K) \to
\Cl(\w_{\goth K})$.

We formulate the function field analog of Theorem 13.2.6 only for the case when
\medskip
{\bf (*)} There exists only one
point over $\infty\in P^1(\n F_q)$ in the extension $\goth K^{ref}/\n F_q(T)$.
\medskip
In this case the field
$\goth K^{ref \ Hilb \ \infty}$ and the ring $\w_{\goth K^{ref}}$ are naturally defined, and we have
an
isomorphism
$\theta^{ref}: \Gal(\goth K^{ref \ Hilb \ \infty}/\goth K^{ref}) \to \Cl(\w_{\goth K^{ref}})$.

Let $M$ be an uniformizable
t-motive of rank $r$ and dimension $n$ having complete multiplication
by $\w_\goth K$, and $\Phi$ its CM-type. $\Cl(M)$ is defined like $\Cl(A)$ in the number
field case,
it is $\Cl(L,\Phi)$ of 13.1.
\medskip
{\bf Conjecture 13.2.7.} If (*) holds, then $M$ is defined over $\goth K^{ref\ Hilb \
\infty}$, and for any
$\gamma\in\Gal(\goth K^{ref\ Hilb \ \infty}/\goth K^{ref})$
we have $\Cl(\gamma(M))=\det_{Cl} \Phi^{ref}\circ \theta^{ref} (\gamma)^{-1}\Cl(M)$.
\medskip

Now we can formulate the main theorem of this section.
\medskip
{\bf Theorem 13.2.8.} If conjecture 13.2.7 is true for $M$ then it is true for $M'$.
\medskip
{\bf Proof.} It follows immediately from the functional analogs of 13.2.2 -- 13.2.4 that
$$(\goth K,\Phi')^{ref}=(\goth K^{ref},(\Phi^{ref})')\eqno{(13.2.9)}$$
Further,
$$\hbox{det}_{Cl} {\Phi'}^{ref}=(\hbox{det}_{Cl} \Phi^{ref})^{-1}\eqno{(13.2.10)}$$
Really, $\det \Phi^{ref}(x) \cdot \det (\Phi^{ref})'(x) = N_{\goth
K^{ref}/\n F_q(T)}(x)\in\n F_q(T)^\times$, hence gives the trivial class of ideals (we use
here (13.2.9).
Finally, for $\gamma\in\Gal(\goth K^{ref})$ we have
$$(\gamma(M))'=\gamma(M')\eqno{(13.2.11)}$$
The theorem follows immediately from 13.1, 13.2.10, 13.2.11 (recall that $\Cl(M)$ is
$\Cl(L,\Phi)$ of 13.1). $\square$
\medskip
{\bf 13.3. Some explicit formulas.} We give
here an elementary explicit proof of the theorem 12.7 in two simple
cases: $\goth K=\n F_{q^r}(T)$ and $\n F_q(T^{1/r})$. By the way,
since the extension $\n F_{q^r}(T)/\n F_{q}(T)$ is not absolutely
irreducible, formally this case is not covered by the theorem 12.7.
\medskip
{\bf Case $\w_\goth K= \n F_{q^r}[T]$.} Let $\alpha_i$, where $i=0, ... ,r-1$, be inclusions $\goth K \to \p$. For $\omega \in\n F_{q^r}$
 we have $\alpha_i(\omega )=\omega^{q^{i}}$. Let $$0 \le
i_1 < i_2 < ... <
i_n \le r-1\eqno{(13.3.0)}$$ be numbers such that $\Phi=\{\alpha_{i_j}\}$, $j=1, ... , n$.
We consider the following t-motive $M=M(\goth K, \Phi)$. Let $e_1, ... , e_n$ be a basis of $M_{\p[\tau]}$ such that $\goth
m_\omega(e_j)=\omega^{q^{i_j}}e_j$ and the multiplication by $T$ is
defined by formulas
$$Te_1=\theta e_1 + \tau^{i_1-i_n+r}e_n\eqno{(13.3.1)} $$
$$Te_j=\theta e_j + \tau^{i_j-i_{j-1}}e_{j-1}, \ \ j=2, ... , n
\eqno{(13.3.2)} $$

It is easy to check that $M$ has complete multiplication by $\w_\goth K$, and its CM-type is $\Phi$.
\medskip
{\bf Remark.} It is possible to prove that $M(\goth K, \Phi)$ is the only t-motive having these properties; we omit the proof.
\medskip
{\bf Proposition 13.3.3.} For $\w_\goth K= \n F_{q^r}[T]$ we have: $M(\goth K, \Phi)'=M(\goth K, \Phi')$.
\medskip
{\bf Proof.} Elements $\tau^je_k$ for $k=1, ... ,n$, $j=0, ..., i_{k+1}-i_{k}-1$ for $k<n$
and $j=0,
..., i_1-i_n+r-1$ for $k=n$ form a basis of $M_{\p[T]}$. Let us
arrange these elements in
the lexicographic order ($\tau^{j_1}e_{k_1}$ precedes to
$\tau^{j_2}e_{k_2}$ if $k_1 <
k_2$) and make a cyclic shift of them by $i_1$ denoting $e_1$ by $f_{i_1+1}$,
$\tau^{i_2-i_1-1}e_1$ by $f_{i_2}$ etc. until
$\tau^{i_1-i_n+r-1}e_n=f_{i_1}$. Formulas
13.3.1, 13.3.2 become
$$\tau(f_i)=f_{i+1} \hbox{ if } i\not\in \{i_1, ..., i_n\} $$
$$\tau(f_i)=(T-\theta)f_{i+1} \hbox{ if } i\in \{i_1, ..., i_n\} $$
($i \mod r$, i.e. $f_{r+1}=f_1$). Formula 1.10.1 shows that in the
dual basis $f'_*$ we
have
$$\tau(f'_i)=f'_{i+1} \hbox{ if } i\in \{i_1, ..., i_n\} $$
$$\tau(f'_i)=(T-\theta)f'_{i+1} \hbox{ if } i\not\in \{i_1, ..., i_n\} $$
which proves the proposition. $\square$
\medskip
{\bf Case $\w_\goth K=\n F_q[T^{1/r}]$, $(r,q)=1$.} In order to define $M(\goth K, \Phi)$ we need more notations. We denote
$\theta^{1/r}$ and $T^{1/r}$ by $\goth s$ and $S$ respectively, and
let $\zeta_r$ be a primitive $r$-th root of 1. Let $\alpha_i$, $i_1 < i_2 < ... < i_n$ and $\Phi$ be the
same as in the case $\w_\goth K= \n F_{q^r}[T]$. We
have $\alpha_i(S)=\zeta_r^iS$. Further, we consider an overring $\p[S,\tau]$ of $\p[T,\tau]$ ($S$ is in the
center of this ring), and we consider the category of modules
over $\p[S,\tau]$ such that the condition 1.9.2 is changed by a
weakened condition 13.3.4 (here $A_{S,0}\in M_n(\p)$ is defined by the formula $Se_*=A_Se_*$, where $A_S\in M_n(\p)[\tau]$, $A_S=\sum_{i=0}^*A_{S,i}\tau^i$):
$$A_{S,0}^r=\theta I_n + N\eqno{(13.3.4)}$$
Let $\bar M$ be a $\p[S,\tau]$-module such that $\dim \bar
M_{\p[S]}=1$, $f_1$ the only element of a basis of $\bar M_{\p[S]}$ and
$$\tau f_1=(S-\zeta_r^{i_1}\goth s)\cdot ... \cdot
(S-\zeta_r^{i_n}\goth s)f_1$$
By definition, $M=M(\goth K, \Phi)$ is the restriction of scalars from
$\p[S,\tau]$ to $\p[T,\tau]$ of $\bar M$. Like in the case $\w_\goth K= \n F_{q^r}[T]$, it is easy to check that $M$ has complete multiplication by $\w_\goth K$ with CM-type $\Phi$, and it is possible to prove that it is the only t-motive having these properties.
\medskip
{\bf Proposition 13.3.5.} For $\w_\goth K= \n F_q[T^{1/r}]$, $(r,q)=1$ we have: $M(\goth K, \Phi)'=M(\goth K, \Phi')$.
\medskip
{\bf Proof.} For $i=1, ..., r$ we denote $f_i=S^{i-1}f_1$. These $f_*=f_*(\Phi)$
form a basis of $M_{\p[T]}$, and the matrix $Q=Q(f_*,\Phi)$ of
multiplication of $\tau$ in this basis has the following description.
We denote by $\sigma_k(\Phi)$ the elementary symmetric polynomial
$\sigma_k(\zeta_r^{i_1}, ... , \zeta_r^{i_n})$.

The first line of $Q$ is
$$\sigma_n(\Phi)\goth s^n \ \ \ \sigma_{n-1}(\Phi)\goth s^{n-1} \ \ \ ...
\ \ \  \sigma_1(\Phi)\goth s \ \ \ 1 \ \ \ 0 \ \ \ ... \ \ \ 0$$
and its $i$-th line is obtained from the first line by 2 operations:

1. Cyclic shift of elements of the first line by $i-1$ positions to the right;

2. Multiplication of the first $i+n-r$ elements of the obtained line by $T$.

We consider another basis $g_*=g_*(\Phi)$ of $M_{\p[T]}$ obtained by
inversion of order of $f_i$, i.e. $g_i=f_{r+1-i}$. The elements of
$Q(g_*)$ are obtained
by reflection of positions of elements of $Q(f_*)$ respectively the
center of the matrix.

The theorem for the present case follows from the formula
$$Q(f_*,\Phi)Q(g_*,\Phi')^t=(T-\theta)I_r$$
whose proof is an elementary exercise: let $\Phi'=\{j_1, ... ,
j_{r-n}\}$; we apply
equality
$$\sigma_k(x_1, ... x_r)=\sum_l \sigma_l(x_{i_1}, ... ,
x_{i_n})\sigma_{k-l}(x_{j_1}, ...
, x_{j_{r-n}})$$
to $1, \zeta_r, ... , \zeta_r^{r-1}$. $\square$
\medskip
{\bf 13.4. Reduction.} Recall notations of 1.16.  Let $L$ be a finite extension of $\n F_q(\theta)$, $\goth p$ a valuation of $L$ over a valuation $P\ne\infty$ of $\n F_q(\theta)$, and we denote $\iota^{-1}(P)\subset \w$ by $\Cal P$. Let $M$ be a t-motive defined over $L$ having a good ordinary reduction $\tilde M$ at $\goth p$ and such that the dual $M'$ exists. According 1.15.1, the $L$-structure on $M'$ is well-defined. We denote by $M_{\Cal P,0}$ the kernel of the reduction map
$M_\Cal P \to \tilde M_\Cal P$. Condition of ordinarity means that  $M_{\Cal P,0}=(\w/\Cal P)^n$.
\medskip
{\bf Conjecture 13.4.1.} For the above $M$, $M'$ we have:
\medskip
$M_{\Cal P,0}$ and $M'_{\Cal P,0}$ are mutually dual with respect to the pairing of Remarks 4.2, 5.1.6 (recall that conjecturally $M'$ also has good ordinary reduction at $\goth p$).
\medskip
{\bf Proof for a particular case:} $M$ is a Drinfeld module, $\Cal P=T$.
\medskip
(1.9.1) for $M$ has a form
$$Te=\theta e+a_1\tau e+... a_{r-1}\tau^{r-1}e+\tau^r e$$
Condition of good ordinary reduction means $a_i\in L$, $\ord_\goth p(a_i)\ge 0$,
$\ord_\goth p a_1=0$. Let $x\in
M_T$, $y\in M'_T$; we can consider $x$ (resp. $y$) as an element of
$\p$ (resp. $\p^{r-1}$) satisfying some polynomial equation(s). Considering
Newton polygon of these polynomials we get immediately (1) for both
$M$, $M'$. Let $y=(y_1, ... ,y_{r-1})$ be the coordinates of $y$; explicit
formula (5.3.5) for the present case has the form
$$<x,y>_M=\Xi(xy_{r-1}^q+x^qy_1+x^{q^2}y_2+ ... +x^{q^{r-1}}y_{r-1})$$
The same consideration of the Newton polygon of the above polynomials
shows that for $x\in M_{T,0}$, $y\in M'_{T,0}$ we have
$\ord_\goth p x$, $\ord_\goth p y_i \ge 1/(q-1)$. Since $\ord_\goth p \Xi = -1/(q-1)$ we get
that $\ord_\goth p (<x,y>_M)>0$ and hence (because $<x,y>_M\in \n F_q$) we have
$<x,y>_M=0$. Dimensions of $M_{T,0}$, $M'_{T,0}$ are
complementary, hence they are mutually dual. $\square$
\medskip
{\bf Remark 13.4.2.} Analogous explicit proof exists for any standard-3 $M$ of
Section 11.8.
\medskip
\medskip
{\bf References}

\nopagebreak
\medskip
[A] Anderson Greg W. $t$-motives. Duke Math. J. Volume 53, Number 2 (1986), 457-502.

[BH] Matthias Bornhofen, Urs Hartl. Pure Anderson Motives and Abelian $\tau$-Sheaves. arXiv:0709.2809

[F] Faltings, Gerd. Group schemes with strict $\Cal O$-action. Mosc.
Math. J.  2
(2002),  no. 2, 249--279.

[G] Goss, David. Basic structures of function
field arithmetic.

[GL17] Grishkov A., Logachev, D. Lattice map for Anderson t-motives: first approach.
J. Number Theory 180 (2017), 373–402. http://arxiv.org/pdf/1109.0679.pdf

[GL18] Grishkov A., Logachev, D. Duality of Anderson t-motives having $N\ne0$.
https://arxiv.org/pdf/1812.11576.pdf

[H] Urs Hartl. Uniformizing the Stacks of Abelian Sheaves.

http://arxiv.org/abs/math.NT/0409341

[L09] Logachev, D. Anderson t-motives are analogs of abelian varieties with multiplication by imaginary quadratic fields.
http://arxiv.org 0907.4712.

[L] Logachev, D. Reductions of Hecke correspondences on Anderson varieties. In preparation.

[P] Richard Pink, Hodge structures over function fields. Universit\"at Mannheim.
Preprint. September 17, 1997.

[Sh63] Shimura, Goro. On analytic families of polarized abelian
varieties and automorphic functions. Annals of Math., 1 (1963), vol.
78, p. 149 -- 192

[Sh71]  Shimura, Goro. Introduction to the
arithmetic theory of
automorphic functions.

[Sh98]  Shimura, Goro. Abelian varieties with
complex multiplication
and modular functions. Princeton Mathematical Series, 46.

[Tae] Taelman Lenny, Artin t-motifs. J. Number Theory, 129 (2009), 142 - 157

[T] Taguchi, Yuichiro. A duality for finite $t$-modules.
J. Math. Sci. Univ. Tokyo 2 (1995), no. 3, 563--588.

\enddocument

and $\mu$ a number such that the $\mu$-dual ${M'}^{\mu}$ exists. Let $m$ be the minimal number such that $N^m=0$.
\medskip
{\bf Conjecture 6B1.} 1. $\mu \ge m$, i.e. $N^\mu=0$.
\medskip
2. ${M'}^{\mu}$ is uniformizable. There exists a canonical perfect $\w$-valued pairing between $L_T(M)$ and $L_T({M'}^{\mu})$.
\medskip
This pairing extends by $\p[[T-\theta]]$-linearity to the $\p[[T-\theta]]$-valued pairing between $L_T(M)\underset{\w}\to{\otimes}\p[[T-\theta]]$, $L_T({M'}^{\mu})\underset{\w}\to{\otimes}\p[[T-\theta]]$. We denote it by $<.,.>$.
\medskip
{\bf Conjecture 6B2.} We have:
$$\goth q({M'}^{\mu})=\{x\in L_T({M'}^{\mu})\underset{\w}\to{\otimes}\p[[T-\theta]]  \hbox{ such that }  \forall y\in \goth q(M) $$ $$\hbox{ we have } <x,y>\in (T-\theta)^\mu\p[[T-\theta]]\} $$

Let us consider uniformizable $M$ having $N=0$ (i.e. $m=1$) such that $M'$ --- the dual of $M$ --- exists and has $N'=0$.
\medskip
{\bf Corollary 6B3.} For this case Conjecture 6B2 is true, i.e.
$$\goth q({M'})=\{x\in L_T({M'})\underset{\w}\to{\otimes}\p[[T-\theta]]  \hbox{ such that }  \forall y\in \goth q(M) $$ $$\hbox{ we have } <x,y>\in (T-\theta)\p[[T-\theta]]\} $$
\medskip
This follows immediately from Theorem 5 and equivalence of Definition 2.3 and Property 2.4 (we consider the reduction of the above pairing $<.,.>$ modulo the maximal ideal $ (T-\theta) \p[[T-\theta]]$ of $\p[[T-\theta]]$; this reduction coincides with the pairing of Section 5).

\newpage
{\bf 1.16.4. $\Cal P$-rank of $M$ and of $M'$.} By analogy with the number field case, the $\Cal P$-rank of $M$ is the dimension of the $\Cal P$-torsion points of $E(M)$ over $\w/\Cal P$; it varies from $r-n$ (ordinary $M$) to 0 (completely suresingular $M$).
\medskip
{\bf Conjecture 1.16.5.} $M$ is ordinary $\iff M'$ is ordinary.
\medskip
For standard-3 t-motives  apparently this can be shown by explicit calculations. See 13.4.1 for the proof for the case of Drinfeld modules.
\medskip
{\bf Question 1.16.6.} What are possible values of $\Cal P$-rank of $M$, $\Cal P$-rank of $M'$? Is it true that the pair of numbers ($\Cal P$-rank of $M$, $\Cal P$-rank of $M'$) characterizes completely (in some meaning) the type of $M$? Particularly, whether the $\Cal P$-ranks of other tensor operations of $M$ (exterior powers etc.) are defined completely by ($\Cal P$-rank of $M$, $\Cal P$-rank of $M'$), or not?
\medskip
{\bf Example 1.16.7.} Let us consider $M$ defined by (8.2.2), entries of $A$ belong to $\bar \n F_q$, and let $\Cal P=T$. According Lang's Theorem, $T$-rank of $M$ is $n \iff \det A\ne 0$. An explicit calculation for the case $n=2$ shows that if the entries of $A$ belong to $\n F_q$ then the $T$-rank of $M$ is equal to the quantity of the non-zero eigenvalues of $A$ (if the entries of $A$ do not belong to $\n F_q$ then the formula for the $T$-rank of $M$ is more complicated). The dual $M'$ is defined by the same (8.2.2) with $A$ replaced by $-A^t$; this means that  if the entries of $A$ belong to $\n F_{q}$ then the $T$-rank of $M$ is equal to the $T$-rank of $M'$. It is easy to check that the same is true if the entries of $A$ belong to $\n F_{q^2}$ but not to a larger field:  for example, if $\gamma\in \bar \n F_q-\n F_{q^2} $ then for $A=\left(\matrix \gamma & 1 \\ -\gamma^{q+1}& -\gamma^{q}\endmatrix \right)$ we have $T$-rank of $M$ is 0, $T$-rank of $M'$ is 1.
\medskip
{\bf Example 1.16.8.} For pure non-ordinary standard-3 $M$ having $r=5$, $n=2$ we get easily that the case $T$-rank of $M$ is 2, $T$-rank of $M'$ is 0 cannot be realized; all other possible cases can be realised.
\medskip
{\bf Remark 13.4.3: Case $N\ne0$.} I do not know a definition of ordinarity for $M$ in finite characteristic for the case $N\ne0$; we can expect that these are those $M$ whose $\Cal P$-rank is the maximal possible in some family of these $M$. A version of Conjecture 13.4.1 should hold for this case; instead of the dual $M'$ we should consider the $m$-th dual ${M'}^m$ where $m$ satisfies $M^m=0$ (or $m$ is the minimal number with this property?)
\medskip
{\bf Example 13.4.3.1.} Let $M$ be given by the following modified formula (8.2.2):
$$Te_* =( \theta+N) e_* + A \tau e_* + \tau^2 e_*\eqno{(13.4.3.2)}$$
where $n=2$, $N=\left(\matrix 0&1\\0&0 \endmatrix\right)$, $A= \left(\matrix a_{11}&a_{12} \\a_{21}&a_{22} \endmatrix\right)$. For the case of finite characteristic (i.e. $a_{ij}\in \bar\n F_q$) and $\Cal P=T$ we have $M$ is ordinary iff $a_{21}\ne0$, in this case the $\Cal P$-rank of $M$ is 3. We have $m=2$, ${M'}^2$ is given by the formula
$$Te'_*=( \theta+N') e'_* + A' \tau e'_*$$
where $N'=\left(\matrix 0&-1&0&0&0&0\\0&0&0&0&0&0\\0&0&0&1&0&0\\0&0&0&0&0&0
\\0&0&0&0&0&1\\0&0&0&0&0&0\endmatrix\right)$,  $A'=\left(\matrix 0&0&0&0&1&0 \\0&0&1&0&0&0\\0&0&0&0&0&0\\0&1&-a_{11}&0&-a_{21}&0\\0&0&0&0&0&0
\\1&0&-a_{12}&0&-a_{22}&0 \endmatrix\right)$.

An explicit calculation shows that the $\Cal P$-rank of $M$ is $\le 1$ and it is 1 iff $a_{21}\ne0$. So, an analog of 13.4.1 for $M$ is the following:
\medskip
{\bf Conjecture 13.4.3.3.} Let $M$ be given by 13.4.3.2 in the generic characteristic, and let $M$ have good ordinary reduction at $\Cal P=T$.  We have: $M_{T,0}$ and $({M'}^2)_{T,0}$ are mutually dual with respect to the pairing.
\medskip
This conjecture can be easily checked by explicit calculation.

\enddocument

More exactly, if $\vf: M \to N$ is a map of abelian
t-motives and both ${M_{\goth C}'}^{\mu}$, ${N_{\goth C}'}^{\mu}$ exist then there exists
the map of rational pr\'e-t-motives ${\vf_{\goth C}'}^{\mu}: {N_{\goth C}'}^{\mu} \to
{M_{\goth C}'}^{\mu}$.


$\sigma: \p \to \p$ the Frobenius automorphism of $\p$, i.e. $\sigma(z)=z^q$.

We extend $\sigma$ to $K_C$, resp. $\w_C$ by the formula
$\sigma_{K}(k\otimes z)=k\otimes \sigma(z)$, $z\in\p$, $k\in K$, resp. $k\in
\w$.

: we consider $M\mapsto M^{(1)}$ as a
functor; if $M$ is free $\p[T]$-module and $\alpha: M \to M$ is a map whose matrix in
some basis $\{m_*\}$ is $C$, then the matrix of $\alpha^{(1)}: M^{(1)} \to M^{(1)}$ in the
base $\{m_*\otimes1\}$ is $C^{(1)}$

We have an inclusion $\w \hookrightarrow K_\infty$,
we prolonge $\iota$ to the inclusion $K_\infty \hookrightarrow \p$, and we denote
the image $\iota(K_\infty)\subset \p$ by $K_\infty$ as well.
zarhin@math.psu.edu
Dorogoj Yura,

mozhet byt', Vy znaete chto-nibud' po takim voprosam:

1. Rassmotrim ob'ekty: $C$-prostranstvo $V$, $Z$-reshyotka $L$ v nyom i ermitova forma $H$ na $V$, kotoraya ne >>0 i takaya, chto $im H$ celochislenna na $L$. Kakov kriterij togo, chto eti ob'ekty yavlyayutsya realizaciej motiva?

My znaem, chto esli $H >> 0$, to vsegda. Esli net, to po krajnej mere inogda, naprimer esli $X$ - trifold takoj, chto ego $J^3(X)$ - ne abelevo mnogoobrazie, to $H^3(X)$ - kak raz takogo tipa.

Otsyuda vytekaet takoj vopros

2. My znaem, chto abelevo mnogoobrazie - bolee obshchij ob'ekt, chem jacobian krivoj. Rassmotrim trifold $X$ takoj, chto ego $J^3(X)$ - ne abelevo mnogoobrazie, rassmotrim ego motiv $M=H^3(X)$ i sprosim: motivy kakogo tipa obobshchayut $M$ tak zhe, kak abelevo mnogoobrazie obobshchaet jacobian krivoj?

Mozhet byt', lyudi dumali nad etim 30 - 40 let nazad i uvideli, chto na eti voprosy net razumnyh otvetov?

Prois|hozhdenie etih voprosov takoe. God nazad ya poslal Vam stat'yu ob analogii mezhdu t-motivami s $N=0$ i abelevymi mnogoobraziyami s umnozheniem na mnimoe kvadratichnoe pole. Vopros: a kakov analog t-motivov s $n \ne 0$? Mozhet byt', eto motiv? My legko mozhem opredelit' vysheupomyanutye $V$, $L$ i $H$ dlya vneshnej stepeni abelevyh mnogoobrazij s umnozheniem na mnimoe kvadratichnoe pole. Esli by my znali, chto oni - realizaciya motiva, to kak raz eto i byl by motiv, sootvetstvuyushchij t-motivu s $n \ne 0$.

Yours, Dima

 We consider the case $P\ne \infty$, i.e. we can identify $P$ with an irreducible polynomial in $\n Z_\infty$ (up to a scalar factor).

C-MATRICY                 POZOR !

{\bf 8.2. Infinitesimal lattice map and C-lattices.} For $r=2n$, \ $n>1$ the degree of the lattice map is infinite, hence in order to get a 1 - 1 correspondence it is necessary to introduce a new object called a C-lattice. There is a C-lattice map from a neighbourhood of a t-motive $M_0$ to a neighbourhood of a C-lattice $L_{C,0}$ (see below for the definition of $M_0$, $L_{C,0}$) which is 1 - 1.

For the case $(r,n)=1$ I do not see an inicial t-motive $M_0$ such that in its neighbourhood the degree of the lattice map is not 1 (see Remark 8.2.14 for more details), hence we can ask whether in this case the answer to 8.1.2 is yes. I think that this is few likely.
\medskip
{\bf Definition 8.2.1.} A C-lattice is a pair $(L, e_*)$ where
\medskip
$L \subset V=\p^n$ is a lattice and
\medskip
$e_1, ... , e_r$ is a $\n Z_\infty$-basis of $L$.
\medskip
Further, two such pairs $(L; e_1, ... , e_r)$ and $(L'; e'_1, ... ,
e'_r)$ are called equivalent if there exists a $\p$-linear map $\psi:
V\to V$ such that $\psi(L; e_1, ... , e_r) = (L'; e'_1, ... , e'_r)$,
or $L=L'$ and the matrix of the change of basis from $e_1, ... , e_r$
to $e'_1, ... , e'_r$ (which a priori belongs to $GL_r(\n Z_\infty)$ ) belongs to $GL_r(\n F_q)$.
\medskip
The functor of forgetting the basis from C-lattices to lattices is
denoted by $\goth i$. The notion of a Siegel matrix for a C-lattice is
the same as for a lattice.

We consider infinitesimal degree of the lattice map in a neighbourhood
of some distinguished (having many endomorphisms) t-motive $M_0$. We
consider the case $M_0=\goth C_2^{\oplus n}$ (recall that $\goth C_2$ is the Carlitz module over $\n F_{q^2}[T]$; in 7.7 $M_0$ is denoted by $M$).
\medskip
{\bf Remark.} It is easy to see that $M_0$ is the (only) t-motive with complete multiplication by $\n F_{q^{2n}}[\theta]$ with CM-type $Id, \fr^2, \fr^4, ..., \fr^{2n-2}$, see 12.5.3 and 13.3, first case.
\medskip
Let $\omega \in \n F_{q^2} - \n F_q$ be a
fixed element. A Siegel matrix of $M_0$ is $\omega I_n$. We denote the lattice (resp. the C-lattice) corresponding to $\omega I_n$ by $L_0$ (resp. $L_{C,0}$). $L_0$ is the lattice of $M_0$. We consider 5 sets $S_1, ... , S_5$:
\medskip
$S_1$. The set of $n \times n$ matrices $A$.
\medskip
$S_2$. The set of t-motives $M$ given by the equation (see 1.9.1)

$$Te_* = \theta e_* + A \tau e_* + \tau^2 e_*\eqno{(8.2.2)}$$ where $A
\in S_1$, $e_* = (e_1,...,e_n)^t$ (notations of 1.9.1: $\goth A_1=A$, $\goth A_2=I_n$).
\medskip
$S_3$. The set of $n \times n$ Siegel matrices $Z$.
\medskip
$S_4$. The set of C-lattices of rank $r=2n$ in $C_\infty ^n$.
\medskip
$S_5$. The set of lattices of rank $r=2n$ in $C_\infty ^n$.
\medskip
We consider initial elements 0, $M_0$, $\omega I_n$, $L_{C,0}$, $L_0$ of $S_1, ...,S_5$ respectively and open neighbourhoods $U_i\subset S_i$ of these initial elements.
\medskip
{\bf Proposition 8.2.3.} There exist neighbourhoods $U_2$, $U_4$, $U_5$ such that

(a) The restriction of $\goth i$ to $U_4$ gives us an epimorphism $U_4\to U_5$.

(b) there exists a 1 -- 1 map $\mu_{24}$ from $U_2$ to $U_4$ such that $\mu_{25}:=\goth i\circ\mu_{24}$ (see the below diagram 8.2.4) is the lattice map from uniformizable t-motives to lattices. Particularly, for $n>1$ the fibre of $\mu_{25}$ is discrete infinite.
\medskip
{\bf Proof.} It is sufficient to prove that there exists a commutative diagram $$\matrix U_1 &
\overset{\mu_{12}}\to{\to} & U_2 && \\ \\ \mu_{13} \downarrow &
&\mu_{24}\downarrow &
\overset{\mu_{25}}\to{\searrow} &\\ \\ U_3 &
\overset{\mu_{34}}\to{\to} & U_4 & \overset{\goth i}\to{\to} & U_5
\endmatrix \eqno{(8.2.4)}$$ where $\mu_{12}(A)$ is the t-motive defined by
8.2.2, $\mu_{34}(Z)$ is the C-lattice corresponding to a Siegel matrix
$Z$ and $\mu_{13}$ is defined as
follows. We identify $\Lie(\goth C_2)$ with $\p$ and hence $\Lie(M_0)$ with $\p^n$.
We consider the following basis
$l_{0,1}, ... , l_{0,2n}$ of $L_0\subset \Lie(M_0)=\p^n$: $l_{0,i}=(0,...,0,1,0,...,0)$ (1 at the $i$-th place), $l_{0,n+i}=\omega l_{0,i}$, $i=1,...,n$. The Siegel
matrix of this basis is $\omega I_n$. Let $A\in U_1$, $M$ its $\mu_{12}$-image and $L$ its
lattice, i.e. its $\mu_{25}\circ\mu_{12}$-image. For any $l_0\in L_0$
there exists a well-defined $l\in L$ which is close to $l_0$ (because
entries of $A$ are near 0). So, we consider a basis $l_{1}, ... ,
l_{2n}$ of $L$ where any $l_i$ is near the corresponding $l_{0,i}$.
The Siegel matrix corresponding to $l_{1}, ... , l_{2n}$ is exactly
the $\mu_{13}$-image of $M$. By definitions, the outer quadrangle is
commutative.

We denote by $d_{\alpha\beta}$ the degree of
$\mu_{\alpha\beta}$ at a
generic point near the initial element of $U_\alpha$. We must prove that there exists a map $\mu_{24}$ preserving
commutativity, and that $d_{24}=1$.
\medskip
{\bf Lemma 8.2.5.} $d_{13}=1$.
\medskip
{\bf Proof.} We consider maps of the diagram 5.1.3 for $M$ given by 8.2.2. For $Z\in E=\p^n$ we have $m_T(Z)=\theta Z + A Z^{(1)} +
Z^{(2)}$, and for $Z\in \Lie(M)=\p^n$ we have $\Exp(Z)=\Exp_A(Z)=\sum_{i=0}^\infty C_iZ^{(i)}$ where
$C_i=C_i(A)$, $C_0=1$.

For the reader's convenience, we consider only the case $n=1$ (the
general case does not require any new ideas), hence $Z$, $A$ will be
denoted by $z$, $a$ respectively. We denote
$\theta_{ij}=\theta^{q^i}-\theta^{q^j}$. Recall that the exponent for $\goth C_2$ has the form
$$\Exp_0(z)=z+\frac{1}{\theta_{20}}z^{q^2}+\frac{1}{\theta_{42}\theta_{40}}z^{q^4}+...\eqno{(8.2.6)}$$
($C_{2i}(0)=\frac{1}{\prod_{j=0}^{i-1}\theta_{2i,2j}}$). We denote by
$y_0\in \p$ a nearest-to-zero root to $\Exp_0(z)=0$ (this is $\xi$ of $\goth C_2$ in notations of [G]). It is defined up to
multiplication by elements of $\n F_{q^2}^*$, and it generates over $\n F_{q^2}[\theta]$ the
lattice of the Carlitz module $\goth C_2$. We fix one such $y_0$.
If $a$ is sufficiently small then there is the only root to
$\Exp_a(z)=0$ near $y_{0}$, and there is the only root to
$\Exp_a(z)=0$ near $\omega y_{0}$. We denote these roots
by $z=z(a)$, $z'=z'(a)$ respectively, and we denote $z=y_{0}+\delta$,
$z'=\omega y_{0}+\delta'$. $\delta$ (resp. $\delta'$) is a root to the
power series
$$\sum_{i,j=0}^\infty d_{ij}a^i\delta^j=0 \eqno{(8.2.7)}$$
$$\hbox{resp. }\sum_{i,j=0}^\infty d'_{ij}a^i{\delta'}^j=0 \eqno{(8.2.7')}$$
where
$$d_{00}=0, \ \
d_{10}=\frac{y_{0}^{q}}{\theta_{10}}+\frac{y_{0}^{q^3}}{\theta_{31}\theta_{30}}
+ \frac{y_{0}^{q^5}}{\theta_{53}\theta_{51}\theta_{50}}+
\frac{y_{0}^{q^7}}{\theta_{75}\theta_{73}\theta_{71}\theta_{70}}+..., \ \ d_{01}=1$$
$$d'_{00}=0, \ \ d'_{10}=\frac{(\omega
y_{0})^{q}}{\theta_{10}}+\frac{(\omega
y_{0})^{q^3}}{\theta_{31}\theta_{30}} + \frac{(\omega
y_{0})^{q^5}}{\theta_{53}\theta_{51}\theta_{50}}+
\frac{(\omega y_{0})^{q^7}}{\theta_{75}\theta_{73}\theta_{71}\theta_{70}}+..., \ \ d'_{01}=1$$
Moreover, $\delta$, $\delta'$ are the nearest-to-0 roots to (8.2.7),
$(8.2.7')$ respectively. It is clear that $d_{10}, d'_{10}\ne 0$. This
means that the approximate value of $\delta$, $\delta'$ is
$$-d_{10}a, \ \ -d'_{10}a\eqno{(8.2.8)}$$
respectively. Exactly, both $\delta$, $\delta'$ are power series in
$a$ whose first term is given by (8.2.8). This means that
$$z=y_0-d_{10}a+\sum_{i=2}^\infty k_i
a^i$$
$$z'=\omega y_0-d'_{10}a+\sum_{i=2}^\infty k'_i
a^i$$ It is easy to see that $d'_{10}\ne \omega d_{10}$,
hence the Siegel matrix $\goth Z=z^{-1}z'$ (which is a number because $n=1$) is given by the formula $$\goth Z=\omega
+ \sum_{i=1}^\infty l_i a^i\eqno{(8.2.9)}$$ and $l_1\ne 0$. Since (for $n=1$)
$d_{13}$ is the minimal $i$ such that $l_i\ne 0$ we get that
$d_{13}=1$.
\medskip
For $n > 1$ the calculation is the same (this is the main part of the proof, because for the case $n=1$ the result is known). Analog of (8.2.9) is $$\goth Z=\omega I_n + l_1 A + P_{\ge 2}(A)\eqno{(8.2.10)}$$ where $P_{\ge 2}(A)$ is a power series of entries of $A$ such that all its terms have degree $\ge 2$.
Condition $l_1\ne 0$ implies $d_{13}=1$. $\square$
\medskip
In order to simplify notations, we identify until the end of the proof $L_T(M)$ and $L(M)$ via $\goth a$, and $\n Z_\infty$ and $\w$ via $\iota$. We denote the monodromy group of $\mu_{12}$, resp. $\mu_{35}:=\goth i \circ \mu_{34}$ by $\Cal M_{12}$, resp. $\Cal M_{35}$. We have $\Cal M_{12}=GL_n(\n F_{q^2})/\n F_q^*$. Really, the automorphism group of $M_0$ is $GL_n(\n F_{q^2}[T])$: an element $g$ of this group acts on the basis $e_*$ of (8.2.2) if $A=0$. But if $A$ is a generic matrix, then for $g\in GL_n(\n F_{q^2}[T]) - GL_n(\n F_{q^2})$ the result of the action of $g$ on (8.2.2) becomes a more complicated equation: terms having higher powers of $\tau$ appear, so these $g$ do not belong to $\Cal M_{12}$. Factorization by $\n F_q^*$ is obvious. Further, obviously $\Cal M_{35}=\{\gamma\in PGL_{2n}(\n Z_\infty)|\gamma(\omega I_n)=\omega I_n\}$. The outer quadrangle of 8.2.4 defines a map from $\Cal M_{12}$ to $\Cal M_{35}$ which we denote by $\alpha$.
\medskip
{\bf Lemma 8.2.11.} $\alpha$ is injective, and $\im \alpha= \{\gamma\in PGL_{2n}(\n F_{q})|\gamma(\omega I_n)=\omega I_n\}=\Cal M_{35} \cap PGL_{2n}(\n F_{q}) \subset PGL_{2n}(\n Z_\infty)$.
\medskip
{\bf Proof.} $\alpha$ is defined by the condition: for any $A\in S_1$, $\gamma \in \Cal M_{12}$ we have
$$\mu_{13}(\gamma(A))=(\alpha(\gamma))(\mu_{13}(A))\eqno{(8.2.12)}$$ The explicit formula for $\alpha$ is the following. For odd $q$ we fix $\omega$ satisfying $\omega^2\in \n F_q^*$ (it is an easy exercise to find analog of the below formula for even $q$). Let $\gamma=U+\omega V$, \ \ $U,V\in GL_n(\n F_q)$ (we consider cleary a representative of $\gamma\in GL_n(\n F_{q^2})/\n F_q^*$ in $GL_n(\n F_{q^2})$). Then $$\alpha(\gamma)=\left(\matrix U&-\omega^2V\\ -V&U \endmatrix\right)^{-1}\eqno{(8.2.13)}$$ ($\alpha$ is an antihomomorphism, because the functor of lattice is contravariant). It is checked immediately that (8.2.12) holds. We see that $\im \alpha= \{\gamma\in PGL_{2n}(\n F_{q})|\gamma(\omega I_n)=\omega I_n\}$. $\square$ 
\medskip
Proposition 8.2.3 follows immediately from these lemmas. $\square$
\medskip
{\bf Remark 8.2.14.} We see that if $r=2n$, $n>1$ then the lattice map $\mu_{25}$ is not a local isomorphism near $M_0$. The origin of this phenomenon is reducibility of $M_0$ which implies that the monodromy group of $\mu_{35}$ is much bigger then the one of $\mu_{12}$. For other values of $r$, $n$ a natural analog of $M_0$ is a t-motive with complete multiplication. Apparently if $(r,n)=1$ then all pure t-motives with complete multiplication are irreducible (example: CM-field is $\n F_{q^r}[\theta]$), and analog of Proposition 8.2.3 for this case shows that $\mu_{25}$ is a local isomorphism near this t-motive.
\medskip
{\bf Remark 8.2.15.} From the first sight, for $n=1$ Proposition 8.2.3 contradicts to a result of Drinfeld about 1 -- 1 correspondence between Drinfeld modules and lattices. Really, there is no contradiction: if $n=1$ then $\im \alpha=\Cal M_{35}$, and --- although $\goth i: S_4 \to S_5$ clearly is not an isomorphism --- its restriction to $U_4$ is an isomorphism $U_4 \to U_5$.
\medskip
{\bf 8.3. Duality of C-lattices.} We have no analog of 2.2 for C-lattices, so we use an analog of 3.2 as a definition of duality. Namely, if $Z$ is a Siegel matrix of a C-lattice $(L, e_*)$ then its dual $(L, e_*)'$ is a C-lattice whose Siegel matrix is $-Z^t$ (or $Z^t$ which is the same, but more convenient for further calculations, because $(-\omega I_n)^t\ne \omega I_n$). Equality 3.8.2 shows that this notion is well-defined (entries of $A$, $B$, $C$, $D$ belong to $\n F_q$).

An analog of Theorem 5 holds for C-lattices:
\medskip
{\bf Theorem 8.3.1.} Let $M\in U_2$. Then $\mu_{24}(M')=\mu_{24}(M)'$.
\medskip
{\bf Proof} is completely analogous to the proof of Theorem 5, so we omit it. Alternatively, we can show that the exact form of (8.2.10) is $$\goth Z=\omega I_n + \sum_{k=1}^\infty \sum_{d_1,...,d_k}l_{d_1,...,d_k}A^{(d_1)}\cdot...\cdot A^{(d_k)}$$ where coefficients $l_{d_1,...,d_k}$ satisfy $$l_{d_1,...,d_k}=l_{d_k,...,d_1}$$ This obviously implies the theorem. $\square$
\medskip

\newpage
{\bf Theorem.} ${M'}^{m}$ is uniformizable, and $\underline{H}({M'}^{m})= {\underline{H}(M)'}^{m}$.
\medskip
We shall need several elementary lemmas.
Let $w$ be a number such that $N^w=0$. We define numbers $k_i=k_i(M)$, $i=2, \dots,w+1$, as follows:
$$k_{i}:=\dim \Ker N^{i-1}/ \Ker N^{i-2}- \dim \Ker N^{i}/ \Ker N^{i-1} $$
$$=\dim \im N^{i-2} / \im N^{i-1} - \dim \im N^{i-1} / \im N^{i}\eqno{(9.3)}$$
Equivalently, let $n=d_1+...+d_\al$, where $d_1\ge d_2\ge...\ge d_\al>0$, be a partition of $n$ corresponding to the Jordan form of $N$, i.e. the Jordan form of $N$ consists of $\al$ 0-Jordan blocks of sizes $d_1, d_2,...,d_\al$. We have $w\ge d_1$, $\al\le r$. We shall call a 0-partition of length $\g r$ of a number $\g n$ a representation of $\g n$ as a sum
$$\g n = \g d_1+\g d_2+...+\g d_{\g r}$$ where $\g d_i\in \n Z$ and $\g d_1\ge\g d_2\ge...\ge\g d_{\g r}\ge0$, i.e. a 0-partition is a partition plus several zeroes at its end. We extend the partition $n=d_1+...+d_\al$ to a 0-partition of length $r$ denoted by $\g p=\g p(M)$:
$n=d_1+...+d_\al+d_{\al+1}+...+d_r$ where $d_{\al+1}=...=d_r=0$. Let $n=c_1+...+c_{d_1}+c_{d_1+1}+...+c_w$ be the 0-partition of length $w$ dual to $\g p$ (the definition of the dual 0-partition of a given length is clear). We have $\al=c_1\ge c_2\ge...\ge c_{d_1}>0$, $c_{d_1+1}=...=c_w=0$. We have $\dim \Ker N^{i}=c_1+...+c_i$, hence $k_i=c_{i-1}-c_i\ge0$ (for $i=w+1$ we let $c_{w+1}=0$), and
$$n=\sum_{i=1}^{w}ik_{i+1}\eqno{(9.4)}$$ We have $\al=c_1=\sum_{i=2}^{w+1}k_{i}$, hence $r\ge \sum_{i=2}^{w+1}k_{i}$. We let $k_1:=r- \sum_{i=2}^{w+1}k_{i}$.
\medskip
According ??, elements $T^il_j$, $i=0,1,\dots$, $j=1,\dots,r$, generate $\Lie(M)$ as a $\p$-vector space. Hence, elements $N^il_j$, $i=0,\dots,w-1$, $j=1,\dots,r$, also generate $\Lie(M)$ as a $\p$-vector space. Using this fact we arrange elements $l_1,\dots,l_r$ in $w+1$ segments as follows. First, elements $N^{w-1}l_j$, $j=1,\dots,r$, generate $N^{w-1}\Lie(M)$ as a $\p$-vector space. Its dimension is $k_{w+1}$, hence (first step of a process) we can choose $k_{w+1}$ elements from $l_1,\dots,l_r$ (we denote them by $l_{w+1,1}, \dots, l_{w+1,k_{w+1}}$ respectively) such that
\medskip
(9.6) $N^{w-1}(l_{w+1,i})$, $i=1, \dots, k_{w+1}$, form a $\p$-basis of $N^{w-1}\Lie(M)$.
\medskip
Further (second step), elements $N^{w-2}l_j$, $N^{w-1}l_j$, $j=1,\dots,r$, generate $N^{w-2}\Lie(M)$ as a $\p$-vector space. Elements $N^{w-2}(l_{w+1,i})$, $N^{w-1}(l_{w+1,i})$, $i=1, \dots, k_{w+1}$, are linearly independent over $\p$. Really, let $$\sum_{\al_2=1}^{k_{w+1}}c_{\al_2}N^{w-2}(l_{w+1,\al_2})+ \sum_{\al_1=1}^{k_{w+1}}c_{\al_1}N^{w-1}(l_{w+1,\al_1})=0\eqno{(9.7)}$$
be a non-trivial dependence relation. Applying $N$ to (9.7) we get $$\sum_{\al_2=1}^{k_{w+1}}c_{\al_2}N^{w-1}(l_{w+1,\al_2})=0$$ that contradicts (9.6). Hence, all $c_{\al_2}$ are 0 that again contradicts (9.6).
\medskip
We have $\dim N^{w-1}\Lie(M)=2k_{w+1}+k_{w}$, this follows immediately from (9.3). Hence, we get:
\medskip
Among $l_1,\dots,l_r$ there exist $k_w$ elements (we denote them by $l_{w,1}, \dots, l_{w,k_{w}}$ respectively) such that their intersection with $l_{w+1,1}, \dots, l_{w+1,k_{w+1}}$ is empty and such that
\medskip
(9.8) $N^{w-2}(l_{w,i})$, $i=1, \dots, k_{w}$, $N^{w-2}(l_{w+1,i})$, $i=1, \dots, k_{w+1}$, $N^{w-1}(l_{w+1,i})$, $i=1, \dots, k_{w+1}$, form a $\p$-basis of $N^{w-2}\Lie(M)$.
\medskip
Third step of the process: elements $N^{w-3}l_j$, $N^{w-2}l_j$, $N^{w-1}l_j$, $j=1,\dots,r$, generate $N^{w-3}\Lie(M)$ as a $\p$-vector space. Elements $N^{w-3}(l_{w,i})$, $N^{w-2}(l_{w,i})$, $i=1, \dots, k_{w}$, and $N^{w-3}(l_{w+1,i})$, $N^{w-2}(l_{w+1,i})$, $N^{w-1}(l_{w+1,i})$, $i=1, \dots, k_{w+1}$, are linearly independent over $\p$ (the proof is as the one above for the second step). Hence, we get:
\medskip
Among $l_1,\dots,l_r$ there exist $k_{w-1}$ elements (we denote them by $l_{w-1,1}, \dots,$ $ l_{w-1,k_{w-1}}$ respectively) such that their intersection with $l_{w,1}, \dots, l_{w,k_{w}}$, $l_{w+1,1}, \dots, l_{w+1,k_{w+1}}$ is empty and such that
\medskip
(9.9) $N^{w-3}(l_{w-1,i})$, $i=1, \dots, k_{w-1}$, $N^{w-3}(l_{w,i})$, $N^{w-2}(l_{w,i})$, $i=1, \dots, k_{w}$, and $N^{w-3}(l_{w+1,i})$, $N^{w-2}(l_{w+1,i})$, $N^{w-1}(l_{w+1,i})$, $i=1, \dots, k_{w+1}$, form a $\p$-basis of $N^{w-3}\Lie(M)$.
\medskip
Continuing this process, we represent $r$ as an ordered partition
$$r=k_1+...+k_{w+1}\eqno{(9.10)}$$
(recall that some $k_*$ can be 0) and we represent the set $\{l_1,\dots,l_{r}\}$ as a union of segments $$\{l_1,\dots,l_{r}\}=\bigcup_{u=1}^{w+1}\{l_{u1}, \dots, l_{uk_u}\}\eqno{(9.11)}$$ (the union is ordered and disjoint) such that $\forall \ u=0,\dots,w-1$ we have:
\medskip
$(9.12)$ A $\p$-basis of $N^u\Lie(M)$ is formed by elements $N^\al(l_{\be\ga})$, where $\al\in[u,\dots, w-1]$, $\be\in[\al+2,\dots, w+1]$, $\ga\in [1,\dots, k_\be]$.
\medskip
This implies that for any
$$u\in [1,w], \ \ z\in [u-1, w-1], \ \  y\in [z+2,w+1], \ \ \eqno{(9.13)}$$
$$v=u-1\eqno{(9.14)}$$
there exist matrices $S_{uvyz}$ of size $k_u\times k_y$ with entries in $\p$ (analogues of the Siegel matrix for the case $w=1$) such that $\forall \ u=1,\dots, w, \forall \ i=1,\dots, k_u$ the following holds:
$$N^{u-1}l_{ui}=-\sum_{z=u-1}^{w-1}\sum_{y=z+2}^{w+1} \sum_{j=1}^{k_y}(S_{uvyz})_{ij}N^zl_{yj}\eqno{(9.15)}$$
(if some $k_*$ are 0 then the corresponding $S_{****}$ do not exist).
\medskip
To simplify formulas, below for any $\al$ we consider $\hat l_\al:=l_{\al*}$ as matrix columns. (9.15) becomes a matrix equality
$$N^{u-1}\hat l_{u}=-\sum_{z=u-1}^{w-1}\sum_{y=z+2}^{w+1}S_{uvyz}N^z\hat l_{y}\eqno{(9.16)}$$

{\bf Remark 9.17.} Since always $v=u-1$, really matrices $S_{uvyz}$ depend on 3 parameters $u,y,z$. Number $v$ indicates the exponent of $N$ in the left hand side of (9.15), by analogy with $z$ which indicates the exponent of $N$ in the right hand side of (9.15). This notation is convenient to define a symmetry between $A_*$ and $P_*$, see below.
\medskip
{\bf 9.18.} Example for $w=3$, $u=1$:
$$\matrix l_{1i}=-(\sum_{j=1}^{k_2}(S_{1020})_{ij}l_{2j} &+& \sum_{j=1}^{k_3}(S_{1030})_{ij}l_{3j} &+& \sum_{j=1}^{k_4}(S_{1040})_{ij}l_{4j} \\ \\
&+& \sum_{j=1}^{k_3}(S_{1031})_{ij}N(l_{3j}) &+& \sum_{j=1}^{k_4}(S_{1041})_{ij} N(l_{4j}) \\ \\ &&&+& \sum_{j=1}^{k_4}(S_{1042})_{ij} N^2(l_{4j}) \ )\endmatrix $$
(terms of a fixed column of this formula correspond to a fixed $y$ and different $z$ of (9.15), and terms of a fixed row of this formula correspond to a fixed $z$ and different $y$ of (9.15) ).
\medskip
Applying powers of $N$ to (9.16), for any $v\in [0,\dots, w-1]$, $u\in [1,\dots, w+1]$ we can represent $N^{v}(\hat l_{u})$ as a linear combination of $N^z(\hat l_{y})$ where for a fixed $v$ the numbers $z, \ y$ satisfy $$z\in [v,\dots,w-1], \ \ y\in [z+2,\dots,w+1], \ \ \eqno{(9.19)}$$ Namely, there exist polynomials in $S_{****}$ denoted by $P_{uvyz}$ such that (matrix notations)

$$N^{v}\hat l_{u}=-\sum_{z=v}^{w-1}\sum_{y=z+2}^{w+1}P_{uvyz}N^z\hat l_{y}\eqno{(9.20)}$$
Clearly for $v=u-1$ we have $P_{uvyz}=S_{uvyz}$.

{\bf 9.21.} The domain $v\ge u-1 \ \ \wedge \ \ \{z, \ y$ satisfy (9.19)\} is called the non-trivial domain of the definition of $P_{****}$.

For $v< u-1$ (trivial domain) we have:

$$ P_{u,v,y,z}=-1, \hbox{ resp. } P_{u,v,y,z}=0  \eqno{(9.22)}$$ for $y$, $z$ satisfying (9.19), $(y,z)=(u,v)$, resp. $(y,z)\ne(u,v)$.

\medskip
{\bf 9.23.} Example for $w=3$:

\medskip

$N^2\hat l_2=(S_{2131}S_{3242}-S_{2141})N^2\hat l_4$, i.e. $P_{2242}=-S_{2131}S_{3242}+S_{2141}$;

\medskip

$N^2\hat l_1=(-S_{1020}S_{2131}S_{3242}+S_{1020}S_{2141} +S_{1030}S_{3242}-S_{1040})N^2\hat l_4$, i.e.

\medskip

$P_{1242}=S_{1020}S_{2131}S_{3242}-S_{1020}S_{2141} -S_{1030}S_{3242}+S_{1040}$;

\medskip

$N\hat l_1=(S_{1020}S_{2131}-S_{1030})N\hat l_3+(S_{1020}S_{2141}-S_{1040})N\hat l_4+$

\medskip

$+(S_{1020}S_{2142}+S_{1031}S_{3242}-S_{1041})N^2\hat l_4$, i.e.

\medskip

$P_{1131}=-S_{1020}S_{2131}+S_{1030}$, $P_{1141}=-S_{1020}S_{2141}+S_{1040}$,

\medskip

$P_{1142}=-S_{1020}S_{2142}-S_{1031}S_{3242}+S_{1041}$.
\medskip
{\bf Remark 9.24.} Although we do not need this fact, let us give a formula for some $P_{****}$. Let us define a block unitriangular matrix $\g S$ whose $(i,j)$-th block is $S_{i,i-1,j,i-1}$ for $j>i$, $I_{k_i}$ for $i=j$ and 0 for $j<i$. Further, we define a block unitriangular matrix $\g P$ whose $(i,j)$-th block is $-P_{i,j-2,j,j-2}$ for $j>i$, $I_{k_i}$ for $i=j$ and 0 for $j<i$. We have $\g P=\g S^{-1}$ (see the the below propositions).
\medskip
Some $P_{****}$ that enter in the below formula for $\bar \g D$ are not of the form of the elements of the inverse unitriangular matrix, for example $P_{1142}$, $w=3$.
\medskip
For the proof of Lemmas 9.32, 23, we need
\medskip
{\bf Lemma 9.25.} For all $i,\ j,\ \psi,\ \xi$ (domain?)
$$(\sum_{\be=0}^{j+\xi-w} \sum_{\al=i+\be}^{w+1-j+\be} S_{i-1,i-2,\al,i-2+\be} P_{\al,w-j+\be,\psi,\xi})-$$ $$-S_{i-1,i-2,\psi,i-2+\xi+j-w}+P_{i-1,w-j,\psi,\xi}=0\eqno{(9.25.1)}$$
(A recurrent formula for $P_{****}$).
\medskip
{\bf Proof.} First, we rewrite (9.16) for $u=i-1$:

$$N^{i-2}\hat l_{i-1}=-\sum_{z=i-2}^{w-1}\sum_{y=z+2}^{w+1}S_{i-1,i-2,y,z}N^z\hat l_{y}\eqno{(9.25.2)}$$

Now, for any $$j=1,\dots,w-i+2\eqno{(9.25.3)}$$ we apply $N^{w-i+2-j}$ to (9.25.2):

$$N^{w-j}\hat l_{i-1}=-\sum_{z=i-2}^{i+j-3}\sum_{y=z+2}^{w+1}S_{i-1,i-2,y,z}N^{z+w-i+2-j}\hat l_{y}\eqno{(9.25.4)}$$

(since $N^w=0$, we get that $z\le i+j-3$ ).

We change a summation variable: $z=i-2+\be$, and $y \to \al$, we get

$$N^{w-j}\hat l_{i-1}=-\sum_{\be=0}^{j-1}\sum_{\al=i+\be}^{w+1}S_{i-1,i-2,\al,i-2+\be}N^{w+\be-j}\hat l_{\al}\eqno{(9.25.5)}$$

Now we use (9.20), we make the following variable change:

$$u \to \al\ \ \ \ \ y\to \psi$$
$$v\to w-j+\be\ \ \ \ \ z\to \xi$$
we get
$$N^{w-j+\be}\hat l_\al=-\sum_{\xi=w-j+\be}^{w-1} \sum_{\psi=\xi+2}^{w+1} P_{\al,w-j+\be,\psi,\xi}N^{\xi}\hat l_\psi\eqno{(9.20a)}$$

We substitute (9.20a) in (9.25.5):

$$N^{w-j}\hat l_{i-1}=\sum_{\be=0}^{j-1}\sum_{\al=i+\be}^{w+1} \sum_{\xi=w-j+\be}^{w-1}\sum_{\psi=\xi+2}^{w+1} S_{i-1,i-2,\al,i-2+\be} P_{\al,w-j+\be,\psi,\xi}N^{\xi}\hat l_\psi\eqno{(9.25.6)}$$

We change the order of summation in (9.25.6):

$$N^{w-j}\hat l_{i-1}=\sum_{\xi=w-j}^{w-1} \sum_{\psi=\xi+2}^{w+1} (\sum_{\be=0}^{j+\xi-w} \sum_{\al=i+\be}^{w+1} S_{i-1,i-2,\al,i-2+\be} P_{\al,w-j+\be,\psi,\xi}) N^{\xi}\hat l_\psi\eqno{(9.25.7)}$$

We rewrite (9.20) making changes:

$$u\to i-1\ \ \ \ \ \ y \to \psi$$
$$v\to w-j\ \ \ \ \ \  z\to \xi$$
we get
$$N^{w-j}\hat l_{i-1}=-\sum_{\xi=w-j}^{w-1}\sum_{\psi=\xi+2}^{w+1} P_{i-1,w-j,\psi,\xi}N^\xi \hat l_{\psi}\eqno{(9.25.8)}$$

For $\psi\ge \xi+2$ elements $N^\xi l_{\psi i}$, $i=1,\dots,k_\psi$, are linearly independent over $\p$. Hence, (9.25.7), (9.25.8) imply

$$P_{i-1,w-j,\psi,\xi}=-\sum_{\be=0}^{j+\xi-w} \sum_{\al=i+\be}^{w+1} S_{i-1,i-2,\al,i-2+\be} P_{\al,w-j+\be,\psi,\xi}\eqno{(9.25.9)}$$

Here the domain of $\xi$, $\psi$ is:

$$\xi\in [w-j,\dots,w-1], \ \ \ \psi\in[\xi+2,\dots, w+1]$$
Taking into consideration (9.22) we can rewrite (9.25.9) as follows:

$$P_{i-1,w-j,\psi,\xi}=-(\sum_{\be=0}^{j+\xi-w} \sum_{\al=i+\be}^{w+1-j+\be} S_{i-1,i-2,\al,i-2+\be} P_{\al,w-j+\be,\psi,\xi})+$$ $$+S_{i-1,i-2,\psi,i-2+\xi+j-w}\eqno{(9.25.10)}$$ with the same domain of $\xi$, $\psi$. This is (9.25.1). Because of (9.25.3), this formula is valid for $w-j\ge i-2$ (the non-trivial case of the definition of $P_{****}$). $\square$
\medskip
Let us consider a symmetry $\g s: \n Z^4\to\n Z^4$ defined as follows: $\g s(\al,\be,\ga,\de)=(w+2-\ga,w-1-\de,w+2-\al,w-1-\be)$.
\medskip
{\bf Remark 9.26.} $\g s$ has the following geometric interpretation. Let us consider a matrix $NL$ whose $(i,j)$-th entry is a symbol $N^{i-1}\hat l_{j}$. We interpret a quadruple $(\al,\be,\ga,\de)$ as a vector from $N^\be \hat l_\al$ to $N^\de \hat l_\ga$ in $NL$. $\g s$ is the reflection of this vector with respect to the center of $NL$ and the inversion of its direction.
\medskip
{\bf Definition 9.27.} $\bar S_{uvyz}:=-P_{\g s(uvyz)}^t$ (defined if $P_{\g s(uvyz)}$ has meaning).
\medskip
We shall consider block $r\times r$-matrices having the following block structure: their block size is $(w+1)\times(w+1)$, quantities of columns in blocks are $k_{w+1},k_w,\dots,k_1$ (counting from the left to the right), and quantities of lines in blocks are $k_{1},k_2,\dots,k_{w+1}$ (counting from up to down). Hence, the $(\al,\be)$-th block of this matrix is a $k_{\al}\times k_{w+2-\be}$-matrix. These matrices will be called skew $k_*$-block matrices.

\medskip
$\forall \ i = 0,\dots, w$ we define skew $k_*$-block matrices $C_i=C_i(S_{****})$ as follows:
\medskip
The $(\al, \be)$-th block of $C_i$ is $S^t_{w+2-\be,w+1-\be,\al,i}$ if the quadruple $(w+2-\be,w+1-\be,\al,i)$ satisfies (9.13, 9.14)\footnotemark \footnotetext{Really, it always satisfies (9.14).} (i.e. if $S_{w+2-\be,w+1-\be,\al,i}$ exists);
\medskip
The $(i+1,w+1-i)$-th block of $C_i$ is $I_{k_{i+1}}$, all other blocks of $C_i$ are 0:
$$(C_i)_{\al\be}=S^t_{w+2-\be,w+1-\be,\al,i}\eqno{(9.28.1)}$$
$$(C_i)_{i+1,w+1-i}=I_{k_{i+1}}\eqno{(9.28.2)}$$
Particularly, the $(\al,\be)$-th block of $C_i$ is a $(k_\al\times k_{w+2-\be})$-th matrix.
\medskip
$\forall \ i = 0,\dots, w$ we define skew $k_*$-block matrices $\bar C_i=\bar C_i(S_{****})$ as follows:
\medskip
The $(\al,\be)$-th block of $\bar C_i$ is given by the formula

$$(\bar C_i)_{\al,\be}=-P_{\al,w-1-i,w+2-\be,w-\be}=\bar S^{t}_{\be,\be-1,w+2-\al,i}\eqno{(9.29)}$$
if the quadruple $(\al,w-1-i,w+2-\be,w-\be)$ belongs to the non-trivial domain of $P_{****}$;
\medskip
For $i=0,\dots,w$
$$(\bar C_i)_{w+1-i,i+1}=I_{k_{w+1-i}}\eqno{(9.30)}$$
other block entries of $\bar C_i$ are 0.
\medskip
{\bf Remark 9.31.} Formula (9.30) is concordant with (9.29), if we consider $P_{****}$ from (9.22). Nevertheless, some 0-blocks of $\bar C_i$ correspond to $P_{**yz}$ where $(y,z)$ do not satisfy (9.19), and hence this $P_{****}$ is not defined.
\medskip
Finally, we define elements $B(S_{****}):=\sum_{i=0}^w C_iN^i\in M_r(\p)[N]$ and $\bar B(S_{****}):=\sum_{i=0}^w \bar C_iN^i\in M_r(\p)[N]$.
\medskip
Example for $w=3$:
$$B(S_{****})=\left(\matrix 0&0&0&I_{k_1}\\0&0&0& S^t_{1020}\\0&0&0& S^t_{1030}\\0&0&0& S^t_{1040} \endmatrix \right)+
\left(\matrix 0&0&0&0\\0&0&I_{k_2}& 0\\0&0&S^t_{2131}& S^t_{1031}\\0&0&S^t_{2141}& S^t_{1041} \endmatrix \right)N+$$ $$+ \left(\matrix 0&0&0&0\\0&0&0& 0\\0&I_{k_3}&0& 0\\0& S^t_{3242}&S^t_{2142}& S^t_{1042} \endmatrix \right)N^2 + \left(\matrix 0&0&0&0\\0&0&0&0\\0&0&0& 0\\I_{k_4}&0&0&0 \endmatrix \right)N^3$$
\medskip
$$\bar B(S_{****})=\left(\matrix \bar S^{t}_{1040}&0&0&0&\\ \bar S^{t}_{1030}&0&0&0\\ \bar S^{t}_{1020}&0&0&0\\I_{k_4}&0&0&0  \endmatrix \right)+
\left(\matrix \bar S^{t}_{1041}&\bar S^{t}_{2141}&0&0\\ \bar S^{t}_{1031}&\bar S^{t}_{2131}&0& 0\\0&I_{k_3}&0&0 \\0&0&0&0  \endmatrix \right)N+$$ $$+ \left(\matrix \bar S^{t}_{1042}&\bar S^{t}_{2142}& \bar S^{t}_{3242}&0\\0&0&I_{k_2}& 0\\0&0&0& 0\\0&0&0& 0  \endmatrix \right)N^2 + \left(\matrix 0&0&0&I_{k_1}\\0&0&0&0\\0&0&0& 0\\0&0&0&0 \endmatrix \right)N^3=$$
\medskip
$$=\left(\matrix -P_{1242}&0&0&0&\\-P_{2242}&0&0&0\\-P_{3242}&0&0&0\\I_{k_4}&0&0&0  \endmatrix \right)+
\left(\matrix -P_{1142}&-P_{1131}&0&0\\-P_{2142}&-P_{2131}&0& 0\\0&I_{k_3}&0&0 \\0&0&0&0  \endmatrix \right)N+$$ $$+ \left(\matrix -P_{1042}&-P_{1031}& -P_{1020}&0\\0&0&I_{k_2}& 0\\0&0&0& 0\\0&0&0& 0  \endmatrix \right)N^2 + \left(\matrix 0&0&0&I_{k_1}\\0&0&0&0\\0&0&0& 0\\0&0&0&0 \endmatrix \right)N^3$$
\medskip
{\bf Lemma 9.32.} $B(S_{****})^t \cdot \bar B(S_{****})=I_r N^w\in M_r(\p)[N]$.
\medskip
{\bf Proof.} We denote $B(S_{****})^t \cdot \bar B(S_{****})$ by $\sum_\mu \Cal C_\mu N^\mu$. The fact that $\Cal C_w=I_{r}$ is obvious: the only non-zero factors that enter in the sum $\sum_{\ga=0}^{w} C_\ga^t \bar C_{w-\ga}$ are products of blocks of $C_*$, $\bar C_*$ containing $I_*$, and they form $I_r$. Also it is obvious that for $\mu>w$ we have $\Cal C_\mu=0$, because all products whose sum is $\Cal C_\mu$, have at least one factor 0. We need to consider $\Cal C_\mu$ for $\mu<w$. (9.28.1), (9.28.2), (9.29), (9.30) give us (here and below $(C^t_\ga)_{\nu\de}$ is the $(\nu\de)$-th block of $C^t_\ga$, i.e. $(C^t_\ga)_{\nu\de}=((C_\ga)_{\de\nu})^t$ )
$$(\Cal C_\mu)_{\nu\pi}=\sum_{\ga=0}^\mu \sum_{\de=1}^{w+1} (C^t_\ga)_{\nu\de}(\bar C_{\mu-\ga})_{\de\pi}\eqno{(9.32.1.1)}$$ $$=-\sum_{\ga,\de} S_{w+2-\nu, w+1-\nu, \de, \ga}P_{\de,w-1-\mu+\ga,w+2-\pi,w-\pi}+ \eqno{(9.32.1.2)}$$
$$-P_{w+2-\nu,2w-\mu-\nu,w+2-\pi,w-\pi}+S_{w+2-\nu,w+1-\nu,w+2-\pi,\mu+1-\pi}\eqno{(9.32.1.3)}$$
where (9.32.1.2) corresponds to the products $(C^t_\ga)_{\nu\de}(\bar C_{\mu-\ga})_{\de\pi}$ where both terms $\ne 0, \ I_*$, and (9.32.1.3) corresponds to the products where one of the terms is $I_*$.

Let us find the relations satisfied by $\mu,\ \nu, \ \pi$ and the domain of summation by $\ga, \ \de$ in (9.32.1.2). We have
$$\matrix (C^t_\ga)_{\nu\de}\ne0, I_{k_*} \ \iff \ \nu\ge w+1-\ga \ \ \wedge \ \ \de\ge\ga+2 \\
(\bar C_{\mu-\ga})_{\de\pi}\ne0, I_{k_*} \ \iff \ \de\le w-(\mu-\ga) \ \ \wedge \ \ \pi \le \mu-\ga+1\endmatrix \eqno{(9.32.2)}$$
The set of $\ga, \ \de$ is non-empty $\iff \ \mu\le w-2$ and $\mu +\nu -\pi\ge w$. In this case the conditions (9.32.2) on $\ga, \ \de$ become
$$\mu+1-\pi\ge\ga\ge w+1-\nu\eqno{(9.32.3.1)}$$
$$w+\ga-\mu\ge\de\ge\ga+2\eqno{(9.32.3.2)}$$
Now we use Proposition 9.25 for
$$\matrix i=w+3-\nu && \mu=j+i-3 \\
j=\mu+\nu-w &\iff & \nu=w-i+3 \\
\psi=w-\pi+2 && \pi=w-\xi=w-\psi+2\\
\xi=w-\pi \endmatrix\eqno{(9.32.4)}$$
and summation variables $\al, \ \be$ in (9.25.1) are
$$\matrix \al=\de \\ \be=\ga-i+2 \endmatrix\eqno{(9.32.5)}$$
Under this variable change, (9.32.1.2) becomes the double sum in (9.25.10), and (9.32.1.3) becomes
$$-P_{i-1,w-j,\psi,\xi}+S_{i-1,i-2,\psi,i-2+\xi+j-w}$$
hence the desired. $\square$
\medskip
Let for $i=0,\dots,w-1$ $X_i$ be skew $k_*$-matrices having the following property:
\medskip
If $(\al, \ \be)$ are such that the $(\al,\ \be)$-block of $\bar C_i$ is 0 or $I_*$ then $(X_i)_{\al\be}=(\bar C_i)_{\al\be}$;

If $(\al, \ \be)$ are such that the $(\al,\ \be)$-block of $\bar C_i$ is $\ne 0, \ I_*$ then $(X_i)_{\al\be}$ is arbitrary.
\medskip
We denote $X:=\sum_{i=0}^{w-1}X_iN^i$.
\medskip
{\bf Lemma 9.33.} If $B(S_{****})^t \cdot X\in N^wM_r(\p[N])$ then $X=\bar B(S_{****})$.
\medskip
{\bf Proof.} For any fixed $\mu, \ \nu, \ \pi$ (9.32.1.2), (9.32.1.3) become
$$\sum_{\ga,\de} S_{w+2-\nu, w+1-\nu, \de, \ga} (X_{\ga-\mu})_{\de,\pi}\eqno{(9.33.1a)}$$ $$
+ (X_{\mu+\nu-w-1})_{w+2-\nu,\pi}+S_{w+2-\nu,w+1-\nu,w+2-\pi,\mu+1-\pi}=0\eqno{(9.33.1b)}$$
This is system of linear equations with unknowns $(X_i)_{\al\be}$ where $$\matrix 0\le i \le w-1 \\ \\ 1\le \be\le i+1\\ \\ 1\le\al \le w-i \endmatrix \eqno{(9.33.2)}$$ (for other values of $i,\ \al, \ \be$ we have $(X_i)_{\al\be}=0$). We arrange $(X_i)_{\al\be}$ in decreasing order of $i+\al$ (for $(X_i)_{\al\be}$ having equal $i+\al$ their ordering is arbitrary), and we arrange equations (9.33.1) in decreasing order of $\mu$ (for equations having equal $\mu$ the order of equations corresponds to the order of $(X_i)_{\al\be}$ having equal $i+\al$). Under this arrangement of unknowns and equations, the matrix of the system (9.33.1) becomes unitriangular. Really, for any $i,\ \al, \ \be$ satisfying (9.33.2) there is exactly one values of $\mu, \ \nu, \ \pi$ --- namely, $$\matrix \mu=i+\al-1 \\ \\ \nu=w+2-\al\\ \\ \pi=\be \endmatrix $$ such that the first term of (9.33.1b) is $(X_i)_{\al\be}$. Other terms of (9.33.1) for these $\mu, \ \nu, \ \pi$ --- namely, the terms that enter in (9.33.1a) --- contain $(X_{\ga-\mu})_{\de,\pi}$ such that $\ga-\mu+\de>i+\al$ hence unitriangularity.
\medskip
Lemma 9.32 affirms that
$$(X_i)_{\de,\pi}=-P_{\de,w-1-i,w+2-\pi,w-\pi}$$
is a solution to this system. Unitriangularity implies that this solution is unique. $\square$
\medskip
Let us consider $\g q_H$ from above (9.1).
\medskip
{\bf Lemma 9.34.} $\forall \ u=1,\dots, w+1$, for $v=u-1$, $\forall \ i=1,\dots, k_u$ the elements

$$\om_{ui}:=N^{v}l_{ui}+\sum_{y=u+1}^{w+1}\sum_{z=u-1}^{y-2}\sum_{j=1}^{k_y}(S_{uvyz})_{ij}N^zl_{yj}\eqno{(9.34.1)}$$

form a basis of $\g q_H$. $\square$
\medskip
We denote the set of elements $\om_{ui}$ ($u$ is fixed, $i$ varies) by $\hat \om_{u}$ (matrix columns). So, (9.34.1) becomes ($v=u-1$)
$$\hat \om_{u}=N^{v}\hat l_{u}+\sum_{z=u-1}^{w-1}\sum_{y=z+2}^{w+1}S_{uvyz}N^z\hat l_{y}\eqno{(9.35)}$$

Let $\vf_i$, $i=1,...,r$ be the basis of $L'$ dual to $l_i$, i.e. $\vf_i(l_j)=\delta_{ij}$. We shall need the dual numbers $k'_i:=k_{w+2-i}$ (inverse order of $k_*$). We consider the analogous two-subscript notation of $\vf_i$, but the order of segments of the partition of $\vf_i$ is opposite, namely:
$$(\vf_{w+1,1},\dots,\vf_{w+1,k'_{w+1}},\ \ \vf_{w,1},\dots,\vf_{w,k'_{w}},\ \ \ \dots \ \ \ \vf_{11},\dots,\vf_{1,k'_{1}}):= (\vf_1,\dots, \vf_r)\eqno{(9.35.1)}$$
(order of elements $\vf_*$ is the same in both sides of this equality).
\medskip
{\bf Lemma 9.36.} $\forall \ u=1,\dots, w+1$, for $v=u-1$, $ \forall \ i=1,\dots, k'_u$ the elements

$$\chi_{ui}:=N^{v}\vf_{ui}+\sum_{y=u+1}^{w+1}\sum_{z=u-1}^{y-2}\sum_{j=1}^{k'_y}(\bar S_{uvyz})_{ij}N^z\vf_{yj}\eqno{(9.36.1)}$$

form a basis of ${\g q'_H}^w$.

{\bf Proof.} As above we denote the set of elements $\vf_{ui}$, resp. $\chi_{ui}$ ($u$ is fixed, $i$ varies) by $\hat \vf_{u}$, resp. $\hat \chi_{u}$ (matrix columns). (9.35), (9.36.1) can be written in terms of blocks of $C_i$, $\bar C_i$:

$$\hat \om_{u}=\sum_{z=0}^{w}\sum_{y=1}^{w+1} (C_z^t)_{w+2-u,y}N^z\hat l_{y}$$
$$\hat \chi_{u}=\sum_{z=0}^{w}\sum_{y=1}^{w+1} (\bar C_z^t)_{uy}N^z\hat \vf_{w+2-y}$$
We must prove that $\forall \ u_1, \ u_2$ we have $\hat \om_{u_1}\hat \chi_{u_2}^t=\de_{u_1}^{u_2}I_{k_{u_1}}N^w$ (product is pairing). This is immediate:

$$\hat \om_{u_1}\hat \chi_{u_2}^t=\sum_{z_1=0}^{w}\sum_{y_1=1}^{w+1}\sum_{z_2=0}^{w}\sum_{y_2=1}^{w+1} (C_{z_1}^t)_{w+2-u_1,y_1}N^{z_1}\hat l_{y_1} N^{z_2}\hat \vf_{w+2-y_2}^t (\bar C_{z_2})_{y_2,u_2}$$
We have $\hat l_{y_1} \hat \vf_{w+2-y_2}^t=\de_{y_1}^{y_2}I_{k_{y_1}}$, hence
$$\hat \om_{u_1}\hat \chi_{u_2}^t=\sum_{z_1=0}^{w}\sum_{y=1}^{w+1}\sum_{z_2=0}^{w} (C_{z_1}^t)_{w+2-u_1,y}N^{z_1+z_2} (\bar C_{z_2})_{y,u_2} = \sum_{z=0}^{2w}(\Cal C_z)_{w+2-u_1,u_2}N^z$$
Lemma 9.32 implies the desired. $\square$
\medskip
{\bf Corollary 9.37.} Matrices $S_{uvyz}(L')$ for the dual lattice $L'$ are $\bar S_{uvyz}(L)$ (order of segments of $\vf_*$, and hence of numbers $k_*$, is inverse).
\medskip
We define $D_i$ ($i=1,\dots, w$ ) as $-<N^{i-1}(l_*),f_*>$ where the ordering of $l_*$ is $l_{11},\dots, l_{w+1, k_{w+1}}$. Hence, $D_i$ is a union of $w+1$ matrices $D_{ij}$, where $D_{ij}$ is a $r\times k_j$-matrix, and $(D_{ij})_{\al\be}=-<N^{i-1}(l_{j\be}),f_\al>$. This means that there are relations between $D_{ij}$ coming from (9.20), namely:

$$D_{v+1,u}= -\sum_{z=v}^{w-1}\sum_{y=z+2}^{w+1}D_{z+1,y} P_{uvyz}^t\eqno{(9.38)}$$

Like in (9.21), for $(z,y)$ satisfying $y\ge z+1$ (resp. $y< z+1$) we shall call the corresponding $D_{zy}$ as belonging to the trivial (resp. non-trivial) domain.
\medskip
{\bf Lemma 9.39.} $\Psi_N B(S_{****})\in M_r(\p[[N]])$.
\medskip
{\bf Proof.} For any $\mu=0,\dots, w-1$ we must prove that $\g C_\mu:=\sum_{\de =0}^\mu D_{w-\de}C_{\mu-\de}$ is 0. The $\nu$-th block ($\nu=1,\dots,w+1$) of this matrix is
$$(\g C_\mu)_\nu:=\sum_{\de =0}^\mu \sum_{\ga=1}^{w+1} D_{w-\de,\ga}(C_{\mu-\de})_{\ga \nu}\eqno{(9.39.1)}$$
We have
$$(C_{\mu-\de})_{\ga \nu}\ne0, I_{k_*} \ \ \iff \de\le \nu+\mu-w-1 \ \ \wedge \ \ \ga\ge \mu-\de+2$$
$$(C_{\mu-\de})_{\ga \nu}=I_{k_*} \ \ \iff \de=\nu+\mu-w-1 \ \ \wedge \ \ \ga=w+2-\nu$$
hence (9.39.1) becomes
$$(\g C_\mu)_\nu=\sum_{\de =0}^{\nu+\mu-w-1}\sum_{\ga=\mu-\de+2}^{w+1}D_{w-\de,\ga}S^t_{w+2-\nu,w+1-\nu,\ga,\mu-\de}+\eqno{(9.39.2.1)}$$
$$+D_{2w+1-\nu-\mu,w+2-\nu}\eqno{(9.39.2.2)}$$
where (9.39.2.1) is non-empty if $\nu+\mu\ge w+1$.

Terms $D_{2w+1-\nu-\mu,w+2-\nu}$ always belong to the non-trivial domain. We separate the terms of (9.39.2.1) in terms of trivial and non-trivial domain:
$$(\g C_\mu)_\nu=\sum_{\de =0}^{\nu+\mu-w-1}\sum_{\ga=\mu-\de+2}^{w-\de} D_{w-\de,\ga}  S^t_{w+2-\nu,w+1-\nu,\ga,\mu-\de}+\eqno{(9.39.3.1)}$$
$$+\sum_{\de =0}^{\nu+\mu-w-1}\sum_{\ga=w-\de+1}^{w+1}  D_{w-\de,\ga}  S^t_{w+2-\nu,w+1-\nu,\ga,\mu-\de}+\eqno{(9.39.3.2)}$$
$$+D_{2w+1-\nu-\mu,w+2-\nu}\eqno{(9.39.3.3)}$$

Now we substitute non-trivial $D_{**}$ by linear combinations of the trivial ones, using (9.38):
$$(\g C_\mu)_\nu=-\sum_{\de =0}^{\nu+\mu-w-1}\sum_{\ga=\mu-\de+2}^{w-\de}   \sum_{z=w-\de-1}^{w-1}\sum_{y=z+2}^{w+1}D_{z+1,y} P_{\ga,w-\de-1,y,z}^t       S^t_{w+2-\nu,w+1-\nu,\ga,\mu-\de}+\eqno{(9.39.4.1)}$$
$$+\sum_{\de =0}^{\nu+\mu-w-1}\sum_{\ga=w-\de+1}^{w+1}D_{w-\de,\ga}S^t_{w+2-\nu,w+1-\nu,\ga,\mu-\de}-\eqno{(9.39.4.2)}$$
$$-\sum_{z=2w-\nu-\mu}^{w-1}\sum_{y=z+2}^{w+1}D_{z+1,y} P_{w+2-\nu,2w-\nu-\mu,y,z}^t \eqno{(9.39.4.3)}$$
Now we change variables in (9.39.4.2):
$$\de=w-z-1$$
$$\ga=y$$
interchange the order of summation and transpose:
$$(\g C_\mu)_\nu=\sum_{z=2w-\nu-\mu}^{w-1}\sum_{y=z+2}^{w+1} \g K(\mu,\nu,z,y) D^t_{z+1,y}$$

where $$\g K(\mu,\nu,z,y)=-( \sum_{\de =w-z-1}^{\nu+\mu-w-1}\sum_{\ga=\mu-\de+2}^{w-\de} S_{w+2-\nu,w+1-\nu,\ga,\mu-\de}    P_{\ga,w-\de-1,y,z}) +$$
$$+S_{w+2-\nu,w+1-\nu,y,\mu-w+z+1}-P_{w+2-\nu,2w-\nu-\mu,y,z} \eqno{(9.39.5)}$$
Change of variables in (9.39.5):
$$\matrix y=\psi & \nu=w+3-i & \ga=\al \\ z=\xi & \mu=i+\al-3 & \de=\al-\be-1 \endmatrix $$
transforms (9.39.5) to (9.25.1), hence all $\g K(\mu,\nu,z,y)$ are 0. $\square$
\medskip
Let $M'=M^{\prime w}$ be the $w$-dual of $M$.
\medskip
{\bf Lemma 9.40.} For $i=1,\dots, w+1$ we have $k_i(M')=k'_{i}(M)$.
\medskip
{\bf Proof.} We use Lemma 10.2. $\mu$ of 10.2 is $w$, $m_i$ of 10.2 are $d_{w+2-i}$. (10.2.5) means that $d_i(M')=w-d_i$. The result follows immediately from the properties of dual 0-partitions. $\square$

We shall another basis $\hat \eta$ of $L(M')$ obtained by a permutation matrix from the basis (9.35.1). Namely, we let $\eta_{ij}:=\vf_{ij}$ ($i=1, \dots, w+1, \ j=1,\dots, k'_i$), but the order of elements $\eta_{ij}$ is the following ($\hat \eta$ is a matrix column):

$$\hat \eta=(\eta_{11},\dots, \eta_{1,k'_1}, \eta_{21},\dots, \eta_{2,k'_2}, \dots, \eta_{w+1,1},\dots, \eta_{w+1,k'_{w+1}})^t$$

We shall consider only $M$ satisfying the following (compare with 9.12)
\medskip
{\bf Condition 9.41.} For any $u$ elements $N^\al(\eta_{\be\ga})$, where $\al\in[u,\dots, w-1]$, $\be\in[\al+2,\dots, w+1]$, $\ga\in [1,\dots, k'_\be]$, are linearly independent over $\p$.
\medskip
According the general principle that almost all $n$-uples $(v_1,\dots,v_n)$ of vectors in $n$-dimensional vector space form a basis of this space, we can guess that almost all $M$ satisfy 9.41. Really, Lemma 9.40 affirms that the dimension of $N^u \Lie(M')$ is exactly the quantity of elements $N^\al(\eta_{\be\ga})$ mentioned in 9.41. Again by Lemma 9.40, we see that Condition 9.41 implies that these elements are a basis of $N^u \Lie(M')$.
\medskip
We can use ideas of [GL] in order to show that a large class of Anderson t-motives $M$ satisfies Condition 9.41. Namely,
\medskip
Finally, we have
\medskip
{\bf Conjecture 9.42.} All $M$ satisfy Condition 9.41.
\medskip
{\bf Theorem 9.43.} Let $M$ satisfy Condition 9.41. Then
\medskip
{\bf Proof.} We denote the basis $\vf_*$ from (9.35.1) by $\hat \vf$ (matrix column). We have $\hat \eta= \g I\cdot \hat \vf$ where $\g I$ is a matrix of the change of bases. It is a skew $k'_*$-block anti-identity matrix, i.e. its antidiagonal block entries are identity matrices: the $(w+2-i,i)$-block is $I_{k_i}$, and other block entries are 0.

We denote by $\Psi'_N$, resp. $\Psi'_{N,\eta}$ the $N$-$\Psi$-series of $M'$ in bases $\hat \vf$, resp. $\hat \eta$ (as a base of $M'$ over $\p[T]$ we use $\hat f'$ in both cases). We have $\Psi'_N=\Psi_N^{t-1}\Xi_N^{-w}$, $\Psi'_{N,\eta}=\Psi'_N\g I^t$. Since $M$ satisfies Condition 9.41, there exists a set of Siegel matrices for $M'$ with respect to the basis $\hat \eta$. We denote it by $U_{****}$. It defines $B(U_{****})$ --- the corresponding $B$ in $M_r(\p)[N]$.
\medskip
We denote $\Psi_NB(A(M)_{****})$, resp. $\Psi'_{N,\eta}B(U_{****})$ by $\g Z(M)$, resp. $\g Z(M^d)$. We have $\g Z(M)\in M_r(\p[[N]])$, $\g Z(M^d)\in M_r(\p[[N]])$ (Lemma 9.39), hence
$$\g Z(M)^t\g Z(M^d)=B(A(M)_{****})^t\Psi_N^t \Psi'_N \g I^t B(U_{****})=$$ $$B(A(M)_{****})^t \g I^t B(U_{****}) \Xi_N^{-w}\in M_r(\p[[N]])$$ and hence
$$  B(A(M)_{****})^t \g I^t B(U_{****}) \g I^t  \in \Xi_N^{w}M_r(\p[[N]])$$
We have  $\g I^t B(U_{****}) \g I^t$ is of the form $X$ of Lemma 9.33. Further, $\Xi_N^{w}\in N^wM_r(\p[[N]])$, hence $\g I^t B(U_{****}) \g I^t=\bar B(A(M)_{****})$ (Lemma 9.33). This means that $U_{uvyz}=\bar A_{uvyz}$. The theorem follows from Lemma 9.36. $\square$
\medskip

\newpage
{\bf Case of linearly dependent columns.}

We shall consider for simplicity the case $w=2$. We use notations: $\Psi:=\Psi_N(M)$, $\Psi^d:=\Psi_N(M^d)$, $\Psi^d=D^d_2N^{-2}+D^d_1N^{-1}+...$ in the basis $\vf_1,...,\vf_r$ dual to $l_1,...,l_r$. We denote $D^d_2=D_{21}^d | D_{22}^d | D_{23}^d $ where $D_{2i}^d$ are $r\times k_i$ matrices.  Let some of the columns of $D_{21}^d$ are linearly dependent. Interchanging the order of columns we can assume that $\exists \ u>0$ such that the first $k_1-u$ columns of $D_{21}^d$ are linearly independent and the last $u$ columns of $D_{21}^d$ are their linear combinations. Further, we assume that the remaining $u$ columns of $D_{21}^d$ forming (together with the first $k_1-u$ columns of $D_{21}^d$) a basis of $NV$ are the first $u$ columns of $D_{22}^d$. This means that we can represent $D_{21}^d$, $D_{22}^d$ as unions:

$D_{21}^d=(D_{211}^d | D_{212}^d)$, $D_{22}^d=(D_{221}^d | D_{222}^d)$, where the quantities of rows in all these matrices are $r$ and the quantities of columns in $D_{211}^d$, resp. $ D_{212}^d$, $D_{221}^d, \ D_{222}^d$ are $k_1-u$, resp. $u, \ u, \ k_2-u$. We have: $\exists$ a $(k_1-u)\times u$-matrix $Y$ such that $D_{212}^d=D_{211}^dY$.

We denote by $\tilde \Psi^d$ the matrix obtained from $\Psi^d$ by permutation of (12)-block and (21)-block. We assume that (1,2,3)-blocks of $\tilde \Psi^d$ are "good", i.e. there exists $\tilde D^d$ such that

$$\tilde \Psi^d\cdot\tilde D^d\in M_r(\p[[N]])$$

We denote by $\hat \Psi^d$ the matrix obtained from $\Psi^d$ by elimination of (12)-block. It is a $5\times 4$-block matrix, and we have

$$\Psi^d=\hat \Psi^d U_3$$ where $U_3=\left(\matrix 1&Y&0&0&0\\ 0&0&1&0&0\\ 0&0&0&1&0\\ 0&0&0&0&1\endmatrix \right)$ (block structure; sizes of blocks: $k_1-u; \ u; \ k_2-u; \ k_3$ for columns; $k_1-u; \ u; \ u; \ k_2-u; \ k_3$ for lines).

$$\tilde \Psi^d=\hat \Psi^d U_2$$ where $U_2=\left(\matrix 1&0&Y&0&0\\ 0&1&0&0&0\\ 0&0&0&1&0\\ 0&0&0&0&1\endmatrix \right)$ (block structure; sizes of blocks: $k_1-u; \ u; \ k_2-u; \ k_3$ for columns; $k_1-u; \ u; \ u; \ k_2-u; \ k_3$ for lines).
\medskip
$U_3$ has a right inverse $U_3^{-1}=\left(\matrix 1&0&0&0\\ 0&0&0&0\\ 0&1&0&0\\ 0&0&1&0\\ 0&0&0&1\endmatrix \right)$ such that we have

$$\hat \Psi^d = \Psi^dU_3^{-1}$$

hence $B^t\Psi^t\tilde \Psi^d\cdot\tilde B^d\in M_r(\p[[N]])$

and $\tilde \Psi^d=\Psi^dU_3^{-1}U_2$

As earlier we have $$\Psi^t\Psi^d=\Xi^{-w}$$

hence $B^t U_3^{-1}U_2 B^d\in N^2 M_r(\p[[N]])$
\medskip
\medskip

Generalizations.
\medskip
1. To prove that $\g q_{M_1\otimes M_2}=\g q_{M_1} \otimes \g q_{M_2}$.
\medskip
2. It is known that if $A_1$, $A_2$ are abelian varieties then $A_1\otimes A_2$ is a mixed motive. Is it possible to get an analog of (1) for it?
\medskip
3. Analog of the main theorem for groups other than $GL_r$. We have a theorem that any full sublattice of $\n Z^r$ is described by its elementary divisors $d_1 | d_2| ... | d_r$. There is an analog of this theorem for any reductive group. We should find the corresponding generalization of the main theorem of the present paper.
\medskip
This is an analog of Theorem 5 for arbitrary $N$, and of Theorem 6 for the operator of duality. Its proof, as well as its generalization to the case of non-pure $M$, can be easily obtained using the same ideas of the proofs of Theorems 5, 6 (this explains the terminology: this is not a conjecture, but a theorem whose proof is not written yet).
\medskip
For $m=1$ Formula 9.1 and Theorem 5 imply immediately
\medskip
{\bf Proposition 9.4.} Result 9.3 is proved for $m=1$. $\square$
\medskip

\input amstex
\documentstyle{amsppt}
\magnification1200
\tolerance=10000
\overfullrule=0pt
\def\n#1{\Bbb #1}
\def\p{\Bbb C_{\infty}}
\def\fr{\hbox{fr}}

\def\im{\hbox{im }}
\def\invlim{\hbox{invlim}}
\def\tr{\hbox{tr }}

\def\Gal{\hbox{Gal }}

\def\Exp{\hbox{Exp}}

\def\Hom{\hbox{Hom}}

\def\End{\hbox{End}}
\def\Prin{\hbox{Prin}}
\def\Ker{\hbox{Ker }}
\def\Lie{\hbox{Lie}}

\def\Div{\hbox{Div}}
\def\Pic{\hbox{Pic}}

\def\ord{\hbox{ord}}

\def\Id{\hbox{Id}}
\def\Cl{\hbox{Cl}}
\def\Supp{\hbox{ Supp }}
\def\Spec{\hbox{ Spec }}

\def\diag{\hbox{ diag }}
\def\Diag{\hbox{ Diag }}

\def\e11{E_{11}}

\def\ga{\goth A}
\def\w{\hbox{\bf A}}
\def\x{\hbox{\bf K}}
\def\ve{\varepsilon}
\def\vf{\varphi}

\def\ve{\varepsilon}

\def\vf{\varphi}

\def\de{\delta}

\def\ga{\gamma}

\def\be{\beta}

\def\al{\alpha}

\def\om{\omega}
\def\g{\goth }

\topmatter
\title
Duality of Anderson $t$-motives
\endtitle
\author
A. Grishkov, D. Logachev\footnotemark \footnotetext{E-mails: shuragri{\@}gmail.com; logachev94{\@}gmail.com (corresponding author)\phantom{*******************}}
\endauthor
\thanks Thanks: The authors are grateful to FAPESP, S\~ao Paulo, Brazil for a financial support (process No. 2017/19777-6). The first author is grateful to SNPq, Brazil, to RFBR, Russia, grant 16-01-00577a (Secs. 1-4), and to Russian Science Foundation, project 16-11-10002 (Secs. 5-8) for a financial support. The second author is grateful to Gilles Lachaud for invitation to IML and to Laurent Lafforgue for invitation to
IHES where this paper was started, and to Vladimir
Drinfeld for invitation to the University of Chicago where
this paper was continued. Discussions with Greg Anderson, Vladimir
Drinfeld, Laurent Fargues, Alain Genestier, David Goss, Richard Pink, Yuichiro
Taguchi, Dinesh Thakur on
the subject of this paper were very important. Particularly, Alain
Genestier informed me about the paper of Taguchi where the notion of
the dual of a Drinfeld module is defined. Further, Richard Pink indicated me an important reference (see Section 6 for details); proof of the main theorem of the present paper grew from it. Finally, Vladimir
Drinfeld indicated me the proof of the Theorem 12.6 and Jorge Morales
gave me a reference on classification of quadratic forms over $\n
F_q[T]$ (Remark 7.8).
\endthanks
\NoRunningHeads
\address
First author: Departamento de Matem\'atica e estatistica
Universidade de S\~ao Paulo. Rua de Mat\~ao 1010, CEP 05508-090, S\~ao Paulo, Brasil, and Omsk State University n.a. F.M.Dostoevskii. Pr. Mira 55-A, Omsk 644077, Russia.
\medskip
Second author: Departamento de Matem\'atica, Universidade Federal do Amazonas, Manaus, Brasil
\endaddress

\abstract Let $M$ be a t-motive. We introduce the notion of duality for $M$. Main results of the paper (we consider uniformizable $M$ over $\n F_q[T]$ of rank $r$, dimension $n$, whose nilpotent operator $N$ is 0):

1. Algebraic duality implies analytic duality (Theorem 5). Explicitly, this means that the lattice of the dual of $M$ is the dual of the lattice of $M$, i.e. the transposed of a Siegel matrix of $M$ is a Siegel matrix of the dual of $M$.

2. Let $n=r-1$. There is a 1 -- 1 correspondence between pure t-motives (all they are uniformizable), and lattices of rank $r$ in $\p^{n}$ having dual (Corollary 8.4).

3. Pure t-motives have duals which are pure t-motives as well (Theorem 10.3).

4. For a self-dual uniformizable $M$ a polarization form on its lattice $L(M)$ is defined. For some $M$ this form is skew symmetric like in the number field case, but for some other $M$ it is symmetric. An example is given.

5. Some explicit results are proved for $M$ having complete
multiplication. The CM-type of the dual of $M$ is the complement of
the CM-type of $M$. Moreover, for $M$ having multiplication by a
division algebra there exists a simple formula for the
CM-type of the dual of $M$ (Section 12).

6. We construct a class of non-pure t-motives (t-motives having the
completely non-pure row echelon form) for which duals
are explicitly calculated (Theorem 11.5). This is the first step
of the problem of description of all t-motives having duals.

7. If $M$ has good ordinary reduction then the kernels of reduction maps on groups of  torsion points for $M$ and its dual are complementary with respect to a natural pairing (proof is given for a particular case, Conjecture 13.4.1).
\endabstract
\endtopmatter
\document
{\bf 0. Introduction.}
\nopagebreak
\medskip
t-motives are the function field analogs of abelian varieties (more exactly, of abelian varieties with multiplication by an imaginary quadratic field, see [L1]). Main references for t-motives are [A], [G]. Nevertheless, function field analogs of some basic results in the theory of abelian varieties are not known yet.

The present paper contains an analog of such result. Namely, we introduce the notion of duality for a t-motive $M$ (this is not the duality in a Tannakian category!), and we prove some properties of this notion, see the abstract. Particularly, if $M$ is uniformizable and has dual then the lattice of the dual of $M$ is the dual of the lattice of $M$ (Theorem 5). An immediate corollary of the above theorem and the result of Drinfeld on 1 -- 1 correspondence between Drinfeld modules and lattices in $\p$ (here $\p$ is the function field analog of $\n C$) is Corollary 8.4: there is a 1 -- 1 correspondence between pure t-motives of dimension $r-1$ and rank $r$, and lattices of rank $r$ in $\p^{r-1}$ having dual (not all such lattices have dual).

Let us give more details on the contents of the paper. For simplicity, most results are proved for t-motives over the ring $\hbox{\bf A}=\n F_q[T]$, and we consider, with few exceptions, only the case $N=0$. The main definition of duality of t-motives (definition 1.8 --- case $\hbox{\bf A}=\n F_q[T]$ and definition 1.13 --- general case) is given in Section 1.\footnotemark \footnotetext{A version of the definition of duality is obtained independently in [Tae], 2.2.} Section 2 contains the definition of the dual lattice. Section 3 contains explicit formulas for the dual lattice. Section 5 contains the statement and the proof of the main theorem 5 --- coincidence of algebraic and analytic duality for the case $\hbox{\bf A}=\n F_q[T]$ (section 4 contains the statement of the corresponding conjecture for the case of general $\hbox{\bf A}$). Section 6 contains the theorem 6 describing the lattice of the tensor product of two t-motives (case $N=0$; the proof for the general case was obtained, but not published, by Anderson). Section 7 contains the notion of self-dual t-motives and polarization form on them. Some examples are given. We discuss in Section 8 the problem of correspondence between uniformizable t-motives and lattices. Section 9 gives the statement of the main result for the case $N\ne 0$ without proof and a reformulation of the theorem 5 in terms of Pink-Hodge structures of constant weight.

Further, Lemma 1.10 gives the explicit matrix form of the definition of duality of t-motives. Since Taguchi in [T] gave a definition of dual to a Drinfeld module, we prove in Proposition 1.12.3 that the definition of the present paper is equivalent to the original definition of Taguchi. Section 1.14 contains a definition of duality for abelian $\tau$-sheaves ([BH], Definition 2.1), but we do not develop this subject. We prove in Section 10 that pure t-motives have duals which are pure
t-motives as well, and some related results (a proof that the dual of an abelian
$\tau$-sheaf is also an abelian $\tau$-sheaf can be obtained using ideas
of Section 10). In Section 11 we consider t-motives having the
completely non-pure row echelon form, and we give an explicit formula
for their duals. In Section 12 we consider t-motives with complete multiplication, and we give for them a very simple version of the proof of the first part of the main theorem.
Section 13 contains some explicit formulas for t-motives of complete multiplication. In 13.1 we describe the dual lattice, in 13.2 we show that the results of Section 12 are compatible with (the first form of) the main theorem of complete multiplication. Section 13.3 contains an explicit proof of the main theorem for t-motives with complete multiplication by two types of simplest fields.
Section 13.4 gives us an application of the notion of duality to the
reduction of t-motives (subject in development, see [L2]).
\medskip
{\bf Notations. }
\medskip
$q$ is a power of a prime $p$;
\medskip
{\it Case of $M$ over $\n F_q[T]$}:
\medskip
$\n Z_\infty:=\n F_q[\theta]$, $\n R_{\infty}:=\n F_q((1/\theta))$, $\p$ is the completion of its algebraic closure
($\n Z_\infty$, $\n R_{\infty}$, $\p$ are the function field analogs of $\n Z$, $\n R$, $\n C$ respectively);

$\w:=\n F_q[T]$, $\hbox{\bf K}:=\n F_q((1/T))$;

$\iota: \hbox{\bf A} \to \p$ ($\iota(T)=\theta$) is the standard map of generic characteristic (with one exception (1.16), we shall not consider the case of finite characteristic);
we extend $\iota$ to $\hbox{\bf K}$, and we have $\n Z_\infty=\iota(\w)\subset \p$, $\n R_\infty=\iota(\hbox{\bf K})\subset \p$.

$\goth C$ (resp. $\goth C_2$) is the Carlitz module over $\w=\n F_q[T]$ (resp. over $\n F_{q^2}[T]$).
\medskip
{\it Case of $M$ over an extension of $\n F_q[T]$}:
\medskip
$\n Q_\infty$ is a finite separable extension of $\n F_q(\theta)$;

$\infty$ is a fixed valuation of $\n Q_\infty$ over the infinity of $\n F_q(\theta)$;

${\n Z_\infty \subset \n Q_\infty}$ is the subring of elements which are regular outside $\infty$;

$\n R_\infty$ is the completion of $\n Q_\infty$ at infinity, and $\p$ --- the completion of its algebraic closure --- is the same as of the case of $M$ over $\n F_q[T]$.

$\w\supset \n F_q[T]$, $\hbox{\bf K}\supset \n F_q((1/T))$ are defined by the condition that $\iota: \w \to \n Z_\infty$, $\iota: \hbox{\bf K} \to \n R_\infty$ are isomorphisms.

$\w_C:=\w\underset{\n F_q}\to{\otimes} \p$ (i.e. $\w_C=\p[T]$ for the case of $M$ over $\n F_q[T]$).

$\goth C$ is a Drinfeld module of rank 1 over {\bf A}.
\medskip
If $P=\frac{\sum a_iT^i}{\sum b_iT^i}\in \p(T)$ then $P^{(k)}:=\frac{\sum a_i^{q^k}T^i}{\sum b_i^{q^k}T^i}$. For $x\in \w_C$, $x=a\otimes z$, $a\in \w$, $z\in \p$ we let $x^{(k)}:= a\otimes z^{q^k}$.

$M_r$ is the set of $r\times r$ matrices. If $C=\{c_{ij}\}$ is a matrix with entries $c_{ij}\in \p(T)$ then $C^{(k)}:=\{c_{ij}^{(k)}\}$, $C^t$ is the  transposed of $C$, $C^{(k)\ -1}=(C^{(k)})^{-1}$, $C^{t-1}=(C^t)^{-1}$.

If $M$ is an $\w_C$-module, we define $M^{(1)}$ as the tensor
product $M\otimes_{\w_C,*^{(k)}}\w_C$ with respect to the map $*^{(k)}: \w_C\to \w_C$ (this notation is concordant in the obvious sense with the above notation $C^{(1)}$).

For a t-motive $M$ we denote by $E=E(M)$ the corresponding t-module (see [G], Theorem 5.4.11; Goss uses the inverse functor $E\mapsto M=M(E)$).

$\Lie(M)$ is $\Lie(E(M))$ ([G], 5.4).

$E_k$ is the unit matrix of size $k$.

Throughout the whole paper the word "canonical" will mean "canonical up to multiplication by elements of $\n F_q^*$".
\medskip
{\bf 1. Definitions. }

\nopagebreak
\medskip
If otherwise is not explicitly stated, throughout the whole paper we consider the case of t-motives $M$ over the ring $\hbox{\bf A}=\n F_q[T]$ such that $N=N(M)=0$. Exceptions: case of arbitrary $\w$ is treated in Sections 1.13, 1.14, 2, 4, 5.2. Case of arbitrary $N$ is treated in Sections 1, 10 and in statements of some results of Anderson in Sections 5, 6.
\medskip
In the present section we consider $M$ such that $N(M)$ is arbitary.
\medskip
Let $\p[T,\tau]$ be the Anderson ring, i.e. the ring of non-commutative polynomials satisfying the following relations (here
$a \in \p$):
$$Ta=aT, \ T\tau = \tau T, \ \tau a = a^q \tau \eqno{(1.1)}$$
We need also an extension of $\p[T,\tau]$ --- the ring $\p(T)[\tau]$
which is the ring of
non-commutative polynomials in $\tau$ over the field of rational
functions $\p(T)$ with
the same relations (1.1). For a left $\p[T,\tau]$-module $M$ we denote by $M_{\p[T]}$ the
same
$M$ treated as a
$\p[T]$-module with respect to the natural inclusion
$\p[T]\hookrightarrow\p[T,\tau]$.
Analogously, we define $M_{\p[\tau]}$; we shall use similar notations
also for the left
$\p(T)[\tau]$-modules.
\medskip
Obviously we have:
\medskip
{\bf (1.2)} For $C\in M_{r}(\p(T))$ operations $C^t$, $C^{-1}$
and $C^{(i)}$ commute.
\medskip
{\bf Definition 1.3.} ([G], 5.4.2, 5.4.12, 5.4.10). A t-motive $M$ is a left
$\p[T, \tau]$-module which is free and finitely generated as both $\p[T]$-,
$\p[\tau]$-module and such that
$$ \exists m=m(M) \ \hbox{such that}\ (T-\theta)^m M/\tau M=0\eqno{(1.3.1)}$$
\medskip
{\bf Remark.} The above object is called "abelian t-motive" (resp. "t-motive") in [G] (resp. [A]), while the name "t-motive" is used in [G] for a more general object ([G], Definition 5.4.2). Since we shall not use objects defined in [G], 5.4.2, I prefer to use a shorter name for the above $M$.
\medskip
t-motives are main objects of the present paper. If we affirm
that an object exists this means that it exists as a t-motive if otherwise
is not stated. We denote dimension of $M$ over $\p[\tau]$ (resp. $\p[T]$) by $n$ (resp.
$r$), these numbers are called dimension and rank of $M$. Morphisms of abelian
t-motives are morphisms of left $\p[T, \tau]$-modules.

To define a left $\p[T, \tau]$-module $M$ is the same as to define a left $\p[T]$-module
$M_{\p[T]}$
endowed by an action of $\tau$ satisfying $\tau(Pm)=P^{(1)}\tau(m)$, $P\in \p[T]$. In
this situation
we can also treate $\tau$ as a $\p[T]$-linear map $M^{(1)}\to M$. This interpretation is
necessary if
we consider the general case $\w\supset \n F_q[T]$.

We need two categories which are larger than the category of abelian
t-motives.
\medskip
{\bf Definition 1.4.} A pr\'e-t-motive is a left $\p[T,
\tau]$-module which is free
and finitely generated as $\p[T]$-module, and satisfies (1.3.1).
\medskip
{\bf Definition 1.5.} A rational pr\'e-t-motive is a left
$\p(T)[\tau]$-module which is
free and finitely generated as $\p(T)$-module.
\medskip
{\bf Remark 1.6.} An analog of (1.3.1) does not exist for them.
\medskip
There is an obvious functor from the category of t-motives to
the category of
pr\'e-t-motives which is fully faithful, and an obvious functor from the category of
pr\'e-t-motives to the
category of rational pr\'e-t-motives. We denote these functors by
$i_1$, $i_2$
respectively. It is easy to see (Remark 10.2.3) that if $M$ is a
pr\'e-t-motive then the
action of $\tau$ on $i_2(M)$ is invertible.

Let $M_1$, $M_2$ be rational pr\'e-t-motives such
that the action of
$\tau$ on $(M_1)_{\p(T)}$ is invertible.
\medskip
{\bf Definition 1.7.} $\Hom(M_1,M_2)$ is a rational pr\'e-t-motive such that

$$\Hom(M_1,M_2)_{\p(T)}=\Hom_{\p(T)}((M_1)_{\p(T)}, (M_2)_{\p(T)})$$
and the action of $\tau$ is defined by the usual manner: for
$\varphi:M_1 \to M_2$, $m\in M_1$
$$(\tau\varphi)(m)=\tau(\varphi(\tau^{-1}(m)))$$
\medskip
{\bf Definition 1.8.} Let $M$ be a t-motive and
$\mu$ a positive
number. A t-motive
$M'={M'}^{\mu}$ is called the $\mu$-dual of $M$ (dual if $\mu=1$) if
$M'=\Hom(M,\goth
C^{\otimes
\mu})$ as a rational pr\'e-t-motive, i.e. $$i_2\circ
i_1(M')=\Hom(i_2\circ i_1(M),\goth
C^{\otimes
\mu})\eqno{(1.8.1)}$$

{\bf Remark.} This definition generalizes the original one of Taguchi
([T], Section 5), see 1.12 below. A similar definition is in [F].
\medskip
{\bf 1.9.} We shall need the explicit matrix description of the above
objects. Let
$e_*=(e_1, ..., e_n)^t$ be the vector column of elements of a basis
of $M$ over $\p[\tau]$. There exists a matrix $\goth A\in M_n(\p[\tau])$ such that

$$T e_* = \goth A e_*, \ \ \goth A = \sum_{i=0}^l \goth A_i \tau^i \hbox{ where } \goth A_i
\in M_n(\p)\eqno{(1.9.1)}$$
Condition (1.3.1) is equivalent to the condition
$$\goth A_0=\theta E_n + N\eqno{(1.9.2)}$$
where $N$ is a
nilpotent matrix, and the
condition
$m(M)=1$ is equivalent to the condition $N=0$.
\medskip
Let $f_*=(f_1, ..., f_r)^t$ be the vector column of elements of a basis
of $M$ over
$\p[T]$. There exists a matrix $Q=Q(f_*)\in M_r(\p[T])$ such that
$$\tau f_* = Q f_*\eqno{(1.9.3)}$$

{\bf Lemma 1.10.} Let $M$ be as above. A t-motive
$M'$ is the $\mu$-dual of $M$ iff there exists a basis
$f'_*=(f'_1, ..., f'_r)^t$ of $M'$ over $\p[T]$ such that its matrix
$Q'=Q(f'_*)$ satisfies
$$Q'=(T-\theta)^{\mu}Q^{t-1} \eqno{(1.10.1)}$$
$\square$
\medskip
{\bf 1.10.2.} For further applications we shall need the following
lemma. The above
$f_*$, $f'_*$ are the dual bases (i.e. if we consider $f'_i$ as
elements of $\Hom(M,\goth
C)$ then $f'_i(f_j)=\delta^i_j\goth f$, where $\goth f$ is canonically defined by the condition that it generates $\goth C_{\p[T]}$ and satisfies $\tau \goth
f=(T-\theta)\goth f$). Let $\gamma$ be an endomorphism of $M$ and $D$ its
matrix in the basis
$f_*$ (i.e. $\gamma(f_*)=Df_*$). Let $\gamma'$ be the dual endomorphism.
\medskip
{\bf Lemma 1.10.3.} The matrix of $\gamma'$ in the basis $f'_*$ is $D^t$.
$\square$
\medskip
{\bf Remark 1.11.1.} For any $M$ having dual there exists a canonical homomorphism $\delta: \goth
C \to M\otimes M'$.
This is a well-known theorem of linear algebra. Really, in the above notations we have
$\goth f \mapsto \sum_i f_i \times f'_i$. It is obvious that $\delta$ is well-defined,
canonical and compatible
with the action of $\tau$.
\medskip
{\bf Remark 1.11.2.} The $\mu$-dual of $M$ --- if it exists --- is
unique, i.e. does not
depend on base change. This follows immediately from Definition 1.8,
but can be deduced
easily from 1.10.1. Really, let $g_*=(g_1, ..., g_r)^t$ be another basis
of $M$ over
$\p[T]$ and $C\in GL_r(\p[T])$ the matrix of base change (i.e. $g_*= C
f_*$). Then
$Q(g_*)=C^{(1)}QC^{-1}$. Let $g'_*=(g'_1, ..., g'_r)^t$ be a basis of
$M'$ over $\p[T]$
satisfying $g'_*= C^{t-1} f'_*$. Elementary calculation shows that
matrices $Q(g_*)$,
$Q(g'_*)$ satisfy (1.10.1).
\medskip
{\bf Remark 1.11.3.} The operation $M \mapsto {M'}^{\mu}$ is
obviously contravariant functorial. I leave as an exercise to the reader to give an exact definition of the corresponding category such that the functor of duality is defined on it, and is involutive (recall that not all t-motives have duals, and the dual of a map of t-motives is a priori a map of rational pr\'e-t-motives).
\medskip
{\bf 1.12.} The original definition of duality ([T], Definition 4.1; Theorem 5.1) from
the first
sight seems to be more
restrictive than the definition 1.8 of the present paper, but really
they are
equivalent. We recall some notations and definitions of [T] in a slightly less
general setting (rough statements; see [T] for the exact statements). Let $G$ be a finite affine group scheme over $\p$, i.e. $G=\Spec R$
where $R$ is a finite-dimensional $\p$-algebra. Let $\mu:R \to R\otimes R$ be the
comultiplication of $R$. Such group $G$ is called a finite $v$-module ([T], Definition
3.1) if there is a homomorphism $\psi: \w \to \End_{gr. \ sch.}(G)$ satisfying some natural conditions (for
example, an analog of 1.3.1). Further, let $\Cal E_G$ be a $\p$-subspace of $R$ defined as
follows: $$\Cal E_G=\{x\in R \ \vert \ \mu(x)=x\otimes 1+1\otimes x\}$$
The map $x\mapsto x^q$ is a $\p$-linear map $\fr: \Cal E_G^{(1)}\to \Cal E_G$. Further,
the map $\psi(T): G \to G$ can be defined on $\Cal E_G$. Let $v: \Cal E_G\to\Cal
E_G^{(1)}$ be a map satisfying $\fr \circ v= \psi(T)-\theta$.

We consider two finite $v$-modules $G$, $H$, the above objects $\fr$, $v$ etc. will carry
the respective subscript. Let * be the dual in the meaning of linear algebra.
\medskip
{\bf Definition 1.12.1} ([T], 4.1). Two finite $v$-modules $G$, $H$
are called dual if there exists an
isomorphism $\alpha: \Cal E^*_H\to\Cal E_G$ such that if we denote by
$\goth v: \Cal E_G\to \Cal E^{(1)}_G$ a map which enters in the commutative diagram
$$ \matrix \Cal E^*_H & \overset{\fr^*_H}\to{\longrightarrow} & \Cal E^{*(1)}_H \\ & & \\
\alpha \downarrow & & \alpha ^{(1)}\downarrow \\ & & \\ \Cal E_G & \overset{\goth
v}\to{\longrightarrow} & \Cal E^{(1)}_G \endmatrix $$
then we have:
$$\fr_G\circ \goth v= \psi_G(T)-\theta \eqno{(1.12.2)}$$
i.e. $\goth v=v_G$.
\medskip
Let $M$ be a t-motive having $m(M)=1$, $E=E(M)$ the corresponding t-module and
$a\in \hbox{\bf
A}$. We denote $E_a$ --- the set of $a$-torsion elements of $E$ --- by $M_a$. It is a
finite
$v$-module.
\medskip
{\bf Proposition 1.12.3.} Let $M$, $M'$ be t-motives which are
dual in the meaning
of Definition 1.8. Then $\forall a\in \hbox{\bf A}$, $a\ne 0$ we have: $M_a$,
$M'_a$ are dual in the meaning of 1.12.1 = [T], Definition 4.1.
\medskip
{\bf Proof.} Condition $a\in \n F_q[T]$ implies that multiplication by $\tau$ is
well-defined on $M/aM$.
\medskip
{\bf Lemma 1.12.3.1.} We have canonical isomorphisms $i: M/aM \to
\Cal E_{M_a}$,
$i^{(1)}: M/aM \to \Cal E^{(1)}_{M_a}$ such that the following
diagrams are commutative:
$$ \matrix M/aM & \overset{\tau}\to{\longrightarrow} & M/aM &&& M/aM &
\overset{T}\to{\longrightarrow} & M/aM
\\ & & &&&&&
\\ i^{(1)}\downarrow & & i\downarrow &&& i\downarrow & & i\downarrow \\ & & &&&&&
\\ \Cal E^{(1)}_{M_a} & \overset{fr}\to{\longrightarrow} & \Cal E_{M_a} &&&
\Cal E_{M_a} & \overset{\psi_T}\to{\longrightarrow} & \Cal E_{M_a}
\endmatrix $$
\medskip
{\bf Proof.} Let $R$ be a ring such that $\Spec R=M_a$. The pairing between $M$ and $E$
shows that there exists a map $M\to R$ which is obviously factorized via an inclusion
$M/aM\to R$. It is easy to see that the image of this inclusion is contained in $\Cal
E_{M_a}$, i.e. we get $i$. Since $\dim_{\p}(M/aM)=\deg a \cdot r(M)$ and
$\dim_{\p}(R)=q^{\deg a \cdot r(M)}$ we get from [T], Definition 1.3 that $i$ is an
isomorphism. Other statements of the lemma are obvious. $\square$
\medskip
This lemma means that we can rewrite Definition 1.12.1 for the case $G=M_a$, $H=N_a$ by
the following way:\footnotemark \footnotetext{Here and below a t-motive $N$ should not be confused with $N$ of 1.9.2.}
\medskip
{\bf 1.12.3.2.} Two finite $v$-modules $M_a$, $N_a$ are dual if there
exists an
isomorphism $\alpha: (N/aN)^* \to M/aM $ such that after identification
via $\alpha$ of
$\tau^*: (N/aN)^* \to (N/aN)^*$ with a map $\goth v: M/aM \to M/aM$ we
have on $M/aM$:
$$\tau\circ \goth v= t-\theta \eqno{(1.12.3.3)}$$
We need a
\medskip
{\bf Lemma 1.12.3.4.} For $i=1,2$ let $N_i$ be a free $\p[T]$-module of
dimension $r$ with a
base $f_{i*} =(f_{i1}, ... , f_{ir} )$, let $\varphi_i:N_i \to N_i$ be
$\p[T]$-linear
maps having matrices $\goth Q_i$ in $f_{i*} $ such that $\goth
Q_2=\goth Q_1^t$, and let
$a$ be as above. Let, further, $\varphi_{i,a}: N_i/aN_i \to N_i/aN_i $
be the natural
quotient of $\varphi_i$. Then there exist $\p$-bases $\tilde f_{i*} $
of $N_i/aN_i $ such
that the matrix of $\varphi_{1,a}$ in the base $\tilde f_{1*} $ is
transposed to the
matrix of $\varphi_{2,a}$ in the base $\tilde f_{2*} $.
\medskip
{\bf Proof.} We can identify elements of $N_2$ with $\p[T]$-linear
forms on $N_1$
(notation: for $x\in N_2$ the corresponding form is denoted by
$\chi_x$) such that
$\chi_{\varphi_2(x)}=\chi_x\circ\varphi_1$. Any $\p[T]$-linear form
$\chi$ on $N_i$
defines a $\p[T]/a\p[T]$-linear form on $N_i/aN_i$ which is denoted by
$\chi_a$. Let now
$x\in N_2/aN_2$, $\bar x$ its lift on $N_2$, then
$\chi_{x,a}=(\chi_{\bar x})_a$ is a
well-defined $\p[T]/a\p[T]$-linear form on $N_1/aN_1$. For $x\in
N_2/aN_2$ we have
$$\chi_{\vf_{2,a}(x),a}= \chi_{x,a}\circ \vf_{1,a}$$
Further, let $\lambda: \p[T] \to \p$ be a $\p$-linear map such that
\medskip
{\bf 1.12.3.5.} Its kernel does not contain any non-zero ideal of
$\p[T]/a\p[T]$.
\medskip
(such $\lambda$ obviously exist.) For $x\in N_2/aN_2$ we denote $\lambda \circ
\chi_{x,a}$ by $\psi_x$, it is a $\p$-linear form on $\p$-vector space
$N_1/aN_1$.
Obviously condition (1.12.3.5) implies that the map $x\mapsto \psi_x$
is an isomorphism
from $N_2/aN_2$ to the space of $\p$-linear forms on $\p$-vector space
$N_1/aN_1$, and we
have
$$\psi_{\vf_{2,a}(x)}=\psi_x\circ\vf_{1,a}$$
which is equivalent to the statement of the lemma. $\square$
\medskip
Finally, the proposition follows immediately from this lemma
multiplied by $T-\theta$,
formula 1.10.1 and 1.12.3.2. $\square$
\medskip
{\bf Remark.} Let $a=\sum_{i=0}^k g_iT^i$, $g_i\in \n F_q$, $g_k=1$.
Taguchi ([T], proof
of 5.1 (iv)) uses the following $\lambda$: $\lambda(T^j)=0$ for $j<k-1$,
$\lambda(T^{k-1})=1$. It is easy to check that for $x=(T^i+
T^{i-1}g_{k-1}+T^{i-2}g_{k-2}+ ...
+g_{k-i})f_{2j}$ for this $\lambda$ we have: $\psi_x(T^i f_{1j})=1$,
$\psi_x(T^{i'}
f_{1j'})=0$ for other $i'$, $j'$.
\medskip
{\bf 1.13.} We consider in Sections 1.13, 1.14 the case of arbitrary $\w\supset \n F_q[T]$.
\medskip
A t-motive over {\bf A} is defined for example in [BH], p.1. Let us reproduce this
definition for the case of characteristic 0. Let $J$ be an ideal of
$\w_C$ generated by the elements $a\otimes 1 - 1 \otimes \iota(a)$ for all $a\in \w$. The ring $\w_C[\tau]$ is defined by the formula $\tau\cdot (a\otimes z)=(a\otimes z^q)\cdot \tau$, $a\in \w$, $z\in \p$.
\medskip
{\bf Definition 1.13.1.} A t-motive $M$ over {\bf A} is a pair $(M, \tau)$
where $M$ is
a locally free $\hbox{\bf A}_C$-module and $\tau$ is an $\hbox{\bf A}_C$-linear map
$M^{(1)}\to M$ satisfying the following analog of 1.3.1, 1.9.2:
$$\exists m \hbox{ such that } J^m(M/\tau (M^{(1)}))=0\eqno{(1.13.2)}$$

{\bf Remark 1.13.3.} We can consider $M$ as an $\w_C[\tau]$-module using the following formula for the product $\tau \cdot m$:
$$\tau \cdot m = \tau(m\otimes 1)$$
where $m\in M$, $m\otimes 1\in M^{(1)}$.
\medskip
The rank of $M$ as a locally free $\hbox{\bf A}_C$-module is called the rank of the
corresponding t-motive $(M, \tau)$. If $\w=\n F_q[T]$ then $M^{(1)}$ is isomorphic to
$M$, we can
consider $M$ as a $\p[T,\tau]$-module, and it is possible to show that in this case
1.13.2 implies
that $M_{\p[\tau]}$ is a free $\p[\tau]$-module. In the general case, the dimension $n$
of $(M, \tau)$
is defined as $\dim_{\p}(M/\tau (M^{(1)}))$.

Let us fix $\goth C=(\goth C, \tau_{\goth C})$ --- a t-motive of rank 1 over
{\bf A}.
For a t-motive $M=(M,\tau_M)$ a t-motive ${M_\goth C'}$ --- the
$\goth C$-dual of $M$
--- is defined as follows. We put
$M_\goth C'=\Hom_{\w_C}(M,\goth C)$.
Since for any locally free $\hbox{\bf A}_C$-modules $M_1$, $M_2$ we have
$$\Hom_{\w_C}(M_1,M_2)^{(1)}=\Hom_{\w_C}(M_1^{(1)},M_2^{(1)})$$
we can define $\tau(M_\goth C')$ by the following formula:

$$\hbox{ For } \vf\in \Hom_{\w_C}(M,\goth C)^{(1)} \hbox{ we have }
\tau(M_\goth C')(\vf)=\tau_{\goth C}\circ \vf \circ \tau_M^{-1}$$
\medskip
{\bf 1.14. Duality for abelian $\tau$-sheaves.} We use notations of [BH],
Definition 2.1 if they do not differ from the notations of the present
paper; otherwise we continue to use notations of the present paper
(for example, $d$ (resp. $\sigma^*(\goth X)$ for any object $\goth X$)
of [BH] is
$n$ (resp. $\goth X^{(1)}$) of the present paper). For any abelian
$\tau$-sheaf $\underline{\Cal F}$ we denote its $\Pi_i$, $\tau_i$ by
$\Pi_i(\underline{\Cal F})$, $\tau_i(\underline{\Cal F})$
respectively. If $M$, $N$ are invertible sheaves on $X$ and $\rho: M
\to N$ a rational map then we denote by $\rho^{inv}: N \to M$ the
rational map which is inverse to $\rho$ with respect to the composition.
We define
$\tau_{\goth r,i-1}(\underline{\Cal F})$ (the rational $\tau_i$) as
the composition map $\tau_{i-1}(\underline{\Cal F})
\circ {\Pi_{i-1}^{(1)}}^{inv}(\underline{\Cal F})$, it is a rational
map from $\Cal F_i^{(1)}$ to $\Cal F_i$.

Let $\underline{\Cal O}$ be a fixed abelian
$\tau$-sheaf having $r=n=1$. The $\underline{\Cal O}$-dual abelian
$\tau$-sheaf $\underline{\Cal F}'= \underline{\Cal
F}'_{\underline{\Cal O}}$ is defined by the formulas
$$\Cal F'_0=\Hom_{X}(\Cal F_0, \Cal O_0)$$ where Hom is the sheaf's
one, and the map $\tau_{\goth r,-1}(\underline{\Cal F}'): {\Cal
F'_0}^{(1)}\to \Cal F'_0$ is defined as follows. We have ${\Cal
F'_0}^{(1)}= \Hom_{X}({\Cal F_0}^{(1)}, {\Cal O_0}^{(1)})$. Let
$\gamma\in\Hom_{X}({\Cal F_0}^{(1)}, {\Cal O_0}^{(1)})(U)$ where $U$
is a sufficiently small affine subset of $X_{\p}$, such that $\gamma:
{\Cal F_0}^{(1)}(U) \rightarrow {\Cal O_0}^{(1)}(U)$.
\medskip
{\bf 1.14.1.} We define: $[[\tau_{\goth r,-1}(\underline{\Cal
F}')](U)](\gamma)$ is the following composition map:
$$\Cal F_0(U)\overset{[\tau_{\goth r,-1}^{inv}(\underline{\Cal
F})](U)}\to{\longrightarrow}
\Cal F_0^{(1)}(U)\overset{\gamma}\to{\to}\Cal O_0^{(1)}(U)
\overset{[\tau_{\goth r,-1}(\underline{\Cal
O})](U)}\to{\longrightarrow}\Cal O_0(U)\in \Hom_{X}({\Cal F_0}, {\Cal
O_0})(U)$$

Clearly that this definition and the definitions 1.8, 1.13 are
compatible with the forgetting functor $\underline{M}(\underline{\Cal
F})$ from abelian $\tau$-sheaves to pure Anderson t-motives of [BH],
Section 3, page 8.
\medskip
{\bf 1.15. Duality over fields.} Let $L\supset \n F_q(\theta)$ be a field extension of $\n F_q(\theta)$, and $M$ a t-motive over $L$ (i.e. a pair ($M$, an $L$-structure on $M$)). Obviously we have
\medskip
{\bf Proposition 1.15.1.} The notion of duality for $M$ over $L$ is well-defined. $\square$
\medskip
Similarly, we have a proposition for Galois action:
\medskip
{\bf Proposition 1.15.2.} Let $M$ be defined over $\overline{ \n F_q(\theta)}$ and $\gamma\in \Gal( \n F_q(\theta))$. Then $(\gamma(M))'=\gamma(M')$. $\square$
\medskip
\medskip
{\bf 1.16. Case of finite characteristic.} Let $\iota: \w\to \bar \n F_q$ be a map of finite characteristic, we denote $\Ker \iota$ by $\Cal P$. The definition of t-motive for this case is similar to 1.3, see [G] for the details. The definition of duality also is similar to the one of the case of generic characteristic. Duality commutes with reduction. Namely, let $M$ be from 1.15, $\goth P$ a prime of $L$ not over the infinity of $\n F_q(\theta)$, $\Cal P\subset \w$ is $\iota^{-1}(\goth P\cap \n F_q[\theta]$) --- the finite characteristic. We consider the case of good reduction of $M$ at $\goth P$, we denote it by $\tilde M$. It is a t-motive in characteristic $\Cal P$.  Let $M$ have dual $M'$.
\medskip
{\bf Proposition 1.16.1.} $\tilde M$ has dual iff $M'$ has good reduction at $\goth P$; in this case they coincide. $\square$
\medskip
{\bf Remark 1.16.2.} Apparently if $M$ has good reduction and dual, then $M'$ also has good reduction (in this case 1.16.1 means that $M'$ exists implies $(\tilde M)'$ exists). For standard-3 t-motives (this is a simple tipe of t-motives, see 11.8.1) apparently this can be shown by explicit calculations.
\medskip
{\bf Remark 1.16.3.} Clearly 1.16.1 is true for the case of bad reductions. I do not give exact definitions for this case.
\medskip
{\bf 1.16.4. Ordinarity.} Let $M$ be of finite characteristic. By analogy with the number field case, $M$ is called ordinary if its Newton polygon consists of 2 segments. If $N=0$ then the Newton polygon of $M'$ is the dual of the one of $M$ (the notion of duality of polygons is clear; apparently the condition $N=0$ can be omitted). So, we have
\medskip
{\bf Proposition 1.16.5.} $M$ is ordinary $\iff M'$ is ordinary. $\square$
\medskip
See 13.4.1 for a more exact result.
\medskip
{\bf 2. Analytic duality.}
\medskip
We consider in the present section the case of arbitrary $\w\supset \n F_q[T]$ (and $N=0$ as usually).
\medskip
Condition $N=0$ implies that an element $a\in \w$ acts on $\Lie(M)$ by multiplication by $\iota(a)$. Hence, we have a
\medskip
{\bf Definition 2.1.} Let $V$ be the space $\p^n$. A locally free $r$-dimensional
$\n Z_\infty$-submodule
$L$ of $V$ is called a lattice if

(a) $L$ generates $V$ as a $\p$-module and

(b) The $\n R_\infty$-linear envelope of $L$ has dimension $r$ over
$\n R_\infty$.
\medskip
Numbers $n$, $r$ are called the dimension and the rank of $L$
respectively. Attached to $(L,V)$ is the tautological inclusion $\vf=\vf(L,V):
L \to V$. We shall consider the category of triples $(\vf, L, V)$; a map $\psi:
(\vf, L,V)\to(\vf_1, L_1,V_1)$ is a pair $(\psi_L, \psi_V)$ where $\psi_L: L \to L_1$ is a $\n Z_\infty$-linear map, $\psi_V: V\to V_1$ is a $\p$-linear map such that $\vf_1\circ \psi_L=\psi_V \circ \vf$.

Inclusion $\vf$ can be extended to a map
$L\underset{\n Z_\infty}\to{\otimes}\p \to V$ (which is surjective
by 2.1 a), we denote it by $\vf=\vf(L,V)$ as well.
We can also attach to $(L,V)$ an exact sequence
$$0\to \Ker \vf \to L\underset{\n Z_\infty} \to{\otimes}\p \overset{\vf}\to{\to} V \to
0\eqno{(2.2)}$$

Let $\Cal I\in \Cl(\w)$ be a class of ideals; we shall use the same notation $\Cal I$ to
denote a representative in the $\iota$-image of this class. Let $(\vf', L', V')$ be another lattice and $D$ a structure of a perfect $\Cal I$-pairing $<* , * >_D$ between $L$
and $L'$. Let us fix an isomorphism $$\alpha: \Cal I \underset{\n Z_\infty} \to{\otimes} \p \to \p\eqno{(2.2')}$$ $D$ extends via $\alpha$ to a perfect $\p$-pairing between $L\underset{\n Z_\infty}\to{\otimes}\p$
and $L'\underset{\n Z_\infty}\to{\otimes}\p$, we denote this pairing by $D_{\alpha, \infty}$.
\medskip
{\bf Definition 2.3.} Two lattices $(\vf, L, V)$ and $(\vf', L', V')$ are called $(\alpha, \Cal I)$-dual if there exists a perfect $\Cal I$-pairing $D$ between $L$
and $L'$ such that $\Ker \vf \subset L\underset{\n Z_\infty}\to{\otimes}\p$, $\Ker \vf' \subset L'\underset{\n Z_\infty}\to{\otimes}\p$ are mutually orthogonal with respect to $D_{\alpha, \infty}$.
\medskip
Let $(n,r)$, $(n',r')$ be the dimension and rank of $(\vf, L, V)$ and $(\vf', L', V')$ respectively. If they are $(\alpha, \Cal I)$-dual then $r'=r$, $n'=r-n$. There exists the following reformulation of the definition of duality. $D_{\alpha, \infty}$ induces an isomorphism
$\gamma_{\alpha, D}: (L\underset{\n Z_\infty}\to{\otimes}\p)^* \to
L'\underset{\n Z_\infty}\to{\otimes}\p$ (here and below for any object $W$ we
denote $W^*=\Hom_{\p}(W,\p)$ ).
\medskip
{\bf Property 2.4.} $(\vf, L, V)$ and $(\vf', L', V')$ are $(\alpha, \Cal I)$-dual iff there exists an isomorphism from $(\Ker \vf)^*$
to $V'$ making the following diagram commutative: $$\matrix 0 & \to &
V^* & \overset{\vf^*}\to{\to} &  (L\underset{\n Z_\infty}
\to{\otimes}\p)^* & \to & (\Ker \vf)^* & \to & 0\\
&&&&&&&& \\
& & \downarrow & & \gamma_{\alpha, D} \downarrow && \downarrow \\ &&&&&&&& \\ 0
& \to & \Ker \vf' & \rightarrow  & L'\underset{\n Z_\infty}\to{\otimes}\p &
\overset{\vf'}\to{\to} & V' & \to &
0\endmatrix \eqno{(2.5)}$$

Further, this property is equivalent to the following two conditions:
\medskip
{\bf 2.6.} $\dim V'=r-n$;
\medskip
{\bf 2.7.} The composition map $\vf'\circ \gamma_D \circ \vf^*:
V^* \to V'$ is 0.
\medskip
Both 2.4 and (2.6, 2.7) are obvious.
\medskip
{\bf Remark 2.8.} It is
easy to see that the functor $(\vf, L, V) \mapsto (\vf', L', V')$ is well-defined
on a subcategory
(not all lattices have duals, see
below) of the category of the triples $(\vf, L, V)$, it is contravariant and
involutive.
\medskip
{\bf 3. Explicit formulas for analytic duality.}
\medskip
Here we consider the case $\w=\n
F_q[T]$. In this case $\Cl(\w)=0$, and $(\alpha, \Cal I)$-dual is called simply
dual. The coordinate
description of the dual lattice is the following. Let
$e_1, ..., e_r$ be a $\n Z_\infty$-basis of
$L$ such that $\vf(e_1), ..., \vf(e_n)$ form a $\p$-basis of $V$. Like in the
theory of abelian
varieties, we denote by $Z=(z_{ij})$ the Siegel matrix whose lines are
coordinates of $\vf(e_{n+1}), ..., \vf(e_r)$ in the basis $\vf(e_1), ..., \vf(e_n)$, more
exactly, the size of $Z$ is $(r-n)\times n$ and
$$\forall i =1,..., r-n \ \ \ \ \vf(e_{n+i})=\sum_{j=1}^n z_{ij}\vf(e_j)\eqno{(3.1)}$$ $Z$
defines $L$, we denote $L$ by $\goth L(Z)$.
\medskip
{\bf Proposition 3.2.} $[\goth L(Z)]'=\goth L(-Z^t)$, i.e. a Siegel matrix of the dual
lattice is the minus transposed Siegel matrix.
\medskip
{\bf Proof.} Follows immediately from the definitions. Really, let $f_1, ..., f_r$ be a
basis of $L'$, we define the pairing by the formula $$<e_i,
f_j>=\delta_i^j\eqno{(3.3)}$$ and the map $\vf'$ by the formula
$$\forall i =1,... ,n \ \ \ \ \vf'(f_{i})=\sum_{j=1}^{r-n} -z_{ji}\vf'(f_{n+j})$$
(minus transposed Siegel matrix).
$\Ker \vf$ is generated by elements $$v_i=e_{n+i}-\sum_{j=1}^n z_{ij}e_j, \ \ \ \ i
=1,..., r-n$$ and $\Ker \vf'$ is generated by elements $$w_i=f_{i}+\sum_{j=1}^{r-n}
z_{ji}f_{n+j}, \ \ \ \ i =1,..., n\eqno{(3.4)}$$ It is sufficient to check that $\forall i,j$ we have
$<v_i, w_j>=0$; this follows immediately from 3.3. $\square$
\medskip
{\bf Remark 3.5.} $L'$ exists not for all $L$.
Trivial counterexample: case $n=r=1$. To get another counterexamples, we use that for
$n=1$ (lattices of Drinfeld modules) a Siegel matrix is a column matrix $Z=\left(\matrix
z_1&...&z_{r-1} \endmatrix \right)^t$ and
$$ \goth L(Z) \hbox{ is not a lattice } \iff 1,z_1, ... ,z_{r-1} \hbox{ are linearly
dependent over } \n R_\infty \eqno{(3.6)}$$ while for $n=r-1$ a Siegel matrix is a
row matrix $Z=\left(\matrix -z_1&...&-z_{r-1} \endmatrix \right)$ and
$$ \goth L(Z) \hbox{ is not a lattice } \iff \forall i \ \ z_i\in \n R_\infty
\eqno{(3.7)}$$ Since condition (3.7) is strictly stronger than (3.6) we see
that all lattices having $n=1$, $r>1$ have duals while not all lattices having $n=r-1$,
$r>2$ have duals.

It is clear that almost all matrices have duals. Here "almost all" has the same meaning
that as "Almost all matrices $Z$ are a Siegel matrice of a lattice", i.e. if we choose an
(infinite) basis of $\p/\n R_\infty$, then coordinates of the entries of $Z$ in this
basis must satisfy some polynomial relations in order that $Z$ is not a Siegel matrice of
a lattice.
\medskip
{\bf Remark 3.8.} The coordinate proof of the theorem that the notion
of the dual lattice is well-defined, is the following. Two Siegel matrices
$Z$, $Z_1$ are called equivalent iff there exists an isomorphism of their pairs $(\goth
L(Z), V)$, $(\goth L(Z_1),V_1)$. Like in the classical theory of modular forms, $Z$,
$Z_1$ are equivalent iff there exists a matrix $\gamma \in GL_r(\n Z_\infty)=\left(\matrix A&B\\ C&D \endmatrix \right)$ ($A,B,C,D$ are the ($n\times
n$), ($n\times r-n$), ($r-n\times n$), ($r-n\times r-n$)-blocks of $\gamma$ respectively;
we shall call this block structire by the $(n, r-n)$-block structure) such that
$$C+DZ=Z_1(A+BZ)\eqno{(3.8.1)}$$

Let $A_1,B_1, C_1, D_1$ be the $(n, r-n)$-block structure of the matrix $\gamma^{-1}$.
The equality
$$-C_1^t+A_1^tZ^t={Z_1}^t(D_1^t-B_1^tZ^t)\eqno{(3.8.2)}$$
shows that if $Z$, $Z_1$ are equivalent
then $-Z^t$, $-Z_1^t$ are equivalent. [Proof of (3.8.2): (3.8.1) implies $Z_1=(C+DZ)(A+BZ)^{-1}$; substituting this value of $Z_1$ to the transposed (3.8.2), we get $-C_1+ZA_1=(D_1-ZB_1)(C+DZ)(A+BZ)^{-1}$, or $(-C_1+ZA_1)(A+BZ)=(D_1-ZB_1)(C+DZ)$. This formula follows immediately from $\left(\matrix A_1&B_1\\C_1&D_1\endmatrix \right)\left(\matrix A&B\\C&D\endmatrix \right)=\left(\matrix E_n&0\\0&E_{r-n}\endmatrix \right)$].  

Further, let $\alpha: (L_1\subset \p^n) \to
(L_2\subset \p^n)$ be a map of lattices. If $L'_1$, $L'_2$
exist, then the map
$\alpha': (L'_2\subset \p^{r-n}) \to (L'_1\subset \p^{r-n})$ is defined by the
following formulas. Let $Z_i$ be the Siegel matrices of $L_i$ in the
bases
$e_{i1}, ... e_{ir}$ of $L_i$ ($i=1,2$). Let us consider the matrix
$\goth M=(m_{ij})\in
M_{r}(\n Z_\infty)$ of $\alpha$ in the bases $e_{i1}, ...,
e_{ir}$ (i.e.
$\alpha(e_{1i})=\sum_j m_{ij} e_{2j}$). Let $f_{i1}, ..., f_{ir}$ be the dual
base of $L'_i$ (see 3.3) and $e'_{i1}, ... e'_{ir}$ another base of $L'_i$ defined by $$e'_{ij}=f_{i,j+n}, \ \ \ \ j+n \mod r\eqno{(3.8.3)}$$ Formulas (3.8.3), (3.4) show that an analog of 3.1 is satisfied for both bases $e'_{i1}, ..., e'_{ir}$, their Siegel matrices are $-Z_i^t$.

Let
$$\goth M=\left(\matrix \goth M_{11} & \goth M_{12} \\ \goth M_{21} &
\goth M_{22}
\endmatrix \right)$$
be the $(n, r-n)$-block structure of $\goth M$. The matrix of $\alpha'$ in the bases $f_{i1}, ..., f_{ir}$ is $\goth M^t$, and using the matrix 3.8.3 of change of base, we get that $\goth M'$ --- the matrix of $\alpha'$ in
the bases $e'_{i1}, ..., e'_{ir}$ --- has the following $(r-n,
n)$-block structure:
$$\goth M'=\left(\matrix \goth M_{22}^t & \goth M_{12}^t \\
\goth M_{21}^t & \goth
M_{11}^t \endmatrix \right)\eqno{(3.8.4)} $$
The property that $\goth M$ comes from a $\p$-linear map $\p^n \to \p^n$ implies
that $\goth M'$ comes from a $\p$-linear map $\p^{r-n} \to
\p^{r-n}$. This follows immediately from the definition of dual lattice, or can be easily checked algebraically.
\medskip
{\bf Remark 3.9.} Taking $\gamma =\left(\matrix 1&0\\ 0&-1 \endmatrix \right)$ we get
that $Z$ is equivalent to $-Z$, hence $Z'$ is also a Siegel matrix of the dual lattice.
\medskip
{\bf 4. Main conjecture for arbitrary $\w$}.
\medskip
The main result of the paper is the following Theorem 5 on coincidence of algebraic and
analytic duality. We formulate it as a conjecture 4.1 for any $\w$, but we prove it only for the case $\w=\n
F_q[T]$.  Let $M$ be a uniformizable t-motive. Its lattice $L(M)$ is really a lattice in
the meaning of Definition 2.1, because [A],
Corollary 3.3.6 (resp. [G], Lemma 5.9.12) means that it satisfies 2.1a (resp. 2.1b);
recall that we consider the case $N=0$,
i.e. the action of $T$ on $\Lie(M)$ is simply multiplication by $\theta$.
Let us fix (like in 1.13) $\goth C=(\goth C, \tau_{\goth C})$ --- a t-motive of rank 1 over
{\bf A}, and let $L(\goth C)$ be its lattice. It is a $\n Z_\infty$-module. $\Omega=\Omega(\w)$ is an $\w$-module, we consider a $\n Z_\infty$-module $\iota^{-1}(\Omega)$. There exists the   notion of the $L(\goth C)\otimes \iota^{-1}(\Omega)$-duality.
\medskip
{\bf Conjecture 4.1.} Let $M$ be a uniformizable t-motive having $N=0$ such
that its $\goth C$-dual $M'$ exists. Then $M'$ is uniformizable, it has $N':=N(M')=0$, and $(L(M), \Lie(M))$ and $(L(M'), \Lie(M'))$ are $\alpha, L(\goth C)\otimes \iota^{-1}(\Omega)$-dual for some $\alpha$ from $2.2'$ (it can be explicitly described).
\medskip
We prove in Section 5 the first step of the proof of this conjecture.
\medskip
{\bf Remark 4.2.} It is possible to generalize the above conjecture to the case of non-uniformizable $M$, $M'$. The pairing is defined between $\Hom_{\w_C[\tau]}(M,Z_1)$ and $\Hom_{\w_C[\tau]}(M',Z_1)$ (see (5.2.1a) for the definition of $Z_1$), or, the same, between $M_a$ and $M'_a$ for any $a\in \w$ (see 5.1.6).
\medskip
{\bf 5. Main theorem.}
\medskip
Recall that the word "canonical" means "canonical up to multiplication by elements of $\n F_q^*$".
\medskip
{\bf Theorem 5.}\footnotemark \footnotetext{The proof of this theorem was inspired by a result of Anderson, see Section 6 for details.} Let $M$ be a uniformizable t-motive over $\w=\n F_q[T]$ having $N=0$ such
that its dual $M'$ exists and has $N':=N(M')=0$. Then $M'$ is uniformizable, and $(L(M), \Lie(M))$ and $(L(M'), \Lie(M'))$ are dual.
\medskip
{\bf Remark 5A.} Condition $N'=0$ holds for pure $M$ (Theorem 10.3) and for a large class of non-pure $M$ (Theorem 11.5). Most likely, a modification of the end of the proof of the present theorem will permit us to prove that $N'=0$ holds for all $M$ having $N=0$ and having dual.
\medskip
{\bf Remark 5B.} A reformulation of this theorem in terms of Pink-Hodge structures, as well as its statement for $N\ne0$ (without proof), are given in Section 9 (Result 9.3 and Proposition 9.4 respectively).
\medskip
{\bf Corollary 5.1.1.} If $\w=\n F_q[T]$ then a Siegel matrix of $M'$ is the minus transposed
of a Siegel matrix of $M$.
\medskip
In the section 8 below we give a corollary of this theorem and some conjectures related to the
problem of 1 -- 1 correspondence between t-motives and lattices.
\medskip
{\bf 5.1.2. Some definitions.} Recall that $E=E(M)$ is isomorphic to $\p^n$. There is a structure of $\w$-module on $E$; multiplication by $T$ is denoted by $m_T$, and this operator $m_T$ is defined in coordinates by the formula
$$m_T(x)=\sum_{i=0}^l\goth A_ix^{(i)}$$ where $x\in E=\p^n$ is a vector column, $\goth A_i$ are from 1.9.1. There is a map $\exp: \Lie(M) \to E$ making the following diagram commutative:
$$\matrix \Lie(M) &
\overset{\Exp}\to{\to} & E \\ \\ \theta \downarrow & & m_T
\downarrow \\ \\ \Lie(M) & \overset{\Exp}\to{\to} & E \endmatrix \eqno{(5.1.3)}$$
By definition, $L(M)=\Ker \Exp$.

We need another space $\Lie_T(M)$ together with an isomorphism $\goth a:\Lie_T(M) \to \Lie(M)$ and a structure of $\w$-module on $\Lie_T(M)$ such that the multiplication by $T$ on $\Lie_T(M)$ is simply the multiplication by $\theta$ on $\Lie(M)$, i.e. $$\goth a(Tx)=\theta\cdot(\goth a(x))\eqno{(5.1.4)}$$ where $x\in\Lie_T(M)$. Commutativity of 5.1.3 means that $\Exp\circ\goth a:\Lie_T(M)\to E$ is a map of $\w$-modules.
\medskip
{\bf 5.1.5.} We shall work merely with $L_T(M):=\Ker (\Exp\circ\goth a)\subset \Lie_T(M)$ rather than $L(M)$. Clearly $L_T(M)$ is an $\w$-module, $\goth a:L_T(M)\to L(M)$ is an isomorphism satisfying 5.1.4 for $x\in L_T(M)$.
\medskip
The proof of Theorem 5 consists of two steps. We formulate and prove Step 1 for the case of arbitrary $\w$.
\medskip
{\bf Step 1.} For the above $M$, $M'$ we have:
\medskip
(A) Uniformizability of $M$ implies uniformizability of $M'$.
\medskip
(B) There exists a canonical $\w$-linear $L_T(\goth C)\otimes
\Omega$-valued perfect pairing $<* , * >_M$
between $L_T(M)$ and $L_T(M')$ (by 5.1.5, this is the same as the $\n Z_\infty$-linear pairing between $L(M)$ and $L(M')$, which, in its turn, is $D$ of Definition 2.3). It is functorial.
\medskip
{\bf Remark 5.1.6.} Practically, (B) comes from [T], Theorem 4.3 (case
$\w=\n F_q[T]$). Really, to define a pairing between $L(M)$ and
$L(M')$ it is sufficient to define (concordant) pairings between $L(M)
/aL(M)$ and $L(M') /aL(M')$ for any $a\in \w$. Since $M_a:=E(M)_a=L(M)/aL(M)$ and because of Proposition 1.12.3 which affirms that
$M_a$ and $M'_a$ are Taguchi-dual, we see that [T], Theorem 4.3 gives
exactly the desired pairing.
\medskip
We give two versions of the proof of Step 1: the first one --- for the general case of arbitrary $\w$ and the second one --- for the case $\w=\n F_q[T]$ --- is based on explicit calculations, it is used for the proof of Step 2.
\medskip
{\bf 5.2. Proof: Step 1, Version 1.} Here we consider the general case of arbitrary $\w$.
Let $\Omega=\Omega(\w/\n
F_q)$ be the module of differential forms; we can consider it as an
element of $\Cl(\w)$. We use formulas and notations of [G], Section 5.9
modifying them to the case of arbitrary $\w$. For example, {\bf
A} (resp. {\bf K}) of [G], 5.9.16 is {\bf A} (resp. {\bf K})
of the present paper (recall that $\bar K$ (resp. $\bar K[T,\tau]$) of [G] is
$\p$ (resp. $\w_C[\tau]$, see 1.13) of the present paper). Hence, we denote $\bar K\{T\}$ of
[G], Definition 5.9.10 by $\p\{T\}$. For the general case it must be
replaced by a ring $Z_0$ defined by the formula
$$Z_0:=\w\underset{\n F_q[T]}\to{\otimes}\p\{T\}\eqno{(5.2.1)}$$
$Z_0$ is a $\w_C[\tau]$-module, i.e. $\tau$ acts on $Z_0$, and $Z_0^\tau=\w$.

$Z_1$ for the present case is defined by the same formula [G], 5.9.22. Explicitly,
$$Z_1:=\Hom^{cont}_{\w}(\x/\w,\p)\eqno{(5.2.1a)}$$
It is a locally free $Z_0$-module of dimension 1 (the module structure
is compatible with the action of $\tau$; see [G], p. 168, lines 3 - 4
for the case $\w=\n F_q[T]$). We have: $Z_1^\tau$ is a
$Z_0^\tau$-module ( = $\w$-module) which is isomorphic to $\Omega(\w)$
(see the last lines of the proof of [G], Corollary 5.9.35 for the case
$\w=\n F_q[T]$), and $Z_1$ is isomorphic to $Z_0\otimes_{\w}\Omega(\w)$.

We shall consider $M$ as a $\w_C[\tau]$-module, like in 1.13.3. We denote
$M\{T\}:=M\underset{\w_C}\to{\otimes}Z_0$ ( = [G],
Definition 5.9.11.1 for the case $\w=\n F_q[T]$) and
$H^1(M):=M\{T\}^\tau$ like in [G], Definition 5.9.11.2.
Analogous to [G], Corollary 5.9.25 we get that for the present case
$$H_1(M):=\Hom_{\w_C[\tau]}(M,Z_1)=L_T(M)$$
($H_1(M)=H_1(E)$ of [G], 5.9). Particularly, for $M=\goth C$ we have
$$L_T(\goth C)=\Hom_{\w_C[\tau]}(\goth C,Z_1)$$

{\bf Lemma 5.2.2. } $H_1(M')=H^1(M)\underset{\w}\to{\otimes}L_T(\goth C)$.
\medskip
{\bf Proof.} By definition,
$\Hom_{\w_C} (M',Z_1)=\Hom_{\w_C} (\Hom_{\w_C} (M, \goth C), Z_1)$. Further,
$$\Hom_{\w_C} (\Hom_{\w_C} (M, \goth C), Z_1)=(M\underset{\w_C}
\to{\otimes}Z_0)\underset{Z_0}\to{\otimes} (\Hom_{\w_C}(\goth C,
Z_1))\eqno{(5.2.3)}$$ (an equality of  linear algebra). In order to show that we can
consider $\tau$-invariant subspaces, we need the following objects.
Let $I$ be an ideal of $\w$, $\Cal M_0=IZ_0$. It is clear that $\Cal M_0^\tau=I$.
Further, let $\Cal M_1$ be a locally free $Z_0$-module. We have a formula: $$(\Cal
M_0\underset{Z_0}\to{\otimes} \Cal M_1)^\tau=\Cal M_0^\tau \underset{\w}\to{\otimes} \Cal
M_1^\tau \eqno{(5.2.4)}$$
Really, $\Cal M_0\underset{Z_0}\to{\otimes} \Cal M_1=I\Cal M_1$, and $$(I\Cal
M_1)^\tau=I\Cal M_1^\tau\eqno{(5.2.5)}$$ where this formula is true by the following
reason. Obviously $(I\Cal M_1)^\tau\supset I\Cal M_1^\tau$. Let $J$ be an ideal of $\w$
such that $IJ$ is a principal ideal. We have $(IJ(J^{-1}\Cal M_1))^\tau=IJ(J^{-1}\Cal
M_1)^\tau$ and $(IJ(J^{-1}\Cal M_1))^\tau\supset I(J(J^{-1}\Cal M_1))^\tau\supset
IJ(J^{-1}\Cal M_1)^\tau$, hence all these objects are equal and we get 5.2.5 and hence
5.2.4.

The action of $\tau$ on both sides of 5.2.3 coincide.  Considering $\tau$-invariant
elements of both sides of 5.2.3 and taking into consideration 5.2.4 ($\Cal
M_0=\Hom_{\w_C}(\goth C,
Z_1)$ and $\Cal M_1=M\underset{\w_C}
\to{\otimes}Z_0$) we get the lemma. $\square$
\medskip
This lemma proves (A) of Step 1.
\medskip
{\bf Lemma 5.2.6.} Let $\Cal M_i$ ($i=0,1$) be two locally free $Z_0$-modules with $\tau$-action satisfying $\tau(cm)=\tau(c) \tau(m)$ ($c\in Z_0$, $m\in \Cal M_i$), and $\psi: \Cal M_0\otimes_{Z_0}\Cal M_1\to Z_1$ a perfect pairing of $Z_0$-modules with $\tau$-action. Let, further, both $\Cal M_i$ satisfy $\Cal M_i^\tau\otimes_{\w} Z_0=\Cal M_i$. Then the restriction of $\psi$ to $\Cal M_0^\tau\otimes_{\w}\Cal M_1^\tau\to \Omega$ is a perfect pairing as well.
\medskip
{\bf Proof.} Let $\alpha: \Cal M_0^\tau \to \Omega$ be an $\w$-linear map. We prolonge it to a map $\bar \alpha: \Cal M_0 \to Z_1$ by $Z_0$-$\tau$-linearity. By perfectness of $\psi$, there exists $m_1\in \Cal M_1$ such that $\bar \alpha(m_0)= \psi(m_0 \otimes m_1)$. It is easy to see that $m_1$ is $\tau$-invariant (we use the fact that $\tau: Z_0 \to Z_0$ is surjective). $\square$
\medskip
{\bf Lemma 5.2.7.} There is a natural perfect $\w$-linear
$\Omega$-valued pairing
between $H_1(M)$ and $H^1(M)$:
$H_1(M)\underset{\w}\to{\otimes}H^1(M)\to \Omega$.
\medskip
{\bf Proof.} For the case $\w=\n F_q[T]$ this is [G], Corollary
5.9.35. General case: we have a perfect $Z_0$-pairing
$$\Hom_{\w_C}(M,Z_1)\underset{Z_0}\to{\otimes} (M\underset{\w_C}
\to{\otimes} Z_0) \to Z_1$$
Now we take $\Cal M_0=\Hom_{\w_C}(M,Z_1)$, $\Cal M_1=M\underset{\w_C} \to{\otimes}Z_0$ and we apply Lemma 5.2.6. $\square$
\medskip
Step 1 of the theorem follows from these lemmas. 
\medskip
{\bf Remark 5.2.8.} The pairing can be defined also as the composition of
$$\matrix H_1(M)\underset{\w} \to{\otimes} H_1(M')=\Hom_{\w_C[\tau]}(M,Z_1)
\underset{\w}
\to{\otimes} \Hom_{\w_C[\tau]}(M',Z_1) \\
\to \Hom_{\w_C[\tau]}(M\underset{\w_C} \to{\otimes}M'
,Z_1\underset{Z_0} \to{\otimes}Z_1) \to
\Hom_{\w_C[\tau]}(\goth C,Z_1\underset{Z_0} \to{\otimes}Z_1) = L_T(\goth C) \underset{\w} \to{\otimes}\Omega\endmatrix \eqno{(5.2.9)}$$
where the second map comes from a canonical map $\delta: \goth C \to M\underset{\w_C}
\to{\otimes}M'$ of Remark 1.11.1 (more exactly, of its analog for arbitrary $\w$).
\medskip
{\bf Remark 5.2.10.} Recall that the explicit formula for functoriality is
the following. Let $\alpha: M_1 \to M_2$ be a map of t-motives,
$\alpha': M'_2 \to M'_1$
the dual map and $L_T(\alpha): L_T(M_2) \to L_T(M_1)$, $L_T(\alpha'): L_T(M'_1)
\to L_T(M'_2)$ the
corresponding maps on lattices. For any $l_1'\in L_T(M_1')$, $l_2\in L_T(M_2)$ we
have:
$$<L_T(\alpha)(l_2), l_1'>_{M_1}=<l_2, L_T(\alpha')(l_1')>_{M_2}\eqno{(5.2.11)}$$

{\bf 5.3. Proof: Step 1, Version 2.} Case $\w=\n F_q[T]$. We identify $Z_1$ of [G], p.168, lines 3 -- 4 with $\p\{T\}$ (see [G], Definition 5.9.10) and $\w$ with $\Omega$. Like above, we have an isomorphism of $\w$-modules (recall that $\w$ is the center of $\p[T,\tau]$):
$$L_T(M)=\Hom_{\p[T,\tau]}(M,Z_1)\eqno{(5.3.1)}$$ ([G], first terms of 5.9.25, 5.9.19). Let $\vf: M \to Z_1$, $\vf':
M' \to Z_1$ be elements of $L_T(M)$, $L_T(M')$ respectively, and let
$f_*$, $f'_*$, $Q$, $Q'$ be from
1.9.3, 1.10. We denote $$\vf(f_*)=v_*\eqno{(5.3.2)}$$ where $v_*\in
(Z_1)^r$ is a vector column (it is a column of the scattering matrix ([A], p. 486) of
$M$, see 5.4.1 below). The same notation for
the dual:
$\vf'(f'_*)=v'_*$. Condition
that $\vf$, $\vf'$ are
$\tau$-homomorphisms is equivalent to $$Qv_*=v_*^{(1)}, \ \
Q'v'_*={v'}_*^{(1)}\eqno{(5.3.3)}$$ (analog of the formula for scattering matrices [A],
(3.2.2)). Let us consider $\Xi=\sum_{i=0}^\infty
a_iT^i\in\p\{T\}\subset\p[[T]]$ of [G], p. 172, line 1; recall that it is the only
element (up to
multiplication by $\n F_q^*$) satisfying
$$\Xi=(T-\theta)\Xi^{(1)}, \ \ \ \lim_{i\to\infty}a_i=0, \ \ \ |a_0|>|a_i| \ \ \forall
i>0\eqno{(5.3.4)}$$ (see [G], p. 171,
(*); there is a formula $\Xi=a_0\prod_{i\ge0}(1-T/\theta^{q^i})$ where $a_0$ satisfies
$a_0^{q-1}=-1/\theta$). Finally, we define
$$<\vf, \vf'>=\Xi v_{*}^tv'_{*}\eqno{(5.3.5)}$$
Obviously $<\vf, \vf'>$ does not depend on a choice of a basis $f_*$.
\medskip
{\bf Lemma 5.3.6.} $<\vf, \vf'>\in \w$.
\medskip
{\bf Proof.} Firstly, this element belongs to $\n F_q[[T]]$, because $$\Xi
v_{*}^tv'_{*} -(\Xi v_{*}^tv'_{*})^{(1)}=\Xi (v_{*}^tv'_{*}-(T-\theta)^{-1}
v_{*}^{(1)t}{v'}_{*}^{(1)})= \Xi v_{*}^t(E_r-(T-\theta)^{-1}
Q^tQ')v'_{*}$$ because of (5.3.3). But we have (see (1.10.1) --- the
definition of $Q'$)
$$E_r-(T-\theta)^{-1} Q^tQ'=0$$
Secondly, let $<\vf, \vf'>=\sum_{i=0}^\infty c_iT^i$. Since
coefficients of all factors of (5.3.5): $\Xi$, $v_*$ and $v'_*$ ---
tend to 0, we get
that $c_i$ also tend to 0. But $c_i\in\n F_q$, i.e. they are almost
all 0. $\square$
\medskip
{\bf Lemma 5.3.7.} The above pairing is perfect.
\medskip
{\bf Proof.} We have an isomorphism (here $M\{T\}=M\otimes_{\p[T]}\p\{T\}$ with
the natural action of $\tau$ (see [G], Definition 5.9.11))
$$\alpha: \Hom_{\p[T,\tau]}(M,Z_1)\to\Hom_{\w}(M\{T\}^\tau,\w)\eqno{(5.3.8)}$$ defined as the composition of the maps
$$\Hom_{\p[T,\tau]}(M,Z_1)=\Hom_{\p[T]}(M,Z_1)^\tau\overset{\beta'}\to{\to}\Hom_{\p\{T\}}(M\{T\},\p\{T\})^\tau$$
$$\overset{\gamma}\to{\to} \Hom_{\w}(M\{T\}^\tau,\w)$$ where
$\beta: \Hom_{\p[T]}(M,Z_1)\to \Hom_{\p\{T\}}(M\{T\},\p\{T\} )$ is the natural map and
$\beta'$
is the restriction of $\beta$ to $\tau$-invariant elements. Using the Anderson's
criterion of uniformizability
of $M$ (see, for example, [G], 5.9.14.3 and 5.9.13) we get immediately that both
$\gamma$, $\beta$, and hence
$\beta'$, and hence $\alpha$ are isomorphisms. Further, let us consider a homomorphism
$$i: \Hom_{\p[T,\tau]}(M',Z_1)\to M\{T\}^\tau\eqno{(5.3.9)}$$ defined as follows. Let
$\vf'$, $f'_*$, $v'_*$ be as above. We set
$$i(\vf')=\Xi{v'}^t_*f_*\in M\underset{\p[T]}\to{\otimes}\p[[T]]$$ Since $\Xi\in \p
\{T\}$, we get that $\Xi{v'}^t_*f_*\in M\{T\}$. A simple calculation
(like in the Lemma 5.3.6, but simpler) shows that $i(\vf')$ is $\tau$-invariant, hence
$i$ really defines a map from $\Hom_{\p[T,\tau]}(M',Z_1)$ to $M\{T\}^\tau$. Obviously it
is an inclusion. Let us prove that $i$ is surjective. Really, let $c_*\in (Z_1)^r$ be a
column
vector such that $c_*^tf_*\in M\{T\}^\tau$. An analog of the above calculation shows that
if we define $\vf'$ by the formula $\vf'(f'_*)=\Xi^{-1}c_*$ then $\vf'\in
\Hom_{\p[T,\tau]}(M',Z_1)$, and $i(\vf')=c_*^tf_*\in M\{T\}^\tau$. Finally, the
combination of isomorphisms (5.3.8) and (5.3.9) corresponds to the
pairing (5.3.5). $\square$
\medskip
{\bf 5.4. Step 2 -- End of the proof of Theorem 5.} It is easy to see that the converse of the Corollary 5.1.1 (taking into consideration
Proposition 3.2) is also true, i.e. in order to prove Theorem 5 it is sufficient
to prove that a Siegel matrix of $M'$ is $-Z^t$ where $Z$ is a Siegel matrix of $M$. Let
us consider a basis $l_1,...,l_r$ of $L_T(M)$ and for each $l_i$ we consider the
corresponding (under identification 5.3.1) $\vf_i\in \Hom_{\p[T,\tau]}(M,Z_1)$. Let
$\Psi$ be the scattering matrix of $M$ ([A], p. 486) with respect to the bases $l_1,...,l_r$,
$f_1,...,f_r$, and we denote $\vf_i(f_*)$ by $v_{i*}$ (notations of 5.3.2).
\medskip
{\bf Lemma 5.4.1.} $v_{i*}$ is the $i$-th column of $\Psi$ ($Z_1$ is identified with $\p\{T\}$, see the proof).
\medskip
{\bf Proof.} Follows from the definitions. Recall that $\hbox{\bf K}=\n F_q((1/T))$. The isomorphism 5.3.1 is the composition of
2 isomorphisms $i_1: L_T(M) \to \Hom_{\w}^c(\hbox{\bf K}/\w, E)$ ([G], 5.9.19) and $i_2:
\Hom_{\w}^c(\hbox{\bf K}/\w, E) \to \Hom_{\p[T,\tau]}(M,\Hom^c(\hbox{\bf K}/\w, \p)$
([G], 5.9.24; recall that $Z_1=\Hom^c(\hbox{\bf K}/\w, \p)$). For $l_i\in L_T(M)$ we have
$(i_1(l_i))(T^{-k})=\exp(\theta^{-k}l_i)$ ([G], line above the lemma 5.9.18) and
$$((i_2\circ i_1(l_i))(f_j))(T^{-k})=<f_j,\exp(\theta^{-k}l_i)>$$ ([G], two lines above the
lemma 5.9.24). Using the identification of $Z_1$ and $\p\{T\}$ ([G], p. 168, lines 3 - 4)
and the definition of $\Psi$ ([A], p. 486, first formula of 3.2) we get immediately the
lemma. $\square$
\medskip
Let $l'_1,...,l'_r$ be a basis of $L_T(M')$ which is dual to a basis $l_1,...,l_r$ of
$L_T(M)$ with respect to the pairing 5.3.5.
\medskip
{\bf Lemma 5.4.2.} The scattering matrix of $M'$ with respect to the bases
$l'_1,...,l'_r$, $f'_1,...,f'_r$ (denoted by $\Psi'$) is $\Xi^{-1}\Psi^{t-1}$.
\medskip
{\bf Proof.} Follows immediately from 5.4.1 applied to both $M$, $M'$, and formula 5.3.5.
$\square$
\medskip
{\bf Remark 5.4.3.} An alternative proof for the case of pure $M$ (for $some$ basis of $L_T(M')$)
is the following. We denote $\Xi^{-1}\Psi^{t-1}$ by $\Psi_1$. It satisfies
$\Psi_1^{(1)}=(T-\theta) Q^{t-1}\Psi_1$ and other conditions of [A], 3.1. According [A],
Theorem 5, p. 488, there exists a pure uniformizable t-motive $M_1$ with
$\sigma$-structure such that its scattering matrix is $\Psi_1$. Since $\Psi_1$ satisfies
$$\Psi_1^{(1)}=Q'\Psi_1$$
we get that $Q(M_1)=Q'$, i.e. $M_1=M'$. $\square$
\medskip
Let us recall the statement of the crucial proposition 3.3.2 of [A]. Here we consider the
case of those $M$ whose $N$ is not necessarily 0. Let $\Psi$ be a scattering matrix of
$M$. We consider the $(T-\theta)$-Laurent series for $\Psi$ (here $k(M)<0$ is a number,
and $D_{-i}\in M_{r}(\p)$): $$\Psi=\sum_{i=k(M)}^\infty D_{-i}(T-\theta)^i$$ We consider its
negative part $$\Psi^-:=\sum_{i=k(M)}^{-1} D_{-i}(T-\theta)^i$$ as an element of $M_{r}(\p)((T-\theta))/M_{r}(\p)[[T-\theta]]$.
\medskip
We consider the space $(T-\theta)^{k(M)}\p[[T-\theta]]/\p[[T-\theta]]$ as a $\p$-vector space endowed by the action of $\w$, and we denote by $\goth V$ its $r$-th direct sum written as vector columns of length $r$. Obviously $$k(M)=-1 \iff \hbox{ the action of $T$ on $\goth V$ coincides with multiplication by $\theta$} \eqno{(5.4.3a)}$$
We denote the $i$-th column of $\Psi^-$ by $\Psi^-_{i*}$, it belongs to $\goth V$. Further, we denote by $\Prin(M)$ (resp. by $\Prin_0(M)$) the $\p[T]$-linear envelope (resp. the $\w$-linear envelope) of all $\Psi^-_{i*}$ in $\goth V$. Finally, we obviously extend the definition of $\Lie_T(M)$, $L_T(M)$ to the case $N\ne 0$; formula 5.1.4 becomes
$$\goth a(Tx)=(\theta+N)(\goth a(x))\eqno{(5.4.3b)}$$
\medskip
{\bf Proposition 3.3.2, [A]} (see also Remark 5.5 below). There exists a $\p[T]$-linear isomorphism $\psi_E: \Lie_T(M)
\to \Prin(M)$ such that its restriction to $L_T(M)\subset \Lie_T(M)$ defines an isomorphism
$L_T(M) \to \Prin_0(M)$ (denoted by $\psi_E$ as well). $\square$
\medskip
{\bf Corollary 5.4.4.} $N=0 \iff k(M)=-1$ (because $N=0 \iff $ the action of $T$ on both $\Lie_T(M)$, $\goth V$ coincides with multiplication by $\theta$, by 5.4.3a). $\square$
\medskip
We return to the case $N=0$.
\medskip
Let us consider the $(T-\theta)$-Laurent series for $\Psi'$ and $\Xi^{-1}$:
$$\Psi'=\sum_{i=k(M')}^\infty D'_{-i}(T-\theta)^i,  \ \ \ \Xi^{-1}=\sum_{i=k(\xi)}^\infty
a_i(T-\theta)^i$$
Since for both $M$, $M'$ we have $N=N'=0$, we get $k(M)=k(M')=-1$. An
elementary calculation shows that $k(\xi)$ is also $-1$. Hence, equality
$\Psi'\Psi^t=\Xi^{-1}$ (Lemma 5.4.2) implies that $D'_{1}D^t_{1}=0$.

Further, there exist $n$ columns of $D_{1}$ which are
$\p$-linerly independent (they are $\psi_E$-images of elements of $L_T(M)$ which form a
$\p$-basis of $\Lie_T(M)$) and all other columns of $D_{1}$ are their linear combinations.
Interchanging columns of $D_{1}$ if necessary we can assume that these columns are the
last $n$
columns. We denote by $D_{12}$ (resp. $D_{11}$ ) the $r\times n$ (resp. $r\times
(r-n)$ ) matrix formed by the last $n$ (resp. the first $r-n$) columns of $D_{1}$. There
exists a matrix $S$ such that $D_{11}=D_{12}S^t$. Again according Proposition 3.3.2,
[A], we have:
$$S \hbox{ is a Siegel matrix of $L(M)$ }\eqno{(5.4.5)}$$
(see also Remark 5.5 below).

Analogous objects are defined for $D'_{1}$. We denote by $D'_{12}$ (resp.
$D'_{11}$) the $r\times n$- (resp. $r\times (r-n)$)-matrix formed by the last $n$
(resp. the first $r-n$) columns of $D'_{1}$. Since
$D'_{1}D^t_{1}=D'_{11}D_{11}^t+D'_{12}D_{12}^t$ we get that
$D'_{12}D_{12}^t+D'_{11}SD_{12}^t=0$. Since $D_{12}^t$ is a $n\times
r$-matrix of rank $n$, it is not a zero-divisor from the right, so $$D'_{12}=
-D'_{11}S\eqno{(5.4.6)}$$ Since the rank of $D'_{1}$ is $r-n$ and $D'_{11}$ is a
$r\times (r-n)$ matrix, (5.4.6) implies that columns of $D'_{11}$ are linearly independent,
and by (5.4.6) and Proposition 3.3.2, [A] we get that $-S$ is a Siegel matrix of
$M'$. $\square$
\medskip
{\bf Remark 5.5.} Since the notations of [A] differ from the ones of the present
paper, for the reader's convenience we give here a sketch of the proof for the case $N=0$
of two  crucial facts: Corollary 5.4.4 and 5.4.5 ([A], Theorem 3.3.2).

Let $\alpha: \Lie(M)\to E(M)$ be a linear isomorphism which is the first term of the
series for $\exp: \Lie(M)\to E(M)$, and let $l\in \Lie(M)$, $f \in M$ be arbitrary. We
consider the $(T-\theta)$-Laurent series $\sum_{i=k}^\infty b_i(T-\theta)^i$ of
$\sum_{j=0}^\infty <\exp(\frac{1}{\theta^{j+1}}l),f>T^j$.
\medskip
{\bf Lemma 5.6.} If $N=0$ then $k=-1$, and $b_{-1}=-<\alpha(l),f>$ (this is [A],
3.3.4).
\medskip
{\bf Sketch of the proof.} For $z\in \Lie(M)$ we denote $\exp(z)- \alpha(z)$ by $\ve(z)$,
hence
$\sum_{j=0}^\infty <\exp(\frac{1}{\theta^{j+1}}l),f>T^j=\underline{A} +\underline{E}$, where
$$\underline{A}=\sum_{j=0}^\infty <\alpha(\frac{1}{\theta^{j+1}}l),f>T^j; \ \ \ \ \underline{E}=\sum_{j=0}^\infty <\ve(\frac{1}{\theta^{j+1}}l),f>T^j$$
We consider their $(T-\theta)$-Laurent series:
$$\underline{A}=\sum_{i=k(\underline{A})}^\infty \underline{a}_i(T-\theta)^i; \ \ \ \ \underline{E}=\sum_{i=k(\underline{E})}^\infty \underline{e}_i(T-\theta)^i$$
Since we have $\exp(z)=\sum_{i=0}^{\infty}C_iz^{(i)}$ where $C_0=E_n$ we get that
$\ve(z)=\sum_{i=1}^{\infty}C_iz^{(i)}$. This means that for large $j$ the element
$\ve(\frac{1}{\theta^{j+1}}l)$ is small, and hence $k(\underline{E})=0$, because finitely many
terms having small $j$ do not contribute to the pole of the $(T-\theta)$-Laurent series
of $\underline{E}$ (the reader can prove easily the exact estimations himself, or to look [A],
p. 491). Since $\alpha$ is $\p$-linear, equality $\sum_{j=0}^\infty
\frac{1}{\theta^{j+1}}T^j= -(T-\theta)^{-1}$ implies that $k(\underline{A})=-1$ and $\underline{a}_{\ -1}=-<\alpha(l),f>$ (and other $\underline{a}_i=0$), hence the lemma. $\square$
\medskip
This lemma obviously implies Corollary 5.4.4. Further, elements $f_1,...,f_r$ generate
the $\p$-space $M/\tau M$, because multiplication by $T$ on $M/\tau M$ coincides with
multiplication by $\theta$, hence the fact that $f_1,...,f_r$ \ \ $\p[T]$-generate $M/\tau
M$ implies that they $\p$-generate $M/\tau M$.

Let $l_1,...,l_n$ form a $\p$-basis of $\Lie(M)$ (here we identify $\Lie_T(M)$ and $\Lie(M)$ via $\goth a$). Since the pairing $<*,*>$ between
$E(M)$ and $M/\tau M$ is non-degenerate and $\alpha$ is an isomorphism, we get that
columns $<\alpha(l_1),f_*>,..., <\alpha(l_n),f_*>$ are linearly independent. Again since
$\alpha$ is an isomorphism and the pairing with $f_*$ is linear, we get that
$$(<\alpha(l_{n+1}),f_*> \ ... \ <\alpha(l_r),f_*>)=(<\alpha(l_{1}),f_*> \ ... \
<\alpha(l_n),f_*>)Z^t $$ Applying the lemma 5.6 to this formula we get immediately 5.4.5.
\medskip
{\bf 6. Tensor products.}
\medskip
There exists an analog of the Theorem 5 for the case of tensor
products of
t-motives. It describes the lattice $L(M_1\otimes M_2)$ in terms of $L(M_1)$,
$L(M_2)$. This is a theorem of Anderson; it is formulated in [P], end of page 3, but its
proof was not published. We recall its statement for the case of arbitrary $N\ne 0$, and
we give its proof for the case $N=0$ (case of arbitrary $N$ can be obtained easily using
the same ideas).

Let $M$ be an uniformizable t-motive whose $N$ is not necessarily 0. Since $N$ is nilpotent, formula 5.4.3b shows that $\Lie_T(M)$ is a $\p[[T-\theta]]$-module. There exists an epimorphism of
$\p[[T-\theta]]$-modules
$$L_T(M)\underset{\w}\to{\otimes}\p[[T-\theta]]\to \Lie_T(M)$$
whose kernel $\goth q=\goth q(M)$ carries information on the pair $(L(M), \Lie(M))$.
\medskip
{\bf Theorem 6} (Anderson). Let $M$, $\bar M$ be any two uniformizable abelian
t-motives. Then $$\goth q(M\otimes \bar M)=\goth q(M)\underset{\p[[T-\theta]]}
\to{\otimes}\goth q(\bar M)\eqno{(6.1)}$$
\medskip
{\bf Remark 6A.} $M\otimes \bar M$ is a uniformizable t-motive ([G], Corollary 5.9.38).
\medskip
{\bf Proof of Theorem 6 (case $N=0$).} We define notations for $M$, and all notations for $\bar M$
will carry bar. Let $e_i$ and $Z$ be from the beginning of Section 3. We denote $\goth a^{-1}(e_i)\in \Lie_T(M)$ by $e_i$ (there is no possibility of confusion). So, $\{e_i\}$ is a
$\p[[T-\theta]]$-basis of $L_T(M)\underset{\w}\to{\otimes}\p[[T-\theta]]$. Elements
$b_i:=(T-\theta)e_i$, $i=1,...,n$ and $b_{n+i}:=e_{n+i}-\sum_{j=1}^n z_{ij}e_j$,
$i=1,...,r-n$ form a $\p[[T-\theta]]$-basis of $\goth q$. We need a
\medskip
{\bf Lemma 6.2.} $\Psi(M\otimes \bar M)=\Psi(M)\otimes\Psi(\bar M)$ where $\Psi(M)$
(resp. $\Psi(\bar M)$; $\Psi(M\otimes \bar M)$) is taken with respect to bases $e_*$ of
$L_T(M)$, $f_*$ of $M_{\p[T]}$ (resp. $\bar e_*$ of $L_T(\bar M)$, $\bar f_*$ of $\bar
M_{\p[T]}$; $e_*\otimes \bar e_*$ of $L_T(M\otimes \bar M)$, $f_*\otimes \bar f_*$ of
$(M\otimes \bar M)_{\p[T]}$) (see the proof for the notations).
\medskip
{\bf Proof.} We consider a map $$\alpha:\Hom_{\p[T]}(M,Z_1)^\tau
\underset{\w}\to{\otimes} \Hom_{\p[T]}(\bar M,Z_1)^\tau\to \Hom_{\p[T]}(M\otimes \bar
M,Z_1)^\tau$$ defined as follows: for $\vf\in \Hom_{\p[T]}(M,Z_1)^\tau$, $\bar \vf\in
\Hom_{\p[T]}(\bar M,Z_1)^\tau$ we let $[\alpha(\vf \otimes \bar \vf)](f\otimes \bar
f)=\vf(f)\cdot\bar \vf(\bar f)$ (it is obvious that $\alpha(\vf \otimes \bar \vf)$ is
$\tau$-stable). Since $e_1, ..., e_r$ (resp. $\bar e_1, ..., \bar e_{\bar r}$) is a basis
of $\Hom_{\p[T]}(M,Z_1)^\tau$ (resp. $\Hom_{\p[T]}(\bar M,Z_1)^\tau$; we identify $L_T(M)$,
resp. $L_T(\bar M)$ with $\Hom_{\p[T]}(M,Z_1)^\tau$ (resp. $\Hom_{\p[T]}(\bar M,Z_1)^\tau$)
we get (using Lemma 5.4.1) that $\Psi(M)$, $\Psi(\bar M)$ are non-degenerate. Since their
product is also non-degenerate, we get $\alpha(e_i\otimes \bar e_{\bar i})$ are linearly
independent and hence a basis of $\Hom_{\p[T]}(M\otimes \bar M,Z_1)^\tau$. Applying once
again Lemma 5.4.1 we get the lemma. $\square$
\medskip
If $A$, $B$ are two matrices then columns of $A\otimes B$ are indexed by pairs $(k,l)$
where $k$ (resp. $l$) is the number of a column of $A$ (resp. $B$). We denote by $A_k$,
$B_l$, $A\otimes B_{(k,l)}$ the respective columns. Obviosly we have: $A\otimes
B_{(k,l)}=A_k\otimes B_l$ (tensor product of column matrices).
\medskip
Let us prove that for $i=1,...,r-n$, $\bar i=1,...,\bar r-\bar n$ the element
$b_{n+i}\otimes \bar b_{\bar n+\bar i}\in \goth q(M\otimes \bar M)$. According [A],
Proposition 3.3.2, it is sufficient to prove that the corresponding linear combination
(see 6.3 below) of the columns of the matrix $\Psi^-_{M\otimes \bar M}$ is 0. Since
$$b_{n+i}\otimes \bar b_{\bar n+\bar i}=\sum_{j,\bar j}z_{ij}\bar z_{\bar i\bar
j}e_j\otimes \bar e_{\bar j}- \sum_{j}z_{ij}e_j\otimes \bar e_{\bar n+\bar i}
-\sum_{\bar j}\bar z_{\bar i\bar j}e_{n+i}\otimes \bar e_{\bar j}+ e_{n+i}\otimes \bar
e_{\bar n+\bar i}$$ we get the explicit form of this linear combination: it is sufficient
to prove that for all $i$, $\bar i$ we have
$$\sum_{j,\bar j}z_{ij}\bar z_{\bar i\bar j} (\Psi^-_{M\otimes \bar M})_{(j,\bar j)} -
\sum_{j}z_{ij}(\Psi^-_{M\otimes \bar M})_{(j, \bar n+\bar i)} $$ $$ -\sum_{\bar j}\bar
z_{\bar i\bar j}(\Psi^-_{M\otimes \bar M})_{(n+i,\bar j)} +(\Psi^-_{M\otimes \bar
M})_{(n+i,\bar n+\bar i)} =0\eqno{(6.3)}$$
Further, 6.2 implies that $$(\Psi^-_{M\otimes \bar M})_{(k,\bar
k)}=\frac{A_{-1,k}\otimes \bar A_{-1,\bar k}}{(T- \theta)^2}+\frac{A_{-1,k}\otimes \bar
A_{0,\bar k}+A_{0,k}\otimes \bar A_{-1,\bar k}}{T- \theta}$$ hence 6.3 becomes
$$\sum_{j,\bar j}z_{ij}\bar z_{\bar i\bar j} (\frac{A_{-1,j}\otimes \bar A_{-1,\bar
j}}{(T- \theta)^2}+\frac{A_{-1,j}\otimes \bar A_{0,\bar j}+A_{0,j}\otimes \bar A_{-1,\bar
j}}{T- \theta}) $$ $$- \sum_{j}z_{ij}(\frac{A_{-1,j}\otimes \bar A_{-1,\bar n+\bar i}}{(T-
\theta)^2}+\frac{A_{-1,j}\otimes \bar A_{0,\bar n+\bar i}+A_{0,j}\otimes \bar A_{-1,\bar
n+\bar i}}{T- \theta}) $$ $$-\sum_{\bar j}\bar z_{\bar i\bar j}(\frac{A_{-1,n+i}\otimes
\bar A_{-1,\bar j}}{(T- \theta)^2}+\frac{A_{-1,n+i}\otimes \bar A_{0,\bar
j}+A_{0,n+i}\otimes \bar A_{-1,\bar j}}{T- \theta}) $$ $$+\frac{A_{-1,n+i}\otimes \bar
A_{-1,\bar n+\bar i}}{(T- \theta)^2}+\frac{A_{-1,n+i}\otimes \bar A_{0,\bar n+\bar
i}+A_{0,n+i}\otimes \bar A_{-1,\bar n+\bar i}}{T- \theta}=0\eqno{(6.4)}$$
It is easy to see that 6.4 follows immediately from the equalities
$$A_{-1,n+i}=\sum_j z_{ij}A_{-1,j}\eqno{(6.5)}$$ $$ \bar A_{-1,\bar n+\bar
i}=\sum_{\bar j} \bar z_{\bar i\bar j}\bar A_{-1,\bar j}$$ For example, the left hand
side of (6.4) has 2 terms containing $\bar A_{0,\bar j}$ (in the middle of the first
and the third lines of (6.4)). Multiplying (6.5) by $\bar z_{\bar i\bar j}\bar
A_{0,\bar j}$ we get that the sum of these 2 terms of (6.4) is 0. For other pairs of
terms of (6.4) the situation is the same.
\medskip
The proof that for $i=1,...,r-n$, $\bar i=1,...,\bar n$ the element $b_{n+i}\otimes \bar
b_{\bar i}\in \goth q(M\otimes \bar M)$ is analogous but simpler. We have
$$b_{n+i}\otimes \bar b_{\bar i}=(T-\theta) (-\sum_{j} z_{ij} e_j\otimes \bar e_{\bar i}
+ e_{n+i}\otimes \bar e_{\bar i})$$ The analog of (6.3)) is $$(T-\theta) (-
\sum_{j}z_{ij}(\Psi^-_{M\otimes \bar M})_{(j, \bar i)}  +(\Psi^-_{M\otimes \bar
M})_{(n+i,\bar i)}) =0$$ and the analog of (6.4)) is $$-
\sum_{j}z_{ij}\frac{A_{-1,j}\otimes \bar A_{-1,\bar i}}{T- \theta}
+\frac{A_{-1,n+i}\otimes \bar A_{-1,\bar i}}{T- \theta}=0$$ This equality follows
immediately from (6.5).
\medskip
Finally, elements $b_{i}\otimes \bar b_{\bar i}$ ($i=1,...,n$, $\bar i=1,...,\bar n$)
obviously belong to $ \goth q(M\otimes \bar M)$.
\medskip
So, we proved that $\goth q(M)\underset{\p[[T-\theta]]}\to{\otimes}\goth q(\bar M)
\subset \goth q(M\otimes \bar M)$. Since the $\p$-codimension of both subspaces in
$L_T(M)\underset{\w}\to{\otimes} L_T(\bar M)\underset{\w}\to{\otimes}
\p[[T-\theta]]$ is $n\bar n$, they are equal. $\square$
\medskip
{\bf 7. Self-dual t-motives.}
\medskip
{\bf Case $\w=\n F_q[T]$.} A uniformizable t-motive $M$ is called self-dual if there exists an isogeny $\alpha: M \to M'$. It defines an $\w$-valued, $\w$-bilinear form $<*,*>_\alpha$ on $L_T(M')$ as follows:
$$<\vf_1, \vf_2>_\alpha=<L_T(\alpha)(\vf_1), \vf_2>_M$$
5.2.11 implies that if $\alpha'=-\alpha$ (resp. $\alpha'=\alpha$) then
$<*,*>_\alpha$ is skew symmetric (resp. symmetric). $M$ is called positively (resp. negatively) self-dual if $\alpha$ satisfies $\alpha'=\alpha$ (resp. $\alpha'=-\alpha$). Hence, we have an
\medskip
{\bf Analogy 7a.} The number field case analog of a pair: $\{$negatively self-dual t-motive of rank $2n$, dimension $n$; negative $\alpha: M \to M'\}$ is a (generic) abelian variety of dimension $n$ with a fixed polarization form.
\medskip
For example, like in the number field case, we can define the Rosati involution $I_\alpha$ on  $\End_0(M):=\End(M)\otimes \n F_q(T)$ by the same formula $I_\alpha(\vf)=\alpha^{-1}\circ \vf'\circ \alpha$.
\medskip
Further, we have a
\medskip
{\bf Conjecture 7b.} The dimension of the moduly variety of negatively self-dual t-motives (if it exists) is $n(n+1)/2$.
\medskip
{\bf Examples.} Let $e_*$ be from 1.9, and let $M=M(A)$ given by the equation (here $A\in M_n(\p)$ is $\goth A_1$ of 1.9.1)
$$Te_*= \theta e_*+ A\tau e_* + \tau^2 e_*\eqno{(7.1)}$$
be a t-motive of dimension $n$ and rank $2n$. Elements $f_i=e_i$,
$f_{n+i}=\tau e_i$ $(i=1,...,n)$ form a $\p[T]$-basis of $M$. We have
(see, for example, Section 11): $M'$ is given by the equation $$Te'_*=
\theta e'_*- A^t\tau e'_* + \tau^2 e'_*$$ and if we define $$f'_i=\tau
e'_i, \ \ f'_{n+i}= e'_i\eqno{(7.2)}$$ then bases $f_*$, $f'_*$ are dual in the
meaning of Lemma 1.10.

Let $\alpha: M \to M'$ be given by the formula $\alpha(e_*)=De'_*$
where $D\in M_n(\p)$ (we impose this essential restriction only in order to simplify exposition. In the general case $D\in M_n(\p[\tau])$, $D_f\in M_{2n}(\p[T])$, $D_f$ from 7.4). Condition that $\alpha$ is a $\p[T,\tau]$-map is
equivalent to
$$D^{(2)}=D, \ \ AD^{(1)}=-DA^t\eqno{(7.3)}$$
Further, we have $$\alpha(f_*)=D_ff'_*\eqno{(7.4)}$$ where
$D_f=\left(\matrix 0 & D\\ D^{(1)}&0 \endmatrix \right)$, hence
$$\alpha'=\pm \alpha\iff
D_f^t=\pm D_f\iff D^{(1)}=\pm D^t\eqno{(7.5)}$$
Let us fix $\varepsilon_0\in \n F_{q^2}$ satisfying $\varepsilon_0^{q-1}=-1$. Then $D=\varepsilon_0E_n$ satisfies 7.5 with the sign minus, and the set of $A$ satisfying 7.3 with this $D$ is the set of symmetric matrices. This justifies 7b, because the set of $A_1\in  M_n(\p)$ such that $M(A)=M(A_1)$ is conjecturally discrete.

For $D=E_n$ the sign in 7.5 is plus and hence a skew symmetric $A$ defines a positively self-dual $M(A)$.
\medskip
{\bf Remark 7.6.} The below statements are conjectures based on arguments similar to the ones which justify the below Conjecture 11.8.3. Since they are of secondary importance, we do not give any details of justification here.
\medskip
{\bf 7.6.1. Conjecture.} If $n\ge 3$ then for a
generic skew symmetric $A$ we have: $\End(M(A))=\w$.
\medskip
{\bf 7.6.2. Corollary.} Conjecture 7.6.1 implies that the "minimal" $\alpha: M \to M'$ is defined uniquely up to an element of
$\n F_q^*$, and hence the symmetric pairing $<*,*>_\alpha$ is also defined uniquely up to an element of $\n F_q^*$.
\medskip
{\bf 7.6.3. Conjecture.} If $n=2$, $\alpha'=\alpha$ then $\End (M)$ is strictly larger than $\w$.
\medskip
Other examples of a self-dual t-motive are $M\oplus M'$ where $M$ is any t-motive, but they do not give interesting examples of pairings.
\medskip
{\bf 7.6.4. Conjecture.} There exist other (distinct from the ones defined by 7.1) self-dual t-motives $M$ having $\End (M)=\w$ (we can use a version of standard t-motives of Section 11).
\medskip
{\bf Example 7.7.} Case $A=0$, $D=E_n$.
\medskip
In this case we can find explicitly the matrix of the symmetric form
$<*,*>_\alpha$ in some basis of $L_T(M')$. Let $\goth C_2$ be the
Carlitz module over the field $\n F_{q^2}$ considered as a rank 2
Drinfeld module over $\n F_{q}$ given by the equation $$Te=\theta e
+ \tau^2e$$ We have $M=\goth C_2^{\oplus n}$. Let $\goth T_T(\goth
C_2)$ be the convergent $T$-Tate module of $\goth C_2$, i.e. the set
of elements $\{z_i\}\in E(\goth C_2)=\p$ $(i\ge -1, \ \ z_{-1}=0)$ such that
$$\hbox{ $T z_i=z_{i-1}$ for $i \ge 0$ (i.e. $z_i^{q^2}+\theta
z_i=z_{i-1}$) and $z_i\to 0$}$$ It is a free 1-dimensional module over
$\n F_{q^2}[T]$. We choose and fix its generator; its $\{z_i\}$ satisfy (like in 5.3.4)
$|z_0|>|z_i| \ \ \forall i>0$. We denote $\sum_{k=0}^\infty z_kT^k$ by $\goth Z$.

Let $c$ be a fixed element of $\n F_{q^2}-\n F_{q}$. Formulas (5.3.3)
show that the following elements $\vf_i$, $\vf'_i$ ($i=1, ... , 2n$)
form bases of $L(M)$, $L(M')$ respectively ($j=1, ..., n$; clearly that thanks to 7.2
we have $\vf'_i(f'_j)=\vf_i(f_{n+j})$, $n+j$ mod $2n$):
$$i\le n: \ \ \vf_i(f_j)=\goth Z\delta_i^j, \ \
\vf_i(f_{n+j})=\goth Z^{(1)}\delta_i^j $$

$$i>n: \ \ \vf_i(f_j)=c\goth Z\delta_{i-n}^j, \ \
\vf_i(f_{n+j})=c^q\goth Z^{(1)}\delta_{i-n}^j $$

$$i\le n: \ \ \vf'_i(f'_j)=\goth Z^{(1)}\delta_i^j, \ \
\vf'_i(f'_{n+j})=\goth Z\delta_i^j $$

$$i>n: \ \ \vf'_i(f'_j)=c^q\goth Z^{(1)}\delta_{i-n}^j, \
\ \vf'_i(f'_{n+j})=c\goth Z\delta_{i-n}^j $$
(by the way, it is clear that the same relation between elements of $\goth T_T(M)$ and
$\Hom_{\p[T,\tau]}(M,Z_1)$ holds for all $M$). Formula 7.4 shows
that $\alpha'(\vf'_i)=\vf_{i+n}$, where $i+n \mod 2n$. Let us denote
$\Xi\cdot \goth Z\cdot \goth Z^{(1)}\in\n F_q^*$
by $\gamma$. The above definitions and formulas show that the matrix
of $<*,*>_\alpha$ in the basis $\vf_{1}, \vf_{n+1}, ... ,\vf_{n},
\vf_{2n}$ consists of $n$ \ \ $(2\times 2)$-blocks (trace and norm of $\n
F_{q^2}/\n F_q$)
$$\gamma\left(\matrix \tr(1) & \tr(c)\\ \tr(c)&\tr(N(c))  \endmatrix
\right)=\gamma\left(\matrix 2 & c+c^q \\ c+c^q & 2c^{q+1} \endmatrix
\right)$$ The determinant of
this block is $-(c-c^q)^2\gamma^2$; it belongs to $\n F_q^{*2}\iff
q\equiv 3 \mod 4$ or $q$ is even. Since we have $n$ blocks, we have:

$$\hbox{det $<*,*>_\alpha\not\in\n F_q^{*2}\iff q\equiv 1 \mod 4$ and
$n$ is odd.}$$
\medskip
{\bf Remark 7.8 (Jorge Morales).}
There is a theorem of Harder (see
e.g. W. Scharlau, "Quadratic and Hermitian forms", Springer-Verlag,
Berlin, 1985, Chapter 6, Theorem 3.3) that states that a unimodular form over
$k[X]$ \ \ --- \ \ $k$ being any field of characteristic not 2 --- is the
extension of a form over $k$, i.e. there is a basis in which all the
entries of the associated symmetric matrix are constant. This means
that the classification of the above quadratic forms over $\n F_q[T]$
($q$ odd) is very simple.
\medskip
{\bf Remark 7.9.} Let $M$ be a t-motive which is both negatively and positively self-dual. There is a natural idea 7.9.2 to define an analog of Hodge structure on $M$. Nevertheless, this idea fails. Namely, the exact sequence
$$0\to \Ker \vf \to L(M)\otimes \p \overset{\vf}\to{\to} \Lie(M)\to 0$$
is the functional field analog of an exact sequence for an abelian variety $A$:
$$ 0 \to H^{0,-1}(A) \to H^{-1}(A) \to (H^{1,0})^*(A) \to 0 $$
Hence, we can define $H^{0,-1}(M):= \Ker \vf$, and the problem is to define an analog of $H^{-1,0}(M)$.

Let us fix a negative isogeny $\alpha: M\to M'$, and let us extend the skew form $<*,*>_\alpha$ to $L(M)\otimes \p$ by $\p$-linearity. It is easy to check that $\Ker \vf$ is isotropic with respect to this form (there is an analogy with the number field case). Let us consider the following elementary lemma of linear algebra:
\medskip
{\bf Lemma 7.9.1.} Let $W$ be a vector space of dimension $2n$ over a field of characteristic $\ne 2$, $B^+$ (resp. $B^-$) a symmetric (resp. skew symmetric) non-degenerate bilinear form on $W$, and $W_0\subset W$ a subspace of dimension $n$ which is isotropic with respect to both $B^+$, $B^-$. Then almost always there exists the only $W_1\subset W$ of dimension $n$ having properties:
$$W_0\cap W_1=0; \ \ \ W_1 \hbox{  is isotropic with respect to both } B^+, \  B^-$$
where almost always means that entries of the matrices of $B^+$, $B^-$ in a basis of $W$ must not satisfy (at least one of) polynomial relations. $\square$
\medskip
If $\End_0(M)\ne \n F_q(T)$ and the action of $I_\alpha$ on $\End_0(M)$ is not identical, then there exists a positive isogeny $\beta: M\to M'$ and hence the symmetric form $<*,*>_\beta$ on $L(M)\otimes \p$. $\Ker \vf$ is isotropic with respect to  $<*,*>_\beta$. Let us fix $\beta$.
\medskip
{\bf Idea 7.9.2.} To apply Lemma 7.9.1 to this situation ($W=L(M)\otimes \p$, $W_0=\Ker \vf$, $B^+=<*,*>_\beta$, $B^-=<*,*>_\alpha$) in order to get a canonical subspace of $L(M)\otimes \p$ which is complementary to $\Ker \vf$ and hence can be considered as an analog of $H^{-1,0}(M)$.
\medskip
Clearly there is no complete analogy with the number field case. But the situation is even worse:
\medskip
{\bf Proposition 7.9.3.} For all $M$, $\alpha$, $\beta$ the "almost always" condition of Lemma 7.9.1 is not satisfied. $\square$
\medskip
{\bf 8. Relations between lattices and t-motives.}
\medskip
We have\footnotemark \footnotetext{I am grateful to Urs
Hartl
who indicated me this reference.}
\medskip
{\bf Theorem 8.1.} ([H],  Theorem 3.2). The dimension of the moduli set of pure t-motives of
dimension $n$ and rank $r$ is $n(r-n)$. $\square$
\medskip
{\bf Remark.} A tuple
$(e_1,...,e_r)$ of
integers entering in the statement of this theorem in [H] is
(0,...,0,1,...,1) with 0 repeated
$r-n$ times and 1 repeated $n$ times for the case under
consideration.
\medskip
Since this number $n(r-n)$ is equal to the dimension of the set of lattices of rank $r$
and dimension $n$, we can state an
\medskip
{\bf Open question 8.2.} Let $r$, $n$ be given. Let us consider the lattice map from the set of the pure uniformizable
t-motives of rank $r$ and dimension $n$ to the set
of lattices of rank $r$ and dimension $n$. Is it true that its image is open and the fibre at a generic point is discrete? If yes, what is the fibre?
\medskip
{\bf Remark.} Results of [L3] give some evidence that for the case $r=2n$ in a "neighborhood" of the $n$-th power of the rank 2 Carlitz module the fibre consists of 1 point.
\medskip
Theorem 5 implies that for $n=r-1$ the answer to 8.2 is yes (the below Proposition 11.8.5 shows that most likely the condition of purity is essential):
\medskip
{\bf Corollary 8.4.} All pure t-motives of
dimension $r-1$ and rank $r$ having $N=0$ are uniformizable. There is a 1 -- 1 functorial
correspondence between pure t-motives of dimension $r-1$ and rank $r$ having $N=0$ ($r\ge
2$), and lattices of rank $r$ in $\p^{r-1}$ having dual.
\medskip
{\bf Proof.} Let $L$ be a lattice of rank $r$ in $\p^{r-1}$ having dual $L'$. There
exists the only Drinfeld module $M'$ such that $L(M')=L'$, and let $M$ be its dual.
Theorem 5 implies that $L(M)=L$. If there exists another pure t-motive $M_1$ of
dimension $r-1$ and rank $r$ having $N=0$ such that $L(M_1)=L$ then by Corollary 10.4 (its proof is logically independent: there is no vicious circle) the
dual $M'_1$ is a Drinfeld module, according Theorem 5 it satisfies $L(M'_1)=L'$, hence
$M'_1=M'$ and hence $M_1=M$. $\square$
\medskip
{\bf Remark 8.5.} Recall that lattices of rank $r$ in $\p^{r-1}$ having dual are
described in 3.5 (formulas 3.6, 3.7).
We see that for the case $n=r-1$, $N=0$ purity implies uniformizability. We have
\medskip
{\bf Question 8.5a.} Do exist non-uniformizable t-motives having $n=r-1$, $N=0$?
\medskip
{\bf Question 8.5b.} Do exist uniformizable t-motives having $n=r-1$, $N=0$ such that its lattice has no dual? (Clearly this is a subquestion of 8.2).
\medskip
{\bf Remark 8.6.} Clearly for any $r$, $n$ we have: if a lattice $L$ of rank $r$ and dimension $n$ has no dual then $L\ne L(M)$ for any pure uniformizable $M$. I do not know whether Theorem 6 (which is an analog of Theorem 5 for another tensor operation) imposes a more strong similar restriction on the property of $L$ to be the $L(M)$ of some pure uniformizable $M$, or not.
\medskip
Further, for any uniformizable t-motive $M$ we have a
\medskip
{\bf Corollary 8.7.} If the dual of $(L(M), \Lie(M))$ does not exist then the
dual of $M$ does not exist. Example: the Carlitz module.
\medskip

{\bf 9. Main theorem in terms of Pink-Hodge structure.}
\medskip
Let us consider a version of a special case of the general definition of  Pink-Hodge structure ([P], 0.2; 9.1).
\medskip
{\bf Definition.} A Pink-Hodge structure of constant weight and complete dimension is a pair $\underline{H}=(H,\goth q_H)$ where $H$ is a free finite dimensional $\w$-module and $\goth q_H$ is a $\p[[T-\theta]]$-lattice in $H\underset{\w}\to{\otimes}\p[[T-\theta]]$ such that the dimension of $\goth q_H$ over $\p[[T-\theta]]$ is equal to the dimension of $H$ over $\w$ (condition of complete dimension).
\medskip
Let $\vf: L \hookrightarrow \p^n$ be a lattice. It defines a Pink-Hodge structure $\underline{H}=\underline{H}(L)$ of constant weight and complete dimension. Firstly, instead of a $\n F_q[\theta]$-module $L$ we consider an isomorphic $\w$-module $H$ formally defined by the property $H\underset{\w}\to{\otimes}\n F_q[\theta]=L$ where the map $\w \to \n F_q[\theta]$ is $\iota$. We denote the isomorphism $H \to L$ by $\iota$ as well; the composition $\vf \circ \iota: H \to \p^n$ is a map of $\w$-modules where $T\in \w$ acts on $ \p^n$ by multipication by $\theta$. Further, $\vf \circ \iota$ extends to a surjection of $\p[[T-\theta]]$-modules $H\underset{\w}\to{\otimes}\p[[T-\theta]] \to \p^n$ denoted by $\vf \circ \iota$ as well.  Finally, $\goth q_H$ is defined as $\Ker \vf \circ \iota$.

If $M$ is a pure uniformizable t-motive then we associate it a Pink-Hodge structure of constant weight and complete dimension $\underline{H}(M)=\underline{H}(L(M))$.

Let $m=m(\underline{H})$ be the minimal number such that $\goth q_H \supset (T-\theta)^m H\underset{\w}\to{\otimes}\p[[T-\theta]]$. For $\mu\ge m$ we define the $\mu$-dual structure ${\underline{H}'}^{\mu}=({H'}^{\mu}, \goth q_{{H'}^{\mu}})$ as follows:

$${H'}^{\mu}=H^*, \ \  \goth q_{{H'}^{\mu}}=\{\chi\in H^*\underset{\w}\to{\otimes}\p[[T-\theta]]  $$ $$\hbox{ such that }  \forall y\in \goth q_H  \hbox{ we have } \chi(y)\in (T-\theta)^\mu\p[[T-\theta]]\} $$ It is obvious that it is really a Pink-Hodge structure of constant weight and complete dimension.
\medskip
If $\underline{H}=\underline{H}(L)$ for a lattice $L$ then $m=1$ and if $L$ has dual then $${\underline{H}'}^1=\underline{H}(L')\eqno{(9.1)}$$ this is easy to prove.
\medskip
{\bf Remark 9.2.} And if $L$ has no dual? Really, $\underline{H}(L)$ exists even if $L$ does not satisfy Definition 2.1 (b). If $L$ is a lattice having no dual this means that $L'$  does not satisfy Definition 2.1 (b). Nevertheless, equality ${\underline{H}'}^1=\underline{H}(L')$ is meaningful and holds. We are not interested in these lattices because they cannot be lattices of uniformizable t-motives having dual.
\medskip
An analog of Theorem 6 for dual t-motives is the following. Let $M$ be as above. Obviously $m=m(\underline{H}(M))$ is the minimal number such that $N^m=0$. According Theorem 10.3, ${M'}^{m}$ exists, and it is pure.
\medskip

{\bf 10. Duals of pures, and other elementary results. }

\nopagebreak
\medskip
We consider in this section the case of arbitrary $N$ (i.e. not necessarily $N=0$), and $\w=\n F_q[T]$. The definition 1.8 extends to the case of pr\'e-t-motives, and
remarks 1.11 hold for this case.
\medskip
{\bf Lemma 10.2. } Let $M$ be a pr\'e-t-motive, $m=m(M)$ from its
(1.3.1), and $\mu\ge m$.
Then $M'$ --- the $\mu$-dual of $M$ --- exists as a pr\'e-t-motive,
and $m(M')\le\mu$. If
$M'$ is a t-motive then $\dim M'=r\mu - \dim M$ ($r$ is the rank of $M$).
\medskip
{\bf Proof.} We must check that $Q'$ has no denominators, and the
condition (1.3.1). The
module $\tau M$ is a $\p[T]$-submodule of $M$ (because $a \tau x =
\tau a^{1/q} x$ for
$x\in M$),
hence there are $\p[T]$-bases $f_*=(f_1, ... f_r)^t$, $g_*=(g_1, ...
g_r)^t$ of $M$,
$\tau M$ respectively such that $g_i=P_if_i$, where $P_1 | P_2 | ... |
P_r$, $P_i \in
\p[T]$. Condition (1.3.1) means that $\forall i$ \ $(T-\theta)^m f_i
\in \tau M$, i.e.
$P_i|(T-\theta)^m$, i.e. $\forall i$ \ $P_i=(T-\theta)^{m_i}$ where $0
\le m_i \le
m_{i+1} \le m$. There exists a matrix $\goth Q=\{\goth q_{ij}\}\in
M_r(\p[T])$ such that
$$\tau f_i = \sum_{j=1}^r \goth q_{ij}g_j= \sum_{j=1}^r \goth q_{ij} P_j
f_j\eqno{(10.2.1)}$$
Although $\tau$ is not a linear operator, it is easy to see that
$\goth Q\in GL_r(\p[T])$
(really, there exists $C=\{c_{ij}\}\in M_r(\p[T])$ such that $g_i=P_i
f_i=\tau (\sum_{j=1}^r
c_{ij}f_j)$, we have $C^{(1)}\goth Q =E_r$).

We denote the matrix $\diag(P_1, P_2, ... ,P_r)$ by $\goth P$, so
(10.2.1) means that
$$Q= \goth Q \goth P \eqno{(10.2.2)}$$

{\bf Remark 10.2.3.} Since
$\goth Q \goth P \in
GL_r(\p(T))$, we get that the action of $\tau$ on $i_2(M)$ is invertible.
\medskip
It is clear that if $M$ is a t-motive then
$$\dim M = \sum_{j=1}^r m_j\eqno{(10.2.4)}$$
(because $\dim M = \dim_{\p}(M/\tau M)$. Further, (10.2.2) implies that
for $Q'=Q({M'})$
we have
$$Q'=\goth Q^{t-1}\diag((T-\theta)^{\mu-m_1}, ...,
(T-\theta)^{\mu-m_r})\eqno{(10.2.5)}$$
This means that elements of $Q'$ have no denominators. The condition
(1.3.1) for $M'$
follows easily from (10.2.5) (because $\goth Q^{t-1}\in GL_r(\p[T])$),
and the dimension
formula (for the case $M'$ is a t-motive) follows immediately from
(10.2.4) applied to
$M'$. $\square$
\medskip
A definition of a pure t-motive can be found in [G] ((5.5.2),
(5.5.6) of [G] + formula (1.3.1) of the present paper).
\medskip
{\bf Theorem 10.3.} Let $M$ be a pure t-motive and $m=m(M)$ from
(1.3.1). Then (if
$rm-n>0$) its $m$-dual $M'$ exists, and it is pure.
\medskip
{\bf Proof.} The definition of pure ([G], (5.5.2)) is valid for
pr\'e-t-motives. We use
its following
matrix form. We denote $T^{-1}$ by $S$ and for any $C$ we let
$$C^{[i]}=C^{(i-1)}\cdot C^{(i-2)}\cdot ... \cdot C^{(1)}\cdot C$$

{\bf Lemma 10.3.1.} Let $Q\in M_r(\p[T])$ be a matrix such that formula
(1.9.3) defines an
t-motive $M$. Then it is pure iff there exists $C\in
GL_r(\p((S)) \ )$ such that
for some $\goth q$, $s>0$
$$S^\goth q C^{(s)}Q^{[s]}C^{-1}\in GL_r(\p[[S]])$$
i.e. iff $S^\goth q C^{(s)}Q^{[s]}C^{-1}$ is $S$-integer and its inicial
coefficient is
invertible.
\medskip
{\bf Proof.} Elementary matrix calculations. We take $C$ as a matrix
of base change of
$f_*$ to a $\p[[S]]$-basis of $W$ of (5.5.2) of [G]. $\square$
\medskip
{\bf Lemma 10.3.2.} Let $\mu = m$. We have: $M'={M'}^\mu$ of Lemma 10.2 is a pure
pr\'e-t-motive.
\medskip
{\bf Proof.} Let $\goth q$, $s$ and $C$ be from Lemma 10.3.1. We have
$${Q'}^{[s]}=((T-\theta)^{[s]})^\mu Q^{[s]\ t-1}$$
(we use (1.2)). We take $C'=C^{t-1}$. We have
$$S^{s\mu-\goth q} {C'}^{(s)}{Q'}^{[s]}{C'}^{-1}=$$
$$=S^{s\mu-\goth q}C^{(s)\ t-1}Q^{[s]\ t-1}((\frac1S-\theta)^{[s]})^\mu C^t$$
$$=((1-S\theta)^{[s]})^\mu S^{-\goth q}C^{(s)\ t-1}Q^{[s]\ t-1}C^t$$
$$=((1-S\theta)^{[s]})^\mu (S^\goth q C^{(s)}Q^{[s]}C^{-1})^{t-1}$$
We have: $\goth q/s=n/r$ ([G], (5.5.6)), hence $(s\mu-\goth q)/s=(r\mu-n)/r$ and
$s\mu-\goth q>0$. Further,
$((1-S\theta)^{[s]})^\mu\in GL_r(\p[[S]])$, and the result follows from
Lemma 10.3.1. $\square$
\medskip
{\bf Remark.} This result holds also for $\mu>m$.
\medskip
The theorem 10.3 follows from Lemma 10.2, the above lemmas and the
proposition that a pure
pr\'e-t-motive satisfying (1.3.1) is a t-motive ([G], (5.5.6),
(5.5.7)). $\square$
\medskip
{\bf Corollary 10.4.} Let $M$ be a t-motive such that $m=1$,
$n=r-1$. Then $M$
has dual $\iff$ $M$ is pure $\iff$ $M$ is dual to a Drinfeld module.
\medskip
{\bf Proof.} Dimension formula shows that $M'$ (if it exists) is a
Drinfeld module, and they are all pure. $\square$
\medskip
{\bf Example 10.5.} Let $M$ be given by (notations of 1.9.1)
$$\goth A_0=\theta E_2, \ \ \goth A_{1}= \left(\matrix a_{111} & 0
\\ a_{121} & 1 \endmatrix \right), \ \ \goth A_2= \left(\matrix 1 & 0
\\0 & 0 \endmatrix \right)$$
This $M$ has $m=1$, $n=2$, $r=3$, and it is easy to see that it has no dual. Really, for
this $M$ we have (notations of 1.9) $f_1=e_1$, $f_2=\tau e_1$, $f_3=e_2$,
$$Q=\left(\matrix 0&1&0\\T-\theta&-a_{111} & 0
\\ 0&-a_{121} & t-\theta \endmatrix \right), \ \ Q'=\left(\matrix a_{111} & t-\theta
& a_{121} \\ 1 &0&0\\0&0&1\endmatrix \right)$$ The last line of $Q'$ means that $\tau
f_3'=f_3'$. This is a contradiction to the property that $M'_{\p[\tau]}$ is free. It is
possible also to show (Proposition 11.3.4) that $M$ is not pure, and to use 10.4 in order to
prove that it has no dual.
\medskip
Later (Section 11) we shall construct examples of non-pure abelian
t-motives which have dual. Considerations of 11.8 predict that there is enough such t-motives.
\medskip
{\bf Theorem 10.6. } For any t-motive $M$ there exists $\mu_0$ such that for all
$\mu\ge\mu_0$ the object
${M'}^{\mu}$ exists as a t-motive. For these $\mu$ we have
$${M'}^{\mu+1}={M'}^{\mu}\otimes \goth
C\eqno{(10.6.1)}$$

{\bf Proof.} (10.6.1) holds at the level of pr\'e-t-motives, because $Q(\goth
C)=(T-\theta)E_1$. According [G], Lemma 5.4.10 it is sufficient to prove that
${M'}^{\mu}$ is finitely generated as a $\p [\tau]$-module. We shall
use notations of
Lemma 10.2. We take $$\mu_0=1+\hbox{ \{the maximum of the degrees of
entries of $\goth
Q(M)$ as polynomials in $T$\} } $$ $$+ \max(m_k)$$
Let $f'_1, ...f'_r$ be the basis of ${M'}^{\mu}$ over $\p [T]$ dual to
$f_1, ...f_r$. It
is sufficient to prove the
\medskip
{\bf Lemma 10.6.2.} Let $i_0=\mu-\min(m_k)$. Then elements $T^if'_j$,
$i<i_0$, $j=1,
....,r$, generate ${M'}^{\mu}$ as a $\p [\tau]$-module.
\medskip
{\bf Proof of the lemma.} By induction, it is sufficient to show that
for all $\alpha\ge
i_0$ the equation
$$\tau x = (T-\theta)^\alpha f'_j\eqno{(10.6.3)}$$
(equality in ${M'}^{\mu}$) has a solution
$$x=\sum_{k=1}^r C_kf'_k$$
where $C_k \in \p [T]$, $\deg(C_k)< \alpha$. According (10.2.5), the
solution to (10.6.3)
is given by
$$(C_1^{(1)}, ... , C_r^{(1)})=(0,... 0,
(T-\theta)^{\alpha-\mu+m_j},0,... 0)\goth Q^t$$
(the non-0 element of the row matrix is at the $j$-th place).
Unequalities satisfied by
$\mu$ and $\alpha$ show that all $C_k^{(1)}$ are polynomials of degree
$<\alpha$. Since
$c \mapsto c^q$ is surjective on $\p$, we get the desired. $\square$
\medskip
{\bf 10.7. Virtual t-motives.} \footnotemark \footnotetext{This
notion was
indicated me by
Taguchi.} We need two elementary lemmas.

\nopagebreak
\medskip
{\bf Lemma 10.7.0.}\footnotemark \footnotetext{Anderson proved (not published) that the tensor product of any t-motives is also a t-motive.} If $M$ is a t-motive then $M\otimes \goth C$ is also a t-motive.
\medskip
{\bf Proof.} Let $f_{j}$ ($j=1,...,r$) be a $\p[T]$-basis of
$M_{\p[T]}$ and $\goth f$
from 1.10.2, so $f_{j}\otimes \goth f$ is a $\p[T]$-basis of $(M\otimes \goth
C)_{\p[T]}$. It is sufficient to prove that $(M\otimes \goth
C)_{\p[\tau]}$ is finitely
generated. Since $M_{\p[\tau]}$ is finitely generated, it is easy to
see that there
exists $a$ such that elements
$$(T-\theta)^i f_j, \ \ i=0, ... , a, \ \ j=1, ... ,r$$
generate $M_{\p[\tau]}$. This means that $\forall j=1, ... ,r$ there
exist $c_{ijkl}\in
\p$ such that
$$(T-\theta)^{a+1} f_j=\sum_{i=0}^a \sum_{k=0}^\gamma\sum_{l=1}^r
c_{ijkl}(T-\theta)^i
\tau^k f_l\eqno{(10.7.0.1)}$$
where $\gamma$ is a number.

Let us multiply (10.7.0.1) by $(T-\theta)^\gamma$. Taking into
consideration the formula
of the action of $\tau$ on $M\otimes \goth C$ we get that the result
gives us the
following formula in $M\otimes \goth C$:
$$(T-\theta)^{a+\gamma+1} f_j\otimes \goth f=\sum_{i=0}^a
\sum_{k=0}^\gamma\sum_{l=1}^r
c_{ijkl} (T-\theta)^{i+\gamma-k} \tau^k \cdot (f_l\otimes \goth f)\eqno{(10.7.0.2)}$$
This proves that for all $j$ the element $(T-\theta)^{a+\gamma+1} f_j\otimes \goth f$ is
a linear combination of
$$(T-\theta)^i f_l\otimes \goth f, \ \ i=0, ... , a+\gamma, \ \ l=1, ...
,r\eqno{(10.7.0.3)}$$
in $(M\otimes \goth C)_{\p[\tau]}$.
Multiplying (10.7.0.2) by consecutive powers of $T-\theta$ we get by
induction that
elements of 10.7.0.3 generate $(M\otimes \goth C)_{\p[\tau]}$. $\square$
\medskip
{\bf Lemma 10.7.1.} If $M_1\otimes \goth C$ is isomorphic to
$M_2\otimes \goth C$ then
$M_1$ is isomorphic to $M_2$.
\medskip
{\bf Proof.} Let $f_{i*}$ ($i=1,2$) be a $\p[T]$-basis of
$(M_i)_{\p[T]}$, $Q_i$ from
1.9.3, $\alpha: M_1\otimes \goth C\to M_2\otimes \goth C$ an
isomorphism and $C\in
GL_r(\p[T])$ the matrix of $\alpha$ in $f_{1*} \otimes \goth f $,
$f_{2*}\otimes \goth f
$. The matrix of the action of $\tau$ on $M_i\otimes \goth C$ in the
base $f_{i*} \otimes
\goth f $ is $(T-\theta)Q_i$, and the condition that $\alpha$ commutes with
multiplication by $\tau$ is
$$(T-\theta)Q_1 C = C^{(1)}(T-\theta)Q_2$$
Dividing this equality by $T-\theta$ we get that the map $\alpha_0$
from $M_1$ to $M_2$
having the same matrix $C$ in the bases $f_{i*}$, commutes with
$\tau$, i.e. defines an
isomorphism from $M_1$ to $M_2$. $\square$
\medskip
Using Lemma 10.7.1 we can state the following
\medskip
{\bf Definition.} A virtual t-motive is an object $M\otimes \goth C^{\otimes
\mu}$ where $M$ is
a t-motive
and $\mu\in \n Z$, with the standard equivalence relation (here
$\mu_1\ge\mu_2$):
$$M_1\otimes \goth C^{\otimes \mu_1}=M_2\otimes \goth C^{\otimes \mu_2}\iff
M_2=M_1\otimes \goth C^{\otimes (\mu_1-\mu_2)}$$
$$\iff \exists \mu \hbox{ such that $\mu+\mu_1\ge 0$, $\mu+\mu_2\ge 0$
and } M_1\otimes
\goth C^{\otimes (\mu+\mu_1)}=M_2\otimes \goth C^{\otimes (\mu+\mu_2)}$$
Lemma 10.7.1 shows that these conditions are really equivalent.
\medskip
{\bf Corollary 10.7.2.} The $\mu$-dual of a virtual t-motive is
well-defined and
always exists as a virtual t-motive. $\square$
\medskip
{\bf Proposition 10.8.} The following formula is valid at the level of
pr\'e-t-motives: for
any $\mu_1$, $\mu_2$, if ${M_1'}^{\mu_1}$, ${M_2'}^{\mu_2}$ exist then
${(M_1\otimes
M_2)'}^{(\mu_1+\mu_2)}$ exists and
$${(M_1\otimes M_2)'}^{(\mu_1+\mu_2)}={M_1'}^{\mu_1} \otimes {M_2'}^{\mu_2}$$
{\bf Proof.} This is a functorial equality; also we can check it by
means of elementary
matrix calculations. $\square$
\medskip
{\bf Proposition 10.9.} Let $P\in \hbox{\bf A}$ be an irreducible
element. The Tate
module
$T_P({M'}^{\mu})$ is equal to $$T_P(\goth C)^{\otimes \mu}\otimes
\widehat{T_P(M)}$$
(equality of Galois modules) where $\widehat{T_P(M)}$ is the dual
Galois module.
\medskip
{\bf Proof.} It is completely analogous to the proof of the corresponding theorem for
tensor products
([G], Proposition 5.7.3, p. 157).
All modules in the below proof will be the Galois modules, and equalities of modules will
be equalities of Galois
modules. Recall that $E=E(M)$. Since $T_P(M)=\invlim_n E_{P^n}$, it is sufficient to prove that for any
$a\in \w$ we have $E({M'}^{\mu})_a=E(\goth C^{\otimes \mu})_a\otimes \hat E_a$, where $\hat E_a$ is the dual of $E_a$ in the meaning of [T], Definition 4.1.
We have the following sequence of equalities of modules:
$${M'}^{\mu}/a{M'}^{\mu}=\Hom_{\p[T]}(M/aM, \goth C^{\otimes \mu}/a\goth C^{\otimes
\mu})\eqno{(10.9.2)}$$
such that the action of $\tau$ on both sides of this equality coincide (to define the
action of $\tau$ on the
right and side of (10.9.2) we need the action of $\tau^{-1}$ on $M/aM$; it is
well-defined, because the determinant
of the action of $\tau$ on $M$ is a power of $T-\theta$, hence its image in
$\p[T]/a\p[T]$ is invertible). 10.9.2
follows immediately from the definition of ${M'}^{\mu}$;
$$({M'}^{\mu}/a{M'}^{\mu})^\tau=\Hom_{\n F_q[T]}((M/aM)^\tau,
(\goth C^{\otimes \mu}/a\goth C^{\otimes \mu})^\tau)\eqno{(10.9.3)}$$
This follows from 10.9.2 and the Lang's theorem
$$\goth M/a\goth M=(\goth M/a\goth M)^\tau\underset{\n F_q[T]/a\n
F_q[T]}\to{\otimes}\p[T]/a\p[T]$$
applied to both $\goth M=M$, $\goth M={M'}^{\mu}$ (we use that both $M$, ${M'}^{\mu}$ are
free $\p[T]$-modules).
Finally, we have a formula
$$E(\goth M)_a=\Hom_{\n F_q}((\goth M/a\goth M)^\tau, \n F_q)$$
([G], p. 152, last line of the proof of Proposition 5.6.3). Applying this formula to 10.9.3 we get the desired. $\square$
\medskip
{\bf 11. An explicit formula.}

\nopagebreak
\medskip
We return to the case $N=0$. Let $e_*$, $\goth A$, $\goth A_i$, $l$, $n$ be from (1.9). We consider in the
present section two simple types of t-motives (called standard-1 and standard-2
t-motives respectively) whose $\goth A_i$ have a row echelon form, and we give an
explicit formula for the dual of some standard-1
t-motives. Analogous formula can be easily obtained for more general types of
t-motives. These results are the
first step of the problem of description of all t-motives
having duals.
\medskip
{\bf 11.1.} For the reader's convenience, we give here the definition of standard-1
t-motives
for the case $n=2$ (here $\lambda_1$ and $\lambda_2 $ satisfying
$\lambda_1=l$,
$l>\lambda_2 \ge 2$ are parameters):

$$\goth A_0=\theta E_2, \hbox{ for } 0 < i < \lambda_2 \ \ \goth A_i \hbox{ is
arbitrary, } $$
$$\goth A_{\lambda_2}= \left(\matrix * & 0 \\ * & 1 \endmatrix \right),
\hbox{ for } \lambda_2
< i < l \ \ \goth A_i= \left(\matrix * & 0 \\ * & 0 \endmatrix \right), \ \
\goth A_l= \left(\matrix
1 & 0
\\0 & 0 \endmatrix \right) $$

{\bf 11.2.} To define standard-2 t-motives of dimension $n$, we need to fix
\medskip
1. A permutation $\vf\in S_n$, i.e. a 1 -- 1 map $\vf: (1, ..., n) \to (1, ..., n)$;
\medskip
2. A function $k: (1, ..., n) \to \n Z^+$ where $\n Z^+$ is the set of integers $\ge 1$.
\medskip
{\bf Definition.} A standard-2 t-motive of the type $(\vf, k)$ is an abelian
t-motive of dimension $n$ given by the formulas ($i=1,...,n$):
$$Te_{\vf(i)}= \theta e_{\vf(i)}
+\sum_{\alpha=1}^n\sum_{j=1}^{k(\alpha)-1}a_{j,\vf(i),\alpha}\
\tau^je_\alpha + \tau^{k(i)}e_i\eqno{(11.2.1)}$$
where $a_{j,\vf(i),\alpha}\in\p$ is the $(\vf(i),\alpha)$-th entry of the matrix $\goth A_j$.
\medskip
{\bf Proposition 11.2.2.} Formula 11.2.1 really defines a t-motive denoted by
$M=M(\vf,k)=M(\vf,k,a_{***})$. Its rank is $\sum_{\alpha=1}^n k(\alpha)$ and elements
$X_{\alpha j}:=\tau^je_\alpha$, $\alpha=1, ..., n$, $j=0, ..., k(\alpha)-1$, form its
$\p[T]$-basis.
$\square$
\medskip
The group $S_n$ acts on the set of types $(\vf, k)$ and on the set of the above $M$;
clearly for any $\psi\in S_n$ we have $\psi(M)$ is isomorphic to $M$. Particularly, we
can consider only $\vf$ of the following form of the product of $i$ cycles ($\alpha_0=0,
\ \alpha_i=n$):
$$\vf=(\alpha_0+1, ..., \alpha_1)(\alpha_1+1,...,\alpha_2)
...(\alpha_{i-1}+1,...,\alpha_i)\eqno{(11.2.3)}$$
(standard notation of the theory of permutations, for $\gamma\ne \alpha_j$ we have
$\vf(\gamma)=\gamma+1$, for $\gamma= \alpha_j$ we have  $\vf(\alpha_j)=\alpha_{j-1}+1$).
\medskip
{\bf Example 11.2.4.} Let $\vf$ be defined by 11.2.3, the quantity of cycles $i$ is equal
to $1$ and all $a_{***}=0$. Then the corresponding $M$ is of complete multiplication by a
CM-field $\n F_{q^r}(T)$ and its CM-type $\Phi$ is $\{\Id, \fr^{k(1)}, \fr^{k(1)+k(2)},
..., \fr^{k(1)+k(2)+...+k(n-1)}\}$ where $\fr$ is the Frobenuis homomorphism $\n
F_{q^r}\to \bar \n F_q$ (see 13.3, first case: formulas 13.3.1, 13.3.2 coinside with 11.2.1
for the given $\vf$ and $a_{***}=0$; $i_j$ of 13.3.0 is $k(1)+k(2)+...+k(j-1)$ of the
present notations).
\medskip
{\bf Definition 11.3.} A standard-1 t-motive is a standard-2 t-motive
whose $\vf$ is the identical permutation $Id$.
\medskip
{\bf 11.3.0.} Let $M=M(Id,k)$ be a standard-1 t-motive. Acting
by $\psi\in S_n$ we can consider only the case of non-increasing $k(j)$. We introduce a
number $\goth m\ge1$ --- the quantity of jumps of $k(j)$, and two sequences
$$0=\gamma_0< \gamma_1<...< \gamma_\goth m=n$$
(sequence of arguments of points of jumps of the function $k$) and
$$0=\lambda_{\goth m+1}< \lambda_{\goth m}<...< \lambda_2 <\lambda_1=l$$
(sequence of values of $k$ on segments $[\gamma_{i-1}+1, ...,\gamma_i]$) by the formulas
$$\matrix k(1)=...=k(\gamma_1)=\lambda_1 \\ \\
k(\gamma_1+1)=...=k(\gamma_2)=\lambda_2\\ ... \\
k(\gamma_{\goth m-1}+1)=...=k(\gamma_\goth m)=\lambda_\goth m\endmatrix \eqno{(11.3.1)}$$
\medskip
{\bf Example 11.3.2.} The t-motive $M$ of 11.1 is a standard-1 having $\goth m=2$,
$\gamma_1=1$, $\gamma_2=2$ and $\lambda_1$, $\lambda_2$ as in 11.1. Its rank $r=\lambda_1+\lambda_2$.
\medskip
{\bf Conjecture 11.3.3.} A standard-2 t-motive of the type $(\vf, k)$ (notations of
11.2.3) is pure iff $\forall j=1, ..., i$ we have:
$$\frac{\alpha_j-\alpha_{j-1}}{\sum_{\gamma=\alpha_{j-1}+1}^{\alpha_{j}}k(\gamma)} = \frac{n}{r}$$
This conjecture is obviously true if all $a_{***}$ are 0.
\medskip
To simplify exposition, we prove here only the following particular case of this
conjecture.
\medskip
{\bf Proposition 11.3.4.} Let $M$ be a standard-1 t-motive having $\goth m>1$,
defined over $\n F_q(\theta)$, having a good reduction at a point of degree 1 of $\n
F_q(\theta)$ (i.e. a point $\theta+c$, $c\in \n F_q$). Then $M$ is not pure.
\medskip
{\bf Proof.} Let $M$ be defined by 11.2.1, we use notations of 11.3.1. We consider the
action of Frobenius on $\tilde M$ --- the reduction of
$M$ at $\theta+c$. According [G], Theorem 5.6.10, it is sufficient to prove that orders
of the roots of the characteristic polynomial of Frobenius over $\w$ are not equal. More
exactly, we consider the valuation infinity on $\w$ (defined by the condition
$\ord(T)=-1$); the order corresponds to a continuation of this valuation to $\End(\tilde
M)$. The
action of Frobenius on $\tilde M$ coincides with multiplication by $\tau$, because the
degree of the reduction point is 1.

A basis $f_*$ of $M_{\p[T]}$ is the set of $X_{\alpha j}:=
\tau^je_{\alpha}$ of 11.2.2. The matrix
$Q(M)$ is defined by the following formulas for the action of
$\tau$ on $X_{\alpha j}$:
$$\tau(X_{\alpha j})=X_{\alpha,j+1} \hbox{ if } j <
k(\alpha)-1\eqno{(11.3.4.1)}$$
$$\tau(X_{\alpha,k(\alpha)-1})=TX_{\alpha,0} -
\sum_{\delta=1}^\goth m
\sum_{d=\lambda_{\delta+1}}^{\lambda_{\delta}-1}
\sum_{c=1}^{\gamma_{\delta}} a_{d\alpha c}X_{cd} \eqno{(11.3.4.2)}$$
This means that if we arrange $X_{\alpha j}$ in lexicographic order ($X_{\alpha_1 j_1}$ precedes to $X_{\alpha_2 j_2}$ if $\alpha_1 <
\alpha_2$) then the matrix
$Q(M)$ has the block form: $$Q(M)=(C_{ij})\ \ \ (i,j=1,...,n)$$ where $C_{ij}$ is a
$k(i)\times k(j)$-matrix of the form
$$C_{ii}=\left(\matrix 0&1&0&...&0\\0&0&1&...&0\\ ...&...&...&...&... \\ 0&0&0&...&1 \\
T-\theta & *&*&...&*\endmatrix \right), \ \ C_{ij}=\left(\matrix 0&0&...&0\\
...&...&...&...\\ 0&0&...&0 \\ 0& *&...&*\endmatrix \right) (i\ne j)$$
where asterisks mean elements $a_{***}$ (in some order). We consider the characteristic polynomial $P(X)\in (\p[T])[X]$ of $Q(M)$. We
have $$C_{ii}-XE_{k(i)}=\left(\matrix -X&1&0&...&0\\0&-X&1&...&0\\ ...&...&...&...&... \\
0&0&0&...&1 \\ t-\theta & *&*&...&*-X\endmatrix \right)$$

A subset of the set of entries of a matrix is called (following N.N.Luzin) a lightning if each row and each column of the matrix contains exactly one element of this subset. The product of elements of a lightning is called the value of this lightning (i.e. the determinant is the alternating sum of the values of all lightnings).
\medskip
{\bf Lemma 11.3.4.3.} If a non-zero lightning of $C_{ii}-XE_{k(i)}$ contains the term
$T-\theta$, then it does not contain any term containing $X$. $\square$
\medskip
Let $J$ be a subset of the set $1,...,n$ and $J'$ its complement.
\medskip
{\bf Corollary 11.3.4.4.} If a non-zero lightning of $Q(M)-XE_{r}$ contains terms
$T-\theta$ of blocks $C_{\al}$, $j\in J$, then its value is a polynomial in $X$ of degree
$\le \sum_{j'\in J'}k(j')$, and there exists exactly one such lightning (called the
principal $J$-lightning) whose value is a polynomial in $X$ of degree $\sum_{j'\in
J'}k(j')$. $\square$

Since the characteristic polynomial of Frobenius of $\tilde M$ is $\tilde P$
(respectively the valuation infinity of $\p[T]$), it is sufficient to prove that the
Newton polygon of $P(X)$ is not reduced to the segment $((0,-n); (r,0))$ defined by its
extreme terms $(T-\theta)^n$ and $X^r$. To do it, it is sufficient to find a point on its
Newton polygon which is below this segment. We consider $J_{min}=$ the set of all
$\gamma_\goth m - \gamma_{\goth m-1}$ diagonal blocks $C_{ii}$ ($i=\gamma_{\goth m-1}+1,
..., \gamma_\goth m$) of $Q(M)$ of minimal size $\lambda_\goth m$. The value of the
principal $J_{min}$-lightning is $(T-\theta)^{\gamma_\goth m - \gamma_{\goth m-1}}$ times
polynomial in $X$ of degree $d:=r-(\gamma_\goth m - \gamma_{\goth m-1})\lambda_\goth m$.
Corollary 11.3.4.4 implies that if the value of any other lightning of $Q(M)-XE_r$ contains a
term whose $X$-degree is equal to $d$, then the $T$-degree of this term is strictly less
than $\gamma_\goth m - \gamma_{\goth m-
 1}$. This means that if we write $P(X)=\sum_{i=0}^r C_iX^i$, $C_i\in \p[T]$, then
$\ord_\infty(C_d)= -(\gamma_\goth m - \gamma_{\goth m-1})$, i.e. the point with
coordinates $[-(\gamma_\goth m - \gamma_{\goth m-1}), d]$ belongs to the Newton diagram
of $P(X)$, i.e. it is above (really, at) the Newton polygon of $P(X)$. This point is
below the segment $((0,-n); (r,0))$. $\square$
\medskip
{\bf Remark 11.3.4.5.} It is easy to see that the Newton polygon of $P(X)$ coincides with
the Newton polygon of the direct sum
of trivial Drinfeld modules of ranks $\lambda_*$, i.e. with the Newton polygon of the
polynomial
$$\prod_{i=1}^{\goth m} (X^{\lambda_i}-T)^{\gamma_i - \gamma_{i-1}}$$
\medskip
{\bf 11.4.} To formulate the below theorem 11.5 we need some notations. Let $M$ be a
standard-1 t-motive defined by formulas 11.2.1, 11.3.1. We impose the condition
$\lambda_\goth m \ge 3$. Theorem 11.5 affirms that it has dual. To find explicitly the
dual of $M$, we need to choose an arbitrary function
$\nu: (i,j) \to \nu(i,j)$ which is a 1 - 1 map from the set of pairs
$(i,j)$ such that

$$1 \le i \le n; \ \ 1 \le j \le k(i)-2 \eqno{(11.4.1)}$$
to the set $[n+1, ..., r-n]$ (recall that $r=\sum_{i=1}^n k(i)= \sum_{i=1}^\goth m
(\gamma_{i}-\gamma_{i-1})\lambda_i$).
\medskip
Let the $(r-n)\times(r-n)$-matrices $B_1$, $B_2$ be defined by the
following formulas
(here and until the end of the proof of 11.5 we have $i ,\alpha = 1, ... , n$; \ \
$b_{\beta\gamma\delta}$
is the $(\gamma\delta)$-th entry of $B_\beta$, all
entries of $B_1$,
$B_2$ that are not in the below list are 0):
\medskip
{\bf 11.4.2.} $b_{1i\alpha}= - a_{k(i)-1,\alpha,i}$;
\medskip
$b_{1,\nu(i,j),\alpha} = - a_{j,\alpha,i}$ for $1 \le j \le k(i)-2$;
\medskip
$b_{1,\nu(i,j+1),\nu(i,j)}=1$ for $1 \le j \le k(i)-3$;
\medskip
$b_{1,i,\nu(i,k(i)-2)} =1$;
\medskip
$b_{2,\nu(i,1),i}=1$.
\medskip
We let $B=\theta E_{r-n}+B_1\tau+B_2\tau^2$ and consider a t-motive $M(B)$ (see 11.5.1
below). Formulas
11.4.2 mean that $M(B)$ is standard-2, its $\vf=\vf_B$ is a product of $n$ cycles
$$i\overset{\vf_B}\to{\to}\nu(i,1)\overset{\vf_B}\to{\to}\nu(i,2)\overset{\vf_B}\to{\to}...
\overset{\vf_B}\to{\to}\nu(i,k(i)-2)\overset{\vf_B}\to{\to}i$$ and its $k=k_B$ is defined
by the
formulas $k_B(\gamma)=2$ for $\gamma \in [1, ...,n]$, $k_B(\gamma)=1$ for $\gamma \in
[n+1, ...,r-n]$.
\medskip
{\bf Theorem 11.5.} Let $M$ be from 11.4 (i.e. a standard-1 t-motive having
$\lambda_\goth m \ge 3$). Then $M'=M(B)$.
\medskip
{\bf Proof.}\footnotemark \footnotetext{This proof is a generalization
of the corresponding proof of
Taguchi; we keep his notations.} Let $e'_*=(e'_1, ... e'_{r-n})^t$ be the vector column
of elements of a basis
of $M(B)$ over $\p[\tau]$ satisfying
$$T e'_*=Be'_*\eqno{(11.5.1)}$$
Let us consider the set of pairs $(j,\goth k)$ such that either $j=1,...,n$,
$\goth k=0,1$ or
$j=n+1,...,r-n$, $\goth k=0$. For each pair $(j,\goth k)$ of this set we let (as
in [T], p. 580)
$Y_{j \goth k} = \tau^{\goth k}e'_j$. Formulas (11.4.2) show that these $Y_{**} $
form a basis of
$M(B)_{\p[T]}$, and the action of $\tau$ on this basis is given by the
following formulas
(here $j=1, ... , k(i)-2$):
$$\tau(Y_{i,0})=Y_{i,1} \eqno{(11.5.2.1)}$$
$$\tau(Y_{i,1})= (T-\theta) Y_{\nu(i,1),0} + \sum_{\gamma=1}^n
a_{1\gamma i}Y_{\gamma,1} \eqno{(11.5.2.2)}$$
$$\tau(Y_{\nu(i,j),0})= (T-\theta) Y_{\nu(i,j+1),0} + \sum_{\gamma=1}^n
a_{j+1,\gamma,i}Y_{\gamma,1} \hbox{ if }j<k(i)-2 \eqno{(11.5.2.3)} $$
$$\tau(Y_{\nu(i,k(i)-2),0})= (T-\theta) Y_{i,0} + \sum_{\gamma=1}^n
a_{k(i)-1,\gamma,i}Y_{\gamma,1} \eqno{(11.5.2.4)}$$
Let $X'_{**}$ be the dual basis to the basis $X_{**}$ of 11.2.2.
\medskip
{\bf 11.5.3.} Let us consider the following correspondence between
$X'_{**}$ and $Y_{**}$:
\medskip
$X'_{ij}$ corresponds to $Y_{\nu(i,j),0}$ for
the pair $(i,j)$ like in (11.4.1),
\medskip
$X'_{i0}$ corresponds to $Y_{i1}$ for $1 \le i \le n$;
\medskip
$X'_{i, k(i)-1}$ corresponds to $Y_{i0}$ for $1 \le i \le n$. 
\medskip
Therefore, in order to prove the Theorem 11.5 we must check that
matrices defined by the dual to (11.3.4.*) and by (11.5.2.*) satisfy (1.10.1) under
identification (11.5.3). This is an elementary exercise. $\square $
\medskip
{\bf Remark 11.6.} Clearly it is possible to generalize the Theorem 11.5 to a larger class
of t-motives --- some subclass of standard-3 t-motives, see Definition 11.8.1.
The below example of the proof of Proposition 11.8.7 shows that probably the condition
$\lambda_\goth m \ge 3$ of
the Theorem 11.5 can be changed by
$\lambda_\goth m \ge 2$: it is necessary
to modify slightly formulas 11.4.2. From another side, a standard-1
t-motive of the Example 2.5 shows that this condition cannot
be changed to $\lambda_\goth m\ge1$.
\medskip
{\bf 11.7. An elementary transformation.} To formulate the proposition
11.7.3, we change slightly notations in 1.9.1, namely, instead of $\goth A =
\sum_{i=0}^l \goth A_i \tau^i$ we consider polynomials $P_k(M)$ of $x_1, ...
,x_n$ ($k=1, ...,n$) defined by the formula
$$P_k(M)= \sum_{i=0}^l\sum_{j=1}^n a_{ikj} x^{q^i}_j \eqno{(11.7.1)} $$
Particularly, if $E$ is the t-module associated to $M$ (see [G], 5.4.5),
$x_*=(x_1, ..., x_n)^t$ an element of $E$ then 11.7.1 is equivalent to
$Tx_*=P_*(x_*)$ where $P_*=(P_1(M),...,P_n(M))^t$ is the vector column.
For a standard-1 t-motive $M$ (we use notations of 11.3.0) having $\goth m\ge 2$ we denote vector columns $\goth
P_1(M)=(P_1(M),...,P_{\gamma_1}(M))^t$, $\goth
P_2(M)=(P_{\gamma_1+1}(M),...,P_{\gamma_2}(M))^t$. We use similar
notations for $M'$.
\medskip
{\bf 11.7.2.} Let $M$ be as above, we consider the case
$\lambda_2=\lambda_1-1$. Let $C$ be a fixed $\gamma_1\times (\gamma_
2-\gamma_1)$-matrix. We define a transformed t-motive $M_1$ by the
formulas

$$\goth P_1(M_1)= \goth P_1(M)+C\goth P_2(M)^q$$

$$P_i(M_1)=P_i(M) \hbox{ for } i>\gamma_1$$
\medskip
{\bf Proposition 11.7.3.} For $M$, $C$, $M_1$ of 11.7.2 the dual $M'_1$
of $M_1$ is described by the following formulas:
$$\goth P_2(M'_1)= \goth P_2(M')-C^t\goth P_1(M')^q$$
$$P_i(M'_1)=P_i(M')\hbox{ for }i\not\in[\gamma_1+1, ... , \gamma_2]$$

{\bf Proof} is similar to the proof of the Theorem 11.5, it is omitted.
$\square$
\medskip
\medskip
{\bf 11.8. Non-pure t-motives.} Most results of this subsection are conditional. We shall show that under some natural conjecture the condition of purity in 8.2 and 8.4 is essential, and that for non-pure t-motives the notion of algebraic duality is richer than the notion of analytic duality.

We generalize slightly the definition 11.2.1 as follows. Let
$\succ$ be a linear ordering on the set $[1,...,n]$, and let $\vf$, $k$ be as in 11.2.
\medskip
{\bf Definition 11.8.1.} A standard-3 t-motive of the type $(\vf, k,\succ)$ is
a t-motive of
dimension
$n$ given by the formulas
$$Te_{\vf(i)}= \theta e_{\vf(i)} +\sum_{j=1}^n\sum_{l=1}^{k(j)-1}a_{l,\vf(i),j}\
\tau^le_j +\sum_{j\succ i} a_{k(j),\vf(i),j}\ \tau^{k(j)}e_j +
\tau^{k(i)}e_i\eqno{(11.8.2)}$$
where $a_{***}\in\p$ are coefficients (the only difference with 11.2.1 is the term
$\sum_{j\succ i} a_{k(j),\vf(i),j}\ \tau^{k(j)}e_j$). We denote it by $M(a_{***})$.

Let $M_1=M(a_{1***})$, $M_2=M(a_{2***})$ be two isomorphic standard-3 t-motives of the same type $(\vf,
k,\succ)$ with
$\p[\tau]$-bases $e_{1*}$, $e_{2*}$ respectively (we use notations of 11.8.2 for both
$M_1$, $M_2$). There exists $C\in M_n(\p[\tau])$ such that the formula defining an isomorphism
between $M_1$ and $M_2$ is the following: $e_{2*}=Ce_{1*}$.
\medskip
{\bf Conjecture 11.8.3.} For a generic set of $a_{1***}$ there exists only a countable set of $a_{2***}$ such that $M_2$ is isomorphic to $M_1$.
\medskip
This conjecture is based on calculations in some explicit cases. Particularly, it is proved if $M_1$, $M_2$ are given by the below formula 11.8.5.1 and entries of $C$ are polynomials in $\tau$ of degree $\le 1$.

We denote by $\Cal M_{u}(r,n)$ the moduli space of uniformizable t-motives of the rank
$r$ and dimension $n$,
by $\Cal L(r,n)$ the moduli space of lattices of the rank $r$ and dimension $n$ and by
$\goth L:\Cal M_{u}(r,n) \to \Cal L(r,n)$ the functor of lattice associated to an uniformizable
t-motive.
\medskip
{\bf Proposition 11.8.5.} Conjecture 11.8.3 implies that the dimension of the fibers of
$\goth L$ is $> 0$ for $r=3$, $n=2$. Particularly, we cannot omit condition of purity in the
statement of 8.2.
\medskip
{\bf Proof.} We consider standard-3 t-motives of the type $n=2$,
$\vf=Id$, $k(1)=2$,
$k(2)=1$, $2\succ 1$. Such $M_1=M_1(a_{111}, a_{112}, a_{121})$ is given by
$$\goth A_0=\theta E_2, \ \ \goth A_{1}= \left(\matrix a_{111} & a_{112}
\\ a_{121} & 1 \endmatrix \right), \ \ \goth A_2= \left(\matrix 1 & 0
\\0 & 0 \endmatrix \right)\eqno{(11.8.5.1)}$$ (notations of Example 10.5).
It has $r=3$, it is not pure, hence it has no dual.
Conjecture 11.8.3 implies
that the dimension of the moduli space of these t-motives is 3 (because there are 3
coefficients
$a_{111}, a_{112}, a_{121}$). Uniformizable t-motives form an open subset of this moduli
space, while
the moduli space of lattices of $n=2$ and $r=3$ has dimension 2. $\square$
\medskip
{\bf Remark.} Similar calculations are valid for any sufficiently large $r$, $n$.
\medskip
Standard-3 t-motives of the above type have not dual. The following proposition
shows that the same
phenomenon holds for t-motives having dual. We denote by $\Cal M_{u,d}(r,n)$ the
moduli space of uniformizable t-motives of the rank $r$ and dimension $n$ having dual, by
$\Cal L_d(r,n)$ the moduli space of lattices of the rank $r$ and dimension $n$ having
dual, by $\goth L_d:\Cal M_{u,d}(r,n) \to \Cal L_d(r,n)$ the functor of lattice and by $D_M: \Cal
M_{u,d}(r,n)\to \Cal M_{u,d}(r,r-n)$, $D_L: \Cal L_d(r,n)\to \Cal L_d(r,r-n)$ the
functors of duality on t-motives and lattices respectively. Practically, Theorem 5
means that the following diagram is commutative:
$$\matrix \Cal M_{u,d}(r,n)&\overset{D_M}\to{\to}&\Cal M_{u,d}(r,r-n)
\\ \\ \goth L_d\downarrow&&\goth L_d\downarrow\\  \\ \Cal L_d(r,n)&\overset{D_L}\to{\to}&\Cal
L_d(r,r-n)\endmatrix\eqno{(11.8.6)}$$

{\bf Proposition 11.8.7.} Conjecture 11.8.3 implies that the dimension of the fibers of
$\goth L_d$ in the diagram (11.8.6) is $> 0$ for $r=5$, $n=2$.
\medskip
Practically, this means that the notion of algebraic duality is "richer" than the notion
of analytic duality.
\medskip
{\bf Proof.} We consider standard-3 t-motives of the type $n=2$,
$\vf=Id$, $k(1)=3$,
$k(2)=2$, $2\succ 1$, $r=5$. Such $M$ is given by
$$\goth A_0=\theta E_2, \ \ \goth A_{1}= \left(\matrix a_{111} & a_{112}
\\ a_{121} & a_{122}  \endmatrix \right), \ \ \goth A_{2}= \left(\matrix a_{211} & a_{212}
\\ a_{221} & 1 \endmatrix \right), \ \ \goth A_3= \left(\matrix 1 & 0
\\0 & 0 \endmatrix \right)$$ (notations of Example 10.5). It has dual. Really, we denote by
$A_{i*j}$ the $j$-th column of $\goth A_i$, and we denote by $ (C_1 | C_2 )$ the matrix formed
by union of columns
$C_1$, $C_2$. Then $M'=M(B)$ is also a standard-3 t-motive, where
\medskip
$B_1= \left(\matrix - \det \goth A_2 & -a_{221} & 1
\\ - \det ( A_{1*2} | A_{2*2} ) & -a_{122} & 0
\\ - \det ( A_{1*1} | A_{2*2} ) & -a_{121} & 0
\endmatrix \right)$,
$B_2 = \left(\matrix 0 & 0 & 0
\\ -a_{212}^q & 1 & 0
\\ 1 & 0 & 0 \endmatrix \right)$

The same arguments as in the proof of Proposition 11.8.5 show that the conjecture 11.8.3
implies that the dimension of the moduli space of these t-motives is 7, while
the moduli space of lattices of $n=2$ and $r=5$ has dimension 6. $\square$
\medskip
As above, similar calculations are valid for any sufficiently large $r$, $n$; clearly the dimension of fibers of $\goth L_d$ becomes larger as $r$, $n$ grow.
\medskip
Let us mention two open questions related to the functor $\goth L$. Firstly, let $L$ be a self-dual lattice such that $L\in \goth L(\Cal M_{u,d}(2n,n))$. This means that $D_M: \goth L_d^{-1}(L) \to \goth L_d^{-1}(L)$ is defined.
\medskip
{\bf Open question 11.8.8.} What can we tell on this functor, for example, what is the dimension of its stable elements?
\medskip
Secondly, let us consider $M_1$, $M_2$ of CM-type with CM-field $\n F_{q^r}(T)$, see 13.3.
\medskip
{\bf Open question 11.8.9.} Let the CM-types $\Phi_1$, $\Phi_2$ of the above $M_1$, $M_2$ satisfy $\Phi_1\ne \alpha \Phi_2$, where $\alpha\in \Gal(\n F_{q^r}(T)/\n F_{q}(T))$. Are lattices $L(M_1)$, $L(M_2)$ non-isomorphic?
\medskip
Clearly the negative answer to this question implies the negative answer to the Question 8.2.
\medskip
For any given $M_1$, $M_2$ the answer can be easily found by computer calculation.
Really, let $M$ be one of $M_1$, $M_2$, $c_1,...,c_r$ a basis of $\n F_{q^r}/\n F_{q}$ and
$\alpha_{1},...,\alpha_{n} \subset \Gal(\n F_{q^r}(\theta)/\n F_{q}(\theta))$ the CM-type of $M$.
We define matrices $\Cal M$, $\Cal N$ as follows: $(\Cal M)_{ij}=\alpha_j(c_i)$ $(i,j=1,...,n$),
$(\Cal N)_{ij}=\alpha_{j}(c_{n+i})$, $j=1,...,n$, $i=1,...,r-n$. The Siegel matrix $Z(M)$ is obviously $\Cal N\Cal
M^{-1}$. So, we can find explicitly $Z(M_1)$, $Z(M_2)$ for both $M_1$, $M_2$. To check
whether $Z(M_1)$, $Z(M_2)$ are equivalent or not, it is sufficient to find a solution to
3.8.1 such that the entries of $A$, $B$, $C$, $D$ are in $M_{*,*}(\n F_q)$ (this is
obvious: the condition $\exists \gamma \in GL_r(\n Z_\infty)$ is equivalent to the
condition $\exists \gamma \in GL_r(\n F_q)$, because entries of $Z(M_1)$, $Z(M_2)$ are in $\n F_{q^r}$). The equation 3.8.1 is linear with respect
to $A$, $B$, $C$, $D$, and we can check whether its solution satisfying $\det
\gamma\ne 0$ exists or not.
\medskip
For the case $q=2$, $r=4$, $n=2$, CM-types of $M_1$, $M_2$ are $(Id, Fr)$, $(Id, Fr^2)$
respectively, a calculation shows that the answer is positive: lattices $L(M_1)$, $L(M_2)$ are not isomorphic.
\medskip
{\bf 12. t-motives having multiplications.}
\medskip
Let $\goth K$ be a separable extension of $\n F_q(T)$ such that $\goth K_C:=\goth K\underset{\n F_q}\to{\otimes} \p$ is also a field, $\pi: X\to P^1(\p)$ the projection of curves over $\p$ corresponding to $\p(T)\subset \goth K_C$. Let $\goth K$, $X$ satisfy the condition: $\infty\in X$ is the only point
on $X$ over $\infty\in P^1(\p)$. Let $\w_{\goth K}$ be the subring of
$\goth K$ consisting of elements regular outside of infinity. We
denote $g=\dim \goth K/\n F_q(T)$ and
$ \alpha_1, ... , \alpha_g: \goth K\to\p$ --- inclusions over $\iota: \n F_q(T)\to\p$
(recall that $\iota(T)=\theta$). Let $\Cal W$ be a central simple algebra
over $\goth K$ of dimension $\goth q^2$. Each $ \alpha_i: \goth K \to
\p$ can be extended to a representation $\chi_i: \Cal W \to M_\goth
q(\p)$.
\medskip
{\bf 12.1. Analytic CM-type.} Let $(L, V)$ be as in Section 2 (recall that
$\w=\n F_q[T]$) such that there exists an inclusion $i: \Cal W \to
\End^0(L, V)$, where $\End^0(L, V)=\End(L, V)\underset{\w}\to{\otimes} \n F_q(T)$.
It defines a representation of $\Cal W$
on $V$ denoted by $\Psi$ which is isomorphic to $\sum_{i=1}^g\goth r_i\chi_i$ where
$\{\goth r_i\}$ are some multiplicities (the CM-type of the action of
$\Cal W$ on $(L, V)$). [Proof: restriction of $\Psi$ on $\goth K$ is a sum of one-dimensional representations, i.e. $V=\oplus_{i=1}^g V_i$ where $k\in \goth K$ acts on $V_i$ by multiplication by $\alpha_i(k)$. Spaces $V_i$ are $\Psi$-invariant. We consider an isomorphism $\Cal W\otimes_\goth K\p=M_\goth q(\p)$ where the inclusion of $\goth K$ in $\p$ is $\alpha_i$. We extend $\Psi|_{V_i}$ to $\Cal W\otimes_\goth K\p$ by $\p$-linearity using the inclusion $\alpha_i$ of $\goth K$ in $\p$. It remains to show that a representation of $M_\goth q(\p)$ is a direct sum of its $\goth q$-dimensional standard representations. We consider the corresponding representation of Lie algebra $\goth s\goth l_\goth q(\p)$. It is a sum of irreducible representations. Let $\omega$ be the highest weight of any of these irreducible representations. $\omega$ is extended uniquely to the set of diagonal matrices of $M_\goth q(\p)$, because $\omega$ is identical on scalars. Since our representation is not only of Lie algebra but of algebra $M_\goth q(\p)$, we get that $\omega$ is a ring homomorphism $\Diag(M_\goth q(\p))\to \p$. There exists the only such $\omega$ corresponding to the $\goth q$-dimensional standard representation].

Further, we
denote $m=\dim_\Cal W L \otimes \n F_q(T)$ ($g$, $\goth q$, $\Psi$, $\goth
r_i$, $m$ are analogs of $g$, $q$, $\Phi$, $r_i$, $m$ of [Sh63] respectively).
Clearly we have
$$n=\goth q\sum_{i=1}^g\goth r_i, \ \ r=mg\goth q^2\eqno{(12.2)}$$
By functoriality, we have the dual inclusion $i': \Cal W^{op} \to
\End^0(L',V')$ where $\Cal W^{op}$ is the opposite algebra.
\medskip
{\bf Remark.} A construction of Hilbert-Blumental modules ([A], 4.3, p. 498) practically is a particular case
of the present
construction: for Hilbert-Blumental modules we have $\goth q=1$, i.e. $\goth K=\Cal W$, and all $\goth r_i=1$. Anderson considers the case when $\infty$ splits completely; this difference with the present case is not essential.
\medskip
{\bf Proposition 12.3.} If the dual pair $(L',V')$ exists then the
CM-type of the dual inclusion is $\{m\goth
q - \goth r_i\}$, $i=1, ... ,g$.
\medskip
{\bf Proof.} We have $L\underset{\n Z_\infty}\to{\otimes}\p$ is
isomorphic to $(\Cal W\underset{\n F_q(\theta)}\to{\otimes}\p)^m$ as a
$\Cal W$-module. Since the natural representation of $\Cal W$ on $\Cal
W\underset{\n F_q(\theta)}\to{\otimes}\p$ is isomorphic to $\goth
q\sum_{i=1}^g\chi_i$ we get that $L\underset{\n Z_\infty}\to{\otimes}\p$ is isomorphic to $m\goth q\sum_{i=1}^g\chi_i$ as a $\Cal W$-module. Consideration of the exact
sequence $0\to {V'}^* \to L\underset{\n Z_\infty}\to{\otimes}\p \to V\to0$ gives us the desired. $\square$
\medskip
{\bf Remark 12.4.1.} This result is an analog of the corresponding
theorem in the number field case. We use notations of [Sh63], Section 2. Let
$A$ be an abelian variety having endomorphism algebra of type IV, and
$(r_\nu, s_\nu)=(r_\nu(A), s_\nu(A))$ are from [Sh63], Section 2, (8).
Then
$$r_\nu(A')=mq-r_\nu(A)=s_\nu(A), \ s_\nu(A')=mq-s_\nu(A)=r_\nu(A)$$
By the way, Shimura writes that the CM-types of $A$ and $A'$ coincide
([Sh98], 6.3, second line below (5), case $A$ of CM-type). We see that
his affirmation is not natural: he considers the complex conjugate
action of the endomorphism ring on $A'$. It is necessary to take into
consideration this difference of notations comparing formulas of 12.3 and 13.2
with the corresponding formulas of Shimura.
\medskip
{\bf Remark 12.4.2.} According [L1], a t-motive $M$ is an analog of an abelian variety $A$ with multiplication by an imaginary quadratic field $K$. The above consideration shows that this analogy holds for $M$ and $A$ having more multiplications. Really, if $A$ has more multiplications then (we use notations of [Sh63], Section 2) $F_0=FK$, and numbers $(r_\nu(A), s_\nu(A))$ satisfy $n(A)=q\sum_{i=1}^g r_\nu(A)$, where $(n(A), \dim(A)-n(A))$ is the signature of $A$ treated as an abelian variety with multiplication by $K$. This is an analog of 12.2.
\medskip
{\bf 12.5. Complete multiplication.} Here we consider the case $\goth
q=m=1$, i.e. $\goth K=\Cal W$ and $g=r$.
\medskip
{\bf Lemma 12.5.1.} In this case the condition $N=0$
implies that the CM-type $$\sum_{i=1}^r\goth r_i\alpha_i\eqno{(12.5.2)}$$ of the action of
$\goth K$ on
on $(L,V)$ has the property: all $\goth r_i$ are 0 or 1.
\medskip
{\bf Proof.} $N=0$ means that
the action of $T\in\w$ on $V$ is simply multiplication by $\theta$. We write the
CM-type $\sum_{i=1}^r\goth r_i\chi_i$ in the form $\sum_{i=1}^n\chi_{\alpha_i}$ where
$\alpha_1,...,\alpha_n\in[1,...,r]$ are not necessarily distinct. Let $l_1$ be an (only)
element of a basis of $L\otimes_{\w_{\goth K}}\goth K$ over $\goth K$ and $e_1,...,e_n$ a
basis of $V$
over $\p$ such that the action of $\goth K$ on $V$ is given by the formulas
$$k(e_i)=\chi_{\alpha_i}(k)e_i, \ \ \ k\in \goth K$$
Multiplying $e_i$ by scalars if necessary, we can assume that $l_1=\sum e_i$. Therefore,
if $\alpha_i=\alpha_j$ (i.e. not all $\goth r_*$ in (12.5.2) are 0, \ 1) then the
$e_{\alpha_i}$-th coordinate of any element of $L$ coincide with its $e_{\alpha_j}$-th
coordinate, hence $L$ does not $\p$-generate $V$ --- a contradiction. $\square$
\medskip
{\bf 12.5.3.} Let $M$ be a t-motive of
rank $r$ and dimension $n$ having multiplication by $\w_{\goth K}$. Recall that we
consider only
the case $N=0$. This means that the character of the action of $\goth K$ on $M/\tau M$
is isomorphic to $\sum_{i=1}^r\goth r_i\alpha_i$. Since $E(M)=(M/\tau M)^*$ we get that
the character of the action of $\goth K$ on $E(M)$ is the same. If
$$\hbox{all $\goth r_i$ are 0 or 1}\eqno{(12.5.4)}$$
we shall use the terminology that $M$ has the CM-type $\Phi\subset \{\alpha_1, ... ,
\alpha_r \}$
where $\Phi$ is defined by the condition $\alpha_i\in \Phi \iff \goth r_i=1$.

It is easy to see that this case occurs for uniformizable $M$. Really, if $M$ is
uniformizable then
the action of $\goth K$ can be prolonged on $(L(M),V(M))$, and the character
of the action of $\goth K$ on $V(M)$ coincides with the one on $E(M)$. The result follows
from Lemma 12.5.1.
\medskip
{\bf Lemma 12.5.5.} There exists a canonical isomorphism $\gamma$ from the set of
inclusions $\alpha_1, ... , \alpha_r$ to the set of points $ \theta_{\alpha_1}, ... ,
\theta_{\alpha_r}$ of $X$ over $\theta\in P^1(\p)$.
\medskip
{\bf Proof.} A point $t\in X$ over $\theta\in P^1(\p)$ defines a function
$\vf_t: \goth K_{C} \to P^1(\p)$ --- the value of
an element $f\in \goth K_{C}$ treated as a function on $X$ at the point $t$. This
function
must satisfy the standard axioms of valuation and the condition $\vf_t(T)=\theta$. Let
$\alpha_i$
be an inclusion of $\goth K$ to $\p$ over $\iota$. It defines a valuation
$\vf_{\alpha_i}: \goth K_{C} \to P^1(\p)$ by the formula $\vf_{\alpha_i}(k\otimes
f)=\alpha_i(k) f(\theta)$,
where $k\in\goth K$, $f\in \p(T)$. We define $\gamma(\alpha_i)$ by the condition
$\vf_{\gamma(\alpha_i)}=\vf_{\alpha_i}$; it is easy to see that $\gamma$ is an
isomorphism. $\square$
\medskip
{\bf Theorem 12.6.} For any above \{$\goth K$, $\Phi$\} there exists an
t-motive $(M,\tau)$ with
complete multiplication by $\goth K$ having CM-type $\Phi$.
\medskip
{\bf Proof (Drinfeld).} We denote the divisor $\sum_{\alpha_i\in\Phi}\gamma({\alpha_i})$
by
$\theta_\Phi$. We construct a $\Cal F$-sheaf $F$ of dimension 1 over
$\goth K$ which will give us $M$.
Let $\fr$ be the Frobenius map on $ \Pic_0(X)$. It is an algebraic
map, and the $\fr-\Id: \Pic_0(X) \to \Pic_0(X)$ is an algebraic map as
well. Since the action of $\fr$ on the tangent space of $ \Pic_0(X)$
at 0 is the zero map, the action of $\fr-\Id$ on the tangent space of
$ \Pic_0(X)$ at 0 is the minus identical map and hence $\fr-\Id$ is an
isogeny of $ \Pic_0(X)$. Particularly, there exists a divisor $D$ of
degree 0 on $X$ such that we have the following equality in $\Pic_0(X)$:
$$\fr(D)-D=-\theta_\Phi+n\infty \eqno{(12.6.0)}$$
This means that if we let $F=F_\Phi=O(D)$ then there exists a rational map
$\tau_X=\tau_{X,\Phi}: F^{(1)}\to F$ such
that
$$\Div(\tau_X)=\theta_\Phi-n\infty\eqno{(12.6.1)}$$
The pair $(F_\Phi, \tau_{X,\Phi})$ is the desired $\Cal F$-sheaf.
\medskip
{\bf Remark.} It is easy to see that if the genus of $X$ is $> 0$ then different CM-types $\Phi_1$, $\Phi_2$ give us different sheaves $F_{\Phi_1}$, $F_{\Phi_2}$, while if the genus of $X$ is 0 then $F_{\Phi_1}=F_{\Phi_2}=\Cal O$, but the maps $\tau_{X,\Phi_1}$, $\tau_{X,\Phi_2}$ are clearly different.
\medskip
Let $U_0=X-\{\infty\}$ be an open part of
$X$. We denote $F(U_0)$ by $\Cal M$, hence $F^{(1)}(U_0)=\Cal M^{(1)}$.
Since the support of the negative part of the right hand side of
12.6.1 is $\{\infty\}$,
we get that the (a priory rational) map $\tau_X(U_0): \Cal M^{(1)}\to \Cal M$ is
really a map of
$\w_{\goth K}$-modules.

Let $M$ be a $\p[T]$-module obtained from $\Cal M$ by restriction of scalars from
$\w_{\goth K}$ to $\p[T]$. Construction $F\mapsto M$ is functorial, and we denote this functor by $\delta$. Further, we denote by $\alpha$ the tautological isomorphism
$\Cal M\to M$. $M$ is a free $r$-dimensional $\p[T]$-module, and (because
$M^{(1)}$ is isomorphic to $M$) the
same restriction of scalars of $\tau_X(U_0)$ defines us a $\p[T]$-skew map from $M$ to
$M$ denoted by $\tau$ (skew means that $\tau(zm)=z^q\tau(m)$, $z\in\p$). $\tau$ is
defined by the formula $\tau(m)=\alpha\circ \tau_X((\alpha^{-1}(m))^{(1)})$.

It is easy to check that
$(M,\tau)$ is the required t-motive. Really, $M$ is a $\w_{\goth K}$-module, and
$\tau$
commutes with this multiplication. The fact that the positive part
of the right hand
side of 12.6.1 is $\theta_\Phi$ means that 1.13.2 holds for $M$ and that the CM-type of
the action of
$\w_{\goth K}$ is $\Phi$.
\medskip
{\bf Remark 12.6.2.} It is easy to prove for this case that $M$ is a free
$\p[\tau]$-module.
Really, it is
sufficient to prove (see [G], Lemma 12.4.10)
that $M$ is finitely generated as a $\p[\tau]$-module. We choose $D$ such that
$\infty\not\in \Supp (D)$. There exists $P\in \goth K_{C}^*$
such that $\tau_X(U_0): \Cal M^{(1)}\to \Cal M$ is multiplication by $P$ (recall that both
$\Cal M^{(1)}$, $\Cal M$ are $\w_{\goth K}$-submodules of $\goth K$). 12.6.0 implies that
$-\ord_\infty(P)=n$.
There
exists a number $n_1$ such that $$\hbox{(a) $h^0(X,\Cal O(D+n_1\infty))>0$; \ \ (b)
for any $k\ge 0$ we have}$$ $$h^0(X,\Cal O(D+(n_1+k)\infty))=h^0(X,\Cal
O(D+n_1\infty))+k\eqno{(12.6.3)}$$ $$h^0(X,\Cal O(D^{(1)}+(n_1+k)\infty))=h^0(X,\Cal
O(D^{(1)}+n_1\infty))+k\eqno{(12.6.4)}$$
It is sufficient to prove that if $g_1,...,g_k$ are elements of a basis of $H^0(X,\Cal
O(D+(n_1+n)\infty))$, then
for any $Q\in \Cal M$ the element $\alpha(Q)\in M$ is generated by
$\alpha(g_1),...,\alpha(g_k)$ over
$\p[\tau]$. We prove it by induction by $n_2:=-\ord_\infty(Q)$. If $n_2\le n_1+n$ the
result is trivial. If not
then 12.6.3, 12.6.4 imply that the multiplication by $P$ defines an isomorphism
$$H^0(X, \Cal O(D^{(1)}+(n_2-n)\infty))/H^0(X, \Cal O(D^{(1)}+(n_2-n-1)\infty))\to$$
$$\to H^0(X,\Cal O(D+n_2\infty))/H^0(X, \Cal O(D+(n_2-1)\infty))$$
This means that $\exists Q_1\in H^0(X,\Cal O(D^{(1)}+(n_2-n)\infty))$,
$-\ord_\infty(Q_1)=n_2-n$ such
that
$-\ord_\infty(Q-PQ_1)\le n_2-1$. An element $Q_1^{(-1)}\in \Cal M$ exists;
since $\alpha(Q)=\tau(\alpha(Q_1^{(-1)}))+\alpha(Q-PQ_1)$, the result follows by
induction.
$\square$
\medskip
If $\goth K$ and $\Phi$ are given then the construction of the Theorem 12.6 defines $F$ uniquely up to tensoring by $O(D)$ where $D\in \Div (X(\goth K))$. We denote the set of these $F$ by $F$(\{$\goth K$, $\Phi$\}), and we denote by $M$(\{$\goth K$, $\Phi$\}) the set $\delta(F$(\{$\goth K$, $\Phi$\})). Further, we denote by $\Phi'=\{\alpha_1, ... ,
\alpha_r\}-\Phi$ the complementary CM-type.
\medskip
{\bf Theorem 12.7.} Let $M\in M(\{\goth K,\Phi\})$. Then $M'$ exists, and $M'\in M(\{\goth K,\Phi'\})$. More exactly, if $F\in F$(\{$\goth K$, $\Phi$\}) then $F^{-1}\otimes \Cal D^{-1}\in F$(\{$\goth K$, $\Phi'$\}) where $\Cal D$ is the different
sheaf on $X$, and if $M=\delta(F)$ then $M'=\delta(F^{-1}\otimes \Cal D^{-1})$.
\medskip
{\bf Proof.} Let $G$ be any invertible sheaf on $X$. We have a
\medskip
{\bf Lemma 12.7.0.} There exists the canonical isomorphism $\varphi_G:
\pi_*(G^{-1}\otimes \Cal D^{-1}) \to \Hom_{P^1}(\pi_*(G), \Cal O)$.
\medskip
{\bf Proof.} At the level of affine open sets $\varphi_G$ comes from the trace bilinear form of field extension $\goth K/\n F_q(T)$. Concordance with glueing is obvious. $\square$
\medskip
We need the relative version of this lemma. Let $G_1$, $G_2$ be invertible sheaves on $X$,
$\rho: G_1 \to G_2$ any rational map. Obviously there exists a rational map
$\rho^{-1}: G_1^{-1} \to G_2^{-1}$. Recall that we denote by
$\rho^{inv}: G_2 \to G_1$ the rational map which is inverse to $\rho$
respectively the composition. The map $\pi_*(\rho^{-1}\otimes \Cal D^{-1}): \pi_*(G_1^{-1}\otimes \Cal D^{-1})\to \pi_*(G_2^{-1}\otimes \Cal D^{-1})$ is obviously defined. The map (denoted by $\beta(\rho)$) from $\Hom_{P^1}(\pi_*(G_1), \Cal O)$ to $\Hom_{P^1}(\pi_*(G_2),
\Cal O)$ is defined as follows at the level of affine open sets: let
$\gamma\in\Hom_{P^1}(\pi_*(G_1), \Cal O)(U)$ where $U$ is a sufficiently
small affine subset of $P^1$, such that we have a map $\gamma(U): \pi_*(G_1)(U)
\rightarrow \Cal O(U)$. Then $(\beta(\gamma))(U)$ is the composition
map $\gamma(U)\circ \pi_*(\rho^{inv})(U)$:
$$\pi_*(G_2)(U)\overset{\pi_*(\rho^{inv})(U)}\to{\longrightarrow}
\pi_*(G_1)(U)\overset{\gamma(U)}\to{\to}\Cal O(U)$$

{\bf Lemma 12.7.1.} The above maps form a commutative diagram:
$$\matrix
\pi_*(G_1^{-1}\otimes \Cal D^{-1})&\overset{\pi_*(\rho^{-1}\otimes \Cal
D^{-1})}\to{\longrightarrow}&\pi_*(G_2^{-1}\otimes \Cal D^{-1})&&\\&&&&\\
\varphi_{G_1}\downarrow&&\varphi_{G_2}\downarrow &&\\&&&&\\
\Hom_{P^1}(\pi_*(G_1), \Cal
O)&\overset{\beta(\rho)}\to{\longrightarrow}&\Hom_{P^1}(\pi_*(G_2),
\Cal O)&&\square
\endmatrix$$

We apply this lemma to the case $\{\rho: G_1 \to G_2\}=\{\tau_{X, \Phi}:
F^{(1)}\to F\}$. We have:
$$\Div(\tau_{X, \Phi}^{-1}\otimes \Cal D^{-1})=-\Div(\tau_{X, \Phi})=-\theta_{\Phi}+n\infty$$
Futher, we multiply $\tau_{X, \Phi}^{-1}\otimes \Cal D^{-1}$ by $T-\theta$. We have:
$$\Div((T-\theta)\tau_{X, \Phi}^{-1}\otimes \Cal
D^{-1})=\Div(T-\theta)+\Div(\tau_{X, \Phi}^{-1}\otimes \Cal D^{-1})=\theta_{\Phi'}-(r-n)\infty$$
i.e. $(T-\theta)\tau_{X, \Phi}^{-1}\otimes \Cal D$ is one of
$\tau_{X,\Phi'}$, i.e. $F^{-1}\otimes \Cal D^{-1}\in F$(\{$\goth K$, $\Phi'$\}). Further, $(T-\theta)\beta(\tau_{X, \Phi})$ is the map which is used in the definition of duality of
$M$. This means that the lemma 12.7.1 implies the theorem. $\square$
\medskip
{\bf Remark 12.8.} There exists a simple proof of the second part of the Theorem 5 for uniformizable abelian
t-motives $M$ with complete multiplication by $\w_\goth K\subset \goth K$. Recall that this second part is the proof of 2.7 for $M$. Really, let us
consider the diagram 2.5. The CM-types
of action of $\goth K$ on $\Lie(M)$ and on $E(M)$ coincide,
and the CM-types of action of $\goth K$ on a vector space and on its
dual space coincide. This means that the CM-type of $V^*$ is $\Phi$
and the CM-type of $V'$ is $\Phi'$. Further, $\gamma_D$ of 2.5 commutes with complete multiplication: this follows immediately for example from a description of $\gamma_D$ given in Remark 5.2.8. Really, all homomorphisms of 5.2.9 commute with complete multiplication. For example, this condition for $\delta$ of 1.11.1 is written as follows: if $k\in \goth K$, $\goth m_k(M)$, resp. $\goth m_k(M')$ is the map of complete multiplication by $k$ of $M$, resp. $M'$, then $(\goth m_k(M)\otimes Id) \circ \delta= (Id \otimes \goth m_k(M'))\circ \delta$ --- see any textbook on linear algebra.

Finally, since
$\Phi\cap \Phi'=\emptyset$ and the map $\vf'\circ \gamma_D \circ
\vf^*$ commutes with complete multiplication, we get that it must be
0.
\medskip
{\bf 13. Miscellaneous.}

\nopagebreak
\medskip
Let now $(L,V)$ be from 12.1, case $\goth q=m=1$, i.e. $\goth K=\Cal W$
and $r=g$, and let the ring of complete multiplication be the maximal
order $\w_\goth K$. We identify $\w$ and $\n Z_\infty$ via $\iota$, i.e. we consider $\goth K$ as an extension of $\n F_q(\theta)$. Let $\Phi$ be the CM-type of the action of $\goth
K$ on $V$.
This means that --- as an
$\w_\goth K$-module --- $L$ is isomorphic to $I$ where $I$ is an ideal of
$\w_\goth K$. The class of $I$ in $\Cl(\w_\goth K)$ is defined by $L$
and $\Phi$ uniquely; we denote it by $\Cl(L,\Phi)$.
\medskip
{\bf Remark.} $\Cl(L,\Phi)$ depends on $\Phi$, because the action of $\w_\goth K$ on $V$ depends on $\Phi$. Really, let $a\in L \subset V$, $a=(a_1,...,a_n)$ its coordinates, $\Phi=\{\alpha_{i_1},...,\alpha_{i_n}\}\subset \{\alpha_{1},...,\alpha_{r}\}$ and $k\in \w_\goth K$. Then $ka$ has coordinates $(\alpha_{i_1}(k)a_1,...,\alpha_{i_n}(k)a_n)$, i.e. depends de $\Phi$. Particularly, the $\w_\goth K$-module structure on $L$ depends on $\Phi$, and hence $\Cl(L,\Phi)$ depends on $\Phi$. For example, if $n=1$, $r=2$, $\Phi_1=\{\alpha_{1}\}$, $\Phi_2=\{\alpha_{2}\}$, then $\Cl(L,\Phi_2)$ is the conjugate of $\Cl(L,\Phi_1)$.
\medskip
{\bf Theorem 13.1.} $\Cl(L',\Phi')=(\Cl(\goth d))^{-1}
(\Cl(L,\Phi))^{-1}$ where
$\goth d$ is the different ideal of the ring extension $\w_\goth K/\w$.
\medskip
{\bf Proof.} This theorem follows from the above results;
nevertheless, I give here an explicit elementary proof. Let $a_*=(a_1, ..., a_r)^t$ be a
basis (considered
as a vector column) of $\goth K$ over $\n F_q(\theta)$ and $b_*=(b_1,
..., b_r)^t$ the dual basis. Recall that it satisfies 2 properties:
$$(1) \ \ \forall i\ne j \ \ \alpha_i(a_*)^t\alpha_j(b_*)=0 \ \ \
(\hbox{ i.e. } \sum_{k=1}^r \alpha_i(a_k) \alpha_j(b_k)=0
)\eqno{(13.1.1)} $$
(2) For $x\in \goth K$ let $\goth m_{x,a_*}$ (resp $\goth m_{x,b_*}$)
be the matrix of multiplication by $x$ in the basis $a_*$ (resp.
$b_*$). Then for all $x\in \goth K$ we have
$$\goth m_{x,a_*}=\goth m_{x,b_*}^t\eqno{(13.1.2)} $$

We define $\goth E_{n,r-n}$
as an $r\times r$ block matrix $\left(\matrix 0&E_{r-n}\\ -E_n&0
\endmatrix\right)$, and we define a new basis $\tilde b_*=(\tilde b_1,
..., \tilde b_r)^t$ by
$$\tilde b_*=\goth E_{n,r-n}b_*\eqno{(13.1.3)} $$
(explicit formula: $(\tilde b_1, ..., \tilde b_r)=(b_{n+1},
..., b_r, -b_1, ..., -b_n)$).

We can assume that $\Phi=\{\alpha_1, ... , \alpha_n \}$. Since $L$ has multiplication by $\w_\goth K$ and the CM-type of this
multiplication is $\Phi$, it is possible to choose $a_*$ such that $L\subset
\p^{n}$ is generated over $\n Z_\infty$ by $e_1, ... , e_r$ where
$$e_i=(\alpha_1(a_i), ... , \alpha_n(a_i))\eqno{(13.1.4)}$$ Let $\hat L\subset \p^{r-n}$ be
generated over $\n Z_\infty$ by $\hat e_1, ... , \hat e_r$ where
$$\hat e_i=(\alpha_{n+1}(\tilde b_i), ... , \alpha_r(\tilde b_i))\eqno{(13.1.5)}$$

{\bf Lemma 13.1.6. } $L'=\hat L$.
\medskip
{\bf Proof. } Let $A$ (resp. $B$) be a matrix whose lines are the lines of
coordinates of $e_1, ... , e_n$ (resp. $e_{n+1}, ... ,
e_r$) in 13.1.4, and $C$ (resp. $D$) a matrix whose lines are the lines of
coordinates of $\hat e_1, ... , \hat e_{r-n}$ (resp. $\hat e_{r-n+1}, ... ,
\hat e_r$) in 13.1.5. By definition of Siegel matrix, we have $L=\goth L(BA^{-1})$,
$\hat L=\goth L(DC^{-1})$ ($\goth L$ is defined in 3.1, 3.2). So, it is sufficient to prove that
$(BA^{-1})^t=DC^{-1}$, i.e. $A^tD=B^tC$. This follows immediately from
the definition of $A,B,C,D$ and (13.1.1). $\square$
\medskip
For $x\in \w_\goth K$ we denote by $\goth
M_x(L)$ the matrix of multiplication by $x$ in the basis $e_*$ (see
the notations of Remark 3.8). Obviously $\goth M_x(L)=\goth
m_{x,a_*}$.

Let now $\w_\goth K$ acts on $\p^{r-n}$ (the ambient space of $L'$) by CM-type
$\Phi'$. According (13.1.2) and
(13.1.3), the matrix of the action of $x\in \w_\goth K$ in the basis
$\tilde b_*$ is
$$\goth E_{n,r-n}\goth m_{x,a_*}^t\goth E_{n,r-n}^{-1}\eqno{(13.1.7)} $$
Let $\goth M$, $\goth M'$ be from Remark 3.8. Formula 3.8.4
shows that
$$\goth M'=\goth E_{n,r-n}\goth M^t\goth E_{n,r-n}^{-1}\eqno{(13.1.8)} $$
Formulas (13.1.7) and (13.1.8) --- because of Lemma 1.10.3 --- prove the theorem. $\square$
\medskip
{\bf 13.2. Compatibility with the weak form of the main theorem of
complete multiplication.}

\nopagebreak
\medskip
The reader can think that Theorem 13.1 is incompatible with the main theorem of complex
multiplication, because of the $-1$-th power in its statement. The reason is a bad choice
of notations of Shimura, he affirms that the CM-type of an abelian variety $A$ over a number field coincides with the CM-type
of $A'$, while we see that it is really the complement. Since an analog of even the weak form of the main theorem of complex multiplication --- Theorem 13.2.6 --- for the function field case is not proved yet, the main result of the present section --- Theorem 13.2.8 --- is conditional: it affirms that if this weak form of the main
theorem --- Conjecture 13.2.7 --- is true for a t-motive with complete multiplication $M$, then it is true for $M'$ as well. By the way, even if it will turn out that the statement of the Conjecture 13.2.7 is not correct, the proof of 13.2.8 will not be affected, because the main
ingredient of the proof is the formula 13.2.10 "neutralizing" the $-1$-th power of the
Theorem 13.1.

Let us recall some definitions of [Sh71], Section 5.5. We consider an abelian variety
$A=\n C^n/L$ with
complex multiplication by $K$. The set $\Hom(K, \bar \n Q)$ consists of $n$ pairs of
mutually
conjugate inclusions $\{\vf_1,\bar \vf_1, ..., \vf_n,\bar\vf_n\}$. $\Phi$ is a subset of
the set
$\Hom(K, \bar \n Q)$ such that $\forall i=1,...,n$ we have:
$$\hbox{$\Phi\cap\{\vf_i,\bar \vf_i\}$ consists of one element.}\eqno{(13.2.1)}$$
It is defined by the condition that the action of complex multiplication
of $K$ on $\n C^n$ is isomorphic to the direct sum of the elements of $\Phi$. Let $F$ be
the
Galois envelope of $K/\n Q$,
$$G:=\Gal (F/\n Q), \ \ H:=\Gal(F/K), \ \
S:=\bigcup_{\alpha\in\Phi}H\alpha\eqno{(13.2.2)}$$
(the elements of Galois group act on $x\in F$ from the right, i.e. by the formula
$x^{\alpha\beta}=(x^\alpha)^\beta$; for $\alpha\in\Phi$ we denote by $\alpha$ also a
representative in $G$ of the coset $\alpha$). We denote
$$H^{ref}:=\{\gamma\in G\vert S\gamma=S\}\eqno{(13.2.3)}$$
and let $K^{ref}$ be the subfield of $F$ corresponding to $H^{ref}$. We have:
$$H^{ref}S^{-1}=S^{-1}\eqno{(13.2.4)}$$
i.e. $S^{-1}$ is an union of cosets of $H^{ref}$ in $G$. We can identify these cosets
with elements of $\Hom(K^{ref}, \bar \n Q)$. $\Phi^{ref}\subset \Hom(K^{ref}, \bar \n Q)$
is, by definition, the set of these cosets. There is a map $\det
\Phi^{ref}:{K^{ref}}^\times \to K^\times$ defined as follows:
$$\det \Phi^{ref}(x):=\prod_{\alpha\in\Phi}\alpha(x)\eqno{(13.2.5)}$$
(it follows easily from the above formulas and definitions that $\det \Phi^{ref}(x)$
really belongs to $K^\times$). It can be extended to the group of ideles and factorized to
the group of classes of ideals, we denote this map
by $\det_{Cl} \Phi^{ref}: \Cl(K^{ref}) \to \Cl(K)$. Finally, let
$\theta^{ref}: \Gal(K^{ref\ Hilb}/K^{ref}) \to \Cl(K^{ref})$ be an isomorphism defined by
the Artin reciprocity law.

We consider the case $\End(A)=O_K$. In this case $L$ is isomorphic to an ideal of $O_K$,
its class
is well-defined by the class of isomorphism of $A$, we denote it by $\Cl(A)$.
\medskip
{\bf Theorem 13.2.6.} $A$ is defined over $K^{ref\ Hilb}$;

For any $\gamma\in\Gal(K^{ref\ Hilb}/K^{ref})$ we have

$\Cl(\gamma(A))=\det_{Cl} \Phi^{ref}\circ \theta^{ref} (\gamma)^{-1}(\Cl(A))$. $\square$
\medskip
This is a weak form of [SH71], Theorem 5.15 --- the main theorem of
complex multiplication.
\medskip
Now we define analogous objects for the function field case in order to formulate a
conjectural
analog of Theorem 13.2.6. Let $\goth K$, $\Phi$ be from 12.5.3. $\goth K^{ref}$,
$\Phi^{ref}$,
$\det \Phi^{ref}$ are defined by the same formulas 13.2.2 -- 13.2.5 like in the number
field case
($\n Q$ must be replaced by $\n F_q(T)$). The facts that 13.2.1 has no meaning in the
function field case and
that the order of $S$ is not necessarily the half of the order of $G$ do not affect the
definitions.

The $\infty$-Hilbert class field of $\goth K$ (denoted by $\goth K^{Hilb \ \infty}$) is an abelian
extension of
$\goth K$ corresponding to the subgroup
$$\goth K_\infty^*\cdot \prod_{v\ne\infty}O_{\goth K_v}^*\cdot \goth K^*$$
of the idele group of $\goth K$. We have an isomorphism $\theta: \Gal(\goth K^{Hilb \ \infty}/\goth K) \to
\Cl(\w_{\goth K})$.

We formulate the function field analog of Theorem 13.2.6 only for the case when
\medskip
{\bf (*)} There exists only one
point over $\infty\in P^1(\n F_q)$ in the extension $\goth K^{ref}/\n F_q(T)$.
\medskip
In this case the field
$\goth K^{ref \ Hilb \ \infty}$ and the ring $\w_{\goth K^{ref}}$ are naturally defined, and we have
an
isomorphism
$\theta^{ref}: \Gal(\goth K^{ref \ Hilb \ \infty}/\goth K^{ref}) \to \Cl(\w_{\goth K^{ref}})$.

Let $M$ be an uniformizable
t-motive of rank $r$ and dimension $n$ having complete multiplication
by $\w_\goth K$, and $\Phi$ its CM-type. $\Cl(M)$ is defined like $\Cl(A)$ in the number
field case,
it is $\Cl(L,\Phi)$ of 13.1.
\medskip
{\bf Conjecture 13.2.7.} If (*) holds, then $M$ is defined over $\goth K^{ref\ Hilb \
\infty}$, and for any
$\gamma\in\Gal(\goth K^{ref\ Hilb \ \infty}/\goth K^{ref})$
we have $\Cl(\gamma(M))=\det_{Cl} \Phi^{ref}\circ \theta^{ref} (\gamma)^{-1}\Cl(M)$.
\medskip

Now we can formulate the main theorem of this section.
\medskip
{\bf Theorem 13.2.8.} If conjecture 13.2.7 is true for $M$ then it is true for $M'$.
\medskip
{\bf Proof.} It follows immediately from the functional analogs of 13.2.2 -- 13.2.4 that
$$(\goth K,\Phi')^{ref}=(\goth K^{ref},(\Phi^{ref})')\eqno{(13.2.9)}$$
Further,
$$\hbox{det}_{Cl} {\Phi'}^{ref}=(\hbox{det}_{Cl} \Phi^{ref})^{-1}\eqno{(13.2.10)}$$
Really, $\det \Phi^{ref}(x) \cdot \det (\Phi^{ref})'(x) = N_{\goth
K^{ref}/\n F_q(T)}(x)\in\n F_q(T)^\times$, hence gives the trivial class of ideals (we use
here (13.2.9).
Finally, for $\gamma\in\Gal(\goth K^{ref})$ we have
$$(\gamma(M))'=\gamma(M')\eqno{(13.2.11)}$$
The theorem follows immediately from 13.1, 13.2.10, 13.2.11 (recall that $\Cl(M)$ is
$\Cl(L,\Phi)$ of 13.1). $\square$
\medskip
{\bf 13.3. Some explicit formulas.} We give
here an elementary explicit proof of the theorem 12.7 in two simple
cases: $\goth K=\n F_{q^r}(T)$ and $\n F_q(T^{1/r})$. By the way,
since the extension $\n F_{q^r}(T)/\n F_{q}(T)$ is not absolutely
irreducible, formally this case is not covered by the theorem 12.7.
\medskip
{\bf Case $\w_\goth K= \n F_{q^r}[T]$.} Let $\alpha_i$, where $i=0, ... ,r-1$, be inclusions $\goth K \to \p$. For $\omega \in\n F_{q^r}$
 we have $\alpha_i(\omega )=\omega^{q^{i}}$. Let $$0 \le
i_1 < i_2 < ... <
i_n \le r-1\eqno{(13.3.0)}$$ be numbers such that $\Phi=\{\alpha_{i_j}\}$, $j=1, ... , n$.
We consider the following t-motive $M=M(\goth K, \Phi)$. Let $e_1, ... , e_n$ be a basis of $M_{\p[\tau]}$ such that $\goth
m_\omega(e_j)=\omega^{q^{i_j}}e_j$ and the multiplication by $T$ is
defined by formulas
$$Te_1=\theta e_1 + \tau^{i_1-i_n+r}e_n\eqno{(13.3.1)} $$
$$Te_j=\theta e_j + \tau^{i_j-i_{j-1}}e_{j-1}, \ \ j=2, ... , n
\eqno{(13.3.2)} $$

It is easy to check that $M$ has complete multiplication by $\w_\goth K$, and its CM-type is $\Phi$.
\medskip
{\bf Remark.} It is possible to prove that $M(\goth K, \Phi)$ is the only t-motive having these properties; we omit the proof.
\medskip
{\bf Proposition 13.3.3.} For $\w_\goth K= \n F_{q^r}[T]$ we have: $M(\goth K, \Phi)'=M(\goth K, \Phi')$.
\medskip
{\bf Proof.} Elements $\tau^je_k$ for $k=1, ... ,n$, $j=0, ..., i_{k+1}-i_{k}-1$ for $k<n$
and $j=0,
..., i_1-i_n+r-1$ for $k=n$ form a basis of $M_{\p[T]}$. Let us
arrange these elements in
the lexicographic order ($\tau^{j_1}e_{k_1}$ precedes to
$\tau^{j_2}e_{k_2}$ if $k_1 <
k_2$) and make a cyclic shift of them by $i_1$ denoting $e_1$ by $f_{i_1+1}$,
$\tau^{i_2-i_1-1}e_1$ by $f_{i_2}$ etc. until
$\tau^{i_1-i_n+r-1}e_n=f_{i_1}$. Formulas
13.3.1, 13.3.2 become
$$\tau(f_i)=f_{i+1} \hbox{ if } i\not\in \{i_1, ..., i_n\} $$
$$\tau(f_i)=(T-\theta)f_{i+1} \hbox{ if } i\in \{i_1, ..., i_n\} $$
($i \mod r$, i.e. $f_{r+1}=f_1$). Formula 1.10.1 shows that in the
dual basis $f'_*$ we
have
$$\tau(f'_i)=f'_{i+1} \hbox{ if } i\in \{i_1, ..., i_n\} $$
$$\tau(f'_i)=(T-\theta)f'_{i+1} \hbox{ if } i\not\in \{i_1, ..., i_n\} $$
which proves the proposition. $\square$
\medskip
{\bf Case $\w_\goth K=\n F_q[T^{1/r}]$, $(r,q)=1$.} In order to define $M(\goth K, \Phi)$ we need more notations. We denote
$\theta^{1/r}$ and $T^{1/r}$ by $\goth s$ and $S$ respectively, and
let $\zeta_r$ be a primitive $r$-th root of 1. Let $\alpha_i$, $i_1 < i_2 < ... < i_n$ and $\Phi$ be the
same as in the case $\w_\goth K= \n F_{q^r}[T]$. We
have $\alpha_i(S)=\zeta_r^iS$. Further, we consider an overring $\p[S,\tau]$ of $\p[T,\tau]$ ($S$ is in the
center of this ring), and we consider the category of modules
over $\p[S,\tau]$ such that the condition 1.9.2 is changed by a
weakened condition 13.3.4 (here $A_{S,0}\in M_n(\p)$ is defined by the formula $Se_*=A_Se_*$, where $A_S\in M_n(\p)[\tau]$, $A_S=\sum_{i=0}^*A_{S,i}\tau^i$):
$$A_{S,0}^r=\theta E_n + N\eqno{(13.3.4)}$$
Let $\bar M$ be a $\p[S,\tau]$-module such that $\dim \bar
M_{\p[S]}=1$, $f_1$ the only element of a basis of $\bar M_{\p[S]}$ and
$$\tau f_1=(S-\zeta_r^{i_1}\goth s)\cdot ... \cdot
(S-\zeta_r^{i_n}\goth s)f_1$$
By definition, $M=M(\goth K, \Phi)$ is the restriction of scalars from
$\p[S,\tau]$ to $\p[T,\tau]$ of $\bar M$. Like in the case $\w_\goth K= \n F_{q^r}[T]$, it is easy to check that $M$ has complete multiplication by $\w_\goth K$ with CM-type $\Phi$, and it is possible to prove that it is the only t-motive having these properties.
\medskip
{\bf Proposition 13.3.5.} For $\w_\goth K= \n F_q[T^{1/r}]$, $(r,q)=1$ we have: $M(\goth K, \Phi)'=M(\goth K, \Phi')$.
\medskip
{\bf Proof.} For $i=1, ..., r$ we denote $f_i=S^{i-1}f_1$. These $f_*=f_*(\Phi)$
form a basis of $M_{\p[T]}$, and the matrix $Q=Q(f_*,\Phi)$ of
multiplication of $\tau$ in this basis has the following description.
We denote by $\sigma_k(\Phi)$ the elementary symmetric polynomial
$\sigma_k(\zeta_r^{i_1}, ... , \zeta_r^{i_n})$.

The first line of $Q$ is
$$\sigma_n(\Phi)\goth s^n \ \ \ \sigma_{n-1}(\Phi)\goth s^{n-1} \ \ \ ...
\ \ \  \sigma_1(\Phi)\goth s \ \ \ 1 \ \ \ 0 \ \ \ ... \ \ \ 0$$
and its $i$-th line is obtained from the first line by 2 operations:

1. Cyclic shift of elements of the first line by $i-1$ positions to the right;

2. Multiplication of the first $i+n-r$ elements of the obtained line by $T$.

We consider another basis $g_*=g_*(\Phi)$ of $M_{\p[T]}$ obtained by
inversion of order of $f_i$, i.e. $g_i=f_{r+1-i}$. The elements of
$Q(g_*)$ are obtained
by reflection of positions of elements of $Q(f_*)$ respectively the
center of the matrix.

The theorem for the present case follows from the formula
$$Q(f_*,\Phi)Q(g_*,\Phi')^t=(T-\theta)E_r$$
whose proof is an elementary exercise: let $\Phi'=\{j_1, ... ,
j_{r-n}\}$; we apply
equality
$$\sigma_k(x_1, ... x_r)=\sum_l \sigma_l(x_{i_1}, ... ,
x_{i_n})\sigma_{k-l}(x_{j_1}, ...
, x_{j_{r-n}})$$
to $1, \zeta_r, ... , \zeta_r^{r-1}$. $\square$
\medskip
{\bf 13.4. Reduction.} Recall notations of 1.16.  Let $L$ be a finite extension of $\n F_q(\theta)$, $\goth p$ a valuation of $L$ over a valuation $P\ne\infty$ of $\n F_q(\theta)$, and we denote $\iota^{-1}(P)\subset \w$ by $\Cal P$. Let $M$ be a t-motive defined over $L$ having a good ordinary reduction $\tilde M$ at $\goth p$ and such that the dual $M'$ exists. According 1.15.1, the $L$-structure on $M'$ is well-defined. We denote by $M_{\Cal P,0}$ the kernel of the reduction map
$M_\Cal P \to \tilde M_\Cal P$. Condition of ordinarity means that  $M_{\Cal P,0}=(\w/\Cal P)^n$.
\medskip
{\bf Conjecture 13.4.1.} For the above $M$, $M'$ we have:
\medskip
$M_{\Cal P,0}$ and $M'_{\Cal P,0}$ are mutually dual with respect to the pairing of Remarks 4.2, 5.1.6 (recall that conjecturally $M'$ also has good ordinary reduction at $\goth p$).
\medskip
{\bf Proof for a particular case:} $M$ is a Drinfeld module, $\Cal P=T$.
\medskip
(1.9.1) for $M$ has a form
$$Te=\theta e+a_1\tau e+... a_{r-1}\tau^{r-1}e+\tau^r e$$
Condition of good ordinary reduction means $a_i\in L$, $\ord_\goth p(a_i)\ge 0$,
$\ord_\goth p a_1=0$. Let $x\in
M_T$, $y\in M'_T$; we can consider $x$ (resp. $y$) as an element of
$\p$ (resp. $\p^{r-1}$) satisfying some polynomial equation(s). Considering
Newton polygon of these polynomials we get immediately (1) for both
$M$, $M'$. Let $y=(y_1, ... ,y_{r-1})$ be the coordinates of $y$; explicit
formula (5.3.5) for the present case has the form
$$<x,y>_M=\Xi(xy_{r-1}^q+x^qy_1+x^{q^2}y_2+ ... +x^{q^{r-1}}y_{r-1})$$
The same consideration of the Newton polygon of the above polynomials
shows that for $x\in M_{T,0}$, $y\in M'_{T,0}$ we have
$\ord_\goth p x$, $\ord_\goth p y_i \ge 1/(q-1)$. Since $\ord_\goth p \Xi = -1/(q-1)$ we get
that $\ord_\goth p (<x,y>_M)>0$ and hence (because $<x,y>_M\in \n F_q$) we have
$<x,y>_M=0$. Dimensions of $M_{T,0}$, $M'_{T,0}$ are
complementary, hence they are mutually dual. $\square$
\medskip
{\bf Remark 13.4.2.} Analogous explicit proof exists for any standard-3 $M$ of
Section 11.8.
\medskip
\medskip
{\bf References}

\nopagebreak
\medskip
[A] Anderson Greg W. $t$-motives. Duke Math. J. Volume 53, Number 2 (1986), 457-502.
\medskip
[BH] Matthias Bornhofen, Urs Hartl, Pure Anderson Motives and Abelian $\tau$-Sheaves. arXiv:0709.2809
\medskip
[F] Faltings, Gerd, Group schemes with strict $\Cal O$-action. Mosc.
Math. J.  2
(2002),  no. 2, 249--279.
\medskip
[G] Goss, David Basic structures of function
field arithmetic.
\medskip
[H] Urs Hartl, Uniformizing the Stacks of Abelian Sheaves.

http://arxiv.org/abs/math.NT/0409341
\medskip
[L1] Logachev, Anderson t-motives are analogs of abelian varieties with multiplication by imaginary quadratic fields.
http://arxiv.org 0907.4712.
\medskip
[L2] Logachev, Reductions of Hecke correspondences on Anderson varieties. In preparation.
\medskip
[L3] Logachev, Lattice map for Anderson t-motives: first approach. http://arxiv.org/pdf/1109.0679.pdf
\medskip
[P] Richard Pink, Hodge structures over function fields. Universit\"at Mannheim.
Preprint. September 17, 1997.
\medskip
[Sh63] Shimura, Goro On analytic families of polarized abelian
varieties and automorphic functions. Annals of Math., 1 (1963), vol.
78, p. 149 -- 192
\medskip
[Sh71]  Shimura, Goro Introduction to the
arithmetic theory of
automorphic functions.
\medskip
[Sh98]  Shimura, Goro Abelian varieties with
complex multiplication
and modular functions. Princeton Mathematical Series, 46.
\medskip
[Tae] Taelman Lenny, Artin t-motifs. J. Number Theory, 129 (2009), 142 - 157
\medskip
[T] Taguchi, Yuichiro A duality for finite $t$-modules.
J. Math. Sci. Univ. Tokyo 2 (1995), no. 3, 563--588.
\medskip
E-mail: logachev94{\@}gmail.com

\enddocument

and $\mu$ a number such that the $\mu$-dual ${M'}^{\mu}$ exists. Let $m$ be the minimal number such that $N^m=0$.
\medskip
{\bf Conjecture 6B1.} 1. $\mu \ge m$, i.e. $N^\mu=0$.
\medskip
2. ${M'}^{\mu}$ is uniformizable. There exists a canonical perfect $\w$-valued pairing between $L_T(M)$ and $L_T({M'}^{\mu})$.
\medskip
This pairing extends by $\p[[T-\theta]]$-linearity to the $\p[[T-\theta]]$-valued pairing between $L_T(M)\underset{\w}\to{\otimes}\p[[T-\theta]]$, $L_T({M'}^{\mu})\underset{\w}\to{\otimes}\p[[T-\theta]]$. We denote it by $<.,.>$.
\medskip
{\bf Conjecture 6B2.} We have:
$$\goth q({M'}^{\mu})=\{x\in L_T({M'}^{\mu})\underset{\w}\to{\otimes}\p[[T-\theta]]  \hbox{ such that }  \forall y\in \goth q(M) $$ $$\hbox{ we have } <x,y>\in (T-\theta)^\mu\p[[T-\theta]]\} $$

Let us consider uniformizable $M$ having $N=0$ (i.e. $m=1$) such that $M'$ --- the dual of $M$ --- exists and has $N'=0$.
\medskip
{\bf Corollary 6B3.} For this case Conjecture 6B2 is true, i.e.
$$\goth q({M'})=\{x\in L_T({M'})\underset{\w}\to{\otimes}\p[[T-\theta]]  \hbox{ such that }  \forall y\in \goth q(M) $$ $$\hbox{ we have } <x,y>\in (T-\theta)\p[[T-\theta]]\} $$
\medskip
This follows immediately from Theorem 5 and equivalence of Definition 2.3 and Property 2.4 (we consider the reduction of the above pairing $<.,.>$ modulo the maximal ideal $ (T-\theta) \p[[T-\theta]]$ of $\p[[T-\theta]]$; this reduction coincides with the pairing of Section 5).

\newpage
{\bf 1.16.4. $\Cal P$-rank of $M$ and of $M'$.} By analogy with the number field case, the $\Cal P$-rank of $M$ is the dimension of the $\Cal P$-torsion points of $E(M)$ over $\w/\Cal P$; it varies from $r-n$ (ordinary $M$) to 0 (completely suresingular $M$).
\medskip
{\bf Conjecture 1.16.5.} $M$ is ordinary $\iff M'$ is ordinary.
\medskip
For standard-3 t-motives  apparently this can be shown by explicit calculations. See 13.4.1 for the proof for the case of Drinfeld modules.
\medskip
{\bf Question 1.16.6.} What are possible values of $\Cal P$-rank of $M$, $\Cal P$-rank of $M'$? Is it true that the pair of numbers ($\Cal P$-rank of $M$, $\Cal P$-rank of $M'$) characterizes completely (in some meaning) the type of $M$? Particularly, whether the $\Cal P$-ranks of other tensor operations of $M$ (exterior powers etc.) are defined completely by ($\Cal P$-rank of $M$, $\Cal P$-rank of $M'$), or not?
\medskip
{\bf Example 1.16.7.} Let us consider $M$ defined by (8.2.2), entries of $A$ belong to $\bar \n F_q$, and let $\Cal P=T$. According Lang's Theorem, $T$-rank of $M$ is $n \iff \det A\ne 0$. An explicit calculation for the case $n=2$ shows that if the entries of $A$ belong to $\n F_q$ then the $T$-rank of $M$ is equal to the quantity of the non-zero eigenvalues of $A$ (if the entries of $A$ do not belong to $\n F_q$ then the formula for the $T$-rank of $M$ is more complicated). The dual $M'$ is defined by the same (8.2.2) with $A$ replaced by $-A^t$; this means that  if the entries of $A$ belong to $\n F_{q}$ then the $T$-rank of $M$ is equal to the $T$-rank of $M'$. It is easy to check that the same is true if the entries of $A$ belong to $\n F_{q^2}$ but not to a larger field:  for example, if $\gamma\in \bar \n F_q-\n F_{q^2} $ then for $A=\left(\matrix \gamma & 1 \\ -\gamma^{q+1}& -\gamma^{q}\endmatrix \right)$ we have $T$-rank of $M$ is 0, $T$-rank of $M'$ is 1.
\medskip
{\bf Example 1.16.8.} For pure non-ordinary standard-3 $M$ having $r=5$, $n=2$ we get easily that the case $T$-rank of $M$ is 2, $T$-rank of $M'$ is 0 cannot be realized; all other possible cases can be realised.
\medskip
{\bf Remark 13.4.3: Case $N\ne0$.} I do not know a definition of ordinarity for $M$ in finite characteristic for the case $N\ne0$; we can expect that these are those $M$ whose $\Cal P$-rank is the maximal possible in some family of these $M$. A version of Conjecture 13.4.1 should hold for this case; instead of the dual $M'$ we should consider the $m$-th dual ${M'}^m$ where $m$ satisfies $M^m=0$ (or $m$ is the minimal number with this property?)
\medskip
{\bf Example 13.4.3.1.} Let $M$ be given by the following modified formula (8.2.2):
$$Te_* =( \theta+N) e_* + A \tau e_* + \tau^2 e_*\eqno{(13.4.3.2)}$$
where $n=2$, $N=\left(\matrix 0&1\\0&0 \endmatrix\right)$, $A= \left(\matrix a_{11}&a_{12} \\a_{21}&a_{22} \endmatrix\right)$. For the case of finite characteristic (i.e. $a_{ij}\in \bar\n F_q$) and $\Cal P=T$ we have $M$ is ordinary iff $a_{21}\ne0$, in this case the $\Cal P$-rank of $M$ is 3. We have $m=2$, ${M'}^2$ is given by the formula
$$Te'_*=( \theta+N') e'_* + A' \tau e'_*$$
where $N'=\left(\matrix 0&-1&0&0&0&0\\0&0&0&0&0&0\\0&0&0&1&0&0\\0&0&0&0&0&0
\\0&0&0&0&0&1\\0&0&0&0&0&0\endmatrix\right)$,  $A'=\left(\matrix 0&0&0&0&1&0 \\0&0&1&0&0&0\\0&0&0&0&0&0\\0&1&-a_{11}&0&-a_{21}&0\\0&0&0&0&0&0
\\1&0&-a_{12}&0&-a_{22}&0 \endmatrix\right)$.

An explicit calculation shows that the $\Cal P$-rank of $M$ is $\le 1$ and it is 1 iff $a_{21}\ne0$. So, an analog of 13.4.1 for $M$ is the following:
\medskip
{\bf Conjecture 13.4.3.3.} Let $M$ be given by 13.4.3.2 in the generic characteristic, and let $M$ have good ordinary reduction at $\Cal P=T$.  We have: $M_{T,0}$ and $({M'}^2)_{T,0}$ are mutually dual with respect to the pairing.
\medskip
This conjecture can be easily checked by explicit calculation.

\enddocument